# VOLODYMYR SHCHEDRYK

# ARITHMETIC OF MATRICES OVER RINGS

$$G_\Phi = \{H \in GL_n(R) \mid \exists K \in GL_n(R) : H\Phi = \Phi K\}$$

$$P_A A Q_A = \mathrm{diag}(\varepsilon_1, \varepsilon_2, \ldots, \varepsilon_n) = H$$

$$\begin{Vmatrix} h_{11} & h_{12} \\ \frac{\varphi_2}{\varphi_1} h_{21} & h_{22} \\ \vdots & \vdots \\ \frac{\varphi_n}{\varphi_1} h_{n1} & \frac{\varphi_n}{\varphi_2} h_{n2} \end{Vmatrix}$$

NATIONAL ACADEMY
OF SCIENCES OF UKRAINE
PIDSTRYHACH INSTITUTE FOR APPLIED
PROBLEMS OF MECHANICS
AND MATHEMATICS
OF THE NAS OF UKRAINE

НАЦІОНАЛЬНА
АКАДЕМІЯ НАУК УКРАЇНИ
ІНСТИТУТ ПРИКЛАДНИХ ПРОБЛЕМ
МЕХАНІКИ І МАТЕМАТИКИ
ім. Я.С. ПІДСТРИГАЧА НАН УКРАЇНИ

ВОЛОДИМИР
ЩЕДРИК

# АРИФМЕТИКА МАТРИЦЬ НАД КІЛЬЦЯМИ

VOLODYMYR
SHCHEDRYK

# ARITHMETIC
# OF MATRICES
# OVER RINGS

**Reviewers:**

*V.M. PETRYCHKOVYCH*, Doctor of Physical and Mathematical Sciences, Head of Algebra Department of Pidstryhach Institute for Applied Problems of Mechanics and Mathematics of National Academy of Sciences of Ukraine, Professor

*V.M. BONDARENKO*, Doctor of Physical and Mathematical Sciences, Leading Researcher of Algebra and Topology Department of Institute of Mathematics of National Academy of Sciences of Ukraine, Professor

*Approved to press by the Scientific Council of Pidstryhach Institute for Applied Problems of Mechanics and Mathematics of the NAS of Ukraine (September 10, 2020, Protocol No. 6)*

**The publication was funded within the framework of the Targeted Complex Program of the NAS of Ukraine "Scientific Bases of Functioning and Providing for Conditions for the Development of the Scientific and Publishing Complex of the NAS of Ukraine"**

**Shchedryk V.P.**

S53      Arithmetic of matrices over rings / Volodymyr Shchedryk; Pidstryhach Institute for Applied Problems of Mechanics and Mathematics of the NAS of Ukraine. — Kyiv: Akademperiodyka, 2021. — 278 p.

     ISBN 978-966-360-430-5

     The book is devoted to investigation of arithmetic of the matrix rings over certain classes of commutative finitely generated principal ideals domains. We mainly concentrate on constructing of the matrix factorization theory. We reveal a close relationship between the matrix factorization and specific properties of subgroups of the complete linear group and the special normal form of matrices with respect to unilateral equivalence. The properties of matrices over rings of stable range 1.5 are thoroughly studied.

     The book is intended for experts in the ring theory and linear algebra, senior and post-graduate students.

         **UDC 517.98**

# CONTENTS

CHAPTER **4**

## DIVISIBILITY AND ASSOCIATIVITY OF MATRICES

CHAPTER **5**

## FACTORIZATION OF MATRICES

CHAPTER **6**

## INVARIANTS OF PRIMITIVE MATRICES WITH RESPECT TO ZELISKO GROUP ACTION

CHAPTER **7**

## ONE-SIDED EQUIVALENCE OF MATRICES

# PREFACE

The object of study of the present monograph are matrices over commutative elementary divisor domains at view points of their decomposability into multipliers and related problems establishing relationships between certain subgroups and subsets of a general linear group and reduction of matrices by one-sided transformations to the canonical forms.

I. Kaplansky [1] introduced the concept of the *elementary divisor rings* $R$ as a ring over which any matrix has the canonical diagonal reduction property. That is, for each matrix $A$ over $R$ there exist invertible matrices $P$ and $Q$ such that

$$PAQ = \text{diag}\,(d_1, ..., d_r, 0, ..., 0),$$

in which

$$Rd_{i+1}R \subseteq d_iR \cap Rd_i, \quad (i = 1, ..., r - 1). \qquad (0.1)$$

If $R$ is a commutative ring then the condition $(0.1)$ is equivalent to $d_i|d_{i+1}$, (divides) $i = 1, ..., \mathrm{r} - 1$. The first results of this type concerning integer matrices were obtained by H. Smith [2]. He established that, by elementary transformations of rows and columns, every matrix reduced to a diagonal form, and each diagonal element is a divisor of the following one. Therefore, the diagonal form of the matrix with the condition of divisibility of the diagonal elements is called the Smith form. L. Dixon [3], J. Wedderburn [4, 5], B. van der Waerden [6], N. Jacobson [7] extended Smith's result to different classes of commutative and noncommutative Euclidean rings and commutative principal ideal rings without zero divisors. The ring of analytic functions in the complex plane [4], the ring of continuous real-valued functions on a completely regular Hausdorf space [10], the ring of all integer algebraic numbers [12], the ring of formal power series over the field of rational numbers with free integer term [11], the adequate rings [8] are examples of elementary divisor rings.

In order to study of elementary divisor rings, I. Kaplansky introduced the right Hermitian ring concept. A ring $R$ is a right Hermitian ring provided that, each $1 \times 2$ matrix has he property of a canonical diagonal reduction. I. Kaplansky showed that a right Hermitian ring is a right Bezout ring (finitely generated right principal ideal ring). In the case of commutative domains, these rings coincide [13].

Studying commutative Hermitian rings and elementary divisor rings with zero divisors, L. Gillman, M. Henriksen [10, 11], gave an example of a commutative Bezout ring, which is not a Hermitian ring, and a Hermitian ring, which is not an elementary divisor ring. They raised the question of coincidence of these rings in the case of commutative domains, which is now known as the *problem of elementary divisor rings*. It is no exaggeration to say that this is one of the most pressing problem of modern algebra. Among the large number of articles devoted to this subject, we highlight the studies of P. Kohn [14], M. Larsen, W. Lewis, T. Shores [15], W. Mc Govern [16], who solved it by the methods of classical ring theory and by module theory methods.

Promising studies of this problem are based on the concept of the stable range of ring. This notion was first introduced [9] so as to study stabilization in algebraic K-theory. B. Zabavskii obtained a number of structural theorems that deeply characterize Bezout rings of finite stable range [17–19]. In particular, he showed that the elementary divisor rings have a stable range not more than 2 [17]. Using the notion of stable range, the widely studied commutative clean rings [20], exchange rings [21], neat rings [22] are elegantly introduced. In-depth studies of this concept have led to the introduction of an idempotent stable range [20, 22] and a unit stable range [23].

The interest to the elementary divisor rings, no doubt, is due not only to the fact that the majority of classical rings have the property of canonical diagonal reduction. Their perfect internal structure makes it possible to obtain deep results in related fields, in particular in the theory of modules and in algebraic K-theory. I. Kaplansky [1], showed that over an elementary divisor ring arbitrary finitely represented module can be decomposed into a direct sum of cyclic modules. In the case of commutative rings the reverse statement is true, namely: if finitely represented module over a ring can be decomposed into a direct sum of cyclic modules, then this ring is an elementary divisor ring [15].

The great potential of these rings was manifested in the possibility of constructing a unified theory of matrix factorization over them. Note that the problems of matrix factorization began to be interested in the second half of the XIX century. In particular, A. Cayley [24], M. Sylvester [25, 26], F. Frobenius [27] considered the problem of numerical matrix representation as the product of identical matrices.

The most thoroughly factorization problems were investigated in polynomial matrix rings. It was stimulated by the fact that according to generali-zed Bezout Theorem [28], separated of linear monic multiplier, i.e., a multiplier of the form $Ix - B$, where $I$ is a identity matrix, from a polynomial matrix is equivalent to finding the root $B$ of the corresponding one-sided matrix equation. M. Ingraham [28, 29], W. Roth [30-32], P. Lancaster [33-35], H. Langer [51], H. Wimmer [37], J. Dennis, J. Traub, R. Weber [38], I. Gohberg, L. Rodman [39], A. Markus, I. Mereutsa [40], A. Malyshev [41] etc. (see papers from the bibliography) were engaged in such tasks.

Such studies were popular in the second half of the twentieth century, when research teams won the lead in solving the problem of separating the regular factor from the polynomial matrix. P. Kazimirsky was the first to propose a constructive method for its solution [55-57]. It should be noted that the impulse for such research was a result of the article of the academician Ja. Lopatinsky (teacher of P. Kazimirsky), which considered systems of differential equations with constant coefficients and received a purely algebraic result [58], which was concerned with separating the regular factor from the polynomial matrix.

Significant results concerning the matrix factorization were obtained by V. Zelisko [63, 64] and V. Petrychkovych [65-71]. In particular, the group introduced by V. Zelisko is the cornerstone of the matrix factorization theory constructed in this monograph. V. Petrychkovych, based on the concepts introduced by him together with P. Kazimirsky, of the semi-scalar and generalized equivalence and a certain triangular form of matrices with invariant factors on the main diagonal, significantly expanded the class of absolutely decomposable matrix polynomials, that is, those polynomials that decompose into the product of linear multipliers. He proposed necessary and sufficient conditions of uniqueness, up to associativity, of divisors with a given Smith form for matrices over adequate domains. He introduced the concept of matrix factorization which parallel the factorization of its Smith form, and described the classes of matrices that have such a property.

Studies of matrix factorization were not limited to the consideration of polynomial matrices. Matrix factorization were also studied over skew fields, universal algebras [59, 60], and Dedekind domains [61].

A special highlight of the paper by Z. Borevich [62], where it was proposed to classify the divisors of matrices over commutative principal ideal domains by their Smith forms and to describe them up to associativity. In the monograph, the approach proposed by Z. Borevich is implemented for matrices over elementary divisor domains, over which such a problem has not yet been considered. Such research has become possible through the use of a fundamentally new approach based on studying the set of invertible matrices $\mathbf{L}(E, \Phi)$, which generates all left divisors with given Smith form $\Phi$ and the studying of the Zelisko group $\mathbf{G}_\Phi$, which is responsible for the associativity of matrices.

We will briefly present the content of the book with an overview of its sections. In accordance with the logical development of the subject, the mate-rials divided into seven parts. The book is written so that the reader who is acquainted with university level linear algebra, can comprehend the material presented, without referring to other sources. For this purpose, in the first two sections, in addition to the original results, are already known. They relate to the properties of Bezout ring elements, the complementary of primitive rows to invertible matrices, normal forms of matrices with respect to one- and two-sided equivalences. Particular attention in the first section is paid to the general properties of the largest common divisors and the least common multiple of matrices.

In the second section, the properties of the transforming matrices and the Zelisko group, which is their generator, are studied. The concept of a ring of stable range 1.5 is introduced and it is proved that the general linear group of this ring is presented as the product of the Zelisko group, the groups of lower and upper unitriangular matrices. It is also shown that over the commutative Bezot domain $R$ the second order matrix rings has a stable range 1.5 if and only if $R$ has a similar stable range.

The third section is technical in nature. It investigates the properties of minors and invariant factors of block-triangular matrices and their diagonal blocks. These results were used in the presentation of subsequent sections.

In the fourth section the concept of the generating set is introduced and on its basis the criteria of divisibility and associativity of matrices are formulated. The structure of its elements is investigated and the relationship of this set with the Zelisko group is established.

The fifth section is devoted to the study of factorization of matrices. In particular, the conditions of uniqueness (up to associativity) of divisors with a given Smith form are indicated. The notion of the Kazimirskii set is introduced and the conditions under which this set generates all unassociated divisors with a given Smith form are indicated. Special attention is paid to the factorization of polynomial matrices. The result was the construction of effective methods for finding the roots of one-sided matrix equations and finding the Jordan normal form. The two final subsections of this section are devoted to finding Smith forms and the transforming matrices of the greatest common divisors and the least common multiple of matrices.

The sixth chapter investigates the act of the Zelisko group on primitive matrices. In this case, the invariants with respect to such action are indicated.

The final section is devoted to the study of one-sided equivalence of matrices with a fixed Smith form. It is shown that the search for divisors that have a given Smith form is parallelized and is reduced to the description of certain matrices with given $\Phi$-skeletons. The consequence of these results was a description of all unassociated matrices with a standard $\Phi$-skeleton, and all irreducible matrix divisors. The obtained results are applied to the description of all nilpotent and idempotent matrices of the third order.

# BEZOUT RINGS
# AND MATRICES OVER THEM

*The section describes the basic properties of elements of com-mutative Bezout domains, specifies the normal form of matri-ces for one-side transformations from a general linear group (Hermite form). Particular attention is paid to the Bezout equality.*

## 1.1. Properties of Bezout rings elements

The notion of greatest common divisor (g.c.d.) is one of the foundations in mathematics. For non-commutative rings are consider (if they exist) right and left greatest common di-visors. It should be noted that there are rings in which the concept of the greatest common divisor is incorrect. Examples of such rings are infinite number of variables polynomial rings over fields. Generally speaking, in these rings in the set of common divisors, there is no such polynomial that is divi-sible by the remainder. In addition, it does not follow that the greatest common divisor is defined uniquely up to associ-ativity. We will denote the greatest common divisor of two elements $a, b$ as $(a, b)$. Particularly "appreciated" if $(a, b)$ is linearly expressed through them. That is

$$(a, \ b) = au + bv. \tag{1.1}$$

This equation is called **Bezout's identity** or **Bezout's equality.**

Note that the existence of g.c.d. does not guarantee its representation as (1.1). In particular, in the ring of polyno-mials of two variables over a field $F$ the elements $x, y$ are relatively prime, that is $(x, y) = 1$, but in $F[x, y]$ there are no $u, v$, such that $xu + yv = 1$.

If the condition (1.1) are satisfied for all elements of a ring $R$, than every ideal generated by two elements is principal. So every finitely generated ideal $R$ also principal.

**Definition 1.1.** A ring $R$ is called Bezout if every finitely generated ideal $I$ of $R$ is principal.

Note that a non-commutative ring is called Bezout if each left and right finitely generated ideal is principal.

Examples of commutative Bezout rings are the integer ring, a ring of polynomials from one variable over fields, an integer analytic ring.

Obviously, the commutative principal ideal rings are Bezout rings. The inverse statement is incorrect. For example, the adequate rings discussed below are Bezout rings, however, they are not principal ideal rings. It follows that the elements of Bezout rings cannot be represented as the product of the irreducible elements. That is, Bezout rings are, in general, non-factorial. Matrix rings over commutative Bezout rings are non-commutative Bezout rings.

All rings considered below unless specifically stated, are commutative Bezout domains, that is, a commutative Bezout rings without zero divisors with $1 \neq 0$.

**Property 1.1.** *If* $(a, b) = 1$ *and* $a|bc$*, then* $a|c$.

**Proof.** There are $u, v$ such that

$$au + bv = 1 \Rightarrow acu + bcv = c.$$

Since $a|bc$, then $bc = am$. So,

$$acu + amv = c \Rightarrow c = a(cu + mv).$$

That is $a|c$. $\qquad\qquad\square$

**Property 1.2.** *If* $(a, b) = 1$ *and* $(a, c) = 1$*, then* $(a, bc) = 1$.

**Proof.** There are $u, v, m, n$ such that

$$au + bv = 1,$$
$$am + cn = 1.$$

Multiplying the right and the left sides of these equations, we get

$$au + bvam + bcvn = 1 \Rightarrow a(u + bvm) + bc(vn) = 1.$$

It follows that $(a, bc) = 1$. $\qquad\qquad\square$

**Property 1.3.** *If* $(a, b) = 1$*, then* $(a, bc) = (a, c)$.

**Proof.** Let $(a, bc) = \alpha$. Whereas $\alpha|bc$, then $(\alpha, bc) = \alpha$. Considering that $a = \alpha a_1$ and

$$1 = (a, b) = (\alpha a_1, b),$$

we get $(\alpha, b) = 1$. Taking into account Property 1, we receive $\alpha|c$. So $(a, bc)|(a, c)$. Since $(a, c)|(a, bc)$, we get $(a, bc) = (a, c)$. $\qquad\qquad\square$

**Property 1.4.** *If* $(a, b) = 1$*,* $a|c$*, and* $b|c$*, then* $ab|c$.

**Proof.** Since $c = a\mu$ and $c = b\nu$, there exist $u, v$ such that

$$au + bv = 1 \Rightarrow (a\mu)u + b\mu v = \mu \Rightarrow$$
$$\Rightarrow cu + b\mu v = \mu \Rightarrow b\nu u + b\mu v = \mu.$$

Then $b(\nu u + \mu v) = \mu$. Thus,

$$c = a\mu = ab(\nu u + \mu v). \qquad \square$$

**Property 1.5.** *If* $(a, mn) = (b, \; mn)$, *then* $(a, m) = (b, m)$.
**Proof**. Since
$$(a, m)|(a, mn) = (b, mn),$$

then $(a, m)|b$, and $(a, m)|m$. Hence $(a, m)|(b, m)$. By a similar scheme we show
that $(b, m)|(a, m)$. Therefore, $(a, m) = (b, m)$. $\qquad \square$

**Property 1.6.** *Equality is fulfilled*

$$\left(\frac{a}{(a, c)}, \frac{b}{(b, c)}\right) = \frac{(a, b)}{(a, b, c)},$$

*where* $ab$ *and* $c$ *do not equal to zero simultaneously.*

**Proof.** We have

$$\left(\frac{a}{(a, c)}, \frac{b}{(b, c)}\right) = \frac{(a, b)}{(a, b, c)} \left(\frac{(a^2, ab, ac)}{(a^2, ab, c(a, b))}, \frac{(b^2, ab, bc)}{(a^2, ab, c(a, b))}\right).$$

Since

$$\frac{a}{(a, b)} \frac{\left(a^2, ab, c(a, b)\right)}{(a^2, ab, ac)} = \frac{\left(a^3, a^2 b, ac(a, b)\right)}{(a^2(a, b), ab(a, b), (a, b)ac)} \in R,$$

then

$$\frac{(a^2, ab, ac)}{(a^2, ab, c(a, b))} \left|\frac{a}{(a, b)}\right..$$

Similarly, we show that

$$\frac{(b^2, ab, bc)}{(b^2, ab, c(a, b))} \left|\frac{b}{(a, b)}\right..$$

Taking into account that

$$\left(\frac{a}{(a, b)}, \frac{b}{(a, b)}\right) = 1,$$

we get

$$\left(\frac{(a^2, ab, ac)}{(a^2, ab, c(a, b))}, \frac{(b^2, ab, bc)}{(b^2, ab, c(a, b))}\right) = 1.$$

It means that

$$\left(\frac{a}{(a, c)}, \frac{b}{(b, c)}\right) = \frac{(a, b)}{(a, b, c)}. \qquad \square$$

**Property 1.7.** *Equality is fulfilled*

$$\left(\frac{a}{(a, b)}, \frac{a}{(a, c)}\right) = \frac{a}{(a, [b, c])},$$

*where* $a$ *and* $bc$ *do not equal to zero simultaneously.*

**Proof.** We have

$$\left(\frac{a}{(a,b)}, \frac{a}{(a,c)}\right) = \left(\frac{ac}{(ac,bc)}, \frac{ab}{(ab,bc)}\right).$$

Using Property 1.6, we get

$$\left(\frac{ac}{(ac,bc)}, \frac{ab}{(ab,bc)}\right) = \frac{(ab,ac)}{(ab,ac,bc)} = \frac{a(b,c)}{(a(b,c),bc)} =$$
$$= \frac{a(b,c)}{(b,c)\,(a,[b,c])} = \frac{a}{(a,[b,c])}.$$  □

**Property 1.8.** *Equality is fulfilled*

$$\left[\frac{a}{(a,b)}, \frac{a}{(a,c)}\right] = \frac{a}{(a,(b,c))},$$

*where a and bc do not equal to zero simultaneously.*

**Proof.** Consider the product

$$\frac{a}{(a,[b,c])}\frac{a}{(a,(b,c))} = \frac{a^2}{(a^2, a\,([b,c],(b,c))\,,bc)} =$$
$$= \frac{a^2}{(a^2, a(b,c), bc)} = \frac{a^2}{(a,b)(a,c)}.$$

Then

$$\left[\frac{a}{(a,b)}, \frac{a}{(a,c)}\right] = \frac{\frac{a^2}{(a,b)(a,c)}}{\left(\frac{a}{(a,b)}, \frac{a}{(a,c)}\right)}.$$

By virtue of Property 1.7, we receive

$$\frac{\frac{a^2}{(a,b)(a,c)}}{\left(\frac{a}{(a,b)}, \frac{a}{(a,c)}\right)} = \frac{\frac{a}{(a,[b,c])}\,\frac{a}{(a,(b,c))}}{\frac{a}{(a,[b,c])}} = \frac{a}{(a,(b,c))}.$$  □

**Property 1.9.** *The equality is fulfilled*

$$\left[\frac{a}{(a,c)}, \frac{b}{(b,c)}\right] = \frac{[a,b]}{([a,b],c)},$$

*where ab and c do not equal to zero simultaneously.*

**Proof.** The equality

$$\left[\frac{a}{(a,c)}, \frac{b}{(b,c)}\right] = \left[\frac{ab}{(ab,cb)}, \frac{ab}{(ab,ac)}\right].$$

is fulfilled. According to Property 1.8, we obtain

$$\left[\frac{ab}{(ab,cb)}, \frac{ab}{(ab,ac)}\right] = \frac{ab}{(ab,(cb,ac))} =$$

$$= \frac{ab}{(ab,c(b,a))} = \frac{ab}{(a,b)\,([a,b],c)} = \frac{[a,b]}{([a,b],c)},$$

as required. $\qquad\qquad\square$

**Property 1.10.** *If $(a_1,...,a_n) = 1$, $n \geq 2$, and $\varepsilon_1|\varepsilon_2|...|\varepsilon_k$, $\varepsilon_k \neq 0$, $1 \leq k < n$, then*

$$(a_1\varepsilon_k, a_2\varepsilon_{k-1}, ..., a_k\varepsilon_1, a_{k+1}, ..., a_n) = (\varepsilon_k, a_2\varepsilon_{k-1}, ..., a_k\varepsilon_1, a_{k+1}, ..., a_n).$$

**Proof.** Set

$$(a_1\varepsilon_k, a_2\varepsilon_{k-1}, ..., a_k\varepsilon_1, a_{k+1}, ..., a_n) = \delta_k.$$

Since

$$(\varepsilon_k, a_2\varepsilon_{k-1}, ..., a_k\varepsilon_1, a_{k+1}, ..., a_n)|\delta_k,$$

then to prove our assertion it is suffices to show that $\delta_k|\varepsilon_k$.

If $k = 1$, then using Property 1.3, we get

$$\delta_1 = (a_1\varepsilon_1, a_2, ..., a_n) = (a_1\varepsilon_1, (a_2, ..., a_n)) = (\varepsilon_1, (a_2, ..., a_n)),$$

i.e., $\delta_1|\varepsilon_1$. Therefore, our assertion is   correct if $k = 1$.

Let $k \geq 2$ and suppose that this statement is true for all $m < k$. Then

$$\delta_k = (a_1\varepsilon_k, a_2\varepsilon_{k-1}, ..., a_k\varepsilon_1, a_{k+1}, ..., a_n) =$$

$$= \left(\varepsilon_1\left(a_1\frac{\varepsilon_k}{\varepsilon_1}, ..., a_{k-1}\frac{\varepsilon_2}{\varepsilon_1}, a_k\right), a_{k+1}, ..., a_n\right).$$

Since

$$\frac{\varepsilon_2}{\varepsilon_1}\left|\frac{\varepsilon_3}{\varepsilon_1}\right|...\left|\frac{\varepsilon_k}{\varepsilon_1},\right.$$

by induction assumption, we have

$$\left(a_1\frac{\varepsilon_k}{\varepsilon_1}, ..., a_{k-1}\frac{\varepsilon_2}{\varepsilon_1}, a_k, a_{k+1}, ..., a_n\right) = d_1\left|\frac{\varepsilon_k}{\varepsilon_1}\right..$$

It means that

$$\delta_k = d_1\left(\varepsilon_1\frac{\left(a_1\frac{\varepsilon_k}{\varepsilon_1}, ..., a_{k-1}\frac{\varepsilon_2}{\varepsilon_1}, a_k\right)}{d_1}, \left(\frac{a_{k+1}}{d_1}, ..., \frac{a_n}{d_1}\right)\right) =$$

$$= d_1\left(\varepsilon_1, \left(\frac{a_{k+1}}{d_1}, ..., \frac{a_n}{d_1}\right)\right) = d_1 d_2,$$

**15**

where $d_2|\varepsilon_1$. Hence,

$$\delta_k = d_1 d_2 \left| \frac{\varepsilon_k}{\varepsilon_1} \varepsilon_1 = \varepsilon_k. \right. \qquad \square$$

**Property 1.11.** *Let $a, b, a_1, b_1, \varphi$ be nonzero elements in $R$ and*

$$(a, b, \varphi) = (a_1, b_1, \varphi) = 1, \varphi| \det \left\| \begin{matrix} a & a_1 \\ b & b_1 \end{matrix} \right\|.$$

*Then for all $r \in R$ the equality*

$$(ar + b, \varphi) = (a_1 r + b_1, \varphi)$$

*is executed.*

**Proof.** Set $(a, \varphi) = \alpha$. Since $ab_1 - a_1 b = \varphi t$, then $\alpha|a_1 b$. Noting that $(\alpha, b) = 1$, we get $\alpha|a_1$. Thus, $\alpha|(a_1, \varphi) = \alpha_1$. Similarly we show that $\alpha_1|\alpha$. It means that

$$(a, \varphi) = (a_1, \varphi) = \alpha. \qquad (1.2)$$

Let $r \in R$. Consider

$$(a_1(ar + b), \varphi) = (a_1 ar + a_1 b, \varphi) = (a_1 ar + a_1 b + \varphi t, \varphi) = d.$$

Since $ab_1 = a_1 b + \varphi t$, then

$$d = (a_1 ar + ab_1, \varphi) = (a(a_1 r + b_1), \varphi),$$

i.e.,

$$(a_1(ar + b), \varphi) = (a(a_1 r + b_1), \varphi).$$

Dividing both parts of this equality by $\alpha$, we obtain

$$\left( \frac{a_1}{\alpha}(ar + b), \frac{\varphi}{\alpha} \right) = \left( \frac{a}{\alpha}(a_1 r + b_1), \frac{\varphi}{\alpha} \right).$$

From equality (1.2), it follows that

$$\left( \frac{a}{\alpha}, \frac{\varphi}{\alpha} \right) = \left( \frac{a_1}{\alpha}, \frac{\varphi}{\alpha} \right) = 1.$$

Consequently,

$$\left( ar + b, \frac{\varphi}{\alpha} \right) = \left( a_1 r + b_1, \frac{\varphi}{\alpha} \right) = \delta.$$

Since $\alpha|a, \alpha|a_1$, and

$$(\alpha, b) = (\alpha, b_1) = 1,$$

then for all $r$ from $R$ the equality

$$(ar + b, \alpha) = (a_1 r + b_1, \alpha) = 1$$

holds. This yields

$$\delta = \left( ar + b, \frac{\varphi}{\alpha} \right) = \left( ar + b, \frac{\varphi}{\alpha} \alpha \right) = (ar + b, \varphi).$$

Similarly, we show that $\delta = (a_1 r + b_1, \varphi)$. $\qquad \square$

**Property 1.12.** *Let*

$$(a_1, b_1, \varphi) = 1,$$

$$\varphi \mid \det \left\| \begin{matrix} a & a_1 \\ b & b_1 \end{matrix} \right\|, \tag{1.3}$$

*where* $aba_1b_1 \neq 0$. *If*

$$(a_1 r + b_1, \varphi) = \alpha,$$

*then*

$$(ar + b, \varphi) = \alpha \left( a, b, \frac{\varphi}{\alpha} \right).$$

**Proof.** Set $(a, b, \varphi) = \delta$, then

$$\left( \frac{a}{\delta}, \frac{b}{\delta}, \frac{\varphi}{\delta} \right) = 1.$$

Since

$$\varphi \left| \left| \begin{matrix} a & a_1 \\ b & b_1 \end{matrix} \right| = \delta \left| \begin{matrix} \frac{a}{\delta} & a_1 \\ \frac{b}{\delta} & b_1 \end{matrix} \right|,$$

then

$$\frac{\varphi}{\delta} \left| \left| \begin{matrix} \frac{a}{\delta} & a_1 \\ \frac{b}{\delta} & b_1 \end{matrix} \right|.$$

In view of (1.3), we have

$$\left( a_1, b_1, \frac{\varphi}{\delta} \right) = 1.$$

According to Property 1.11, for all $r \in R$ the condition

$$\left( a_1 r + b_1, \frac{\varphi}{\delta} \right) = \left( \frac{a}{\delta} r + \frac{b}{\delta}, \frac{\varphi}{\delta} \right).$$

holds. Multiplying this equality by $\delta$, we get

$$(ar + b, \varphi) = (\delta(a_1 r + b_1), \varphi) = \alpha \left( \delta \frac{a_1 r + b_1}{\alpha}, \frac{\varphi}{\alpha} \right) = \alpha \left( \delta, \frac{\varphi}{\alpha} \right).$$

Consequently,

$$(ar + b, \varphi) = \alpha \left( a, b, \varphi, \frac{\varphi}{\alpha} \right) = \alpha \left( a, b, \frac{\varphi}{\alpha} \right),$$

which had to be proved. $\square$

**Property 1.13.** *Let*

$$(a_1 r + b_1, \varphi) = 1, \varphi \mid \det \left\| \begin{matrix} a & a_1 \\ b & b_1 \end{matrix} \right\|,$$

*where* $ba_1b_1 \neq 0$. *Then*

$$(ar + b, \varphi) = (a, b, \varphi).$$

**17**

**Proof.** If $a = 0$, there is nothing to prove. So let it $a \neq 0$. From the equality $(a_1 r + b_1, \varphi) = 1$ it follows that $(a_1, b_1, \varphi) = 1$. By Property 1.12,

$$(ar + b, \varphi) = (a, b, \varphi).$$

Which must to be proved. $\qquad\square$

**Property 1.14.** *Let* $(m, a_i) = (n, a_i)$, $i = 1, ..., k$. *Then*

$$(m, [a_1, a_2, ..., a_k]) = (n, [a_1, a_2, ..., a_k]),$$

*i.e.,*

$$\left( m, \frac{a_1 a_2 \ldots a_k}{(a_2 \ldots a_k, a_1 a_3 \ldots a_k, ..., a_1 a_2 \ldots a_{k-1})} \right) =$$

$$= \left( n, \frac{a_1 a_2 \ldots a_k}{(a_2 \ldots a_k, a_1 a_3 \ldots a_k, ..., a_1 a_2 \ldots a_{k-1})} \right).$$

**Proof.** First we treat the case $i = 2$. Since

$$(m, a_i) = (n, a_i) \Rightarrow (m, a_i) | n,$$

Then $[(m, a_1), (m, a_2)] | n$. That is

$$[(m, a_1), (m, a_2)] = \frac{(m, a_1)(m, a_2)}{(m, a_1, a_2)} = \frac{(m^2, m(a_1, a_2), a_1 a_2)}{(m, a_1, a_2)} =$$

$$= \left( m \left( \frac{m}{(m, a_1, a_2)}, \frac{(a_1, a_2)}{(m, a_1, a_2)} \right), \frac{a_1 a_2}{(m, a_1, a_2)} \right) | n. \qquad (1.4)$$

Note that

$$\left( m, \frac{a_1 a_2}{(a_1, a_2)} \right) \Big| m$$

and

$$\left( m, \frac{a_1 a_2}{(a_1, a_2)} \right) \Big| \frac{a_1 a_2}{(a_1, a_2)} \Big| \frac{a_1 a_2}{(m, a_1, a_2)},$$

in view of (1.4), we obtain

$$\left( m, \frac{a_1 a_2}{(a_1, a_2)} \right) \Big| n.$$

Hence,

$$\left( m, \frac{a_1 a_2}{(a_1, a_2)} \right) \Big| \left( n, \frac{a_1 a_2}{(a_1, a_2)} \right).$$

Since the case is quite symmetrical, then

$$\left( n, \frac{a_1 a_2}{(a_1, a_2)} \right) \Big| \left( m, \frac{a_1 a_2}{(a_1, a_2)} \right).$$

Consequently,

$$\left(m, \frac{a_1 a_2}{(a_1, a_2)}\right) = \left(n, \frac{a_1 a_2}{(a_1, a_2)}\right) \Leftrightarrow (m, [a_1, a_2]) = (n, [a_1, a_2]).$$

Now we proceed to the general case. Let

$$(m, a_1) = (n, a_1), (m, a_2) = (n, a_2), ..., (m, a_k) = (n, a_k).$$

Then

$$\begin{cases} (m, a_1) = (n, a_1), \\ (m, a_2) = (n, a_2) \end{cases} \Rightarrow (m, [a_1, a_2]) = (n, [a_1, a_2]),$$

also

$$\begin{cases} (m, [a_1, a_2]) = (n, [a_1, a_2]), \\ (m, a_3) = (n, a_3) \end{cases} \Rightarrow$$

$$\Rightarrow (m, [[a_1, a_2], a_3]) = (n, [[a_1, a_2], a_3]) \Rightarrow$$

$$\Rightarrow (m, [a_1, a_2, a_3]) = (n, [a_1, a_2, a_3]).$$

Continuing by an analogy, we get

$$(m, [a_1, a_2, ..., a_k]) = (n, [a_1, a_2, ..., a_k]). \qquad \square$$

**Property 1.15.** *Let* $(\alpha_1, \alpha_2, ..., \alpha_k) = 1$*, moreover* $\alpha_i | a, \alpha_i, a \neq 0$ *,* $i = 1, ..., k$*. Then*

$$\left[\frac{a}{\alpha_1}, \frac{a}{\alpha_2}, ..., \frac{a}{\alpha_k}\right] = a.$$

**Proof.** We get

$$\left[\frac{a}{\alpha_1}, \frac{a}{\alpha_2}, ..., \frac{a}{\alpha_k}\right] =$$

$$= \frac{\frac{a}{\alpha_1} \frac{a}{\alpha_2} \cdots \frac{a}{\alpha_k}}{\left(\frac{a}{\alpha_2} \cdots \frac{a}{\alpha_k}, \frac{a}{\alpha_1} \frac{a}{\alpha_3} \cdots \frac{a}{\alpha_k}, ..., \frac{a}{\alpha_1} \frac{a}{\alpha_2} \cdots \frac{a}{\alpha_{k-1}}\right)} =$$

$$= \frac{a^k}{(a^{k-1}\alpha_1, a^{k-1}\alpha_2, ..., a^{k-1}\alpha_k)} = \frac{a^k}{a^{k-1}(\alpha_1, \alpha_2, ..., \alpha_k)} = a. \qquad \square$$

**Property 1.16.** *Let*

$$\left(m, \frac{a}{\alpha_i}\right) = \left(n, \frac{a}{\alpha_i}\right),$$

*moreover* $(\alpha_1, \alpha_2, ..., \alpha_k) = 1$*,* $\alpha_i, a \neq 0$*,* $i = 1, ..., k$*. Then*

$$(m, a) = (n, a).$$

**Proof.** The Property 1.14 implies that

$$\left(m, \left[\frac{a}{\alpha_1}, \frac{a}{\alpha_2}, ..., \frac{a}{\alpha_k}\right]\right) = \left(n, \left[\frac{a}{\alpha_1}, \frac{a}{\alpha_2}, ..., \frac{a}{\alpha_k}\right]\right).$$

Based on Property 1.15, we have

$$\left[\frac{a}{\alpha_1}, \frac{a}{\alpha_2}, ..., \frac{a}{\alpha_k}\right] = a.$$

Hence, $(m, a) = (n, a)$. □

**Property 1.17.** *Let*

$$\varphi | \det \left\|\begin{matrix} a & a_i \\ b & b_i \end{matrix}\right\|, (a_i r + b_i, \varphi) = (a_i, b_i, \varphi),$$

*where* $b, a_i, b_i \neq 0, i = 1, ..., k,$ *moreover*

$$((a_1, b_1), (a_2, b_2), ..., (a_k, b_k), \varphi) = 1.$$

*Then*
$$(ar + b, \varphi) = (a, b, \varphi).$$

**Proof.** If $a = 0$ there is nothing to prove. Let now $a \neq 0$. Set $(a_i, b_i, \varphi) = \delta_i, i = 1, ..., k.$ Then

$$\left(\frac{a_i}{\delta_i} r + \frac{b_i}{\delta_i}, \frac{\varphi}{\delta_i}\right) = 1,$$

and

$$\frac{\varphi}{\delta_i} \left|\begin{matrix} a & \dfrac{a_i}{\delta_i} \\ b & \dfrac{b_i}{\delta_i} \end{matrix}\right|.$$

Taking into account Property 1.13, we receive

$$\left(ar + b, \frac{\varphi}{\delta_i}\right) = \left(a, b, \frac{\varphi}{\delta_i}\right),$$

$i = 1, ..., k.$ By assumption $(\delta_1, \delta_2, ..., \delta_k) = 1.$ Having considered Property 1.16, we get $(ar + b, \varphi) = (a, b, \varphi).$ □

**Definition 1.2.** [8] A ring $R$ is adequate if $R$ is a commutative Bezout domain and for every $a \neq 0, b$ in $R$ one can write $a = cd$, with $(c, b) = 1$ and with $(d_i, b) \neq 1$ for every nonunit divisor $d_i$ of $d$.

Commutative $p.i.$-domains, commutative regular rings with identity are examples of adequate rings.

Adequate rings have a property that distinguishes them among the rest of Bezout rings.

**Property 1.18.** *Let $a_1, a_2, a_3$ be relatively prime elements of an adequate ring $R$, where $a_3 \neq 0$. Then there is $r \in R$ such that*

$$(a_1 + a_2 r, a_3) = 1.$$

**Proof.** Write $a_3$ as $a_3 = rt$, where $(r, a_1) = 1$ and any nontrivial divisor $\delta$ of $t$ satisfies the condition

$$(\delta, a_1) \neq 1. \tag{1.5}$$

Let

$$(a_1 + ra_2, a_3) = \tau.$$

Denote $(r, \tau) = \tau_1$. Since $\tau_1 | (a_1 + ra_2)$ and $\tau_1 | r$, then $\tau_1 | a_1$. It means that

$$\tau_1 | (a_1, r) = 1 \Rightarrow \tau = 1.$$

Therefore $(r, \tau) = 1$. Since $\tau_1 | a_3 = rt$, then $\tau | t$. Suppose that $\tau \notin U(R)$. By (1.5),

$$(\tau, a_2) = \tau_2 \neq 1.$$

Furthermore $\tau_2 | (a_1 + ra_2)$. It follows that $\tau_2 | a_2$. Thus,

$$\tau_2 | (a_1, a_2, a_3) = 1,$$

which contradicts our assumption. Thus, $\tau = 1$. Hence, $(a_1 + ra_2, a_3) = 1$. $\square$

The Property 1.18 makes it possible to estimate the stable range of adequate rings, which is important characteristics of ring.

**Definition 1.3.** [9] The **stable range of a ring** $R$ is the smallest positive integer $n$ provided that

$$a_1 R + ... + a_{n+1} R = R, \tag{1.6}$$

implies that there exist $b_1, ..., b_n \in R$ such that

$$(a_1 + a_{n+1} b_1) R + ... + (a_n + a_{n+1} b_n) R = R.$$

In particular, the ring $R$ has stable range 2 if for each triple $a_1, a_2, a_3 \in R$ the equality

$$a_1 R + a_2 R + a_3 R = R$$

implies

$$(a_1 + a_3 r_1) R + (a_2 + a_3 r_2) R = R$$

for some $r \in R$.

The ring $R$ has stable range 1 if for each pair $a_1, a_2 \in R$ the equality $aR + bR = R$ implies $(a + bs) R = R$ for some $s \in R$.

The ring $F[[x]]$ of formal series power over the field $F$ has stable range 1. The integers, a principal ideal rings, the Bezout domains has stable range 2.

The Property 1.18 shows that there is another class of rings between stable rings 1 and 2.

**Definition 1.4.** The ring $R$ has *stable range* 1.5 if for each $a, b \in R$ and $c \in R \setminus \{0\}$ with the property $(a, b, c) = 1$ there exists $r \in R$ such that

$$(a + br, c) = 1.$$

**Property 1.19.** *Let $a_1, ..., a_n$ be relatively prime elements of Bezout ring of stable range 1.5, where $a_n \neq 0$ then there are $u_2, ..., u_{n-1} \in R$ such that*

$$(a_1 + u_2 a_2 + ... + u_{n-1} a_{n-1}, a_n) = 1.$$

**Proof**. Set $(a_2, ..., a_{n-1}) = \gamma$. Then there are elements $v_2, ..., v_{n-1}$ such that
$$v_2 a_2 + ... + v_{n-1} a_{n-1} = \gamma.$$
Since $(a_1, \gamma, a_n) = 1$, there is an element $r \in R$, for which $(a_1 + r\gamma, a_n) = 1$. Thus

$$(a_1 + (rv_2)a_2 + ... + (rv_{n-1})a_{n-1}, a_n) = 1. \qquad \square$$

Note that not all commutative Bezout rings have a stable range 1.5.

**Example 1.1.** [11] Consider the ring

$$Q = \left\{ a_0 + \sum_{i=1}^{\infty} a_i x^i \mid a_0 \in \mathbb{Z}, \ a_i \in \mathbb{Q}, \ i \in \mathbb{N} \right\}.$$

Make sure that $Q$ is a Bezout ring which stable range is not equal to 1.5. The group $U(Q)$ consists of all elements of the form

$$\pm 1 + \sum_{i=1}^{\infty} a_i x^i.$$

It follows that the elements of the ring $Q$ are represented as:

1) if $a_0 \neq 0$, then
$$\alpha = a_0 e_\alpha,$$
where
$$a_0 \in \mathbb{Z}^*, \ \ e_\alpha = \left( 1 + \sum_{i=1}^{\infty} \frac{a_i}{a_0} x^i \right) \in U(Q);$$

2) if $a_0 = 0$, then

$$\beta = \sum_{i=k}^{\infty} \frac{p_i}{q_i} x^i = \frac{p_k}{q_k} x^k \left( 1 + \sum_{i=k+1}^{\infty} \frac{p_i q_k}{q_i p_k} x^i \right) = \frac{p_k}{q_k} x^k e_\beta,$$

where $p_k, q_k \in \mathbb{Z}^*$, $e_\beta \in U(R)$.

Let us prove that $Q$ is a Bezout ring.

Let $\mu = me_\mu, \nu = ne_\nu$, where $m, n \in \mathbb{Z}$ and $e_\mu, e_\nu \in U(Q)$. There are $u, v \in \mathbb{Z}$ such that
$$mu + nv = (m, \ n).$$
Then
$$me_\mu(e_\mu^{-1}u) + ne_\nu(e_\nu^{-1}v) = (m, n) = (\mu, \nu).$$
Suppose that
$$\alpha = \frac{s_l}{t_l}x^l e_\alpha, \ \ \beta = \frac{p_k}{q_k}x^k e_\beta, \ \ e_\alpha, e_\beta \in U(Q),$$
where $0 \leqslant l < k$. Then
$$\beta = \alpha \frac{p_k t_l}{q_k s_l}x^{k-l}e_\beta e_\alpha^{-1}.$$
Consequently,
$$(\alpha, \beta) = \alpha \Rightarrow \alpha 1 + \beta 0 = \alpha.$$

It remains to consider the case
$$\alpha = \frac{s_k}{t_k}x^k e_\alpha, \ \beta = \frac{p_k}{q_k}x^k e_\beta,$$
where $e_\alpha, e_\beta \in U(Q)$, and $k \geqslant 1$. Set
$$(s_k q_k, \ t_k p_k) = d$$
so
$$(s_k q_k)u + (\ t_k p_k)v = d$$
for some $u, v$. Then
$$(\alpha, \beta) = \frac{d}{t_k q_k}x^k \left( \frac{s_k q_k}{d}e_\alpha, \ \frac{p_k t_k}{d}e_\beta \right) = \frac{d}{t_k q_k}x^k.$$
Moreover
$$\alpha u e_\alpha^{-1} + \beta v e_\beta^{-1} = (\alpha, \beta).$$
Therefore, $Q$ is a Bezout ring.

Consider the triple of relatively prime elements: 3, 7, $x$. Note that $x$ is relatively prime only with units of the ring $Q$. Consider the equation
$$3 + 7r = \pm 1 + b_1 x + b_2 x^2 + ...,$$
where $r$ is variable. Its solutions are
$$r_1 = -\frac{4}{7} + c_1 x + c_2 x^2 + ...,$$
$$r_2 = -\frac{2}{7} + d_1 x + d_2 x^2 + ....$$
However, $r_1, r_2$ are not elements of the ring $Q$. It means that $(7 + 5k, x) \neq 1$ for all $k \in \mathbb{Z}$. Consequently, $Q$ is not the ring of stable range 1.5. $\qquad \diamondsuit$

## 1.2. Complementary of a primitive
## row to an invertible matrix

The matrix is called **primitive** provided that g.c.d. of its maximal order minors is equal to 1. So, a row (column) that consists of relatively prime elements is called a **primitive row (column)**.

The problem of a complementary of a primitive row (column) to an invertible matrix arises while solving a lot of matrix problems. It turns out that the choice of the ring significantly influences the appearance of the desired invertible matrix.

**Theorem 1.1.** *Let $R$ be a Bezout ring and $\|a_1 \ldots a_n\|$ be a row over $R$. There is a matrix*

$$A = \left\| \begin{matrix} u_{11} & u_{12} & \ldots & u_{1.n-2} & u_{1.n-1} & u_{1n} \\ 0 & u_{22} & \ldots & u_{2.n-2} & u_{2.n-1} & u_{2n} \\ \ldots & \ldots & \ldots & \ldots & \ldots & \ldots \\ 0 & 0 & & u_{n-2.n-2} & u_{n-2.n-1} & u_{n-2.n} \\ 0 & 0 & \ldots & 0 & u_{n-1.n-1} & u_{n-1.n} \\ a_1 & a_2 & \ldots & a_{n-2} & a_{n-1} & a_n \end{matrix} \right\|, \tag{1.7}$$

*such that* $\det A = (a_1, ..., a_n) = \alpha$.

**Proof**. Suppose that $n = 2$ and

$$a_2 u_{11} - a_1 u_{12} = (a_1, \ a_2).$$

The matrix

$$\left\| \begin{matrix} u_{11} & u_{12} \\ a_1 & a_2 \end{matrix} \right\|$$

will be desired. Suppose that our statement is true for all matrices of an order less than $n$. Then there exists a matrix

$$\left\| \begin{matrix} B \\ \hline a_2 & \ldots & a_n \end{matrix} \right\| = D$$

with the determinant $(a_2, ..., a_n)$. Since

$$\det D = \sum_{i=2}^{n} (-1)^{i+n-1} a_i d_i,$$

where $d_i$ is the corresponding minor of order $n - 2$ of the matrix $B$, then

$$(a_2, ..., a_n)(d_2, ..., d_n)|(a_2, ..., a_n).$$

Hence, $(d_2, ..., d_n) = 1$. It means that there is an invertible matrix of the form

$$\left\| \begin{matrix} c_2 & \ldots & c_n \\ \hline B \end{matrix} \right\|.$$

There are $p, q$, that

$$a_1 q + (a_2, ..., a_{n-1})p = (a_1, ..., a_{n-1}, a_n).$$

Then the desired matrix will be

$$\left\|\begin{array}{c|ccc} p & (-1)^{n-1}qc_2 & ... & (-1)^{n-1}qc_n \\ \hline \mathbf{0} & & B & \\ \hline a_1 & a_2 & ... & a_n \end{array}\right\|. \qquad \square$$

This result was obtained for the first time by Sh. Hermite [72].

**Corollary 1.1.** *Let* $\|a_1 ... a_n\|$ *be a primitive row. Then there is an invertible matrix of the form* (1.7). $\qquad \square$

**Remark**. By permuting of the rows or transposing and permuting the columns of the matrix (1.7), we can construct an invertible matrix where the given primitive row is in an arbitrary position.

The complementary of a primitive row to an invertible matrix is greatly simplified over Bezout rings of stable range 1.5.

**Theorem 1.2.** *Let* $\|a_1 ... a_n\|$ *be a primitive row over Bezout ring of stable range* 1.5, *where* $a_1 \neq 0$, *then there is an invertible matrix of the form*

$$\left\|\begin{array}{cccccc} u_n & 0 & ... & 0 & 0 & u_1 \\ 0 & 1 & & 0 & 0 & u_2 \\ ... & ... & \ddots & & ... & ... \\ 0 & 0 & & 1 & 0 & u_{n-2} \\ 0 & 0 & ... & 0 & 1 & u_{n-1} \\ a_1 & a_2 & ... & a_{n-2} & a_{n-1} & a_n \end{array}\right\| = U.$$

**Proof.** According to Property 1.19, there are $u_1, ..., u_{n-1}$ such that

$$(a_n - u_{n-1}a_{n-1} - ... - u_2 a_2, a_1) = 1.$$

So there exist $u_1$, $u_n$ that

$$u_n(a_n - u_{n-1}a_{n-1} - ... - u_2 a_2) - u_1 a_1 = 1.$$

Then

$$\det U =$$

$$= u_n \left\|\begin{array}{ccccc} 1 & & 0 & 0 & u_2 \\ & \ddots & & & \vdots \\ 0 & & 1 & 0 & u_{n-2} \\ 0 & ... & 0 & 1 & u_{n-1} \\ a_2 & ... & a_{n-2} & a_{n-1} & a_n \end{array}\right\| + (-1)^{n+1}a_1 \left\|\begin{array}{ccccc} 0 & ... & 0 & 0 & u_1 \\ 1 & & & & u_2 \\ & \ddots & & & \vdots \\ 0 & & 1 & 0 & u_{n-2} \\ 0 & ... & 0 & 1 & u_{n-1} \end{array}\right\| =$$

$$= u_n \left( (-1)^n a_2 \left\| \begin{matrix} 0 & \dots & 0 & 0 & u_2 \\ 1 & & 0 & 0 & u_3 \\ & \ddots & & & \vdots \\ 0 & & 1 & 0 & u_{n-2} \\ 0 & \dots & 0 & 1 & u_{n-1} \end{matrix} \right\| + (-1)^{n+1} a_3 \left\| \begin{matrix} 1 & \dots & 0 & 0 & u_2 \\ 0 & \dots & 0 & 0 & u_3 \\ 0 & 1 & 0 & 0 & u_4 \\ \vdots & & \ddots & & \vdots \\ 0 & \dots & 0 & 1 & u_{n-1} \end{matrix} \right\| \right) +$$

$$+ \dots + u_n \left( (-1)^{2(n-1)-1} a_{n-1} \left\| \begin{matrix} 1 & & 0 & u_2 \\ & \ddots & & \vdots \\ 0 & & 1 & u_{n-2} \\ 0 & \dots & 0 & u_{n-1} \end{matrix} \right\| + a_n \right) + (-1)^{2n+1} a_1 u_1 =$$

$$= u_n(a_n - u_{n-1}a_{n-1} - \dots - u_2 a_2) - u_1 a_1 = 1. \qquad \square$$

Now we can establish one of the important properties of matrices over the Bezout rings.

**Theorem 1.3.** *Set* $(a_1, ..., a_n) = \alpha$. *Then there exists an invertible matrix* $U$ *such that*

$$U \left\| a_1 \dots a_n \right\|^T = \left\| \alpha \quad 0 \quad \dots \quad 0 \right\|^T.$$

**Proof.** Since $(a_1, ..., a_n) = \alpha$, then there exist such elements $u_1, ..., u_n$ that

$$u_1 a_1 + \dots + u_n a_n = \alpha.$$

It follows that

$$(u_1, ..., u_n) = 1.$$

It means that the row $\| u_1 \dots u_n \|$ is primitive. Based on Theorem 1.2, it can be complemented to an invertible matrix of the form

$$\left\| \begin{matrix} v_{11} & \dots & v_{1n} \\ \dots & \dots & \dots \\ \mathbf{0} & \dots & v_{n-1.n} \\ u_1 & \dots & u_n \end{matrix} \right\|.$$

Then

$$\left\| \begin{matrix} v_{11} & \dots & v_{1n} \\ \dots & \dots & \dots \\ \mathbf{0} & \dots & v_{n-1.n} \\ u_1 & \dots & u_n \end{matrix} \right\| \left\| a_1 \dots a_n \right\|^T = \left\| b_1 \quad \dots \quad b_{n-1} \quad \alpha \right\|^T,$$

where

$$b_i = v_{i1}a_1 + \dots + v_{in}a_n, \quad i = 1, ..., n-1.$$

Noting that $\alpha$ is a divisor of all the elements of a column $\|a_1 \dots a_n\|^T$ we conclude that $\alpha | b_i$, $i = 1, ..., n-1$. Therefore,

$$\left\| \begin{matrix} 0 & 0 & & 0 & 1 \\ 1 & 0 & & 0 & -\dfrac{b_1}{\alpha} \\ & & \ddots & & \vdots \\ 0 & 0 & & 1 & -\dfrac{b_{n-1}}{\alpha} \end{matrix} \right\| \left\| b_1 \quad \dots \quad b_{n-1} \quad \alpha \right\|^T = \left\| \alpha \quad 0 \quad \dots \quad 0 \right\|^T.$$

Therefore, the matrix

$$U = \left\| \begin{matrix} 0 & 0 & & 0 & 1 \\ 1 & 0 & & 0 & -\dfrac{b_1}{\alpha} \\ & & \ddots & & \vdots \\ 0 & 0 & & 1 & -\dfrac{b_{n-1}}{\alpha} \end{matrix} \right\| \left\| \begin{matrix} v_{11} & \dots & v_{1n} \\ \dots & \dots & \dots \\ \mathbf{0} & \dots & v_{n-1.n} \\ u_1 & \dots & u_n \end{matrix} \right\|.$$

will be desired. $\qquad \square$

Using the matrix transposition operation, we obtain the following result.

**Corollary 1.2.** Set $(a_1, ..., a_n) = \alpha$. Then there exists such an invertible matrix $V$, that

$$\|a_1 \dots a_n\| V = \left\| \alpha \quad 0 \quad \dots \quad 0 \right\|. \qquad \square$$

## 1.3. Hermite normal form

Matrix transformations change their form. But there are quantities that do not change. That is, they are invariants of such transformations. One of the important invariants of the matrices with respect to acting from $\mathrm{GL}_n(R)$ is g.c.d. of fixed-order minors.

**Theorem 1.4.** *Let $A$ be an $m \times n$ matrix, $V \in \mathrm{GL}_m(R)$, $U \in \mathrm{GL}_n(R)$. Then g.c.d. of $k$th order minors of the matrix $A$ coincides with g.c.d. of the corresponding minors of the matrix $VAU$, $k = 1, 2, ..., \min(m, n)$.*

**Proof.** We put $m \leqslant n$ for definiteness. Denote by $\delta_k$ g.c.d. of $k$th order minors of the matrix $A$, and $\Delta_k$ of the matrix $AU$, respectively.

Consider

$$AU = \left\| \begin{matrix} a_{11} & \dots & a_{1m} & \dots & a_{1n} \\ \dots & \dots & \dots & \dots & \dots \\ a_{m1} & \dots & a_{mm} & \dots & a_{mn} \end{matrix} \right\| \left\| \begin{matrix} u_{11} & \dots & u_{1n} \\ \dots & \dots & \dots \\ u_{n1} & \dots & u_{nn} \end{matrix} \right\| = \left\| \begin{matrix} b_{11} & \dots & b_{1n} \\ \dots & \dots & \dots \\ b_{m1} & \dots & b_{mn} \end{matrix} \right\| = B.$$

According to the Binet—Cauchy formula, every $k$th order minor of this matrix

has the form

$$
\begin{vmatrix} b_{i_1j_1} & ... & b_{i_1j_k} \\ ... & ... & ... \\ b_{i_kj_1} & ... & b_{i_kj_k} \end{vmatrix} = \sum_{1 \leqslant l_1 < ... < l_k \leqslant m} \begin{vmatrix} a_{i_1l_1} & ... & a_{i_1l_k} \\ ... & ... & ... \\ a_{i_kl_1} & ... & a_{i_kl_k} \end{vmatrix} \begin{vmatrix} u_{l_1j_1} & ... & u_{l_1j_n} \\ ... & ... & ... \\ u_{l_kj_1} & ... & u_{l_kj_n} \end{vmatrix},
$$

$1 \leqslant i_1 < ... < i_k \leqslant m$, $1 \leqslant j_1 < ... < j_k \leqslant n$. Since $\delta_k$ divides all $k$th order minors of the matrix $A$, it follows from this equality that $\delta_k$ will also be the divisor of all $k$th order minors of the $AU$ matrix. It means that $\delta_k | \Delta_k$, $k = 1, ..., m$.

On the other hand, $A = BV$, where $V = U^{-1}$. Repeating the above considerations, we show that $\Delta_k | \delta_k$, $k = 1, ..., m$. Hence, $\delta_k = \Delta_k$, $k = 1, ..., m$.

It is similarly proved that $\delta_k$ coincides with the $k$th order minors of the matrix $VA$. Finally, taking into account the associativity of multiplication of a matrices, we complete the proof. □

Recall that the **rank of** $A$ is the largest order of any non-zero minor in $A$. We say that an $m \times n$ matrix $A$ has a maximal rank if rank $A = \min(m, n)$.

**Corollary 1.3.** *Let $A$ be an $m \times n$ matrix and $V \in \mathrm{GL}_m$, $U \in \mathrm{GL}_n(R)$. Then*

$$
\mathrm{rank}\, A = \mathrm{rank}\, VAU. \qquad \square
$$

The matrices $A, B$ are **right associates** if there exists an invertible matrix $V$ such that $A = BV$. The relation of being right associates is an equivalence relation. Naturally, the question of choosing in $A\,\mathrm{GL}_n(R)$ a matrix of the simplest structure arises.

In this subsection, $R$ is the Bezout ring.

**Lemma 1.1.** *Let $A$, $B$ be an $n \times n$ lower triangular matrices over $R$, and the matrix $A$ is nonsingular. If $AC = B$, then $C$ is also a lower triangular matrix.*

**Proof.** Suppose that

$$
AC = \begin{Vmatrix} a_{11} & 0 & ... & 0 \\ a_{21} & a_{22} & & 0 \\ & & \ddots & \\ a_{n1} & ... & a_{n.n-1} & a_{nn} \end{Vmatrix} \begin{Vmatrix} c_{11} & c_{12} & ... & c_{1n} \\ c_{21} & c_{22} & ... & c_{2n} \\ ... & ... & ... & ... \\ c_{n1} & c_{n2} & ... & c_{nn} \end{Vmatrix} =
$$

$$
= \begin{Vmatrix} b_{11} & 0 & ... & 0 \\ b_{21} & b_{22} & & 0 \\ & & \ddots & \\ b_{n1} & ... & b_{n.n-1} & b_{nn} \end{Vmatrix} = B.
$$

Since the matrix $A$ is nonsingular, then $a_{ii} \neq 0$, $i = 1, ..., n$. Therefore, from the equals $a_{11}c_{1j} = 0$, $j = 2, ..., n$ it implies that $c_{1j} = 0$. Also from the equals $a_{22}c_{2k} = 0$, $k = 3, ..., n$ we get $c_{2k} = 0$. Continuing these considerations, we are convinced of the correctness of our assertion. □

Let $a \in R$. Denote by $K(a)$ the complete residue system modulo $a$. We also denote by $Z(R)$ the set of unassociated elements of the ring $R$.

**Theorem 1.5.** *Let* $A = \|a_{ij}\|$ *be an* $n \times m$ *matrix* $(n \leqslant m)$ *and* $\operatorname{rank} A = n$. *Then there exists such an invertible matrix* $U$ *that*

$$AU = \|H \quad \mathbf{0}\|, \quad H = \left\|\begin{array}{cccc} \alpha_1 & 0 & ... & 0 \\ b_{21} & \alpha_2 & & 0 \\ & & \ddots & \\ b_{n1} & ... & b_{n.n-1} & \alpha_n \end{array}\right\|, \tag{1.8}$$

*where* $\alpha_i \in Z(R)$, $i = 1, ..., n$, $b_{ij} \in K(\alpha_i)$ *for all* $i > j$. *And in the class* $A\operatorname{GL}_m(R)$ *the matrix of the form* (1.8) *is unique.*

**Proof.** Let

$$(a_{11}, ..., a_{1m}) = \alpha_1,$$

where $\alpha_1 \in Z(R)$. By Corollary 1.2, there exists such an invertible matrix $U_1$ that

$$\|a_{11} ... a_{1m}\| U_1 = \|\alpha_1 \quad 0 \quad ... \quad 0\|.$$

Thus,

$$AU_1 = \left\|\begin{array}{cccc} \alpha_1 & 0 & ... & 0 \\ c_{21} & c_{22} & ... & c_{2m} \\ ... & ... & ... & ... \\ c_{n1} & c_{n2} & ... & c_{nm} \end{array}\right\|.$$

Let

$$(c_{22}, ..., c_{2m}) = \alpha_2,$$

where $\alpha_2 \in Z(R)$. There exists an invertible matrix $U_2'$ such that

$$\|c_{22} ... c_{2m}\| U_2' = \|\alpha_2 \quad 0 \quad ... \quad 0\|.$$

It means that

$$AU_1 U_2' = \left\|\begin{array}{ccccc} \alpha_1 & 0 & 0 & ... & 0 \\ c_{21} & \alpha_2 & 0 & ... & 0 \\ c_{31} & d_{32} & d_{33} & ... & d_{3m} \\ ... & ... & ... & ... & ... \\ c_{n1} & d_{n2} & d_{n3} & ... & d_{nm} \end{array}\right\|,$$

where $U = \left\|\begin{array}{cc} 1 & \mathbf{0} \\ \mathbf{0} & U_2' \end{array}\right\|$. Continuing our considerations, we obtain the following invertible matrix $V$, which is the product of invertible matrices $U_1, U_2, ..., U_{n-1}$ and

$$AV = \|H_1 \ \mathbf{0}\| = A_1, \quad H_1 = \left\|\begin{array}{cccc} \alpha_1 & 0 & ... & 0 \\ f_{21} & \alpha_2 & & 0 \\ & & \ddots & \\ f_{n1} & ... & f_{n.n-1} & \alpha_n \end{array}\right\|,$$

where $\alpha_i \in Z(R)$, $i = 1, ..., n$. Since $\operatorname{rank} A = n$, is based on Corollary 1.3, all $\alpha_i \neq 0$.

In the set $K(\alpha_2)$, there exists an element $b_{21}$ such that

$$f_{21} \equiv b_{21} (\operatorname{mod} \alpha_2).$$

That is $f_{21} = b_{21} + \alpha_2 r_{21}$, where $r_{21} \in R$. Then

$$A_1 \underbrace{\begin{Vmatrix} 1 & 0 & 0 & & 0 \\ -r_{21} & 1 & 0 & & 0 \\ 0 & 0 & 1 & & 0 \\ & & & \ddots & \\ 0 & 0 & 0 & & 1 \end{Vmatrix}}_{V_2} = \begin{Vmatrix} H_2 & \mathbf{0} \end{Vmatrix}, \quad H_2 = \begin{Vmatrix} \alpha_1 & 0 & 0 & & 0 \\ b_{21} & \alpha_2 & 0 & & 0 \\ f'_{31} & f_{32} & \alpha_3 & & 0 \\ & & & \ddots & \\ f'_{n1} & f_{n2} & ... & f_{n.n-1} & \alpha_n \end{Vmatrix}.$$

In the set $K(\alpha_3)$ there are such elements $b_{31}$, $b_{32}$, that

$$f'_{31} = b_{31} + \alpha_3 r_{31}, \quad f_{32} = b_{32} + \alpha_3 r_{32},$$

where $r_{31}$, $r_{32} \in R$. Then

$$A_1 V_2 \begin{Vmatrix} 1 & 0 & 0 & 0 & & 0 \\ 0 & 1 & 0 & 0 & & 0 \\ -r_{31} & -r_{32} & 1 & 0 & & 0 \\ 0 & 0 & 0 & 1 & & 0 \\ & & & & \ddots & \\ 0 & 0 & 0 & 0 & & 1 \end{Vmatrix} = \begin{Vmatrix} H_3 & \mathbf{0} \end{Vmatrix},$$

$$H_3 = \begin{Vmatrix} \alpha_1 & 0 & 0 & 0 & & 0 \\ b_{21} & \alpha_2 & 0 & 0 & & 0 \\ b_{31} & b_{32} & \alpha_3 & 0 & & 0 \\ f''_{41} & f'_{42} & f_{43} & \alpha_4 & & 0 \\ & & & & \ddots & \\ f''_{n1} & f'_{n2} & f_{n3} & & & \alpha_n \end{Vmatrix}.$$

Continuing this process, we get the required invertible matrix $U$. Thus, we proved that the class $A \operatorname{GL}_n(R)$ has a matrix of the form (1.8). Now we will show that such a matrix is unique.

Suppose that the class $A \operatorname{GL}_n(R)$ has a matrix of the form $\begin{Vmatrix} B & \mathbf{0} \end{Vmatrix}$, where

$$B = \begin{Vmatrix} \beta_1 & 0 & ... & 0 \\ c_{21} & \beta_2 & & 0 \\ & & \ddots & \\ c_{n1} & ... & c_{n.n-1} & \beta_n \end{Vmatrix},$$

moreover $\beta_i \in Z(R)$, $i = 1, ..., n$, $c_{ij} \in K(\beta_i)$, $i > j$. There exists an invertible matrix $L$ such that $B = AL$. Based on Lemma 1.1, $L$ is a lower triangular matrix. That is

$$\begin{Vmatrix} \alpha_1 & 0 & ... & 0 \\ b_{21} & \alpha_2 & & 0 \\ & & \ddots & \\ b_{n1} & ... & b_{n.n-1} & \alpha_n \end{Vmatrix} \begin{Vmatrix} e_1 & 0 & ... & 0 \\ l_{21} & e_2 & & 0 \\ & & \ddots & \\ l_{n1} & ... & l_{n.n-1} & e_n \end{Vmatrix} =$$

$$= \begin{Vmatrix} \beta_1 & 0 & ... & 0 \\ c_{21} & \beta_2 & & 0 \\ & & \ddots & \\ c_{n1} & ... & c_{n.n-1} & \beta_n \end{Vmatrix}. \tag{1.9}$$

Since the matrix $L$ is invertible, then $e_i \in U(R)$, $i = 1, ..., n$. Consequently, $\beta_i = \alpha_i e_i$, $i = 1, ..., n$. It means that $\beta_i, \alpha_i$ are representatives of the same class of associating elements of the ring $R$. Therefore, $\beta_i = \alpha_i$, $i = 1, ..., n$. So $e_i = 1$, $i = 1, ..., n$. From the equality (1.9), we get

$$c_{i.i-1} = b_{i.i-1} + \alpha_i l_{i.i-1}, \quad i = 2, ..., n.$$

That is

$$c_{i.i-1} \equiv b_{i.i-1} (\text{mod } \alpha_i).$$

Hence, $c_{i.i-1}$, $b_{i.i-1} \in K(\alpha_i)$. Thus, $c_{i.i-1} = b_{i.i-1}$, $l_{i.i-1} = 0$, $i = 2, ..., n$. Continuing similar considerations, we obtain that $A = B$. □

The matrix (1.8) is called the **Hermite form** of the matrix $A$. It was named after the famous French mathematician Ch. Hermite [72].

**Theorem 1.6.** *Let $A = \|a_{ij}\|$ be an $m \times n$ matrix and the rank $A = r < \min(m, n)$. Then there are an invertible matrices $U$, $V$ such that*

$$UAV = \begin{Vmatrix} M & \mathbf{0} \\ \mathbf{0} & \mathbf{0} \end{Vmatrix},$$

*where $M$ is a nonsingular matrix of the order $r$.*

**Proof.** For definiteness we put $m \leqslant n$. If the matrix $A$ has the maximal rank, then $U = I$, and the matrix $V$ by the Theorem 1.5, is the matrix that reduces the matrix $A$ to its Hermite form.

Suppose that the rank $A = r < m$. Then the matrix $A$ contains a nonsingular submatrix $A_{11}$ of the order $r$. By permuting columns and rows, which is equivalent to multiplying by a nonsingular matrix $U_1, V_1,$, we get the matrix $U_1 A V_1$ has the form

$$U_1 A V_1 = \begin{Vmatrix} A_{11} & A_{12} \\ A_{21} & A_{22} \end{Vmatrix}.$$

Since the matrix $\left\| A_{11} \; A_{12} \right\|$ has the maximal rank, so according to Theorem 1.5, there exists an invertible matrix $V_2$ such that

$$\left\| A_{11} \; A_{12} \right\| V_2 = \left\| \begin{array}{ccc|c} \alpha_1 & 0 & 0 & \mathbf{0} \\ & \ddots & 0 & \mathbf{0} \\ * & * & \alpha_r & \mathbf{0} \end{array} \right\|$$

is the Hermite form of the matrix $\left\| A_{11} \; A_{12} \right\|$, where $\alpha_1 \ldots \alpha_r \neq 0$.

Consider the matrix

$$U_1 A (V_1 V_2) = \left\| \begin{array}{ccc|ccc} \alpha_1 & 0 & 0 & 0 & \ldots & 0 \\ & \ddots & 0 & 0 & \ldots & 0 \\ * & * & \alpha_r & 0 & \ldots & 0 \\ \hline b_{r+1.1} & \ldots & b_{r+1.r} & b_{r+1.r+1} & \ldots & b_{r+1.n} \\ \ldots & \ldots & \ldots & \ldots & \ldots & \ldots \\ b_{m1} & \ldots & b_{mr} & b_{m.r+1} & \ldots & b_{mn} \end{array} \right\| = \left\| \begin{array}{cc} B_{11} & \mathbf{0} \\ B_{21} & B_{22} \end{array} \right\|.$$

Since all $r + 1$th order minors of this matrix are zero, then

$$\alpha_1 \ldots \alpha_r b_{ij} = 0, \quad i = r+1, r+2, \ldots, m, \quad j = r+1, r+2, \ldots, n.$$

It follows that $b_{ij} = 0$, where $i = r+1, r+2, \ldots, m, \; j = r+1, r+2, \ldots, n$. Consequently, $B_{22} = \mathbf{0}$. The matrix $\left\| \begin{array}{c} B_{11} \\ B_{21} \end{array} \right\|$ has the maximal rank. Again, there is such an invertible matrix $U_2$ that

$$U_2 \left\| \begin{array}{c} B_{11} \\ B_{21} \end{array} \right\| = \left\| \begin{array}{ccc} \beta_1 & 0 & 0 \\ & \ddots & 0 \\ * & * & \beta_r \\ \hline \mathbf{0} & \mathbf{0} & \mathbf{0} \end{array} \right\| = \left\| \begin{array}{c} C \\ \mathbf{0} \end{array} \right\|,$$

where $\det C = \beta_1 \ldots \beta_r \neq 0$. This yields

$$(U_2 U_1) A (V_1 V_2) = \left\| \begin{array}{cc} C & \mathbf{0} \\ \mathbf{0} & \mathbf{0} \end{array} \right\|.$$

The theorem is proved. $\qquad \square$

From the method of reducing the matrix to its Hermite form, it follows that $\alpha_1$ is g.c.d. of the first row elements of the matrix $A$. Using the Theorem 1.5, we can also find other diagonal elements of a nonsingular matrix $A$. Indeed, from the elements of the first two rows of the matrix $AU$ we can construct only one nonzero minor, namely $\alpha_1 \alpha_2$. That is, $\alpha_1 \alpha_2$ is g.c.d. of all the second order minors constructed from the elements of the first two rows of the matrix $AU$. By Theorem 1.5, the elements of the first two rows of $A$ have the same property. Thus, the product of the first two diagonal elements of the matrix $AU$ is equal to g.c.d. of the second order minors constructed from the elements of the first two rows of the matrix $A$. All other products $\alpha_1 \dots \alpha_i$, $i = 1, \dots, n$, can be found in much the same way.

**Theorem 1.7.** *Suppose that $AU$ is the Hermite form of a nonsingular $n \times n$ matrix $A$ and $\alpha_i$, $i = 1, \dots, n$, are diagonal elements of $AU$. Then $\alpha_1 \dots \alpha_i$ is equal to g.c.d. of maximal order minors of the matrix $A$ that are constructed on its first $i$ rows.*                                                                         □

## 1.4. Bezout's identity and its properties

Suppose that $a_1, a_2$ are relatively prime elements of a Bezout ring $R$, then there exist $u_1, u_2 \in R$, that

$$a_1 u_1 + a_2 u_2 = 1.$$

It follows from this equality, that for an arbitrary $n$th of relatively prime elements $a_1, a_2, \dots, a_n$ there exist $v_1, v_2, \dots, v_n$ that

$$a_1 v_1 + a_2 v_2 + \dots + a_n v_n = 1. \tag{1.10}$$

So

$$\left\| a_1 \quad a_2 \quad \dots \quad a_n \right\| \left\| v_1 \quad v_2 \quad \dots \quad v_n \right\|^T = 1.$$

We say that the elements of the row $\left\| v_1 \quad v_2 \quad \dots \quad v_n \right\|$ satisfy equality (1.10). Denote by $\mathbf{U}_{a_1, a_2, \dots, a_n}$ the set of all rows which elements satisfy equality (1.10).

Since $\left\| a_1 \quad a_2 \quad \dots \quad a_n \right\|$ is a primitive row, by Theorem 1.1, there is an invertible matrix of the form

$$A = \left\| \begin{matrix} a_1 \\ a_2 \\ \vdots \\ a_n \end{matrix} \quad * \quad \right\|. \tag{1.11}$$

We suggest the method of finding all rows that satisfy the equality (1.10).

**Theorem 1.8.** *Let $(a_1, a_2, \dots, a_n) = 1$, $n \geq 2$, and $A$ be invertible matrix of the form (1.11). Then*

$$\mathbf{U}_{a_1, a_2, \dots, a_n} = \{ \left\| 1 \quad x_2 \quad \dots \quad x_n \right\| A^{-1} \},$$

*where $x_i$ independently acquire values in $R$, $i = 2, 3, \dots, n$.*

**Proof.** Consider the set

$$\mathbf{V} = \{ \left\| 1 \quad x_2 \quad \dots \quad x_n \right\| A^{-1} | x_i \in R, \ i = 2, 3, \dots, n \}.$$

Let $\left\| v_1 \quad v_2 \quad \dots \quad v_n \right\| \in \mathbf{V}$, i.e.,

$$\left\| v_1 \quad v_2 \quad \dots \quad v_n \right\| = \left\| 1 \quad b_2 \quad \dots \quad b_n \right\| A^{-1},$$

where $b_i \in R, i = 2, ..., n$. Then

$$\left\|v_1 \quad v_2 \quad ... \quad v_n\right\|\left\|a_1 \quad a_2 \quad ... \quad a_n\right\|^T =$$

$$= \left\|1 \quad b_2 \quad ... \quad b_n\right\| A^{-1} \left\|a_1 \quad a_2 \quad ... \quad a_n\right\|^T =$$

$$= \left\|1 \quad b_2 \quad ... \quad b_n\right\|\left\|1 \quad 0 \quad ... \quad 0\right\|^T = 1.$$

Hence, $\mathbf{V} \subseteq \mathbf{U}_{a_1, a_2, ..., a_n}$.

On the contrary, let $\left\|u_1 \quad u_2 \quad ... \quad u_n\right\| \in \mathbf{U}_{a_1, a_2, ..., a_n}$ and $A^{-1} = \left\|b_{ij}\right\|_1^n$. Consider the matrix

$$\left\|\begin{matrix} u_1 & u_2 & ... & u_n \\ b_{21} & b_{22} & ... & b_{2n} \\ ... & ... & ... & ... \\ b_{n1} & b_{n2} & ... & b_{nn} \end{matrix}\right\| = U.$$

Then

$$UA = \left\|\begin{matrix} 1 & c_2 & ... & c_n \\ 0 & 1 & ... & 0 \\ \vdots & & \ddots & \\ 0 & 0 & ... & 1 \end{matrix}\right\|.$$

It follows that

$$\left\|u_1 \quad ... \quad u_n\right\| = \left\|1 \quad c_2 \quad ... \quad c_n\right\| A^{-1}.$$

Therefore, $\mathbf{U}_{a_1, a_2, ..., a_n} \subseteq \mathbf{V}$. Hence, $\mathbf{U}_{a_1, a_2, ..., a_n} = \mathbf{V}$. $\quad\square$

**Corollary 1.4.** *If $a$, $b \in R$ and*

$$au + bv = (a, \ b),$$

*then $\mathbf{U}_{a,b} = \{(u + br, \ v - ar)\}$, where $r \in R$.* $\quad\square$

Bezout Equality in rings of stable range 1.5 has additional properties that do not hold in rings of stable range 2.

**Theorem 1.9.** *Let $R$ be Bezout rings of stable range $1.5$ and*

$$(a_1, ..., a_n) = 1,$$

*$n \geq 2$ and $\psi$ is an arbitrary fixed nonzero element in $R$. Then there are $u_1, ..., u_n$, that satisfy simultaneously equalities:*

1) $u_1 a_1 + ... + u_n a_n = 1$;

2) $(u_1, ..., u_i) = 1$ for an arbitrary fixed $i$, where $2 \leq i \leq n$;

3) $(u_i, \psi) = 1$ for an arbitrary fixed $i$, where $2 \leq i \leq n$.

**Proof.** At first we will show that the column $\left\|a_1 \quad a_2 \quad ... \quad a_n\right\|^T$ can be complemented to an invertible matrix $A$ such that $A^{-1} = \left\|b_{ij}\right\|_1^n$, where

$$b_{3i} = ... = b_{ni} = 0.$$

Consider an arbitrary invertible matrix of the form

$$A = \left\| \begin{matrix} a_1 \\ a_2 \\ ... & * \\ a_n \end{matrix} \right\|.$$

Suppose that among the elements $\bar{b}_{2i}, ..., \bar{b}_{ni}$ of the matrix $A_1^{-1} = \left\| \bar{b}_{ij} \right\|_1^n$ there is at least one nonzero. Then there is a matrix $D \in \mathrm{GL}_{n-1}(R)$ such that

$$D \left\| \bar{b}_{2i} \quad ... \quad \bar{b}_{ni} \right\|^T = \left\| \gamma \quad 0 \quad ... \quad 0 \right\|^T.$$

Then the matrix

$$((1 \oplus D)A_1^{-1})^{-1} = A_1(1 \oplus D^{-1}) = A$$

will be desirable.

Consider a matrix consisting of the firsts $i$ columns of the matrix $A^{-1}$, which is written in a block form:

$$\left\| \begin{matrix} b_{11} & ... & b_{1,i-1} & b_{1i} \\ b_{21} & ... & b_{2,i-1} & \gamma \\ b_{31} & ... & b_{3,i-1} & 0 \\ ... & ... & ... & ... \\ b_{i1} & ... & b_{i,i-1} & 0 \\ \hline b_{i+1,1} & ... & b_{i+1,i-1} & 0 \\ ... & ... & ... & ... \\ b_{n1} & ... & b_{n,i-1} & 0 \end{matrix} \right\| = \left\| \begin{matrix} M \\ N \end{matrix} \right\|.$$

By Theorem 1.8, the elements $u_1, ..., u_n$, which satisfy the condition 1) can be written as following:

$$\left\| u_1 \quad ... \quad u_n \right\| = \left\| 1 \quad x_2 \quad ... \quad x_n \right\| A^{-1},$$

where $x_i \in R, i = 2, ..., n$. So, to prove this statement it is sufficient to show that there are elements $x_2, ..., x_n$ such that

$$\left\| 1 \quad x_2 \quad ... \quad x_n \right\| \left\| \begin{matrix} M \\ N \end{matrix} \right\| = \left\| q_1 \quad ... \quad q_i \right\|,$$

where

$$(q_1, ..., q_i) = (q_i, \psi) = 1.$$

Suppose that $\gamma = 0$. The matrix $\left\| \begin{matrix} M \\ N \end{matrix} \right\|$ is primitive. It follows that $b_{1i} \in U(R)$, i.e., $u_1 = b_{11}, ..., u_n = b_{1n}$.

Suppose that $\gamma \neq 0$, $N \neq \mathbf{0}$ and $b_{tj}$ is a nonzero element of this matrix, $i+1 \leq t \leq n, 1 \leq j \leq i-1$. Since $(b_{1i}, \gamma) = 1$, then $(b_{1i}, \gamma, \psi b_{tj}) = 1$. By virtue of Property 1.19, there is an element $l$ such that

$$(b_{1i} + \gamma l, \psi b_{tj}) = 1. \tag{1.12}$$

Set $d_{1i} = b_{1i} + \gamma l$. Since $\psi \neq 0$, then $d_{1i} \neq 0$. The following equations follow from (1.12):

$$\begin{aligned} &i)\ (d_{1i}, \psi) = 1, \\ &ii)\ (d_{1i}, b_{tj}) = 1. \end{aligned} \tag{1.13}$$

From the last equality we get

$$(d_{1j}, b_{tj}, d_{1i}) = 1,$$

where $d_{1j} = b_{1j} + b_{2j}l$. According to Property 1.19, there is $m$ such that

$$(d_{1j} + b_{tj}m, d_{1i}) = 1.$$

Taking into account equality (1.13), we conclude that the elements of the first row of the matrix

$$\left( \left\| \begin{matrix} 1 & l & 0 & ... & 0 & m \\ 0 & 1 & 0 & ... & 0 & 0 \\ ... & ... & ... & ... & ... & ... \\ 0 & 0 & 0 & ... & 0 & 1 \end{matrix} \right\| \oplus I_{n-t} \right) A^{-1}$$

satisfy all the requirements of our statement.

Let $N = \mathbf{0}$ or $i = n$ (in this case the matrix $N$ is empty). The invertibility of the matrix $A$ implies that $M \in \mathrm{GL}_i(R)$. Therefore, $(b_{1i}, \gamma) = 1$. It follows that $(b_{1i}, \gamma, \psi) = 1$. As in in the previous case, there is $r$ that $(b_{1i} + \gamma r, \psi) = 1$. Thus, the elements of the first row of the matrix

$$\left( \left\| \begin{matrix} \in & 1 & r & \in \\ & 0 & 1 & \end{matrix} \right\| \oplus I_{n-2} \right) A^{-1}$$

will be desirable. The theorem is proved. $\qquad\square$

## 1.5. Greatest common divisor of matrices

If $A = BC$, then the matrix $B$ is called a left divisor of the matrix $A$ and the matrix $A$ is called a right multiple of the matrix $B$. If $A = DA_1$ and $B = DB_1$, then the matrix $D$ is called a left common divisor of the matrices $A$ and $B$. In addition, if the matrix $D$ is a right multiple of each left common divisor of the matrices $A$ and $B$, then the matrix $D$ is called the left greatest common divisor of the matrices $A$ and $B$ and denoted by $(A, B)_l$. Lets prove the correctness of this concept in a matrices rings over Bezout rings.

Let $A, B \in M_n(R)$. Consider an $n \times 2n$ matrix $\|A \quad B\|$. According to Theorem 1.5, there is an invertible $2n \times 2n$ matrix $U$ such that

$$\|A \quad B\| U = \|D \quad \mathbf{0}\|,$$

where $D$ is an $n \times n$ matrix. Based on this result, we propose a method for finding of the left g.c.d. of the matrices $A, B$.

**Theorem 1.10.** *[52] Let $A, B \in M_n(R)$ and*

$$\|A \quad B\| U = \|D \quad \mathbf{0}\|,$$

*where $U$ is an invertible matrix. Then $D = (A, B)_l$ moreover there exist matrices $P, Q$, that*

$$(A, B)_l = AP + BQ.$$

**Proof.** Write the matrix $U$ in a block form

$$U = \left\| \begin{matrix} P & S \\ Q & T \end{matrix} \right\|,$$

where each block has order $n$. Then

$$\|A \quad B\| U = \|A \quad B\| \left\| \begin{matrix} P & S \\ Q & T \end{matrix} \right\| = \|D \quad \mathbf{0}\|.$$

It follows that

$$D = AP + BQ. \tag{1.14}$$

Set

$$U^{-1} = \left\| \begin{matrix} P_1 & S_1 \\ Q_1 & T_1 \end{matrix} \right\|.$$

Then

$$\|A \quad B\| = \|D \quad \mathbf{0}\| U^{-1} = \|D \quad \mathbf{0}\| \left\| \begin{matrix} P_1 & S_1 \\ Q_1 & T_1 \end{matrix} \right\|.$$

It means that $A = DP_1$, $B = DQ_1$. That is, the matrix $D$ is the left common divisor of the matrices $A, B$. Suppose that $D_1$ is another left common divisor of the matrices $A, B$, i.e., $A = D_1 M$, $B = D_1 N$. In view of (1.14), we have

$$D = D_1 MP + D_1 NQ = D_1(MP + NQ).$$

That is, the matrix $D_1$ is the left divisor of the matrix $D$. Therefore, $D$ is the left g.c.d. of the matrices $A, B$. $\square$

In order to prove that the left g.c.d. of matrices is defined uniquely up to the right associativity the following result is required.

**Theorem 1.11.** *Let the matrices $A, B$ are left divisors of each other. Then they are right associated.*

**Proof.** Let $\det A \neq 0$. Noting that $A = BC$, we conclude that $\det B \neq 0$. Since $B = AD$, then

$$A = BC = A(DC) \Rightarrow \det A = \det A \det(DC).$$

Therefore, $\det(DC) = 1$. That is, the matrices $DC$, $D$, $C$ are invertible. Thus, the matrices $A, B$ are right associated.

Suppose that $\det A = 0$. Since $B = AD$, so the matrix $B$ is also singular. Consider the second order matrices. By Theorem 1.5, there are such invertible matrices $U, V$ that

$$AU = \left\|\begin{matrix} a_1 & 0 \\ a_2 & 0 \end{matrix}\right\|, \quad BV = \left\|\begin{matrix} b_1 & 0 \\ b_2 & 0 \end{matrix}\right\|.$$

Since $A = BC$, then

$$AU = (BV)(V^{-1}CU),$$

i.e.,

$$\left\|\begin{matrix} a_1 & 0 \\ a_2 & 0 \end{matrix}\right\| = \left\|\begin{matrix} b_1 & 0 \\ b_2 & 0 \end{matrix}\right\| \left\|\begin{matrix} t_{11} & t_{12} \\ t_{21} & t_{22} \end{matrix}\right\|, \tag{1.15}$$

where $\|t_{ij}\| = V^{-1}CU$. It follows that $b_i t_{12} = 0$, and $a_i = b_i t_{11}$, $i = 1, 2$. On the other hand, $B = AM$. This implies that $b_i = a_i k_{11}$, $i = 1, 2$. It means that there is $e \in U(R)$ such that $a_i = b_i e$. The matrix $BV$ is nonzero, therefore, there are at least one nonzero among the elements $b_1, b_2$. Consequently, $t_{12} = 0$. Thus, the equality (1.15) has the form

$$\left\|\begin{matrix} a_1 & 0 \\ a_2 & 0 \end{matrix}\right\| = \left\|\begin{matrix} b_1 & 0 \\ b_2 & 0 \end{matrix}\right\| \left\|\begin{matrix} e & 0 \\ t_{21} & t_{22} \end{matrix}\right\|.$$

In particular, the equality

$$\left\|\begin{matrix} a_1 & 0 \\ a_2 & 0 \end{matrix}\right\| = \left\|\begin{matrix} b_1 & 0 \\ b_2 & 0 \end{matrix}\right\| \left\|\begin{matrix} e & 0 \\ 0 & 1 \end{matrix}\right\|$$

is satisfied. Multiplying this equality on the right by the matrix $U^{-1}$, we get

$$A = B \left( V \left\|\begin{matrix} e & 0 \\ 0 & 1 \end{matrix}\right\| U^{-1} \right).$$

Therefore, the matrices $A, B$ are right associated.

Suppose that our assertion is true for all matrices of an order less than $n$. And let $A, B$ are matrices of the order $n$. By transposing the rows (which is equivalent to left-multiplying by some an invertible matrix $T$) we make the first row of the matrix $TA$ nonzero. Let $TA = \|a_{ij}\|_1^n$, $TB = \|b_{ij}\|_1^n$. By virtue of Theorem 1.3, there are invertible matrices $P, Q$ such that

$$\left\| \begin{matrix} a_{11} & \dots & a_{1n} \end{matrix} \right\| P = \left\| \begin{matrix} a & 0 & \dots & 0 \end{matrix} \right\|,$$

$$\left\| b_{11} \quad ... \quad b_{1n} \right\| Q = \left\| b \quad 0 \quad ... \quad 0 \right\|,$$

where $a \neq 0$, $b \neq 0$. Reasoning similarly to the above, we show that

$$T(AP) = \left\| \begin{array}{c|c} a & \mathbf{0} \\ \hline A_{21} & A_{22} \end{array} \right\| = \left\| \begin{array}{c|c} b & \mathbf{0} \\ \hline B_{21} & B_{22} \end{array} \right\| \left\| \begin{array}{c|c} f & \mathbf{0} \\ \hline F_{21} & F_{22} \end{array} \right\|,$$

$$T(BQ) = \left\| \begin{array}{c|c} b & \mathbf{0} \\ \hline B_{21} & B_{22} \end{array} \right\| = \left\| \begin{array}{c|c} a & \mathbf{0} \\ \hline A_{21} & A_{22} \end{array} \right\| \left\| \begin{array}{c|c} f & \mathbf{0} \\ \hline F_{21} & F_{22} \end{array} \right\|,$$

where $f \in U(R)$. These implies that the matrices $A_{22}$ and $B_{22}$ are left divisors of each other. According to the induction assumption $A_{22} = B_{22}K$, where $K \in \mathrm{GL}_{n-1}(R)$. Hence,

$$TAP = \left\| \begin{array}{cc} a & \mathbf{0} \\ A_{21} & A_{22} \end{array} \right\| = \left\| \begin{array}{cc} b & \mathbf{0} \\ B_{21} & B_{22} \end{array} \right\| \left\| \begin{array}{cc} e & \mathbf{0} \\ F_{21} & K \end{array} \right\|.$$

Left-multiplying this equality by $T^{-1}$ and right-multiplying by $P^{-1}$, we get

$$A = B \left( Q \left\| \begin{array}{cc} e & \mathbf{0} \\ F_{21} & K \end{array} \right\| P^{-1} \right).$$

Noting that

$$Q \left\| \begin{array}{cc} e & \mathbf{0} \\ F_{21} & K \end{array} \right\| P^{-1} \in \mathrm{GL}_n(R),$$

we complete the proof. $\qquad\square$

Similarly, if the matrices $A, B$ are right divisors of each other, they are left associated.

**Theorem 1.12.** *The left g.c.d. of matrices over a Bezout ring is defined uniquely up to right associativity.*

**Proof.** Suppose that $D_1, D_2$ are left g.c.d. of the matrices $A, B$. Then they are left dividers of each other. Based on Theorem 1.11, the matrices $D_1$ and $D_2$ are right associated. $\qquad\square$

The proved theorem states that $M_n(R)$ is a left Bezout ring. Having considered the matrix $\left\| \begin{array}{c} A \\ B \end{array} \right\|$, we similarly show that $M_n(R)$ is a right Bezout ring. Thus, we have proved the Theorem.

**Theorem 1.13.** *If $R$ is a Bezout ring, so is $M_n(R)$.* $\qquad\square$

Note that left and right g.c.d. of matrices, in general, are not related. The reason for this lies in the structure of the left and right divisors of matrices, which will be covered in Chapter 5. Now let's just show it by an example.

**Example 1.2.** Consider the matrices

$$A = \left\| \begin{array}{cc} 1 & 0 \\ 0 & 6 \end{array} \right\|, \quad B = \left\| \begin{array}{cc} 1 & 0 \\ 1 & 8 \end{array} \right\|.$$

Then

$$\left\| A \quad B \right\| U = \left\| \begin{array}{cc|cc} 1 & 0 & 0 & 0 \\ 0 & 1 & 0 & 0 \end{array} \right\| = \left\| I \quad \mathbf{0} \right\|,$$

where

$$U = \left\| \begin{array}{cc|cc} 1 & 6 & -1 & 8 \\ 0 & 1 & 0 & 0 \\ \hline 0 & -6 & 1 & -8 \\ 0 & 0 & 0 & 1 \end{array} \right\|$$

is an invertible matrix. That is, the matrices $A, B$ are relatively prime on the left. However, they have the right g.c.d. $\mathrm{diag}(1, 2)$. $\diamond$

## 1.6. Zero divisors in matrices rings

Zero divisors play an important role in the matrix rings. A nonzero matrix $A$ is called a left zero divisor if there exists $M \neq \mathbf{0}$ such that $AM = \mathbf{0}$. Similarly, matrix $B$ is called a right zero divisor if there exists $N \neq \mathbf{0}$ such that $NB = \mathbf{0}$. A matrix that is a left and a right zero divisor is simply called a zero divisor.

**Theorem 1.14.** *In order that the matrix $A$ to be a zero divisor it is necessary and sufficient that $\det A = 0$.*

**Proof. Necessity.** Let $\det A \neq 0$ and $AM = \mathbf{0}$. There is an invertible matrix $U$ such that

$$AU = \left\| \begin{array}{cccc} \alpha_1 & 0 & ... & 0 \\ b_{21} & \alpha_2 & & 0 \\ & & \ddots & \\ b_{n1} & ... & b_{n.n-1} & \alpha_n \end{array} \right\|$$

is the Hermite form of the matrix $A$. The matrix $A$ is nonsingular. It follows that $\alpha_i \neq 0$, $i = 1, ..., n$. Then

$$\mathbf{0} = AM = (AU)(U^{-1}M) = \left\| \begin{array}{cccc} \alpha_1 & 0 & ... & 0 \\ b_{21} & \alpha_2 & & 0 \\ & & \ddots & \\ b_{n1} & ... & b_{n.n-1} & \alpha_n \end{array} \right\| \left\| v_{ij} \right\|_1^n, \qquad (1.16)$$

where $U^{-1}M = \left\| v_{ij} \right\|_1^n$. The equality (1.16) implies that

$$\alpha_1 \left\| v_{11} \quad ... \quad v_{1n} \right\| = \mathbf{0}.$$

Hence,

$$\left\| v_{11} \quad ... \quad v_{1n} \right\| = \mathbf{0}.$$

Considering sequentially the 2nd, etc. $n$th row of Product (1.16), we get $U^{-1}M = \mathbf{0}$. It follows that $M = \mathbf{0}$. Therefore, the matrix $A$ is not a zero divisor.

**Sufficiency**. Let $\det A = 0$. Hence, rang $A = r < n$. By Theorem 1.6, there are such invertible matrices $P, Q$, that

$$PAQ = \left\| \begin{array}{cc} D & \mathbf{0} \\ \mathbf{0} & \mathbf{0} \end{array} \right\|,$$

where $D$ is a nonsingular matrix of the order $r$. Then for an arbitrary $S$ matrix

$$S = \left\|\begin{matrix} \mathbf{0} & \mathbf{0} \\ S_1 & S_2 \end{matrix}\right\|,$$

where $S_1$ is a $(n-r) \times r$ matrix, the equality

$$PAQS = \mathbf{0} \Rightarrow A(QS) = \mathbf{0}$$

is fulfilled. That is, the matrix $A$ is a left zero divisor.

Similarly, we show that the condition $NA = \mathbf{0}$, where $N \neq \mathbf{0}$ is equivalent to the condition $\det A = 0$. Consequently, the matrix is at the same time a left and a right zero divisor or a nonsingular matrix. □

Denote by $Ann^r(D)$ the set of right annihilators of the matrix $D$:

$$Ann^r(D) = \{Q \in M_n(R) | DQ = \mathbf{0}\}.$$

Obviously, $Ann^r(D)$ is the right ideal of the ring $M_n(R)$. By virtue of Theorem 1.14, if the matrix $D$ is nonsingular then $Q = \mathbf{0}$. That is $Ann^r(D) = \{\mathbf{0}\}$. We show the form of matrices of the set $Ann^r(D)$, when the matrix $D$ is singular.

**Theorem 1.15.** *Let $D$ be the $n \times n$ matrix of rank $r < n$. Then $Ann^r(D)$ consists of all matrices of the form*

$$U \left\|\begin{matrix} \mathbf{0} \\ B \end{matrix}\right\|, \text{ where } DU = \left\|D_1 \quad \mathbf{0}\right\|,$$

*$D_1$ is the $n \times r$ matrix of rank $r$, $B$ is an arbitrary $(n-r) \times n$ matrix.*

**Proof.** Since

$$D = \left\|D_1 \quad \mathbf{0}\right\| U^{-1},$$

then

$$DU \left\|\begin{matrix} \mathbf{0} \\ B \end{matrix}\right\| = \left\|D_1 \quad \mathbf{0}\right\| U^{-1} U \left\|\begin{matrix} \mathbf{0} \\ B \end{matrix}\right\| = \mathbf{0}.$$

Hence, all matrices of the form $U \left\|\begin{matrix} \mathbf{0} \\ B \end{matrix}\right\|$ belong to $Ann^r(D)$.

Let $DQ = \mathbf{0}$. That is

$$\left\|D_1 \quad \mathbf{0}\right\| U^{-1} Q = \left\|D_1 \quad \mathbf{0}\right\| V,$$

where $V = U^{-1}Q$. Write the matrix $V$ in a block form: $V = \left\|\begin{matrix} V_1 \\ V_2 \end{matrix}\right\|$, where $V_1$ is $r \times n$ matrix. Since $D_1$ has a maximum rank, it contains a nonsingular submatrix of order $r$. Without loss of generality, we can assume that this matrix consists of the first $r$ of its rows. We denote it by $B_1$. Then

$$DQ = \left\|\begin{matrix} B_1 & \mathbf{0} \\ B_2 & \mathbf{0} \end{matrix}\right\| \left\|\begin{matrix} V_1 \\ V_2 \end{matrix}\right\| = \mathbf{0}.$$

**41**

Hence, $B_1 V_1 = \mathbf{0}$. The matrix $B_1$ is nonsingular, therefore $V_1 = \mathbf{0}$. Obviously, the equality

$$\left\| \begin{matrix} B_1 & \mathbf{0} \\ B_2 & \mathbf{0} \end{matrix} \right\| \left\| \begin{matrix} \mathbf{0} \\ V_2 \end{matrix} \right\| = \mathbf{0}$$

is executed without any restrictions on the matrix $V_2$. That is, $V_2$ is an arbitrary $(n-r) \times n$ matrix. Noting that

$$Q = UV = U \left\| \begin{matrix} \mathbf{0} \\ V_2 \end{matrix} \right\|,$$

we complete the proof. $\qquad\qquad\qquad\qquad\qquad\qquad\qquad\qquad\qquad\qquad \square$

## 1.7. Matrix Bezout's identity

We say that the pair $(U, V)$ reduces the matrices $A, B$ to $(A, B)_l$ provide that
$$AU + BV = (A, B)_l = D.$$

Denote the set of all such pairs by $\mathbf{U}^r_{A,B}$.

**Theorem 1.16.** *Let*
$$\left\| A \quad B \right\| \left\| \begin{matrix} K_{11} & K_{12} \\ K_{21} & K_{22} \end{matrix} \right\| = \left\| D \quad \mathbf{0} \right\|,$$

*where* $\left\| \begin{matrix} K_{11} & K_{12} \\ K_{21} & K_{22} \end{matrix} \right\|$ *is an invertible matrix. The set* $\mathbf{U}^r_{A,B}$ *consists of all pairs of the form*
$$(K_{11} + K_{11}Q + K_{12}P, \ K_{21} + K_{21}Q + K_{22}P),$$

*where* $Q \in Ann^r(D)$, *and* $P \in M_n(R)$.

**Proof**. Denote by $\mathbf{V}^r_{A,B}$ the set

$$\{(K_{11} + K_{11}Q + K_{12}P, \ K_{21} + K_{21}Q + K_{22}P)\},$$

where $Q \in Ann^r(D)$, and $P \in M_n(R)$. Then

$$A(K_{11} + K_{11}Q + K_{12}P) + B(K_{21} + K_{21}Q + K_{22}P) =$$
$$= (AK_{11} + BK_{21}) + (AK_{11} + BK_{21})Q + (AK_{12} + BK_{22})P =$$
$$= D + DQ + \mathbf{0}P = D.$$

Hence, $\mathbf{V}^r_{A,B} \subseteq \mathbf{U}^r_{A,B}$.

Suppose that $(U, V) \in \mathbf{U}^r_{A,B}$ and

$$\left\| \begin{matrix} K_{11} & K_{12} \\ K_{21} & K_{22} \end{matrix} \right\|^{-1} = \left\| \begin{matrix} L_{11} & L_{12} \\ L_{21} & L_{22} \end{matrix} \right\|.$$

Then

$$\left\| A \quad B \right\| = \left\| D \quad \mathbf{0} \right\| \left\| \begin{matrix} L_{11} & L_{12} \\ L_{21} & L_{22} \end{matrix} \right\|.$$

It follows that
$$A = DL_{11}, \ B = DL_{12}.$$

Consider the product
$$\begin{Vmatrix} L_{11} & L_{12} \\ L_{21} & L_{22} \end{Vmatrix} \begin{Vmatrix} U & K_{12} \\ V & K_{22} \end{Vmatrix} = \begin{Vmatrix} T & \mathbf{0} \\ P & I \end{Vmatrix}, \tag{1.17}$$

where $I$ is an identity matrix and
$$T = L_{11}U + L_{12}V, \ P = L_{21}U + L_{22}V.$$

Then
$$DT = (DL_{11})U + (DL_{12})V = AU + BV = D,$$

e.i.,
$$DT = D \Rightarrow D(T - I) = \mathbf{0} \Rightarrow T - I = Q \in Ann^r(D).$$

Therefore, $T = I + Q$. The equality (1.17) implies that
$$\begin{Vmatrix} U & K_{12} \\ V & K_{22} \end{Vmatrix} = \begin{Vmatrix} K_{11} & K_{12} \\ K_{21} & K_{22} \end{Vmatrix} \begin{Vmatrix} I + Q & \mathbf{0} \\ P & I \end{Vmatrix} =$$
$$= \begin{Vmatrix} K_{11} + K_{11}Q + K_{12}P & K_{12} \\ K_{21} + K_{21}Q + K_{22}P & K_{22} \end{Vmatrix}.$$

Hence,
$$U = K_{11} + K_{11}Q + K_{12}P, \ V = K_{21} + K_{21}Q + K_{22}P.$$

It means that $\mathbf{U}^r_{A,B} \subseteq \mathbf{V}^r_{A,B}$. Therefore, $\mathbf{U}^r_{A,B} = \mathbf{V}^r_{A,B}$. $\qquad \square$

**Theorem 1.17.** *Let* $A, B \in M_n(R)$. *Then there is an invertible matrix* $\begin{Vmatrix} U & M \\ V & N \end{Vmatrix}$ *such that*
$$\begin{Vmatrix} A & B \end{Vmatrix} \begin{Vmatrix} U & M \\ V & N \end{Vmatrix} = \begin{Vmatrix} D & \mathbf{0} \end{Vmatrix},$$

*moreover*
$$Ann^r(D) \subseteq Ann^r(V). \tag{1.18}$$

**Proof.** Let $D = (A, B)_l$ be a nonsingular matrix. Then $Ann^r(D) = \mathbf{0}$ and the inclusion (1.18) is correct.

Suppose that $(A, B)_l$ be a singular matrix. It follows that $A, B$ are singular matrices. Let rang $B = k < n$. Then there are invertible matrices $P_B, Q_B$ such that
$$P_B B Q_B = \begin{Vmatrix} B_{11} & \mathbf{0} \\ \mathbf{0} & \mathbf{0} \end{Vmatrix},$$

where $B_{11}$ is $k \times k$ nonsingular matrix.

Let rang $A = r < n$. Hence, rang $P_B A = r$. Then there exists an invertible matrix $Q_A$ such that
$$(P_B A) Q_A = \begin{Vmatrix} A_{11} & \mathbf{0} \\ A_{21} & \mathbf{0} \end{Vmatrix},$$

where
$$\operatorname{rang} \left\| \begin{matrix} A_{11} \\ A_{21} \end{matrix} \right\| = r,$$

and $A_{11}$ is $k \times r$ matrix.

Suppose that $\operatorname{rang} A_{21} = t \leqslant r$. There exists an invertible matrix $P_{21}$ such that
$$P_{21} A_{21} = \left\| \begin{matrix} A'_{21} \\ \mathbf{0} \end{matrix} \right\|,$$

where $A'_{21}$ is $t \times r$ matrix, and $\operatorname{rang} A'_{21} = t$. Then

$$\left\| A \quad B \right\| = P_B^{-1} \left\| \begin{matrix} I_k & \mathbf{0} \\ \mathbf{0} & P_{21}^{-1} \end{matrix} \right\| \underbrace{\left\| \begin{matrix} A_{11} & \mathbf{0} & B_{11} & \mathbf{0} \\ A'_{21} & \mathbf{0} & \mathbf{0} & \mathbf{0} \\ \mathbf{0} & \mathbf{0} & \mathbf{0} & \mathbf{0} \end{matrix} \right\|}_{F} \left\| \begin{matrix} Q_A^{-1} & \mathbf{0} \\ \mathbf{0} & Q_B^{-1} \end{matrix} \right\|.$$

The matrix $B_{11}$ is nonsingular and g.c.d. $t$th order minors of the matrix $A'_{21}$ is not equal to zero. It follows that
$$\operatorname{rang} \left\| A \quad B \right\| = \operatorname{rang} B_{11} + \operatorname{rang} A'_{21}.$$

That is
$$\operatorname{rang} \left\| A \quad B \right\| = k + t.$$

Complement the matrices $\left\| \begin{matrix} A_{11} \\ A'_{21} \end{matrix} \right\|$, $B_{11}$ by null blocks to matrices of the order $k + t$ of the forms
$$\left\| \begin{matrix} A_{11} & \mathbf{0} \\ A'_{21} & \mathbf{0} \end{matrix} \right\| = A', \quad \left\| \begin{matrix} B_{11} & \mathbf{0} \\ \mathbf{0} & \mathbf{0} \end{matrix} \right\| = B'.$$

Let $(A', B')_l = D'$. Since
$$\operatorname{rang} \left\| A' \quad B' \right\| = \operatorname{rang} \left\| A \quad B \right\| = k + t,$$

then the matrix $D'$ is nonsingular. Then there is an invertible matrix $\left\| \begin{matrix} U' & M' \\ V' & N' \end{matrix} \right\|$ such that
$$\left\| A' \quad B' \right\| \left\| \begin{matrix} U' & M' \\ V' & N' \end{matrix} \right\| = \left\| D' \quad \mathbf{0} \right\|.$$

Write the matrix $F$ in the form
$$F = \left\| \begin{matrix} A' & \mathbf{0} & B' & \mathbf{0} \\ \mathbf{0} & \mathbf{0} & \mathbf{0} & \mathbf{0} \end{matrix} \right\|.$$

Then
$$F \underbrace{\left\| \begin{matrix} U' & \mathbf{0} & M' & \mathbf{0} \\ \mathbf{0} & I_{n-k-t} & \mathbf{0} & \mathbf{0} \\ V' & \mathbf{0} & N' & \mathbf{0} \\ \mathbf{0} & \mathbf{0} & \mathbf{0} & I_{n-k-t} \end{matrix} \right\|}_{G} = \left\| \begin{matrix} D' & \mathbf{0} & \mathbf{0} & \mathbf{0} \\ \mathbf{0} & \mathbf{0} & \mathbf{0} & \mathbf{0} \end{matrix} \right\|,$$

where $I_{n-k-t}$ is an identity matrix of the order $n-k-t$. It is obvious that the matrix $G$ is invertible. Set

$$P = P_B^{-1} \begin{Vmatrix} I_k & \mathbf{0} \\ \mathbf{0} & P_{21}^{-1} \end{Vmatrix}.$$

Therefore,

$$\|A \quad B\| \left( \begin{Vmatrix} Q_A & \mathbf{0} \\ \mathbf{0} & Q_B \end{Vmatrix} \begin{Vmatrix} U' & \mathbf{0} & M' & \mathbf{0} \\ \mathbf{0} & I_{n-k-t} & \mathbf{0} & \mathbf{0} \\ V' & \mathbf{0} & N' & \mathbf{0} \\ \mathbf{0} & \mathbf{0} & \mathbf{0} & I_{n-k-t} \end{Vmatrix} \right) =$$

$$= P \begin{Vmatrix} D' & \mathbf{0} & \mathbf{0} & \mathbf{0} \\ \mathbf{0} & \mathbf{0} & \mathbf{0} & \mathbf{0} \end{Vmatrix}. \tag{1.19}$$

Hence,

$$(A, B)_l = P \begin{Vmatrix} D' & \mathbf{0} \\ \mathbf{0} & \mathbf{0} \end{Vmatrix} = D.$$

The matrix $D'$ is nonsingular. It follows that $Ann^r(D)$ consists of all matrices of the form $\begin{Vmatrix} \mathbf{0} \\ T \end{Vmatrix}$, where $T$ is $(n-k-t) \times n$ matrix. The equality (1.19) implies that

$$A \underbrace{\left( Q_A \begin{Vmatrix} U' & \mathbf{0} \\ \mathbf{0} & I_{n-k-t} \end{Vmatrix} \right)}_{U} + B \underbrace{\left( Q_B \begin{Vmatrix} V' & \mathbf{0} \\ \mathbf{0} & \mathbf{0} \end{Vmatrix} \right)}_{V} = D.$$

Let the matrix $V'$ is nonsingular. Then the form of the matrix $V$ implies that $Ann^r(D) = Ann^r(V)$. Otherwise we have $Ann^r(D) \subseteq Ann^r(V)$. The theorem is proved. □

If $M = AP = BQ$, then the matrix $M$ is called a right common multiple of the matrices $A$ and $B$. Moreover, if the matrix $M$ is a left divisor of every common right multiple of the matrices $A$ and $B$, then the matrix $M$ is called the right least common multiple **l.c.m.** of the matrices $A$ and $B$ (and we write $[A, B]_r$). By analogy with the statement of the Theorem 1.12, right l.c.m. of a matrices over a Bezout ring is defined uniquely up to associativity.

**Theorem 1.18.** *Let $A, B \in M_n(R)$. Then there exists an invertible matrix $F$ such that*

$$\begin{Vmatrix} A & B \\ \mathbf{0} & B \end{Vmatrix} F = \begin{Vmatrix} (A, B)_l & \mathbf{0} \\ * & [A, B]_r \end{Vmatrix}.$$

**Proof.** Let $(A, B)_l = D$. By Theorem 1.17, there exists an invertible matrix $\begin{Vmatrix} U & M \\ V & N \end{Vmatrix}$ such that

$$\|A \quad B\| \begin{Vmatrix} U & M \\ V & N \end{Vmatrix} = \|D \quad \mathbf{0}\|, \tag{1.20}$$

moreover $Ann^r(D) \subseteq Ann^r(V)$. The equality (1.20) implies that $AM + BN = \mathbf{0}$. Therefore, $S = BN = -AM$ is the right common multiple of matrices

**45**

$A, B$. Let $S_1 = AA_1 = BB_1$ is any other right common multiple of matrices $A, B$. We claim that $S$ is the left divisor of the matrix $S_1$.

Consider a pair of matrices $(U - A_1, V + B_1)$. Then

$$A(U - A_1) + B(V + B_1) = (AU + BV) + (BB_1 - AA_1) = D.$$

It means that $(U - A_1, V + B_1) \in \mathbf{U}^r_{A,B}$. By virtue of Theorem 1.16, there exists $Q \in Ann^r(D)$ and $P \in M_n(R)$ such that

$$\begin{aligned} U - A_1 &= U + UQ + MP, \\ V + B_1 &= V + VQ + NP. \end{aligned} \tag{1.21}$$

Since $VQ = \mathbf{0}$ the equality (1.21) implies that $B_1 = NP$. Consequently,

$$S_1 = BB_1 = B(NP) = (BN)P = SP.$$

This yields, $S = BN = [A, B]_r$. Taking into account the equality (1.20), we receive

$$\left\| \begin{matrix} A & B \\ \mathbf{0} & B \end{matrix} \right\| \left\| \begin{matrix} U & M \\ V & N \end{matrix} \right\| = \left\| \begin{matrix} D & \mathbf{0} \\ BV & BN \end{matrix} \right\| = \left\| \begin{matrix} (A,B)_l & \mathbf{0} \\ BV & [A,B]_r \end{matrix} \right\|.$$

The theorem is proved. □

**Corollary 1.5.** *The equality*

$$\det(AB) = \det(A,B)_l \det[A,B]_r.$$

*is fulfilled.* □

The existence of zero divisors in matrix rings leads to the fact that some properties of the Bezout's identity, g.c.d., and l.c.m., which are executed in the ring $R$, are not inherited by the matrix rings over $R$. In particular, the equality

$$\left\| \begin{matrix} A & B \end{matrix} \right\| \left\| \begin{matrix} U & M \\ V & N \end{matrix} \right\| = \left\| \begin{matrix} D & \mathbf{0} \end{matrix} \right\|,$$

where $\left\| \begin{matrix} U & M \\ V & N \end{matrix} \right\|$ is an invertible matrix, it does not follow that the matrix $S = BN$ (see Theorem 1.18) is a right l.c.m. of the matrices $A, B$.

Also, the equality $AU + BV = (A, B)_l$ it does not imply that there exist an invertible matrix of the form $\left\| \begin{matrix} U & M \\ V & N \end{matrix} \right\|$.

Also it does not follow from the equality $A = (A, B)_l A_1$, $B = (A, B)_l B_1$ that $(A_1, B_1)_l = I$.

Let us demonstrate these specific properties in an example.

**Example 1.3.** Consider the matrices

$$A = \left\| \begin{matrix} \alpha & 0 \\ 0 & 0 \end{matrix} \right\|, \quad B = \left\| \begin{matrix} \beta & 0 \\ 0 & 0 \end{matrix} \right\|,$$

where $\alpha u + \beta v = 1$. Then

$$\left\| \begin{matrix} \alpha & 0 \\ 0 & 0 \end{matrix} \middle| \begin{matrix} \beta & 0 \\ 0 & 0 \end{matrix} \right\| \left\| \begin{matrix} u & -\beta \\ 0 & 0 \end{matrix} \middle| \begin{matrix} -\beta r & 0 \\ 1 & 0 \end{matrix} \\ \begin{matrix} v & \alpha \\ 0 & 0 \end{matrix} \middle| \begin{matrix} \alpha r & 0 \\ 0 & 1 \end{matrix} \right\| = \left\| \begin{matrix} 1 & 0 \\ 0 & 0 \end{matrix} \middle| \begin{matrix} 0 & 0 \\ 0 & 0 \end{matrix} \right\|.$$

Hence, $(A, B)_l = \mathrm{diag}(1, 0)$. However, the matrix

$$S = BN = \left\| \begin{matrix} \beta & 0 \\ 0 & 0 \end{matrix} \right\| \left\| \begin{matrix} \alpha r & 0 \\ 0 & 1 \end{matrix} \right\| = \left\| \begin{matrix} \alpha \beta r & 0 \\ 0 & 0 \end{matrix} \right\|$$

will be the right l.c.m. of $A$, $B$ if and only if $r$ is an invertible element of the ring $R$.

Consider the Bezout equation:

$$\left\| \begin{matrix} \alpha & 0 \\ 0 & 0 \end{matrix} \right\| \left\| \begin{matrix} u & 0 \\ 0 & 0 \end{matrix} \right\| + \left\| \begin{matrix} \beta & 0 \\ 0 & 0 \end{matrix} \right\| \left\| \begin{matrix} v & 0 \\ 0 & 0 \end{matrix} \right\| = \left\| \begin{matrix} 1 & 0 \\ 0 & 0 \end{matrix} \right\|.$$

However, the matrix

$$\left\| \begin{matrix} u & 0 \\ 0 & 0 \\ \hline v & 0 \\ 0 & 0 \end{matrix} \right\|$$

is not complemented to the fourth order invertible matrix.

Taking into account form of $(A, B)_l$, we receive $A = (A, B)_l A$, $B = (A, B)_l B$, i.e., $A_1 = A$, $B_1 = B$. Consequently, $(A_1, B_1)_l \neq I$. $\diamond$

We show that the left fractions of the matrix $A, B$ by $(A, B)_l$ can still be chosen so that they will be left relatively prime.

**Lemma 1.2.** *Let $A, B$ be $m \times n$ left relatively prime, $m < n$. Then this matrices can be complement to square left relatively prime matrices of the form*

$$\left\| \begin{matrix} A \\ M \end{matrix} \right\|, \quad \left\| \begin{matrix} B \\ N \end{matrix} \right\|.$$

**Proof.** There are invertible matrices $U, V$ such that

$$AU = \left\| A_1 \quad \mathbf{0} \right\|, \quad BV = \left\| B_1 \quad \mathbf{0} \right\|,$$

where $A_1, B_1$ are $m \times m$ matrices. Since $(A, B)_l = I_m$, there is an invertible matrix $T$ that

$$\left\| A \quad B \right\| T = \left\| I_m \quad \mathbf{0} \right\|.$$

On the other hand

$$\left( \left\| A \quad B \right\| \left\| \begin{matrix} U & \mathbf{0} \\ \mathbf{0} & V \end{matrix} \right\| \right) \left( \left\| \begin{matrix} U^{-1} & \mathbf{0} \\ \mathbf{0} & V^{-1} \end{matrix} \right\| T \right) = \left\| I_m \quad \mathbf{0} \right\|.$$

That is
$$\|A_1 \quad \mathbf{0} \quad B_1 \quad \mathbf{0}\| \, T_1 = \|I_m \quad \mathbf{0}\|,$$
where $T_1 = \left\| \begin{matrix} U^{-1} & \mathbf{0} \\ \mathbf{0} & V^{-1} \end{matrix} \right\| T$. Thus, $(A_1, B_1)_l = I_m$. Consequently, there is an invertible matrix $P = \|P_{ij}\|_1^2$, such that
$$\|A_1 \quad B_1\| \, P = \|I_m \quad \mathbf{0}\|.$$

Thus
$$\|A_1 \quad \mathbf{0} \quad B_1 \quad \mathbf{0}\| \underbrace{\left\| \begin{matrix} P_{11} & \mathbf{0} & P_{12} & \mathbf{0} \\ \mathbf{0} & I_{n-m} & \mathbf{0} & \mathbf{0} \\ P_{21} & \mathbf{0} & P_{22} & \mathbf{0} \\ \mathbf{0} & \mathbf{0} & \mathbf{0} & I_{n-m} \end{matrix} \right\|}_{Q} = \|I_m \quad \mathbf{0}\|.$$

Consider the matrix
$$\left\| \begin{matrix} A_1 & \mathbf{0} & B_1 & \mathbf{0} \\ \mathbf{0} & I_{n-m} & \mathbf{0} & \mathbf{0} \end{matrix} \right\|.$$

Then
$$\left\| \begin{matrix} A_1 & \mathbf{0} & B_1 & \mathbf{0} \\ \mathbf{0} & I_{n-m} & \mathbf{0} & \mathbf{0} \end{matrix} \right\| Q = \left\| \begin{matrix} I_m & \mathbf{0} & \mathbf{0} & \mathbf{0} \\ \mathbf{0} & I_{n-m} & \mathbf{0} & \mathbf{0} \end{matrix} \right\|.$$

It follows that
$$\left( \left\| \begin{matrix} A_1 & \mathbf{0} & B_1 & \mathbf{0} \\ \mathbf{0} & I_{n-m} & \mathbf{0} & \mathbf{0} \end{matrix} \right\| \left\| \begin{matrix} U^{-1} & \mathbf{0} \\ \mathbf{0} & V^{-1} \end{matrix} \right\| \right) \left( \left\| \begin{matrix} U & \mathbf{0} \\ \mathbf{0} & V \end{matrix} \right\| Q \right) = \|I_m \quad \mathbf{0}\|.$$

That is
$$\left\| \begin{matrix} A & B \\ * & * \end{matrix} \right\| \left( \left\| \begin{matrix} U & \mathbf{0} \\ \mathbf{0} & V \end{matrix} \right\| Q \right) = \|I_m \quad \mathbf{0}\|. \qquad \square$$

**Theorem 1.19.** *Let $D = (A, B)_l$. There are $M$, $N$ such that $A = DM$, $B = DN$, where $(M, N)_l = I$.*

**Proof.** Let $A = DM$, $B = DN$, moreover $\det D \neq 0$ i $(M, N)_l = D_1$. Then $M = D_1 M_1$, $N = D_1 N_1$. Hence, $A = DD_1 M_1$, $B = DD_1 N_1$. Since $D$ is the left g.c.d. of matrices $A, B$, then $DD_1$ is its left divisor:
$$D = DD_1 L \Rightarrow D(I - D_1 L) = \mathbf{0}.$$

The matrix $D$ is not a divisor of zero, so $D_1 L = I$. That is $D_1 \in \mathrm{GL}_n(R)$. Hence, $(M, N)_l = I$.

Let $\det D = 0$ and $\mathrm{rang}\, D = r < n$. There are invertible matrices $U, V$ such that
$$UDV = \left\| \begin{matrix} K & \mathbf{0} \\ \mathbf{0} & \mathbf{0} \end{matrix} \right\|,$$

where $K$ is the nonsingular matrix of an order $r$. The equalities $A = DA_1$, $B = DB_1$ imply that

$$UA = UDV(V^{-1}A_1) = \begin{Vmatrix} K & \mathbf{0} \\ \mathbf{0} & \mathbf{0} \end{Vmatrix} \begin{Vmatrix} A_{11} & A_{12} \\ A_{21} & A_{22} \end{Vmatrix}, \qquad (1.22)$$

$$UB = UDV(V^{-1}B_1) = \begin{Vmatrix} K & \mathbf{0} \\ \mathbf{0} & \mathbf{0} \end{Vmatrix} \begin{Vmatrix} B_{11} & B_{12} \\ B_{21} & B_{22} \end{Vmatrix}. \qquad (1.23)$$

Suppose that matrices $\begin{Vmatrix} A_{11} & A_{12} \end{Vmatrix}$, $\begin{Vmatrix} B_{11} & B_{12} \end{Vmatrix}$ have a left common divisor $T$: $\begin{Vmatrix} A_{11} & A_{12} \end{Vmatrix} = TL_1$, $\begin{Vmatrix} B_{11} & B_{12} \end{Vmatrix} = TL_2$. It follows that the matrix $\begin{Vmatrix} KT & \mathbf{0} \\ \mathbf{0} & \mathbf{0} \end{Vmatrix}$ is the left common divisor of matrices $UA$ and $UB$, having the left g.c.d. $\begin{Vmatrix} K & \mathbf{0} \\ \mathbf{0} & \mathbf{0} \end{Vmatrix}$. Hence, $K = KTS$. Since $K$ is not singular, so $TS = I$. That is, $T$ is an invertible matrix. Thus, the matrices $\begin{Vmatrix} A_{11} & A_{12} \end{Vmatrix}$, $\begin{Vmatrix} B_{11} & B_{12} \end{Vmatrix}$ are left relatively prime. On the basis of Lemma 1.2, there are left relatively prime square matrices

$$\begin{Vmatrix} A_{11} & A_{12} \\ S_1 & S_2 \end{Vmatrix} = A_1', \quad \begin{Vmatrix} B_{11} & B_{12} \\ L_1 & L_2 \end{Vmatrix} = B_1'.$$

The equalities (1.22) and (1.23) imply that the equalities

$$UA = UDV(V^{-1}A_1) = \begin{Vmatrix} K & \mathbf{0} \\ \mathbf{0} & \mathbf{0} \end{Vmatrix} \underbrace{\begin{Vmatrix} A_{11} & A_{12} \\ S_1 & S_2 \end{Vmatrix}}_{A_1'} = (UDV)A_1',$$

$$UB = UDV(V^{-1}B_1) = \begin{Vmatrix} K & \mathbf{0} \\ \mathbf{0} & \mathbf{0} \end{Vmatrix} \underbrace{\begin{Vmatrix} B_{11} & B_{12} \\ B_{21} & B_{22} \end{Vmatrix}}_{B_1'} = (UDV)B_1'$$

are correct. Left-multiply these equations by $U^{-1}$, we get $A = D(VA_1')$, $B = D(VB_1')$. Since $(A_1', B_1')_l = I$, then $(VA_1', VB_1')_l = I$. The Theorem is proved. $\qquad \square$

**Theorem 1.20.** *The matrix $(A, B)_l$ can be chosen such that*

$$[A, B]_r = BA_1, \ \ where \ A = (A, B)_l A_1.$$

**Proof.** According to Theorem 1.18, there exists an invertible matrix $U = \|U_{ij}\|$ such that

$$\begin{Vmatrix} A & B \\ \mathbf{0} & B \end{Vmatrix} \begin{Vmatrix} U_{11} & U_{12} \\ U_{21} & U_{22} \end{Vmatrix} = \begin{Vmatrix} (A, B)_l & \mathbf{0} \\ * & [A, B]_r \end{Vmatrix}.$$

Thus, $BU_{22} = [A, B]_r = M$. Set $U^{-1} = V = \|V_{ij}\|_1^2$. Then

$$\left\| \begin{matrix} A & B \\ 0 & B \end{matrix} \right\| = \left\| \begin{matrix} D & \mathbf{0} \\ * & M \end{matrix} \right\| \left\| \begin{matrix} V_{11} & V_{12} \\ V_{21} & V_{22} \end{matrix} \right\|,$$

where $D = (A, B)_l$. It follows that $A = DV_{11}$. Since $U^{-1} = V$, on the basis of Corollary 3.4* (c. 103) the matrices $U_{22}$ and $V_{11}$ are equivalent. That is $U_{22} = KV_{11}S$, where $K, S \in \mathrm{GL}_n(R)$. Then

$$M = BU_{22} = BKV_{11}S \Rightarrow MS^{-1} = M_1 = B(KV_{11}) = BA_1,$$

where $A_1 = KV_{11}$. The equalities

$$A = DV_{11} = D(K^{-1}K)V_{11} = (DK^{-1})(KV_{11}) = D_1A_1$$

are also fulfilled. To finish the proof it is enough to note that the matrices $M, M_1$ and $D, D_1$ are right associated. That is, $M_1$ is the right l.c.m. and $D_1$ is the left g.c.d. of the matrices $A, B$. □

# ELEMENTARY DIVISOR RINGS

*In the previous chapter, one-sided transformations of matrices over Bezout rings were studied. In this case, a canonical form was established for such transformations. Further research will focus on finding the same form for bilateral matrix transformations.*

## 2.1. Smith normal form

Matrices $A$ and $B$ are called **equivalent** provide that $A = = PBQ$ for some invertible $P$ and $Q$.

By a **diagonal matrix** we mean a matrix in which all elements outside the main diagonal are zeros. If each previous diagonal element divides the next one, then we will call it the **d-matrix**. This matrix can also be rectangular.

**Theorem 2.1.** *Let $R$ be a Bezout ring of stable range $1.5$. Then each $m \times n$ matrix $A$ over $R$ is equivalent to the d-matrix $\Phi = \mathrm{diag}(\varphi_1, ..., \varphi_k)$, $k = \min(m, n)$, moreover $\varphi_i \in Z(R)$, $i = 1, ..., k$, and the matrix $\Phi$ is unique in the class matrices that equivalent to $A$.*

**Proof**. At first, suppose that $A$ be a nonsingular $n \times n$ matrix. Denote by $B_1$ a submatrix of the matrix $A$, consist of its first two rows. By Theorem 1.5, there exists an invertible matrix $V_1$ such that

$$B_1 V_1 = \left\| \begin{matrix} b_{11} & 0 & 0 & ... & 0 \\ b_{21} & b_{22} & 0 & ... & 0 \end{matrix} \right\|,$$

where $b_{11}$, $b_{22} \neq 0$. There is $r$ such that

$$(b_{21} + b_{11}r, b_{22}) = (b_{21}, b_{11}, b_{22}) = \beta.$$

Then

$$\underbrace{\left\| \begin{matrix} 1 & 0 \\ r & 1 \end{matrix} \right\|}_{U_1} B_1 V_1 = \left\| \begin{matrix} b_{11} & 0 & 0 & ... & 0 \\ b'_{21} & b_{22} & 0 & ... & 0 \end{matrix} \right\|,$$

where $b'_{21} = b_{21} + b_{11}r$. Since $(b'_{21}, b_{22}) = \beta$, there are $u, v$ such that

$$ub'_{21} + vb_{22} = \beta.$$

Therefore,

$$U_1 B_1 V_1 \left( \left\| \begin{matrix} u & -\dfrac{b_{22}}{\beta} \\ v & \dfrac{b'_{21}}{\beta} \end{matrix} \right\| \oplus I_{n-2} \right) = \left\| \begin{matrix} b_{11}u & -b_{11}\dfrac{b_{22}}{\beta} & 0 & ... & 0 \\ \beta & 0 & 0 & ... & 0 \end{matrix} \right\| = B_2,$$

where $I_{n-2}$ is an identity matrix of an order $n-2$. Thus,

$$\left\| \begin{matrix} 0 & 1 \\ -1 & \dfrac{b_{11}}{\beta}u \end{matrix} \right\| B_2 = \left\| \begin{matrix} \beta & 0 & 0 & ... & 0 \\ 0 & b_{11}\dfrac{b_{22}}{\beta} & 0 & ... & 0 \end{matrix} \right\| = B_3.$$

Since $\beta | b_{11}$, so $B_3$ is $d$-matrix. By Theorem 1.4, multiplication of matrix by invertible one does not change g.c.d. of its elements. Therefore, $\beta$ is g.c.d. all elements of the first two rows of the $A$ matrix. It is obvious that

$$A \sim \left\| \begin{matrix} \beta & 0 & 0 & ... & 0 \\ 0 & b_{11}\dfrac{b_{22}}{\beta} & 0 & ... & 0 \\ c_{31} & c_{32} & c_{33} & ... & c_{3n} \\ ... & ... & ... & ... & ... \\ c_{n1} & c_{n2} & c_{n3} & ... & c_{nn} \end{matrix} \right\| = A_1.$$

Then

$$\left( \left\| \begin{matrix} 1 & 0 & 0 \\ 0 & 0 & 1 \\ 0 & 1 & 0 \end{matrix} \right\| \oplus I_{n-3} \right) A_1 \sim \left\| \begin{matrix} \beta & 0 & 0 & ... & 0 \\ c_{31} & c_{32} & c_{33} & ... & c_{3n} \\ 0 & 0 & b_{11}\dfrac{b_{22}}{\beta} & & 0 \\ c_{41} & c_{42} & c_{43} & ... & c_{4n} \\ ... & ... & ... & ... & ... \\ c_{n1} & c_{n2} & c_{n3} & ... & c_{nn} \end{matrix} \right\| = A_2.$$

Denote by $C_1$ the first two rows of this matrix. Reasoning similarly to the above $C_1 \sim \text{diag}(\gamma_1, \gamma_2)$, where $\gamma_1 = (\beta, c_{31}, ..., c_{3n})$ and $\gamma_1 | \gamma_2$. It means that

$$A \sim \left\| \begin{matrix} \gamma_1 & 0 & 0 & ... & 0 \\ 0 & \gamma_2 & 0 & ... & 0 \\ d_{31} & d_{32} & d_{33} & ... & d_{3n} \\ ... & ... & ... & ... & ... \\ d_{n1} & d_{n2} & d_{n3} & ... & d_{nn} \end{matrix} \right\| = A_3.$$

Note that the third row of $A_3$ is obtained from the third row of $A_2$ by multiplying by some invertible matrix. Hence

$$(d_{31}, ..., d_{3n}) = b_{11}\frac{b_{22}}{\beta}.$$

Since $\gamma_1|\beta$ and $\beta|b_{11}$, then $\gamma_1|b_{11}\frac{b_{22}}{\beta}$. This means that $\gamma_1$ is g.c.d. of the first three rows of the matrix $A_3$. Continuing the process described, we eventually get

$$A \sim \begin{Vmatrix} \varphi_1 & 0 & ... & 0 \\ f_{21} & f_{22} & ... & f_{2n} \\ ... & ... & ... & ... \\ f_{n1} & f_{n2} & ... & f_{nn} \end{Vmatrix} = A',$$

where $\varphi_1 \in Z(R)$ and equal to g.c.d. of all matrix $A'$ elements. Then

$$\begin{Vmatrix} 1 & 0 & ... & 0 \\ -\frac{f_{21}}{\varphi_1} & 1 & ... & 0 \\ \vdots & & \ddots & \\ -\frac{f_{n1}}{\varphi_1} & 0 & ... & 1 \end{Vmatrix} A' = \begin{Vmatrix} \varphi_1 & 0 & ... & 0 \\ 0 & f_{22} & ... & f_{2n} \\ ... & ... & ... & ... \\ 0 & f_{n2} & ... & f_{nn} \end{Vmatrix}.$$

Having made the same reasoning with the matrix

$$F = \begin{Vmatrix} f_{22} & ... & f_{2n} \\ ... & ... & ... \\ f_{n2} & ... & f_{nn} \end{Vmatrix},$$

we get

$$F \sim \begin{Vmatrix} \varphi_2 & 0 & ... & 0 \\ 0 & g_{33} & ... & g_{3n} \\ ... & ... & ... & ... \\ 0 & g_{n3} & ... & g_{nn} \end{Vmatrix},$$

where $\varphi_2 \in Z(R)$ and equal to g.c.d. of all matrix $F$. Therefore,

$$A \sim \begin{Vmatrix} \varphi_1 & 0 & 0 & ... & 0 \\ 0 & \varphi_2 & 0 & ... & 0 \\ 0 & 0 & g_{33} & ... & g_{3n} \\ ... & ... & ... & ... & ... \\ 0 & 0 & g_{n3} & ... & g_{nn} \end{Vmatrix}.$$

Since $\varphi_1$ is a divisor of all elements of the matrix $F$, then $\varphi_1|\varphi_2$. Continuing the process described, we show that $A \sim \mathrm{diag}(\varphi_1, ..., \varphi_n)$, $\varphi_i|\varphi_{i+1}$, $i = 1, ..., n-1$.

Suppose that rang $A = k \leq \min(m, n)$. By Theorem 1.6, there are invertible matrices $U, V$ such that

$$UAV = \begin{Vmatrix} M & \mathbf{0} \\ \mathbf{0} & \mathbf{0} \end{Vmatrix},$$

where $M$ is a nonsingular $k \times k$ matrix. According to the just proved, the matrix $M$ is equivalent to $d$-matrix $\mathrm{diag}(\mu_1, ..., \mu_k)$, where $\mu_i \in Z(R)$, $i = 1, ..., k$. Thus,

$$A \sim \mathrm{diag}(\mu_1, ..., \mu_k, 0, ..., 0) = \mathrm{M}_1.$$

Consequently, in the class of matrices equivalent to $A$ there is a $d$-matrix which diagonal elements belong to $Z(R)$.

Suppose

$$A \sim \mathrm{diag}(\nu_1, ..., \nu_p, \ 0, ..., 0) = \mathrm{N}, \quad \nu_i | \nu_{i+1}, \quad \nu_j \in Z(R).$$

Since N is the $d$-matrix, then $\nu_1 \nu_2 ... \nu_i$ divides all $i$th order minors of the matrix N. Hence, $\nu_1 \nu_2 ... \nu_i$ is g.c.d. of $i$th order minors of this matrix. Since $A \sim \mathrm{N}$ and $A \sim \mathrm{M}_1$, then $\mathrm{N} \sim \mathrm{M}_1$. It means that g.c.d. of corresponding minors of these matrices coincide. It follows that

$$\mathrm{rang}\,\mathrm{M}_1 = \mathrm{rang}\,\mathrm{N}.$$

Therefore, $k = p$. Taking into account $\mu_1, \ \nu_1 \in Z(R)$, we get $\mu_1 = \nu_1$. The equality $\mu_1 \mu_2 = \nu_1 \nu_2$ implies that $\mu_2 = \nu_2$. Continuing this process, we will make sure that in the class of matrices equivalent to the matrix $A$, there is only one $d$-matrix with elements in $Z(R)$. $\qquad \square$

The $d$-matrix obtained in this Theorem is called the Smith form of the matrix. It is named after the British mathematician H.J. Smith [2].

While proving Theorem 2.1, we substantially relied on the fact that the ring $R$ has a stable range 1.5. We show that every matrix over the ring $Q$ (see Example 1.1 (p. 22)) of stable rank 2 is also reduced to the Smith form. That is, the class of rings over which matrices have the Smith form is wider than the class of Bezout rings of stable range 1.5.

We will say that the matrix $A$ has the property of canonical diagonal reduction if there exist invertible matrices $P, Q$ such that

$$PAQ = \mathrm{diag}(\varepsilon_1, ..., \varepsilon_k, 0, ..., 0), \tag{2.1}$$

where $\varepsilon_k \neq 0$, $\varepsilon_i | \varepsilon_{i+1}, i = 1, ..., k-1$. The elements $\varepsilon_i$ are called **invariant factors**, and $P, Q$ are **transforming matrices** of the matrix $A$. The set of all matrices $P$, that satisfies the equality (2.1), will be denoted by $\mathbf{P}_A$.

**Theorem 2.2.** *The first invariant factor $\varepsilon_1$ of the Smith form of the matrix $A$ is equal to g.c.d. of its elements. The remaining invariant factors are found by the formula*

$$\varepsilon_i = \frac{\delta_i}{\delta_{i-1}},$$

*where $\delta_i$ is g.c.d. of ith order minors of the matrix $A$, $i = 2, ..., k$, where $k$ is the order of the last nonzero minor.*

**Proof.** Let $A$ be an $m \times n$ matrix and $P, Q$ its transforming matrices. By Theorem 1.3, g.c.d. of the $k$th order minors of the matrix $A$ and $PAQ$ coincide, $i = 1, ..., m$. Since $\varepsilon_1$ is g.c.d. of all elements of the matrix $PAQ$, then it is equal to g.c.d. of all elements of the matrix $A$. G.c.d. of the $i$th order minors of the matrix $PAQ$ is equal to $\varepsilon_1\varepsilon_2...\varepsilon_i$ and coincide with $\delta_i$ be g.c.d. of the corresponding minors of the matrix $A$. It follows that

$$\frac{\delta_i}{\delta_{i-1}} = \frac{\varepsilon_1\varepsilon_2...\varepsilon_i}{\varepsilon_1\varepsilon_2...\varepsilon_{i-1}} = \varepsilon_i,$$

where $i = 2, ..., k$. $\qquad\square$

We say that the matrix $A$ **admits a diagonal reduction** if $A$ is equivalent to the $d$-matrix. If every matrix over $R$ admits a diagonal reduction, we call (due to I. Kaplansky [1]) $R$ **an elementary divisor ring**.

The following Theorem specifies conditions under which the ring is an elementary divisor ring.

**Theorem 2.3.** *[1] If all $2 \times 2$ matrices over $R$ admit a diagonal reduction, then all matrices admit a diagonal reduction and $R$ is an elementary divisor ring.*

**Proof.** Let $A$ be an $m \times n$ matrix. For definiteness, we put $m \leqslant n$. At first consider the case $m = 2$. By Theorem 1.5, there exists an invertible matrix $V_1$ such that

$$AV_1 = \begin{Vmatrix} a & 0 & 0 & ... & 0 \\ b & c & 0 & ... & 0 \end{Vmatrix} = \begin{Vmatrix} A_1 & \mathbf{0} \end{Vmatrix}.$$

The matrix $A_1$ admits a diagonal reduction. Hence, there exist an invertible matrix $P_1, Q_1$ such that

$$P_1 A_1 Q_1 = \mathrm{diag}(\alpha_1, \alpha_2), \quad \alpha_1 | \alpha_2, \quad \alpha_1, \alpha_1 \in Z(R).$$

Then

$$P_1 A_1 V_1 (Q_1 \oplus I_{n-2}) = \begin{Vmatrix} \alpha_1 & 0 & 0 & ... & 0 \\ 0 & \alpha_2 & 0 & ... & 0 \end{Vmatrix}.$$

Consequently, our statement is correct for two-row matrices. Assume its correctness for all $k < m$ and consider the $m \times n$ matrix $A$. Write the matrix $A$ in a block form: $A = \begin{Vmatrix} B_1 \\ B_2 \end{Vmatrix}$, where $B_1$ is the first row of $A$. Since $B_2$ is an $(m-1) \times n$ matrix, there are such invertible matrices $P_2$, $Q_2$ that

$$P_2 B_2 Q_2 = \mathrm{diag}(\beta_1, ..., \beta_{m-1})$$

is the Smith form of the matrix $B_2$ (some $\beta_i$ may be zeros). Then

$$(1 \oplus P_2) \left\| \begin{matrix} B_1 \\ B_2 \end{matrix} \right\| Q_2 = \left\| \begin{matrix} B_1 Q_2 \\ P_2 B_2 Q_2 \end{matrix} \right\| = \left\| \begin{matrix} * & * & * & ... & * \\ \beta_1 & 0 & 0 & ... & 0 \\ 0 & \beta_2 & 0 & ... & 0 \\ ... & ... & ... & ... & ... \\ 0 & 0 & \beta_{m-1} & & 0 \end{matrix} \right\| = B.$$

Write this matrix in a block form: $B = \left\| \begin{matrix} C_1 \\ C_2 \end{matrix} \right\|$, where $C_1$ is the first two rows of the matrix $B$. As shown above, there are such invertible matrices $P_3, Q_3$ that

$$(P_3 \oplus I_{m-2}) B Q_3 = \left\| \begin{matrix} \gamma_1 & 0 & 0 & ... & 0 \\ 0 & \gamma_2 & 0 & ... & 0 \\ * & * & * & ... & * \\ ... & ... & ... & ... & ... \\ * & * & * & ... & * \end{matrix} \right\| = C,$$

where $\gamma_1 | \gamma_2$, $\gamma_1 \in Z(R)$. Noting that $\beta_1$ is g.c.d. of the matrix $P_2 B_2 Q_2$ elements, and the row $\left\| \beta_1 \quad 0 \quad ... \quad 0 \right\|$ is the second row of the matrix $C_1$, we conclude that $\gamma_1$ is g.c.d. elements of the matrix $C$. By elementary transformations of rows of the matrix $C$ we reduce it to the form

$$\left\| \begin{matrix} \gamma_1 & \mathbf{0} \\ \mathbf{0} & D \end{matrix} \right\| = F \sim C.$$

According to the induction assumption, there are such invertible matrices $P_4, Q_4$ that

$$P_4 D Q_4 = \mathrm{diag}(\delta_1, ..., \delta_{m-1})$$

is the Smith form of the matrix $D$. So,

$$C \sim (1 \oplus P_4) F Q_4 = \mathrm{diag}(\gamma_1, \delta_1, ..., \delta_{m-1}).$$

Since $\gamma_1$ is g.c.d. of elements of $C$, then $\gamma_1$ is the divisor of all elements of the resulting matrix. Therefore, $\mathrm{diag}(\gamma_1, \delta_1, ..., \delta_{m-1})$ is the Smith form of the matrix $A$. The theorem is proved. □

Using the results of this theorem, we can formulate other conditions under which $R$ will be the elementary divisor ring.

**Theorem 2.4.** *[1] Necessary and sufficient conditions for $R$ to be an elementary divisor ring are*

*1) all finitely generated ideals of $R$ are principal, i.e., $R$ is a Bezout ring,*

*2) if $(a, b, c) = l$, there exist $p$ and $q$ such that*

$$(pa + qb, \ qc) = 1.$$

**Proof. Necessity.** Any row $\|a \quad b\|$ over $R$ admits reduction to $\|(a,b) \quad 0\|$. This is equivalent to the fact that $R$ is a Bezout ring.

Consider the matrix $A = \left\|\begin{matrix} a & 0 \\ b & c \end{matrix}\right\|$. Since $(a,b,c) = 1$, than $A \sim \mathrm{diag}(1, ac)$. There are invertible matrices $P = \|p_{ij}\|_1^2$, $Q = \|q_{ij}\|_1^2$ such that

$$PAQ = \mathrm{diag}(1, ac).$$

It means that

$$PAQ = \left\|\begin{matrix} p_{11}a + p_{12}b & p_{12}c \\ p_{21}a + p_{22}b & p_{22}c \end{matrix}\right\| \left\|\begin{matrix} q_{11} & q_{12} \\ q_{21} & q_{22} \end{matrix}\right\| = \left\|\begin{matrix} 1 & 0 \\ 0 & ac \end{matrix}\right\|.$$

Whence

$$(p_{11}a + p_{12}b)q_{11} + (p_{12}c)q_{21} = 1.$$

It follows that

$$(p_{11}a + p_{12}b, p_{12}c) = 1.$$

**Sufficiency.** Based on Theorem 2.3, it suffices to show that every $2 \times 2$ matrix admits a diagonal reduction. Let $B$ be an arbitrary $2 \times 2$ matrix over $R$. According to Theorem 1.5, there exists an invertible matrix $U$ such that

$$BU = \left\|\begin{matrix} b_1 & 0 \\ b_2 & b_3 \end{matrix}\right\| = \delta \left\|\begin{matrix} \dfrac{b_1}{\delta} & 0 \\ \dfrac{b_2}{\delta} & \dfrac{b_3}{\delta} \end{matrix}\right\|,$$

where $\delta = (b_1, b_2, b_3)$. Then

$$\left(p\frac{b_1}{\delta} + q\frac{b_2}{\delta}, q\frac{b_3}{\delta}\right) = 1$$

for some $p, q$. It follows that

$$(pb_1 + qb_2, qb_3) = \delta.$$

The equality $(p, q) = 1$ implies existence of an invertible matrix of the form $\left\|\begin{matrix} p & q \\ u & v \end{matrix}\right\|$. Then

$$\left\|\begin{matrix} p & q \\ u & v \end{matrix}\right\| \left\|\begin{matrix} b_1 & 0 \\ b_2 & b_3 \end{matrix}\right\| = \left\|\begin{matrix} pb_1 + qb_2 & qb_3 \\ * & * \end{matrix}\right\| = B_1.$$

There are $m, n$ such that

$$(pb_1 + qb_2)m + (qb_3)n = \delta.$$

So

$$B_1 \left\|\begin{matrix} m & -\dfrac{b_3}{\delta} \\ n & p\dfrac{b_1}{\delta} + q\dfrac{b_2}{\delta} \end{matrix}\right\| \sim \left\|\begin{matrix} \delta & 0 \\ * & \dfrac{b_1 b_3}{\delta} \end{matrix}\right\| \sim \left\|\begin{matrix} \delta & 0 \\ 0 & \dfrac{b_1 b_3}{\delta} \end{matrix}\right\|.$$

The theorem is proved. $\qquad\qquad\square$

**Corollary 2.1.** *Bezout ring $R$ to be an elementary divisor ring if and only if for each triplet of relatively prime elements $a, b, c$ there are $p, q, u, v$ such that*

$$a(pu) + b(qu) + c(qv) = 1. \qquad \square$$

**Remark.** By Theorem 2.1, a Bezout ring of stable range 1.5 is an elementary divisor ring. Proving this fact was cumbersome. Using Theorem 2.4, this process is simplified. Indeed, if $(b, a, c) = 1$, $c \neq 0$, there is $r$ such that $(b + ar, c) = 1$. That is, $p = r$, $q = 1$ (denoting of Theorem 2.4). If $c = 0$, then $p, q$ are obtained from the equality $ap + bq = 1$. Therefore, according to Theorem 2.4, a Bezout ring of stable range 1.5 is an elementary divisor ring.

**Theorem 2.5.** *[11] The ring*

$$Q = \left\{ a_0 + \sum_{i=1}^{\infty} a_i x^i \,|\, a_0 \in \mathbb{Z}, \ a_i \in \mathbb{Q}, \ i \in \mathbb{N} \right\}$$

*is an elementary divisor ring.*

**Proof.** Example 1.1 shows that $Q$ is a Bezout ring. Therefore, the first condition of Theorem 2.4 is held.

The element $x^k e$, $e \in U(R)$ are relatively prime only with units of the ring $Q$. Hence, only such variants of triple of relatively prime elements are possible:

1) $ae_1, be_2, ce_3$, where $a, b, c$ are relatively prime integers, $e_1, e_2, e_3 \in U(R)$;
2) $x^k e_1, be_2, ce_3$, where $b, c$ are relatively prime integers;
3) $x^k e_1, x^p e_2, e_3$.

First we treat the case 1). Without loss of generality, we can assume that $c \neq 0$. An integer ring $\mathbb{Z}$ is an adequate ring. Therefore, based on Property 1.18, $\mathbb{Z}$ has a stable range of 1.5. It means that there exists $r$ such that $(a + rb, c) = 1$ whence

$$((re_1 e_2^{-1})be_2 + ae_1, ce_3) = 1.$$

So $p = re_1 e_2^{-1}$, $q = 1$.

2). Since $(be_2, ce_3) = 1$ there are $p, q$ such that

$$pbe_2 + qce_3 = 1.$$

Thus,

$$(pbe_2 + qce_3, qx^k e_1) = 1.$$

In the third case, it is enough to choose $q$ in the group $U(Q)$. According to Theorem 2.4, a ring $Q$ is an elementary divisor ring. $\qquad \square$

## 2.2. Solution of the matrix equation $AX = B$ over an elementary divisor rings

Each matrix over elementary divisor ring can be represented as $A = P^{-1}\Phi Q^{-1}$, where $\Phi$ is the Smith form, and $P, Q$ are the left and the right transforming matrices for the matrix $A$. We use this fact to find the roots of the linear matrix equation $AX = B$.

**Theorem 2.6.** *Linear matrix equation $AX = B$ is solvable if and only if all the corresponding invariant factors of the matrices $A$ and $\|A \quad B\|$ coincide.*

**Proof. Necessity.** Suppose that the equation $AX = B$ has the solution $C$. That is, $AC = B$. The equality

$$\|A \quad AC\| \left\| \begin{matrix} I & -C \\ \mathbf{0} & I \end{matrix} \right\| = \|A \quad \mathbf{0}\|$$

is fulfilled. Hence, the matrices $\|A \quad B\|$, $\|A \quad \mathbf{0}\|$ are right associates and equivalent. Therefore, the corresponding invariant factors of these matrices coincide.

**Sufficiency.** Let $\Phi = \mathrm{diag}(\varphi_1, ..., \varphi_r, 0, ..., 0)$ is the Smith form of the matrix $A$. Then the Smith form of the matrix $\|A \quad B\|$ will be a matrix $\Phi'$, which differs from the matrix $\Phi$ only by null rows or columns. There are invertible matrices $P, P', Q, Q'$ such that

$$PAQ = \Phi, \quad P' \|A \quad B\| Q' = \Phi'.$$

The equations $AX = B$ and $(PAQ)(Q^{-1}X) = PB$ are equivalent. Consider the equation $\Phi Y = B_1$, where $Y = Q^{-1}X$, $B_1 = PB$. The equalities

$$P \|A \quad B\| \left\| \begin{matrix} Q & \mathbf{0} \\ \mathbf{0} & I \end{matrix} \right\| = \|PAQ \quad PB\| = \|\Phi \quad B_1\|.$$

hold. It follows that the matrices $\|A \quad B\|$, $\|\Phi \quad B_1\|$ are equivalent. Therefore,

$$\|\Phi \quad \mathbf{0}\| \sim \|A \quad \mathbf{0}\| \sim \|A \quad B\| \sim \|\Phi \quad B_1\|.$$

Consequently,

$$\left\| \begin{matrix} \varphi_1 & & & & 0 & ... & 0 \\ & \ddots & & & ... & ... & ... \\ & & \varphi_r & & 0 & ... & 0 \\ & & & 0 & 0 & ... & 0 \\ & & & & \ddots & ... & ... & ... \\ & & & 0 & 0 & ... & 0 \end{matrix} \right\| \sim$$

$$\sim \begin{Vmatrix} \varphi_1 & & & & b_{11} & ... & b_{1n} \\ & \ddots & & & ... & ... & ... \\ & & \varphi_r & & b_{r1} & ... & b_{rn} \\ & & & 0 & b_{r+1,1} & ... & b_{r+1,n} \\ & & & \ddots & ... & ... & ... \\ & & & 0 & b_{m1} & ... & b_{mn} \end{Vmatrix} = \begin{Vmatrix} \Phi & B_1 \end{Vmatrix}.$$

It follows that

$$\begin{Vmatrix} b_{r+1,1} & ... & b_{r+1,n} \\ ... & ... & ... \\ b_{m1} & ... & b_{mn} \end{Vmatrix} = \mathbf{0}.$$

The element $\varphi_1$ is the divisor of all elements of the matrix $\begin{Vmatrix} \Phi & B_1 \end{Vmatrix}$. Hence, $b_{1i} = \varphi_1 d_{1i}$, $i = 1, ..., n$. So

$$\begin{Vmatrix} \Phi & B_1 \end{Vmatrix} \begin{Vmatrix} 1 & & -d_{11} & ... & -d_{1n} \\ & \ddots & & & \\ 0 & & 1 & & 0 \\ & & & \ddots & \\ 0 & & 0 & & 1 \end{Vmatrix} =$$

$$= \begin{Vmatrix} \varphi_1 & & & & 0 & ... & 0 \\ & \varphi_2 & & & b_{21} & ... & b_{2n} \\ & & \ddots & & ... & ... & ... \\ & & & \varphi_r & b_{r1} & ... & b_{rn} \\ & & & 0 & 0 & ... & 0 \\ & & & \ddots & ... & ... & ... \\ & & & 0 & 0 & ... & 0 \end{Vmatrix} = \begin{Vmatrix} \Phi & B_1' \end{Vmatrix}.$$

This matrix has the 2nd order minors of the form $\varphi_1 b_{2i}$, $i = 1, ..., n$. Considering that g.c.d. of 2nd order minors of this matrix is $\varphi_1 \varphi_2$, we conclude that $b_{2i} = \varphi_2 d_{1i}$, $i = 1, ..., n$. So

$$\begin{Vmatrix} \Phi & B_1' \end{Vmatrix} \begin{Vmatrix} 1 & 0 & & 0 & ... & 0 \\ 0 & 1 & & -d_{21} & ... & -d_{2n} \\ & & \ddots & & & \\ 0 & 0 & & 1 & & 0 \\ & & & & \ddots & \\ 0 & 0 & & 0 & & 1 \end{Vmatrix} =$$

$$= \left\| \begin{matrix} \varphi_1 & & & & & & 0 & ... & 0 \\ & \varphi_2 & & & & & 0 & ... & 0 \\ & & \varphi_3 & & & & b_{31} & ... & b_{3n} \\ & & & \ddots & & & ... & ... & ... \\ & & & & \varphi_r & & b_{r1} & ... & b_{rn} \\ & & & & & 0 & 0 & ... & 0 \\ & & & & & & \ddots & ... & ... & ... \\ & & & & & & 0 & ... & 0 \end{matrix} \right\|.$$

By a similar scheme we show that $b_{ij} = \varphi_i d_{ij}$, $i = 3, ..., r$, $j = 1, ..., n$. Thus, the matrix $B_1$ has the form $B_1 = \Phi D$. It follows that

$$PB = B_1 = \Phi D \Rightarrow B = P^{-1}\Phi D.$$

Then,

$$A(QD) = P^{-1}\Phi Q^{-1}(QD) = P^{-1}\Phi D = B.$$

Consequently, $QD$ is the desired root. $\qquad\square$

## 2.3. Transforming matrices and Zelisko group

According to Theorem 2.2, search of the Smith forms of matrices is related to finding of g.c.d. of corresponding minors. However, in many problems, including matrix factorization, we must know the transforming matrices. Further studies focus on the research of these matrices.

Let $\Phi$ be a $d$-matrix. Consider the set of matrices

$$\mathbf{G}_\Phi = \{ H \in \mathrm{GL}_n(R) \mid \exists \, K \in \mathrm{GL}_n(R) : \ H\Phi = \Phi K \}.$$

**Property 2.1.** *The set $\mathbf{G}_\Phi$ is a multiplicative group.*
**Proof.** It is obvious that $I \in \mathbf{G}_\Phi$. Let $H_1, H_2 \in \mathbf{G}_\Phi$. So

$$H_1\Phi = \Phi K_1, H_2\Phi = \Phi K_2, K_1, K_2 \in \mathrm{GL}_n(R).$$

Then

$$H_2 H_1 \Phi = H_2 \Phi K_1 = \Phi K_1 K_2.$$

That is, the set $\mathbf{G}_\Phi$ is multiplicatively closed. Moreover, the equality $H\Phi = \Phi K$ implies that $H^{-1}\Phi = \Phi K^{-1}$. Thus, $H^{-1} \in \mathbf{G}_\Phi$. $\qquad\square$

The group $\mathbf{G}_\Phi$ is called the **Zelisko group** [63].
**Property 2.2.** *If $B = P^{-1}\Phi Q^{-1}$, then $\mathbf{P}_B = \mathbf{G}_\Phi P$.*
**Proof.** Let $P_1 \in \mathbf{P}_B$ and $B = P_1^{-1}\Phi Q_1^{-1}$. Hence,

$$P^{-1}\Phi Q^{-1} = P_1^{-1}\Phi Q_1^{-1}.$$

It follows that

$$P_1 P^{-1}\Phi = \Phi Q_1^{-1} Q.$$

So $P_1 P^{-1} = H \in \mathbf{G}_\Phi$. That is $P_1 = HP$. Consequently, $\mathbf{P}_B \subseteq \mathbf{G}_\Phi P$.

On the contrary, if $P \in \mathbf{P}_B$ and $H \in \mathbf{G}_\Phi$, then

$$HPB = HPP^{-1}\Phi Q^{-1} = H\Phi Q^{-1} = \Phi K Q^{-1}.$$

Hence, $(HP)B(QK^{-1}) = \Phi$. So $HP \in \mathbf{P}_B$. Therefore, $\mathbf{P}_B \supseteq \mathbf{G}_\Phi P$. This yields $\mathbf{P}_B = \mathbf{G}_\Phi P$.  □

Consequently, the set $\mathbf{P}_B$ is the right coset of $\mathbf{G}_\Phi$ in $\mathrm{GL}_n(R)$. We characterize the structure of matrices from the group $\mathbf{G}_\Phi$.

**Theorem 2.7.** *The group* $\mathbf{G}_\Phi$, *where* $\Phi = \mathrm{diag}(\varphi_1, ..., \varphi_t, 0, ..., 0)$ *consists of all invertible matrices of the form*

$$H = \left\| \begin{matrix} H_1 & * \\ \mathbf{0} & H_2 \end{matrix} \right\|, \tag{2.2}$$

*where*

$$H_1 = \left\| \begin{matrix} h_{11} & h_{12} & ... & h_{1.t-1} & h_{1t} \\ \dfrac{\varphi_2}{\varphi_1}h_{21} & h_{22} & ... & h_{2.t-1} & h_{2t} \\ ... & ... & ... & ... & ... \\ \dfrac{\varphi_t}{\varphi_1}h_{t1} & \dfrac{\varphi_t}{\varphi_2}h_{t2} & ... & \dfrac{\varphi_t}{\varphi_{t-1}}h_{t.t-1} & h_{tt} \end{matrix} \right\|, \quad H_2 \in \mathrm{GL}_{n-t}(R).$$

**Proof.** Let $H = \|p_{ij}\|_1^n \in \mathbf{G}_\Phi$. It means that there is an invertible matrix $K = \|k_{ij}\|_1^n$ such that

$$\left\| \begin{matrix} \varphi_1 p_{11} & ... & \varphi_t p_{1t} & 0 & ... & 0 \\ ... & ... & ... & ... & ... & ... \\ \varphi_1 p_{t1} & ... & \varphi_t p_{tt} & 0 & ... & 0 \\ \varphi_1 p_{t+1.1} & ... & \varphi_t p_{t+1.t} & 0 & ... & 0 \\ ... & ... & ... & ... & ... & ... \\ \varphi_1 p_{n1} & ... & \varphi_t p_{nt} & 0 & ... & 0 \end{matrix} \right\| = \left\| \begin{matrix} \varphi_1 k_{11} & ... & \varphi_1 k_{1n} \\ ... & ... & ... \\ \varphi_t k_{t1} & ... & \varphi_t k_{tn} \\ 0 & ... & 0 \\ ... & ... & ... \\ 0 & ... & 0 \end{matrix} \right\|. \tag{2.3}$$

Since, $\varphi_1, ..., \varphi_t \neq 0$, then

$$\left\| \begin{matrix} p_{t+1.1} & ... & p_{t+1.t} \\ ... & ... & ... \\ p_{n1} & ... & p_{nt} \end{matrix} \right\| = \mathbf{0}.$$

In view of (2.3), we have $\varphi_i | \varphi_j p_{ij}$, $i,j = 1, ..., t$. If $i \leq j$ then $\varphi_i | \varphi_j$. Thus, elements $p_{ij}$, where $i \leq j$ do not have any restriction. If $i > j$, then $\varphi_j | \varphi_i$. It follows that the equality

$$\varphi_j p_{ij} = \varphi_i k_{ij}$$

is equivalent to

$$p_{ij} = \dfrac{\varphi_i}{\varphi_j} k_{ij}.$$

Thus, the matrix $H$ has the form (2.3).

Suppose that $H$ is an invertible matrix of the form (2.2). The equality

$$\left\| \begin{matrix} H_1 & * \\ \mathbf{0} & H_2 \end{matrix} \right\| \Phi = \Phi \left\| \begin{matrix} K_1 & \mathbf{0} \\ * & H_2 \end{matrix} \right\|,$$

where

$$K_1 = \left\| \begin{matrix} h_{11} & \dfrac{\varphi_2}{\varphi_1}h_{12} & ... & \dfrac{\varphi_{t-1}}{\varphi_1}h_{1.t-1} & \dfrac{\varphi_t}{\varphi_1}h_{1t} \\ h_{21} & h_{22} & ... & \dfrac{\varphi_{t-1}}{\varphi_2}h_{2.t-1} & \dfrac{\varphi_t}{\varphi_2}h_{2t} \\ ... & ... & ... & ... & ... \\ h_{t1} & h_{t2} & ... & h_{t.t-1} & h_{tt} \end{matrix} \right\|$$

is fulfilled. In view of

$$\det H = \det H_1 \det H_2,$$

we have $H_1, H_2$ are invertible matrices. The matrices $H_1, K_1$ satisfy the equality

$$H_1 \mathrm{diag}(\varphi_1, ..., \varphi_t) = \mathrm{diag}(\varphi_1, ..., \varphi_t)K_1.$$

Taking into account $\mathrm{diag}(\varphi_1, ..., \varphi_t)$ is a nonsingular matrix, $K_1$ is an invertible matrix. Consequently, the matrix $\left\| \begin{matrix} K_1 & \mathbf{0} \\ * & H_2 \end{matrix} \right\|$ is an invertible. Hence, $H \in \mathbf{G}_\Phi$. $\qquad\square$

**Corollary 2.2.** *The group of invertible upper triangular matrices is a subgroup of any Zelisko group.* $\qquad\square$

**Property 2.3.** *The element $\dfrac{\varphi_i}{\varphi_j}$ is a product of fractions of the first subdiagonal of the matrix $H_1$, being above and to the right:*

$$\begin{matrix} h_{jj} & * & * & & * \\ \dfrac{\varphi_{j+1}}{\varphi_j} & h_{j+1.j+1} & * & & * \\ \vdots & \dfrac{\varphi_{j+2}}{\varphi_{j+1}} & h_{j+2.j+2} & & * \\ \vdots & & \ddots & \ddots & \\ \dfrac{\varphi_i}{\varphi_j} & ... & ... & \dfrac{\varphi_i}{\varphi_{i-1}} & h_{ii} \end{matrix}.$$

**Proof.** Indeed,

$$\frac{\varphi_i}{\varphi_j} = \frac{\varphi_{j+1}}{\varphi_j}\frac{\varphi_{j+2}}{\varphi_{j+1}}...\frac{\varphi_i}{\varphi_{i-1}}. \qquad\square$$

**Corollary 2.3.** *If*

$$\left(\sigma, \frac{\varphi_i}{\varphi_j}\right) = 1,$$

*then the element $\sigma$ is relatively prime with all elements $\dfrac{\varphi_p}{\varphi_q}$ of the matrix $H_1$, delineated by a right triangle with vertices $(j+1, j)$, $(i, i-1)$, $(i, j)$.* $\qquad\square$

**Property 2.4.** *The element $\dfrac{\varphi_i}{\varphi_j}$, $i > j$, is a divisor of all elements of the matrix delineated by a rectangle with vertices $(i,1)$, $(i,j)$, $(t,j)$, $(t,1)$:*

$$\frac{\varphi_i}{\varphi_1}h_{i1} \qquad \frac{\varphi_i}{\varphi_2}h_{i2} \qquad ... \qquad \frac{\varphi_i}{\varphi_j}h_{ij}$$

$$\frac{\varphi_{i+1}}{\varphi_1}h_{i+1.1} \quad \frac{\varphi_{i+1}}{\varphi_2}h_{i+1.2} \quad ... \quad \frac{\varphi_{i+1}}{\varphi_j}h_{i+1.j}$$

$$\begin{array}{cccc} ... & ... & ... & ... \\ \dfrac{\varphi_t}{\varphi_1}h_{t1} & \dfrac{\varphi_t}{\varphi_2}h_{t2} & ... & \dfrac{\varphi_t}{\varphi_j}h_{tj}. \end{array}$$

*It means that*

$$\frac{\varphi_i}{\varphi_j}\left|\frac{\varphi_{i+p}}{\varphi_{j-q}}h_{i+p.j-q},\ p = 0,1,...,n-i,\ q = 0,1,...,j-1.\right.$$

**Proof.** Considering Property 2.3, the proving of this statement does not cause any difficulties. $\square$

**Property 2.5.** *Let*

$$\sigma\left|\frac{\varphi_i}{\varphi_j}\right.,\quad \left(\sigma, \frac{\varphi_i}{\varphi_{j+1}}\right) = 1,$$

*Then*

*1) the element $\sigma$ is a divisor of all elements of the matrix delineated by a rectangle with vertices $(j+1,1)$, $(j+1,j)$, $(n,j)$, $(n,1)$. That is*

$$\sigma\left|\frac{\varphi_{j+1+p}}{\varphi_{j-q}},\ p = 0,1,...,n-j-1,\ q = 0,1,...,j-1,\ j+1+p > j-q.\right.$$

*2) the element $\sigma$ is relatively prime with all fractions $\dfrac{\varphi_p}{\varphi_q}$ of the matrix delineated by a right triangle with vertices $(j+2,j+1)$, $(i,i-1)$, $(i,j+1)$. That is*

$$\left(\sigma, \frac{\varphi_{i-s}}{\varphi_{j+t}}\right) = 1,\ s = 0,1,...,i-j-2,\ t = 1,2,...,i-j-1,\ i-s > j+t,$$

$$\begin{array}{ccccccc} \dfrac{\varphi_{j+1}}{\varphi_1} & ... & \dfrac{\varphi_{j+1}}{\varphi_j} & * & * & ... & * \\[2mm] \dfrac{\varphi_{j+2}}{\varphi_1} & ... & \dfrac{\varphi_{j+2}}{\varphi_j} & \dfrac{\varphi_{j+2}}{\varphi_{j+1}} & * & ... & * \\[2mm] \dfrac{\varphi_{j+3}}{\varphi_1} & ... & \dfrac{\varphi_{j+3}}{\varphi_j} & \dfrac{\varphi_{j+3}}{\varphi_{j+1}} & \dfrac{\varphi_{j+3}}{\varphi_{j+2}} & ... & * \\[2mm] \vdots & & \vdots & \vdots & \vdots & \ddots & \\[1mm] \dfrac{\varphi_i}{\varphi_1} & ... & \dfrac{\varphi_i}{\varphi_j} & \dfrac{\varphi_i}{\varphi_{j+1}} & \dfrac{\varphi_i}{\varphi_{j+2}} & ... & \dfrac{\varphi_i}{\varphi_{i-1}} \\[2mm] \vdots & & \vdots & & & & \\[1mm] \dfrac{\varphi_{j+1}}{\varphi_1} & ... & \dfrac{\varphi_{j+1}}{\varphi_j} & * & * & ... & * \end{array}.$$

**Proof.** Since

$$\frac{\varphi_i}{\varphi_j} = \frac{\varphi_i}{\varphi_{j+1}} \frac{\varphi_{j+1}}{\varphi_j}$$

and

$$\sigma \left| \frac{\varphi_i}{\varphi_j}, \quad \left( \sigma, \frac{\varphi_i}{\varphi_{j+1}} \right) = 1, \right.$$

we claim that $\sigma \left| \dfrac{\varphi_{j+1}}{\varphi_j} \right.$. Using Property 2.4, the case 1) is complete. The Case 2) follows from Corollary 2.3. $\square$

Combining the results of the last two statements, we obtain.

**Corollary 2.4.** *Let*

$$\sigma \left| \frac{\varphi_i}{\varphi_j}, \quad \left( \sigma, \frac{\varphi_{i-1}}{\varphi_j} \right) = \left( \sigma, \frac{\varphi_i}{\varphi_{j+1}} \right) = 1. \right.$$

*Then $j = i - 1$, i.e., $\sigma \left| \dfrac{\varphi_i}{\varphi_{i-1}} \right.$.* $\square$

## 2.4. Transforming lower unitriangular matrices

As noted in the previous subsection, the set $\mathbf{P}_B$ of left transforming matrices of the matrix $B = P^{-1} \Phi Q^{-1}$ is the right coset of $\mathbf{G}_\Phi$ in $\mathrm{GL}_n(R)$. Let us state the conditions under which there exists a lower unitriangular matrix in this set.

By $A_{(i)}$ denote a submatrix of the matrix $A = \|a_{ij}\|_1^n \in M_n(R)$ of the form

$$A_{(i)} = \left\|
\begin{array}{cccc}
a_{ii} & a_{i.i+1} & ... & a_{in} \\
a_{i+1.i} & a_{i+1.i+1} & ... & a_{i+1.n} \\
... & ... & ... & ... \\
a_{ni} & a_{n.i+1} & ... & a_{nn}
\end{array}
\right\|,$$

$i = 1, ..., n$.

**Lemma 2.1.** *Suppose that $P = \|p_{ij}\|_1^n$ is an invertible matrix and $HP = Q$, where $H \in \mathbf{G}_\Phi$, $\det \Phi \neq 0$. Then*

$$\left( \frac{\varphi_{i+1}}{\varphi_i}, |P_{(i+1)}| \right) = \left( \frac{\varphi_{i+1}}{\varphi_i}, |Q_{(i+1)}| \right), \ i = 1, ..., n-1.$$

**Proof**. By Theorem 2.7, a submatrix of the matrix $H$, consisting of $m$ its last rows has the form

$$K_s = \left\|
\begin{array}{ccc|cccc}
\frac{\varphi_s}{\varphi_1} h_{s1} & ... & \frac{\varphi_s}{\varphi_{s-1}} h_{s.s-1} & h_{ss} & ... & h_{s.n-1} & h_{sn} \\
\frac{\varphi_{s+1}}{\varphi_1} h_{s+1.1} & ... & \frac{\varphi_{s+1}}{\varphi_{s-1}} h_{s+1.s-1} & \frac{\varphi_{s+1}}{\varphi_s} h_{s+1.s} & ... & h_{s+1.n-1} & h_{s+1 n} \\
... & ... & ... & ... & ... & ... & ... \\
\frac{\varphi_n}{\varphi_1} h_{n1} & ... & \frac{\varphi_n}{\varphi_{s-1}} h_{n.s-1} & \frac{\varphi_n}{\varphi_s} h_{ns} & ... & \frac{\varphi_n}{\varphi_{n-1}} h_{n.n-1} & h_{n n}
\end{array}
\right\| =$$

$$= \left\| H'_s \quad H_s \right\|,$$

**65**

where $s = n - m + 1$. Since

$$\frac{\varphi_s}{\varphi_{s-1}} \left| \frac{\varphi_{s+k}}{\varphi_{s-1-l}}, \right.$$

so all elements of the matrix $H'_s$ are divisible by $\frac{\varphi_s}{\varphi_{s-1}}$. Hence, all $m$ order minors of the matrix $K_s$, except the minor $|H_s|$, are divisible by $\frac{\varphi_s}{\varphi_{s-1}}$. This implies that

$$\left( |H_s|, \frac{\varphi_s}{\varphi_{s-1}} \right) = 1. \tag{2.4}$$

Since

$$Q_{(s)} = K_s \begin{Vmatrix} p_{1s} & p_{1.s+1} & ... & p_{1n} \\ p_{2s} & p_{2.s+1} & ... & p_{2n} \\ ... & ... & ... & ... \\ p_{ns} & p_{n.s+1} & ... & p_{nn} \end{Vmatrix},$$

by using the Binet—Cauchy formula, we can represent the determinant $|Q_{(s)}|$ in the form

$$|Q_{(s)}| = \sum_j |H_j||P_j| + |H_s||P_{(s)}|,$$

where $\sum_j |H_j||P_j|$ is the sum of products of all possible $m$ order minors of the matrix $K_s$ (except the minor $|H_s|$) by the corresponding minors of the matrix $P$. Since all minors $|H_j|$ are divisible by $\frac{\varphi_s}{\varphi_{s-1}}$, we find

$$|Q_{(s)}| = \frac{\varphi_s}{\varphi_{s-1}} d + |H_s||P_{(s)}|,$$

where $d \in R$. Thus, in view of the equality (2.4), we obtain

$$\left( \frac{\varphi_s}{\varphi_{s-1}}, |Q_{(s)}|, \right) = \left( \frac{\varphi_s}{\varphi_{s-1}}, |H_s||P_{(s)}| \right) = \left( \frac{\varphi_s}{\varphi_{s-1}}, |P_{(s)}| \right). \qquad \square$$

**Theorem 2.8.** *Suppose that $S \in \mathrm{GL}_n(R)$ and $\Phi$ be a nonsingular d-matrix. In order that the group $\mathbf{G}_\Phi$ has a matrix $H$ such that $HS$ is a lower unitriangular matrix, it is necessary and sufficient that*

$$\left( \frac{\varphi_{i+1}}{\varphi_i}, |S_{(i+1)}| \right) = 1, \quad i = 1, ..., n-1. \tag{2.5}$$

**Proof. Necessity.** Let $HS = T = \|t_{ij}\|_1^n$ be a lower unitriangular matrix. Then, by virtue of Lemma 2.1, we conclude that

$$\left( \frac{\varphi_{i+1}}{\varphi_i}, |S_{(i+1)}| \right) = \left( \frac{\varphi_{i+1}}{\varphi_i}, \begin{vmatrix} 1 & 0 & & 0 \\ t_{i+2.i+1} & 1 & & 0 \\ ... & ... & & \\ t_{n.i+1} & t_{n.i+2} & ... & t_{n.n-1} & 1 \end{vmatrix} \right) = 1, \ i = 1, ..., n-1.$$

**Sufficiency.** Let $S = \|s_{ij}\|_1^2$ be an invertible matrix over $R$. Since $(s_{12}, s_{22}) = 1$ and

$$\left(\frac{\varphi_2}{\varphi_1}, s_{22}\right) = 1,$$

then

$$\left(\frac{\varphi_2}{\varphi_1} s_{12}, s_{22}\right) = 1.$$

There exist $u_1$, $u_2$ such that

$$s_{12} \frac{\varphi_2}{\varphi_1} u_1 + s_{22} u_2 = 1.$$

This implies that the matrix

$$K = \left\| \begin{matrix} s_{22} & -s_{12} \\ \dfrac{\varphi_2}{\varphi_1} u_1 & u_2 \end{matrix} \right\|$$

is an element of the group $\mathbf{G}_\Phi$. Hence,

$$KS = \left\| \begin{matrix} d & 0 \\ c & 1 \end{matrix} \right\|.$$

The invertibility of the matrix $KS$ implies that $d \in U(R)$. Therefore, $H = \mathrm{diag}(d^{-1}, 1)K$ is just the required matrix. Thus, the theorem is true for matrices of order 2.

Further, we assume that the theorem is true for all matrices of order lower than $n$ and consider an invertible matrix $S = \|s_{ij}\|_1^n$. By using equalities (2.5), we get

$$\left(\frac{\varphi_{i+1}}{\varphi_i}, \ (s_{n.i+1}, s_{n.i+2}, ..., s_{nn})\right) = 1,$$

$i = 1, ..., n - 1$. In particular, we find

$$\left(\frac{\varphi_2}{\varphi_1}, \ (s_{n2}, s_{n3}, ..., s_{nn})\right) = 1.$$

Since

$$(s_{n1}, (s_{n2}, s_{n3}, ..., s_{nn})) = 1,$$

we obtain

$$\left(\frac{\varphi_2}{\varphi_1} s_{n1}, (s_{n2}, s_{n3}, ..., s_{nn})\right) = \left(\left(\frac{\varphi_2}{\varphi_1} s_{n1}, s_{n2}\right), (s_{n3}, ..., s_{nn})\right) = 1.$$

In view of the fact that

$$\left(\frac{\varphi_3}{\varphi_2}, (s_{n3}, ..., s_{nn})\right) = 1,$$

we can write

$$\left(\left(\frac{\varphi_3}{\varphi_1}s_{n1}, \frac{\varphi_3}{\varphi_2}s_{n2}\right), (s_{n3}, ..., s_{nn})\right) = 1.$$

We continue the outlined process step by step and, finally, get

$$\left(\frac{\varphi_n}{\varphi_1}s_{1n}, ..., \frac{\varphi_n}{\varphi_{n-2}}s_{n-2.n}, \ \frac{\varphi_n}{\varphi_{n-1}}s_{n-1.n}, s_{nn}\right) = 1.$$

There exist $u_{n1}, u_{n2}, ..., u_{nn} \in R$ such that

$$\frac{\varphi_n}{\varphi_1}u_{n1}s_{1n} + ... + \frac{\varphi_n}{\varphi_{n-1}}u_{n.n-1}s_{n-1.n} + u_{nn}s_{nn} = 1.$$

This implies that the row

$$\left\|\frac{\varphi_n}{\varphi_1}u_{n1} \quad ... \quad \frac{\varphi_n}{\varphi_{n-1}}u_{n.n-1} \quad u_{nn}\right\|$$

is unimodular. By using Theorem 1.1, we complement this row to an invertible matrix of the form

$$H_1 = \left\|\begin{array}{cccccc} u_{11} & u_{12} & ... & u_{1.n-2} & u_{1.n-1} & u_{1n} \\ 0 & u_{22} & ... & u_{2.n-2} & u_{2.n-1} & u_{2n} \\ ... & ... & ... & ... & ... & ... \\ 0 & 0 & & u_{n-2.n-2} & u_{n-2.n-1} & u_{n-2.n} \\ 0 & 0 & ... & 0 & u_{n-1.n-1} & u_{n-1.n} \\ \frac{\varphi_n}{\varphi_1}u_{n1} & \frac{\varphi_n}{\varphi_2}u_{n2} & ... & \frac{\varphi_n}{\varphi_{n-2}}u_{n.n-2} & \frac{\varphi_n}{\varphi_{n-1}}u_{n.n-1} & u_{nn} \end{array}\right\|,$$

which is an element of the group $\mathbf{G}_\Phi$. Thus,

$$H_1 S = \left\|\begin{array}{cc} S_1 & \mathbf{q} \\ \mathbf{p} & 1 \end{array}\right\|,$$

where $\mathbf{p} = \left\|\begin{array}{cccc} s'_{n1} & s'_{n2} & ... & s'_{n.n-1}\end{array}\right\|$, $\mathbf{q} = \left\|\begin{array}{cccc} s'_{1n} & s'_{2n} & ... & s'_{n-1.n}\end{array}\right\|^T$. Hence,

$$\underbrace{\left\|\begin{array}{cc} I & -\mathbf{q} \\ \mathbf{0} & 1 \end{array}\right\|}_{H_2} H_1 S = \left\|\begin{array}{cc} S_2 & \mathbf{0} \\ \mathbf{p} & 1 \end{array}\right\|,$$

where $S_2 = \left\|s''_{ij}\right\|_1^n$. Since $H_2 H_1 \in \mathbf{G}_\Phi$, by Lemma 2.1, we get

$$\left(\frac{\varphi_{i+1}}{\varphi_i}, \left|\begin{array}{cccc} s''_{i+1.i+1} & ... & s''_{i+1.n-1} & 0 \\ ... & ... & ... & ... \\ s''_{n-1.i+1} & ... & s''_{n-1.n-1} & 0 \\ s'_{n.i+1} & ... & s'_{n.n-1} & 1 \end{array}\right|\right) = 1,$$

$i = 1, ..., n - 1$. This implies that

$$\left( \frac{\varphi_{i+1}}{\varphi_i}, \begin{vmatrix} s''_{i+1.i+1} & \cdots & s''_{i+1.n-1} \\ \cdots & \cdots & \cdots \\ s''_{n-1.i+1} & \cdots & s''_{n-1.n-1} \end{vmatrix} \right) = 1, \;\; i = 1, ..., n - 2.$$

Consider the matrix $\Phi_{n-1} = \mathrm{diag}\,(\varphi_1, ..., \varphi_{n-1})$. By the induction hypothesis, the group $\mathbf{G}_{\Phi_{n-1}}$ contains a matrix $K$ such that $KS_2$ is a lower unitriangular matrix. Thus,

$$\begin{Vmatrix} K & \mathbf{0} \\ \mathbf{0} & 1 \end{Vmatrix} H_2 H_1 S$$

is also a lower unitriangular matrix. To complete the proof, note that $\begin{Vmatrix} K & \mathbf{0} \\ \mathbf{0} & 1 \end{Vmatrix}$ is an element of the group $\mathbf{G}_\Phi$. $\qquad\square$

**Theorem 2.9.** *Suppose that the nonsingular matrix $B$ has the Smith form $\Phi$ and $S \in \mathbf{P}_B$. The set $\mathbf{P}_B$ contains a lower unitriangular matrix if and only if the matrix $S$ satisfies the equality* (2.5).

**Proof** follows from Property 2.2 and Theorem 2.8. $\qquad\square$

## 2.5. Structure of transforming matrices of one class of matrices

By Theorem 2.3, the second order matrices are playing the main role in establishing the fact that $R$ is an elementary divisor ring. Let us show the explicit form of their transforming matrices.

**Theorem 2.10.** *Assume that $A = \begin{Vmatrix} b & c \\ a & 0 \end{Vmatrix}$ is the matrix over an elementary divisor ring, where $(a, b, c) = 1$. Then*

$$P = \begin{Vmatrix} \beta & m \\ -\alpha a & \alpha b + cn \end{Vmatrix} \in \mathbf{P}_A, \quad Q = \begin{Vmatrix} \alpha & c\beta \\ n & -b\beta - am \end{Vmatrix} \in \mathbf{Q}_A,$$

*where*

$$a(m\alpha) + b(\alpha\beta) + c(\beta n) = 1. \tag{2.6}$$

**Proof.** Since $(a, b, c) = 1$, based on Corollary 2.1, there are $m, n, \alpha, \beta \in R$ such that equality (2.6) is hold. Then

$$\begin{Vmatrix} b & c \\ a & 0 \end{Vmatrix} = \underbrace{\begin{Vmatrix} cn + \alpha b & -m \\ \alpha a & \beta \end{Vmatrix}}_{M} \begin{Vmatrix} 1 & 0 \\ 0 & ac \end{Vmatrix} \underbrace{\begin{Vmatrix} \beta b + am & \beta c \\ n & -\alpha \end{Vmatrix}}_{N}.$$

To complete the proof it suffices to verify that $M^{-1} = P$ and $N^{-1} = Q$. $\quad\square$

**69**

**Corollary 2.5.** *If $c = 0$ and $am + cn = 1$, then*

$$P = \left\| \begin{matrix} 1 & m \\ -a & cn \end{matrix} \right\| \in \mathbf{P}_A, \quad Q = \left\| \begin{matrix} 1 & c \\ n & -am \end{matrix} \right\| \in \mathbf{Q}_A. \qquad \square$$

Naturally, the question arises whether the matrix $A$ has a transforming matrices with a different structure. An answer to this question is negative.

**Theorem 2.11.** *All matrices from the set $\mathbf{P}_A$ have the form*

$$\left\| \begin{matrix} \beta_1 & m_1 \\ -\alpha_1 a & \alpha_1 b + c n_1 \end{matrix} \right\|,$$

*where*

$$a(m_1 \alpha_1) + b(\alpha_1 \beta_1) + c(\beta_1 n_1) = e \in U(R).$$

**Proof.** According to Property 2.2, we have $\mathbf{P}_A = \mathbf{G}_\Phi P$, where $P$ is an arbitrary matrix from $\mathbf{P}_A$, and $\mathbf{G}_\Phi$ consists of all invertible matrices of the form

$$\left\| \begin{matrix} k_{11} & k_{12} \\ ack_{21} & k_{22} \end{matrix} \right\|.$$

As $P$ we choose the matrix obtained in Theorem 2.10, i.e.,

$$P = \left\| \begin{matrix} \beta & m \\ -\alpha a & \alpha b + cn \end{matrix} \right\|,$$

where

$$a(m\alpha) + b(\alpha\beta) + c(\beta n) = 1.$$

Let $P_1$ belongs to $\mathbf{P}_A$. Then the group $\mathbf{G}_\Phi$ contains the matrix

$$H = \left\| \begin{matrix} h_{11} & h_{12} \\ ach_{21} & h_{22} \end{matrix} \right\|$$

such that $P_1 = HP$. Thus,

$$P_1 = \left\| \begin{matrix} h_{11} & h_{12} \\ ach_{21} & h_{22} \end{matrix} \right\| \left\| \begin{matrix} \beta & m \\ -\alpha a & \alpha b + cn \end{matrix} \right\| = \left\| \begin{matrix} \beta_1 & m_1 \\ -\alpha_1 a & \alpha_1 b + c n_1 \end{matrix} \right\|,$$

where

$$\alpha_1 = c\beta h_{21} + \alpha h_{22}, \ \beta_1 = \beta h_{11} - \alpha\beta a h_{22},$$
$$m_1 = m h_{11} + h_{22}(\alpha b + cn), \ n_1 = n h_{22} + h_{21}(\beta b + am).$$

Since $\det H = e \in U(R)$, we have $\det P_1 = e$. Consequently,

$$a(m_1 \alpha_1) + b(\alpha_1 \beta_1) + c(\beta_1 n_1) = e. \qquad \square$$

A similar statement can be made regarding the structure of right transforming matrices of the matrix $A$.

## 2.6. Some properties
## of elements of elementary divisor rings

It was shown in Theorem 1.1 that, each primitive row $\|a \quad b \quad c\|$ over a Bezout ring $R$ can be complemented to an invertible matrix of the form

$$\begin{Vmatrix} a & b & c \\ 0 & * & * \\ * & * & * \end{Vmatrix}.$$

Over an elementary divisor ring this statement can be improved.

**Theorem 2.12.** *A commutative Bezout domain $R$ is a commutative elementary divisor domain if and only if each primitive row $\|a \quad b \quad c\|$ complemented to an invertible matrix of the form*

$$\begin{Vmatrix} a & b & c \\ 0 & * & * \\ * & * & 0 \end{Vmatrix}.$$

**Proof.** Since $(a, b, c) = 1$, there exist $m, n, \alpha, \beta$ such that

$$a(m\alpha) + b(\alpha\beta) + c(\beta n) = 1.$$

Then the required matrix has the form

$$\begin{Vmatrix} a & b & c \\ 0 & -n & \alpha \\ \beta & -m & 0 \end{Vmatrix}.$$

The converse argument is obvious. □

Note that, according to Property 1.18, every primitive row $\|a \quad b \quad c\|$, $c \neq \neq 0$, can be complemented to an invertible matrix of the form

$$\begin{Vmatrix} a & b & c \\ * & 1 & 0 \\ * & 0 & * \end{Vmatrix}.$$

The following Theorem formulates the property of the elements of the Bezout ring $R$ under which $R$ becomes the elementary divisor ring. To prove this result, we need the following statement.

**Lemma 2.2.** *Let $A = \|a_{ij}\|_1^3$, $\det A = 1$ and*

$$B = \begin{Vmatrix} b_{11} & b_{12} & b_{13} \\ a_{21} & a_{22} & a_{23} \\ a_{31} & a_{32} & a_{33} \end{Vmatrix}.$$

*In order that $\det B = 1$, necessary and sufficient there exist elements $s_2, s_3$ such that*

$$\|b_{11} \quad b_{12} \quad b_{13}\| = \|1 \quad s_2 \quad s_3\| A.$$

**71**

**Proof. Necessity**. Since

$$\det B = b_{11} \det \begin{Vmatrix} a_{22} & a_{23} \\ a_{32} & a_{33} \end{Vmatrix} - b_{12} \det \begin{Vmatrix} a_{21} & a_{23} \\ a_{31} & a_{33} \end{Vmatrix} + b_{13} \det \begin{Vmatrix} a_{21} & a_{22} \\ a_{31} & a_{32} \end{Vmatrix} = 1,$$

we have

$$BA^{-1} = \begin{Vmatrix} 1 & s_2 & s_3 \\ 0 & 1 & 0 \\ 0 & 0 & 1 \end{Vmatrix}.$$

Thus

$$B = \begin{Vmatrix} 1 & s_2 & s_3 \\ 0 & 1 & 0 \\ 0 & 0 & 1 \end{Vmatrix} A.$$

Consequently,

$$\begin{Vmatrix} b_{11} & b_{12} & b_{13} \end{Vmatrix} = \begin{Vmatrix} 1 & s_2 & s_3 \end{Vmatrix} A.$$

**Sufficiency**. Since

$$\begin{Vmatrix} b_{11} & b_{12} & b_{13} \end{Vmatrix} A^{-1} = \begin{Vmatrix} 1 & s_2 & s_3 \end{Vmatrix} AA^{-1} =$$
$$= \begin{Vmatrix} 1 & s_2 & s_3 \end{Vmatrix} E = \begin{Vmatrix} 1 & s_2 & s_3 \end{Vmatrix},$$

we get

$$BA^{-1} = \begin{Vmatrix} 1 & s_2 & s_3 \\ 0 & 1 & 0 \\ 0 & 0 & 1 \end{Vmatrix}.$$

This implies that $\det B = 1$. □

**Theorem 2.13.** *In order that a commutative Bezout domain $R$ to be a commutative elementary divisor domain it is necessary and sufficient for every four elements $a_1, a_2, b_1, b_2$ such that*

$$(a_1, a_2) = (b_1, b_2) = 1$$

*there exists an element $r$ such that*

$$b_1 + rb_2 = \alpha\beta,$$

*where*

$$(\alpha, \beta) = (a_1, \alpha) = (a_2, \beta) = 1.$$

**Proof. Necessity.** The ring $R$ contains $v_1, v_2, u_1, u_2$ such that

$$a_1 v_1 - a_2 v_2 = 1,$$
$$b_1 u_1 - b_2 u_2 = 1.$$

Thus, the matrix

$$\begin{Vmatrix} u_2 v_2 & b_1 & u_2 v_1 \\ u_1 v_2 & b_2 & u_1 v_1 \\ a_1 & 0 & a_2 \end{Vmatrix} = A$$

is unimodular, i.e.,

$$\det A = a_2 b_2 (u_2 v_2) - \begin{vmatrix} u_1 v_2 & u_1 v_1 \\ a_1 & a_2 \end{vmatrix} b_1 - a_1 b_2 (u_2 v_1) = 1.$$

Denote

$$a = u_2 v_2, \ b = \begin{vmatrix} u_1 v_2 & u_1 v_1 \\ a_1 & a_2 \end{vmatrix}, \ c = v_1 u_2.$$

It is obvious that

$$(a, b, c) = 1.$$

There exist elements $\alpha, \beta, m, n$ such that

$$a(m\alpha) - b(\alpha\beta) - c(\beta n) = 1.$$

It follows that

$$(\alpha, \beta) = 1.$$

Thus, the matrix

$$\begin{Vmatrix} \alpha m & \alpha\beta & \beta n \\ u_1 v_2 & b_2 & u_1 v_1 \\ a_1 & 0 & a_2 \end{Vmatrix} = B$$

is unimodular. By virtue of Lemma 2.2, there exist elements $s_2, s_3$ such that

$$\begin{Vmatrix} \alpha m & \alpha\beta & \beta n \end{Vmatrix} = \begin{Vmatrix} 1 & s_2 & s_3 \end{Vmatrix} A,$$

i.e.,

$$b_1 + s_2 b_2 = \alpha\beta.$$

Since

$$1 = \det B = a_1(-\beta n b_2) + a_1(\alpha\beta u_1 v_1) + \alpha \left( \begin{vmatrix} m & \beta \\ u_1 v_2 & b_2 \end{vmatrix} a_2 \right),$$

we have

$$(\alpha, a_1) = 1.$$

By analogy, we show that

$$(\beta, a_2) = 1.$$

Setting $r = s_2$ we complete the proof of the necessity part.

**Sufficiency.** Let $a, b, c$ be relatively prime elements of the ring $R$ and, moreover, $(a, c) \neq 0$. Consider the matrix

$$\begin{Vmatrix} -\dfrac{bu}{(a,c)} & (a,c) & \dfrac{bv}{a} \\ \dfrac{}{(a,c)} & 0 & \dfrac{a}{(a,c)} \end{Vmatrix},$$

where

$$\frac{a}{(a,c)} u + \frac{c}{(a,c)} v = -1.$$

For $a, b, c$ there exist $m_1, m_2, m_3$ such that

$$a m_1 + b m_2 + c m_3 = 1.$$

Thus,

$$\det \underbrace{\left\| \begin{array}{ccc} m_1 & m_2 & m_3 \\ bu & (a,c) & bv \\ -\dfrac{c}{(a,c)} & 0 & \dfrac{a}{(a,c)} \end{array} \right\|}_{A} = \det A = 1.$$

Set

$$A = \left\| \begin{array}{ccc} a_{11} & a_{12} & a_{13} \\ a_{21} & a_{22} & a_{23} \\ a_{31} & 0 & a_{33} \end{array} \right\|.$$

It is obvious that

$$(a_{12}, a_{22}) = (a_{31}, a_{33}) = 1.$$

By assumption, there exists $r$ such that

$$a_{12} + r a_{22} = \alpha\beta,$$

where

$$(\alpha, \beta) = (a_{31}, \alpha) = (a_{33}, \beta) = 1. \tag{2.7}$$

Then

$$\left\| \begin{array}{ccc} 1 & r & 0 \\ 0 & 1 & 0 \\ 0 & 0 & 1 \end{array} \right\| \left\| \begin{array}{ccc} a_{11} & a_{12} & a_{13} \\ a_{21} & a_{22} & a_{23} \\ a_{31} & 0 & a_{33} \end{array} \right\| =$$

$$= \left\| \begin{array}{ccc} a_{11} + r a_{21} & \alpha\beta & a_{13} + r a_{23} \\ a_{21} & a_{22} & a_{23} \\ a_{31} & 0 & a_{33} \end{array} \right\| = \left\| \begin{array}{ccc} a'_{11} & \alpha\beta & a'_{13} \\ a_{21} & a_{22} & a_{23} \\ a_{31} & 0 & a_{33} \end{array} \right\|.$$

Consider the system of congruences

$$\begin{cases} a_{31} x \equiv a'_{11} (\mathrm{mod}\,\alpha), \\ a_{33} x \equiv a'_{13} (\mathrm{mod}\,\beta). \end{cases}$$

According to (2.7), this system has a solution $x = s$. Therefore,

$$a_{31} s - a'_{11} = -\alpha m,$$

$$a_{33} s - a'_{13} = -\beta n,$$

whence

$$a'_{11} + a_{31}(-s) = \alpha m,$$

$$a'_{13} + a_{33}(-s) = \beta n.$$

Thus,

$$\left\| \begin{array}{ccc} 1 & 0 & -s \\ 0 & 1 & 0 \\ 0 & 0 & 1 \end{array} \right\| \left\| \begin{array}{ccc} a'_{11} & \alpha\beta & a'_{13} \\ a_{21} & a_{22} & a_{23} \\ a_{31} & 0 & a_{33} \end{array} \right\| = \left\| \begin{array}{ccc} m\alpha & \alpha\beta & \beta n \\ a_{21} & a_{22} & a_{23} \\ a_{31} & 0 & a_{33} \end{array} \right\| = B.$$

Taking into account

$$1 = \det B = a(m\alpha) + b(\alpha\beta) + c(\beta n)$$

and using Corollary 2.1, we complete the analysis of this case.

If $a = c = 0$, then $b \in U(R)$. In this case, it is sufficient to set

$$\alpha = b^{-1}, \quad \beta = m = n = 1.$$

The Theorem is proved. □

## 2.7. Zelisko group and stable range of rings

The representation of matrices in the form of products of two and more factors with given properties is one of the most efficient methods for the solution and simplification of matrix problems. It is known that each matrix over a left Hermitian ring is the product of invertible and triangular matrices and each matrix over a ring of elementary divisors is the product of invertible, diagonal, and invertible matrices. Note that the possibility of representation of a matrix in the form of the product of given factors is closely connected to a stable range of the ring over which these matrices are considered.

In [74-76], it was shown that every invertible matrix over a commutative ring is the product of upper triangular, lower and upper unitriangular matrices (triangular matrices with 1 on the diagonal) if and only if $R$ is an Hermitian ring of stable range 1. We show similar connections in the case of commutative rings of stable range 1.5.

Recall that a ring $R$ has stable range 1.5 if for any $(a, b, c) = 1$, $a, b, c \in R$, $c \neq 0$ there exists $r \in R$ such that

$$(a + br, c) = 1.$$

By $U_n^{up}(R)$ and $U_n^{lw}(R)$ we denote the groups of upper and lower unitriangular $n \times n$ matrices over $R$, respectively.

**Theorem 2.14.** *Suppose that $R$ is a commutative Bezout domain. The following conditions are equivalent:*

*1) $R$ has the stable range 1.5;*

*2) $\mathrm{GL}_2(R) = \mathbf{G}_\Phi U_2^{lw}(R) U_2^{up}(R)$ for all nonsingular d-matrices $\Phi = \mathrm{diag}(\varphi_1, ..., \varphi_1)$;*

*3) $\mathrm{GL}_n(R) = \mathbf{G}_\Phi U_n^{lw}(R) U_n^{up}(R)$ for all nonsingular d-matrices $\Phi$ and any $n \geq 2$.*

**Proof.** 1) ⇒ 2). Let $A = \|a_{ij}\|_1^2 \in \mathrm{GL}_2(R)$. Then

$$(a_{21}, a_{22}) = 1 \implies \left(a_{21}, a_{22}, \frac{\varphi_2}{\varphi_1}\right) = 1.$$

**75**

In the ring $R$, there exists $r$ such that

$$\left(a_{22} + a_{21}r, \frac{\varphi_2}{\varphi_1}\right) = 1.$$

Consider a matrix $AU = \|a_{ij}'\|_1^2$, where

$$U = \left\|\begin{matrix} 1 & r \\ 0 & 1 \end{matrix}\right\|.$$

Since

$$\left(a_{22}', \frac{\varphi_2}{\varphi_1}\right) = 1,$$

by Theorem 2.8, the group $\mathbf{G}_\Phi$ contains a matrix $H$ such that $HAU = V$ is a lower unitriangular matrix. Hence, $A = H^{-1}VU^{-1}$, therefore,

$$\mathrm{GL}_2(R) = \mathbf{G}_\Phi U_2^{lw}(R)U_2^{up}(R).$$

2) $\Rightarrow$ 1). Let $(a, b, c) = 1, \ abc \neq 0$. We have

$$a = (a, b)a_1, \ \ b = (a, b)b_1, \ \ (a_1, b_1) = 1.$$

There exist $u, v$ such that

$$a_1 u + b_1 v = 1.$$

Hence, the matrix

$$\left\|\begin{matrix} u & -v \\ b_1 & a_1 \end{matrix}\right\| = A$$

is invertible. Consider a $d$-matrix

$$\left\|\begin{matrix} 1 & 0 \\ 0 & c \end{matrix}\right\| = \Phi.$$

According to the condition of the theorem, the matrix $A$ is the product of three matrices: $A = HUV$, where $H \in \mathbf{G}_\Phi$, $U \in U_2^{lw}(R)$, $V \in U_2^{up}(R)$. Note that

$$H^{-1} = \left\|\begin{matrix} h_{11} & h_{12} \\ ch_{21} & h_{22} \end{matrix}\right\|, \quad U = \left\|\begin{matrix} 1 & 0 \\ u & 1 \end{matrix}\right\|, \quad V^{-1} = \left\|\begin{matrix} 1 & r \\ 0 & 1 \end{matrix}\right\|.$$

This yields

$$U = \left\|\begin{matrix} 1 & 0 \\ u & 1 \end{matrix}\right\| = H^{-1}AV^{-1} = H^{-1}\left(\left\|\begin{matrix} u & -v \\ b_1 & a_1 \end{matrix}\right\|\left\|\begin{matrix} 1 & r \\ 0 & 1 \end{matrix}\right\|\right) =$$

$$= \left\|\begin{matrix} h_{11} & h_{12} \\ ch_{21} & h_{22} \end{matrix}\right\|\left\|\begin{matrix} u & ur - v \\ b_1 & b_1r + a_1 \end{matrix}\right\| = \left\|\begin{matrix} * & * \\ * & ch_{21}(ur - v) + h_{22}(b_1r + a_1) \end{matrix}\right\|.$$

Hence,

$$ch_{21}(ur - v) + h_{22}(b_1r + a_1) = 1.$$

This implies that
$$(b_1r + a_1, c) = 1.$$
In view of the fact that
$$((a, b), c) = 1,$$
we obtain
$$((a, b)(b_1r + a_1), c) = (a + br, c) = 1.$$

This means that the ring $R$ has the stable range 1.5. The case where $a = 0$ or $b = 0$ is obvious. The implication 3) $\Rightarrow$ 2) is obvious.

2) $\Rightarrow$ 3). We prove this implication by induction. As shown above, this statement is true for the matrices of order 2 (see proof of condition 2)). Assume that it is true for the matrices which order is lower than $n$. Since 2) $\Leftrightarrow$ 1), $R$ is a commutative Bezout domain of stable range 1.5. Let $A = \|a_{ij}\|_1^n \in \mathrm{GL}_n(R)$. Then

$$(a_{n1}, ..., a_{nn}) = 1 \Rightarrow \left(a_{n1}, ..., a_{nn}, \frac{\varphi_n}{\varphi_1}\right) = 1.$$

According to Lemma 1.19, there exist $r_1, ..., r_{n-1}$ such that

$$\left(a_{nn} + a_{n.n-1}r_{n-1} + ... + a_{n1}r_1, \frac{\varphi_n}{\varphi_1}\right) = 1.$$

Consider a matrix $AU_n = \|a'_{ij}\|_1^n$, where

$$U_n = \left\| \begin{matrix} 1 & & 0 & r_1 \\ & \ddots & 0 & \vdots \\ 0 & & 1 & r_{n-1} \\ 0 & ... & 0 & 1 \end{matrix} \right\|.$$

Since
$$((a'_{n1}, ..., a'_{n.n-1}), a'_{nn}) = 1$$
and
$$\left(\frac{\varphi_n}{\varphi_1}, a'_{nn}\right) = 1,$$
we conclude that
$$\left(\frac{\varphi_n}{\varphi_1}(a'_{n1}, ..., a'_{n.n-1}), a'_{nn}\right) = 1.$$
This implies that
$$\left(\frac{\varphi_n}{\varphi_1}a'_{n1}, \frac{\varphi_n}{\varphi_2}a'_{n2}, ..., \frac{\varphi_n}{\varphi_{n-1}}a'_{n.n-1}, a'_{nn}\right) = 1.$$

From the proof of sufficiency of Theorem 2.8, it follows that the group $\mathbf{G}_\Phi$ has a matrix $H_n$ such that
$$H_n A U_n = \left\| \begin{matrix} A_1 & \mathbf{0} \\ A_2 & 1 \end{matrix} \right\|.$$

Consider the $d$-matrix $\Phi_{n-1} = \operatorname{diag}(\varphi_1, ..., \varphi_{n-1})$. Since $A_1 \in \operatorname{GL}_{n-1}(R)$, by the assumption, we get

$$A_1 = H_{n-1}V_{n-1}U_{n-1},$$

where $H_{n-1} \in \mathbf{G}_{\Phi_{n-1}}$, $V_{n-1} \in U_{n-1}^{lw}(R)$, $U_{n-1} \in U_{n-1}^{up}(R)$. Hence,

$$H_n A U_n = \left\| \begin{matrix} A_1 & \mathbf{0} \\ A_2 & 1 \end{matrix} \right\| = \left\| \begin{matrix} H_{n-1}V_{n-1}U_{n-1} & \mathbf{0} \\ A_2 & 1 \end{matrix} \right\| =$$

$$= \underbrace{\left\| \begin{matrix} H_{n-1} & \mathbf{0} \\ \mathbf{0} & 1 \end{matrix} \right\|}_{M} \underbrace{\left\| \begin{matrix} V_{n-1} & \mathbf{0} \\ A_2 U_{n-1}^{-1} & 1 \end{matrix} \right\|}_{S} \underbrace{\left\| \begin{matrix} U_{n-1} & \mathbf{0} \\ \mathbf{0} & 1 \end{matrix} \right\|}_{N} = MSN.$$

Thus, the matrix $A$ can be rewritten in the form

$$A = (H_n^{-1}M)S(NU_n^{-1}).$$

Note that $H_n^{-1}M \in \mathbf{G}_\Phi$, $S \in U_n^{lw}(R)$, $NU_n^{-1} \in U_n^{up}(R)$, we have

$$\operatorname{GL}_n(R) = \mathbf{G}_\Phi U_n^{lw}(R)U_n^{up}(R).$$

Theorem is proved. □

**Corollary 2.6.** *Suppose that $\Phi$ be nonsingular, $E$ be arbitrary $d$-matrices, and $A \sim E$, $B \sim \Phi$. The set $\mathbf{P}_B\mathbf{P}_A^{-1}$ contains a lower unitriangular matrix.*

**Proof.** Let $P_A \in \mathbf{P}_A, P_B \in \mathbf{P}_B$. By Property 2.2, $\mathbf{P}_A = \mathbf{G}_E P_A$, $\mathbf{P}_B = \mathbf{G}_\Phi P_B$. Consequently,

$$\mathbf{P}_B\mathbf{P}_A^{-1} = \mathbf{G}_\Phi P_B P_A^{-1} \mathbf{G}_E. \tag{2.8}$$

By virtue of Theorem 2.14,

$$P_B P_A^{-1} = LUV,$$

where $L \in \mathbf{G}_\Phi$, $U \in U_n^{lw}(R)$, $V \in U_n^{up}(R)$. It follows that

$$L^{-1}(P_B P_A^{-1})V^{-1} = U.$$

By Corollary 2.2, $U_n^{up}(R) \subseteq \mathbf{G}_E$. Therefore, $V^{-1} \in \mathbf{G}_E$. Taking into account (2.8), we receive $U \in \mathbf{P}_B\mathbf{P}_A^{-1}$. □

**Remark.** If $\det \Phi = \det E = 0$, then the statement of Lemma 2.6 is false.

**Example 2.1.** Let $\Phi = \operatorname{diag}(\alpha, 0)$, $E = \operatorname{diag}(\varepsilon, 0)$ be nonzero matrices over $\mathbb{N}$. In this case, the groups $\mathbf{G}_\Phi$, $\mathbf{G}_E$ coincide and consist of all invertible matrices of the form

$$\left\| \begin{matrix} e_1 & a \\ 0 & e_2 \end{matrix} \right\|, e_1, e_2 = \pm 1, a \in \mathbb{N}.$$

It is easy to check that the set $\mathbf{P}_B\mathbf{P}_A^{-1}$, where

$$P_B P_A^{-1} = \left\| \begin{matrix} 7 & 3 \\ 5 & 2 \end{matrix} \right\|$$

does not contain a lower unitriangular matrix. ◇

## 2.8. Matrix rings of stable range 1.5

All rings of stable range 1.5 considered above were commutative rings without zero divisors. It will be shown that the second order matrix rings over commutative Bezout rings of stable range 1.5 have a similar stable range.

Let $A, B$ be $2 \times 2$ matrices over an elementary divisor ring $R$, which have the Smith forms

$$E = \operatorname{diag}(\varepsilon_1, \varepsilon_2), \quad \Delta = \operatorname{diag}(\delta_1, \delta_2),$$

respectively. Let $P_A \in \mathbf{P}_A$, $P_B \in \mathbf{P}_B$, and $P_B P_A^{-1} = \|s_{ij}\|_1^2 = S$.

**Lemma 2.3.** *Element $((\varepsilon_2, \delta_2), s_{21}[\varepsilon_1, \delta_1])$ is invariant with respect to the choice of transforming matrices $P_B$ and $P_A$.*

**Proof.** 1). Let $A, B$ be nonsingular matrices, and $F_A$ and $F_B$ be the others their left transforming matrices, i.e., $F_A \in \mathbf{P}_A$, $F_B \in \mathbf{P}_B$. There are such $H_A \in \in \mathbf{G}_E$ and $H_B \in \mathbf{G}_\Delta$, such $F_A = H_A P_A$, $F_B = H_B P_B$. Set $F_B F_A^{-1} = \|s'_{ij}\|_1^2$. To prove the Lemma it is necessary to show that

$$((\varepsilon_2, \delta_2), s_{21}[\varepsilon_1, \delta_1]) = ((\varepsilon_2, \delta_2), s'_{21}[\varepsilon_1, \delta_1]).$$

Consider the product of matrices

$$F_B F_A^{-1} = H_B P_B (H_A P_A)^{-1} = H_B P_B P_A^{-1} H_A^{-1} = H_B S H_A^{-1},$$

where $S = P_B P_A^{-1}$. Set $H_B S = \|k_{ij}\|_1^2$. Since

$$H_B = \left\| \begin{matrix} h_{11} & h_{12} \\ \dfrac{\delta_2}{\delta_1} h_{21} & h_{22} \end{matrix} \right\|,$$

than

$$k_{21} = \left\| \begin{matrix} \dfrac{\delta_2}{\delta_1} h_{21} & h_{22} \end{matrix} \right\| \left\| \begin{matrix} s_{11} \\ s_{21} \end{matrix} \right\| = \frac{\delta_2}{\delta_1} h_{21} s_{11} + h_{22} s_{21}.$$

Consider

$$(k_{21}[\varepsilon_1, \delta_1], (\varepsilon_2, \delta_2)) = \left( \left( \frac{\delta_2}{\delta_1} h_{21} s_{21} + h_{22} s_{21} \right) [\varepsilon_1, \delta_1], (\varepsilon_2, \delta_2) \right).$$

Since

$$\frac{\delta_2}{\delta_1} \cdot \frac{[\varepsilon_1, \delta_1]}{(\varepsilon_2, \delta_2)} = \frac{\delta_2 [\varepsilon_1, \delta_1]}{\delta_1 (\varepsilon_2, \delta_2)} \in R,$$

than

$$(\varepsilon_2, \delta_2) \left| \frac{\delta_2}{\delta_1} [\varepsilon_1, \delta_1]. \right.$$

Hence,

$$\left( \left( \frac{\delta_2}{\delta_1} h_{21} s_{21} + h_{22} s_{21} \right) [\varepsilon_1, \delta_1], (\varepsilon_2, \delta_2) \right) = (h_{22} s_{21}[\varepsilon_1, \delta_1], (\varepsilon_2, \delta_2)).$$

**79**

It follows from the invertibility of the matrix $H_B$ that

$$\left( h_{22}, \frac{\delta_2}{\delta_1} \right) = 1.$$

Consequently,

$$(h_{22}, (\varepsilon_2, \delta_2)) = 1.$$

That is

$$(k_{21}[\varepsilon_1, \delta_1], (\varepsilon_2, \delta_2)) = (s_{21}[\varepsilon_1, \delta_1], (\varepsilon_2, \delta_2)).$$

Consider the matrix $SH_A^{-1} = \|t_{ij}\|_1^2$. Since

$$H_A^{-1} = \left\| \begin{matrix} v_{11} & v_{12} \\ \dfrac{\varepsilon_2}{\varepsilon_1} v_{21} & v_{22} \end{matrix} \right\|,$$

so

$$t_{21} = \left\| \begin{matrix} s_{21} & s_{22} \end{matrix} \right\| \left\| \begin{matrix} \nu_{11} \\ \dfrac{\varepsilon_2}{\varepsilon_1} \nu_{21} \end{matrix} \right\| = s_{21}\nu_{11} + \frac{\varepsilon_2}{\varepsilon_1} \nu_{21} s_{22}.$$

Noting that

$$(\varepsilon_2, \delta_2) \left| \frac{\varepsilon_2}{\varepsilon_1} [\varepsilon_1, \delta_1] \right.$$

and thinking similarly, we get

$$(t_{21}[\varepsilon_1, \delta_1], (\varepsilon_2, \delta_2)) = (s_{21}[\varepsilon_1, \delta_1], (\varepsilon_2, \delta_2)).$$

By virtue of the associativity of the ring $M_2(R)$, we conclude consideration of this case.

2). Assume that the matrix $A$ is nonsingular and $B \sim \operatorname{diag}(\delta_1, \, 0)$, then

$$((\varepsilon_2, \delta_2), s_{21}[\varepsilon_1, \delta_1]) = (\varepsilon_2, s_{21}[\varepsilon_1, \delta_1]).$$

Let $F_A$ and $F_B$ be the other left transforming matrices of the matrices $A$ and $B$. Then there are $H_A \in \mathbf{G}_E$ and

$$H_B = \left\| \begin{matrix} e_1 & h_{12} \\ 0 & e_2 \end{matrix} \right\|$$

such that $F_A = H_A P_A$, $F_B = H_B P_B$. Set $F_B F_A^{-1} = \|s'_{ij}\|_1^2$. To prove the Lemma it is necessary to show that

$$(\varepsilon_2, s_{21}[\varepsilon_1, \delta_1]) = (\varepsilon_2, s'_{21}[\varepsilon_1, \delta_1]).$$

Consider the product of matrices

$$F_B F_A^{-1} = H_B P_B (H_A P_A)^{-1} = H_B P_B P_A^{-1} H_A^{-1} = H_B S H_A^{-1},$$

where $S = P_B P_A^{-1}$. Set $H_B S = \|k_{ij}\|_1^2$. Since

$$H_B = \begin{Vmatrix} e_1 & h_{12} \\ 0 & e_2 \end{Vmatrix},$$

so

$$k_{21} = \begin{Vmatrix} 0 & e_2 \end{Vmatrix} \begin{Vmatrix} s_{11} \\ s_{21} \end{Vmatrix} = e_2 s_{21}.$$

Consider

$$(k_{21}[\varepsilon_1, \delta_1], \varepsilon_2) = (e_2 s_{21}[\varepsilon_1, \delta_1], \varepsilon_2).$$

Since $e_2$ is an invertible ring element, then

$$(k_{21}[\varepsilon_1, \delta_1], \varepsilon_2) = (s_{21}[\varepsilon_1, \delta_1], \varepsilon_2).$$

Consider the matrix $S H_A^{-1} = \|t_{ij}\|_1^2$. Since

$$H_A^{-1} = \begin{Vmatrix} v_{11} & v_{12} \\ \dfrac{\varepsilon_2}{\varepsilon_1} v_{21} & v_{22} \end{Vmatrix},$$

then

$$t_{21} = \begin{Vmatrix} s_{21} & s_{22} \end{Vmatrix} \begin{Vmatrix} \nu_{11} \\ \dfrac{\varepsilon_2}{\varepsilon_1} \nu_{21} \end{Vmatrix} = s_{21} \nu_{11} + \frac{\varepsilon_2}{\varepsilon_1} \nu_{21} s_{22}.$$

Noting that

$$\varepsilon_2 \left| \frac{\varepsilon_2}{\varepsilon_1} [\varepsilon_1, \delta_1] \right.$$

and considering the same as in the case 1), we obtain

$$(t_{21}[\varepsilon_1, \delta_1], \varepsilon_2) = (s_{21}[\varepsilon_1, \delta_1], \varepsilon_2).$$

By virtue of the associativity of the ring $M_2(R)$, we conclude consideration of this case.

3). Let $B$ be nonsingular and $A \sim \mathrm{diag}(\varepsilon_1, 0)$, then

$$((\varepsilon_2, \delta_2), s_{21}[\varepsilon_1, \delta_1]) = (\delta_2, s_{21}[\varepsilon_1, \delta_1]).$$

Let $F_A$ and $F_B$ be the other left transforming matrices of the matrices $A$ and $B$. Then there are

$$H_A = \begin{Vmatrix} e_1 & l_{12} \\ 0 & e_2 \end{Vmatrix}, \quad H_B \in \mathbf{G}_\Delta,$$

such that $F_A = H_A P_A$, $F_B = H_B P_B$. Set $F_B F_A^{-1} = \|s'_{ij}\|_1^2$. To prove the Lemma it is necessary to show that

$$(\delta_2, s_{21}[\varepsilon_1, \delta_1]) = (\delta_2, s'_{21}[\varepsilon_1, \delta_1]).$$

Consider the product of matrices

$$F_B F_A^{-1} = H_B P_B (H_A P_A)^{-1} = H_B P_B P_A^{-1} H_A^{-1} = H_B S H_A^{-1},$$

where $S = P_B P_A^{-1}$. Set $H_B S = \|k_{ij}\|_1^2$. Since

$$H_B = \left\| \begin{matrix} h_{11} & h_{12} \\ \dfrac{\delta_2}{\delta_1} h_{21} & h_{22} \end{matrix} \right\|,$$

then

$$k_{21} = \left\| \dfrac{\delta_2}{\delta_1} h_{21} \quad h_{22} \right\| \left\| \begin{matrix} s_{11} \\ s_{12} \end{matrix} \right\| = \dfrac{\delta_2}{\delta_1} h_{21} s_{11} + h_{22} s_{12}.$$

Noting that

$$\delta_2 \left| \dfrac{\delta_2}{\delta_1} [\varepsilon_1, \delta_1] \right.$$

and considering the same as in the case 1), we obtain

$$(k_{21}[\varepsilon_1, \delta_1], \delta_2) = (s_{21}[\varepsilon_1, \delta_1], \delta_2).$$

Consider $SH_A^{-1} = \|t_{ij}\|_1^2$. Since

$$H_A^{-1} = \left\| \begin{matrix} e_1 & -l_{12} \\ 0 & e_2 \end{matrix} \right\|,$$

then

$$t_{21} = \left\| s_{21} \quad s_{22} \right\| \left\| \begin{matrix} e_1 \\ 0 \end{matrix} \right\| = e_1 s_{21}.$$

Consider

$$(t_{21}[\varepsilon_1, \delta_1], \delta_2) = (e_1 s_{21}[\varepsilon_1, \delta_1], \delta_2).$$

Since $e_1$ is an invertible ring element, so

$$(t_{21}[\varepsilon_1, \delta_1], \delta_2) = (s_{21}[\varepsilon_1, \delta_1], \delta_2).$$

By virtue of the associativity of the ring $M_2(R)$, we conclude consideration of this case.

4). Suppose that

$$A \sim \left\| \begin{matrix} \varepsilon_1 & 0 \\ 0 & 0 \end{matrix} \right\|, \quad B \sim \left\| \begin{matrix} \delta_1 & 0 \\ 0 & 0 \end{matrix} \right\|$$

and $F_A, F_B$ are the other left transforming matrices of these matrices. Then there are those

$$H_A = \left\| \begin{matrix} e_1' & v_{12} \\ 0 & e_2' \end{matrix} \right\|, \quad H_B = \left\| \begin{matrix} e_1 & h_{12} \\ 0 & e_2 \end{matrix} \right\|$$

such that $F_A = H_A P_A$, $F_B = H_B P_B$. Set $F_B F_A^{-1} = \|s'_{ij}\|_1^2 = S'$. To prove the Lemma it is necessary to show that $s_{21}$ and $s'_{21}$ are associates. Consider the product of matrices

$$S' = F_B F_A^{-1} = H_B P_B (H_A P_A)^{-1} = H_B P_B P_A^{-1} H_A^{-1} = H_B S H_A^{-1},$$

where $S = P_B P_A^{-1}$. Write these matrices explicitly, that is

$$\left\| \begin{matrix} s'_{11} & s'_{12} \\ s'_{21} & s'_{22} \end{matrix} \right\| = \left\| \begin{matrix} e_1 & h_{12} \\ 0 & e_2 \end{matrix} \right\| \left\| \begin{matrix} s_{11} & s_{12} \\ s_{21} & s_{22} \end{matrix} \right\| \left\| \begin{matrix} e'^{-1}_1 & * \\ 0 & e'^{-1}_2 \end{matrix} \right\| =$$

$$= \left\| \begin{matrix} s''_{11} & s''_{12} \\ s_{21} u & s''_{22} \end{matrix} \right\|,$$

where $u = e_2 e'^{-1}_2$ is an invertible element of the ring $R$. Therefore, $s'_{21} = s_{21} u$, which is to be proved. $\qquad\square$

**Theorem 2.15.** *The Smith form of the matrix* $(A, \ B)_l$ *is a matrix*

$$\mathrm{diag}\,((\varepsilon_1, \delta_1), \ \ (\varepsilon_2, \delta_2, [\varepsilon_1, \delta_1] s_{21})).$$

**Proof.** Consider the matrix $\|A \quad B\|$. Based on Theorem 1.10, there exists an invertible matrix $U$, that

$$\|A \quad B\| U = \|D \quad \mathbf{0}\|,$$

moreover $D = (A, B)_l$. Equalities

$$\|A \quad B\| =$$

$$= \|P_A^{-1} \mathrm{E} Q_A^{-1} \quad P_B^{-1} \Delta Q_B^{-1}\| = P_B^{-1} \|P_B P_A^{-1} \mathrm{E} \quad \Delta\| \left\| \begin{matrix} Q_A^{-1} & \mathbf{0} \\ \mathbf{0} & Q_B^{-1} \end{matrix} \right\| =$$

$$= P_B^{-1} \|S\mathrm{E} \quad \Delta\| \left\| \begin{matrix} Q_A^{-1} & \mathbf{0} \\ \mathbf{0} & Q_B^{-1} \end{matrix} \right\|.$$

are fulfilled. Right-multiply this equality by $U$:

$$\|A \quad B\| U = \|D \quad \mathbf{0}\| = \left( P_B^{-1} \right) \|S\mathrm{E} \quad \Delta\| \left( \left\| \begin{matrix} Q_A^{-1} & \mathbf{0} \\ \mathbf{0} & Q_B^{-1} \end{matrix} \right\| U \right).$$

Hence,

$$\|A \quad B\| \sim \|D \quad \mathbf{0}\| \sim \|S\mathrm{E} \quad \Delta\| \sim \left\| \begin{matrix} \nu_1 & 0 & 0 & 0 \\ 0 & \nu_2 & 0 & 0 \end{matrix} \right\|$$

is the Smith form of these matrices. Based on Theorem 2.2, $\nu_1$ is equal to g.c.d. all elements of the matrices $S\mathrm{E}$, $\Delta$ and

$$\nu_2 = \frac{\mu_2}{\nu_1},$$

where $\mu_2$ is g.c.d. of second-order minors of these matrices. It follows that $D = (A, B)_l$ have the Smith form $\mathrm{diag}(\nu_1, \nu_2)$. Since g.c.d. of elements of the matrix $SE$ is $\varepsilon_1$, and $\Delta$ is $\delta_1$, then

$$\nu_1 = (\varepsilon_1, \delta_1).$$

Consider the matrix

$$\left\| SE \quad \Delta \right\| = \left\| \begin{matrix} s_{11}\varepsilon_1 & s_{12}\varepsilon_2 & \delta_1 & 0 \\ s_{21}\varepsilon_1 & s_{22}\varepsilon_2 & 0 & \delta_2 \end{matrix} \right\|.$$

Noting that $s_{ij}$ are elements of the invertible matrix $S$, we get

$$\mu_2 = (\varepsilon_1\varepsilon_2, \delta_1\delta_2, \varepsilon_1\delta_2 s_{11}, \varepsilon_2\delta_2 s_{12}, \varepsilon_1\delta_1 s_{21}, \varepsilon_2\delta_1 s_{22}) =$$

$$= \left( \varepsilon_1\varepsilon_2, \delta_1\delta_2, \varepsilon_1\delta_2 \left( s_{11}, \frac{\varepsilon_2}{\varepsilon_1} s_{12} \right), \varepsilon_1\delta_1 \left( s_{21}, \frac{\varepsilon_2}{\varepsilon_1} s_{22} \right) \right) =$$

$$= \left( \varepsilon_1\varepsilon_2, \delta_1\delta_2, \varepsilon_1\delta_2 \left( s_{11}, \frac{\varepsilon_2}{\varepsilon_1} \right), \varepsilon_1\delta_1 \left( s_{21}, \frac{\varepsilon_2}{\varepsilon_1} \right) \right) =$$

$$= \left( \varepsilon_1\varepsilon_2, \delta_1\delta_2, \varepsilon_1\delta_1 \left( \frac{\delta_2}{\delta_1} \left( s_{11}, \frac{\varepsilon_2}{\varepsilon_1} \right), s_{21}, \frac{\varepsilon_2}{\varepsilon_1} \right) \right) =$$

$$= \left( \varepsilon_1\varepsilon_2, \delta_1\delta_2, \varepsilon_1\delta_1 \left( \frac{\delta_2}{\delta_1} s_{11}, \frac{\delta_2\varepsilon_2}{\delta_1\varepsilon_1}, s_{21}, \frac{\varepsilon_2}{\varepsilon_1} \right) \right) =$$

$$= \left( \varepsilon_1\varepsilon_2, \delta_1\delta_2, \varepsilon_1\delta_1 \left( \left( \frac{\delta_2}{\delta_1} s_{11}, s_{21} \right), \left( \frac{\delta_2\varepsilon_2}{\delta_1\varepsilon_1}, \frac{\varepsilon_2}{\varepsilon_1} \right) \right) \right) =$$

$$= \left( \varepsilon_1\varepsilon_2, \delta_1\delta_2, \varepsilon_1\delta_1 \left( \frac{\delta_2}{\delta_1} s_{21}, \frac{\varepsilon_2}{\varepsilon_1} \right) \right) =$$

$$= (\varepsilon_1\varepsilon_2, \varepsilon_2\delta_1, \delta_1\delta_2, \varepsilon_1\delta_2, \varepsilon_1\delta_1 s_{21}) =$$

$$= (\varepsilon_1, \delta_1) \left( \varepsilon_2 \left( \frac{\varepsilon_1}{(\varepsilon_1, \delta_1)}, \frac{\delta_1}{(\varepsilon_1, \delta_1)} \right), \delta_2 \left( \frac{\varepsilon_1}{(\varepsilon_1, \delta_1)}, \frac{\delta_1}{(\varepsilon_1, \delta_1)} \right), \frac{\varepsilon_1\delta_1}{(\varepsilon_1, \delta_1)} s_{21} \right) =$$

$$= (\varepsilon_1, \delta_1) (\varepsilon_2, \delta_2, [\varepsilon_1, \delta_1] s_{21}).$$

Consequently,

$$\nu_2 = \frac{\mu_2}{\nu_1} = (\varepsilon_2, \delta_2, [\varepsilon_1, \delta_1] s_{21}).$$

Finally, let us look at Lemma 2.3, which shows that $(\varepsilon_2, \delta_2, [\varepsilon_1, \delta_1] s_{21})$ is invariance with respect to the choice of transforming matrices $P_B$ and $P_A$. $\square$

**Corollary 2.7.** *In order that* $(A, B)_l = I$, *it is necessary and sufficient that*

$$(\varepsilon_2, \delta_2, [\varepsilon_1, \delta_1] s_{21}) = 1.$$

**Proof.** In order that $(A, B)_l = I$, it is necessary and sufficient that

$$(\varepsilon_1, \delta_1) (\varepsilon_2, \delta_2, [\varepsilon_1, \delta_1] s_{21}) = 1. \tag{2.9}$$

Since $(\varepsilon_1, \delta_1)$ is a divisor of $(\varepsilon_2, \delta_2, [\varepsilon_1, \delta_1]s_{21})$, the equality (2.9) is equivalent to the condition

$$(\varepsilon_2, \delta_2, [\varepsilon_1, \delta_1]s_{21}) = 1. \qquad \square$$

**Corollary 2.8.** *If $\varepsilon_2 = 0$, $\delta_2 \neq 0$, then $(A, B)_l = I$ if and only if*

$$\delta_1 = 1 \quad and \quad (\delta_2, \varepsilon_1 s_{21}) = 1. \qquad \square$$

**Corollary 2.9.** *If $\varepsilon_2$, $\delta_2 = 0$, then $(A, B)_l = I$ if and only if $\varepsilon_1 \delta_1 s_{21}$ is an invertible element of the ring $R$.* $\qquad \square$

**Corollary 2.10.** *If $(A, B)_l = I$ and $s_{21} = 0$, then the matrices $A$, $B$ are nonsingular and, in addition,*

$$(\det A, \det B) = 1. \qquad \square$$

**Theorem 2.16.** *Let $A, B \in M_2(R)$ and let at least one of the matrices be nonsingular. Then the set $\mathbf{P}_B \mathbf{P}_A^{-1}$ contains a lower unitriangular matrix.*

**Proof.** We preserve the notation of Theorem 2.7. Since $\mathbf{P}_B = \mathbf{G}_\Delta P_B$, $\mathbf{P}_A = \mathbf{G}_E P_A$, we find

$$\mathbf{P}_B \mathbf{P}_A^{-1} = \mathbf{G}_\Delta P_B (\mathbf{G}_E P_A)^{-1} = \mathbf{G}_\Delta P_B P_A^{-1} \mathbf{G}_E = \mathbf{G}_\Delta S \mathbf{G}_E.$$

Left-multiply the matrix $P_B P_A^{-1} = S$ by elements of group $\mathbf{G}_\Delta$ and right-multiply by elements from $\mathbf{G}_E$, we remain inside the set $\mathbf{P}_B \mathbf{P}_A^{-1}$.

Let $\det \Delta \neq 0$. By Theorem 2.14

$$\mathrm{GL}_2(R) = \mathbf{G}_\Delta U_2^{\mathrm{lw}}(R) U_2^{\mathrm{up}}(R).$$

Since $S \in \mathrm{GL}_2(R)$, then

$$S = HUV, \ \ H \in \mathbf{G}_\Delta, \ \ U \in U_2^{\mathrm{lw}}(R), \ \ V \in U_2^{\mathrm{up}}(R).$$

Consequently, $U = H^{-1} S V^{-1}$. Noting that $H^{-1} \in \mathbf{G}_\Delta$, $V^{-1} \in \mathbf{G}_E$, we conclude that $U \in \mathbf{P}_B \mathbf{P}_A^{-1}$.

Let $\det E \neq 0$. The equality

$$\mathrm{GL}_2(R) = \mathbf{G}_E U_2^{\mathrm{lw}}(R) U_2^{\mathrm{up}}(R)$$

is fulfilled. Going to the invertible matrices, we get

$$\mathrm{GL}_2(R) = U_2^{\mathrm{up}}(R) U_2^{\mathrm{lw}}(R) \mathbf{G}_E.$$

That is

$$S = MNK, \ \ M \in U_2^{\mathrm{up}}(R), \ \ N \in U_2^{\mathrm{lw}}(R), \ \ K \in \mathbf{G}_E.$$

Hence, $N = M^{-1} S K^{-1}$. Since $U_2^{\mathrm{up}}(R) \subset \mathbf{G}_\Delta$, then $M^{-1} \in \mathbf{G}_\Delta$. Noting that $K^{-1} \in U_2^{\mathrm{up}}(R)$, we conclude that $N \in \mathbf{P}_B \mathbf{P}_A^{-1}$. The theorem is proved. $\qquad \square$

To prove the following result, we need to use the following theorem. Its proof is cumbersome and therefore not included in the monograph.

**Theorem 2.17.** *Let* $A = P_A^{-1}\mathrm{E}Q_A^{-1}$ *and* $B = P_B^{-1}\Phi Q_B^{-1}$. *In order that*

$$AB \sim \mathrm{E}\Phi$$

*it is necessary and sufficient that* $Q_A^{-1}P_B^{-1} = K_1 K_2$, *where* $K_1 \in \mathbf{G}_{\mathrm{E}}^T$ *and* $K_2 \in \mathbf{G}_\Phi$.

The proof of this Theorem can be found in [103].

**Theorem 2.18.** *Let* $A = P_A^{-1}\mathrm{E}Q_A^{-1}$, $B = P_B^{-1}\Phi Q_B^{-1}$. *If* $AB \sim \mathrm{E}\Delta$, *then* $\boldsymbol{P}_{AB} \subseteq \boldsymbol{P}_A$.

**Proof**. Since,

$$AB = P_A^{-1}\mathrm{E}Q_A^{-1}\, P_B^{-1}\mathrm{E}Q_B^{-1} \sim \mathrm{E}\Delta$$

then by virtue of Theorem 2.17, $Q_A^{-1}P_B^{-1}$ can be represented in the form $Q_A^{-1}P_B^{-1} = K_1 K_2$, where $K_1 \in \mathbf{G}_{\mathrm{E}}^T$. This yields $\mathrm{E}K_1 = L_1\mathrm{E}$, where $L_1 \in \mathbf{G}_{\mathrm{E}}$ and $K_2 \in \mathbf{G}_\Delta$. Then

$$AB = P_A^{-1}\mathrm{E}(Q_A^{-1}P_B^{-1})\Delta Q_B^{-1} = P_A^{-1}\mathrm{E}K_1 K_2 \Delta Q_B^{-1} = (L_1^{-1}P_A)^{-1}\mathrm{E}\Delta(Q_B L_2^{-1})^{-1}.$$

This implies that $L_1^{-1}P_A \in \boldsymbol{P}_{AB}$. Hence, $\mathbf{P}_{AB} = \mathbf{G}_{\mathrm{E}\Delta}(L_1^{-1}P_A)$ and, in turn $\mathbf{P}_A = \mathbf{G}_{\mathrm{E}}P_A$. Noting that $L_1^{-1} \in \mathbf{G}_{\mathrm{E}}$, we get

$$\mathbf{G}_{\mathrm{E}}(L_1^{-1}P_A) = (\mathbf{G}_{\mathrm{E}}L_1^{-1})P_A = \mathbf{G}_{\mathrm{E}}P_A = \mathbf{P}_A.$$

Since $\mathbf{G}_{\mathrm{E}\Delta} \subseteq \mathbf{G}_{\mathrm{E}}$, we arrive at the conclusion that $\boldsymbol{P}_{AB} \subseteq \boldsymbol{P}_A$. □

**Theorem 2.19.** *Let* $A, B$ *be nonzero left relatively prime* $2 \times 2$ *matrices over a Bezout ring* $R$ *of stable range 1.5. If the matrix* $A$ *is singular, then there exists a matrix* $T$ *such that*

$$A + BT \in \mathrm{GL}_2(R).$$

*If the matrix* $A$ *is nonsingular, then, for any fixed* $\varphi \neq 0$ *from* $R$, *there exists a matrix* $F_\varphi$ *such that*

$$A + BF_\varphi \sim \begin{Vmatrix} 1 & 0 \\ 0 & \mu \end{Vmatrix},$$

*moreover* $(\mu, \varphi) = 1$.

To prove the theorem, we need the following assertions.

**Lemma 2.4.** *The matrices* $(A, B)_l$ *and* $P_B^{-1}(P_B P_A^{-1}\mathrm{E}, \Delta)_l$ *are right associated.*

**Proof.** Denote $P_B P_A^{-1} = S$. As indicated in the proof of Theorem 2.15 $(S\mathrm{E}, \Delta)_l = D$, where

$$\begin{Vmatrix} S\mathrm{E} & \Delta \end{Vmatrix} U = \begin{Vmatrix} D & \mathbf{0} \end{Vmatrix},$$

$U \in \mathrm{GL}_4(R)$. Multiplying this equality from the left by $P_B^{-1}$, we obtain

$$\begin{Vmatrix} P_A^{-1}\mathrm{E}Q_A^{-1} & P_B^{-1}\Delta Q_B^{-1} \end{Vmatrix} \left( \begin{Vmatrix} Q_A & \mathbf{0} \\ \mathbf{0} & Q_B \end{Vmatrix} U \right) = \begin{Vmatrix} P_B^{-1}D & \mathbf{0} \end{Vmatrix}.$$

That is $(A, B)_l = P_B^{-1}D = P_B^{-1}(S\mathrm{E}, \Delta)_l$. □

**Lemma 2.5.** *The matrices* $P_B P_A^{-1}\mathrm{E} + \Delta V$ *and* $A + B(Q_B V Q_A^{-1})$ *are equivalent.*

**Proof.** Indeed,

$$A + B(Q_B V Q_A^{-1}) = P_A^{-1}\mathrm{E}Q_A^{-1} + P_B^{-1}\Delta Q_B^{-1}(Q_B V Q_A^{-1}) =$$
$$= P_B^{-1}(P_B P_A^{-1}\mathrm{E} + \Delta V)Q_A^{-1}. \qquad \square$$

**Proof Theorem 2.19.** Assume that $\det A = 0$, i.e., $A = P_A^{-1}\mathrm{E}Q_A^{-1}$, $\mathrm{E} = \mathrm{diag}(\varepsilon_1, 0)$. In view of Corollary 2.8, the matrix $B$ has the form $B = P_B^{-1}\Delta Q_B^{-1}, \Delta = \mathrm{diag}(1, \delta_2)$. Let $V = \|v_{ij}\|_1^2$ be a parametric matrix.

Since $(A, B)_l = I$, it follows from Lemma 2.4 that $(S\mathrm{E}, \Delta)_l = I$.

1). Let $\delta_2 \neq 0$. By Theorem 2.16, it is possible to choose matrices $P_A, P_B$ such that

$$P_B P_A^{-1} = \left\| \begin{matrix} 1 & 0 \\ s & 1 \end{matrix} \right\| = S.$$

In view of the fact that $P_{S\mathrm{E}} = S^{-1}$ and $P_\Delta = I$, we obtain

$$P_{S\mathrm{E}} P_\Delta^{-1} = S^{-1} = \left\| \begin{matrix} 1 & 0 \\ -s & 1 \end{matrix} \right\|.$$

Since $(S\mathrm{E}, \Delta)_l = I$, in view of Corollary 2.8, we find

$$(\delta_2, -s\varepsilon_1) = 1 \Rightarrow (\delta_2, s\varepsilon_1) = 1.$$

Consider the equality

$$S\mathrm{E} + \Delta V = \left\| \begin{matrix} \varepsilon_1 + v_{11} & v_{12} \\ s\varepsilon_1 + v_{21}\delta_2 & v_{22}\delta_2 \end{matrix} \right\|.$$

Set $v_{21}^0 = 0, v_{22}^0 = 1$. There exist $m$ and $n$ such that

$$m\delta_2 - ns\varepsilon_1 = 1.$$

Set $v_{11}^0 = m - \varepsilon_1$, $v_{12}^0 = n$. The matrix $\|v_{ij}^0\|_1^2$ is denoted by $V^0$. The obtained matrix $S\mathrm{E} + \Delta V^0$ is invertible. By Lemma 2.5, we get $A + BU \in \mathrm{GL}_2(R)$, where $U = Q_B V^0 Q_A^{-1}$.

2). Assume that $\delta_2 = 0$ and $P_B P_A^{-1} = \|s_{ij}\|_1^2$. Then

$$S\mathrm{E} + \Delta V = \left\| \begin{matrix} s_{11}\varepsilon_1 + v_{11} & v_{12} \\ s_{21}\varepsilon_1 & 0 \end{matrix} \right\|.$$

Since

$$S^{-1} = e^{-1} \left\| \begin{matrix} s_{22} & -s_{12} \\ -s_{21} & s_{11} \end{matrix} \right\|,$$

where $e = \det S \in U(R)$, ), by using the previous arguments and Corollary 2.9, we get

$$-e\varepsilon_1 s_{21} \in U(R) \Rightarrow \varepsilon_1 s_{21} \in U(R).$$

Set

$$V^0 = \begin{Vmatrix} 0 & 1 \\ 0 & 0 \end{Vmatrix}.$$

Then $SE + \Delta V^0 \in \mathrm{GL}_2(R)$. According to Lemma 2.5, there exists a matrix $T$ such that $A + BT \in \mathrm{GL}_2(R)$.

3). Assume that $\det A \neq 0, \det B \neq 0$. We choose matrices $P_A, P_B$ such that

$$P_B P_A^{-1} = \begin{Vmatrix} 1 & 0 \\ s & 1 \end{Vmatrix} = S.$$

In the parametric matrix $V$, we set $v_{22} = 0$. Consider the equality

$$SE + \Delta V = \begin{Vmatrix} \varepsilon_1 + \delta_1 v_{11} & \delta_1 v_{12} \\ s\varepsilon_1 + v_{21}\delta_2 & \varepsilon_2 \end{Vmatrix}.$$

Since $(\varepsilon_1, \delta_1)|(\varepsilon_2, \delta_2, [\varepsilon_1, \delta_1]s)$, we have $(\varepsilon_1, \delta_1) = 1$. Then

$$(\varepsilon_2, \delta_2, [\varepsilon_1, \ \delta_1]s) = (\varepsilon_2, \delta_2, \varepsilon_1\delta_1 s) = 1. \tag{2.10}$$

Since $\delta_1|\delta_2$, we conclude that $(\varepsilon_2, \delta_1\delta_2, \varepsilon_1\delta_1 s) = 1$. Thus, there exists $u$ such that

$$(\varepsilon_1\delta_1 s + \delta_1\delta_2 u, \varepsilon_2) = 1.$$

Set $v_{21}^0 = u, v_{22}^0 = 0$. Then

$$\begin{Vmatrix} \varepsilon_1 s + \delta_2 v_{21}^0 & \varepsilon_2 \end{Vmatrix} = \begin{Vmatrix} a & \varepsilon_2 \end{Vmatrix} \sim \begin{Vmatrix} 1 & 0 \end{Vmatrix}. \tag{2.11}$$

It follows from (2.10) that $(\varepsilon_1\varepsilon_2, \delta_1) = 1$. Thus, for any $\varphi \neq 0$ from $R$, we obtain $(\varepsilon_1\varepsilon_2, \delta_1, \varphi) = 1$. Hence, there exists $r$ such that

$$(\varepsilon_1\varepsilon_2 + \delta_1 r, \varphi) = 1.$$

There are $p, q$ such that

$$p(\varepsilon_1\varepsilon_2 + \delta_1 r) + q\varphi = 1.$$

In the parametric matrix $V$, we set $v_{21} = v_{22} = 0$. Consider

$$\det \begin{Vmatrix} \varepsilon_1 + \delta_1 v_{11} & \delta_1 v_{12} \\ a & \varepsilon_2 \end{Vmatrix} = \varepsilon_1\varepsilon_2 + \delta_1(\varepsilon_2 v_{11} - av_{12}).$$

Since $(a, \varepsilon_2) = 1$, there exist $v_{11}^0, v_{12}^0$ such that $\varepsilon_2 v_{11}^0 - av_{12}^0 = r$. Denote $\|v_{ij}^0\|_1^2 = V^0$. Then $(\det(SE + \Delta V^0), \varphi) = 1$. Moreover, it follows from (2.11) that

$$SE + \Delta V^0 \sim \mathrm{diag}(1, \mu), \quad \mu = \det(SE + \Delta V^0).$$

To complete the analysis of this case, we use Lemma 2.5.

4). Assume that $\det A \neq 0, \det B = 0$. In this case,

$$A = P_A^{-1} \mathrm{E} Q_A^{-1}, \quad \mathrm{E} = \mathrm{diag}(1, \varepsilon_2), \quad B = P_B^{-1} \Delta Q_B^{-1},$$

$$\Delta = \mathrm{diag}(\delta_1, 0), \quad P_B P_A^{-1} = \left\| \begin{matrix} 1 & 0 \\ s & 1 \end{matrix} \right\| = S.$$

Consider the equality

$$S\mathrm{E} + \Delta V = \left\| \begin{matrix} 1 + \delta_1 v_{11} & \delta_1 v_{12} \\ s & \varepsilon_2 \end{matrix} \right\|.$$

Then

$$\det(S\mathrm{E} + \Delta V) = \varepsilon_2 + (v_{11}\varepsilon_2 - sv_{12})\,\delta_1.$$

Since $(S\mathrm{E}, \Delta)_l = I$, in view of Corollary 2.8, we get

$$(\varepsilon_2, -\delta_1 s) = (\varepsilon_2, \delta_1 s) = 1. \tag{2.12}$$

Consequently, $(\varepsilon_2, \delta_1) = 1$. Then, for any $\varphi \neq 0$ from $R$, we find $(\varepsilon_2, \delta_1, \varphi) = 1$. Thus, there exists $d$ such that

$$(\varepsilon_2 + \delta_1 d, \varphi) = 1.$$

It follows from the equality (2.12) that $(\varepsilon_2, s) = 1$. We choose $v_{11}^0, v_{12}^0$ such that

$$v_{11}^0 \varepsilon_2 - sv_{12}^0 = d.$$

Denote

$$\left\| \begin{matrix} v_{11}^0 & v_{12}^0 \\ 0 & 0 \end{matrix} \right\| = V^0.$$

Then $(\det(S\mathrm{E} + \Delta V^0), \varphi) = 1$. Hence, according to Lemma 2.5

$$A + B(Q_B V^0 Q_A^{-1}) \sim \mathrm{diag}(1, \ \nu),$$

where $(\nu, \varphi) = 1$. $\qquad\square$

**Corollary 2.11.** *Let $A, B$ be nonzero left relatively prime matrices from $M_2(R)$ and let $C$ be a nonsingular matrix. Then there exists a matrix $F_C$ such that*

$$(\det(A + BF_C), \det C) = 1. \qquad\square$$

**Lemma 2.6.** *Let $A, B$ be $2 \times 2$ matrices over an elementary divisor ring. Then there exist $A_1, B_1$ such that*

$$A = (A, B)_l A_1, \ B = (A, B)_l B_1,$$

*moreover $(A_1, B_1)_l = I$.*

**Proof.** Suppose that $(A, B)_l$ be a nonsingular matrix and let $(A, B)_l = D \notin GL_2(R)$. Hence, $(A, B)_l D$ is a left common divisor of the matrices $A, B$. According to the definition of a left common divisor, the matrix $(A, B)_l D$ is a left divisor of the matrix $(A, B)_l$:

$$(A, B)_l = (A, B)_l DU$$

for some matrix $U$, i.e., $(A, B)_l(I - DU) = 0$. Since $(A, B)_l$ is a nonsingular matrix, we get $DU = I$, i.e., $D \in GL_2(R)$. We arrive at a contradiction.

Let $(A, B)_l$ be a singular matrix. Hence, the matrices $A, B$ are also singular and have the form

$$A = P_A^{-1} \mathrm{E} Q_A^{-1}, \ \ \mathrm{E} = \mathrm{diag}(\varepsilon_1, 0), \ \ B = P_B^{-1} \Delta Q_B^{-1}, \ \ \Delta = \mathrm{diag}(\delta_1, 0).$$

It follows from Theorem 2.7 that the matrix $(A, B)_l$ is singular if and only if $s_{21} = 0$, where $P_B P_A^{-1} = S = \|s_{ij}\|_1^2$, i.e., $P_B P_A^{-1} \in \mathbf{G}_\mathrm{E} = \mathbf{G}_\Delta$ is a group of invertible upper triangular matrices. This implies that $\mathbf{P}_B \cap \mathbf{P}_A \neq \varnothing$. Let $P \in \mathbf{P}_B \cap \mathbf{P}_A$. Hence, the matrices $A$ and $B$ have the form $A = P^{-1} \mathrm{E} Q_A^{-1}$, $B = P^{-1} \Delta Q_B^{-1}$. Thus, the following factorization of the matrices:

$$A = \left( P^{-1} \left\| \begin{matrix} (\varepsilon_1, \delta_1) & 0 \\ 0 & 0 \end{matrix} \right\| \right) \left( \left\| \begin{matrix} \dfrac{\varepsilon_1}{(\varepsilon_1, \delta_1)} & 0 \\ u & 0 \end{matrix} \right\| Q_A^{-1} \right),$$

$$B = \left( P^{-1} \left\| \begin{matrix} (\varepsilon_1, \delta_1) & 0 \\ 0 & 0 \end{matrix} \right\| \right) \left( \left\| \begin{matrix} \dfrac{\delta_1}{(\varepsilon_1, \delta_1)} & 0 \\ v & 0 \end{matrix} \right\| Q_B^{-1} \right),$$

where

$$\frac{\varepsilon_1}{(\varepsilon_1, \delta_1)} v - \frac{\delta_1}{(\varepsilon_1, \delta_1)} u = 1,$$

is just the required representation. $\qquad \square$

**Theorem 2.20.** *Let $R$ be a commutative Bezout domain. In order that a ring $M_2(R)$ have the stable range 1.5, it is necessary and sufficient that the ring $R$ has the same stable range 1.5.*

**Proof. Sufficiency.** Let $A, B$, and $C$ be nonzero left relatively prime matrices from $M_2(R)$.

1). Suppose that $\det A = 0$ and $(A, B)_l = D$. Then $A = DA_1$, $B = DB_1$. By virtue of Lemma 2.6, the matrices $A_1, B_1$ can be chosen to guarantee that $(A_1, B_1)_l = I$. If the matrix $D$ is nonsingular, then the matrix $A_1$ is singular. If the matrix $D$ is singular, then the matrix $B$ is singular. By Lemma 2.6, the matrix $A_1$ can also be chosen as singular. So, we will assume that in this case always $\det A_1 = 0$.

By using Theorem 2.19, we determine $T$ such that

$$A_1 + B_1 T = U \in \mathrm{GL}_2(R).$$

Since $P_D = P_{DU}$ and $(D, C)_l = I$, it follows from Theorem 2.7 that $(A + BT, C)_l = I$.

2). Suppose that $\det A \neq 0$, $\det C \neq 0$ and $(A, B)_l = I$. By virtue of Corollary 2.11, there exists a matrix $T$ such that

$$(\det(A + BT), \det C) = 1.$$

Then $(A + BT, C)_l = I$.

Let $(A, B)_l = D \neq I$. Hence, $A = DA_1$, $B = DB_1$, where

$$D = P_D^{-1}\Gamma Q_D^{-1}, \ \ \Gamma = \mathrm{diag}(\gamma_1, \gamma_2), \ \ C = P_C^{-1}\Omega Q_C^{-1}, \ \ \Omega = \mathrm{diag}(\omega_1, \omega_2).$$

By Corollary 2.11, there exists a matrix $L$ such that

$$(\det(A_1 + B_1L), \det D \det C) = 1. \tag{2.13}$$

This implies that $(A_1 + B_1L, C)_l = I$. Consider the matrix

$$D(A_1 + B_1L) = F.$$

Let $\Psi = \mathrm{diag}(\psi_1, \psi_2)$ be the Smith form of the matrix $A_1 + B_1L$. By using (2.13) and Theorem 4.4 (p. 148), we obtain

$$F = D(A_1 + B_1L) \sim \Delta\Psi = \mathrm{diag}(\gamma_1\psi_1, \ \gamma_2\psi_2) = \Gamma.$$

By virtue of Theorem 2.18, $\mathbf{P}_F \subseteq \mathbf{P}_D$. This implies that the matrix $F$ can be rewritten in the form $F = P_D^{-1}\Gamma Q_F^{-1}$. In view of (2.13), we get

$$(\psi_1\psi_2, \gamma_1\gamma_2) = (\psi_1\psi_2, \omega_1\omega_2) = 1.$$

Hence,

$$(\psi_2, \omega_2) = (\psi_1, \omega_1) = (\psi_1, \omega_2) = 1.$$

Let $P_D P_C^{-1} = \|s_{ij}\|_1^2$. Consider

$$(\psi_2\gamma_2, \omega_2, [\psi_1\gamma_1, \omega_1]s_{21}) = \left((\psi_2\gamma_2, \omega_2), \frac{\psi_1\gamma_1\omega_1}{(\psi_1\gamma_1, \omega_1)}s_{21}\right) =$$

$$= \left(\gamma_2, \omega_2, \frac{\psi_1\gamma_1\omega_1}{(\gamma_1, \omega_1)}s_{21}\right) = \left(\gamma_2, \omega_2, \frac{\gamma_1\omega_1}{(\gamma_1, \omega_1)}\psi_1 s_{21}\right) =$$

$$= (\gamma_2, (\omega_2, [\gamma_1, \omega_1]\psi_1)s_{21}) = (\gamma_2, \omega_2, [\gamma_1, \omega_1]s_{21}).$$

Since $(D, C)_l = I$, we conclude that

$$(\gamma_2, \omega_2, [\gamma_1, \omega_1]s_{21}) = 1.$$

Thus,

$$(\psi_2\gamma_2, \omega_2, [\psi_1\gamma_1, \omega_1]s_{21}) = 1,$$

which implies that $(D(A_1 + B_1L), C)_l = I$, i.e., $(A + BL, C)_l = I$.

3). Suppose that $\det A \neq 0$, $\det C = 0$. Hence, $C = P_C^{-1}\Omega Q_C^{-1}$, $\Omega = \mathrm{diag}(\omega_1, 0)$. First we consider the case where $(A, B)_l = D \neq I$, i.e., $D = P_D^{-1}\Gamma Q_D^{-1}$, $\Gamma = \mathrm{diag}(1, \gamma_2)$. Since $(D, C)_l = I$, we obtain

$$(\gamma_2, \omega_1 s) = 1, \tag{2.14}$$

**91**

where $P_D P_C^{-1} = \begin{Vmatrix} 1 & 0 \\ s & 1 \end{Vmatrix}$. Then $s \neq 0$. By virtue of Theorem 2.19, there exists a matrix $K$ such that $A_1 + B_1 K \sim \mathrm{diag}(1, \mu)$, where $(\mu, \omega_1 s \det D) = 1$. Since $(\mu, \det D) = 1$, by analogy with the case 2), we can show that

$$A + BK = D(A_1 + B_1 K) \sim P_D^{-1} \mathrm{diag}(1, \gamma_2 \mu).$$

Since we also have $(\mu, \omega_1 s) = 1$, it follows from (2.14) that $(\mu \gamma_2, \omega_1 s) = 1$. This means that $(A + BK, C)_l = I$.

Let $(A, B)_l = I$. We choose the matrices $P_A$ and $P_B$ such that

$$P_B P_A^{-1} = \begin{Vmatrix} 1 & 0 \\ c & 1 \end{Vmatrix}.$$

First we consider the case where c $c \neq 0$. Denote $P_B P_C^{-1} = \|p_{ij}\|_1^2$. It follows from the equality $(\varepsilon_1, \delta_1) = 1$ that $(\varepsilon_1, \delta_1, \varepsilon_1 \delta_1 c) = 1$. We choose an element $t_{11}$ such that

$$(\varepsilon_1 + \delta_1 t_{11}, \varepsilon_1 \delta_1 c) = 1, \tag{2.15}$$

moreover

$$\det \begin{Vmatrix} \varepsilon_1 + \delta_1 t_{11} & p_{11} \\ \varepsilon_1 c & p_{21} \end{Vmatrix} \neq 0.$$

This is possible because the equation

$$\det \begin{Vmatrix} \varepsilon_1 + \delta_1 x & p_{11} \\ \varepsilon_1 c & p_{21} \end{Vmatrix} = 0$$

has at most one solution in the variable $x$ and the number of elements satisfying equality (2.15) is greater. Hence, each element of the neighboring class $t_{11} + (\varepsilon_1 \delta_1 c) R$ also satisfies equality (2.15).

Since $(A, B)_l = I$, we have $(\varepsilon_2, \delta_2, \varepsilon_1 \delta_1 c) = 1$

$$(\varepsilon_2 + \delta_2 t_{22}, \varepsilon_1 \delta_1 c) = 1.$$

If $\delta_2 = 0$, then we can set $t_{22} = 0$, because $(\varepsilon_2, \varepsilon_1 \delta_1 c) = 1$. Hence, $(ab, \varepsilon_1 \delta_1 c) = 1$, where $a = \varepsilon_1 + \delta_1 t_{11}, b = \varepsilon_2 + \delta_2 t_{22}$. Therefore, we obtain

$$\left( ab, \varepsilon_1 \delta_1 c, \det \begin{Vmatrix} a & p_{11} \\ \varepsilon_1 c & p_{21} \end{Vmatrix} \omega \right) = 1.$$

Thus, there exists $t_{22}$ such that

$$\left( ab - \varepsilon_1 \delta_1 c t_{12}, \det \begin{Vmatrix} a & p_{11} \\ \varepsilon_1 c & p_{21} \end{Vmatrix} \omega \right) = 1.$$

This implies that

$$\begin{Vmatrix} \varepsilon_1 + \delta_1 t_{11} & \delta_1 t_{12} & \omega p_{11} & 0 \\ \varepsilon_1 c & \varepsilon_2 + \delta_2 t_{22} & \omega p_{21} & 0 \end{Vmatrix} \sim \begin{Vmatrix} I & 0 \end{Vmatrix}.$$

Hence,

$$\left(P_B P_A^{-1}\mathrm{E} + \Delta \left\|\begin{matrix} t_{11} & t_{12} \\ 0 & t_{22} \end{matrix}\right\|, P_B P_C^{-1}\Omega\right)_l = I.$$

and, therefore

$$\left(P_A^{-1}\mathrm{E}Q_A^{-1} + P_B^{-1}\Delta Q_B^{-1}\left(Q_B \left\|\begin{matrix} t_{11} & t_{12} \\ 0 & t_{22} \end{matrix}\right\|\right), P_C^{-1}\Omega Q_C^{-1}\right)_l = (A + BU, C)_l = I,$$

where $U = Q_B \left\|\begin{matrix} t_{11} & t_{12} \\ 0 & t_{22} \end{matrix}\right\|$.

In conclusion, it remains to consider the case where $(A, B)_l = I$, and, in addition, $P_B P_A^{-1} = I$, i.e., $P_A = P_B$. By virtue of Corollary 2.10, the matrices $A, B$ are nonsingular and $(\det A, \det B) = 1$. Hence, $(\varepsilon_1\varepsilon_2, \delta_1\delta_2) = 1$. This implies that $(\varepsilon_1, \delta_1, \omega) = 1$, i $(\varepsilon_2, \delta_2, \omega) = 1$. We choose $q_{11}, q_{22}$ such that

$$(\varepsilon_1 + \delta_1 q_{11}, \omega) = 1, (\varepsilon_2 + \delta_2 q_{22}, \omega) = 1.$$

This yields
$$\left\|\begin{matrix} \varepsilon_1 + \delta_1 t_{11} & 0 & \omega p_{11} & 0 \\ 0 & \varepsilon_2 + \delta_2 q_{22} & \omega p_{21} & 0 \end{matrix}\right\| \sim \left\| I \quad 0 \right\|,$$

where $P_B P_C^{-1} = \|p_{ij}\|_1^2$. Therefore, $(A + BV, C)_l = I$, where $V = Q_B \mathrm{diag}(q_{11}, q_{22})$. Sufficiency is proved.

**Necessity.** Let $R$ be a Bezout ring and let $(a, b, c) = 1, c \neq 0$. Consider the matrices $A = \mathrm{diag}(1, a), B = \mathrm{diag}(0, b), C = \mathrm{diag}(0, c)$. It is clear that $(A, B, C)_l = I$. Then there exists a matrix $T = \|t_{ij}\|_1^2$ such that $(A + BT, C)_l = I$. Since

$$A + BT = \left\|\begin{matrix} 1 & 0 \\ bt_{21} & a + bt_{22} \end{matrix}\right\|,$$

we obtain

$$\left\|\begin{matrix} 1 & 0 & 0 & 0 \\ bt_{21} & a + bt_{22} & 0 & c \end{matrix}\right\| \sim \left\| I \quad 0 \right\|.$$

This implies that $(a + bt_{22}, c) = 1$. $\qquad\square$

The stable range of the ring $R$ (we write st.r.$(R)$ ) is closely connected with the stable range of the ring $n \times n$ matrices over this ring. This relationship was established by Vasershtein in [77]:

$$\mathrm{st.r.}(M_n(R)) = 1 - \left[-\frac{\mathrm{st.r.}(R) - 1}{n}\right].$$

Here, the symbol $[*]$ denotes the integral part of the number. According to this formula, if st.r.$(R)$ is equal to 1 or 2, then the ring $M_n(R)$ has a similar stable range. Theorem 2.20 states that the ring of matrices of the second order over a commutative Bezout domain $R$ of stable range 1.5 inherits this property.

At the end of this subsection, we point out another important property of the Bezout rings of stable range 1.5.

**Theorem 2.21.** *Let $R$ be a commutative Bezout domain. The following are equivalent:*

*1) $R$ has stable range 1.5;*

*2) for each $n \times m$ matrix $A$ over $R$, $\operatorname{rank} A > 1$, there is a row $u = \|1 \quad u_2 \quad ... \quad u_n\|$ such that*

$$uA = \|b_1 \quad b_2 \quad ... \quad b_m\|,$$

*where $(b_1, b_2, ..., b_m)$ coincide with the g.c.d. of elements $A$.*

**Proof.** 1) $\Rightarrow$ 2). Without loss of generality, assume that $m \geq n$. By Theorem 2.1, any commutative Bezout domain of stable range 1.5 is an elementary divisor ring. Then there exist $P \in \operatorname{GL}_n(R)$, $Q \in \operatorname{GL}_m(R)$ such that

$$PAQ = \operatorname{diag}(\varphi_1, ..., \varphi_k, 0, ..., 0) = \Phi,$$

$\varphi_i | \varphi_{i+1}$, $i = 1, ..., k-1$. Consider the matrix

$$U = \left\| \begin{matrix} 0 & \mathbf{0} & 1 & \mathbf{0} \\ \mathbf{0} & I_{k-2} & \mathbf{0} & \mathbf{0} \\ 1 & \mathbf{0} & 0 & \mathbf{0} \\ \mathbf{0} & \mathbf{0} & \mathbf{0} & I_{n-k} \end{matrix} \right\|,$$

where $I_j$ is $j \times j$ identity matrix. Then

$$(UP)A(QU) = \operatorname{diag}(\varphi_k, \varphi_2, \varphi_3, ..., \varphi_{k-1}, \varphi_1, 0, ..., 0) = \Phi'.$$

Denote by $P_1 = UP \det(UP)^{-1}$, $Q_1 = QU \det UP$. Then $P_1 A Q_1 = \Phi'$, moreover, $\det P_1 = 1$. Consider the matrix $P_1^{-1} = \|p_{ij}\|_1^n$. We have

$$\det P_1^{-1} = p_{11}\Delta_1 + ... + (-1)^{n+1}p_{1n}\Delta_n = 1,$$

where $\Delta_i$ are the corresponding minors. According to Theorem 1.9 there are $s_{11}, s_{12}, ..., s_{1n}$ such that the equalities

$$s_{11}\Delta_1 + s_{12}\Delta_2 + ... + s_{1n}\Delta_n = 1,$$
$$(s_{11}, s_{12}, ..., s_{1k}) = 1,$$
$$(s_{1k}, \varphi_k) = 1$$

hold. By Theorem 1.8, the elements $s_{11}, s_{12}, ..., s_{1n}$ satisfy the condition

$$\|1 \quad t_2 \quad ... \quad t_n\| P_1^{-1} = \|s_{11} \quad s_{12} \quad ... \quad s_{1n}\|,$$

where $t_2, ..., t_n \in R$. The equality $P_1 A Q_1 = \Phi'$ implies that $A Q_1 = P_1^{-1} \Phi'$. Then

$$\left\| \begin{matrix} 1 & t_2 & ... & t_n \\ 0 & 1 & ... & 0 \\ \vdots & & \ddots & \\ 0 & 0 & ... & 1 \end{matrix} \right\| A Q_1 = \left\| \begin{matrix} s_{11} & s_{12} & ... & s_{1n} \\ p_{21} & p_{22} & ... & p_{2n} \\ ... & ... & ... & ... \\ p_{n1} & p_{n2} & ... & p_{nn} \end{matrix} \right\| \Phi' =$$

$$= \begin{Vmatrix} s_{11}\varphi_k & s_{12}\varphi_2 & ... & s_{1.k-1}\varphi_{k-1} & s_{1k}\varphi_1 & \mathbf{0} \\ p_{21}\varphi_k & p_{22}\varphi_2 & ... & p_{2.k-1}\varphi_{k-1} & p_{2k}\varphi_1 & \mathbf{0} \\ ... & ... & ... & ... & ... & ... \\ p_{n1}\varphi_k & p_{n2}\varphi_2 & ... & p_{n.k-1}\varphi_{k-1} & p_{nk}\varphi_1 & \mathbf{0} \end{Vmatrix}. \qquad (2.16)$$

The g.c.d. of the first row elements of this matrix is equal to

$$(s_{11}\varphi_k, s_{12}\varphi_2, ..., s_{1.k-1}\varphi_{k-1}, s_{1k}\varphi_1) =$$
$$= \varphi_1 \left( s_{11}\frac{\varphi_k}{\varphi_1}, s_{12}\frac{\varphi_2}{\varphi_1}, ..., \frac{\varphi_{k-1}}{\varphi_1}s_{1.k-1}, s_{1k} \right) =$$
$$= \varphi_1 \left( \left(s_{1k}, s_{11}\frac{\varphi_k}{\varphi_1}\right), \left(s_{1k}, s_{12}\frac{\varphi_2}{\varphi_1}\right), ..., \left(s_{1k}, s_{1.k-1}\frac{\varphi_{k-1}}{\varphi_1}\right) \right) =$$
$$= \varphi_1(s_{11}, s_{12}, ..., s_{1k}) = \varphi_1.$$

Taking into account (2.16) and arguing as above we have

$$\|1 \quad t_2 \quad ... \quad t_n\|AQ_1 = \|s_{11} \quad s_{12} \quad ... \quad s_{1n}\|\Phi' \sim \|\varphi_1 \quad 0 \quad ... \quad 0\|.$$

Hence

$$\|1 \quad t_2 \quad ... \quad t_n\|A \sim \|\varphi_1 \quad 0 \quad ... \quad 0\|.$$

2) $\Rightarrow$ 1). Consider the matrix

$$A = \begin{Vmatrix} a & c \\ b & 0 \end{Vmatrix},$$

where $(a, b, c) = 1$, $b, c \neq 0$. There is a row $\|1 \quad r\|$ such that

$$\|1 \quad r\| \begin{Vmatrix} a & c \\ b & 0 \end{Vmatrix} = \|b_1 \quad b_2\|,$$

where

$$(b_1, b_2) = (a, b, c) = 1.$$

Hence, $(a + br, c) = 1$. It follows that st.r.$(R) = 1.5$. $\qquad \square$

Note that a similar result for adequate rings was obtained by O. Helmer [8], specified and summarized by V. Petrychovych [71].

**Theorem 2.22.** *Let $R$ be a commutative Bezout domain of stable range 1.5 and $A = \|a_{ij}\|$ is $n \times m$ matrix over $R$, rank $A > 1$. Then the set $\mathbf{P}_A$ contains an elementary matrix. A similar statement is correct to the set $\mathbf{Q}_A$.*

**Proof.** By Theorem 2.21, there is a row $u = \|1 \quad u_2 \quad ... \quad u_n\|$ such that

$$uA = \|b_1 \quad b_2 \quad ... \quad b_m\|,$$

where $(b_1, b_2, ..., b_m) = \varphi_1$ is g.c.d. of the elements of the matrix $A$. It means that $\varphi_1$ is the first invariant factor of $A$. Then

$$\underbrace{\begin{Vmatrix} 1 & u_2 & ... & u_n \\ 0 & 1 & ... & 0 \\ \vdots & & \ddots & \\ 0 & 0 & ... & 1 \end{Vmatrix}}_{U_1} A = \begin{Vmatrix} b_1 & b_2 & ... & b_m \\ a_{21} & a_{22} & ... & a_{2m} \\ ... & ... & ... & ... \\ a_{n1} & a_{n2} & ... & a_{nm} \end{Vmatrix}.$$

The group $\mathrm{GL}_m(R)$ contain a matrix $V_1$ such that

$$U_1AV_1 = \begin{Vmatrix} \varphi_1 & 0 & ... & 0 \\ a'_{21} & a'_{22} & ... & a'_{2m} \\ ... & ... & ... & ... \\ a'_{n1} & a'_{n2} & ... & a'_{nm} \end{Vmatrix} = A_1.$$

Obviously $\varphi_1$ also is the g.c.d. of elements of the matrix $A_1$. Thus, $\varphi_1 | a'_{i1}$, $i = 2, ..., n$. So

$$\underbrace{\begin{Vmatrix} 1 & 0 & ... & 0 \\ -\dfrac{a'_{21}}{\varphi_1} & 1 & ... & 0 \\ \vdots & & \ddots & \\ -\dfrac{a'_{n1}}{\varphi_1} & 0 & ... & 1 \end{Vmatrix}}_{U_2} A_1 = \begin{Vmatrix} \varphi_1 & 0 & ... & 0 \\ 0 & c_{22} & ... & c_{2m} \\ ... & ... & ... & ... \\ 0 & c_{n2} & ... & c_{nm} \end{Vmatrix} = A_2.$$

Note that $U_1, U_2$ are elementary matrices. Further arguments are obvious. Eventually we find matrices $P, Q$ such that $PAQ = \Phi$ is the Smith form of $A$. Step by step, we multiply from the left the matrix $A$ by an elementary matrix. This means that $P$ is an elementary matrix.

Consider the matrix $A^T$. As it was proved, there is an elementary matrix $P_1$ and $Q_1 \in \mathrm{GL}_n(R)$ such that

$$P_1 A^T Q_1 = \Phi^T = \Phi \;\; \Rightarrow \;\; A^T = P_1^{-1} \Phi Q_1^{-1}.$$

Transposing this matrix, we get $A = (Q_1^{-1})^T \Phi (P_1^{-1})^T$. Note that $(P_1^{-1})^T$ is an elementary matrix and $(P_1^{-1})^T \in \mathbf{Q}_A$. $\qquad \square$

# LINEAR ALGEBRA
# TECHNIQUE IN MATRIX RINGS

*This section is devoted to the research of matrix minors properties and establishing Smith forms for matrix of special form. The results obtained will be applied in further sections.*

## 3.1. Properties of matrix minors over Bezout rings

Denote by $\langle A \rangle_i$ g.c.d. of the $i$th order minors of the matrix $A$, and $\langle A \rangle$ denote g.c.d. of maximal order minors of this matrix.

**Proposition 3.1.** *Let $A = \|a_{ij}\|$ be $n \times n$ matrix and $B$ be its $t \times k$ submatrix. If $t + k > n$, then $\langle B \rangle_1 \,|\, \det A$.*

**Proof.** Without loss of generality, assume that the matrix $A$ has the form

$$A = \left\| \begin{matrix} B & C \\ * & * \end{matrix} \right\|.$$

Otherwise, by a permutation of the rows and columns of the matrix $A$ we will achieve the desired result. According to Theorem 1.4, such transformations do not change g.c.d. of the minors of corresponding order.

First we treat the case $t + k = n + 1$. The matrix $C$ has the size

$$t \times (n - k) = t \times (t - 1).$$

It means, all $t \times t$ submatrices of the matrix $\|B \quad C\|$ contain at least one column of the matrix $B$. Therefore, every $t$ order minors of the matrix $\|B \quad C\|$ is divided into $\langle B \rangle_1$. By virtue of Laplace Theorem, we get that $\langle B \rangle_1 \,|\, \det A$.

Obviously, this statement will be correct if $t+k > n+1$. $\square$

**Corollary 3.1.** *Let $n \times n$ matrix $A$ contains the zero $t \times k$ submatrix, where $t + k > n$. Then $\det A = 0$.* $\square$

**Corollary 3.2.** *Let $n \times n$ matrix $A$ has the form*

$$A = \left\| \begin{matrix} A_{11} & A_{12} & A_{13} \\ \mathbf{0}_{(n-p)\times p} & A_{22} & A_{23} \end{matrix} \right\|,$$

*where $\mathbf{0}_{(n-p)\times p}$ be the zero $(n-p) \times p$ matrix. If $\langle A_{23} \rangle = 0$, then $\det A = 0$.*

**Proof.** The matrix $\|A_{22} \quad A_{23}\|$ is square of order $n - p$. Since $\langle A_{23}\rangle = 0$ by Laplace Theorem, $\det \|A_{22} \quad A_{23}\| = 0$. It means that all maximum order minors of the matrix $\|\mathbf{0}_{(n-p)\times p} \quad A_{22} \quad A_{23}\|$ are equal to zero. It follows that $\det A = 0$. $\qquad\square$

Let $A = \|a_{ij}\|$ be $m \times n$ matrix, $m > n$. Denote by $\delta_k(A)$ g.c.d. of maximum order minors of this matrix, which contain its $k$ last rows, $k \leqslant n$.

**Lemma 3.1.** *If* $U = \|u_{ij}\| \in \mathrm{GL}_{m-k}(R)$, *then*

$$\delta_k(A) = \delta_k((U \oplus I_k)A).$$

**Proof.** Consider the matrix

$$B = (U \oplus I_k)A = \left\|\begin{array}{ccc} b_{11} & ... & b_{1n} \\ ... & ... & ... \\ b_{m-k.1} & ... & b_{m-k.n} \\ a_{m-k+1.1} & ... & a_{m-k+1.n} \\ ... & ... & ... \\ a_{m1} & ... & a_{mn} \end{array}\right\|.$$

Each $n$th order submatrix of the matrix $B$ containing its $k$ last rows can be written as follows:

$$\underbrace{\left\|\begin{array}{ccc} b_{i_11} & ... & b_{i_1n} \\ ... & ... & ... \\ b_{i_{n-k}1} & ... & b_{i_{n-k}n} \\ a_{m-k+1.1} & ... & a_{m-k+1.n} \\ ... & ... & ... \\ a_{m1} & ... & a_{mn} \end{array}\right\|}_{B_{i_1,...,i_{n-k}}} = \underbrace{\left\|\begin{array}{cccc} u_{i_11} & ... & u_{i_1m-k} & \\ ... & ... & ... & \mathbf{0} \\ u_{i_{n-k}1} & ... & u_{i_{n-k}.m-k} & \\ & & 1 & \\ \mathbf{0} & & & \ddots \\ & & & 1 \end{array}\right\|}_{U_{i_1,...,i_{n-k}}} \times$$

$$\times \left\|\begin{array}{ccc} a_{11} & ... & a_{1n} \\ ... & ... & ... \\ a_{m-k.1} & ... & a_{m-k.n} \\ a_{m-k+1.1} & ... & a_{m-k+1.n} \\ ... & ... & ... \\ a_{m1} & ... & a_{mn} \end{array}\right\|,$$

where $1 \leqslant i_1 < ... < i_{n-k} \leqslant m - k$. A necessary condition that the maximum order ($n$th) minor of the matrix $U_{i_1,...,i_{n-k}}$, is nonzero, there is the presence in it the last $k$ rows of this matrix. Write $\det B_{i_1,...,i_{n-k}}$ as the sum of products of $n$th order minors of matrices $U_{i_1,...,i_{n-k}}$ and $A$. From the remark just made implies that each nonzero addition is a product of the minor of the matrix $A$, containing its last $k$ rows, to the corresponding minor of the matrix $U_{i_1,...,i_{n-k}}$. It means $\delta_k(A)|\det B_{i_1,...,i_{n-k}}$. It follows that

$$\delta_k(A)|\delta_k((U \oplus I_k)A).$$

Reasoning similarly to the above, we get

$$\delta_k(B)|\delta_k((U^{-1} \oplus I_k)B) = \delta_k((U^{-1} \oplus I_k)(U \oplus I_k)A) = \delta_k(A).$$

Since g.c.d. is determined up to associativity, so $\delta_k(A) = \delta_k(B)$. Proof is complete. $\qquad\square$

Let $A = \|a_{ij}\|$ be $m \times n$ matrix. Cross out its $i$th row and $j$th column. Denote by $b_{ij}$ g.c.d. of maximum order minors of the resulting matrix. The matrix $A^* = \|b_{ij}\|$ is called **complementary matrix** of $A$.

**Proposition 3.2.** *Let $A = \|a_{ij}\|$ be $n \times n$ matrix. G.c.d. of elements of each $t \times k$, $t + k > n$, submatrices of the matrix $A^*$ are divisor of $\det A$.*

**Proof.** Let

$$\left\langle \begin{matrix} b_{11} & ... & b_{1k} \\ ... & ... & ... \\ b_{t1} & ... & b_{tk} \end{matrix} \right\rangle_1 = \beta,$$

where $t + k = n + 1$. Consider the $t \times (t-1)$ submatrix

$$A_{tk} = \left\| \begin{matrix} a_{1.k+1} & ... & a_{1n} \\ ... & ... & ... \\ a_{t.k+1} & ... & a_{tn} \end{matrix} \right\|,$$

of the matrix $A$. There is $U \in \mathrm{GL}_t(R)$ such that

$$UA_{tk} = \left\| \begin{matrix} 0 & 0 & ... & 0 \\ c_{21} & 0 & ... & 0 \\ c_{31} & c_{32} & & 0 \\ ... & ... & ... & ... \\ c_{t1} & c_{t2} & ... & c_{t.t-1} \end{matrix} \right\|.$$

Hence,

$$(U \oplus I_{n-t})A = \left\| \begin{matrix} f_{11} & ... & f_{1k} & 0 & 0 & ... & 0 \\ f_{21} & ... & f_{2k} & c_{21} & 0 & ... & 0 \\ f_{31} & ... & f_{3k} & c_{31} & c_{32} & ... & 0 \\ ... & ... & ... & ... & ... & ... & ... \\ f_{t1} & ... & f_{tk} & c_{t1} & c_{t2} & ... & c_{t.t-1} \\ a_{t+1.1} & ... & a_{t+1.k} & a_{t+1.k+1} & a_{t+1.k+2} & ... & a_{t+1.n} \\ ... & ... & ... & ... & ... & ... & ... \\ a_{n1} & ... & a_{nk} & a_{n.k+1} & a_{n.k+2} & ... & a_{nn} \end{matrix} \right\| = C.$$

Note that the minors are $b_{11}, b_{21}, ..., b_{t1}$ can be described as all maximum order minors of the matrix

$$\left\| \begin{matrix} a_{12} & ... & a_{1n} \\ a_{22} & ... & a_{2n} \\ ... & ... & ... \\ a_{n1} & ... & a_{nn} \end{matrix} \right\| = F_1,$$

which contain its last $n - t$ rows. That is

$$\delta_{n-t}(F_1) = (b_{11}, b_{21}, ..., b_{t1}).$$

Denote by $L_1$ the matrix obtained from the matrix $C$ by crossing out its first column. Let $C^* = \|d_{ij}\|_1^n$. Since $L_1 = (U \oplus I_{n-t})F_1$, by Lemma 3.1, we get

$$\delta_{n-t}(F_1) = \delta_{n-t}(L_1).$$

Moreover
$$\delta_{n-t}(L_1) = (d_{11}, d_{21}, ..., d_{t1}).$$

Therefore,
$$(b_{11}, b_{21}, ..., b_{t1}) = (d_{11}, d_{21}, ..., d_{t1}).$$

Reasoning similarly to the above

$$(b_{1i}, b_{2i}, ..., b_{ti}) = (d_{1i}, d_{2i}, ..., d_{ti}),$$

$i = 1, ..., k$. This yields

$$\beta \left| \left\langle \begin{array}{ccc} d_{11} & ... & d_{1k} \\ ... & ... & ... \\ d_{t1} & ... & d_{tk} \end{array} \right\rangle_1 \right. ,$$

Hence, $\beta | (d_{11}, d_{12}, ..., d_{1k})$. Note that

$$\det C = f_{11}d_{11} - f_{12}d_{12} + ... (-1)^{k+1} f_{1k}d_{1k},$$

we have $\beta | \det C$. Since $\det A$ is different from $\det C$ by the invertible factor, we conclude that $\beta | \det A$.

If $t \times k$ submatrix of the matrix $A^*$ is not in the top left corner, then by permuting of rows and columns of the matrix $A$, we get matrix $B$, in which desired submatrix lies in the right place. Then $\beta | \det B$. Since $\det A = \pm \det B$, than $\beta | \det A$. Proof is complete. $\qquad \square$

The elements of a mutual matrix $A_*$ and elements of a complementary matrix $A^*$ differ by sign and position. Therefore, the following statement holds.

**Corollary 3.3.** *G.c.d. of elements of each $t \times k$, $t + k > n$, submatrix $A_*$ is a divisor* $\det A$. $\qquad \square$

**Proposition 3.3.** *Let $A = \|a_{ij}\|$ be $m \times n$ matrix. G.c.d. of elements of each $t \times k$, $t + k > \min(m, n)$, submatrix $A^* = \|b_{ij}\|$ is a divisor* $\langle A \rangle$.

**Proof**. Let $m \geqslant n$. The case $m = n$ discussed in Proposition 3.2. Suppose that $m > n$. Denote by

$$\beta = \left\langle \begin{array}{ccc} b_{11} & ... & b_{1k} \\ ... & ... & ... \\ b_{t1} & ... & b_{tk} \end{array} \right\rangle_1 .$$

Consider $n \times n$ submatrix of the matrix $A$:

$$C = \begin{Vmatrix} a_{11} & ... & a_{1n} \\ ... & ... & ... \\ a_{t1} & ... & a_{tn} \\ * & ... & * \\ ... & ... & ... \\ * & ... & * \end{Vmatrix}.$$

Let $C^* = \|d_{ij}\|_1^n$ be a complementary matrix of the matrix $C$. Since $b_{ij}|d_{ij}$, $i = 1, ..., t$, $j = 1, ..., k$, then

$$\beta \begin{Vmatrix} d_{11} & ... & d_{1k} \\ ... & ... & ... \\ d_{t1} & ... & d_{tk} \end{Vmatrix}.$$

If

$$\left\langle \begin{matrix} d_{11} & ... & d_{1k} \\ ... & ... & ... \\ d_{t1} & ... & d_{tk} \end{matrix} \right\rangle_1 = \gamma,$$

then $\beta|\gamma$. By Proposition 3.2, $\gamma|\det C$. Therefore, $\beta|\det C$.

Now consider the submatrix of the matrix $A$:

$$D_j = \begin{Vmatrix} a_{i_1 1} & ... & a_{i_1 n} \\ ... & ... & ... \\ a_{i_n 1} & ... & a_{i_n n} \end{Vmatrix},$$

$1 \leqslant i_1 < ... < i_n \leqslant m$, which does not contain its $j$th row, $1 \leqslant j \leqslant t$. By definition, $b_{j1}$ is equal to g.c.d. of $(n-1)$ order minor of the matrix

$$A_{j1} = \begin{Vmatrix} a_{12} & ... & a_{1n} \\ ... & ... & ... \\ a_{j-1.2} & ... & a_{j-1.n} \\ a_{j+1.2} & ... & a_{j+1.n} \\ ... & ... & ... \\ a_{m2} & ... & a_{mn} \end{Vmatrix}.$$

The matrix

$$\begin{Vmatrix} a_{i_1 2} & ... & a_{i_1 n} \\ ... & ... & ... \\ a_{i_n 2} & ... & a_{i_n n} \end{Vmatrix}$$

is the submatrix of the matrix $A_{j1}$. We claim that

$$b_{j1} \left| \left\langle \begin{matrix} a_{i_1 2} & ... & a_{i_1 n} \\ ... & ... & ... \\ a_{i_n 2} & ... & a_{i_n n} \end{matrix} \right\rangle \right. .$$

**101**

It follows that $b_{j1}|\det D_j$. Since $\beta|b_{j1}$, then $\beta|\det D_j$. It means that $\beta|\langle A\rangle$. Noting that the location of separated $t \times k$ submatrices of the matrix $A^*$ is irrelevant, we finish the proof. □

Let $U$ be a primitive $m \times n$, $m \leqslant n$, matrix. Denote by $U_{i_1\ldots i_k}$ a matrix consisting of $i_1, ..., i_k$ columns of the matrix $U$, $1 \leqslant k \leqslant m$, $1 \leqslant i_1 < ... < i_k \leqslant n$.

**Proposition 3.4.** *G.c.d. of maximal order minors of the matrix $U$, containing the matrix $U_{i_1\ldots i_k}$ is equal to $\langle U_{i_1\ldots i_k}\rangle$.*

**Proof.** If $k = m$, the statement is obvious. Let $1 \leqslant k < m$. There exists an invertible matrix $L$ such that

$$UL = \left\| S \quad U_{i_1\ldots i_k} \right\|.$$

There exists an invertible matrix $K$ such that

$$KU_{i_1\ldots i_k} = \left\|\begin{matrix} u_1 & * & * \\ 0 & \ddots & * \\ 0 & 0 & u_k \\ & \mathbf{0} & \end{matrix}\right\| = \left\|\begin{matrix} U'_{i_1\ldots i_k} \\ \mathbf{0} \end{matrix}\right\|,$$

where $U'_{i_1\ldots i_k}$ be $k \times k$ matrix. Then

$$KUL = K\left\| S \quad U_{i_1\ldots i_k}\right\| = \left\|\begin{matrix} * & U'_{i_1\ldots i_k} \\ T & \mathbf{0} \end{matrix}\right\|.$$

Each maximal order minor of the matrix $KUL$ containing the matrix

$$\left\|\begin{matrix} U'_{i_1\ldots i_k} \\ \mathbf{0} \end{matrix}\right\|,$$

has the form

$$\det\left\|\begin{matrix} * & U'_{i_1\ldots i_k} \\ T_{m-k} & \mathbf{0} \end{matrix}\right\|,$$

where $T_{m-k}$ be $(m - k) \times (n - k)$ — submatrix of the matrix $T$. Thus, g.c.d. of such minors is equal to $u_1... u_k\tau$, where $\tau = \langle T\rangle$. The matrix $KUL$ is primitive. It follows that the matrix $T$ is also primitive. So, $\tau = 1$. Taking into account that any maximal order minor of the matrix $U$ is differ from the corresponding minor of the matrix $KUL$ at most per unit of the ring $R$, namely the determinant of an invertible matrix, we complete the proof. □

Denote by

$$A\begin{pmatrix} i_1 & \ldots & i_p \\ k_1 & \ldots & k_{n-p} \end{pmatrix}$$

the minor

$$\det \begin{Vmatrix} a_{i_1 k_1} & ... & a_{i_1 k_p} \\ ... & ... & ... \\ a_{i_p k_1} & ... & a_{i_p k_p} \end{Vmatrix}$$

of the matrix $A = \|a_{ij}\|$.

**Corollary 3.4.** *Let $U = \|U_{ij}\|_1^2$ be an invertible block matrix, where $U_{ij}$ be $n \times n$ matrices. Let $U^{-1} = \|V_{ij}\|_1^2$ and $\delta_k(U_{ij})$, $\delta_k(V_{ij})$ are g.c.d. of kth order minors of matrices $U_{ij}$, $V_{ij}$, respectively. Then*

$$\delta_k(U_{11}) = \delta_k(V_{22}), \;\; \delta_k(U_{21}) = \delta_k(V_{21}), \;\; \delta_k(U_{12}) = \delta_k(V_{12}), \;\; \delta_k(U_{22}) = \delta_k(V_{11}),$$

$k = 1, ..., n$.

**Proof.** Minors of matrices $U$, $U^{-1}$ are related as follows:

$$U^{-1} \begin{pmatrix} i_1 & ... & i_p \\ k_1 & ... & k_p \end{pmatrix} = \pm (\det U)^{-1} U \begin{pmatrix} k'_1 & ... & k'_{n-p} \\ i'_1 & ... & i'_{n-p} \end{pmatrix}, \tag{3.1}$$

where $i_1 < i_2 < ... < i_p$ together with $i'_1 < i'_2 < ... < i'_{n-p}$ and $k_1 < k_2 < ... < k_p$ together with $k'_1 < k'_2 < ... < k'_{n-p}$ form a complete index system $1, 2, ..., 2n$ (see p. 30 of [79]). Using Property 3.4, and the formula (3.1), we complete the proof. $\square$

Taking into account Theorem 2.2, this assertion we can reformulate also like this.

**Corollary 3.4\*.** *The Smith forms of matrices $U_{11}$ and $V_{22}$, $U_{12}$ and $V_{12}$, $U_{21}$ and $V_{22}$, $U_{22}$ and $V_{11}$ coincide, that is, the corresponding matrices are equivalent.* $\square$

A similar result was obtained in [73].

**Proposition 3.5.** *Let $U = \|u_{ij}\|_1^n$, $U^{-1} = V = \|v_{ij}\|_1^n$, and*

$$U' = \begin{Vmatrix} u_{i_1 k_1} & u_{i_1 k_2} & ... & u_{i_1 k_p} \\ u_{i_2 k_1} & u_{i_2 k_2} & ... & u_{i_2 k_p} \\ ... & ... & ... & ... \\ u_{i_s k_1} & u_{i_s k_2} & ... & u_{i_s k_p} \end{Vmatrix}$$

*and*

$$V' = \begin{Vmatrix} v_{k'_1 i'_1} & v_{k'_1 i'_2} & ... & v_{k'_1 i'_{n-s}} \\ v_{k'_2 i'_1} & v_{k'_2 i'_2} & ... & v_{k'_2 i'_{n-s}} \\ ... & ... & ... & ... \\ v_{k'_{n-p} i'_1} & v_{k'_{n-p} i'_2} & ... & v_{k'_{n-p} i'_{n-s}} \end{Vmatrix}.$$

*If $i_1 < i_2 < ... < i_s$ together with $i'_1 < i'_2 < ... < i'_{n-s}$ i $k_1 < k_2 < ... < k_p$ together with $k'_1 < k'_2 < ... < k'_{n-p}$ form a complete index system $1, 2, ..., n$, then $\langle U' \rangle = \langle V' \rangle$.*

**Proof.** Suppose that $s \geq p$. In view of (3.1), we have that minors

$$U \begin{pmatrix} i_1 & ... & i_p \\ k_1 & ... & k_p \end{pmatrix}, \quad V \begin{pmatrix} k'_1 & ... & k'_{n-p} \\ i'_1 & ... & i'_{n-p} \end{pmatrix}$$

are associates. Moreover, all minors of $;V \begin{pmatrix} k'_1 & ... & k'_{n-p} \\ i'_1 & ... & i'_{n-p} \end{pmatrix}$ can be described as all minors, containing a submatrix $V'$. Using Proposition 3.4, we conclude our consideration of this case. The case of $s < p$ is similar. $\qquad\square$

**Proposition 3.6.** *Let $a_1, ..., a_n$ be elements of a Bezout ring. There exists $(n-1) \times n$ matrix in which $(n-1)th$ order minors are equal to $a_1, ..., a_n$.*

**Proof.** First we treat the case m $n = 3$. If $(a_1, a_2) \neq 0$, then a desirable matrix will be

$$\left\| \begin{matrix} (a_1, a_2) & ua_3 & -va_3 \\ 0 & \frac{a_1}{(a_1,a_2)} & \frac{a_2}{(a_1,a_2)} \end{matrix} \right\|,$$

where $u, v$ satisfy the equality

$$u \frac{a_2}{(a_1, a_2)} + v \frac{a_1}{(a_1, a_2)} = 1.$$

If $a_1 = a_2 = 0$, then the matrix

$$\left\| \begin{matrix} 0 & 1 & 0 \\ 0 & 0 & a_3 \end{matrix} \right\|$$

will be required.

Let $(a_1, ..., a_{n-1}) = \delta$. Suppose that there exists $(n-2) \times (n-1)$ matrix $S$, having

$$\frac{a_1}{\delta}, ..., \frac{a_{n-1}}{\delta}$$

as $(n-2)$ order minors. Since

$$\left( \frac{a_1}{\delta}, ..., \frac{a_{n-1}}{\delta} \right) = 1,$$

there exists $v_1, ..., v_{n-1}$ such that

$$\det \left\| \begin{matrix} v_1 & v_2 & ... & v_{n-1} \\ & & S & \end{matrix} \right\| = 1.$$

Then the desired matrix will be

$$\left\| \begin{matrix} \delta & v_1 a_n & v_2 a_n & ... & v_{n-1} a_n \\ \mathbf{0} & & & S & \end{matrix} \right\|.$$

The proof is complete. $\qquad\square$

**Lemma 3.2.** *Let*

$$A = \begin{Vmatrix} a_{11} & a_{12} & ... & a_{1,n-1} & a_{1n} \\ a_{21} & a_{22} & ... & a_{2,n-1} & 0 \\ ... & ... & ... & ... & ... \\ a_{n1} & 0 & ... & 0 & 0 \\ \alpha_1 & 0 & ... & 0 & 0 \\ 0 & \alpha_2 & & & 0 \\ & & \ddots & & \\ & & & \alpha_{n-1} & 0 \\ 0 & 0 & ... & 0 & \alpha_n \end{Vmatrix}$$

*be a primitive matrix, moreover* $\alpha_i | \alpha_{i+1}$, $i = 1, ..., n-1$. *Then*

$$(a_{n-i+1,i}, \alpha_i) = 1, i = 1, ..., n.$$

**Proof.** The primitive matrix contains a primitive columns. It follows that $(a_{1n}, \alpha_n) = 1$. There are $u$ and $v$ such that $ua_{1n} + v\alpha_n = 1$. Then

$$\begin{Vmatrix} u & \mathbf{0} & v \\ \mathbf{0} & I_{2n-2} & \mathbf{0} \\ -\alpha_n & \mathbf{0} & a_{1n} \end{Vmatrix} A =$$

$$= \begin{Vmatrix} ua_{11} & ua_{12} & ... & ua_{1,n-1} & 1 \\ a_{21} & a_{22} & ... & a_{2,n-1} & 0 \\ ... & ... & ... & ... & ... \\ a_{n1} & 0 & ... & 0 & 0 \\ \alpha_1 & 0 & ... & 0 & 0 \\ 0 & \alpha_2 & & & 0 \\ \vdots & & \ddots & & \vdots \\ 0 & 0 & & \alpha_{n-1} & 0 \\ -\alpha_n a_{11} & -\alpha_n a_{12} & ... & -\alpha_n a_{1,n-1} & 0 \end{Vmatrix} \sim$$

$$\sim \begin{Vmatrix} 0 & 0 & ... & 0 & 1 \\ a_{21} & a_{22} & ... & a_{2,n-1} & 0 \\ ... & ... & ... & ... & ... \\ a_{n1} & 0 & ... & 0 & 0 \\ \alpha_1 & 0 & ... & 0 & 0 \\ 0 & \alpha_2 & & 0 & 0 \\ \vdots & & \ddots & & \vdots \\ 0 & 0 & & \alpha_{n-1} & 0 \\ 0 & 0 & ... & 0 & 0 \end{Vmatrix} = A_1.$$

The matrix $A_1$ is primitive. It follows that

$$(a_{2,n-1}, \alpha_{n-1}) = 1.$$

Further proving is as it was described. $\qquad \square$

**Proposition 3.7.** *Let $R$ be an elementary divisor ring and $A$ be $m \times s$ matrix over $R$ having*

$$\operatorname{diag}(\underbrace{1, ..., 1}_{k}, \alpha_{k+1}, ..., \alpha_s)$$

*as the Smith form, where $\alpha_{k+1} \notin U(R)$, $m \geq s$. The matrix $A$ is a submatrix of some an invertible $n \times n$ matrix if and only if $n \geq m + s - k$.*

**Proof. Necessity.** Suppose that $A$ is a submatrix of an invertible $n \times n$ matrix $U$. Without loss of generality, we can assume that the matrix $U$ has the form

$$U = \begin{Vmatrix} A & B \\ C & D \end{Vmatrix},$$

where $B, C, D$ are matrices of corresponding sizes. Let $P$ and $Q$ are transforming matrices of the matrix $A$, i.e.,

$$PAQ = \begin{Vmatrix} I_k & \mathbf{0} \\ \mathbf{0} & S \\ & \mathbf{0} \end{Vmatrix} = \Phi,$$

where $S = \operatorname{diag}(\alpha_{k+1}, ..., \alpha_s)$. Then

$$\begin{Vmatrix} P & \mathbf{0} \\ \mathbf{0} & I_{n-m} \end{Vmatrix} U \begin{Vmatrix} Q & \mathbf{0} \\ \mathbf{0} & I_{n-s} \end{Vmatrix} = \begin{Vmatrix} \Phi & PB \\ CQ & D \end{Vmatrix}.$$

Let $F$ be $n \times (s-k)$ matrix consisting of $k+1, k+2, ..., s$ columns of the matrix

$$\begin{Vmatrix} \Phi & PB \\ CQ & D \end{Vmatrix}.$$

Then

$$F = \begin{Vmatrix} \mathbf{0} \\ S \\ \mathbf{0} \\ C_k \end{Vmatrix},$$

where $C_k$ is a submatrix of the matrix $C$. Since $\alpha_{k+1} | \alpha_j$ , $j = k+1, ..., s$, then $\alpha_{k+1}$ divides all maximal $s - k$th order minors of the matrix $F$, containing at least one of the first $m$ rows of this matrix. The $F$ matrix is primitive. Therefore, among its maximal order minors there is one that does not contain the first $m$ rows of this matrix. It means that $(n - m) \times (s - k)$ matrix $C_k$ contains at least $s - k$ rows. That is,

$$n - m \geq s - k \Rightarrow n \geq m + s - k.$$

**Sufficiency.** Let

$$U = \left\| \begin{matrix} P^{-1} & \mathbf{0} \\ \mathbf{0} & I_{s-k} \end{matrix} \right\| \left\| \begin{matrix} I_k & \mathbf{0} & \mathbf{0} \\ \mathbf{0} & S & I_{m-k} \\ \mathbf{0} & \mathbf{0} & \\ \mathbf{0} & I_{s-k} & \mathbf{0} \end{matrix} \right\| \left\| \begin{matrix} Q^{-1} & \mathbf{0} \\ \mathbf{0} & I_{m-k} \end{matrix} \right\|.$$

The $U$ matrix is an example of a desired invertible matrix of the order $n = m + s - k$. Proof is complete. $\qquad\square$

**Corollary 3.5.** *Let* $\mathrm{diag}(1, ..., 1, \underbrace{\alpha, ..., \omega}_{l})$ *be the Smith form of the matrix*

*A. Then A is complemented to a primitive matrix by l rows.* $\qquad\square$

**Corollary 3.6.** *Let A is complemented to a primitive matrix by l rows. Then the number of 1's invariant factors does not exceed l.* $\qquad\square$

## 3.2. Invariant factors of block-triangular matrices and their diagonal blocks

In 1952, Roth published a paper [80] that stood out among others due to the fact that not only a complete solution of the formulated problem was given, but also the elegance of the result was remarkable. Namely, it was shown that the matrix equation

$$AX + YB = C \tag{3.2}$$

(respectively, $AX + XB = C$) over a field has a solution if and only if the matrices

$$\left\| \begin{matrix} A & 0 \\ C & B \end{matrix} \right\| \quad \text{and} \quad \left\| \begin{matrix} A & 0 \\ 0 & B \end{matrix} \right\|$$

are equivalent (respectively, similar).

Let us dwell only on some of the papers dealing with this subject. In [81], it was shown that the Roth's result concerning the solvability of equation (3.2) (further, the *Roth property*) remains correct in the case of commutative principal ideal domains. In [82], this result was generalized to the case of arbitrary commutative rings. As a tribute to the Roth results, the paper [83] introduced the notion of a ring with Roth property; these are rings for which the Roth property holds. Commutative elementary divisor domains are examples of such rings. Over such rings, the equivalence of matrices is equivalent to the equality (up to associativity) of the corresponding invariant factors of these matrices. Thus, in this case, the problem is reduced to establishing the interdependence between the invariant factors of a block-triangular matrix and its diagonal blocks. Newman [84] was the first to obtain a result in this direction. He showed that, in the case $(\det A, \det B) = 1$, over the commutative principal ideal domains, the set of elementary divisors of the matrix

$$\begin{Vmatrix} A & 0 \\ C & B \end{Vmatrix}$$

is the union of the elementary divisors of the matrices $A$ and $B$. The present subchapter is devoted to the research of the interdependence between the invariant factors of the matrices

$$\begin{Vmatrix} A & 0 \\ C & B \end{Vmatrix}, \; A, B.$$

Let $R$ be a commutative elementary divisor domain

**Theorem 3.1.** *Suppose that $k \times k$ matrices $D_1$ and $D_3$ have the Smith forms*

$$\Delta_1 = \mathrm{diag}(\alpha_1, \alpha_2, ..., \alpha_k), \quad \Delta_3 = \mathrm{diag}(\beta_1, \beta_2, ..., \beta_k),$$

*respectively. Further, let*

$$D = \begin{Vmatrix} D_1 & \mathbf{0} \\ D_2 & D_3 \end{Vmatrix} \sim \mathrm{diag}(\underbrace{1, ..., 1}_{k}, \underbrace{\varphi, ..., \varphi}_{k}).$$

*Then the invariant factors $\alpha_i$ and $\beta_j$ can be chosen so that*

$$\alpha_t \beta_{k-t+1} = \varphi, \quad i = 1, ..., k.$$

**Proof.** For the matrices $D_1$ and $D_3$, there exist invertible matrices $P_1$, $Q_1$ and $P_3$, $Q_3$ such that

$$P_1 D_1 Q_1 = \Delta_1, \quad P_3 D_3 Q_3 = \Delta_3.$$

Then

$$\begin{Vmatrix} P_1 & \mathbf{0} \\ 0 & P_3 \end{Vmatrix} D \begin{Vmatrix} Q_1 & \mathbf{0} \\ 0 & Q_3 \end{Vmatrix} = \begin{Vmatrix} P_1 D_1 Q_1 & \mathbf{0} \\ P_3 D_2 Q_1 & P_3 D_3 Q_3 \end{Vmatrix} = \begin{Vmatrix} \Delta_1 & \mathbf{0} \\ D_2' & \Delta_3 \end{Vmatrix}.$$

Therefore, without loss of generality, we can assume that the matrix $D$ has the form

$$D = \begin{Vmatrix} \alpha_1 & & & & & & & & & \mathbf{0} & & & & \mathbf{0} \\ & \ddots & & & & & & & & & & & & \\ & & \alpha_{s-1} & & & & & & & & & & & \\ & & & \alpha_s & & & & & & & & & & \\ & & & & \ddots & & & & & & & & & \\ \mathbf{0} & & & & & \alpha_k & & & & & & & & \\ a_{11} & ... & a_{1.s-1} & a_{1s} & ... & a_{1k} & \beta_1 & & & & \mathbf{0} & & & \\ ... & ... & ... & ... & ... & ... & & \ddots & & & & & & \\ a_{t-1.1} & ... & a_{t-1.s-1} & a_{t-1.s} & ... & a_{t-1.k} & & & \beta_{t-1} & & & & & \\ a_{t1} & ... & a_{t.s-1} & a_{ts} & ... & a_{tk} & & & & \beta_t & & & & \\ ... & ... & ... & ... & ... & ... & & & & & \ddots & & & \\ a_{k1} & ... & a_{k.s-1} & a_{ks} & ... & a_{kk} & \mathbf{0} & & & & & \beta_k \end{Vmatrix} =$$

$$= \left\| \begin{matrix} \Delta_1 & \mathbf{0} \\ A & \Delta_3 \end{matrix} \right\|.$$

Suppose that $t + s = k + 1$. First, let us show that $\varphi \mid \alpha_s \beta_t$. Let $(\alpha_s \beta_t, \varphi) = \tau_s$. Then $\varphi = \tau_s \sigma_s$ and

$$\left( \frac{\alpha_s \beta_t}{\tau_s}, \ \sigma_s \right) = 1. \tag{3.3}$$

Let us show that $\sigma_s$ is a divisor of all the minors of $k$th order of the matrix $D$. Cross out in this matrix $t$th row and $s$th column. This matrix is denoted by $D_{st}$:

$$D_{st} =$$

$$= \left\| \begin{matrix}
\alpha_1 & & & & & & \mathbf{0} & & & & & & \mathbf{0} \\
& \ddots & & & & & & & & & & & \\
& & \alpha_{s-1} & & & & & & & & & & \\
& & & \alpha_{s+1} & & & & & & & & & \\
& & & & \ddots & & & & & & & & \\
\mathbf{0} & & & & & \alpha_k & & & & & & & \\
a_{11} & \dots & a_{1.s-1} & a_{1.s+1} & \dots & a_{1k} & \beta_1 & & & & & & \mathbf{0} \\
\dots & \dots & \dots & \dots & \dots & \dots & & \ddots & & & & & \\
a_{t-1.1} & \dots & a_{t-1.s-1} & a_{t-1.s+1} & \dots & a_{t-1.k} & & & \beta_{t-1} & & & & \\
a_{t+1.1} & \dots & a_{t+1.s-1} & a_{t+1.s+1} & \dots & a_{t+1.k} & & & & \beta_{t+1} & & & \\
\dots & \dots & \dots & \dots & \dots & \dots & & & & & \ddots & & \\
a_{k1} & \dots & a_{k.s-1} & a_{k.s+1} & \dots & a_{kk} & \mathbf{0} & & & & & \beta_k
\end{matrix} \right\|.$$

Let $d_{st}^{(k-1)}$ is an arbitrary minor of $(k-1)$th order of the matrix $D_{st}$. Since

$$D \sim \mathrm{diag}(\underbrace{1, ..., 1}_{k}, \varphi, ..., \varphi),$$

it follows that each minor of $(k+1)$th order of the matrix $D$ is divisible by $\varphi$. In particular, the same applies to the minor $\alpha_s \beta_t d_{st}^{(k-1)}$. Therefore,

$$\sigma_s \left| \frac{\alpha_s \beta_t}{\tau_s} d_{st}^{(k-1)} \right..$$

It follows from the relation (3.3) that $\sigma_s \mid d_{st}^{(k-1)}$. Thus, $\sigma_s$ is a divisor of all the minors of $(k-1)$th order of the matrix $D_{st}$ and, in particular, we have

$$\sigma_s | \alpha_1 \beta_1 \beta_2 \dots \beta_{t-1} \beta_{t+1} \dots \beta_{k-1} \ \Rightarrow \ \sigma_s | \alpha_1 \beta_1 \beta_2 \dots \beta_{k-1};$$

$$\sigma_s | \alpha_1 \alpha_2 \beta_1 \beta_2 \dots \beta_{t-1} \beta_{t+1} \dots \beta_{k-2} \ \Rightarrow \ \sigma_s | \alpha_1 \alpha_2 \beta_1 \beta_2 \dots \beta_{k-2};$$

$$\dots\dots\dots\dots\dots\dots\dots\dots\dots\dots\dots\dots\dots\dots\dots\dots\dots\dots\dots\dots\dots\dots\dots$$

$$\sigma_s | \alpha_1 \alpha_2 \dots \alpha_{s-1} \beta_1 \beta_2 \dots \beta_{t-1} \ \Rightarrow \ \sigma_s | \alpha_1 \alpha_2 \dots \alpha_{s-1} \beta_1 \beta_2 \dots \beta_t.$$

Also we have

$$\sigma_s|\alpha_1\alpha_2\ldots\alpha_{s-1}\beta_1\beta_2\ldots\beta_{t-1} \Rightarrow \sigma_s|\alpha_1\alpha_2\ldots\alpha_s\beta_1\beta_2\ldots\beta_{t-1};$$
$$\sigma_s|\alpha_1\alpha_2\ldots\alpha_{s-1}\alpha_{s+1}\beta_1\beta_2\ldots\beta_{t-2} \Rightarrow \sigma_s|\alpha_1\alpha_2\ldots\alpha_{s+1}\beta_1\beta_2\ldots\beta_{t-2};$$
$$\ldots\ldots\ldots\ldots\ldots\ldots\ldots\ldots\ldots\ldots\ldots\ldots\ldots\ldots\ldots\ldots\ldots\ldots\ldots\ldots\ldots$$
$$\sigma_s|\alpha_1\alpha_2\ldots\alpha_{s-1}\alpha_{s+1}\ldots\alpha_{k-1}\beta_1 \Rightarrow \sigma_s|\alpha_1\alpha_2\ldots\alpha_{k-1}\beta_1.$$

Each minor of $k$th order which main diagonal contains only the diagonal elements of the matrix $D$ is of the form $\alpha_{i_1}\ldots\alpha_{i_p}\beta_{j_1}\ldots\beta_{j_{1q}}$, where

$$i_1,...,i_p,j_1,...,j_q \in \{1,...,k\},$$

$p + q = k$. Since $\alpha_1\ldots\alpha_p|\alpha_{i_1}\ldots\alpha_{i_p}$ and $\beta_1\ldots\beta_q|\beta_{j_1}\ldots\beta_{j_{1q}}$, then

$$\alpha_1\ldots\alpha_p\beta_1\ldots\beta_q|\alpha_{i_1}\ldots\alpha_{i_p}\beta_{j_1}\ldots\beta_{j_{1q}}.$$

Then the divisibility relations given above imply that

$$\sigma_s \mid \alpha_{i_1}\ldots\alpha_{i_p}\beta_{j_1}\ldots\beta_{j_{1q}}.$$

Suppose that $l$, $m$ are indices such that $1 \leqslant l \leqslant s$, $1 \leqslant m \leqslant t$. Then $\varphi \mid \alpha_l\beta_m d_{lm}^{(k-1)}$, where $d_{lm}^{(k-1)}$ is an arbitrary minor of $(k-1)$th order of the matrix $D_{lm}$. Since $\alpha_l \mid \alpha_s$, and $\beta_m \mid \beta_t$, it follows that $\varphi \mid \alpha_s\beta_t d_{lm}^{(k-1)}$. Therefore, just as in the previous case, $\sigma_s \mid d_{lm}^{(k-1)}$. Thus, $\sigma_s$ is a divisor of all the minors of $(k-1)$th order of the following matrices:

$$D_{11},...,D_{1s},$$
$$\ldots\ldots\ldots\ldots\ldots$$
$$D_{t1},...,D_{ts}.$$

Consider the matrix

$$\left\|\begin{array}{c} A \\ \Delta_1 \end{array}\right\|,$$

and let

$$\left\|\begin{array}{c} A \\ \Delta_1 \end{array}\right\|^* = \left\|\begin{array}{ccccc} b_{11} & ... & b_{1s} & ... & b_{1k} \\ ... & ... & ... & ... & ... \\ b_{t1} & ... & b_{ts} & ... & b_{tk} \\ ... & ... & ... & ... & ... \\ b_{2k.1} & ... & b_{2k.s} & ... & b_{2k.k} \end{array}\right\|$$

be its complementary matrix. The fact that $\sigma_s$ is a divisor of all the minors of $(k-1)$th order of the matrix $D_{ij}$ implies that $\sigma_s \mid b_{ij}$, $1 \leqslant i \leqslant t$, $1 \leqslant j \leqslant s$. Therefore,

$$\sigma_s \left\|\left\|\begin{array}{ccc} b_{11} & ... & b_{1s} \\ ... & ... & ... \\ b_{t1} & ... & b_{ts} \end{array}\right\|\right\|.$$

Since $t + s = k + 1 > k$, it follows that, by Proposition 3.3, $\sigma_s$ is a divisor of all the minors of $k$th order of the matrix $\left\| \begin{matrix} A \\ \Delta_1 \end{matrix} \right\|$ hence, and the matrices $\left\| \begin{matrix} \Delta_1 \\ A \end{matrix} \right\|$.

Using similar arguments, we can show that $\sigma_s$ divides all the minors of $k$th order of the matrix $\| A \quad \Delta_3 \|$.

Consider the minor of $k$th order $\alpha_1 \dots \alpha_\mu \beta_1 \dots \beta_\nu \det K$, $\mu \geqslant \nu$, where $K$ is a submatrix of the corresponding order of the matrix $A$. As proved above,

$$\sigma_s | \alpha_1 \dots \alpha_\mu \begin{vmatrix} a_{1.\mu+1} & \dots & a_{1k} \\ \dots & \dots & \dots \\ a_{\nu.\mu+1} & \dots & a_{\nu k} \\ a_{\nu+1.\mu+1} & \dots & a_{\nu+1.k} \\ \dots & \dots & \dots \\ a_{k.\mu+1} & \dots & a_{kk} \end{vmatrix}.$$

Then, a fortiori,

$$\sigma_s | \alpha_1 \dots \alpha_\mu \begin{vmatrix} a_{\nu+1.\mu+1} & \dots & a_{\nu+1.k} \\ a_{\nu+2.\mu+1} & \dots & a_{\nu+2.k} \\ \dots & \dots & \dots \\ a_{k.\mu+1} & \dots & a_{kk} \end{vmatrix} e \beta_1 \dots \beta_\nu.$$

Therefore, $\sigma_s \mid \alpha_1 \dots \alpha_\mu \beta_1 \dots \beta_\nu \det K$. But if $\nu > \mu$, then the required result is a consequence of the fact that

$$\sigma_s | \beta_1 \dots \beta_\nu \begin{vmatrix} a_{\nu+1.1} & \dots & a_{\nu+1.k} \\ a_{\nu+2.1} & \dots & a_{\nu+2.k} \\ \dots & \dots & \dots \\ a_{k1} & \dots & a_{kk} \end{vmatrix}.$$

Similarly, we can show that $\sigma_s | \alpha_{i_1} \dots \alpha_{i_p} \beta_{j_1} \dots \beta_{j_{1q}}$, where $i_1, \dots, i_p$, $j_1, \dots, j_q \in \{1, \dots, k\}$, $p + q = k$. Thus, we have shown that $\sigma_s$ is a divisor of all the minors of $k$th order of the matrix $D$. In view of the relation

$$D \sim \mathrm{diag}(\underbrace{1, \dots, 1}_{k}, \varphi, \dots, \varphi), \tag{3.4}$$

we have $\sigma_s \in U(R)$. Recalling that $\varphi = \tau_s \sigma_s$, we conclude that $\varphi$ and $\tau_s$ are associates. Since the g.c.d. of the elements is defined up to associativity, we have $(\alpha_s \beta_t, \varphi) = \varphi$. This implies that $\varphi \mid \alpha_s \beta_t$, $t + s = k + 1$. Therefore, $\alpha_s \beta_t = \varphi u_s$. Then

$$\det D = \alpha_1 \alpha_2 \dots \alpha_k \beta_1 \beta_2 \dots \beta_k = (\alpha_1 \beta_k)(\alpha_2 \beta_{k-1}) \dots (\alpha_k \beta_1) = \varphi^k u_1 u_2 \dots u_k.$$

On the other hand, by (3.4) we have $\det D = \varphi^k e$, where $e \in U(R)$. Thus, $u_1 u_2 \dots u_k = e \in U(R)$. Noting that the invariant factors of the matrix are chosen up to divisors of 1, we conclude that we can choose $\alpha_i$ and $\beta_j$ so that

$$\alpha_1 \beta_k = \alpha_2 \beta_{k-1} = \dots = \alpha_k \beta_1 = \varphi.$$

The theorem is proved. $\qquad \square$

In order to study the general case, we establish several auxiliary statements.

**Lemma 3.3.** *If*

$$A = \begin{Vmatrix} a_{11} & ... & a_{1.k-1} & a_{1k} & a_{1.k+1} & ... & a_{1n} \\ ... & ... & ... & ... & ... & ... & ... \\ a_{s1} & ... & a_{s.k-1} & a_{sk} & a_{s.k+1} & ... & a_{sn} \\ 0 & ... & 0 & 1 & 0 & ... & 0 \\ a_{s+2.1} & ... & a_{s+2.k-1} & a_{s+2.k} & a_{s+2.k+1} & ... & a_{s+2.n} \\ ... & ... & ... & ... & ... & ... & ... \\ a_{m1} & ... & a_{m.k-1} & a_{mk} & a_{m.k+1} & ... & a_{mn} \end{Vmatrix} \sim$$

$$\sim \operatorname{diag}(\alpha_1, \alpha_2, ..., \alpha_k),$$

*where $k = \min(m, n)$, $\alpha_i | \alpha_{i+1}$, $i = 1, ..., k-1$, then*

$$\begin{Vmatrix} a_{11} & ... & a_{1.k-1} & a_{1.k+1} & ... & a_{1n} \\ ... & ... & ... & ... & ... & ... \\ a_{s1} & ... & a_{s.k-1} & a_{s.k+1} & ... & a_{sn} \\ a_{s+2.1} & ... & a_{s+2.k-1} & a_{s+2.k+1} & ... & a_{s+2.n} \\ ... & ... & ... & ... & ... & ... \\ a_{m1} & ... & a_{m.k-1} & a_{m.k+1} & ... & a_{mn} \end{Vmatrix} \sim$$

$$\sim \operatorname{diag}(\alpha_2, \alpha_3, ..., \alpha_k).$$

**Proof.** We can easily verify the following equivalences:

$$A \sim \begin{Vmatrix} a_{11} & ... & a_{1.k-1} & 0 & a_{1.k+1} & ... & a_{1n} \\ ... & ... & ... & ... & ... & ... & ... \\ a_{s1} & ... & a_{s.k-1} & 0 & a_{s.k+1} & ... & a_{sn} \\ 0 & ... & 0 & 1 & 0 & ... & 0 \\ a_{s+2.1} & ... & a_{s+2.k-1} & 0 & a_{s+2.k+1} & ... & a_{s+2.n} \\ ... & ... & ... & ... & ... & ... & ... \\ a_{m1} & ... & a_{m.k-1} & 0 & a_{m.k+1} & ... & a_{mn} \end{Vmatrix} \sim$$

$$\sim \begin{Vmatrix} 1 & 0 & ... & 0 & 0 & ... & 0 \\ 0 & a_{11} & ... & a_{1.k-1} & a_{1.k+1} & ... & a_{1n} \\ ... & ... & ... & ... & ... & ... & ... \\ 0 & a_{s1} & ... & a_{s.k-1} & a_{s.k+1} & ... & a_{sn} \\ 0 & a_{s+2.1} & ... & a_{s+2.k-1} & a_{s+2.k+1} & ... & a_{s+2.n} \\ ... & ... & ... & ... & ... & ... & ... \\ 0 & a_{m1} & ... & a_{m.k-1} & a_{m.k+1} & ... & a_{mn} \end{Vmatrix} =$$

$$= \begin{Vmatrix} 1 & \mathbf{0} \\ \mathbf{0} & A_1 \end{Vmatrix} \sim \begin{Vmatrix} 1 & \mathbf{0} \\ \mathbf{0} & \operatorname{diag}(\beta_2, ..., \beta_k) \end{Vmatrix},$$

where

$$A_1 \sim \mathrm{diag}(\beta_2, ..., \beta_k),$$

$\beta_i | \beta_{i+1}, i = 2, ..., k-1$. This implies that $\alpha_i$ and $\beta_i, i = 2, ..., k$, are associates. □

**Lemma 3.4.** *Let*

$$A = \begin{Vmatrix} A_{11} & A_{12} \\ A_{21} & A_{22} \end{Vmatrix} \sim \mathrm{diag}(\alpha_1, \alpha_2, ..., \alpha_n) = \Delta,$$

*where $\alpha_i | \alpha_{i+1}, i = 1, ..., n-1$. Then*

$$B = \begin{Vmatrix} A_{11} & \mathbf{0} & A_{12} \\ \mathbf{0} & 1 & \mathbf{0} \\ A_{21} & \mathbf{0} & A_{22} \end{Vmatrix} \sim \mathrm{diag}(1, \alpha_1, \alpha_2, ..., \alpha_n).$$

**Proof.** There are invertible matrices $P, Q, U, V$ such that

$$PBQ = \begin{Vmatrix} 1 & \mathbf{0} \\ \mathbf{0} & A \end{Vmatrix} = B_1, UAV = \Delta.$$

Then

$$B \sim B_1 \sim \begin{Vmatrix} 1 & \mathbf{0} \\ \mathbf{0} & U \end{Vmatrix} B_1 \begin{Vmatrix} 1 & \mathbf{0} \\ \mathbf{0} & V \end{Vmatrix} =$$

$$= \begin{Vmatrix} 1 & \mathbf{0} \\ \mathbf{0} & UAV \end{Vmatrix} = \begin{Vmatrix} 1 & \mathbf{0} \\ \mathbf{0} & \Delta \end{Vmatrix} = \mathrm{diag}(1, \alpha_1, \alpha_2, ..., \alpha_n),$$

as required. □

**Theorem 3.2.** *Let*

$$D = \begin{Vmatrix} D_1 & \mathbf{0} \\ D_2 & D_3 \end{Vmatrix} \sim \mathrm{diag}(\overbrace{1, ..., 1}^{p+q}, \underbrace{\varphi, ..., \varphi}_{s}),$$

$\varphi \neq 0$, *and $D_1$, $D_3$ be square matrices with the Smith forms*

$$\Delta_1 = \mathrm{diag}(\alpha_1, \alpha_2, ..., \alpha_q), \Delta_3 = \mathrm{diag}(\beta_1, \beta_2, ..., \beta_p),$$

*respectively. Then the following assertions are valid:*

*1) if $p \leq s \leq q$, then the invariant factors $\alpha_i$ and $\beta_j$ can be chosen so that*

$$\alpha_1 = \alpha_2 = ... = \alpha_{q-s} = 1,$$
$$\alpha_{q-s+1}\beta_p = \alpha_{q-s+2}\beta_{p-1} = ... = \alpha_{q-s+p}\beta_1 = \varphi,$$
$$\alpha_{q-s+p+1} = \alpha_{q-s+p+2} = ... = \alpha_q = \varphi;$$

*2) if $p, q \geq s$, then the invariant factors $\alpha_i$ and $\beta_j$ can be chosen so that*

$$\alpha_1 = \alpha_2 = ... = \alpha_{q-s} = 1,$$
$$\beta_1 = \beta_2 = ... = \beta_{p-s} = 1,$$
$$\alpha_{q-s+1}\beta_p = \alpha_{q-s+2}\beta_{p-1} = ... = \alpha_q\beta_{p-s+1} = \varphi;$$

3) *if $p, q \leq s$, then the invariant factors $\alpha_i$ and $\beta_j$ can be chosen so that*

$$\alpha_1\beta_{q+p-s} = \alpha_2\beta_{q+p-s-1} = ... = \alpha_{q+p-s}\beta_1 = \varphi,$$

$$\alpha_{q+p-s+1} = \alpha_{q+p-s+2} = ... = \alpha_q = \beta_{q+p-s+1} = \beta_{q+p-s+2}... = \beta_p = \varphi.$$

**Proof.** Suppose that $p \leq s \leq q$. By Corollary 3.6, we have

$$\alpha_1 = \alpha_2 = ... = \alpha_{q-s} = 1.$$

Crossing out the first $q - s$ rows and columns in the matrix $D$, we obtain the matrix

$$\left\| \begin{matrix} C_1 & \mathbf{0} \\ C_2 & D_3 \end{matrix} \right\|,$$

having, by Lemma 3.3, the Smith form

$$\text{diag}(\underbrace{1, ..., 1}_{p}, \underbrace{\varphi, ..., \varphi}_{s}).$$

Using Lemma 3.4, we obtain

$$\left\| \begin{matrix} C_1 & \mathbf{0} & \mathbf{0} \\ \mathbf{0} & E_{s-p} & \mathbf{0} \\ C_2 & \mathbf{0} & D_3 \end{matrix} \right\| = \left\| \begin{matrix} C_1' & \mathbf{0} \\ C_2' & D_3' \end{matrix} \right\| \sim \text{diag}(\underbrace{1, ..., 1}_{s}, \underbrace{\varphi, ..., \varphi}_{s}),$$

where

$$C_1' \sim \text{diag}(\underbrace{\alpha_{q-s+1}, \alpha_{q-s+2}, ..., \alpha_q}_{s}),$$

$$D_3' \sim \text{diag}(\underbrace{1, ..., 1}_{s-p}, \beta_1, \beta_2, ..., \beta_p).$$

According to Theorem 3.1,

$$\alpha_{q-s}\beta_p = \alpha_{q-s+1}\beta_{p-1} = ... = \alpha_{q-s+p}\beta_1 = \varphi;$$

$$\alpha_{q-s+p+1} = \alpha_{q-s+p+2} = ... = \alpha_q = \varphi.$$

Suppose that $p, q \geq s$. Again, by Corollary 3.6, we have

$$\alpha_1 = \alpha_2 = ... = \alpha_{q-s} = 1,$$

$$\beta_1 = \beta_2 = ... = \beta_{p-s} = 1.$$

Crossing out the rows and the columns containing these 1's and using Lemma 3.3, we obtain

$$\left\| \begin{matrix} F_1 & \mathbf{0} \\ F_2 & F_3 \end{matrix} \right\| \sim \text{diag}(\underbrace{1, ..., 1}_{s}, \underbrace{\varphi, ..., \varphi}_{s}),$$

where

$$F_1 \sim \text{diag}(\underbrace{\alpha_{q-s+1}, \alpha_{q-s+2}, ..., \alpha_q}_{s}),$$

$$F_3 \sim \text{diag}(\underbrace{\beta_{p-s+1}, \beta_{p-s+2}, ..., \beta_p}_{s}).$$

Then, by Theorem 3.1, we have

$$\alpha_{q-s+1}\beta_p = \alpha_{q-s+2}\beta_{p-1} = ... = \alpha_q\beta_{p-s+1} = \varphi.$$

Suppose that $p, q \leqslant s$. By Lemma 3.4, we can write

$$\left\| \begin{array}{cc|cc} D_1 & \mathbf{0} & \mathbf{0} & \mathbf{0} \\ \mathbf{0} & I_{s-q} & \mathbf{0} & \mathbf{0} \\ \hline D_2 & \mathbf{0} & D_3 & \mathbf{0} \\ \mathbf{0} & \mathbf{0} & \mathbf{0} & I_{s-p} \end{array} \right\| = \left\| \begin{array}{cc} H_1 & \mathbf{0} \\ H_2 & H_3 \end{array} \right\| \sim$$

$$\sim \mathrm{diag}(\underbrace{1, ..., 1}_{s}, \underbrace{\varphi, ..., \varphi}_{s}),$$

here

$$H_1 \sim \mathrm{diag}(\underbrace{1, ..., 1}_{s-q}, \alpha_1, \alpha_2, ..., \alpha_q),$$

$$H_3 \sim \mathrm{diag}(\underbrace{1, ..., 1}_{s-p}, \beta_1, \beta_2, ..., \beta_p).$$

Therefore,

$$\alpha_1\beta_{q+p-s} = \alpha_2\beta_{q+p-s-1} = ... = \alpha_{q+p-s}\beta_1 = \varphi,$$
$$\alpha_{q+p-s+1} = \alpha_{q+p-s+2} = ... = \alpha_q = \beta_{q+p-s+1} = \beta_{q+p-s+2}... = \beta_p = \varphi.$$

The theorem is proved. $\qquad\qquad\qquad\qquad\qquad\qquad\qquad\qquad\square$

**Remark.** The case $q \leqslant s \leqslant p$ is not considered because of its symmetry to the case $p \leqslant s \leqslant q$.

If $R$ is a commutative principal ideal domain, then the obtained results can be used to find the invariant factors of the diagonal blocks of the nonsingular matrix

$$D = \left\| \begin{array}{cc} D_1 & \mathbf{0} \\ D_2 & D_3 \end{array} \right\| \sim \Phi = \mathrm{diag}(\varphi_1, ..., \varphi_n)$$

for the case in which

$$\frac{\varphi_2}{\varphi_1}, \frac{\varphi_3}{\varphi_2}, ..., \frac{\varphi_n}{\varphi_{n-1}}$$

are pairwise coprime. To do this, let us apply the locally global method proposed by Gerstein in [85].

Suppose that $p$ is an indecomposable element of the ring $R$. Let $R_{(p)}$ denote the localization of the ring $R$ at the prime ideal $(p)$. In other words, $R_{(p)}$ is the ring consisting of elements of the form $x = p^\nu \frac{a}{b}$, where $a$ and $b$ are elements of the ring $R$ coprime to $p$ and $\nu \in N \bigcup \{0\}$. Since

$$\Phi = I\varphi_1 \, \mathrm{diag}\left(1, \frac{\varphi_2}{\varphi_1}, ..., \frac{\varphi_2}{\varphi_1}\right) \times$$

$$\times \, \mathrm{diag}\left(1, 1, \frac{\varphi_3}{\varphi_2}, ..., \frac{\varphi_3}{\varphi_2}\right) ... \, \mathrm{diag}\left(1, ..., 1, \frac{\varphi_n}{\varphi_{n-1}}\right),$$

it follows that the matrix $\Phi$ can be written as the product of matrices of the form

$$\text{diag}(1, ..., 1, p^\mu, ..., p^\mu).$$

By Theorem 5.6 from [85], the matrix $D$ in the ring $R_{(p)}$ will have the Smith form $\text{diag}(1, ..., 1, p^\mu, ..., p^\mu)$. Therefore, the invariant factors of the matrices $D_1$ and $D_3$ in the ring $R_{(p)}$ are related by the equalities that were written in Theorem 3.2. Then, in order to find the canonical diagonal forms of the matrices $D_1$ and $D_3$ in the ring $R$, we must find the canonical diagonal forms of the matrices $D_1$ and $D_3$ in all the localizations of the ring $R$ at the indecomposable elements that are divisors of the entry $\varphi_n$ and multiply the corresponding invariant factors of these matrices.

## 3.3. Smith form of some matrices

Let $\Phi = \text{diag}(\varphi_1, ..., \varphi_n)$ be a nonsingular $d$-matrix.

**Lemma 3.5.** *Let $S$ be $n \times m$ matrix and*

$$\Phi_i = \text{diag}\left(\frac{\varphi_i}{\varphi_1}, ..., \frac{\varphi_i}{\varphi_{i-1}}, \underbrace{1, ..., 1}_{n-i+1}\right),$$

$i = 2, ..., n.$ *If $H \in \mathbf{G}_\Phi$, then*

$$\Phi_i H S \overset{l}{\sim} \Phi_i S, \ \ i = 2, ..., n.$$

**Proof.** By Theorem 2.7, the group $\mathbf{G}_\Phi$ consists of all invertible matrices of the form

$$\left\| \begin{array}{ccccc} h_{11} & h_{12} & ... & h_{1.n-1} & h_{1n} \\ \dfrac{\varphi_2}{\varphi_1}h_{21} & h_{22} & ... & h_{2.n-1} & h_{2n} \\ ... & ... & ... & ... & ... \\ \dfrac{\varphi_n}{\varphi_1}h_{n1} & \dfrac{\varphi_n}{\varphi_2}h_{n2} & ... & \dfrac{\varphi_n}{\varphi_{n-1}}h_{n.n-1} & h_{nn} \end{array} \right\|.$$

So the $j$th column of the $H$ matrix has the form

$$h_j = \left\| h_{1j} ... h_{jj} \frac{\varphi_{j+1}}{\varphi_j}h_{j+1.j} ... \frac{\varphi_n}{\varphi_j}h_{nj} \right\|^T, \quad j = 1, ..., n-1.$$

Then

$$\Phi_i h_j =$$

$$= \left\| \frac{\varphi_i}{\varphi_1}h_{1j} ... \frac{\varphi_i}{\varphi_{i-1}}h_{i-1.j} \frac{\varphi_i}{\varphi_j}h_{jj} \frac{\varphi_i}{\varphi_j}h_{j+1.j} ... \frac{\varphi_i}{\varphi_j}h_{ij} \frac{\varphi_{i+1}}{\varphi_j}h_{i+1.j} ... \frac{\varphi_n}{\varphi_j}h_{nj} \right\|^T =$$

$$= \frac{\varphi_i}{\varphi_j} \left\| \frac{\varphi_j}{\varphi_1}h_{1j} ... \frac{\varphi_j}{\varphi_{j-1}}h_{j-1.j}h_{jj} ... h_{ij} \ \frac{\varphi_{i+1}}{\varphi_i}h_{i+1.j} ... \frac{\varphi_n}{\varphi_i}h_{nj} \right\|^T,$$

$i = 2, ..., n, i > j$. That is all elements of the first column of the matrix $\Phi_i H$ are divisible by $\dfrac{\varphi_i}{\varphi_1}$, the second are divisible by $\dfrac{\varphi_i}{\varphi_2}$, etc., $(i-1)$th are divisible by $\dfrac{\varphi_i}{\varphi_{i-1}}$. Therefore,
$$\Phi_i H = K_i \Phi_i,$$

where the matrix $K_i$ is a fraction of the right-division of the matrix $\Phi_i H$ to $\Phi_i$. The matrix $\Phi_i$ is nonsingular and $H \in \mathrm{GL}_n(R)$. Hence, $K_i \in \mathrm{GL}_n(R)$. Thus,

$$K_i^{-1} \Phi_i H = \Phi_i.$$

Right-multiply this equality by $n \times m$ matrix $S$, we get

$$K_i^{-1} \Phi_i H S = \Phi_i S.$$

It follows that
$$\Phi_i S \overset{l}{\sim} \Phi_i H S,$$

$i = 2, ..., n$. □

Denote by
$$F_i = \mathrm{diag}\left(\frac{\varphi_i}{\varphi_{i-1}}, ..., \frac{\varphi_i}{\varphi_{i-1}}, \underbrace{1, ..., 1}_{n-i+1}\right), \ \ i = 2, ..., n.$$

**Lemma 3.6.** *If $H \in \mathbf{G}_\Phi$ then there are invertible matrices $H_i$ such that $F_i H = H_i F_i$, $i = 2, ..., n$.*

**Proof.** Consider the matrix

$$\Delta^i = \mathrm{diag}\left(\varphi_{i-1}, ..., \varphi_{i-1}, \varphi_i, \underbrace{1, ..., 1}_{n-i+1}\right).$$

Denote by
$$F_i = \mathrm{diag}\left(\frac{\varphi_i}{\varphi_{i-1}}, ..., \frac{\varphi_i}{\varphi_{i-1}}, \underbrace{1, ..., 1}_{n-i+1}\right).$$

The proof of Lemma 3.5, implies that there is an invertible matrix $H_i$ such that $F_i H = H_i F_i$, $i = 2, ..., n$. □

Let $P = \|p_{ij}\|_1^n$ be an invertible matrix and

$$P_{ij} = \left\|\begin{matrix} p_{ij} & p_{i.j+1} & ... & p_{in} \\ p_{i+1.j} & p_{i+1.j+1} & ... & p_{i+1.n} \\ ... & ... & ... & ... \\ p_{nj} & p_{n.j+1} & ... & p_{nn} \end{matrix}\right\|, \quad P_j = \left\|\begin{matrix} p_{1j} & p_{1.j+1} & ... & p_{1n} \\ p_{2j} & p_{2.j+1} & ... & p_{2n} \\ ... & ... & ... & ... \\ p_{nj} & p_{n.j+1} & ... & p_{nn} \end{matrix}\right\|$$

its submatrices, $i = 2, ..., n$, $j = 1, ..., n$. The Smith form of the matrix $P_{ij}$, where $i > j$, has the form

$$S(P_{ij}) = \left\|Q_j^i \quad \mathbf{0}\right\|, Q_j^i = \mathrm{diag}(q_{j1}^i, ..., q_{j,n-i+1}^i),$$

**117**

and

$$S(P_{ij}) = \left\|\begin{matrix} Q_j^i \\ \mathbf{0} \end{matrix}\right\|, \quad Q_j^i = \operatorname{diag}(q_{j1}^i, ..., q_{j,n-j+1}^i),$$

if $i \leqslant j$.

**Proposition 3.8.** *The Smith form of the matrix* $F_i P_j$, $i = 2, ..., n$, $j = 1, ..., n$, *is*

$$\left\|\begin{matrix} \mathrm{E}_j^i \\ \mathbf{0} \end{matrix}\right\|,$$

*where*

$$\mathrm{E}_j^i = \begin{cases} \operatorname{diag}\left(\left(\dfrac{\varphi_i}{\varphi_{i-1}}, q_{j1}^i\right), ..., \left(\dfrac{\varphi_i}{\varphi_{i-1}}, q_{j,n-i+1}^i\right)\right) \oplus \dfrac{\varphi_i}{\varphi_{i-1}} E_{i-j}, \text{ where } i > j; \\ \operatorname{diag}\left(\left(\dfrac{\varphi_i}{\varphi_{i-1}}, q_{j1}^i\right), ..., \left(\dfrac{\varphi_i}{\varphi_{i-1}}, q_{j,n-j+1}^i\right)\right), \text{ where } i \leqslant j. \end{cases}$$

**Proof.** First we treat the case $i > j$, $n - i + 1 > j - 1$. There are invertible matrices $U, V$ such that

$$U P_{ij} V = S(P_{ij}) = \left\| Q_j^i \quad \mathbf{0} \right\|$$

is the Smith form of the matrix $P_{ij}$, where

$$Q_j^i = \operatorname{diag}(q_{j1}^i, ..., q_{j,n-i+1}^i).$$

Then

$$(I_{i-1} \oplus U) P_j V = \left\|\begin{matrix} M_1 & M_2 \\ Q_j^i & \mathbf{0} \end{matrix}\right\| = D_1.$$

The matrix $D_1$ is primitive. It follows that

$$\left\|\begin{matrix} M_2 \\ \mathbf{0} \end{matrix}\right\|$$

is primitive. That is $(i-1) \times (i-j)$ matrix $M_2$ is also primitive. Since $i > j$, we get $i - 1 \geqslant i - j$. So, there is an invertible matrix $L$ such that

$$L M_2 = \left\|\begin{matrix} I_{i-j} \\ \mathbf{0} \end{matrix}\right\|.$$

Then

$$(L \oplus I_{n-i+1}) D_1 = \left\|\begin{matrix} K_1 & I_{i-j} \\ K_2 & \mathbf{0} \\ Q_j^i & \mathbf{0} \end{matrix}\right\| = D_2.$$

Consequently,

$$D_2 \left\|\begin{matrix} I_{n-i+1} & \mathbf{0} \\ -K_1 & I_{i-j} \end{matrix}\right\| = \left\|\begin{matrix} \mathbf{0} & I_{i-j} \\ K_2 & \mathbf{0} \\ Q_j^i & \mathbf{0} \end{matrix}\right\| = D_3,$$

where $K_2$ is $(j-1) \times (n-i+1)$ matrix. Consider the primitive matrix $P^i$, consisting of the last $n-i+1$ rows of the matrix $P$:

$$P^i = \left\| \begin{array}{ccc|ccc} p_{i1} & ... & p_{i.j-1} & p_{ij} & ... & p_{in} \\ ... & ... & ... & ... & ... & ... \\ p_{n1} & ... & p_{n.j-1} & p_{nj} & ... & p_{nn} \end{array} \right\| = \left\| \begin{array}{ccc|c} p_{i1} & ... & p_{i.j-1} & \\ ... & ... & ... & P_{ij} \\ p_{n1} & ... & p_{n.j-1} & \end{array} \right\|.$$

Then

$$UP^i(I_{j-1} \oplus V) = \left\| \begin{array}{cccccc} p'_{i1} & ... & p'_{i.j-1} & q^i_{j1} & & 0 \\ ... & ... & ... & & \ddots & \mathbf{0} \\ p'_{n1} & ... & p'_{n.j-1} & 0 & & q^i_{j,n-i+1} \end{array} \right\|.$$

This matrix is primitive. It follows that the matrix

$$\left\| \begin{array}{cccccc} p'_{i1} & ... & p'_{i.j-1} & q^i_{j1} & & 0 \\ ... & ... & ... & & \ddots & \\ p'_{n1} & ... & p'_{n.j-1} & 0 & & q^i_{j,n-i+1} \end{array} \right\|$$

is primitive. So $d$-matrix $\mathrm{diag}(q^i_{j1},...,q^i_{j,n-i+1})$ is complemented to a primitive matrix by $j-1$ columns. According to Corollary 3.6, number of unit invariant factor of the matrix $\mathrm{diag}(q^i_{j1},...,q^i_{j,n-i+1})$ does not exceed $j-1$. It follows that

$$q^i_{j1} = q^i_{j2} = ... = q^i_{jt} = 1,$$

where $t = n-i-j+2$, i.e.,

$$D_3 = \left\| \begin{array}{ccc} \mathbf{0} & \mathbf{0} & I_{i-j} \\ K^1_2 & K^2_2 & \mathbf{0} \\ I_t & \mathbf{0} & \mathbf{0} \\ \mathbf{0} & Q^i_{j,t+1} & \mathbf{0} \end{array} \right\|,$$

where $Q^i_{j,t+1} = \mathrm{diag}(q^i_{j,t+1},...,q^i_{j,n-i+1})$, and $K^2_2$ is $j-1 \times j-1$ matrix. Then

$$\left\| \begin{array}{cccc} I_{i-j} & \mathbf{0} & \mathbf{0} & \mathbf{0} \\ \mathbf{0} & I_{j-1} & -K^1_2 & \mathbf{0} \\ \mathbf{0} & \mathbf{0} & I_t & \mathbf{0} \\ \mathbf{0} & \mathbf{0} & \mathbf{0} & I_{j-1} \end{array} \right\| D_3 = \left\| \begin{array}{ccc} \mathbf{0} & \mathbf{0} & I_{i-j} \\ \mathbf{0} & K^2_2 & \mathbf{0} \\ I_t & \mathbf{0} & \mathbf{0} \\ \mathbf{0} & Q^i_{j,t+1} & \mathbf{0} \end{array} \right\| = D_4.$$

There is a matrix $S \in \mathrm{GL}_{j-1}(R)$ such that

$$SK^2_2 = \left\| \begin{array}{ccccc} k_{1,t+1} & k_{1,t+2} & ... & k_{1,\bar{i}-1} & k_{1\bar{i}} \\ k_{2,t+1} & k_{2,t+2} & ... & k_{2,\bar{i}-1} & 0 \\ ... & ... & ... & ... & ... \\ k_{j-1,t+1} & 0 & ... & 0 & 0 \end{array} \right\| = \bar{K}^2_2, \tag{3.5}$$

where $\bar{i} = n - i + 1$. Hence,

$$(I_{i-j} \oplus S \oplus I_{\bar{i}})D_4 = \left\| \begin{matrix} \mathbf{0} & \mathbf{0} & I_{i-j} \\ \mathbf{0} & \bar{K}_2^2 & \mathbf{0} \\ I_t & \mathbf{0} & \mathbf{0} \\ \mathbf{0} & Q_{j,t+1}^i & \mathbf{0} \end{matrix} \right\| = D_5.$$

Thus, $D_5 = NP_jM$, where the matrices $N$ and $M$ are the product of all left-multiply and right-multiply $P_j$ by invertible matrices. The matrix $N$ has the form

$$N = \left\| \begin{matrix} N_1 & N_2 \\ \mathbf{0} & N_3 \end{matrix} \right\|,$$

where $N_3$ is $n - i + 1 \times n - i + 1$ matrix. Then

$$F_i N = \left\| \begin{matrix} \dfrac{\varphi_i}{\varphi_{i-1}}N_1 & \dfrac{\varphi_i}{\varphi_{i-1}}N_2 \\ \mathbf{0} & N_3 \end{matrix} \right\| = \left\| \begin{matrix} N_1 & \dfrac{\varphi_i}{\varphi_{i-1}}N_2 \\ \mathbf{0} & N_3 \end{matrix} \right\| F_i = \bar{N}F_i.$$

The matrix $F_i$ is nonsingular. Therefore, $\det N = \det \bar{N}$. Consequently, the matrix $\bar{N}$ is invertible. Then

$$F_i D_5 = F_i N P_j M = \bar{N}F_i P_j M.$$

It follows that $F_i P_j \sim F_i D_5$. It is obvious

$$F_i D_5 = \underbrace{\left\| \begin{matrix} 0 & 0 & I_{i-j} \\ 0 & \bar{K}_2^2 \Gamma_{j,t+1}^i & 0 \\ I_t & 0 & 0 \\ 0 & \Delta_{j,t+1}^i & 0 \end{matrix} \right\|}_{D_6} \times$$

$$\times \Big( I_t \oplus \mathrm{diag}((\psi_i, q_{j,t+1}^i), ..., (\psi_i, q_{j\bar{i}}^i)) \oplus \psi_i I_{i-j} \Big),$$

where

$$\psi_i = \frac{\varphi_i}{\varphi_{i-1}}, \ \Gamma_{j,t+1}^i = \mathrm{diag}(\gamma_{t+1}, \gamma_{t+2}, ..., \gamma_{\bar{i}}), \gamma_l = \frac{\psi_i}{(\psi_i, q_{jl}^i)},$$

$$\Delta_{j,t+1}^i = \mathrm{diag}(\delta_{j,t+1}^i, \delta_{j,t+2}^i, ..., \delta_{j\bar{i}}^i), \delta_{jl}^i = \frac{q_{jl}^i}{(\psi_i, q_{jl}^i)},$$

$l = t+1, t+2, ..., \bar{i}$. To complete the proof, we need to show that the matrix $D_6$ is primitive. It sufficiently to prove that the matrix

$$D_7 = \left\| \begin{matrix} \bar{K}_2^2 \Gamma_{j,t+1}^i \\ \Delta_{j,t+1}^i \end{matrix} \right\|$$

is primitive. Primitiveness of the matrix $D_5$ implies that the matrix

$$\left\| \begin{matrix} \bar{K}_2^2 \\ Q_{j,t+1}^i \end{matrix} \right\|.$$

is primitive. By Lemma 3.2, the equalities

$$(q^i_{j,t+1}, k_{j-1,t+1}) = (q^i_{j,t+2}, k_{j-2,t+2}) = ... = (q^i_{j\bar{i}}, k_{1\bar{i}}) = 1 \qquad (3.6)$$

are satisfied. Since $\delta^i_{j\bar{i}}|q^i_{j\bar{i}}$ and taking into account (3.6), we get

$$(q^i_{j\bar{i}}, k_{1\bar{i}}) = 1 \Rightarrow (\delta^i_{j\bar{i}}, k_{1\bar{i}}) = 1.$$

Also the equality

$$(\delta^i_{j\bar{i}}, \gamma_{\bar{i}}) = 1$$

is fulfilled. It follows that

$$(\gamma_{\bar{i}} k_{1\bar{i}}, \delta^i_{j\bar{i}}) = 1.$$

There are $u$ and $v$ such that

$$u\gamma_{\bar{i}} k_{1\bar{i}} + v\delta^i_{j\bar{i}} = 1.$$

Taking into account that $\delta^i_{j,t+1}|\delta^i_{j,t+2}|...|\delta^i_{j\bar{i}}$, we have

$$\left\| \begin{array}{ccc} u & \mathbf{0} & v \\ \mathbf{0} & I_{2j-4} & \mathbf{0} \\ -\delta^i_{j\bar{i}} & \mathbf{0} & \gamma_{\bar{i}} k_{1\bar{i}} \end{array} \right\| D_7 =$$

$$= \left\| \begin{array}{ccccc} * & * & ... & * & 1 \\ \gamma_{t+1} k_{2,t+1} & \gamma_{t+2} k_{2,t+2} & ... & \gamma_{\bar{i}-1} k_{2,\bar{i}-1} & 0 \\ ... & ... & ... & ... & ... \\ \gamma_{t+1} k_{j-1,t+1} & 0 & ... & 0 & 0 \\ \delta^i_{j,t+1} & 0 & ... & 0 & 0 \\ 0 & \delta^i_{j,t+2} & & 0 & 0 \\ \vdots & & \ddots & & \vdots \\ 0 & 0 & & \delta^i_{j,\bar{i}-1} & 0 \\ \delta^i_{j\bar{i}} s_{t+1} & \delta^i_{j\bar{i}} s_{t+2} & ... & \delta^i_{j\bar{i}} s_{\bar{i}-1} & 0 \end{array} \right\| \sim$$

$$\sim \left\| \begin{array}{ccccc} 0 & 0 & ... & 0 & 1 \\ \gamma_{t+1} k_{2,t+1} & \gamma_{t+2} k_{2,t+2} & ... & \gamma_{\bar{i}-1} k_{2,\bar{i}-1} & 0 \\ ... & ... & ... & ... & ... \\ \gamma_{t+1} k_{j-1,t+1} & 0 & ... & 0 & 0 \\ \delta^i_{j,t+1} & 0 & ... & 0 & 0 \\ 0 & \delta^i_{j,t+2} & & 0 & 0 \\ \vdots & & \ddots & & \vdots \\ 0 & 0 & & \delta^i_{j,\bar{i}-1} & 0 \\ 0 & 0 & ... & 0 & 0 \end{array} \right\| = D_8.$$

Continuing the process, we show that the matrix $D_8$ is primitive. So the matrix $D_6$ is primitive. Noting that

$$I_t \oplus \mathrm{diag}((\psi_i, q^i_{j,t+1}), ..., (\psi_i, q^i_{j\bar{i}})) \oplus \psi_i I_{i-j} =$$

$$= \operatorname{diag}\left((\psi_i, q_{j1}^i), ..., (\psi_i, q_{j\bar{i}}^i)\right) \oplus \psi_i I_{i-j} = \mathrm{E}_j^i,$$

we complete this case.

Suppose that $i > j, n - i + 1 \leqslant j - 1$. By right-sided transformations from $\mathrm{GL}_{n-j+1}(R)$ and left-sided transformations that do not change the form of the matrix $F_i P_j$ Smith's form, we reduce the matrix $P_j$ to the form

$$\left\|\begin{array}{cc} \mathbf{0} & I_{i-j} \\ L_1 & \mathbf{0} \\ \mathbf{0} & \mathbf{0} \\ Q_{j,t+1}^i & \mathbf{0} \end{array}\right\|.$$

Similarly if $i < j$, the matrix $P_j$ is equivalent to

$$\left\|\begin{array}{c} \mathbf{0} \\ L_2 \\ \mathbf{0} \\ Q_{j,t+1}^i \end{array}\right\|,$$

where matrices $L_1, L_2$ have the form (3.5). That is, these cases are partial cases of the prime, and therefore they are proved by a similar scheme. The claim is proved. □

Set
$$S(\Phi, P_{ij}) = \operatorname{diag}\left(\left(\frac{\varphi_i}{\varphi_{i-1}}, q_{j1}^i\right), ..., \left(\frac{\varphi_i}{\varphi_{i-1}}, q_{jk}^i\right)\right),$$

where $k = n - i + 1$, if $i > j$, and $k = n - j + 1$, if $i \leqslant j$.

**Theorem 3.3.** *If $H \in \mathbf{G}_\Phi$, then the diagonal elements of the matrices $S(\Phi, P_{ij})$, $S(\Phi, (HP)_{ij},)$ $i = 2, ..., n, j = 1, ..., n$ are associates.*

**Proof.** By Lemma 3.6,
$$F_i H = H_i F_i, \;\; i = 2, ..., n.$$

The matrices $F_i$ are nonsingular. So
$$\det H = \det H_i.$$

It follows that the matrices $H_i$ are invertible. Denote by $P_j$ and $U_j$ the matrices consisting of the last $n - j + 1$ columns of matrices $P$ and $U$, respectively. Noting that $U_j = HP_j$, we get
$$F_i U_j = F_i H P_j = H_i F_i P_j \sim F_i P_j.$$

That is an invariant factor of Smith forms of the matrices $F_i U_j$ and $F_i P_j$ are associates. Thus, a nonzero element of the matrices $S(\Phi, P_{ij}), S(\Phi', U_{ij})$, $i = 2, ..., n, j = 1, ..., n$, has the form
$$\left(\frac{\varphi_i}{\varphi_{i-1}}, q_{js}^i\right), \left(\frac{\varphi_i}{\varphi_{i-1}}, q_{js}^i l_{ijs}\right),$$

where $l_{ijs} \in U(R)$, and obviously, are associates. □

# DIVISIBILITY AND ASSOCIATIVITY OF MATRICES

*This section is devoted to systematic study of the matrix divisors. Whenever in this section R be an elementary divisor ring.*

## 4.1. Structure and properties of generating set and its elements

Let $A$ and $B$ be $n \times n$ matrices over $R$. There exist invertible matrices $P_A, P_B, Q_A, Q_B$ such that

$$P_A A Q_A = \mathrm{diag}(\varepsilon_1, ..., \varepsilon_k, 0, ..., 0) = \mathrm{E},$$
$$P_B B Q_B = \mathrm{diag}(\varphi_1, ..., \varphi_t, 0, ..., 0) = \Phi,$$

where $\varepsilon_k \neq 0$, $\varphi_t \neq 0$, $\varepsilon_i | \varepsilon_{i+1}$, $\varphi_j | \varphi_{j+1}$, $i = 1, ..., k-1$, $j = 1, ..., t-1$.

**Definition 4.1.** The **generating set** is called the set

$$\mathbf{L}(\mathrm{E}, \Phi) = \{L \in \mathrm{GL}_n(R)| \; \exists \, S \in M_n(R) \colon LE = \Phi S\}.$$

Let us describe the structure of generating set elements.

**Theorem 4.1.** *Let* $\mathrm{E} = \mathrm{diag}(\varepsilon_1, ..., \varepsilon_k, 0, ..., 0)$,

$$\Phi = \mathrm{diag}(\varphi_1, ..., \varphi_t, 0, ..., 0),$$

*where* $\varepsilon_k \neq 0$, $\varphi_t \neq 0$, $\varepsilon_i | \varepsilon_{i+1}$, $\varphi_j | \varphi_{j+1}$, $i = 1, ..., k-1$, $j = 1, ..., t-1$, $k, t \leq n$.

*If* $\Phi | \mathrm{E}$, *then the set* $\mathbf{L}(\mathrm{E}, \Phi)$ *consists of all invertible matrices of the form*

$$\begin{Vmatrix} L_1 & * \\ L_2 & * \\ \mathbf{0} & * \end{Vmatrix}, \qquad (4.1)$$

*where*

$$L_1 = \begin{Vmatrix} l_{11} & l_{12} & ... & l_{1.k-1} & l_{1k} \\ \dfrac{\varphi_2}{(\varphi_2, \varepsilon_1)} l_{21} & l_{22} & ... & l_{2.k-1} & l_{2k} \\ ... & ... & ... & ... & ... \\ \dfrac{\varphi_k}{(\varphi_k, \varepsilon_1)} l_{k1} & \dfrac{\varphi_k}{(\varphi_k, \varepsilon_2)} l_{k2} & ... & \dfrac{\varphi_k}{(\varphi_k, \varepsilon_{k-1})} l_{k.k-1} & l_{kk} \end{Vmatrix},$$

$$L_2 = \left\| \begin{array}{ccc} \dfrac{\varphi_{k+1}}{(\varphi_{k+1},\varepsilon_1)}l_{k+1.1} & ... & \dfrac{\varphi_{k+1}}{(\varphi_{k+1},\varepsilon_k)}l_{k+1.k} \\ ... & ... & ... \\ \dfrac{\varphi_t}{(\varphi_t,\varepsilon_1)}l_{t1} & ... & \dfrac{\varphi_t}{(\varphi_t,\varepsilon_k)}l_{tk} \end{array} \right\|.$$

*If* $\Phi \nmid$ E, *then* $\mathbf{L}(E,\Phi) = \varnothing$.

**Proof.** Let $L = \|p_{ij}\|_1^n \in \mathbf{L}(E,\Phi)$. That is $LE = \Phi S'$, where $S = \|s_{ij}\|_1^n \in \in M_n(R)$. So

$$\left\| \begin{array}{cccccc} \varepsilon_1 p_{11} & ... & \varepsilon_k p_{1k} & 0 & ... & 0 \\ ... & ... & ... & ... & ... & ... \\ \varepsilon_1 p_{t1} & ... & \varepsilon_k p_{tk} & 0 & ... & 0 \\ \varepsilon_1 p_{t+1.1} & ... & \varepsilon_k p_{t+1.k} & 0 & ... & 0 \\ ... & ... & ... & ... & ... & ... \\ \varepsilon_1 p_{n1} & ... & \varepsilon_k p_{nk} & 0 & ... & 0 \end{array} \right\| = \left\| \begin{array}{ccc} \varphi_1 s_{11} & ... & \varphi_1 s_{1n} \\ ... & ... & ... \\ \varphi_t s_{t1} & ... & \varphi_t s_{tn} \\ 0 & ... & 0 \\ ... & ... & ... \\ 0 & ... & 0 \end{array} \right\|. \qquad (4.2)$$

It follows that

$$\left\| \begin{array}{ccc} p_{t+1.1} & ... & p_{t+1.k} \\ ... & ... & ... \\ p_{n1} & ... & p_{nk} \end{array} \right\| = \mathbf{0}.$$

That is, the invertible matrix $L$ contains zero $(n-t) \times k$ submatrix. Applying Corollary 3.1, we conclude that the element $p_{t+1.k}$ lies below the main diagonal, i.e., $t \geq k$. From (4.2) we conclude that $\varphi_i | \varepsilon_j p_{ij}$, $i = 1,...,t, j = 1,...,k$. Consequently,

$$\frac{\varphi_i}{(\varphi_i,\varepsilon_j)} \left| \frac{\varepsilon_j}{(\varphi_i,\varepsilon_j)} p_{ij}. \right.$$

It follows that $\dfrac{\varphi_i}{(\varphi_i,\varepsilon_j)} \Big| p_{ij}$, i.e.,

$$p_{ij} = \frac{\varphi_i}{(\varphi_i,\varepsilon_j)} p'_{ij}.$$

Therefore the matrix $L$ has the form

$$L = \left\| \begin{array}{cccccc} \dfrac{\varphi_1}{(\varphi_1,\varepsilon_1)}p'_{11} & ... & \dfrac{\varphi_1}{(\varphi_1,\varepsilon_k)}p'_{1k} & p_{1.k+1} & ... & p_{1n} \\ ... & ... & ... & ... & ... & ... \\ \dfrac{\varphi_k}{(\varphi_k,\varepsilon_1)}p'_{k1} & ... & \dfrac{\varphi_k}{(\varphi_k,\varepsilon_k)}p'_{kk} & p_{k.k+1} & ... & p_{kn} \\ ... & ... & ... & ... & ... & ... \\ \dfrac{\varphi_t}{(\varphi_t,\varepsilon_1)}p'_{t1} & ... & \dfrac{\varphi_t}{(\varphi_t,\varepsilon_k)}p'_{tk} & p_{t.k+1} & ... & p_{tn} \\ 0 & ... & 0 & p_{t+1.k+1} & ... & p_{t+1.n} \\ ... & ... & ... & ... & ... & ... \\ 0 & ... & 0 & p_{n.k+1} & ... & p_{nn} \end{array} \right\|. \qquad (4.3)$$

Consider submatrices

$$L_i = \left\| \begin{matrix} p_{i1} & ... & p_{ii} \\ ... & ... & ... \\ p_{n1} & ... & p_{ni} \end{matrix} \right\|,$$

$i = 1, ..., k$, of the matrix $L$. Since

$$\frac{\varphi_{i+r}}{(\varphi_{i+r}, \varepsilon_{j-l})} = \frac{\varphi_i}{(\varphi_i, \varepsilon_j)} \frac{\left(\varphi_{i+r}, \frac{\varphi_{i+r}}{\varphi_i}\varepsilon_j\right)}{(\varphi_{i+r}, \varepsilon_{j-l})}, \quad l < j, \tag{4.4}$$

then $\frac{\varphi_i}{(\varphi_i, \varepsilon_i)}\Big| \langle L_i \rangle_1$. The matrices $L_i$ have sizes $(n - i + 1) \times i$. Taking into account

$$n - i + 1 + i = n + 1 > n,$$

and using Proposition 3, we get $\langle L_i \rangle_1 | \det L \in U(R)$. Therefore $\frac{\varphi_i}{(\varphi_i, \varepsilon_i)} \in U(R)$, $i = 1, ..., k$. Since g.c.d. is determined up to divisors of unit, we can assume that

$$\frac{\varphi_i}{(\varphi_i, \varepsilon_i)} = 1.$$

That is $\varphi_i = (\varphi_i, \varepsilon_i)$. It follows that $\varphi_i | \varepsilon_i, i = 1, ..., k$. Hence, $\Phi | E$. Consequently, if $\Phi \nmid E$, then $\mathbf{L}(E, \Phi) = \varnothing$.

Let $\Phi | E$ it follows that $\varphi_i | \varepsilon_{i+j}, j = 0, ..., k - i$. Then

$$\frac{\varphi_i}{(\varphi_i, \varepsilon_{i+j})} \in U(R).$$

Thus, elements $p_{i.i+j}$ do not have any restriction. That is, the matrix $L$ has the form (4.1).

Conversely, suppose $\Phi | E$ and the matrix $L$ has the form (4.1). Then it is easy to make sure that

$$LE = \Phi S,$$

where

$$S = \left\| \begin{matrix} M_1 & \mathbf{0} \\ M_2 & \mathbf{0} \\ M_3 & M_4 \end{matrix} \right\|,$$

$$M_1 = \left\| \begin{matrix} \frac{\varepsilon_1}{\varphi_1}l_{11} & \frac{\varepsilon_2}{\varphi_1}l_{12} & ... & \frac{\varepsilon_{k-1}}{\varphi_1}l_{1.k-1} & \frac{\varepsilon_k}{\varphi_1}l_{1k} \\ \frac{\varepsilon_1}{(\varphi_2, \varepsilon_1)}l_{21} & \frac{\varepsilon_2}{\varphi_2}l_{22} & ... & \frac{\varepsilon_{k-1}}{\varphi_2}l_{2.k-1} & \frac{\varepsilon_k}{\varphi_2}l_{2k} \\ ... & ... & ... & ... & ... \\ \frac{\varepsilon_1}{(\varphi_k, \varepsilon_1)}l_{k1} & ... & ... & \frac{\varepsilon_{k-1}}{(\varphi_k, \varepsilon_{k-1})}l_{k.k-1} & \frac{\varepsilon_k}{\varphi_k}l_{kk} \end{matrix} \right\|,$$

$$M_2 = \left\| \begin{matrix} \frac{\varepsilon_1}{(\varphi_{k+1}, \varepsilon_1)}l_{k+1.1} & ... & \frac{\varepsilon_k}{(\varphi_{k+1}, \varepsilon_k)}l_{k+1.k} \\ ... & ... & ... \\ \frac{\varepsilon_1}{(\varphi_t, \varepsilon_1)}l_{t1} & ... & \frac{\varepsilon_k}{(\varphi_t, \varepsilon_k)}l_{tk} \end{matrix} \right\|,$$

$M_3$, $M_4$ are arbitrary matrices of corresponding sizes. The theorem is proved. $\square$

The elements of matrices from the generating set and the Zelisko group have similar properties (see Property 2.4, page 64):

**Property 4.1.** *Element* $\dfrac{\varphi_i}{(\varphi_i, \varepsilon_j)}$, $i > j$, *is a divisor of all elements of the matrices* $L_1$, $L_2$, *which are outlined by a rectangle with vertices* $(i, 1)$, $(i, j)$, $(t, j)$, $(t, 1)$:

$$
\begin{array}{cccc}
f_{i1}h_{i1} & f_{i2}h_{i2} & ... & f_{ij}h_{ij} \\
f_{i+1.1}h_{i+1.1} & f_{i+1.2}h_{i+1.2} & ... & f_{i+1.j}h_{i+1.j} \\
... & ... & ... & ... \\
f_{t1}h_{t1} & f_{t2}h_{t2} & ... & f_{tj}h_{tj},
\end{array}
$$

*where*

$$f_{pq} = \frac{\varphi_p}{(\varphi_p, \varepsilon_q)}.$$

*It means*

$$f_{ij}|f_{i+p.j-q}, \ p = 0, 1, ..., n-i, \ q = 0, 1, ..., j-1.$$

**Proof** follows from equality (4.4). $\square$

We describe the action of the Zelisko group on the generating set.

**Property 4.2.** *Let* $H \in \mathbf{G}_\Phi$. *Then*

$$H\mathbf{L}(\mathrm{E}, \Phi) = \mathbf{L}(\mathrm{E}, \Phi).$$

**Proof.** Let $H \in \mathbf{G}_\Phi$ and $L \in \mathbf{L}(\mathrm{E}, \Phi)$, i.e.,

$$H\Phi = \Phi K, \ LE = \Phi S,$$

where $K \in \mathrm{GL}_n(R)$, $S \in M_n(R)$. Then

$$(HL)\mathrm{E} = H(LE) = H\Phi S = (H\Phi)S = \Phi(KS).$$

It means that $H\mathbf{L}(\mathrm{E}, \Phi) \subseteq \mathbf{L}(\mathrm{E}, \Phi)$.

On the other hand, each element $L$ from $\mathbf{L}(\mathrm{E}, \Phi)$ can be written as $L = HL_1$, where $L_1 = H^{-1}L \in \mathbf{L}(\mathrm{E}, \Phi)$. Therefore, $\mathbf{L}(\mathrm{E}, \Phi) \subseteq H\mathbf{L}(\mathrm{E}, \Phi)$. Consequently, $H\mathbf{L}(\mathrm{E}, \Phi) = \mathbf{L}(\mathrm{E}, \Phi)$. $\square$

**Corollary 4.1.** *The equality*

$$\mathbf{G}_\Phi\mathbf{L}(\mathrm{E}, \Phi) = \mathbf{L}(\mathrm{E}, \Phi).$$

*is fulfilled.* $\square$

**Property 4.3.** *Let* $H \in \mathbf{G}_\mathrm{E}$. *Then*

$$\mathbf{L}(\mathrm{E}, \Phi)H = \mathbf{L}(\mathrm{E}, \Phi).$$

**Proof** is similar to proving Property 4.1. $\square$

**Corollary 4.2.** *The equality*

$$\mathbf{L}(\mathrm{E}, \Phi)\mathbf{G}_\mathrm{E} = \mathbf{L}(\mathrm{E}, \Phi).$$

*is fulfilled.* $\square$

**Property 4.4.** *Let* $\mathbf{L}(\mathrm{E}, \Phi) \neq \varnothing$. *Then*

$$\mathbf{G}_\Phi \mathbf{G}_\mathrm{E} \subseteq \mathbf{L}(\mathrm{E}, \Phi)$$

**Proof.** Let $H \in \mathbf{G}_\Phi$, $K \in \mathbf{G}_\mathrm{E}$. Since $\mathbf{L}(\mathrm{E}, \Phi) \neq \varnothing$ by Theorem 4.1, $\mathrm{E} = \Phi\Delta$. Then the equalities

$$(HK)\mathrm{E} = HEK_1 = H\Phi\Delta K_1 = \Phi(H_1 \Delta K_1)$$

are fulfilled. Hence, $HK \in \mathbf{L}(\mathrm{E}, \Phi)$. Consequently, $\mathbf{G}_\Phi \mathbf{G}_\mathrm{E} \subseteq \mathbf{L}(\mathrm{E}, \Phi)$. $\qquad\square$

## 4.2. Generating set and divisibility of matrices

The generating set plays a basic role in the description of matrix divisors.

**Theorem 4.2.** *The matrix* $B = P_B^{-1}\Phi Q_B^{-1}$ *is left divisor of the matrix* $A = P_A^{-1}\mathrm{E}Q_A^{-1}$, *i.e.,* $A = BC$ *if and only if* $P_B = LP_A$, *where* $L \in \mathbf{L}(\mathrm{E}, \Phi)$.

**Proof. Necessity.** Left-multiplying $A$ by $P_B$, we get

$$P_B A = P_B(BC) = (P_B B)C = (\Phi Q_B^{-1})C = \Phi(Q_B^{-1}C).$$

On the other hand,

$$P_B A = (P_B P_A^{-1})(P_A A) = (P_B P_A^{-1})\mathrm{E}Q_A^{-1}.$$

Hence,

$$(P_B P_A^{-1})\mathrm{E}Q_A^{-1} = \Phi(Q_B^{-1}C).$$

This yields,

$$(P_B P_A^{-1})\mathrm{E} = \Phi S,$$

where $S = Q_B^{-1}CQ_A$. Therefore,

$$P_B P_A^{-1} = L \in \mathbf{L}(\mathrm{E}, \Phi).$$

Consequently, $P_B = LP_A$.

**Sufficiency.** Let

$$(P_B P_A^{-1})\mathrm{E} = \Phi S.$$

Then

$$P_B A = P_B(P_A^{-1}\mathrm{E}Q_A^{-1}) = (P_B P_A^{-1})\mathrm{E}Q_A^{-1} = \Phi SQ_A^{-1}.$$

It means that

$$A = P_B^{-1}\Phi SQ_A^{-1} = (P_B^{-1}\Phi Q_B^{-1})(Q_B SQ_A^{-1}) = BC,$$

where $C = Q_B SQ_A^{-1}$. $\qquad\square$

Combining the results of Theorems 4.1, 4.2, we obtain.

**Theorem 4.3.** The matrix $B = P_B^{-1}\Phi Q_B^{-1}$ is a left divisor of the matrix $A = P_A^{-1}\mathrm{E}Q_A^{-1}$ if and only if $\Phi|\mathrm{E}$ and $P_B = LP_A$, where $L$ is invertible matrix of the form (4.1). $\qquad\square$

The obvious consequence of Theorem 4.2 is the following statement.

**Corollary 4.3.** *All divisors of the matrix $A = P_A^{-1} E Q_A^{-1}$ with the Smith form $\Phi$ have the form $(L P_A)^{-1} \Phi Q^{-1}$, where $L \in \mathbf{L}(E, \Phi), Q^{-1} \in \mathrm{GL}_n(R)$.* □

**Theorem 4.4.** *The set $(\mathbf{L}(E, \Phi) P_A)^{-1} \Phi \mathrm{GL}_n(R)$ is the set of all left divisors of the matrix $A = P_A^{-1} E Q_A^{-1}$ with the Smith form $\Phi$.* □

**Proof**. To prove this statement, we must show that the set $(\mathbf{L}(E, \Phi) P_A)^{-1} \Phi \mathrm{GL}_n(R)$ does not depend on the choice of the transforming matrix $P_A$. Let $P_A' \in \mathbf{P}_A$. Then there exists $H \in \mathbf{G}_E$ such that $P_A' = H P_A$. So

$$(\mathbf{L}(E, \Phi) P_A')^{-1} \Phi \mathrm{GL}_n(R) = (\mathbf{L}(E, \Phi) H P_A)^{-1} \Phi \mathrm{GL}_n(R).$$

By Property 4, $\mathbf{L}(E, \Phi) H = \mathbf{L}(E, \Phi)$. Consequently,

$$(\mathbf{L}(E, \Phi) P_A')^{-1} \Phi \mathrm{GL}_n(R) = (\mathbf{L}(E, \Phi) P_A)^{-1} \Phi \mathrm{GL}_n(R),$$

as required. □

The matrices $A, B$ are right (left) associates if there is an invertible matrix $U$ such that $A = BU$ $(A = UB)$. In notation $A \overset{r}{\sim} B$ $(A \overset{l}{\sim} B)$. Obviously, such matrices are equivalent. Therefore, the associated matrices have the same Smith form. Thus, associated right or left matrices $A, B$ can always be written in the form

$$A = P_A^{-1} \Phi Q_A^{-1}, \ B = P_B^{-1} \Phi V_B^{-1}.$$

**Theorem 4.5.** *Let $A = P_A^{-1} \Phi Q_A^{-1}$, $B = P_B^{-1} \Phi Q_B^{-1}$. The following conditions are equivalent:*

1) $A \overset{r}{\sim} B$;
2) $P_B = H P_A$, where $H \in \mathbf{G}_\Phi$;
3) $\mathbf{P}_B = \mathbf{P}_A$.

**Proof.** We will prove that 1) $\Leftrightarrow$ 2). There is an invertible matrix $U$ such that $A = BU$, i.e.,

$$P_A^{-1} \Phi Q_A^{-1} = P_B^{-1} \Phi Q_B^{-1} U.$$

Hence,

$$P_B P_A^{-1} \Phi = \Phi Q_B^{-1} U Q_A.$$

It means that

$$P_B P_A^{-1} = H \in \mathbf{G}_\Phi.$$

Thus $P_B = H P_A$.

On the contrary, let $P_B = H P_A$, where $H \in \mathbf{G}_\Phi$. Since

$$H\Phi = \Phi K,$$

where $K \in \mathrm{GL}_n(R)$, then

$$H^{-1}\Phi = \Phi K^{-1}.$$

So

$$B = P_B^{-1} \Phi Q_B^{-1} = (H P_A)^{-1} \Phi Q_B^{-1} = P_A^{-1} H^{-1} \Phi Q_B^{-1} =$$
$$= P_A^{-1} \Phi K^{-1} Q_B^{-1} = (P_A^{-1} \Phi Q_A^{-1})(Q_A K^{-1} Q_B^{-1}) = AV,$$

where $V = Q_A K^{-1} Q_B^{-1} \in \mathrm{GL}_n(R)$. Hence, $A \overset{r}{\sim} B$. So conditions 1) and 2) are equivalent.

Now we will prove that 2)⇔ 3). By Property 2.2,

$$\mathbf{P}_A = \mathbf{G}_\Phi P_A, \ \ \mathbf{P}_B = \mathbf{G}_\Phi P_B.$$

If $P_B = H P_A$, where $H \in \mathbf{G}_\Phi$, then

$$\mathbf{P}_B = \mathbf{G}_\Phi P_B = \mathbf{G}_\Phi H P_A = \mathbf{G}_\Phi P_A = \mathbf{P}_A.$$

Let $\mathbf{P}_B = \mathbf{P}_A$, i.e., $\mathbf{G}_\Phi P_B = \mathbf{G}_\Phi P_A$. Since $I \in \mathbf{G}_\Phi$, then $P_B \in \mathbf{G}_\Phi P_A$. Therefore, there is the matrix $H \in \mathbf{G}_\Phi$ such that $P_B = H P_A$. □

Let $B$ be the left divisor of the matrix $A$. Then, for an arbitrary invertible matrix $U$ the equality

$$A = BC = (BU)(U^{-1}C)$$

holds. That is, all matrices that are right associates to $B$ are again the left divisors of the matrix $A$. Therefore, of course, describe only unassociated divisors of $A$.

By Theorem 4.4, the set of all left divisors of $A$ with the Smith form $\Phi$ have the form $(\mathbf{L}(\mathrm{E}, \Phi) P_A)^{-1} \Phi \mathrm{GL}_n(R)$. So our task is to select all unassociated matrices from this set.

The Corollary 4 implies that the set $\mathbf{L}(\mathrm{E}, \Phi)$ can be broken by right cosets $\mathbf{G}_\Phi L$, where $L \in \mathbf{L}(\mathrm{E}, \Phi)$. Denote by $\mathbf{W}(\mathrm{E}, \Phi)$ the set of representatives of these cosets.

**Theorem 4.6.** *The set* $(\mathbf{W}(\mathrm{E}, \Phi) P_A)^{-1} \Phi$ *consists of all right unassociated divisors of the matrix $A$ with the Smith form $\Phi$.*

**Proof.** Let $L_1, L_2 \in \mathbf{W}(\mathrm{E}, \Phi) \subset \mathbf{L}(\mathrm{E}, \Phi)$. According to Theorem 4.3, the matrices $(L_1 P_A)^{-1} \Phi$ and $(L_2 P_A)^{-1} \Phi$ are the left divisors of the matrix $A$. Suppose that matrices $B_1 = (L_1 P_A)^{-1} \Phi$ and $B_2 = (L_2 P_A)^{-1} \Phi$ are right associates. Since

$$L_1 P_A \in \mathbf{P}_{B_1}, \ \ L_2 P_A \in \mathbf{P}_{B_2},$$

by virtue of Theorem 4.5, we get

$$L_2 P_A = H L_1 P_A,$$

where $H \in \mathbf{G}_\Phi$. It follows that $L_2 = H L_1$. Since $L_1, L_2 \in \mathbf{W}(\mathrm{E}, \Phi)$, we have $L_1 = L_2$. It means that $(\mathbf{W}(\mathrm{E}, \Phi) P_A)^{-1} \Phi$ consists of all right unassociated divisors of the matrix $A$ with the Smith form $\Phi$.

Now let $B$ be the left divisor of the $A$ matrix with the Smith form $\Phi$. By Theorem 4.3 $B$ has the form $B = (L P_A)^{-1} \Phi Q^{-1}$, where $L \in \mathbf{L}(\mathrm{E}, \Phi)$. Let $W \in \mathbf{W}(\mathrm{E}, \Phi)$ is a representative of an adjacent class $\mathbf{G}_\Phi L$, i.e., $L = HW$, where $H \in \mathbf{G}_\Phi$. Consider the matrix $B_1 = (W P_A)^{-1} \Phi$ of elements of the set $(\mathbf{W}(\mathrm{E}, \Phi) P_A)^{-1} \Phi$. Since $L P_A = H(W P_A)$, by Theorem 4.5, matrices $B$, $B_1$ are right associate. That is, for every divisor of the matrix $A$ with Smith form $\Phi$, in the set $(\mathbf{W}(\mathrm{E}, \Phi)P_A)^{-1}\Phi$, there exists a matrix right associate with it.

A transforming matrix $P_A$ is chosen ambiguously. Let $P \in \mathbf{P}_A$ and $B =$ $= (WP)^{-1}\Phi$, where $W \in \mathbf{W}(\mathrm{E}, \Phi)$. According to Property 2.2, there is $S \in \mathbf{G}_\mathrm{E}$ such that $P = SP_A$. Thus, $B = (WSP_A)^{-1}\Phi$. By Property 4.2, $WS \in \mathbf{L}(\mathrm{E}, \Phi)$. So $WS$ belongs to some coset $\mathbf{G}_\Phi V$, where $V \in \mathbf{W}(\mathrm{E}, \Phi)$. Hence, there is $H \in$ $\in \mathbf{G}_\Phi$ such that $WS = HV$. Then

$$B = (WP)^{-1}\Phi = (WSP_A)^{-1}\Phi = (HVP_A)^{-1}\Phi =$$
$$= (VP_A)^{-1}H^{-1}\Phi = (VP_A)^{-1}\Phi K^{-1} = B_1 K^{-1},$$

where $B_1 \in (\mathbf{W}(\mathrm{E}, \Phi)P_A)^{-1}\Phi$. That is for any matrix in $(\mathbf{W}(\mathrm{E}, \Phi)P)^{-1}\Phi$ there is a right associate matrix in $(\mathbf{W}(\mathrm{E}, \Phi)P_A)^{-1}\Phi$. The converse is also correct. Therefore, for randomly selected transforming matrix $P_A$ the set $(\mathbf{W}(\mathrm{E}, \Phi)P_A)^{-1}\Phi$ consists of all left unassociated divisors of the matrix $A$ with the Smith form $\Phi$. $\qquad\square$

## 4.3. G.c.d. and l.c.m. of a matrix equation $BX = A$ solutions

The equation $BX = A$ was investigated in subsection 2.2. Necessary and sufficient conditions for the solvability of this equation were obtained. However, the method of finding roots that was proposed in Theorem 2.6 (see p. 59) does not allow to investigate their properties. So let us look at this equation from a slightly different point of view. If the equation $BX = A$ is solvable, then there exists a matrix $C$, such that $A = BC$. That is, the matrix $B$ is the left divisor of the matrix $A$. The inverse reasoning will also be correct. Thus, the equation $BX = A$ is solvable if and only if the matrix $B$ is the left divisor of the matrix $A$. We use this approach to find the roots of this equation and study their properties.

**Theorem 4.7.** *Left g.c.d. and left l.c.m. of roots of a solvable matrix equation $BX = A$ are again their solutions.*

**Proof.** Let
$$A = P^{-1}\mathrm{E}Q^{-1}, \quad \mathrm{E} = \mathrm{diag}(\varepsilon_1, ..., \varepsilon_k, 0, ..., 0),$$
and
$$B = P_B^{-1}\Phi Q_B^{-1}, \quad \Phi = \mathrm{diag}(\varphi_1, ..., \varphi_t, 0, ..., 0),$$

where $\mathrm{E}$ and $\Phi$ are $d$-matrices. By Theorem 4.2, the matrix $B$ is a left divisor of the matrix $A$ if and only if $P_B = LP$, where $L \in \mathbf{L}(\mathrm{E}, \Phi) \neq \varnothing$. According to Corollary 4.4, the matrix $B$ has the form $B = (LP)^{-1}\Phi U^{-1}$, where $U \in \mathrm{GL}_n R$. By definition, the matrix $L$ satisfies the equality

$$L\mathrm{E} = \Phi S \qquad (4.5)$$

for some $S$ in $M_n(R)$. It follows that $\mathrm{E} = L^{-1}\Phi S$. Therefore,

$$A = P^{-1}\mathrm{E}Q^{-1} = P^{-1}(L^{-1}\Phi S)Q^{-1} = (LP)^{-1}\Phi SQ^{-1} =$$

$$= ((LP)^{-1}\Phi U^{-1})(USQ^{-1}) = BC,$$

where $C = USQ^{-1}$. Consequently, $C$ is the solution of the equation $BX = A$.

Let $C_1$ be another root. That is $BC_1 = A$. Then

$$A = BC_1 = BC \quad \Rightarrow \quad B(C_1 - C) = \mathbf{0}.$$

Hence, $C_1 - C = F \in Ann_r(B)$. This yields

$$C_1 = C + F \quad \Rightarrow \quad C_1 \in (C + Ann_r(B)).$$

On the contrary, let $C_2 \in (C + Ann_r(B))$. Then $C_2 = C + N$, where $N \in Ann_r(B)$. Therefore,

$$BC_2 = B(C + N) = BC + BN = A + \mathbf{0} = A.$$

Consequently, the set $C + Ann_r(B)$ consists of solutions of the equation $BX = A$. Thus the set $C + Ann_r(B)$ is the set of all solutions of the equation $BX = A$.

By Theorem 1.15, the set $Ann_r(B)$ consists of all matrices of the form

$$U \left\| \begin{matrix} \mathbf{0}_{t \times n} \\ D \end{matrix} \right\|,$$

where $D$ is an arbitrary matrix in $M_{(n-t) \times n}(R)$.

The matrix $S$, which appears in the equality (4.5), as follows from the proof of Theorem 4.1, has the form

$$S = \left\| \begin{matrix} M_1 & \mathbf{0} \\ M_2 & \mathbf{0} \\ M_3 & M_4 \end{matrix} \right\|,$$

where

$$\left\| \begin{matrix} M_1 \\ M_2 \end{matrix} \right\|$$

is $t \times k$ submatrix of the matrix $S$, and $\left\| M_3 \quad M_4 \right\|$ is an arbitrary matrix in $M_{(n-t) \times n}(R)$. Consequently,

$$C + Ann_r(B) = USQ^{-1} + U \left\| \begin{matrix} \mathbf{0}_{t \times n} \\ D \end{matrix} \right\| =$$

$$= U \left( S + \left\| \begin{matrix} \mathbf{0}_{t \times n} \\ DQ \end{matrix} \right\| \right) Q^{-1} = U \left( S + \left\| \begin{matrix} \mathbf{0}_{t \times n} \\ D_1 \end{matrix} \right\| \right) Q^{-1}.$$

Since $Q$ is an invertible matrix then

$$M_{(n-t) \times n}(R)Q = M_{(n-t) \times n}(R).$$

So $D_1 = DN$ is an arbitrary matrix in $M_{(n-t) \times n}(R)$. That is

$$C + Ann_r(B) = U \left( \left\| \begin{matrix} M_1 & \mathbf{0} \\ M_2 & \mathbf{0} \\ M_3 & M_4 \end{matrix} \right\| + \left\| \begin{matrix} \mathbf{0}_{t \times n} \\ D_1 \end{matrix} \right\| \right) Q^{-1} = U \left\| \begin{matrix} M_1 & \mathbf{0} \\ M_2 & \mathbf{0} \\ T_3 & T_4 \end{matrix} \right\| Q^{-1},$$

where $\left\| T_3 \quad T_4 \right\|$ is an arbitrary matrix in $M_{(n-t) \times n}(R)$. This set contains a matrix

$$U \left\| \begin{matrix} M_1 & \mathbf{0} \\ M_2 & \mathbf{0} \\ \mathbf{0} & I_{n-t} \end{matrix} \right\| Q^{-1} = F,$$

where $I_{n-t}$ is identity matrix of order $n-t$. All elements of the set $C + Ann_r(B)$ are roots of $BX = A$. Hence, $F$ is a root of $BX = A$. Note that each element of the set $C + Ann_r(B)$ can be written as

$$U \left\| \begin{matrix} M_1 & \mathbf{0} \\ M_2 & \mathbf{0} \\ T_3 & T_4 \end{matrix} \right\| Q^{-1} = \left( U \left\| \begin{matrix} M_1 & \mathbf{0} \\ M_2 & \mathbf{0} \\ \mathbf{0} & I_{n-t} \end{matrix} \right\| Q^{-1} \right) \left( Q \left\| \begin{matrix} I_t & \mathbf{0} \\ T_3 & T_4 \end{matrix} \right\| Q^{-1} \right) = FM.$$

That is, the matrix $F$ is the left divisor of all elements of the set $C + Ann_r(B)$. Taking into account that $F \in (C + Ann_r(B))$ we come to the conclusion that $F$ is the left g.c.d. of all elements of this set.

Also the following equation

$$N = U \left\| \begin{matrix} M_1 & \mathbf{0} \\ M_2 & \mathbf{0} \\ \mathbf{0} & \mathbf{0} \end{matrix} \right\| Q^{-1} = \left( U \left\| \begin{matrix} I_t & \mathbf{0} \\ \mathbf{0} & \mathbf{0} \end{matrix} \right\| U^{-1} \right) \left( U \left\| \begin{matrix} M_1 & \mathbf{0} \\ M_2 & \mathbf{0} \\ T_3 & T_4 \end{matrix} \right\| Q^{-1} \right)$$

holds. That is, the matrix $N$ is the left multiple of all elements of the set $C + Ann_r(B)$. Taking into account that $N \in (C + Ann_r(B))$ we come to the conclusion that $V$ is the left l.c.m. of all elements of this set. The theorem is proved. □

**Corollary 4.4.** *Left g.c.d. and left l.c.m. of solutions of the equation $BX = A$ are*

$$U \left\| \begin{matrix} M_1 & \mathbf{0} \\ M_2 & \mathbf{0} \\ \mathbf{0} & I_{n-t} \end{matrix} \right\| Q^{-1}, \;\; U \left\| \begin{matrix} M_1 & \mathbf{0} \\ M_2 & \mathbf{0} \\ \mathbf{0} & \mathbf{0} \end{matrix} \right\| Q^{-1},$$

*respectively.* □

**Corollary 4.5.** *The set*

$$\left( U \left\| \begin{matrix} M_1 & \mathbf{0} \\ M_2 & \mathbf{0} \\ \mathbf{0} & I_{n-t} \end{matrix} \right\| Q^{-1} \right) \left( Q \left\| \begin{matrix} I_t & \mathbf{0} \\ T_3 & T_4 \end{matrix} \right\| Q^{-1} \right),$$

where $\|T_3 \quad T_4\|$ *is an arbitrary matrix in* $M_{(n-t)\times n}(R)$, *is the set of all solutions of the equation* $BX = A$. $\qquad\qquad\Box$

**Example 4.1.** Consider the matrix linear equation over $\mathbb{Z}$: $BX = A$, where

$$B = \begin{Vmatrix} 0 & 1 & 0 & 0 & 0 & 0 & 0 \\ 0 & -2 & 2 & 0 & 0 & 0 & 0 \\ 1 & 0 & 0 & 0 & 0 & 0 & 0 \\ 0 & 0 & 0 & 0 & 0 & 0 & 0 \\ -2 & 0 & -12 & 0 & 12 & 0 & 0 \\ 0 & 0 & 0 & 0 & 0 & 0 & 0 \\ -2 & 0 & -4 & 4 & 0 & 0 & 0 \end{Vmatrix}, \quad A = \operatorname{diag}(1,2,6,0,0,0,0).$$

The matrix $A$ coincides with its Smith form E. Hence, $P, Q = I$. We write the matrix $B$ in the form: $B = P_B^{-1}\Phi$, where

$$P_B^{-1} = \begin{Vmatrix} 0 & 1 & 0 & 0 & 0 & 0 & 0 \\ 0 & -2 & 1 & 0 & 0 & 0 & 0 \\ 1 & 0 & 0 & 0 & 0 & 0 & 0 \\ 0 & 0 & 0 & 0 & 0 & 0 & 0 \\ -2 & 0 & -6 & 0 & 1 & 0 & 0 \\ 0 & 0 & 0 & 0 & 0 & 0 & 0 \\ -2 & 0 & -2 & 1 & 0 & 0 & 0 \end{Vmatrix}, \quad \Phi = \operatorname{diag}(1,1,2,4,12,0,0).$$

Since $\Phi|\mathrm{E}$, then $\mathbf{L}(\mathrm{E}, \Phi) \neq \varnothing$ and consists of all invertible matrices of the form

$$\begin{Vmatrix} * & * & * & * & * & * & * \\ * & * & * & * & * & * & * \\ 2h_{31} & * & * & * & * & * & * \\ \hline 4h_{41} & 2h_{42} & 2h_{43} & * & * & * & * \\ 12h_{51} & 6h_{52} & 2h_{53} & * & * & * & * \\ \hline 0 & 0 & 0 & * & * & * & * \\ 0 & 0 & 0 & * & * & * & * \end{Vmatrix}.$$

From $P = I$, it follows that $P_B = LP = L$. So

$$P_B = \begin{Vmatrix} 0 & 0 & 1 & 0 & 0 & 0 & 0 \\ 1 & 0 & 0 & 0 & 0 & 0 & 0 \\ 2 & 1 & 0 & 0 & 0 & 0 & 0 \\ \hline 4 & 2 & 2 & 0 & 0 & 0 & 1 \\ 12 & 6 & 2 & 0 & 1 & 0 & 0 \\ \hline 0 & 0 & 0 & 0 & 0 & 1 & 0 \\ 0 & 0 & 0 & 1 & 0 & 0 & 0 \end{Vmatrix} = L.$$

It is obviously that $L \in \mathbf{L}(\mathrm{E}, \Phi)$. Therefore, the equation $BX = A$ has a solution. It is easy to see that

$$LE = \Phi \left\| \begin{array}{ccc|cccc} 0 & 0 & 6 & 0 & 0 & 0 & 0 \\ 1 & 0 & 0 & 0 & 0 & 0 & 0 \\ 1 & 1 & 0 & 0 & 0 & 0 & 0 \\ \hline 1 & 1 & 3 & 0 & 0 & 0 & 0 \\ 1 & 1 & 1 & 0 & 0 & 0 & 0 \\ \hline * & * & * & * & * & * & * \\ * & * & * & * & * & * & * \end{array} \right\| = \Phi S.$$

Consequently, left g.c.d. and l.c.m. of the roots of the equation $BX = A$ are:

$$\left\| \begin{array}{ccc|cccc} 0 & 0 & 6 & 0 & 0 & 0 & 0 \\ 1 & 0 & 0 & 0 & 0 & 0 & 0 \\ 1 & 1 & 0 & 0 & 0 & 0 & 0 \\ \hline 1 & 1 & 3 & 0 & 0 & 0 & 0 \\ 1 & 1 & 1 & 0 & 0 & 0 & 0 \\ \hline 0 & 0 & 0 & 0 & 0 & 1 & 0 \\ 0 & 0 & 0 & 0 & 0 & 0 & 1 \end{array} \right\|, \quad \left\| \begin{array}{ccc|cccc} 0 & 0 & 6 & 0 & 0 & 0 & 0 \\ 1 & 0 & 0 & 0 & 0 & 0 & 0 \\ 1 & 1 & 0 & 0 & 0 & 0 & 0 \\ \hline 1 & 1 & 3 & 0 & 0 & 0 & 0 \\ 1 & 1 & 1 & 0 & 0 & 0 & 0 \\ \hline 0 & 0 & 0 & 0 & 0 & 0 & 0 \\ 0 & 0 & 0 & 0 & 0 & 0 & 0 \end{array} \right\|,$$

respectively. Moreover, this matrices are roots of $BX = A$. $\diamond$

## 4.4. Properties of Zelisko group and generating set

In this subsection the relationships between a general linear group, Zelisko group, and the generating set are investigated.

Let E and $\Phi$ be $d$-matrices over an elementary divisor ring $R$, moreover $\Phi | \mathrm{E}$.

Denote by $S_i$ a matrix obtained from $n \times n$ matrix $S$ by crossing out its first $i$ rows and first $i$ columns and $\overline{S}_i$ a matrix obtained by crossing out its first $i$ rows and the last $n - i$ columns.

**Lemma 4.1.** *Let* $\mathrm{E}, \Phi$ *be nonsingular matrices and* $F$ *be the lower unitriangular matrix, written in the form*

$$F = HS, \tag{4.6}$$

*where* $H \in \mathbf{G}_\Phi$, $S \in \mathbf{G}_\mathrm{E}$. *Then there are lower unitriangular matrices* $H_1 \in \mathbf{G}_\Phi$, $S_1 \in \mathbf{G}_\mathrm{E}$ *such that* $F = H_1 S_1$.

**Proof.** In view of (4.6), there is $H \in \mathbf{G}_\Phi$ such that $HS$ is a lower unitriangular matrix. According to Theorem 2.8,

$$\left( \frac{\varphi_{i+1}}{\varphi_i}, \det S_i \right) = 1, \quad i = 1, ..., n - 1.$$

By Property 2.4, $\dfrac{\varepsilon_{i+1}}{\varepsilon_i}\Big|\left\langle \overline{S}_i\right\rangle_1$. Since matrix $\|\overline{S}_iS_i\|$ is primitive, we get

$$\left(\frac{\varepsilon_{i+1}}{\varepsilon_i}, \det S_i\right) = 1, \quad i = 1, ..., n-1.$$

So

$$\left(\frac{\varepsilon_{i+1}\varphi_{i+1}}{\varepsilon_i\varphi_i}, \det S_i\right) = 1, \quad i = 1, ..., n-1. \tag{4.7}$$

Consider $d$-matrix $\Gamma = \text{diag}(\varepsilon_1\varphi_1, ..., \varepsilon_n\varphi_n)$. In view of (4.7) and using Theorem (4.7) we get that there exists $K \in \mathbf{G}_\Gamma$ such that $KS = S_1$ is a lower unitriangular matrix. Since $K \in \mathbf{G}_\text{E}\bigcap\mathbf{G}_\Phi$, then

$$F = HS = (HK^{-1})(KS) = H_1S_1,$$

where $H_1 \in \mathbf{G}_\Phi$, $S_1 \in \mathbf{G}_\text{E}$. The matrices $F$ and $S_1$ are lower unitriangular. So $H_1$ is a lower unitriangular matrix. □

**Lemma 4.2.** *Let* $(a_1, ..., a_n) = 1$ *and* $\left(\dfrac{\varepsilon_n}{\varepsilon_1}, a_n\right) = 1$. *Then the group* $\mathbf{G}_\text{E}$ *contains an invertible matrix having* $\|a_1 ... a_n\|^T$ *as the last column.*

**Proof.** First we treat the case $n = 2$. Since $(a_1, a_2) = 1$ and

$$\left(\frac{\varepsilon_2}{\varepsilon_1}, a_2\right) = 1,$$

then

$$\left(\frac{\varepsilon_2}{\varepsilon_1}a_1, a_2\right) = 1.$$

There are $u, v$ such that

$$\frac{\varepsilon_2}{\varepsilon_1}a_1u + a_2v = 1.$$

So the matrix

$$\left\|\begin{matrix} v & a_1 \\ -\dfrac{\varepsilon_2}{\varepsilon_1}u & a_2 \end{matrix}\right\|$$

is desired.

Suppose this statement is true for all matrices of an order less than $n$. Let $(a_2, ..., a_n) = \delta$. Then

$$\left(\frac{a_2}{\delta}, ..., \frac{a_n}{\delta}\right) = 1.$$

Since $\dfrac{\varepsilon_n}{\varepsilon_2}\Big|\dfrac{\varepsilon_n}{\varepsilon_1}$, i $\dfrac{a_n}{\delta}\Big|a_n$, then

$$\left(\frac{\varepsilon_n}{\varepsilon_2}, \frac{a_n}{\delta}\right) = 1.$$

By induction the group $\mathbf{G}_{\text{diag}(\varepsilon_2, ..., \varepsilon_n)}$ contains the matrix of the form

$$\left\|\begin{matrix} u_{11} & ... & u_{1.n-2} & \dfrac{a_2}{\delta} \\ ... & ... & ... & ... \\ u_{n-1.1} & ... & u_{n-1.n-2} & \dfrac{a_n}{\delta} \end{matrix}\right\| = U_1,$$

**135**

moreover $\det D = 1$. Hence,

$$\det \left\| \begin{matrix} u_{11} & ... & u_{1.n-2} & a_2 \\ ... & ... & ... & ... \\ u_{n-1.1} & ... & u_{n-1.n-2} & a_n \end{matrix} \right\| = \delta.$$

Consider the matrix

$$\left\| \begin{matrix} x_1 & 0 & ... & 0 & a_1 \\ \dfrac{\varepsilon_2}{\varepsilon_1}x_2 & u_{11} & ... & u_{1.n-2} & a_2 \\ ... & ... & ... & ... & ... \\ \dfrac{\varepsilon_n}{\varepsilon_1}x_n & u_{n-1.1} & ... & u_{n-1.n-2} & a_n \end{matrix} \right\| = V,$$

where $x_1, ..., x_n$ are parameters. Then

$$\det V = x_1\delta - \frac{\varepsilon_2}{\varepsilon_1}x_2 a_1 \Delta_1 + ... + (-1)^{n+1}\frac{\varepsilon_n}{\varepsilon_1}x_n a_1 \Delta_{n-1},$$

where

$$\Delta_i = \det \left\| \begin{matrix} u_{11} & ... & u_{1.n-2} \\ ... & ... & ... \\ u_{i-1.1} & ... & u_{i-1.n-2} \\ u_{i+1.1} & ... & u_{i+1.n-2} \\ ... & ... & ... \\ u_{n-1.1} & ... & u_{n-1.n-2} \end{matrix} \right\|.$$

The matrix

$$D = \left\| \begin{matrix} u_{11} & ... & u_{1.n-2} \\ ... & ... & ... \\ u_{n-1.1} & ... & u_{n-1.n-2} \end{matrix} \right\|$$

consists of the first $n-2$ columns of the invertible matrix $U_1$. Consequently,

$$\langle D \rangle = (\Delta_1, ..., \Delta_{n-1}) = 1.$$

Since $\delta | a_n$ and

$$\left( \frac{\varepsilon_n}{\varepsilon_1}, a_n \right) = 1,$$

we get

$$\left( \frac{\varepsilon_n}{\varepsilon_1}, \delta \right) = 1.$$

It follows that

$$\left( \frac{\varepsilon_i}{\varepsilon_1}, \delta \right) = 1, \ \ i = 2, ..., n.$$

Then

$$\left(\delta, \frac{\varepsilon_2}{\varepsilon_1}a_1\Delta_1, ..., \frac{\varepsilon_n}{\varepsilon_1}a_1\Delta_{n-1}\right) = \left(\left(\delta, \frac{\varepsilon_2}{\varepsilon_1}a_1\Delta_1\right), ..., \left(\delta, \frac{\varepsilon_n}{\varepsilon_1}a_1\Delta_{n-1}\right)\right) =$$
$$= (\delta, a_1\Delta_1, ..., a_1\Delta_{n-1}) = (\delta, a_1(\Delta_1, ..., \Delta_{n-1})) = (\delta, a_1) = 1.$$

Therefore, there are $x_1 = v_1, ..., x_n = v_n$ such that $\det V = 1$. $\qquad\square$

**Lemma 4.3.** *Let* $E, \Phi$ *be nonsingular matrices, moreover* $\Phi | E$, *and*

$$\frac{\varphi_i}{(\varphi_i, \varepsilon_j)} = \left(\frac{\varphi_i}{\varphi_j}, \frac{\varepsilon_i}{\varepsilon_j}\right),$$

$i = 2, ..., n$, $j = 1, ..., n-1$, $i > j$. *If* $L$ *is a lower unitriangular matrix, then there are lower unitriangular matrices* $H \in \mathbf{G}_\Phi$, $K \in \mathbf{G}_E$ *such that* $L = HK$.

**Proof.** Let $n = 2$. Consider the matrix

$$L = \left\| \begin{matrix} 1 & 0 \\ f_{21}l & 1 \end{matrix} \right\|,$$

where

$$f_{21} = \frac{\varphi_2}{(\varphi_2, \varepsilon_1)} = \left(\frac{\varphi_2}{\varphi_1}, \frac{\varepsilon_2}{\varepsilon_1}\right).$$

There exist $u, v$ such that

$$f_{21} = \frac{\varphi_2}{\varphi_1}u + \frac{\varepsilon_2}{\varepsilon_1}v.$$

Then

$$L = \left\| \begin{matrix} 1 & 0 \\ f_{21}l & 1 \end{matrix} \right\| = \left\| \begin{matrix} 1 & 0 \\ \frac{\varphi_2}{\varphi_1}ul & 1 \end{matrix} \right\| \left\| \begin{matrix} 1 & 0 \\ \frac{\varepsilon_2}{\varepsilon_1}vl & 1 \end{matrix} \right\| = HK.$$

Therefore, our assertion is correct for the second order matrices.

Suppose the correctness of this assertion for all matrices of order less than $n$. Consider a lower unitriangular matrix

$$L = \left\| \begin{matrix} 1 & 0 & & 0 & 0 \\ f_{21}l_{21} & 1 & & & 0 \\ & & \ddots & & \vdots \\ f_{n-1.1}l_{n-1.1} & f_{n-1.2}l_{n-1.2} & & 1 & 0 \\ f_{n1}l_{n1} & f_{n2}l_{n2} & ... & f_{n.n-1}l_{n.n-1} & 1 \end{matrix} \right\| = \left\| \begin{matrix} L_{11} & \mathbf{0} \\ L_{21} & 1 \end{matrix} \right\|,$$

where $f_{ij} = \frac{\varphi_i}{(\varphi_i, \varepsilon_j)}$. Denote by

$$\Phi_n = \operatorname{diag}(\varphi_1, ..., \varphi_{n-1}), E_n = \operatorname{diag}(\varepsilon_1, ..., \varepsilon_{n-1}),$$
$$\Phi_1 = \operatorname{diag}(\varphi_2, ..., \varphi_n), E_1 = \operatorname{diag}(\varepsilon_2, ..., \varepsilon_n).$$

Since $L_{11} \in \mathbf{L}(\mathrm{E}_n, \Phi_n)$ there are $H_1 \in \mathbf{G}_{\Phi_n}$ and $K_1 \in \mathbf{G}_{\mathrm{E}_n}$ such that $L_{11} = H_1 K_1$. Taking into account that

$$\left\| \begin{matrix} H_1^{-1} & \mathbf{0} \\ \mathbf{0} & 1 \end{matrix} \right\| \in \mathbf{G}_\Phi, \quad \left\| \begin{matrix} K_1^{-1} & \mathbf{0} \\ \mathbf{0} & 1 \end{matrix} \right\| \in \mathbf{G}_\mathrm{E},$$

and using Properties 4.1, 4.2, we get

$$\left\| \begin{matrix} H_1^{-1} & \mathbf{0} \\ \mathbf{0} & 1 \end{matrix} \right\| L \left\| \begin{matrix} K_1^{-1} & \mathbf{0} \\ \mathbf{0} & 1 \end{matrix} \right\| = \left\| \begin{matrix} 1 & 0 & 0 & 0 \\ 0 & 1 & 0 & 0 \\ & & \ddots & \\ 0 & 0 & 1 & 0 \\ f_{n1}l'_{n1} & f_{n2}l'_{n2} & \dots & f_{n.n-1}l'_{n.n-1} & 1 \end{matrix} \right\| =$$

$$= \left\| \begin{matrix} 1 & \mathbf{0} \\ S_{21} & S_{22} \end{matrix} \right\| = S \in \mathbf{L}(\mathrm{E}, \Phi).$$

Since $S_{22} \in \mathbf{L}(\mathrm{E}_1, \Phi_1)$ by induction there exist $H_2 \in \mathbf{G}_{\Phi_1}$ and $K_2 \in \mathbf{G}_{\mathrm{E}_1}$ such that $S_{22} = H_2 K_2$. Then

$$\left\| \begin{matrix} 1 & \mathbf{0} \\ \mathbf{0} & H_2^{-1} \end{matrix} \right\| S \left\| \begin{matrix} 1 & \mathbf{0} \\ \mathbf{0} & K_2^{-1} \end{matrix} \right\| = \left\| \begin{matrix} 1 & 0 & 0 & 0 \\ 0 & 1 & 0 & 0 \\ & & \ddots & \\ 0 & 0 & 1 & 0 \\ f_{n1}a & 0 & \dots & 0 & 1 \end{matrix} \right\|,$$

where $a \in R$. Since

$$f_{n1} = \left( \frac{\varphi_n}{\varphi_1}, \frac{\varepsilon_n}{\varepsilon_1} \right),$$

there are $u_n, v_n$ such that

$$f_{n1} = \frac{\varphi_n}{\varphi_1} u_n + \frac{\varepsilon_n}{\varepsilon_1} v_n.$$

Then

$$\left\| \begin{matrix} 1 & 0 & 0 & 0 \\ 0 & 1 & & 0 \\ & & \ddots & \\ 0 & 0 & 1 & 0 \\ f_{n1}a & 0 & 0 & 1 \end{matrix} \right\| = \underbrace{\left\| \begin{matrix} 1 & 0 & 0 & 0 \\ 0 & 1 & & 0 \\ & & \ddots & \\ 0 & 0 & 1 & 0 \\ \frac{\varphi_n}{\varphi_1}u_n a & 0 & 0 & 1 \end{matrix} \right\|}_{M} \underbrace{\left\| \begin{matrix} 1 & 0 & 0 & 0 \\ 0 & 1 & & 0 \\ & & \ddots & \\ 0 & 0 & 1 & 0 \\ \frac{\varepsilon_n}{\varepsilon_1}v_n a & 0 & 0 & 1 \end{matrix} \right\|}_{N}.$$

So the equality

$$\left\| \begin{matrix} 1 & \mathbf{0} \\ \mathbf{0} & H_2^{-1} \end{matrix} \right\| \left\| \begin{matrix} H_1^{-1} & \mathbf{0} \\ \mathbf{0} & 1 \end{matrix} \right\| L \left\| \begin{matrix} K_1^{-1} & \mathbf{0} \\ \mathbf{0} & 1 \end{matrix} \right\| \left\| \begin{matrix} 1 & \mathbf{0} \\ \mathbf{0} & K_2^{-1} \end{matrix} \right\| = MN$$

is fulfilled. That is

$$L = \left( \begin{Vmatrix} H_1 & \mathbf{0} \\ \mathbf{0} & 1 \end{Vmatrix} \begin{Vmatrix} 1 & \mathbf{0} \\ \mathbf{0} & H_2 \end{Vmatrix} M \right) \left( N \begin{Vmatrix} 1 & \mathbf{0} \\ \mathbf{0} & K_2 \end{Vmatrix} \begin{Vmatrix} K_1 & \mathbf{0} \\ \mathbf{0} & 1 \end{Vmatrix} \right) = HK,$$

where $H \in \mathbf{G}_\Phi$, $K \in \mathbf{G}_\mathrm{E}$. $\qquad\qquad\qquad\qquad\qquad\qquad\Box$

**Property 4.5.** *Let $R$ be a Bezout ring of stable rank* 1.5 *and* E *be a non-singular matrix, moreover* $\Phi|\mathrm{E}$. *In order that*

$$\mathbf{L}(\mathrm{E}, \Phi) = \mathbf{G}_\mathrm{E}\mathbf{G}_\Phi,$$

*necessary and sufficient*

$$\frac{\varphi_i}{(\varphi_i, \varepsilon_j)} = \left( \frac{\varphi_i}{\varphi_j}, \frac{\varepsilon_i}{\varepsilon_j} \right)$$

*for all $i = 2, ..., n$, $j = 1, ..., n - 1$, $i > j$.*

**Proof. Necessity.** Let $n = 2$. The set $\mathbf{L}(\mathrm{E}, \Phi)$ contains the matrix

$$\begin{Vmatrix} 1 & 0 \\ f_{21} & 1 \end{Vmatrix}, \quad f_{21} = \frac{\varphi_2}{(\varphi_2, \varepsilon_1)}.$$

Then $F = HK$, where $H \in \mathbf{G}_\Phi$, $K \in \mathbf{G}_\mathrm{E}$. Moreover, by Lemma 4.1, the matrices $H, K$ can be considered as lower unitriangular. That is

$$\begin{Vmatrix} 1 & 0 \\ f_{21} & 1 \end{Vmatrix} = \begin{Vmatrix} 1 & 0 \\ \frac{\varphi_2}{\varphi_1}h & 1 \end{Vmatrix} \begin{Vmatrix} 1 & 0 \\ \frac{\varepsilon_2}{\varepsilon_1}k & 1 \end{Vmatrix}.$$

Hence,

$$f_{21} = \frac{\varphi_2}{\varphi_1}h + \frac{\varepsilon_2}{\varepsilon_1}k.$$

Since $f_{21}\left|\dfrac{\varphi_2}{\varphi_1}\right.$ and $f_{21}\left|\dfrac{\varepsilon_2}{\varepsilon_1}\right.$, we get

$$f_{21} = \left( \frac{\varphi_2}{\varphi_1}, \frac{\varepsilon_2}{\varepsilon_1} \right).$$

Suppose that statement is correct for all matrices of order less than $n$. The equality

$$\mathbf{L}(\mathrm{E}_n, \Phi_n) = \mathbf{G}_{\Phi_n}\mathbf{G}_{\mathrm{E}_n}$$

implies that

$$\frac{\varphi_i}{(\varphi_i, \varepsilon_j)} = \left( \frac{\varphi_i}{\varphi_j}, \frac{\varepsilon_i}{\varepsilon_j} \right),$$

$i = 2, ..., n - 1$, $j = 1, ..., n - 2$, $i > j$. Also from

$$\mathbf{L}(\mathrm{E}_1, \Phi_1) = \mathbf{G}_{\Phi_1}\mathbf{G}_{\mathrm{E}_1}$$

we have

$$\frac{\varphi_i}{(\varphi_i, \varepsilon_j)} = \left( \frac{\varphi_i}{\varphi_j}, \frac{\varepsilon_i}{\varepsilon_j} \right),$$

$i = 3, ..., n$, $j = 2, ..., n$, $i > j$. To complete the proof you need to show that

$$f_{n1} = \frac{\varphi_n}{(\varphi_n, \varepsilon_1)} = \left(\frac{\varphi_n}{\varphi_1}, \frac{\varepsilon_n}{\varepsilon_1}\right).$$

The set $\mathbf{L}(\mathrm{E}, \Phi)$ contains the matrix

$$L = \left\| \begin{array}{ccccc} 1 & 0 & & 0 & 0 \\ 0 & 1 & & & 0 \\ & & \ddots & & \\ 0 & 0 & & 1 & 0 \\ f_{n1} & 0 & & 0 & 1 \end{array} \right\|.$$

By assumption of Theorem and Lemma 4.1 the matrix $L$ can be written in the form $L = HS$, where $H, S$ are lower unitriangular matrices from $\mathbf{G}_\Phi$, $\mathbf{G}_\mathrm{E}$, respectively. Then

$$H = LS^{-1} = L \left\| \begin{array}{cccc} 1 & 0 & & 0 \\ s_{21} & 1 & & \\ & & \ddots & \\ s_{n1} & s_{n2} & & 1 \end{array} \right\| =$$

$$= \left\| \begin{array}{ccccc} 1 & 0 & & 0 & 0 \\ s_{21} & 1 & & & 0 \\ & & \ddots & & \\ s_{n-1.1} & s_{n-1.2} & & 1 & 0 \\ f_{n1} + s_{n1} & s_{n2} & ... & s_{n.n-1} & 1 \end{array} \right\| \in \mathbf{G}_\Phi.$$

It follows that

$$f_{n1} + s_{n1} = \frac{\varphi_n}{\varphi_1} h_{n1}.$$

Since $s_{n1} = \frac{\varepsilon_n}{\varepsilon_1} s'_{n1}$, we get

$$\left(\frac{\varphi_n}{\varphi_1}, \frac{\varepsilon_n}{\varepsilon_1}\right) \Big| f_{n1}.$$

On the other hand, $f_{n1} | \frac{\varphi_n}{\varphi_1}$ and $f_{n1} \big| \frac{\varepsilon_n}{\varepsilon_1}$. That is

$$f_{n1} \Big| \left(\frac{\varphi_n}{\varphi_1}, \frac{\varepsilon_n}{\varepsilon_1}\right),$$

Therefore,

$$f_{n1} = \left(\frac{\varphi_n}{\varphi_1}, \frac{\varepsilon_n}{\varepsilon_1}\right).$$

**Sufficiency.** Let $L \in \mathbf{L}(\mathrm{E}, \Phi)$. By Theorem 2.14, the matrix $L$ can be written as $L = HVK$, where $H \in \mathbf{G}_\Phi$, and matrices $V$, $K$, are lower and upper unitriangular matrices, respectively. According to Corollary 2.2, the group of upper unitriangular matrices is a subgroup of any Zelisko group. So $K \in \mathbf{G}_E$. By virtue of Properties 4.1, 4.2, the equality

$$V = H^{-1}LK^{-1}$$

implies that $V \in \mathbf{L}(E, \Phi)$. By Lemma 4.3, the matrix $V$ can be written as $V = H_1K_1$, where $H_1 \in \mathbf{G}_\Phi$, $K_1 \in \mathbf{G}_E$. Hence,

$$L = (HH_1)(K_1K),$$

The theorem is proved. □

**Property 4.6.** *The equality*

$$\mathbf{L}(\Phi, \Phi) = \mathbf{G}_\Phi$$

*is fulfilled.*

**Proof**. The proof of this property follows directly from the structure of the matrices consisting the set $\mathbf{L}(\Phi, \Phi)$ (see Theorem 4.1), and the group $\mathbf{G}_\Phi$ (see Theorem 2.7). □

Let us analyze the elements of matrices from the set $\mathbf{L}(E, \Phi)$.

Set

$$d_{ij} = \left( \frac{(\varphi_i, \varepsilon_{i-1})}{\varphi_{i-1}} \left( \frac{(\varphi_{i-1}, \varepsilon_{i-2})}{\varphi_{i-2}} \left( ... \left( \frac{(\varphi_{j+2}, \varepsilon_{j+1})}{\varphi_{j+1}}, \frac{\varepsilon_{j+1}}{\varepsilon_j} \right), ... \right), \frac{\varepsilon_{i-2}}{\varepsilon_j} \right), \frac{\varepsilon_{i-1}}{\varepsilon_j} \right),$$

where $i > j + 1$.

**Lemma 4.4.** *The equality*

$$(\varphi_{j+1}, \varepsilon_j) \frac{(\varphi_{j+2}, \varepsilon_{j+1})}{\varphi_{j+1}} ... \frac{(\varphi_i, \varepsilon_{i-1})}{\varphi_{i-1}} = (\varphi_i, \varepsilon_j d_{ij}) = s_{ij}$$

*is fulfilled.*

**Proof**. We have

$$(\varphi_{j+1}, \varepsilon_j) \frac{(\varphi_{j+2}, \varepsilon_{j+1})}{\varphi_{j+1}} = \left( \varphi_{j+2}, \varepsilon_j \left( \frac{(\varphi_{j+2}, \varepsilon_{j+1})}{\varphi_{j+1}}, \frac{\varepsilon_{j+1}}{\varepsilon_j} \right) \right) =$$

$$= (\varphi_{j+2}, \varepsilon_j d_{j+2.j}) = s_{j+2.j}.$$

Suppose that

$$(\varphi_{j+1}, \varepsilon_j) \frac{(\varphi_{j+2}, \varepsilon_{j+1})}{\varphi_{j+1}} ... \frac{(\varphi_{i-1}, \varepsilon_{i-2})}{\varphi_{i-2}} = (\varphi_{i-1}, \varepsilon_j d_{i-1.j}) = s_{i-1.j}.$$

Then

$$s_{i-1.j} \frac{(\varphi_i, \varepsilon_{i-1})}{\varphi_1} =$$

$$= \left( \varphi_{i-1}, \varepsilon_j \left( \frac{(\varphi_{i-1}, \varepsilon_{i-2})}{\varphi_{i-2}} \left( ... \left( \frac{(\varphi_{i-1}, \varepsilon_{j+2})}{\varphi_{j+2}} \left( \frac{(\varphi_{j+2}, \varepsilon_{j+1})}{\varphi_{j+1}}, \frac{\varepsilon_{j+2}}{\varepsilon_j} \right), \right. \right. \right. \right.$$

$$\left. \left. \left. \left. \frac{\varepsilon_{i-2}}{\varepsilon_j} \right), ... \right), \frac{\varepsilon_{i-2}}{\varepsilon_j} \right) \right) \frac{(\varphi_i, \varepsilon_{i-1})}{\varphi_{i-1}} = (\varphi_i, \varepsilon_j d_{ij}) = s_{ij}.$$

Lemma is proved. □

Set $f_{ij} = \dfrac{\varphi_i}{(\varphi_i, \varepsilon_j)}$, $i > j$.

**Property 4.7.** *The equality*

$$f_{ij} = f_{j+1.j} f_{j+2.j+1} ... f_{i.i-1}(f_{ij}, d_{ij}),$$

*where $i = 3, 4, ..., n$; $j = 1, 2, ..., n-2$; $i > j+1$ is fulfilled.*

**Proof**. We have

$$\frac{f_{ij}}{f_{j+1.j} f_{j+2.j+1} ... f_{i.i-1}} = \frac{\varphi_i}{(\varphi_i, \varepsilon_j)} \frac{(\varphi_{j+1}, \varepsilon_j)}{\varphi_{j+1}} \frac{(\varphi_{j+2}, \varepsilon_{j+1})}{\varphi_{j+2}} ... \frac{(\varphi_i, \varepsilon_{i-1})}{\varphi_i} =$$

$$= \frac{1}{(\varphi_i, \varepsilon_j)}(\varphi_{j+1}, \varepsilon_j)\frac{(\varphi_{j+2}, \varepsilon_{j+1})}{\varphi_{j+1}} ... \frac{(\varphi_i, \varepsilon_{i-1})}{\varphi_{i-1}} = \frac{s_{ij}}{(\varphi_i, \varepsilon_j)} =$$

$$= \frac{(\varphi_i, \varepsilon_j d_{ij})}{(\varphi_i, \varepsilon_j)} = \left(\frac{\varphi_i}{(\varphi_i, \varepsilon_j)}, \frac{\varepsilon_j}{(\varphi_i, \varepsilon_j)}d_{ij}\right) = (f_{ij}, d_{ij}),$$

what needs to be proven. $\square$

**Lemma 4.5.** *The divisibility $(f_{ij}, d_{ij})|(f_{i+k.j-s}, d_{i+k.j-s})$ is fulfilled.*

**Proof**. We have

$$d_{i+k.j-s} = \left(\frac{(\varphi_{i+k}, \varepsilon_{i+k-1})}{\varphi_{i+k-1}}\left(...\left(\frac{(\varphi_{j-s+3}, \varepsilon_{j-s+2})}{\varphi_{j-s+2}}\left(\frac{(\varphi_{j-s+2}, \varepsilon_{j-s+1})}{\varphi_{j-s+1}},\right.\right.\right.\right.$$

$$\left.\left.\left.\frac{\varepsilon_{j-s+1}}{\varepsilon_{j-s}}\right), \frac{\varepsilon_{j-s+2}}{\varepsilon_{j-s}}\right), ...\right), \frac{\varepsilon_{i+k-1}}{\varepsilon_{j-s}} =$$

$$= \left(\frac{(\varphi_{i+k}, \varepsilon_{i+k-1})}{\varphi_{i+k-1}}\left(...\left(\frac{(\varphi_{j-s+3}, \varepsilon_{j-s+2})}{\varphi_{j-s+2}}\left(\frac{(\varphi_{j+2}, \varepsilon_{j+1})}{\varphi_{j+1}}d_{j+1.j-s},\right.\right.\right.\right.$$

$$\left.\left.\left.\frac{\varepsilon_{j+1}}{\varepsilon_{j-s}}\right), \frac{\varepsilon_{j-s+2}}{\varepsilon_{j-s}}\right), ...\right), \frac{\varepsilon_{i+k-1}}{\varepsilon_{j-s}}.$$

Since $\dfrac{\varepsilon_k}{\varepsilon_j}\bigg|\dfrac{\varepsilon_k}{\varepsilon_{j-s}}$, then

$$d_{ij}\left|\left(\frac{(\varphi_i, \varepsilon_{i-1})}{\varphi_{i-1}}\left(...\left(\frac{(\varphi_{j+2}, \varepsilon_{j+1})}{\varphi_{j+1}}d_{j+1.j-s}, \frac{\varepsilon_{j+1}}{\varepsilon_{j-s}}\right), ...\right), \frac{\varepsilon_{i-1}}{\varepsilon_{j-s}}\right).\right.$$

Therefore, $d_{ij}|d_{i+k.j-s}$. Noting that $f_{ij}|f_{i+k.j-s}$, we get $(f_{ij}, d_{ij})|(f_{i+k.j-s}, d_{i+k.j-s})$. Lemma is proved. $\square$

The comparison of the structure of matrix elements from $\mathbf{G}_\Phi$ (see Property 2.3, p. 63) and $\mathbf{L}(\mathrm{E}, \Phi)$ (see Property 4.7, p. 142), shows their apparent similarity.

So naturally, the question arises when there exists $d$-marix $\Delta$ such that $\mathbf{L}(\mathrm{E}, \Phi) = \mathbf{G}_\Delta$. To matrices of order 2, such $d$-matrix $\Delta$ always exists. If $\varphi_2 \neq 0$, then

$$\mathbf{L}(\mathrm{E}, \Phi) = \mathbf{G}_{\Delta_2},$$

where

$$\Delta_2 = \mathrm{diag}\left(1, \frac{\varphi_2}{(\varphi_2, \varepsilon_1)}\right).$$

If $\varphi_2 = 0$, then the set $\mathbf{L}(\mathrm{E}, \Phi)$ is a group of units of the ring of upper triangular matrices. In the case of matrices of an order greater than 2, this statement is not always correct.

**Example 4.2.** Consider integer $d$-matrices $\mathrm{E} = \mathrm{diag}(2, 6, 12)$ and $\Phi = \mathrm{diag}(1, 2, 6)$. By Theorem 4.1, the set $\mathbf{L}(\mathrm{E}, \Phi)$ consists of all invertible matrices of the form

$$\begin{Vmatrix} l_{11} & l_{12} & l_{13} \\ l_{21} & l_{22} & l_{23} \\ 3l_{31} & l_{32} & l_{33} \end{Vmatrix}.$$

Hence,

$$L = \begin{Vmatrix} 1 & 0 & 0 \\ 1 & 1 & 0 \\ 3 & 1 & 1 \end{Vmatrix} \in \mathbf{L}(\mathrm{E}, \Phi).$$

However, the matrix

$$L^2 = \begin{Vmatrix} 1 & 0 & 0 \\ 2 & 1 & 0 \\ 7 & 2 & 1 \end{Vmatrix}$$

does not belong to $\mathbf{L}(\mathrm{E}, \Phi)$. That is, in this case, the set $\mathbf{L}(\mathrm{E}, \Phi)$ does not multiplicatively closed. $\diamond$

**Property 4.8.** *In order that there exists $d$-marix $\Delta$ such that $\mathbf{L}(\mathrm{E}, \Phi) = \mathbf{G}_\Delta$, where*

$$\mathrm{E} = \mathrm{diag}(\varepsilon_1, ..., \varepsilon_k, 0, ..., 0), \quad \Phi = \mathrm{diag}(\varphi_1, ..., \varphi_t, 0, ..., 0),$$

$1 \le k \le n$, $1 \le t \le n$, $n > 2$, $\Phi | \mathrm{E}$, *necessary and sufficient*
*1) if $\det \Phi \ne 0$ then $(f_{n1}, d_{n1}) = 1$, i.e.,*

$$\frac{f_{n1}}{f_{21} f_{32} \cdots f_{n.n-1}} \in U(R);$$

*2) if $\det \Phi = 0$ then $k = t$, $(f_{k1}, d_{k1}) = 1$, i.e.,*

$$\frac{f_{k1}}{f_{21} f_{32} \cdots f_{k.k-1}} \in U(R).$$

**Proof. Necessity.** Let $\Phi$ be a nonsingular matrix. Consider a matrix of the order 3. The set $\mathbf{L}(\mathrm{E}, \Phi)$ contains a matrix

$$A_3 = \begin{Vmatrix} 1 & 0 & 0 \\ f_{21} & 1 & 0 \\ f_{31} & f_{32} & 1 \end{Vmatrix}.$$

Therefore,

$$A_3^2 = \left\| \begin{array}{ccc} 1 & 0 & 0 \\ 2f_{21} & 1 & 0 \\ 2f_{31} + f_{21}f_{32} & 2f_{32} & 1 \end{array} \right\| \in \mathbf{L}(\mathrm{E}, \Phi).$$

Hence,

$$2f_{31} + f_{21}f_{32} = f_{31}l_{31}.$$

According to Property 4.7

$$f_{31} = f_{21}f_{32}(f_{31}, d_{31}).$$

Then

$$f_{21}f_{32}(f_{31}, d_{31})(l_{31} - 2) = f_{21}f_{32}.$$

That is

$$(f_{31}, d_{31})(l_{31} - 2) = 1.$$

Consequently, $(f_{31}, d_{31}) = 1$.

Suppose that our assertion is correct for matrices of less than $n$. Let $\mathbf{L}(\mathrm{E}, \Phi)$ be a multiplicative group of the order $n$. This group contains a subgroup $1 \oplus \oplus \mathbf{L}(\mathrm{E}_1, \Phi_1)$, where

$$\mathrm{E}_1 = \mathrm{diag}(\varepsilon_2, ..., \varepsilon_n), \quad \Phi_1 = \mathrm{diag}(\varphi_2, ..., \varphi_n).$$

It follows that the set $\mathbf{L}(\mathrm{E}_1, \Phi_1)$ is also a multiplicative group. Since the set $\mathbf{L}(\mathrm{E}_1, \Phi_1)$ consists of invertible matrices of the form

$$\left\| \begin{array}{cccc} l_{22} & l_{23} & ... & l_{2n} \\ f_{32}l_{32} & l_{33} & ... & l_{3n} \\ ... & ... & ... & ... \\ f_{n2}l_{n2} & f_{n3}l_{n3} & ... & l_{nn} \end{array} \right\|,$$

by induction, we get

$$(f_{n2}, d_{n2}) = 1. \tag{4.8}$$

Consider the matrix

$$A_n = \left\| \begin{array}{cccccc} 1 & 0 & 0 & & & 0 \\ f_{21} & 1 & 0 & & & 0 \\ 0 & 0 & 1 & & & 0 \\ \vdots & \vdots & & \ddots & & \\ 0 & 0 & 0 & & 1 & 0 \\ f_{n1} & f_{n2} & 0 & ... & 0 & 1 \end{array} \right\| \in \mathbf{L}(\mathrm{E}, \Phi).$$

Since $A_n^2 = \|a_{ij}\|_1^n \in \mathbf{L}(\mathrm{E}, \Phi)$, then

$$a_{n1} = 2f_{n1} + f_{n2}f_{21} = f_{n1}l_{n1}.$$

In view of (4.8), and using Property 4.7, we get

$$f_{n2} = f_{32} ... f_{n.n-1},$$
$$f_{n1} = f_{21}f_{32} ... f_{n.n-1}(f_{n1}, d_{n1}).$$

Hence,

$$f_{21}f_{32}\ldots f_{n.n-1}(f_{n1}, d_{n1})(l_{n1} - 2) = f_{21}f_{32}\ldots f_{n.n-1}.$$

It follows that $(f_{n1}, d_{n1}) = 1$.

Let $E, \Phi$ be nonsingular matrices. By Theorem 4.1, the set $\mathbf{L}(E, \Phi)$ contains a zero $(n - t) \times k$ block. On the other hand, $\mathbf{L}(E, \Phi)$ is a group $\mathbf{G}_\Delta$, which according to Theorem 2.7 contains a zero $(n - s) \times s$ block. Hence, $k = t$. Theorem 2.7 also implies that

$$\mathbf{L}(\mathrm{diag}(\varepsilon_1, ..., \varepsilon_k), \ \mathrm{diag}(\varphi_1, ..., \varphi_k)) = \mathbf{G}_{\Delta_k}.$$

On the basis of just proved $(f_{k1}, d_{k1}) = 1$.

**Sufficiency.** Let $\Phi$ be a nonsingular matrix. By Lemma 4.5,

$$(f_{ij}, d_{ij})|(f_{n1}, d_{n1}) = 1, \ \ i = 3, 4, ..., n; \ \ j = 1, 2, ..., n - 2; \ \ i > j + 1.$$

That is

$$(f_{ij}, d_{ij}) = 1, \ \ i = 3, 4, ..., n; \ \ j = 1, 2, ..., n - 2; \ \ i > j + 1.$$

By virtue of Lemma 4.7,

$$f_{ij} = f_{j+1.j}f_{j+2.j+1}\ldots f_{i.i-1}, \ \ i = 3, 4, ..., n; \ \ j = 1, 2, ..., n - 2; \ \ i > j + 1.$$

So

$$\Delta = \mathrm{diag}(1, f_{21}, f_{21}f_{32}, ..., f_{21}f_{32}\ldots f_{n.n-1})$$

is a desired matrix.

If $\Phi, E$ are singular matrices it is easy to see

$$\Delta = \mathrm{diag}(1, f_{21}, f_{21}f_{32}, ..., f_{21}f_{32}\ldots f_{k.k-1}, 0, ..., 0)$$

is a desired matrix.

Having considered the Property 4.7, the conditions $(f_{n1}, d_{n1}) = 1$ and

$$\frac{f_{n1}}{f_{21}f_{32}\ldots f_{n.n-1}} \in U(R)$$

are equivalent. The theorem is proved. $\qquad\qquad\square$

**Property 4.9.** *The equality*

$$\mathbf{L}(E, \Phi) = \mathbf{G}_\Phi$$

*is fulfilled if and only if*
  *1) if $k = t = n$, then $(\varphi_n, \varepsilon_j) = \varphi_j, j = 1, ..., n - 1$;*
  *2) if $k < n, t = n$, then*
     *i) $\varphi_{k+1} = \varphi_{k+2} = ... = \varphi_n$,*
     *ii) $(\varphi_n, \varepsilon_j) = \varphi_j, j = 1, ..., k$ ;*
  *3) if $k, t < n$, then*
     *iii) $k = t$,*
     *iiii) $(\varphi_k, \varepsilon_j) = \varphi_j, j = 1, ..., k - 1$.*

**Proof.** Suppose that $k = t = n$. The equality of the sets $\mathbf{L}(\mathrm{E}, \Phi)$ and $\mathbf{G}_\Phi$ is equivalent to $L_1 = H_1$. Therefore

$$(\varphi_i, \varepsilon_j) = \varphi_j, \ \ i = 2, ..., n, \ \ j = 1, ..., n-1, \ \ i > j.$$

Specifically,

$$(\varphi_n, \varepsilon_j) = \varphi_j, \ \ j = 1, ..., n-1.$$

Conversely, if

$$(\varphi_n, \varepsilon_j) = \varphi_j, \ \ j = 1, ..., n-1,$$

we have

$$\left( \frac{\varphi_n}{\varphi_j}, \frac{\varepsilon_j}{\varphi_j} \right) = 1,$$

so that

$$\left( \frac{\varphi_i}{\varphi_j}, \frac{\varepsilon_j}{\varphi_j} \right) = 1, \ \ i = j+1, j+2, ..., n.$$

Hence,

$$(\varphi_i, \varepsilon_j) = \varphi_j \left( \frac{\varphi_i}{\varphi_j}, \frac{\varepsilon_j}{\varphi_j} \right) = \varphi_j, \ \ i = 2, ..., n, \ \ j = 1, ..., n-1, \ \ i > j.$$

Case 2). The equality of the sets $\mathbf{L}(\mathrm{E}, \Phi)$ and $\mathbf{G}_\Phi$ in this case is equivalent to

$$H_1 = \left\| \begin{matrix} L_1 & * \\ L_2 & * \end{matrix} \right\|.$$

Hence,

$$\frac{\varphi_i}{\varphi_j} = 1, \ \ i = k+2, k+3, ..., n, \ \ j = k+1, k+2, ..., n-1, \ \ i > j \quad (4.9)$$

and

$$(\varphi_i, \varepsilon_j) = \varphi_j, \ \ i = 2, ..., n, \ \ j = 1, ..., k, \ \ i > j. \quad (4.10)$$

Specifically,

$$\frac{\varphi_{k+2}}{\varphi_{k+1}} = \frac{\varphi_{k+3}}{\varphi_{k+2}} = ... = \frac{\varphi_n}{\varphi_{n-1}} = 1.$$

It means

$$\varphi_{k+1} = \varphi_{k+2} = ... = \varphi_n. \quad (4.11)$$

Having noticed

$$\frac{\varphi_p}{\varphi_q} = \frac{\varphi_p}{\varphi_{p-1}} \frac{\varphi_{p-1}}{\varphi_{p-2}} ... \frac{\varphi_{q+1}}{\varphi_q},$$

where $p > q$, we conclude that equalities (4.9) and (4.11) are equivalent. In the same manner as above we can see that $ii$), and (4.10) are equivalent.

Now consider Case 3). Thus we get

$$\left\| \begin{matrix} L_1 & * \\ L_2 & * \\ \mathbf{0} & * \end{matrix} \right\| = \left\| \begin{matrix} H_1 & * \\ \mathbf{0} & H_2 \end{matrix} \right\|.$$

The sizes of zero submatrices are $(n-t) \times k$ and $(n-t) \times t$ hence, $k = t$. It follows that the matrix $L_2$ is empty. Furthermore the analysis similar to above shows that $(\varphi_k, \varepsilon_j) = \varphi_j$, $j = 1, ..., k-1$. □

**Property 4.10.** *Let* E *be a nonsingular matrix and* $E = \Phi\Delta$. *If* $\mathbf{L}(E, \Phi) = \mathbf{G}_\Phi$, *then* $\Delta$ *is a d-matrix.*

**Proof.** Matrix $\Delta$ has the form

$$\Delta = \text{diag}\left(\frac{\varepsilon_1}{\varphi_1}, \frac{\varepsilon_2}{\varphi_2}, ..., \frac{\varepsilon_n}{\varphi_n}\right).$$

Consider the product

$$\frac{\varepsilon_{i+1}}{\varphi_{i+1}} \frac{\varphi_i}{\varepsilon_i} = \frac{\varepsilon_{i+1}\varphi_i}{\varphi_{i+1}\varepsilon_i} = \frac{\varepsilon_{i+1}\varphi_i}{(\varphi_{i+1}, \varepsilon_i)[\varphi_{i+1}, \varepsilon_i]} = \mu_{i+1.i},$$

$i = 1, ..., n-1$. Since $\mathbf{L}(E, \Phi) = \mathbf{G}_\Phi$, by Property 4.9, we get $(\varphi_{i+1}, \varepsilon_i) = \varphi_i$. So

$$\mu_{i+1.i} = \frac{\varepsilon_{i+1}\varphi_i}{\varphi_i[\varphi_{i+1}, \varepsilon_i]} = \frac{\varepsilon_{i+1}}{[\varphi_{i+1}, \varepsilon_i]}.$$

Since $\varphi_{i+1}|\varepsilon_{i+1}$ and $\varepsilon_i|\varepsilon_{i+1}$, then $[\varphi_{i+1}, \varepsilon_i]|\varepsilon_{i+1}$. That is $\mu_{i+1.i} \in R$ and

$$\frac{\varepsilon_{i+1}}{\varphi_{i+1}} = \frac{\varepsilon_i}{\varphi_i}\mu_{i+1.i}.$$

$i = 1, ..., n-1$. Consequently, $\Delta$ is a $d$-matrix. □

Note that $\Delta$ is a $d$-matrix, it does not follow that $\mathbf{L}(E, \Phi) = \mathbf{G}_\Phi$.

**Example 4.3.** Let $E = \text{diag}(a, a^3)$, $\Phi = \text{diag}(1, a)$, $a \neq 0$. Then $\Delta = \text{diag}(a, a^2)$. In this case $\mathbf{L}(E, \Phi) = \text{GL}_2(R)$, and the group $\mathbf{G}_\Phi$ consists of all invertible matrices of the form

$$\left\|\begin{matrix} h_{11} & h_{12} \\ ah_{21} & h_{22} \end{matrix}\right\|.$$

**Property 4.11.** *In order that*

$$\mathbf{L}(E, \Phi) = \text{GL}_n(R),$$

*it is necessary and sufficient that* $\varphi_n|\varepsilon_1$.

**Proof. Necessity.** By Theorem 4.1, the set $\mathbf{L}(E, \Phi)$ consists of all invertible matrices of the form (4.1). Since

$$\left\|\begin{matrix} 0 & ... & 0 & 1 \\ 0 & ... & 1 & 0 \\ ... & ... & ... & ... \\ 1 & ... & 0 & 0 \end{matrix}\right\| \in \text{GL}_n(R) = \mathbf{L}(E, \Phi),$$

then the matrices from the set $\mathbf{L}(\mathrm{E}, \Phi)$ do not contain a zero block in the lower left corner. This is equivalent to the fact that $\varphi_n \neq 0$, moreover

$$\frac{\varphi_n}{(\varphi_n, \varepsilon_1)} = 1.$$

That is $\varphi_n = (\varphi_n, \varepsilon_1)$. Hence, $\varphi_n | \varepsilon_1$.

**Sufficiency.** Since $\varphi_n | \varepsilon_1$, then $\varphi_n | \varepsilon_j$, $j = 1, ..., n$. It follows that $\varphi_i | \varepsilon_j$, $i = 1, ..., n$, $j = 1, ..., k$. Then

$$\frac{\varphi_i}{(\varphi_i, \varepsilon_j)} = 1.$$

Thus

$$\mathbf{L}(\mathrm{E}, \Phi) = \mathrm{GL}_n(R). \qquad \square$$

Let

$$\Gamma = \mathrm{diag}\,(\gamma_1, \ ..., \ \gamma_n), \ \ \Delta = \mathrm{diag}\,(\delta_1, \ ..., \ \delta_n),$$

be nonsingular $d$-matrices.

**Theorem 4.8.** *In order that* $\mathrm{GL}_n(R) = \mathbf{G}_\Gamma^T \mathbf{G}_\Delta$, *necessary and sufficient*

$$\left(\det \frac{1}{\delta_1}\Delta, \det \frac{1}{\gamma_1}\Gamma\right) = 1.$$

**Proof. Necessity.** Let

$$\left(\det \frac{1}{\delta_1}\Delta, \det \frac{1}{\gamma_1}\Gamma\right) = \sigma \neq 1$$

and $\dfrac{\gamma_r}{\gamma_1}$ be the first diagonal elements of the matrix $\dfrac{1}{\gamma_1}\Gamma$ such that

$$\left(\frac{\gamma_r}{\gamma_1}, \sigma\right) = \sigma_1 \neq 1.$$

Let $\dfrac{\delta_s}{\delta_1}$ is the first diagonal elements of the matrix $\dfrac{1}{\delta_1}\Delta$ such that

$$\left(\frac{\delta_s}{\delta_1}, \sigma_1\right) = \sigma_2 \neq 1.$$

Moreover

$$\left(\frac{\delta_{s-1}}{\delta_1}, \sigma_2\right) = \left(\frac{\gamma_{r-1}}{\gamma_1}, \sigma_2\right) = 1.$$

Since

$$\frac{\delta_s}{\delta_1} = \frac{\delta_{s-1}}{\delta_1}\frac{\delta_s}{\delta_{s-1}},$$

then, considering the previous equality, we get

$$\sigma_2 \left| \frac{\delta_s}{\delta_{s-1}}\right..$$

It follows that $\sigma_2$ is a divisor of all elements of the matrix

$$\left\|\begin{matrix} \dfrac{\delta_s}{\delta_1} & \cdots & \dfrac{\delta_s}{\delta_{s-1}} \\ \cdots & \cdots & \cdots \\ \dfrac{\delta_n}{\delta_1} & \cdots & \dfrac{\delta_n}{\delta_{s-1}} \end{matrix}\right\|.$$

Thus, every matrix from $\mathbf{G}_\Delta$ in the lower left corner contains $(n-s+1)\times(s-1)$ submatrix with elements which are multiple of $\sigma_2$.

Similarly, we show that every matrix from $\mathbf{G}_\Gamma^T$ in upper right corner contains $(r-1)\times(n-r+1)$ submatrix with elements which are multiple of $\sigma_2$. Whereas

$$\mathrm{GL}_n(R) = \mathbf{G}_\Gamma^T\mathbf{G}_\Delta,$$

there are $L \in \mathbf{G}_\Gamma^T$ and $H \in \mathbf{G}_\Delta$ such that

$$LH = \left\|\begin{matrix} 0 & & 1 \\ & \cdot^{\cdot^{\cdot}} & \\ 1 & & 0 \end{matrix}\right\| = T.$$

Then $L = TH^{-1}$. The matrix $H^{-1}$ has the same structure as the matrix $H$. Therefore, the matrix $L$ contains $(n-s+1)\times(s-1)$ submatrix in the left upper corner, and $(r-1)\times(n-r+1)$ submatrix in the right upper corner with elements which are multiple of $\sigma_2$.

If $n-s+1 \leqslant r-1$, then the matrix $L$ contains

$$(n-s+1)\times((s-1)+(n-r+1))$$

submatrix with elements which a multiple of $\sigma_2$. Since

$$(n-s+1)+(s-1)+(n-r+1) = (n+1)+(n-r) \geqslant n+1,$$

by Proposition 3.1, $\sigma_2|\det L$, which contradicts the invertibility of this matrix.

If $n-s+1 > r-1$, then the matrix $L$ contains

$$(r-1)\times((s-1)+(n-r+1))$$

submatrix with elements which are multiple of $\sigma_2$. Whereas

$$(r-1)+(s-1)+(n-r+1) = n+s-1 = (n+1)+(s-2) \geqslant n+1,$$

then and in this case $\sigma_2|\det L$ is contradiction. Consequently, $T \notin \mathbf{G}_\Gamma^T\mathbf{G}_\Delta$. Hence, $\mathbf{G}_\Gamma^T\mathbf{G}_\Delta \neq \mathrm{GL}_n(R)$ is contradiction.

**Sufficiency.** Let $A = \|a_{ij}\|_1^2 \in \mathrm{GL}_2(R)$ and

$$\left(a_{11}, \frac{\delta_2}{\delta_1} a_{12}\right) = \sigma.$$

From Theorem 6.3, it follows that there exists $H \in \mathbf{G}_\Delta$ such that

$$AH = \left\|\begin{matrix} \sigma & b_{12} \\ b_{21} & b_{22} \end{matrix}\right\|.$$

Since

$$\left(\frac{\delta_2}{\delta_1}, \frac{\gamma_2}{\gamma_1}\right) = 1, \quad \text{and} \quad \sigma \left| \frac{\delta_2}{\delta_1}, \right.$$

then

$$\left(\sigma, \frac{\gamma_2}{\gamma_1}\right) = 1.$$

So there is $L \in \mathbf{G}_\Gamma^T$ such that

$$\det(LAH) = 1$$

and

$$LAH = \left\|\begin{matrix} 1 & a \\ b & c \end{matrix}\right\| = \left\|\begin{matrix} 1 & 0 \\ b & 1 \end{matrix}\right\| \left\|\begin{matrix} 1 & a \\ 0 & 1 \end{matrix}\right\|.$$

Therefore,

$$A = \underbrace{\left(L^{-1} \left\|\begin{matrix} 1 & 0 \\ b & 1 \end{matrix}\right\|\right)}_{L_1} \underbrace{\left(\left\|\begin{matrix} 1 & a \\ 0 & 1 \end{matrix}\right\| H^{-1}\right)}_{H_1},$$

where $L_1 \in \mathbf{G}_\Gamma^T$, and $H_1 \in \mathbf{G}_\Delta$. That is, the Theorem is correct for the second order matrices.

Assume that our assumption is correct for matrices of the order $n - 1$. Let $A = \|a_{ij}\|_1^n \in \mathrm{GL}_n(R)$. Reasoning similarity as above, we can prove, that there are $L \in \mathbf{G}_\Gamma^T$ and $H \in \mathbf{G}_\Delta$ such that

$$LAH = \left\|\begin{matrix} 1 & \mathbf{0} \\ \mathbf{0} & A_{n-1} \end{matrix}\right\|.$$

By induction, an invertible matrix $A_{n-1}$ can be written in the form

$$A_{n-1} = L_{n-1} H_{n-1},$$

where

$$L_{n-1} \in \mathbf{G}_{\Gamma_1}^T, \quad H_{n-1} \in \mathbf{G}_{\Delta_1}, \quad \Gamma_1 = \mathrm{diag}(\gamma_2, ..., \gamma_n),$$
$$\Delta_1 = \mathrm{diag}(\delta_2, ..., \delta_n).$$

Hence

$$A = \left(L^{-1} \left\|\begin{matrix} 1 & \mathbf{0} \\ \mathbf{0} & L_{n-1} \end{matrix}\right\|\right)\left(\left\|\begin{matrix} 1 & \mathbf{0} \\ \mathbf{0} & H_{n-1} \end{matrix}\right\| H^{-1}\right).$$

Noting that

$$\left\|\begin{matrix} 1 & \mathbf{0} \\ \mathbf{0} & L_{n-1} \end{matrix}\right\| \in \mathbf{G}_\Gamma^T, \quad \left\|\begin{matrix} 1 & \mathbf{0} \\ \mathbf{0} & H_{n-1} \end{matrix}\right\| \in \mathbf{G}_\Delta,$$

we make sure correctness of our statement. The theorem is proved. $\qquad \square$

# Chapter **5**

# FACTORIZATION OF MATRICES

*This subsection shows the relationship between the properties of the set* $\mathbf{L}(E, \Phi)$ *and matrix divisors generated by this set.*

*Everywhere, unless specifically stated, $R$ is an elementary divisor ring.*

## 5.1. Specific cases of matrix factorizations

**Theorem 5.1.** *The matrix $A = P_A^{-1}EQ_A^{-1}$ has one up to associativity divisor with the Smith form $\Phi$ if and only if*

$$\mathbf{L}(E, \Phi) = \mathbf{G}_\Phi.$$

**Proof.** By Theorem 4.6, the matrix $A$ has one up to associativity divisor with the Smith form $\Phi$ if and only if $\mathbf{W}(E, \Phi) = \{I\}$. This is equivalent to $\mathbf{L}(E, \Phi) = \mathbf{G}_\Phi$. $\square$

Conditions for equality of sets $\mathbf{L}(E, \Phi)$, $\mathbf{G}_\Phi$ in terms of invariant factors are specified in Property 4.9.

Let matrix $A$ has Smith form E and $P \in \mathbf{P}_A$. Suppose that $\Phi$ be a $d$-matrix such that $E = \Phi\Delta$. Then

$$A = P^{-1}EQ^{-1} = (P^{-1}\Phi)(\Delta Q^{-1}) = (P^{-1}\Phi U^{-1})(U\Delta Q^{-1}),$$

where $U \in \mathrm{GL}_n(R)$. It follows that the set $\mathbf{P}_A^{-1}\Phi\mathrm{GL}_n(R)$ is the set of left divisors of the matrix $A$ with the Smith form $\Phi$. The question arises, $\mathbf{P}_A^{-1}\Phi\mathrm{GL}_n(R)$ is the set of all left divisors of the $A$ matrix with the Smith form $\Phi$ or not? The complete answer is given in the following theorem.

**Theorem 5.2.** *The set $\mathbf{P}_A^{-1}\Phi\mathrm{GL}_n(R)$ consists of all left divisors of the matrix $A = P_A^{-1}EQ_A^{-1}$ with the Smith form $\Phi$ if and only if*

$$\mathbf{L}(E, \Phi) = \mathbf{G}_\Phi\mathbf{G}_E.$$

**Proof. Necessity.** According to Corollary 4.4, the set $(\mathbf{L}(E, \Phi)P_A)^{-1}\Phi\mathrm{GL}_n(R)$ is the set of all left divisors of the matrix $A$ with the Smith form $\Phi$. Let

$$(\mathbf{L}(E, \Phi)P_A)^{-1}\Phi\mathrm{GL}_n(R) = \mathbf{P}_A^{-1}\Phi\mathrm{GL}_n(R).$$

**151**

# Chapter **5**

# FACTORIZATION OF MATRICES

*This subsection shows the relationship between the properties of the set* $\mathbf{L}(E, \Phi)$ *and matrix divisors generated by this set.*

*Everywhere, unless specifically stated, $R$ is an elementary divisor ring.*

## 5.1. Specific cases of matrix factorizations

**Theorem 5.1.** *The matrix $A = P_A^{-1}EQ_A^{-1}$ has one up to associativity divisor with the Smith form $\Phi$ if and only if*

$$\mathbf{L}(E, \Phi) = \mathbf{G}_\Phi.$$

**Proof.** By Theorem 4.6, the matrix $A$ has one up to associativity divisor with the Smith form $\Phi$ if and only if $\mathbf{W}(E, \Phi) = \{I\}$. This is equivalent to $\mathbf{L}(E, \Phi) = \mathbf{G}_\Phi$. $\square$

Conditions for equality of sets $\mathbf{L}(E, \Phi)$, $\mathbf{G}_\Phi$ in terms of invariant factors are specified in Property 4.9.

Let matrix $A$ has Smith form E and $P \in \mathbf{P}_A$. Suppose that $\Phi$ be a $d$-matrix such that $E = \Phi\Delta$. Then

$$A = P^{-1}EQ^{-1} = (P^{-1}\Phi)(\Delta Q^{-1}) = (P^{-1}\Phi U^{-1})(U\Delta Q^{-1}),$$

where $U \in \mathrm{GL}_n(R)$. It follows that the set $\mathbf{P}_A^{-1}\Phi\mathrm{GL}_n(R)$ is the set of left divisors of the matrix $A$ with the Smith form $\Phi$. The question arises, $\mathbf{P}_A^{-1}\Phi\mathrm{GL}_n(R)$ is the set of all left divisors of the $A$ matrix with the Smith form $\Phi$ or not? The complete answer is given in the following theorem.

**Theorem 5.2.** *The set $\mathbf{P}_A^{-1}\Phi\mathrm{GL}_n(R)$ consists of all left divisors of the matrix $A = P_A^{-1}EQ_A^{-1}$ with the Smith form $\Phi$ if and only if*

$$\mathbf{L}(E, \Phi) = \mathbf{G}_\Phi\mathbf{G}_E.$$

**Proof. Necessity.** According to Corollary 4.4, the set $(\mathbf{L}(E, \Phi)P_A)^{-1}\Phi\mathrm{GL}_n(R)$ is the set of all left divisors of the matrix $A$ with the Smith form $\Phi$. Let

$$(\mathbf{L}(E, \Phi)P_A)^{-1}\Phi\mathrm{GL}_n(R) = \mathbf{P}_A^{-1}\Phi\mathrm{GL}_n(R).$$

**151**

This is equivalent to that for each matrix $L \in \mathbf{L}(\mathrm{E}, \Phi)$ and $V \in \mathrm{GL}_n(R)$ there are $P \in \mathbf{P}_A$ and $U \in \mathrm{GL}_n(R)$ such that

$$(LP_A)^{-1}\Phi V = P^{-1}\Phi U.$$

Since $\mathbf{P}_A = \mathbf{G}_\mathrm{E} P_A$, there exists $K \in \mathbf{G}_\mathrm{E}$ such that $P = K P_A$. Thus,

$$(LP_A)^{-1}\Phi V = (KP_A)^{-1}\Phi U.$$

It follows that

$$(KL^{-1})\Phi = \Phi(UV^{-1}).$$

It means that $KL^{-1} = H \in \mathbf{G}_\Phi$, i.e., $L = H^{-1}K$. Taking into account that $H^{-1} \in \mathbf{G}_\Phi$, we get $L \in \mathbf{G}_\Phi \mathbf{G}_\mathrm{E}$. Consequently, $\mathbf{L}(\mathrm{E}, \Phi) \subseteq \mathbf{G}_\Phi \mathbf{G}_\mathrm{E}$. By Property 4.4, $\mathbf{G}_\Phi \mathbf{G}_\mathrm{E} \subseteq \mathbf{L}(\mathrm{E}, \Phi)$. So $\mathbf{G}_\Phi \mathbf{G}_\mathrm{E} = \mathbf{L}(\mathrm{E}, \Phi)$.

**Sufficiency.** We have

$$(\mathbf{L}(\mathrm{E}, \Phi)P_A)^{-1}\Phi\mathrm{GL}_n(R) = (\mathbf{G}_\Phi \mathbf{G}_\mathrm{E} P_A)^{-1}\Phi\mathrm{GL}_n(R) = (\mathbf{G}_\Phi(\mathbf{G}_\mathrm{E} P_A))^{-1} \times$$
$$\times \Phi\mathrm{GL}_n(R) = (\mathbf{G}_\Phi \mathbf{P}_A)^{-1}\Phi\mathrm{GL}_n(R) = \mathbf{P}_A^{-1}\mathbf{G}_\Phi\Phi\mathrm{GL}_n(R) = \mathbf{P}_A^{-1}\Phi\mathrm{GL}_n(R).$$

The theorem is proved. $\qquad\square$

Conditions for equality of sets $\mathbf{L}(\mathrm{E}, \Phi)$, $\mathbf{G}_\mathrm{E}\mathbf{G}_\Phi$ over a Bezout ring of stable range 1.5 are formulated in Property 4.5.

The Theorem 5.1 gives the conditions under which the matrix $A$ has only one, up to associativity, divisor with the Smith form $\Phi$. In the next theorem we consider another "extreme" case: all matrices with a given Smith form are divisors of the matrix $A$.

**Theorem 5.3.** *In order that each matrix with the Smith form $\Phi$ be a left divisor of the matrix $A = P_A^{-1}\mathrm{E}Q_A^{-1}$, necessary and sufficient $\Phi|\mathrm{E}$ and*

$$\mathbf{L}(\mathrm{E}, \Phi) = \mathrm{GL}_n(R).$$

**Proof. Necessity** of conditions $\Phi|\mathrm{E}$ is proved in Theorem 4.3.

Let $U$ be an invertible matrix. Then the matrix $(UP_A)^{-1}\Phi$ is the left divisor of the matrix $A$. The set of all left divisors of the matrix $A$ with the Smith form $\Phi$ has the form $(\mathbf{L}(\mathrm{E}, \Phi)P_A)^{-1}\Phi\mathrm{GL}_n(R)$. Therefore, there are $L \in \mathbf{L}(\mathrm{E}, \Phi)$, and $V \in \mathrm{GL}_n(R)$ such that

$$(UP_A)^{-1}\Phi = (LP_A)^{-1}\Phi V^{-1}.$$

It follows that $LU^{-1}\Phi = \Phi V^{-1}$. Consequently, $LU^{-1} = H \in \mathbf{G}_\Phi$. Hence, $U = H^{-1}L$. Therefore, each invertible matrix $U$ can be represented as the product of matrices from $\mathbf{G}_\Phi$ and $\mathbf{L}(\mathrm{E}, \Phi)$. That is,

$$\mathrm{GL}_n(R) \subseteq \mathbf{G}_\Phi\mathbf{L}(\mathrm{E}, \Phi).$$

The inverse inclusion is obvious. So

$$\mathrm{GL}_n(R) = \mathbf{G}_\Phi \mathbf{L}(\mathrm{E}, \Phi).$$

By Property 4.1,

$$\mathbf{G}_\Phi \mathbf{L}(\mathrm{E}, \Phi) = \mathbf{L}(\mathrm{E}, \Phi).$$

Therefore, $\mathbf{L}(\mathrm{E}, \Phi) = \mathrm{GL}_n(R)$.

**Sufficiency.** Let $B = P_B^{-1} \Phi Q_B^{-1}$ be an arbitrary matrix with the Smith form $\Phi$. Since

$$P_B P_A^{-1} \in \mathrm{GL}_n(R) = \mathbf{L}(\mathrm{E}, \Phi),$$

according to Theorem 4.2, the matrix $B$ is a left divisor of $A$. $\qquad\square$

Conditions for equality of sets $\mathbf{L}(\mathrm{E}, \Phi)$, $\mathrm{GL}_n(R)$ in terms of invariant factors are formulated in Properties 4.11.

According to Theorem 1.11, if the matrices $A, B$ are left divisors of each other, then they are right associates. Consider the case when they are both left and right divisors of each other.

**Theorem 5.4.** *Assume that $B$ is a left divisor of $A$, and the matrix $A$ is a right divisor of $B$. Then $A$ and $B$ are right and left associates.*

**Proof.** Let the matrices $A, B$ have Smith forms $\mathrm{E}, \Phi$, respectively. Since $A = BC$, by Theorem 4.2, $\Phi | \mathrm{E}$. The equality $B = DA$ is also fulfilled. It follows that $B^T = A^T D^T$. Therefore, $A^T$ is a left divisor of $B^T$. The operator of transposition does not change the determinants of all submatrices of the matrix $A$. Consequently, by virtue of Theorem 2.2, we get $A^T \sim A \sim \mathrm{E}$. Hence, $\mathrm{E} | \Phi$. It means that the corresponding invariant factors of $\Phi$ and $\mathrm{E}$ are associates. It follows that the matrices $A, B$ can be written as

$$A = P_A^{-1} \Phi Q_A^{-1}, \; B = P_B^{-1} \Phi Q_B^{-1}.$$

Since $A = BC$, by Theorem 4.2, $P_B = L P_A^{-1}$, where $L \in \mathbf{L}(\Phi, \Phi)$. According to Property 4.6, $\mathbf{L}(\Phi, \Phi) = \mathbf{G}_\Phi$. Taking into account Theorem 4.5, we receive $A$ and $B$ are right associates.

For similar reasons, the equality $B^T = A^T D^T$ implies that the matrices $B^T$ and $A^T$ are right associates. Consequently, the matrices $A$ and $B$ are left associates. $\qquad\square$

**Example 5.1.** We suggest an example of matrices that are left and right divisors of each other. Let

$$A = P_1^{-1} \Phi Q_1^{-1} = P_2^{-1} \Phi Q_2^{-1},$$

i.e., $P_1, P_2 \in \mathbf{P}_A$, and $Q_1, Q_2 \in \mathbf{Q}_A$. Then the matrices

$$M = P_1^{-1} \Phi Q_2^{-1}, \; N = P_2^{-1} \Phi Q_1^{-1}$$

**153**

will be desired. Indeed,

$$M = P_1^{-1}\Phi Q_2^{-1} = (P_1^{-1}P_2)P_2^{-1}\Phi Q_2^{-1} =$$
$$= (P_1^{-1}P_2)P_1^{-1}\Phi Q_1^{-1} = (P_1^{-1}P_2)^2 P_2^{-1}\Phi Q_1^{-1} = (P_1^{-1}P_2)^2 N.$$

Analogously, we show that $N = M(Q_2 Q_1^{-1})^2$. $\diamond$

## 5.2. Unassociated matrix divisors and Kazimirskii set

Our further research is aimed at finding the set $\mathbf{W}(\mathrm{E}, \Phi)$. Let $f \in R$. Consider the factor-ring $R/R_f$. Denote by $K(f)$ the set of representatives of this factor-ring cosets. Assume that

$$\mathrm{E} = \mathrm{diag}(\varepsilon_1, ..., \varepsilon_n), \quad \Phi = \mathrm{diag}(\varphi_1, ..., \varphi_n)$$

be nonsingular $d$-matrices, where $\Phi|\mathrm{E}$. Denote by $\mathbf{V}(\mathrm{E}, \Phi)$ the set of lower unitriangular matrices of the form

$$\left\|\begin{array}{ccccc} 1 & 0 & ... & 0 & 0 \\ \dfrac{\varphi_2}{(\varphi_2, \varepsilon_1)}k_{21} & 1 & ... & 0 & 0 \\ ... & ... & ... & ... & ... \\ \dfrac{\varphi_n}{(\varphi_n, \varepsilon_1)}k_{n1} & \dfrac{\varphi_n}{(\varphi_n, \varepsilon_2)}k_{n2} & ... & \dfrac{\varphi_n}{(\varphi_n, \varepsilon_{n-1})}k_{n.n-1} & 1 \end{array}\right\|,$$

where $k_{ij} \in K\left(\dfrac{(\varphi_i, \varepsilon_j)}{\varphi_j}\right)$, $i = 2, ..., n$, $j = 1, ..., n-1$, $i > j$.

The set $\mathbf{V}(\mathrm{E}, \Phi)$ is called the **Kazimirskii set** in honor of a famous Ukrainian algebraist who first considered matrices of this form.

We establish the relationship between the generating set $\mathbf{L}(\mathrm{E}, \Phi)$, the Kazimirskii set $\mathbf{V}(\mathrm{E}, \Phi)$, and groups $\mathbf{G}_\Phi$, $\mathbf{G}_\mathrm{E}$.

**Theorem 5.5.** *Let $R$ be a Bezout ring of stable range* 1.5. *Suppose that $\Phi$ is a non-singular $d$-matrix and $\Phi|\mathrm{E}$, then*

$$\mathbf{L}(\mathrm{E}, \Phi) = \mathbf{G}_\Phi \mathbf{V}(\mathrm{E}, \Phi) U_n^{up}(R).$$

**Proof.** Let $L \in \mathbf{L}(\mathrm{E}, \Phi)$. By Theorem 2.14, there are $H \in \mathbf{G}_\Phi$ a lower unitriangular matrix $S$ and an upper unitriangular matrix $U$ such that $L = HSU$. Consequently $S = H^{-1}LU^{-1}$. According to Corollary 2.2, the group of upper unitriangular matrices is a subgroup of any Zelisko group. So $U \in \mathbf{G}_\mathrm{E}$. Using Properties 4.1, 4.2, we get $S \in \mathbf{L}(\mathrm{E}, \Phi)$. On the basis of Lemma 5.3 in the group $\mathbf{G}_\Phi$ there is the matrix $H_1$ such that $H_1 S = V \in \mathbf{V}(\mathrm{E}, \Phi)$. For this reason

$$L = HSU = (HH_1^{-1})(H_1 S)U = H_2 VU.$$

Noting that $H_2 \in \mathbf{G}_\Phi$, we get

$$\mathbf{L}(\mathrm{E}, \Phi) \subseteq \mathbf{G}_\Phi \mathbf{V}(\mathrm{E}, \Phi) U_n^{up}(R).$$

In accordance with Properties 4.1, 4.2,

$$\mathbf{G}_\Phi \mathbf{L}(\mathrm{E}, \Phi) = \mathbf{L}(\mathrm{E}, \Phi), \quad \mathbf{L}(\mathrm{E}, \Phi)\mathbf{G}_\mathrm{E} = \mathbf{L}(\mathrm{E}, \Phi).$$

Since $U_n^{up}(R) \subset \mathbf{G}_\mathrm{E}$ and $I \in U_n^{up}(R)$, we have

$$\mathbf{L}(\mathrm{E}, \Phi)U_n^{up}(R) = \mathbf{L}(\mathrm{E}, \Phi).$$

Taking into account $\mathbf{V}(\mathrm{E}, \Phi) \subset \mathbf{L}(\mathrm{E}, \Phi)$, we get

$$\mathbf{G}_\Phi \mathbf{V}(\mathrm{E}, \Phi)U_n^{up}(R) \subseteq \mathbf{L}(\mathrm{E}, \Phi).$$

Consequently,

$$\mathbf{L}(\mathrm{E}, \Phi) = \mathbf{G}_\Phi \mathbf{V}(\mathrm{E}, \Phi)U_n^{up}(R).$$

The theorem is proved. $\qquad\qquad\qquad\qquad\qquad\qquad\square$

**Corollary 5.1.** *Let $R$ be a Bezout ring of stable range* 1.5. *The set*

$$(\mathbf{V}(\mathrm{E}, \Phi)\, U_n^{up}(R)P)^{-1}\, \Phi \mathrm{GL}_n(R)$$

*is the set of all left divisors of the matrix $A = P^{-1}\mathrm{E}Q^{-1}$ with the Smith form $\Phi$.*

**Proof.** By Corollary 4.4 the set $(\mathbf{L}(\mathrm{E}, \Phi)P)^{-1}\Phi \mathrm{GL}_n(R)$ is the set of all left divisors of the matrix $A$ with the Smith form $\Phi$. Taking into account Theorem 5.5, we have

$$(\mathbf{L}(\mathrm{E}, \Phi)P)^{-1}\Phi \mathrm{GL}_n(R) = (\mathbf{G}_\Phi \mathbf{V}(\mathrm{E}, \Phi)U_n^{up}(R)P)^{-1}\Phi \mathrm{GL}_n(R) =$$

$$= P^{-1}U_n^{up}(R)\mathbf{V}^{-1}(\mathrm{E}, \Phi)\mathbf{G}_\Phi \Phi \mathrm{GL}_n(R) = P^{-1}U_n^{up}(R)\mathbf{V}^{-1}(\mathrm{E}, \Phi)\Phi \mathrm{GL}_n(R) =$$

$$= (\mathbf{V}(\mathrm{E}, \Phi)U_n^{up}(R)P)^{-1}\Phi \mathrm{GL}_n(R). \qquad\qquad\qquad \square$$

**Theorem 5.6.** *The set $\mathbf{W}(\mathrm{E}, \Phi)$ can be chosen so that*

$$\mathbf{V}(\mathrm{E}, \Phi) \subseteq \mathbf{W}(\mathrm{E}, \Phi).$$

**Proof.** Let $V, V_1$ be matrices from $\mathbf{V}(\mathrm{E}, \Phi)$ with the elements $f_{ij}v_{ij}$, $f_{ij}u_{ij}$, respectively, where $f_{ij} = \dfrac{\varphi_i}{(\varphi_i, \varepsilon_j)}$, $i > j$, and let $HV = V_1$, where $H \in \mathbf{G}_\Phi$. The assertion follows if we prove that $V = V_1$. It is obvious that the matrix $H$ is also a lower unitriangular matrix with the elements $\dfrac{\varphi_i}{\varphi_j}h_{ij}$, $i > j$. We have $H = V_1 V^{-1}$. Putting $n = 2$, we get

$$f_{21}u_{21} - f_{21}v_{21} = \frac{\varphi_2}{\varphi_1}h_{21}.$$

**155**

Notice that
$$\frac{\varphi_2}{\varphi_1} = f_{21}\frac{(\varphi_2, \varepsilon_1)}{\varphi_1},$$
thus
$$u_{21} - v_{21} = \frac{(\varphi_2, \varepsilon_1)}{\varphi_1}h_{21}.$$
Hence,
$$u_{21} \equiv v_{21}\left(\mathrm{mod}\,\frac{(\varphi_2, \varepsilon_1)}{\varphi_1}\right).$$

This means that $u_{21} = v_{21}$, so $V = V_1$.

Suppose that the assumption holds for the matrices of the order $n-1$, we will prove it for $n$. The equality $HV = V_1$ implies that the equalities

$$H'V' = \left\|\begin{matrix} 1 & & & & 0 \\ \dfrac{\varphi_2}{\varphi_1}h_{21} & 1 & & & \\ ... & ... & & & \\ \dfrac{\varphi_{n-1}}{\varphi_1}h_{n-1.1} & \dfrac{\varphi_{n-1}}{\varphi_2}h_{n-1.2} & ... & \dfrac{\varphi_{n-1}}{\varphi_{n-2}}h_{n-1.n-2} & 1 \end{matrix}\right\| \times$$

$$\times \left\|\begin{matrix} 1 & & & & 0 \\ f_{21}v_{21} & 1 & & & \\ ... & ... & & & \\ f_{n-1.1}v_{n-1.1} & f_{n-1.2}v_{n-1.2} & ... & f_{n-1.n-2}v_{n-1.n-2} & 1 \end{matrix}\right\| =$$

$$= \left\|\begin{matrix} 1 & 0 & & & 0 \\ f_{21}u_{21} & 1 & & & \\ ... & ... & & & \\ f_{n-1.1}u_{n-1.1} & f_{n-1.2}u_{n-1.2} & ... & f_{n-1.n-2}u_{n-1.n-2} & 1 \end{matrix}\right\| = V_1',$$

and

$$H''V'' = \left\|\begin{matrix} 1 & & & & 0 \\ \dfrac{\varphi_3}{\varphi_2}h_{21} & 1 & & & \\ ... & ... & & & \\ \dfrac{\varphi_n}{\varphi_2}h_{n2} & \dfrac{\varphi_n}{\varphi_3}h_{n3} & ... & \dfrac{\varphi_n}{\varphi_{n-1}}h_{n.n-1} & 1 \end{matrix}\right\| \times$$

$$\times \left\|\begin{matrix} 1 & & & & 0 \\ f_{32}v_{32} & 1 & & & \\ ... & ... & & & \\ f_{n2}v_{n2} & f_{n3}v_{n3} & ... & f_{n.n-1}v_{n.n-1} & 1 \end{matrix}\right\| =$$

$$= \left\|\begin{matrix} 1 & & & & 0 \\ f_{32}u_{32} & 1 & & & \\ ... & ... & & & \\ f_{n2}u_{n2} & f_{n3}u_{n3} & ... & f_{n.n-1}u_{n.n-1} & 1 \end{matrix}\right\| = V_1''$$

hold. Since $H' \in \mathbf{G}_{\mathrm{diag}(\varphi_1, ..., \varphi_{n-1})}$ and

$$V', V_1' \in \mathbf{V}(\mathrm{diag}(\varepsilon_1, ..., \varepsilon_{n-1}), \mathrm{diag}(\varphi_1, ..., \varphi_{n-1})),$$

by the induction hypothesis, $V' = V_1'$. Analogously, $H'' \in \mathbf{G}_{\mathrm{diag}(\varphi_2, ..., \varphi_n)}$ and $V'', V_1'' \in \mathbf{V}(\mathrm{diag}(\varepsilon_2, ..., \varepsilon_n), \mathrm{diag}(\varphi_2, ..., \varphi_n))$ implies $V'' = V_1''$. It follows that the matrices $V, V_1$ differ from each other by the entry $(n, 1)$ at most. Hence,

$$V_1 V^{-1} = \left\| \begin{matrix} 1 & & & & 0 \\ 0 & 1 & & & \\ \vdots & & \ddots & & \\ 0 & 0 & & 1 & \\ s_{n1} & 0 & ... & 0 & 1 \end{matrix} \right\|,$$

where $s_{n1} = f_{n1}(u_{n1} - v_{n1})$. Thus,

$$u_{n1} \equiv v_{n1} \left( \mathrm{mod} \frac{(\varphi_n, \varepsilon_1)}{\varphi_1} \right).$$

This means that $u_{n1} = v_{n1}$. Consequently, $V = V_1$ and the proof is complete. □

Let $\Phi = \mathrm{diag}(\varphi_1, ..., \varphi_n)$ be a nonsingular $d$-matrix and $2 \leq j_1 < j_2 ... < < j_g \leq n$ is a set of indices where $\varphi_i, \varphi_{i-1}$ elements are non-associative in the ring $R$, $i \in \{j_1, j_2, ..., j_g\}$.

**Theorem 5.7.** *The sets* $\mathbf{V}(\mathrm{E}, \Phi)$ *and* $\mathbf{W}(\mathrm{E}, \Phi)$ *coincide if and only if any divisor of the element* $\dfrac{\varphi_i}{\varphi_{i-1}}$ *has a common divisor with the element* $\dfrac{\varphi_i}{(\varphi_i, \varepsilon_{i-1})}$, $i = j_1, j_2, ..., j_g$.

In order to prove this theorem we establish a series of facts which present interest in their own right.

**Lemma 5.1.** *Let each non-trivial divisor of* $\dfrac{\varphi_i}{\varphi_{i-1}}$ *has a common divisor with* $\dfrac{\varphi_i}{(\varphi_i, \varepsilon_{i-1})}$. *If for some $d$ the condition*

$$\left( \frac{\varphi_i}{(\varphi_i, \varepsilon_{i-1})}, d \right) = 1,$$

*is satisfied, then*

$$\left( \frac{\varphi_i}{(\varphi_{i-1}, d)} \right) = 1.$$

**Proof.** Assume that

$$\left( \frac{\varphi_i}{(\varphi_{i-1}, d)} \right) = \alpha_i \neq 1.$$

It follows that $\alpha_i$ is a divisor of $\dfrac{\varphi_i}{(\varphi_{i-1})}$. According to the assumption of our assumption

$$\left( \frac{\varphi_i}{(\varphi_i, \varepsilon_{i-1})}, \alpha_i \right) \neq 1$$

which is a contradiction. □

**Lemma 5.2.** *Let L be an invertible matrix of the form*

$$L = \left\| \begin{matrix} l_{11} & l_{12} & ... & l_{1.n-1} & l_{1n} \\ \dfrac{\varphi_2}{(\varphi_2,\varepsilon_1)}l_{21} & l_{22} & ... & l_{2.n-1} & l_{2n} \\ ... & ... & ... & ... & ... \\ \dfrac{\varphi_n}{(\varphi_n,\varepsilon_1)}l_{n1} & \dfrac{\varphi_n}{(\varphi_n,\varepsilon_2)}l_{n2} & ... & \dfrac{\varphi_n}{(\varphi_n,\varepsilon_{n-1})}l_{n.n-1} & l_{nn} \end{matrix} \right\| = \|l'_{ij}\|_1^n.$$

*Then*

$$\left( \frac{\varphi_i}{(\varphi_i,\varepsilon_{i-1})}, l'_{ij}, l'_{i+1.j}, ..., l'_{nj} \right) = 1,$$

$i = 2, 3, ..., n, \ j = i, i+1, ..., n.$

**Proof.** Suppose, contrary to our claim, that

$$\left( \frac{\varphi_i}{(\varphi_i,\varepsilon_{i-1})}, l_{ij}, l_{i+1.j}, ..., l_{nj} \right) = \delta_{ij} \neq 1.$$

Let us consider the submatrix

$$L_{ij} = \left\| \begin{matrix} \dfrac{\varphi_i}{(\varphi_i,\varepsilon_1)}l_{i1} & ... & \dfrac{\varphi_i}{(\varphi_i,\varepsilon_{i-1})}l_{i.i-1} & l_{ij} \\ ... & ... & ... & ... \\ \dfrac{\varphi_n}{(\varphi_n,\varepsilon_1)}l_{n1} & ... & \dfrac{\varphi_n}{(\varphi_n,\varepsilon_{i-1})}l_{n.i-1} & l_{nj} \end{matrix} \right\|$$

of the matrix $L$. Since $\dfrac{\varphi_i}{(\varphi_i,\varepsilon_{i-1})} \Big| \dfrac{\varphi_k}{(\varphi_k,\varepsilon_s)}$, $k = i, i+1, ..., n$, $s = 1, ..., i-1$, we have $\delta_{ij}|L_{ij}$. By Proposition 3.1, $\delta_{ij}|\det L$. This contradicts to the fact that $L$ is an invertible matrix. $\square$

**Lemma 5.3.** *Let S be a lower unitriangular matrix from* $\mathbf{L}(E, \Phi)$. *Then there is a matrix* $H \in \mathbf{G}_\Phi$ *such that* $HS \in \mathbf{V}(E, \Phi)$.

**Proof.** According to Theorem 4.1 a matrix $S$ has the form

$$S = \left\| \begin{matrix} 1 & 0 & ... & 0 & 0 \\ \dfrac{\varphi_2}{(\varphi_2,\varepsilon_1)}l_{21} & 1 & & 0 & 0 \\ ... & ... & ... & ... & ... \\ \dfrac{\varphi_n}{(\varphi_n,\varepsilon_1)}l_{n1} & \dfrac{\varphi_n}{(\varphi_n,\varepsilon_2)}l_{n2} & ... & \dfrac{\varphi_n}{(\varphi_n,\varepsilon_{n-1})}l_{n.n-1} & 1 \end{matrix} \right\|.$$

Consider the matrix

$$H_0 = \left\| \begin{matrix} 1 & 0 & ... & 0 & 0 \\ \dfrac{\varphi_2}{\varphi_1}h_{21} & 1 & & 0 & 0 \\ ... & ... & ... & ... & ... \\ \dfrac{\varphi_n}{\varphi_1}h_{n1} & \dfrac{\varphi_n}{\varphi_2}h_{n2} & ... & \dfrac{\varphi_n}{\varphi_{n-1}}h_{n.n-1} & 1 \end{matrix} \right\|,$$

where $h_{ij}$ are parameters, $i = 2, ..., n$, $j = 1, ..., n-1$, $i > j$. If $n = 2$, then

$$H_0 S = \left\| \begin{matrix} 1 & 0 \\ \dfrac{\varphi_2}{\varphi_1} h_{21} & 1 \end{matrix} \right\| \left\| \begin{matrix} 1 & 0 \\ \dfrac{\varphi_2}{(\varphi_2, \varepsilon_1)} s_{21} & 1 \end{matrix} \right\| = \left\| \begin{matrix} 1 & 0 \\ \dfrac{\varphi_2}{(\varphi_2, \varepsilon_1)} \left( \dfrac{(\varphi_2, \varepsilon_1)}{\varphi_1} h_{21} + s_{21} \right) & 1 \end{matrix} \right\| = S_1.$$

Let $s_{21} \equiv k_{21} \left( \bmod \dfrac{(\varphi_2, \varepsilon_1)}{\varphi_1} \right)$, where $k_{21} \in K \left( \dfrac{(\varphi_2, \varepsilon_1)}{\varphi_1} \right)$. It follows that $k_{21} =$

$= s_{21} + \dfrac{(\varphi_2, \varepsilon_1)}{\varphi_1} r_{21}$ for some $r_{21} \in R$. Setting $h_{21} = r_{21}$, we obtain $S_1 \in \mathbf{V}(\mathrm{E}, \Phi)$.

Suppose that the assumption holds for the matrices of the order $n-1$, we will prove it for $n$. The matrix $H_0 S$ is also a lower unitriangular matrix with the elements $d_{ij}$, $i > j$. We have

$$d_{nj} = \left\| \begin{matrix} \dfrac{\varphi_n}{\varphi_1} h_{n1} \, ... \, \dfrac{\varphi_n}{\varphi_{n-1}} h_{n.n-1} & 1 \end{matrix} \right\| \left\| \underbrace{0 ... 0}_{j-1} \quad 1 \quad \dfrac{\varphi_{j+1}}{(\varphi_{j+1}, \varepsilon_j)} s_{j+1.j} \, ... \, \dfrac{\varphi_n}{(\varphi_n, \varepsilon_j)} s_{nj} \right\|^T =$$

$$= \dfrac{\varphi_n}{\varphi_j} h_{nj} + \dfrac{\varphi_n}{(\varphi_{j+1}, \varepsilon_j)} h_{n.j+1} s_{j+1.j} + ... + \dfrac{\varphi_n}{(\varphi_{n-1}, \varepsilon_j)} h_{n.n-1} s_{n-1.j} + \dfrac{\varphi_n}{(\varphi_n, \varepsilon_j)} s_{nj} =$$

$$= \dfrac{\varphi_n}{(\varphi_n, \varepsilon_j)} \left( \dfrac{(\varphi_n, \varepsilon_j)}{\varphi_j} h_{nj} + \dfrac{(\varphi_n, \varepsilon_j)}{(\varphi_{j+1}, \varepsilon_j)} h_{n.j+1} s_{j+1.j} + ... + \right.$$

$$\left. + \dfrac{(\varphi_n, \varepsilon_j)}{(\varphi_{n-1}, \varepsilon_j)} h_{n.n-1} s_{n-1.j} + s_{nj} \right),$$

$j = 1, ..., n-1$. Let $j = n-1$ and $s_{n.n-1} \equiv k_{n.n-1} \left( \bmod \dfrac{(\varphi_n, \varepsilon_{n-1})}{\varphi_{n-1}} \right)$, where

$k_{n.n-1} \in K \left( \dfrac{(\varphi_n, \varepsilon_{n-1})}{\varphi_{n-1}} \right)$. This gives $k_{n.n-1} = s_{n.n-1} + \dfrac{(\varphi_n, \varepsilon_{n-1})}{\varphi_{n-1}} r_{n.n-1}$ for

some $r_{n.n-1} \in R$. We put $h_{n.n-1} = r_{n.n-1}$. Set $j = n-2$ and

$$s_{n.n-2} + \dfrac{(\varphi_n, \varepsilon_{n-2})}{(\varphi_{n-1}, \varepsilon_{n-2})} r_{n.n-1} s_{n-1.n-2} \equiv k_{n.n-2} \left( \bmod \dfrac{(\varphi_n, \varepsilon_{n-2})}{\varphi_{n-2}} \right),$$

for some $k_{n.n-2} \in K \left( \dfrac{(\varphi_n, \varepsilon_{n-2})}{\varphi_{n-2}} \right)$. Then there exists $r_{n.n-2} \in R$ such that

$$k_{n.n-2} = s_{n.n-2} + \dfrac{(\varphi_n, \varepsilon_{n-2})}{(\varphi_{n-1}, \varepsilon_{n-2})} r_{n.n-1} s_{n-1.n-2} + \dfrac{(\varphi_n, \varepsilon_{n-2})}{\varphi_{n-2}} r_{n.n-2}.$$

We set $h_{n.n-2} = r_{n.n-2}$. Continuing the described process we obtain the lower unitriangular matrix

$$H_1 = \left\| \begin{matrix} I_{n-1} & \mathbf{0} \\ h & 1 \end{matrix} \right\|,$$

where $E_{n-1}$ is the identity $(n-1) \times (n-1)$ matrix,

$$h = \left\| \frac{\varphi_n}{\varphi_1} r_{n1} \quad ... \quad \frac{\varphi_n}{\varphi_{n-1}} r_{n.n-1} \right\|,$$

such that

$$H_1 S = \left\| \begin{matrix} S' & \mathbf{0} \\ g & 1 \end{matrix} \right\|,$$

where $S'$ is a lower unitriangular matrix from

$$\mathbf{L}(\mathrm{diag}(\varepsilon_1, ..., \varepsilon_{n-1}), \mathrm{diag}(\varphi_1, ..., \varphi_{n-1})),$$

$$g = \left\| \frac{\varphi_n}{(\varphi_n, \varepsilon_1)} k_{n1} \quad ... \quad \frac{\varphi_n}{(\varphi_n, \varepsilon_{n-1})} k_{n.n-1} \right\|,$$

$k_{nj} \in K\left( \frac{(\varphi_n, \varepsilon_j)}{\varphi_j} \right)$, $j = 1, ..., n - 1$. Thus, by the induction hypothesis, there exists $H' \in \mathbf{G}_{\mathrm{diag}(\varphi_1, ..., \varphi_{n-1})}$ such that

$$H'S' \in \mathbf{V}(\mathrm{diag}(\varepsilon_1, ..., \varepsilon_{n-1}), \mathrm{diag}(\varphi_1, ..., \varphi_{n-1})).$$

Hence, $(H' \oplus 1)H_1 S \in \mathbf{V}(\mathrm{E}, \Phi)$, and this is precisely the assertion of the Lemma. □

**Lemma 5.4.** *Let $\bar{l} = \left\| l_1 \quad ... \quad l_n \right\|^T$ – be a primitive column and*

$$\left( \frac{\varphi_n}{\varphi_1} l_1, ..., \frac{\varphi_n}{\varphi_{n-1}} l_{n-1}, l_n \right) = 1.$$

*Then there exists a matrix $H \in \mathbf{G}_\Phi$ such that $H\bar{l} = \left\| 0 \quad ... \quad 0 \quad 1 \right\|^T$.*

**Proof.** There are elements $u_1, ..., u_n$, such that

$$\frac{\varphi_n}{\varphi_1} l_1 u_1 + ... + \frac{\varphi_n}{\varphi_{n-1}} l_{n-1} u_{n-1} + l_n u_n = 1.$$

It follows that

$$\left( \frac{\varphi_n}{\varphi_1} u_1, ..., \frac{\varphi_n}{\varphi_{n-1}} u_{n-1}, u_n \right) = 1.$$

Using Theorem 1.1, we complement the primitive row

$$\left\| \frac{\varphi_n}{\varphi_1} u_1 \quad ... \quad \frac{\varphi_n}{\varphi_{n-1}} u_{n-1} \quad u_n \right\|$$

to an invertible matrix of the form

$$H_1 = \left\| \begin{matrix} u_{11} & u_{12} & ... & u_{1.n-1} & u_{1n} \\ 0 & u_{22} & ... & u_{2.n-1} & u_{2n} \\ ... & ... & ... & ... & ... \\ 0 & ... & 0 & u_{n-1.n-1} & u_{n-1.n} \\ \frac{\varphi_n}{\varphi_1} u_1 & \frac{\varphi_n}{\varphi_2} u_2 & ... & \frac{\varphi_n}{\varphi_{n-1}} u_{n-1} & u_n \end{matrix} \right\|.$$

Based on Theorem 2.7, $H_1 \in \mathbf{G}_\Phi$. Since

$$H_1 \bar{l} = \left\| t_1 \quad \dots \quad t_{n-1} \quad 1 \right\|^T,$$

then

$$\underbrace{\left\|\begin{matrix} 1 & 0 & \dots & 0 & -t_1 \\ & \ddots & & & \vdots \\ 0 & & & 1 & -t_{n-1} \\ 0 & 0 & \dots & 0 & 1 \end{matrix}\right\|}_{H_2} H_1 \bar{l} = \left\| 0 \quad \dots \quad 0 \quad 1 \right\|^T.$$

The matrix $H_2$ belongs to the group $\mathbf{G}_\Phi$. So, the matrix $H = H_2 H_1$ is desired. The Lemma is proved. □

We proceed the proof of Theorem 5.7.

**Proof. Necessity.** Let $\delta_i$ be a non-trivial divisor of $\dfrac{\varphi_i}{\varphi_{i-1}}$, $i_1 \leqslant i \leqslant i_g$. Suppose, contrary to our claim, that $\left( \dfrac{\varphi_i}{(\varphi_i, \varepsilon_{i-1})}, \delta_i \right) = 1$. Then there exist $u, v \in R$ such that $u \dfrac{\varphi_i}{(\varphi_i, \varepsilon_{i-1})} + v \delta_i = 1$. Consider the matrix

$$L_i = I_{i-2} \oplus \left\| \begin{matrix} v & -u \\ \dfrac{\varphi_i}{(\varphi_i, \varepsilon_{i-1})} & \delta_i \end{matrix} \right\| \oplus I_{n-i}.$$

By Theorem 4.1 $L_i \in \mathbf{L}(\mathrm{E}, \Phi)$. A matrix consisting of the last $n-i+1$ columns of the matrix $L_i$ has the form

$$S_i = \left\| \begin{matrix} 0 & 0 & 0 & \dots & 0 \\ \dots & \dots & \dots & \dots & \dots \\ 0 & 0 & 0 & \dots & 0 \\ -u & 0 & 0 & \dots & 0 \\ \delta_i & 0 & 0 & \dots & 0 \\ 0 & 1 & 0 & \dots & 0 \\ 0 & 0 & 1 & & 0 \\ & & & \ddots & \\ 0 & 0 & 0 & & 1 \end{matrix} \right\|.$$

Then

$$\Phi_i S_i = \left\| \begin{matrix} 0 & 0 & 0 & \dots & 0 \\ \dots & \dots & \dots & \dots & \dots \\ 0 & 0 & 0 & \dots & 0 \\ -u \dfrac{\varphi_i}{\varphi_{i-1}} & 0 & 0 & \dots & 0 \\ \delta_i & 0 & 0 & \dots & 0 \\ 0 & 1 & 0 & \dots & 0 \\ 0 & 0 & 1 & & 0 \\ & & & \ddots & \\ 0 & 0 & 0 & & 1 \end{matrix} \right\|.$$

**161**

Since $\delta_i \left| \dfrac{\varphi_i}{\varphi_{i-1}} \right.$, then

$$\frac{\varphi_i}{\varphi_{i-1}} = \delta_i \gamma_i.$$

So

$$\left(I_{i-2} \oplus \left\| \begin{matrix} 1 & u\gamma_i \\ 0 & 1 \end{matrix} \right\| \oplus I_{n-i}\right) \Phi_i S_i =$$

$$= \left\| \begin{matrix} 0 & 0 & 0 & \dots & 0 \\ \dots & \dots & \dots & \dots & \dots \\ 0 & 0 & 0 & \dots & 0 \\ \delta_i & 0 & 0 & \dots & 0 \\ 0 & 1 & 0 & \dots & 0 \\ 0 & 0 & 1 & & 0 \\ & & & \ddots & \\ 0 & 0 & 0 & & 1 \end{matrix} \right\|.$$

It means that

$$\Phi_i S_i \overset{l}{\sim} \left\| \begin{matrix} 0 & 0 & 0 & \dots & 0 \\ \dots & \dots & \dots & \dots & \dots \\ 0 & 0 & 0 & \dots & 0 \\ \delta_i & 0 & 0 & \dots & 0 \\ 0 & 1 & 0 & \dots & 0 \\ 0 & 0 & 1 & & 0 \\ & & & \ddots & \\ 0 & 0 & 0 & & 1 \end{matrix} \right\|.$$

On the other hand, the last $n-i+1$ columns of any matrix of the set $\mathbf{V}(\mathrm{E}, \Phi)$ have the form

$$M_i = \left\| \begin{matrix} 0 & 0 & \dots & 0 & 0 \\ \dots & \dots & \dots & \dots & \dots \\ 0 & 0 & \dots & 0 & 0 \\ 1 & 0 & \dots & 0 & 0 \\ \dfrac{\varphi_{i+1}}{(\varphi_{i+1}, \varepsilon_i)} l_{i+1.i} & 1 & \dots & 0 & 0 \\ \dots & \dots & \dots & & \\ \dfrac{\varphi_n}{(\varphi_n, \varepsilon_i)} l_{ni} & \dfrac{\varphi_n}{(\varphi_n, \varepsilon_{i+1})} l_{n.i+1} & \dots & \dfrac{\varphi_n}{(\varphi_n, \varepsilon_{n-1})} l_{n.n-1} & 1 \end{matrix} \right\|.$$

Hence,

$$\Phi_i M_i = M_i \overset{l}{\sim} \left\| \begin{matrix} \mathbf{0} \\ I \end{matrix} \right\|.$$

Suppose that $V_i \in \mathbf{V}(\mathrm{E}, \Phi)$ is the representative of the coset $\mathbf{G}_\Phi L_i$. That is, in the group $\mathbf{G}_\Phi$ there exists a matrix $H$, such that $V_i = HL_i$. By Lemma 3.5, $\Phi_i V_i \overset{l}{\sim} \Phi_i L_i$. It follows that $\Phi_i S_i \overset{l}{\sim} \Phi_i M_i$. However, the matrix $\Phi_i S_i$ has the left Hermite form

$$\left\| \begin{array}{ccccc} 0 & 0 & 0 & ... & 0 \\ ... & ... & ... & ... & ... \\ 0 & 0 & 0 & ... & 0 \\ \delta_i & 0 & 0 & ... & 0 \\ 0 & 1 & 0 & ... & 0 \\ 0 & 0 & 1 & & 0 \\ & & & \ddots & \\ 0 & 0 & 0 & & 1 \end{array} \right\|,$$

and the matrix $\Phi_i M_i$:

$$\left\| \begin{array}{c} \mathbf{0} \\ I_{n-i+1} \end{array} \right\|.$$

We came to a contradiction.

**Sufficiency.** Let $L$ be an invertible matrix of the form

$$L = \left\| \begin{array}{ccccc} l_{11} & l_{12} & ... & l_{1.n-1} & l_{1n} \\ \dfrac{\varphi_2}{(\varphi_2, \varepsilon_1)} l_{21} & l_{22} & ... & l_{2.n-1} & l_{2n} \\ ... & ... & ... & ... & ... \\ \dfrac{\varphi_n}{(\varphi_n, \varepsilon_1)} l_{n1} & \dfrac{\varphi_n}{(\varphi_n, \varepsilon_2)} l_{n2} & ... & \dfrac{\varphi_n}{(\varphi_n, \varepsilon_{n-1})} l_{n.n-1} & l_{nn} \end{array} \right\|.$$

In order to prove this statement we only need to show that there exists a matrix $H \in \mathbf{G}_\Phi$ such that $HL \in \mathbf{V}(\mathrm{E}, \Phi)$. The proof will be divided into 2 steps. At the first step we will find a matrix $H_1 \in \mathbf{G}_\Phi$ such that $H_1 L$ has a lower unitriangular form.

In the case $n = 2$ the invertibility of the matrix $L$ implies that

$$\left( \frac{\varphi_2}{(\varphi_2, \varepsilon_1)}, l_{22} \right) = 1.$$

By Lemma 5.1

$$\left( \frac{\varphi_2}{\varphi_1}, l_{22} \right) = 1.$$

Therefore,

$$\frac{\varphi_2}{\varphi_1} l_{12} u_1 + l_{22} u_2 = 1$$

for some $u_1, u_2 \in R$. Clearly, we have

$$\left\| \begin{array}{cc} l_{22} & -l_{12} \\ \dfrac{\varphi_2}{\varphi_1} u_1 & u_2 \end{array} \right\| \left\| \begin{array}{cc} l_{11} & l_{12} \\ \dfrac{\varphi_2}{(\varphi_2, \varepsilon_1)} l_{21} & l_{22} \end{array} \right\| = \left\| \begin{array}{cc} e & 0 \\ \dfrac{\varphi_2}{(\varphi_2, \varepsilon_1)} s_{21} & 1 \end{array} \right\|,$$

where $e \in U(R)$. Hence

$$H_1 = \left\| \begin{matrix} e^{-1}l_{22} & -e^{-1}l_{21} \\ \dfrac{\varphi_2}{\varphi_1}u_1 & u_2 \end{matrix} \right\|$$

is the desired matrix.

Suppose that the assumption holds for the matrices of order $n-1$, we will prove it for $n$. Step by step, using Lemmas 5.2 and 5.1, we obtain

$$\left( \frac{\varphi_n}{\varphi_1}l_{1n}, ..., \frac{\varphi_n}{\varphi_{n-1}}l_{n-1.n}, l_{nn} \right) =$$

$$= \left( \frac{\varphi_n}{\varphi_{n-1}} \left( \frac{\varphi_{n-1}}{\varphi_1}l_{1n}, ..., \frac{\varphi_{n-1}}{\varphi_{n-2}}l_{n-2.n}, l_{n-1.n} \right), l_{nn} \right) =$$

$$= \left( \frac{\varphi_{n-1}}{\varphi_1}l_{1n}, ..., \frac{\varphi_{n-1}}{\varphi_{n-2}}l_{n-2.n}, l_{n-1.n}, l_{nn} \right) =$$

$$= \left( \frac{\varphi_{n-1}}{\varphi_{n-2}} \left( \frac{\varphi_{n-2}}{\varphi_1}l_{1n}, ..., l_{n-2.n} \right), (l_{n-1.n}, l_{nn}) \right) =$$

$$= ... = \left( \frac{\varphi_2}{\varphi_1}l_{1n}, (l_{2n}, ..., l_{nn}) \right) = (l_{1n}, ..., l_{nn}) = 1.$$

By Lemma 5.4, there exists a matrix $H_0 \in \mathbf{G}_\Phi$ such that

$$H_0L = \left\| \begin{matrix} L' & \mathbf{0} \\ g & 1 \end{matrix} \right\|.$$

According to Property 4.2 (see p. 126), $H_0L \in \mathbf{L}(\mathrm{E}, \Phi)$. Therefore,

$$L' \in \mathbf{L}(\mathrm{diag}(\varepsilon_1, ..., \varepsilon_{n-1}), \mathrm{diag}(\varphi_1, ..., \varphi_{n-1})).$$

By the induction hypothesis, there exists a matrix $H' \in \mathbf{G}_{\mathrm{diag}(\varphi_1, ..., \varphi_{n-1})}$ such that the matrix $H'L'$ has lower unitriangular form. Hence the matrix $H_1 = (H' \oplus 1)H_0$ is the desired matrix, i.e., $H_1L \in \mathbf{L}(\mathrm{E}, \Phi)$ and has a lower unitriangular form.

By Lemma 5.3, $H_2H_1L \in \mathbf{V}(\mathrm{E}, \Phi)$ for some matrix $H_2 \in \mathbf{G}_\Phi$. The proof is complete. $\qquad\square$

Combining Theorems 4.6 and 5.7, we obtain the following result.

**Theorem 5.8.** *The set* $(\mathbf{V}(\mathrm{E}, \Phi)P_A)^{-1}\Phi$ *consists of all left non-associative by right divisors of the matrix* $A$, *which have the Smith form* $\Phi$ *if and only if each divisor of elements* $\dfrac{\varphi_i}{\varphi_{i-1}}$ *has a common non-trivial divisor with* $\dfrac{\varphi_i}{(\varphi_i, \varepsilon_{i-1})}$, $i = j_1, j_2, ..., j_g$. $\qquad\square$

**Remark.** If the conditions this Theorem do not hold, the set $(\mathbf{V}(\mathrm{E}, \Phi)P_A)^{-1}\Phi$ contains only part of desired matrix $A$ divisors.

**Example 5.2.** Consider the matrix

$$A = \left\| \begin{matrix} 5 & 0 \\ 0 & 5x^3 \end{matrix} \right\| = \mathrm{E}$$

over the ring

$$Q = \left\{ a_0 + \sum_{i=1}^{\infty} a_i x^i \middle| a_0 \in \mathbb{Z}, \ a_i \in \mathbb{Q}, \ i \in \mathbb{N} \right\}.$$

Let us find all its left non-associative by right divisors with Smith's form

$$\Phi = \left\| \begin{matrix} 1 & 0 \\ 0 & 5x \end{matrix} \right\|.$$

Since

$$x = 5 \left( \frac{1}{5} x \right),$$

then each divisor of the element $\frac{\varphi_2}{\varphi_1} = 5x$ has a common divisor with

$$\frac{\varphi_2}{(\varphi_2, \varepsilon_1)} = \frac{5x}{(5x, 5)} = x.$$

By Theorem 5.8, the required set of divisors has the form

$$\left\{ \left\| \begin{matrix} 1 & 0 \\ xk & 1 \end{matrix} \right\|^{-1} \Phi \right\} = \left\{ \left\| \begin{matrix} 1 & 0 \\ -xk & 5x \end{matrix} \right\| \right\},$$

where $k \in K \left( \frac{(\varphi_2, \varepsilon_1)}{\varphi_1} \right) = K(5) = \{0, 1, 2, 3, 4\}$. $\diamondsuit$

Let $R$ be a principal ideal ring. In this case Theorem 5.8 can be reformulated.

Decompose the elements $\frac{\varphi_i}{\varphi_{i-1}}$, $\frac{\varepsilon_{i-1}}{\varphi_{i-1}}$, $i = j_1, j_2, ..., j_g$, into the product of degrees of irreducible factors:

$$\frac{\varphi_i}{\varphi_{i-1}} = g_{i1}^{k_{i1}} ... g_{il}^{k_{il}}, \quad \frac{\varepsilon_{i-1}}{\varphi_{i-1}} = g_{i1}^{q_{i1}} ... g_{il}^{q_{il}} h_{i1}^{p_{i1}} ... h_{ir}^{p_{ir}}.$$

**Theorem 5.9.** *Let $R$ be a principal ideal ring. The set $(\mathbf{V}(\mathrm{E}, \Phi) P_A)^{-1} \Phi$ consists of all left right unassociated divisors of the matrix $A$, which have the Smith form $\Phi$ if and only if $k_{ij} > q_{ij}$, $i = j_1, j_2, ..., j_g$, $j = 1, ..., l$.*

**Proof.** Write the element $\dfrac{\varphi_i}{(\varphi_i, \varepsilon_{i-1})}$ in the form

$$\frac{\varphi_i}{(\varphi_i, \varepsilon_{i-1})} = \frac{\varphi_i / \varphi_{i-1}}{(\varphi_i / \varphi_{i-1}, \varepsilon_{i-1} / \varphi_{i-1})}, \tag{5.1}$$

**165**

$i = j_1, j_2, ..., j_g$. According to Theorem 5.7, in order that the set $(\mathbf{V}(\mathrm{E}, \Phi)P_A)^{-1}\Phi$ consists of all left right unassociated divisors of the matrix $A$, which have the Smith form $\Phi$ necessary and sufficient all $g_{im}$ divide $\dfrac{\varphi_i}{(\varphi_i, \varepsilon_{i-1})}$. This will be only when

$$g_{im} \left| \frac{g_{im}^{k_{im}}}{(g_{im}^{k_{im}}, g_{im}^{q_{im}})}. \right.$$

This is equivalent to $k_{im} > q_{im}$. $\qquad\square$

**Example 5.3.** Let $A$ be an integer matrix with the Smith form

$$\mathrm{E} = \mathrm{diag}(1,\, 2,\, 2^3 3^1 5^1,\, 2^5 3^2 5^2).$$

We describe the left non-associative by right divisors of the matrix $A$, which have the Smith form $\Phi = \mathrm{diag}(1,\, 1,\, 2^3,\, 2^4 3^2)$. In this case

$$\mathbf{V}(\mathrm{E}, \Phi) = \begin{Vmatrix} 1 & 0 & 0 & 0 \\ 0 & 1 & 0 & 0 \\ 0 & 2^2 r_{32} & 1 & 0 \\ 0 & 2^3 3^2 r_{42} & 2\,3\,r_{43} & 1 \end{Vmatrix}$$

where $r_{32}, r_{42} \in K(2) = \{0, 1\}$, $r_{43} \in K(3) = \{0, 1, 2\}$.

The set of indices for which $\varphi_i$ and $\varphi_{i-1}$ are unassociated, will be $\{3, 4\}$. Then

$$\frac{\varphi_3}{\varphi_2} = 2^3 \Rightarrow k_{31} = 3,$$

$$\frac{\varepsilon_2}{\varphi_2} = 2^1 \Rightarrow q_{31} = 1.$$

Therefore $k_{31} > q_{31}$. So

$$\frac{\varphi_4}{\varphi_3} = 2^1 3^2 \Rightarrow k_{41} = 1, k_{42} = 2,$$

$$\frac{\varepsilon_3}{\varphi_3} = 2^0 3^1 5^1 \Rightarrow q_{41} = 0, q_{42} = 1$$

and

$$k_{41} > q_{41}, k_{42} > q_{42}.$$

Since all the conditions of Theorem 5.9 are satisfied, so the set of all left right unassociated divisors of the matrix $A$, having the Smith form $\Phi$ has the form $(\mathbf{V}(\mathrm{E}, \Phi)P_A)^{-1}\Phi$.

If we look for the divisors of the $A$ matrix with the Smith form

$$\Phi_1 = \mathrm{diag}(1,\, 1,\, 2^3,\, 2^4 3),$$

then

$$\mathbf{V}(\mathrm{E}, \Phi_1) = \begin{Vmatrix} 1 & 0 & 0 & 0 \\ 0 & 1 & 0 & 0 \\ 0 & 2^2 r_{32} & 1 & 0 \\ 0 & 2^3 3 r_{42} & 2 r_{43} & 1 \end{Vmatrix},$$

where $r_{32}, r_{42} \in K(2)$, $r_{43} \in K(3)$. In this case

$$\frac{\varphi_4}{\varphi_3} = 2^1 3^1 \Rightarrow k_{41} = k_{42} = 1,$$

$$\frac{\varepsilon_3}{\varphi_3} = 2^0 3^1 \Rightarrow q_{41} = 0, \ q_{42} = 1.$$

and $k_{42} = q_{42} = 1$. That is, the conditions of Theorem 5.9 are not fulfilled and the set $(\mathbf{V}(\mathrm{E}, \Phi) P_A)^{-1} \Phi$ does not exhaust all left right unassociated divisors of the matrix $A$, with the Smith form $\Phi_1$. In particular, in the set $(\mathbf{V}(\mathrm{E}, \Phi) P_A)^{-1} \Phi$ there is no matrix associated with the matrix

$$P_A^{-1} \begin{Vmatrix} 1 & 0 & 0 & 0 \\ 0 & 1 & 0 & 0 \\ 0 & 0 & 1 & 1 \\ 0 & 0 & 2 & 3 \end{Vmatrix}^{-1} \Phi_1,$$

which is the left divisor of the matrix $A$. $\qquad\qquad\qquad \diamondsuit$

### 5.3. Value of a matrix on the root system of diagonal elements of $d$-matrix

We apply the obtained results to describe the divisors of polynomial matrices. For their presentation, we will need the results of this subsection.

Let $F$ be an algebraically closed field of characteristic zeros and $G(x) = = \|g_{ij}(x)\|$ be an $n \times p$ matrix over $F[x]$. The derived of the matrix $G(x)$ is called $G'(x) = \|g'_{ij}(x)\|$. Derivatives of higher orders will be denoted by $G''(x)$, $G'''(x), ..., G^{(t)}(x)$.

Let

$$\varphi(x) = (x - \alpha_1)^{k_1} (x - \alpha_1)^{k_2} ... (x - \alpha_r)^{k_r}$$

be monic polynomial over $F$.

**Definition 5.1.** [54]The value of the matrix $G(x)$ on the root system of polynomial $\varphi(x)$ is called the numerical matrix

$$M_{G(x)}(\varphi) = \begin{Vmatrix} H_1 \\ H_2 \\ ... \\ H_r \end{Vmatrix}, \quad H_i = \begin{Vmatrix} G(\alpha_i) \\ G'(\alpha_i) \\ ... \\ G^{(k_i-1)}(\alpha_i) \end{Vmatrix},$$

$i = 1, 2, ..., r$.

The form of the matrix $M_{G(x)}(\varphi)$ depends on the representation of the polynomial $\varphi(x)$ as the product of the factors $(x - \alpha_j)^{k_j}$. To emphasize this dependence, we also denote $M_{G(x)}(\varphi)$ by $M_{G(x)}[\alpha_1^{k_1}, ..., \alpha_r^{k_r}]$.

Let $g_i(x)$ be $i$th row of an $n \times m$ matrix $G(x)$, $i = 1, ..., n$, and

$$\Phi(x) = \text{diag} \, (\varphi_1(x), \varphi_2(x), ..., \varphi_n(x)),$$

where $\varphi_i(x)$ is monic polynomials, $i = 1, ..., n$.

**Definition 5.2.** The value of the matrix $G(x)$ on the root system of diagonal elements of the matrix $\Phi(x)$ is called the matrix

$$M_{G(x)}(\Phi) = \begin{Vmatrix} M_{g_1(x)}(\varphi_1) \\ M_{g_2(x)}(\varphi_2) \\ \text{............} \\ M_{g_n(x)}(\varphi_n) \end{Vmatrix},$$

where the block $M_{g_j(x)}(\varphi_j)$ is missing if $\deg \varphi_j(x) = 0$.

**Property 5.1.** *If $L = \|l_{ij}\|$ be a $p \times t$ matrix over $F$, then*

$$M_{G(x)L}(\Phi) = M_{G(x)}(\Phi)L.$$

**Proof.** Since

$$g_i(x)L = \begin{Vmatrix} g_{i1}(x) & ... & g_{ip}(x) \end{Vmatrix} \begin{Vmatrix} l_{11} & ... & l_{1t} \\ ... & ... & ... \\ l_{p1} & ... & l_{pt} \end{Vmatrix} =$$

$$= \begin{Vmatrix} g_{i1}(x)l_{11} + ... + g_{ip}(x)l_{p1} & ... & g_{i1}(x)l_{1t} + ... + g_{ip}(x)l_{pt} \end{Vmatrix},$$

then

$$(g_i(x)L)' = \begin{Vmatrix} g'_{i1}(x)l_{11} + ... + g'_{ip}(x)l_{p1} & ... & g'_{i1}(x)l_{1t} + ... + g'_{ip}(x)l_{pt} \end{Vmatrix} =$$

$$= \begin{Vmatrix} g'_{i1}(x) & ... & g'_{ip}(x) \end{Vmatrix} \begin{Vmatrix} l_{11} & ... & l_{1t} \\ ... & ... & ... \\ l_{p1} & ... & l_{pt} \end{Vmatrix} = g'_i L, \ \ i = 1, ..., n.$$

It follows that

$$(g_i(x)L)^{(j)} = g_i^{(j)}L.$$

Hence,

$$M_{g_i(x)L}(\varphi_i(x)) = M_{g_i(x)}(\varphi_i(x))L.$$

Therefore,

$$M_{G(x)L}(\Phi) = \begin{Vmatrix} M_{g_1(x)L}(\varphi_1) \\ \text{...............} \\ M_{g_n(x)L}(\varphi_n) \end{Vmatrix} = \begin{Vmatrix} M_{g_1(x)}(\varphi_1)L \\ \text{...............} \\ M_{g_n(x)}(\varphi_n)L \end{Vmatrix} = M_{G(x)}(\Phi)L.$$

The proof is complete. $\square$

**Property 5.2.** *The equality*

$$M_{G_1(x)+G_2(x)}(\Phi) = M_{G_1(x)}(\Phi) + M_{G_1(x)}(\Phi)$$

*is fulfilled.*

**Proof** does not cause any difficulties. $\square$

**Property 5.3.** *In order that*

$$M_{G(x)}(\Phi) = \mathbf{0}, \tag{5.2}$$

*necessary and sufficient* $G(x) = \Phi(x)Q(x)$.

**Proof.** In view of (5.2), if $\deg \varphi_i(x) \geqslant 1$, then

$$M_{g_i(x)}(\varphi_i(x)) = \mathbf{0}.$$

This is equivalent to the fact that each elements of the row $g_i(x)$ is a multiple to $\varphi_i(x)$. That is, the matrix $G(x)$ has the form

$$G(x) = \begin{Vmatrix} \varphi_1(x)q_{11}(x) & ... & \varphi_1(x)q_{1p}(x) \\ ... & ... & ... \\ \varphi_n(x)q_{n1}(x) & ... & \varphi_n(x)q_{np}(x) \end{Vmatrix} =$$

$$= \operatorname{diag}(\varphi_1(x), ..., \varphi_n(x)) \begin{Vmatrix} q_{11}(x) & ... & q_{1p}(x) \\ ... & ... & ... \\ q_{n1}(x) & ... & q_{np}(x) \end{Vmatrix} = \Phi(x)Q(x).$$

The proof is complete. $\qquad\qquad\square$

Let $A(x)$ be a nonsingular matrix over $F[x]$. There are an invertible matrices $P(x)$ and $Q(x)$ such that

$$P(x)A(x)Q(x) = \operatorname{diag}(\varphi_1(x), \varphi_2(x), ..., \varphi_n(x)) = \Phi(x),$$

where $\Phi(x)$ is the Smith form of the matrix $A(x)$. Let us research the properties of the group $\mathbf{G}_\Phi$ and the value of the matrix $G(x)$ on the root system of diagonal elements of the matrix $\Phi(x)$.

Let

$$\Phi_\alpha = \operatorname{diag}((x - \alpha)^{s_1}, \ (x - \alpha)^{s_2}, ..., (x - \alpha)^{s_n}) \tag{5.3}$$

be matrix over $F[x]$, moreover $s_i \geqslant 0$, $s_i \leqslant s_{i+1}$, $i = 1, 2, ..., n - 1$. Let

$$s_1 = s_2 = ... = s_{\nu_1}, \ \ s_{\nu_1+1} = s_{\nu_1+2} = ... = s_{\nu_2}, ...,$$
$$s_{\nu_{p-1}+1} = s_{\nu_{p-1}+2} = ... = s_{\nu_p} = s_n,$$

where $s_{\nu_1} < s_{\nu_2} < ... < s_{\nu_p}$. In the group $\mathbf{G}_{\Phi_\alpha}$ we choose an arbitrary matrix $H(x) = \|h_{ij}\|_1^n$. Break it into blocks so that the diagonal blocks are matrices

$$H_{11} = \begin{Vmatrix} h_{11} & ... & h_{1.\nu_1} \\ ... & ... & ... \\ h_{\nu_1.1} & ... & h_{\nu_1.\nu_1} \end{Vmatrix}, \ H_{22} = \begin{Vmatrix} h_{\nu_1+1.\nu_1+1} & ... & h_{\nu_1+1.\nu_2} \\ ... & ... & ... \\ h_{\nu_2.\nu_1+1} & ... & h_{\nu_2\nu_2} \end{Vmatrix}, \ ...,$$

$$H_{pp} = \begin{Vmatrix} h_{\nu_{p-1}+1.\nu_{p-1}+1} & ... & h_{\nu_{p-1}+1.n} \\ ... & ... & ... \\ h_{n.\nu_{p-1}+1} & ... & h_{nn} \end{Vmatrix}.$$

**169**

**Lemma 5.5.** *The inequalities* $\det H_{ii}(\alpha) \neq 0$, $i = 1,\ 2, ..., p$ *are fulfilled.*

**Proof.** By Theorem 2.7, elements of the matrix $H(x)$ have the form $h_{ij} = (x-\alpha)^{s_i - s_j} q_{ij}$ for all $i > j$. Therefore, the matrix $H(\alpha)$ has the form

$$H(\alpha) = \begin{Vmatrix} H_{11}(\alpha) & * & & * \\ \mathbf{0} & H_{22}(\alpha) & & * \\ & & \ddots & \\ \mathbf{0} & \mathbf{0} & & H_{pp}(\alpha) \end{Vmatrix}.$$

Thus,

$$\det H(\alpha) = H_{11}(\alpha) H_{22}(\alpha) ... H_{pp}(\alpha).$$

The matrix $H(x)$ is invertible. So $\det H(\alpha) \neq 0$. Consequently,

$$H_{11}(\alpha) H_{22}(\alpha) ... H_{pp}(\alpha) \neq 0.$$

It follows that $\det H_{ii}(\alpha) \neq 0$, $i = 1, 2, ..., p$. $\qquad\square$

**Lemma 5.6.** *Let $D(x)$ be an $n \times m$ matrix over $F[x]$ and $H(x) \in \mathbf{G}_{\Phi_\alpha}$. Then*

$$M_{H(x)D(x)}(\Phi_\alpha) = F M_{D(x)}(\Phi_\alpha),$$

*where $F$ is an invertible matrix, moreover*

$$\det F = \pm |H_{11}(\alpha)|^{s_{\nu_1}} |H_{22}(\alpha)|^{s_{\nu_2}} ... |H_{pp}(\alpha)|^{s_{\nu_p}}.$$

**Proof.** By definition,

$$M_{H(x)D(x)}(\Phi_\alpha) = \begin{Vmatrix} M_{h_1(x)D(x)}[\alpha^{(s_1)}] \\ M_{h_2(x)D(x)}[\alpha^{(s_2)}] \\ ............ \\ M_{h_n(x)D(x)}[\ \alpha^{(s_n)}] \end{Vmatrix},$$

where $h_i(x)$ is the $i$th row of the matrix $H(x)$. Applying the Leibniz formula, we obtain

$$(h_i(x)D(x))^{(q)} = h_i^{(q)}(x)D(x) + \binom{q}{1} h_i^{(q-1)}(x)D'(x) + ... +$$
$$+ \binom{q}{q-2} h_i''(x)D^{(q-2)}(x) + \binom{q}{1} h_i'(x)D^{(q-1)}(x) + h_i(x)D^{(q)}(x).$$

This yields

$$M_{h_i(x)D(x)}[\alpha^{(s_i)}] =$$

$$= \begin{Vmatrix} h_i(\alpha)D(\alpha) \\ h_i'(\alpha)D(\alpha) + h_i(\alpha)D'(\alpha) \\ ........................................ \\ h_i^{(s_i-1)}(\alpha)D(\alpha) + \binom{s_i-1}{1} h_i^{(s_i-2)}(\alpha)D'(\alpha) + ... + h_i(\alpha)D^{(s_i-1)}(\alpha) \end{Vmatrix} =$$

$$= \left\| \begin{array}{cccc} h_i(\alpha) & \mathbf{0} & ... & \mathbf{0} \\ h_i'(\alpha) & h_i(\alpha) & & \mathbf{0} \\ ... & ... & & \\ h_i^{(s_i-1)}(\alpha) & \begin{pmatrix} s_i-1 \\ 1 \end{pmatrix} h^{(s_i-2)} & ... & h_i(\alpha) \end{array} \right\| M_{D(x)}[\alpha^{(s_i)}] =$$

$$= R_{s_i} M_{D(x)}[\alpha^{(s_i)}], \;\; i = 1, 2, ..., n.$$

Hence,

$$M_{H(x)D(x)}(\Phi_\alpha) = \left\| \begin{array}{cc} R_{s_k} & \mathbf{0} \\ R_{s_{k+1}} & \mathbf{0} \\ ... & \\ R_{s_n} & \end{array} \right\| M_{D(x)}[\alpha^{(s_n)}] = R M_{D(x)}[\alpha^{(s_n)}]. \tag{5.4}$$

Considering the structure of the elements of the blocks of the matrix $H(x)$, which lie below its main diagonal, we see that the matrix $R$ contains zeros columns. Removing them, we get some matrix $L$. Remove in the matrix $M_{D(x)}[\alpha^{(s_n)}] \; M_{D(x)}[\alpha^{(s_n)}]$ the rows which are corresponding to the removed columns of the matrix $R$. We reduce the resulting matrix by a permutation of the rows to the matrix $M_{D(x)}(\Phi_\alpha)$. Thus, equality (5.4) is equivalent to the equality

$$M_{H(x)D(x)}(\Phi_\alpha) = LT M_{D(x)}(\Phi_\alpha),$$

where $T$ is a permutation matrix. For the matrix $L$ there exist permutation matrices $M, N$ such that in the matrix $L_1 = MLN$ the first $s_{\nu_1}$ diagonal blocks will be the matrices $H_{11}(\alpha)$, following $s_{\nu_2}$ is the matrices $H_{22}(\alpha)$, etc. the last $s_{\nu_p}$ is the matrices $H_{pp}(\alpha)$. The structure of the matrix $L_1 = MLN$ is such that

$$\det L_1 = \pm |H_{11}(\alpha)|^{s_{\nu_1}} |H_{22}(\alpha)|^{s_{\nu_2}} ... |H_{pp}(\alpha)|^{s_{\nu_p}}.$$

Since, $L_1 = MLN = M(LT)T^{-1}N$, then $LT = M^{-1}L_1N^{-1}T = F$. Noting that the permutation matrices $M, N, T$ have the determinants $\pm 1$, we get that

$$\det F = \pm |H_{11}(\alpha)|^{s_{\nu_1}} |H_{22}(\alpha)|^{s_{\nu_2}} ... |H_{pp}(\alpha)|^{s_{\nu_p}}.$$

By Lemma 5.5, $\det F \neq 0$. Lemma is proved. $\qquad\square$

**Example 5.4.** Let $\Phi_0(x) = \mathrm{diag}\,(x, x, x, x^2, x^2, x^3)$. The group $\mathbf{G}_\Phi$ consists of all invertible matrices of the form

$$H(x) = \left\| \begin{array}{ccc|cc|c} h_{11}(x) & h_{12}(x) & h_{13}(x) & h_{14}(x) & h_{15}(x) & h_{16}(x) \\ h_{21}(x) & h_{22}(x) & h_{23}(x) & h_{24}(x) & h_{25}(x) & h_{26}(x) \\ h_{31}(x) & h_{32}(x) & h_{33}(x) & h_{34}(x) & h_{35}(x) & h_{36}(x) \\ \hline xh_{41}(x) & xh_{42}(x) & xh_{43}(x) & h_{44}(x) & h_{45}(x) & h_{46}(x) \\ xh_{51}(x) & xh_{52}(x) & xh_{53}(x) & h_{54}(x) & h_{55}(x) & h_{56}(x) \\ \hline x^2h_{61}(x) & x^2h_{62}(x) & x^2h_{63}(x) & xh_{64}(x) & xh_{65}(x) & h_{66}(x) \end{array} \right\|.$$

Zeros columns of the $R$ matrix (see equality (5.4)) are $7, 8, 9, 13\text{-}17$ columns. Removing them and permuting the rows and columns, we get

$$
L_1 = \left\|
\begin{array}{ccc|cc|cc|ccc}
h_{11} & h_{12} & h_{13} & h_{14} & h_{15} & 0 & 0 & 0 & 0 & 0 \\
h_{21} & h_{22} & h_{23} & h_{24} & h_{25} & 0 & 0 & 0 & 0 & 0 \\
h_{31} & h_{32} & h_{33} & h_{34} & h_{35} & 0 & 0 & 0 & 0 & 0 \\
\hline
0 & 0 & 0 & h_{44} & h_{45} & 0 & 0 & h_{46} & 0 & 0 \\
0 & 0 & 0 & h_{54} & h_{55} & 0 & 0 & h_{56} & 0 & 0 \\
\hline
* & * & * & * & * & h_{44} & h_{45} & * & h_{46} & 0 \\
* & * & * & * & * & h_{54} & h_{55} & * & h_{56} & 0 \\
\hline
0 & 0 & 0 & 0 & 0 & 0 & 0 & h_{66} & 0 & 0 \\
0 & 0 & 0 & * & * & 0 & 0 & * & h_{66} & 0 \\
* & * & * & * & * & * & * & * & * & h_{66}
\end{array}
\right\|,
$$

where $h_{ij} = h_{ij}(0)$. It is easy to see that

$$
\det L_1 = \pm
\begin{vmatrix}
h_{11} & h_{12} & h_{13} \\
h_{21} & h_{22} & h_{23} \\
h_{31} & h_{32} & h_{33}
\end{vmatrix}
\begin{vmatrix}
h_{44} & h_{45} \\
h_{54} & h_{55}
\end{vmatrix}^2
h_{66}^3.
\qquad \diamondsuit
$$

We now proceed to the general case. Let

$$
\Phi(x) = \mathrm{diag}(\varphi_1(x), \varphi_2(x), ..., \varphi_n(x))
$$

be a nonsingular $d$-matrix with monic polynomials on the main diagonal and

$$
\det \Phi(x) = (x - \alpha_1)^{t_1}(x - \alpha_2)^{t_2} ... (x - \alpha_q)^{t_q}.
$$

The matrix $\Phi(x)$ can be written as

$$
\Phi(x) = \Phi_{\alpha_1}(x)\Phi_{\alpha_2}(x) ... \Phi_{\alpha_q}(x),
$$

where matrices $\Phi_{\alpha_i}(x)$ have the form (5.3).

**Property 5.4.** *There is an invertible matrix $K$ such that*

$$
K M_{G(x)}(\Phi) = \left\|
\begin{array}{c}
M_{G(x)}(\Phi_{\alpha_1}) \\
M_{G(x)}(\Phi_{\alpha_2}) \\
\text{.................} \\
M_{G(x)}(\Phi_{\alpha_q})
\end{array}
\right\|.
\tag{5.5}
$$

**Proof**. It is easy to see that the matrix $M_{G(x)}(\Phi)$ can be reduced to the form (5.5) by rows permutation. That is, $K$ is a permutation matrix and therefore is invertible. $\qquad \square$

**Theorem 5.10.** *Let $D(x)$ be an $n \times m$ matrix over $F[x]$ and $H(x) \in \mathbf{G}_\Phi$. There is an invertible matrix $N$ such that*

$$M_{H(x)D(x)}(\Phi) = N M_{D(x)}(\Phi).$$

**Proof.** The Property 5.4 implies that there exists an invertible matrix $T$ such that

$$M_{H(x)D(x)}(\Phi) = T \left\| \begin{matrix} M_{H(x)D(x)}(\Phi_{\alpha_1}) \\ M_{H(x)D(x)}(\Phi_{\alpha_2}) \\ \cdots\cdots\cdots \\ M_{H(x)D(x)}(\Phi_{\alpha_q}) \end{matrix} \right\|.$$

By Lemma 5.6, there are invertible matrices $L_{\alpha_i}$ such that

$$M_{H(x)D(x)}(\Phi_{\alpha_i}) = L_{\alpha_i} M_{D(x)}(\Phi_{\alpha_i}),$$

$i = 1, 2, ..., q$. Then

$$M_{H(x)D(x)}(\Phi) = T \left\| \begin{matrix} L_{\alpha_1} M_{D(x)}(\Phi_{\alpha_1}) \\ L_{\alpha_2} M_{D(x)}(\Phi_{\alpha_2}) \\ \cdots\cdots\cdots \\ L_{\alpha_q} M_{D(x)}(\Phi_{\alpha_q}) \end{matrix} \right\| =$$

$$= T \left\| \begin{matrix} L_{\alpha_1} & \mathbf{0} & & \mathbf{0} \\ \mathbf{0} & L_{\alpha_2} & & \mathbf{0} \\ & & \ddots & \\ \mathbf{0} & \mathbf{0} & & L_{\alpha_q} \end{matrix} \right\| \left\| \begin{matrix} M_{D(x)}(\Phi_{\alpha_1}) \\ M_{D(x)}(\Phi_{\alpha_2}) \\ \cdots\cdots\cdots \\ M_{D(x)}(\Phi_{\alpha_q}) \end{matrix} \right\| =$$

$$= \left( T \left\| \begin{matrix} L_{\alpha_1} & \mathbf{0} & & \mathbf{0} \\ \mathbf{0} & L_{\alpha_2} & & \mathbf{0} \\ & & \ddots & \\ \mathbf{0} & \mathbf{0} & & L_{\alpha_q} \end{matrix} \right\| T^{-1} \right) \left( T \left\| \begin{matrix} M_{D(x)}(\Phi_{\alpha_1}) \\ M_{D(x)}(\Phi_{\alpha_2}) \\ \cdots\cdots\cdots \\ M_{D(x)}(\Phi_{\alpha_q}) \end{matrix} \right\| \right) =$$

$$= \underbrace{\left( T \left\| \begin{matrix} L_{\alpha_1} & \mathbf{0} & & \mathbf{0} \\ \mathbf{0} & L_{\alpha_2} & & \mathbf{0} \\ & & \ddots & \\ \mathbf{0} & \mathbf{0} & & L_{\alpha_q} \end{matrix} \right\| T^{-1} \right)}_{N} M_{D(x)}(\Phi).$$

Noting that the matrix $N$ is invertible, we complete the proof of the Theorem. $\square$

**Corollary 5.2.** *If $H(x) \in \mathbf{G}_\Phi$, then*

$$\operatorname{rang} M_{H(x)D(x)}(\Phi) = \operatorname{rang} M_{D(x)}(\Phi). \qquad \square$$

**Corollary 5.3.** *If $H(x) \in \mathbf{G}_\Phi$ and $\det M_{D(x)}(\Phi) \neq 0$, then*

$$\det M_{H(x)D(x)}(\Phi) \neq 0. \qquad \square$$

Denote by

$$J_{\alpha^{(k)}} = \begin{Vmatrix} \alpha & & & & & 0 \\ 1 & \alpha & & & & \\ & 2 & \alpha & & & \\ & & & \ddots & \ddots & \\ 0 & & & & k-1 & \alpha \end{Vmatrix},$$

and

$$J(\varphi) = J_{\alpha_1^{(k_1)}} \oplus ... \oplus J_{\alpha_l^{(k_l)}},$$

where $\varphi(x) = (x - \alpha_1)^{k_1} ... (x - \alpha_l)^{k_l}$.

Also denote by

$$J(\Phi) = J(\varphi_1) \oplus ... \oplus J(\varphi_n),$$

where $\Phi(x) = \mathrm{diag}(\varphi_1(x), ..., \varphi_n(x))$. If $\deg \varphi_i(x) = 0$, then the matrix $J(\varphi_i)$ is empty.

**Property 5.5.** *The equality*

$$M_{xG(x)}(\Phi) = J(\Phi) M_{G(x)}(\Phi)$$

*is fulfilled.*

**Proof.** Let

$$\varphi_i(x) = (x - \alpha_1)^{k_{i1}} ... (x - \alpha_t)^{k_{it}},$$

and $g_i(x)$ is the $i$th row of the matrix $G(x)$, $i = 1, ..., n$. Then

$$M_{xg_i(x)}(x - \alpha_j)^{k_{ij}} = M_{xg_i(x)}[\alpha_j^{(k_{ij})}] =$$

$$= \begin{Vmatrix} \alpha_j g_i(\alpha_j) \\ g_i(\alpha_j) + \alpha_j g_i'(\alpha_j) \\ 2g_i'(\alpha_j) + \alpha_j g_i''(\alpha_j) \\ \cdots\cdots\cdots\cdots\cdots\cdots\cdots\cdots \\ (k_{ij} - 1)g_i^{(k_{ij}-2)}(\alpha_j) + \alpha_j g_i^{(k_{ij}-1)}(\alpha_j) \end{Vmatrix} =$$

$$= \begin{Vmatrix} \alpha_j & & & & 0 \\ 1 & \alpha_j & & & \\ & 2 & \alpha_j & & \\ & & \ddots & \ddots & \\ 0 & & & k_{ij}-1 & \alpha_j \end{Vmatrix} \begin{Vmatrix} g_i(\alpha_j) \\ g_i'(\alpha_j) \\ g_i''(\alpha_j) \\ \cdots \\ g_i^{(k_{ij}-1)}(\alpha_j) \end{Vmatrix} =$$

$$= J_{\alpha_j^{(k_{ij})}} M_{g_i(x)}[\alpha_j^{(k_{ij})}].$$

Therefore,

$$M_{xg_i(x)}(\varphi_i) = \begin{Vmatrix} M_{xg_i(x)}[\alpha_1^{(k_{i1})}] \\ \cdots\cdots\cdots\cdots\cdots \\ M_{xg_i(x)}[\alpha_t^{(k_{it})}] \end{Vmatrix} = \begin{Vmatrix} J_{\alpha_1^{(k_{i1})}} M_{g_i(x)}[\alpha_1^{(k_{i1})}] \\ \cdots\cdots\cdots\cdots\cdots\cdots \\ J_{\alpha_t^{(k_{it})}} M_{g_i(x)}[\alpha_t^{(k_{it})}] \end{Vmatrix} =$$

$$= J(\varphi_i) M_{g_i(x)}(\varphi_i).$$

Thus

$$M_{xG(x)}(\Phi) = \left\| \begin{matrix} M_{xg_1(x)}(\varphi_1) \\ \dots\dots\dots \\ M_{xg_n(x)}(\varphi_n) \end{matrix} \right\| = \left\| \begin{matrix} J(\varphi_1)M_{g_1(x)}(\varphi_1) \\ \dots\dots\dots\dots \\ J(\varphi_n)M_{g_n(x)}(\varphi_n) \end{matrix} \right\| = J(\Phi)M_{G(x)}(\Phi).$$

The proof is complete. □

## 5.4. Regularization of polynomial matrices

Let $A(x)$ be a nonsingular $n \times n$ matrix over $F[x]$. Let us write one as a matrix polynomial:

$$A(x) = A_k x^k + A_{k-1}x^{k-1} + ... + A_0.$$

A matrix $A(x)$ is called a **monic matrix** if $A_k$ is an identity matrix. We say that $A(x)$ is **right regularized** if there exists an invertible matrix $U(x)$ such that

$$A(x)U(x) = Ix^s - D_{s-1}x^{s-1} - ... - D_0.$$

**Lemma 5.7.** *The polynomial matrix $A(x)$ is right regularized uniquely.*

**Proof.** Suppose that there are invertible matrices $U(x)$ and $V(x)$ such that

$$A(x)U(x) = Ix^s - D_{s-1}x^{s-1} - ... - D_0 = A_1(x),$$
$$A(x)V(x) = Ix^k - M_{k-1}x^{k-1} - ... - M_0 = A_2(x).$$

Then $A_2(x) = A_1(x)L(x)$, where $L(x) = U^{-1}(x)V(x)$. Since $L(x)$ is an invertible matrix then $\det L(x)$ is a nonzero element of $F$. Hence, $\deg \det L(x) = 0$. Thus

$$\deg \det A_2(x) = \deg \det(A_1(x)L(x)) =$$
$$= \deg \det A_1(x) + \deg \det L(x) = \deg \det A_1(x).$$

Since $\deg \det A_1(x) = ns$ and $\deg \det A_2(x) = nk$ we conclude that $k = s$. It means that

$$Ix^s - M_{s-1}x^{s-1} - ... - M_0 = (Ix^s - D_{s-1}x^{s-1} - ... - D_0)L(x). \quad (5.6)$$

By writing down the matrix $L(x)$ as a matrix polynomial, multiplying and equating the matrix coefficients in equality (5.6) we get that $L(x) = I$. Consequently, $A_2(x) = A_1(x)$. Lemma is proved. □

There are invertible matrices $P(x)$ and $Q(x)$ such that

$$P(x)A(x)Q(x) = \mathrm{diag}(\varphi_1(x), \varphi_2(x), ..., \varphi_n(x)) = \Phi(x), \quad (5.7)$$

where $\Phi(x)$ is the Smith form of the matrix $A(x)$.

**Theorem 5.11.** *The matrix polynomial $A(x)$ is right regularized if and only if*

*1) $\det A(x) = ns$,*

*2) the equation*

$$M_{P(x)\|Ix^{s-1} \ \dots \ Ix \ I\|}(\Phi)X = M_{P(x)x^s}(\Phi) \tag{5.8}$$

*is solvable.*

**Proof. Necessity**. Suppose that there is an invertible matrix $U(x)$ such that

$$A(x)U(x) = Ix^s - D_{s-1}x^{s-1} - \dots - D_0 = K(x). \tag{5.9}$$

Since

$$\deg \det(A(x)U(x)) = \deg \det(Ix^s - D_{s-1}x^{s-1} - \dots - D_0) = ns$$

and

$$\deg \det(A(x)U(x)) = \deg \det A(x) + \deg \det U(x),$$

where $\deg \det U(x) = 0$, then $\deg \det A(x) = ns$.

Using (5.9) and (5.7), we have

$$P(x)A(x)U(x) = P(x)(Ix^s - D_{s-1}x^{s-1} - \dots - D_0) = \Phi(x)Q^{-1}(x)U(x).$$

By Property 5.3,

$$M_{\Phi(x)Q^{-1}(x)U(x)}(\Phi) = \mathbf{0}.$$

Then

$$M_{P(x)(Ix^s - D_{s-1}x^{s-1} - \dots - D_0)}(\Phi) = \mathbf{0}. \tag{5.10}$$

Taking into account Properties 5.1, 5.2, we receive

$$M_{P(x)(Ix^s - D_{s-1}x^{s-1} - \dots - D_0)}(\Phi) =$$

$$= M_{P(x)x^s}(\Phi) - M_{P(x)x^{s-1}D_{s-1}}(\Phi) - \dots - M_{P(x)D_0}(\Phi) =$$

$$= M_{P(x)x^s}(\Phi) - M_{P(x)x^{s-1}}(\Phi)D_{s-1} - \dots - M_{P(x)}(\Phi)D_0 =$$

$$= M_{P(x)x^s}(\Phi) - \left\| M_{P(x)x^{s-1}}(\Phi) \ \dots \ M_{P(x)x}(\Phi) \ M_{P(x)}(\Phi) \right\| \times$$

$$\times \left\| \begin{matrix} D_{s-1} \\ \dots \\ D_1 \\ D_0 \end{matrix} \right\| = M_{P(x)x^s}(\Phi) - M_{P(x)\|Ix^{s-1} \ \dots \ Ix \ I\|}(\Phi) \left\| \begin{matrix} D_{s-1} \\ \dots \\ D_1 \\ D_0 \end{matrix} \right\|.$$

In view of (5.10), we have

$$M_{P(x)\|Ix^{s-1} \ \dots \ Ix \ I\|}(\Phi) \left\| \begin{matrix} D_{s-1} \\ \dots \\ D_1 \\ D_0 \end{matrix} \right\| = M_{P(x)x^s}(\Phi). \tag{5.11}$$

Consequently, the equation (5.8) is solvable.

To prove necessity, it remains to show that the condition 2) does not depend on the choice of the transforming matrix $P(x)$. Let $P_1(x)$ be another left transforming matrix of the matrix $A(x)$. It means that there exists the matrix $H(x) \in \mathbf{G}_\Phi$ such that $P_1(x) = H(x)P(x)$. By virtue of Theorem 5.10, there is an invertible matrix $N$ such that

$$M_{H(x)D(x)}(\Phi) = NM_{D(x)}(\Phi),$$

where $D(x)$ is an arbitrary $n \times m$ matrix. That is the equation

$$M_{P_1(x)\|Ix^{s-1} \ \dots \ Ix \ I\|}(\Phi)X = M_{P_1(x)x^s}(\Phi)$$

is equivalent to the equation

$$NM_{P(x)\|Ix^{s-1} \ \dots \ Ix \ I\|}(\Phi)X = NM_{P(x)x^s}(\Phi)$$

which in turn is equivalent to the solvable equation (5.8).

**Sufficiency.** Suppose that the equation (5.8) is solvable and $\left\|\begin{array}{c} M_{s-1} \\ \dots \\ M_1 \\ M_0 \end{array}\right\|$ is the solution. Having thought the opposite of what we have just done, we get

$$M_{P(x)(Ix^s - M_{s-1}x^{s-1} - \dots - M_0)}(\Phi) = \mathbf{0}.$$

By Property 5.3,

$$P(x)(Ix^s - M_{s-1}x^{s-1} - \dots - M_0) = P(x)K(x) = \Phi(x)T(x). \qquad (5.12)$$

Since

$$ns = \deg \det \Phi(x) = \deg \det(P(x)K(x)) =$$

$$= \deg \det(\Phi(x)T(x)) = \deg \det \Phi(x) + \deg \ \det T(x),$$

then $\deg \det T(x) = 0$. The matrix $P(x)K(x)$ is nonsingular. It follows that $T(x)$ is an invertible matrix.

Right-multiply (5.12) by an invertible matrix $L(x) = T^{-1}(x)Q^{-1}(x)$, we get

$$P(x)K(x)T^{-1}(x)Q^{-1}(x) = \Phi(x)Q^{-1}(x).$$

Therefore,

$$K(x)L(x) = P^{-1}(x)\Phi(x)Q^{-1}(x) = A(x).$$

Since $A(x)L^{-1}(x) = K(x)$ is a monic matrix, then $A(x)$ is right regularized. The theorem is proved. □

This Theorem can be reformulated as follows.

**Theorem 5.12.** *The matrix polynomial $A(x)$ is right regularized if and only if*

*1) $\deg \det A(x) = ns$,*

*2) $\det M_{P(x)\|I \ Ix \ \dots \ Ix^{s-1}\|}(\Phi) \neq 0$.*

**Proof.** The necessity for condition 1) is proved in Theorem 5.11.

By Lemma 5.7, the matrix polynomial $A(x)$ is right regularized uniquely. It means that equation (5.8) has a single solution. And this is equivalent to fulfilling of the condition 2). □

**Theorem 5.13.** *Suppose that there is an invertible matrix $U(x)$ such that*

$$A(x)U(x) = Ix^s - D_{s-1}x^{s-1} - ... - D_0.$$

*Then the coefficient-matrices $D_i$ are obtained from equality*

$$\left\| \begin{matrix} D_{s-1} \\ ... \\ D_1 \\ D_0 \end{matrix} \right\| = M_{P(x)\|Ix^{s-1} \ ... \ Ix \ I\|}^{-1}(\Phi)M_{P(x)x^s}(\Phi). \tag{5.13}$$

**Proof.** The coefficient-matrices $D_i$ of monic polynomial $A(x)U(x)$ satisfy the equality (5.11), see proving of necessity for Theorem 5.11.

By virtue of Theorem 5.12, the matrix $M_{P(x)\|I \ Ix \ ... \ Ix^{s-1}\|}(\Phi)$ is invertible. Therefore, the equality (5.13) is fulfilled. □

**Corollary 5.4.** *Suppose that*

$$A(x)U(x) = Ix - D.$$

*Then the matrix $D$ is obtained from the equality*

$$D = T^{-1}J(\Phi)T,$$

*where $T = M_{P(x)}(\Phi)$.*

**Proof.** According to Theorem 5.13, and Property 5.5, we get

$$D = M_{P(x)}^{-1}(\Phi)M_{P(x)x}(\Phi) = M_{P(x)}^{-1}(\Phi)J(\Phi)M_{P(x)}(\Phi),$$

The proof is complete. □

## 5.5. One method for finding Jordan normal form

Now apply the obtained results to find the Jordan normal forms (or Jordan form) for a matrices over $F$. Let

$$f(x) = (x - \alpha_1)^{k_1} ... (x - \alpha_m)^{k_m}.$$

Denote by

$$\tilde{M}_{g(x)}(f(x)) = \left\| \begin{matrix} \tilde{M}_{g(x)}[\alpha_1^{(k_1)}] \\ ... \\ \tilde{M}_{g(x)}[\alpha_m^{(k_m)}] \end{matrix} \right\|,$$

where

$$\tilde{M}_{g(x)}[\alpha_i^{(k_i)}] = \left\| \begin{array}{c} 1g(\alpha_i) \\ \dfrac{1}{1!}g'(\alpha_i) \\ \dfrac{1}{2!}g''(\alpha_i) \\ ... \\ \dfrac{1}{(k_i-1)!}g^{(k_i-1)}(\alpha_i) \end{array} \right\|.$$

It is easy to see that

$$\tilde{M}_{g(x)}[\alpha_i^{(k_i)}] = T_{\alpha_i^{(k_i)}}M_{g(x)}[\alpha_i^{(k_i)}],$$

where

$$T_{\alpha_i^{(k_i)}} = \operatorname{diag}\left(1, \frac{1}{1!}, \frac{1}{2!}, ..., \frac{1}{(k_i-1)!}\right).$$

Let $\Phi(x) = \operatorname{diag}(\varphi_1(x),...,\varphi_n(x))$. Denote by

$$\tilde{M}_{G(x)}(\Phi) = \left\| \begin{array}{c} \tilde{M}_{g_1(x)}(\varphi_1) \\ ... ... ... \\ \tilde{M}_{g_n(x)}(\varphi_n) \end{array} \right\|.$$

It is obvious that

$$\tilde{M}_{G(x)}(\Phi) = T_\Phi M_{G(x)}(\Phi),$$

where $T_\Phi$ is a diagonal matrix, which is the direct sum of all matrices $T_{\alpha_{ij}}^{(k_{ij})}$.

**Theorem 5.14.** *Let $A \in M_n(F)$ and $A(x) = Ix - A$ is the characteristic matrix with the Smith form $\Phi(x)$ and $P(x) \in \mathbf{P}_{A(x)}$. Then the matrix*

$$J_*(\Phi) = \tilde{M}_{P(x)}(\Phi)\, A\tilde{M}_{P(x)}^{-1}(\Phi)$$

*is a Jordan form of the matrix A.*

**Proof.** Obviously, the matrix $A(x) = Ex - A$ is right regularized. According to Theorem 5.12,

$$\det M_{P(x)}(\Phi) \neq 0 \;\;\Rightarrow\;\; \det \tilde{M}_{P(x)}(\Phi) \neq 0.$$

Since the matrix $A(x)$ is right regularized uniquely, then on the basis of the Corollary 5.4,

$$A = M_{P(x)}^{-1}(\Phi)J(\Phi)M_{P(x)}(\Phi). \tag{5.14}$$

By direct calculations we check that

$$T_{\alpha^{(k)}} J_{\alpha^{(k)}} T_{\alpha^{(k)}}^{-1} =$$

**179**

$$= \begin{Vmatrix} 1 & & & & 0 \\ & \dfrac{1}{1!} & & & \\ & & \dfrac{1}{2!} & & \\ & & & \ddots & \\ 0 & & & & \dfrac{1}{(k_i-1)!} \end{Vmatrix} \begin{Vmatrix} \alpha & & & & 0 \\ 1 & \alpha & & & \\ & 2 & \alpha & & \\ & & \ddots & \ddots & \\ 0 & & & k-1 & \alpha \end{Vmatrix} \times$$

$$\times \begin{Vmatrix} 1 & & & & 0 \\ & 1! & & & \\ & & 2! & & \\ & & & \ddots & \\ 0 & & & & (k_i-1)! \end{Vmatrix} = \begin{Vmatrix} \alpha & & & & 0 \\ 1 & \alpha & & & \\ & 1 & \alpha & & \\ & & \ddots & \ddots & \\ 0 & & & 1 & \alpha \end{Vmatrix}$$

is a Jordan block corresponding to the elementary divisor $(x-\alpha)^k$. Therefore, $T_\Phi J(\Phi) T_\Phi^{-1}$ is the Jordan matrix.

We write the equality (5.14) in the form

$$A = M_{P(x)}^{-1}(\Phi) J(\Phi) M_{P(x)}(\Phi) =$$
$$= \left( M_{P(x)}^{-1}(\Phi) T_\Phi^{-1} \right) \left( T_\Phi J(\Phi) T_\Phi^{-1} \right) \left( T_\Phi M_{P(x)}(\Phi) \right) =$$
$$= \left( T_\Phi M_{P(x)}(\Phi) \right)^{-1} J_*(\Phi) \left( T_\Phi M_{P(x)}(\Phi) \right) = \tilde{M}_{P(x)}^{-1}(\Phi) J_*(\Phi) \tilde{M}_{P(x)}(\Phi).$$

The theorem is proved. $\qquad\square$

**Example 5.5.** We reduce the matrix

$$A = \begin{Vmatrix} 1 & -3 & 3 \\ -2 & -6 & 13 \\ -1 & -4 & 8 \end{Vmatrix}$$

to the Jordan form.

The Smith form of the characteristic matrix $A(x) = Ex - A$ is

$$P(x)(Ex - A)Q(x) = \mathrm{diag}(1, 1, (x-1)^3) = \Phi(x),$$

where

$$P(x) = \begin{Vmatrix} 0 & 0 & 1 \\ -1 & -4 & x+7 \\ x-2 & 4x-7 & -x^2-5x+12 \end{Vmatrix},$$

$$Q(x) = \begin{Vmatrix} 1 & -4 & -4x^2+5x-4 \\ 0 & 1 & x^2-x-1 \\ 0 & 0 & -1 \end{Vmatrix}.$$

In the second step, we will construct its transforming matrix:

$$\tilde{M}_{P(x)}(\Phi) = \tilde{M}_{\| x-2 \;\; 4x-7 \;\; -x^2-5x+12 \|} (x-1)^3 = \begin{Vmatrix} -1 & -3 & 6 \\ 1 & 4 & -7 \\ 0 & 0 & -1 \end{Vmatrix} = T.$$

Since,

$$T^{-1} = \begin{Vmatrix} -4 & -3 & -3 \\ 1 & 1 & -1 \\ 0 & 0 & -1 \end{Vmatrix},$$

we get

$$TAT^{-1} = \begin{Vmatrix} 1 & 0 & 0 \\ 1 & 1 & 0 \\ 0 & 1 & 1 \end{Vmatrix}$$

is the Jordan form of the matrix $A$. $\diamondsuit$

## 5.6. Determining matrix and its properties

To describe the unassociated divisors of matrices over an elementary divisor rings we used the Kazimirskii set $\mathbf{V}(E, \Phi)$. The determining matrix plays a similar role in finding the monic divisors of polynomial matrices [57]. This subsection is devoted to the study of its properties.

Let

$$E(x) = \mathrm{diag}(\varepsilon_1(x), ..., \varepsilon_n(x)), \quad \Phi(x) = \mathrm{diag}\,(\varphi_1(x), ..., \varphi_n(x))$$

are nonsingular $d$-matrices over $F[x]$, moreover $\Phi(x)|E(x)$, $\deg \det \Phi(x) = nr$.

The **determining matrix** is called the matrix

$$V(E, \Phi) = \begin{Vmatrix} 1 & 0 & ... & 0 & 0 \\ \dfrac{\varphi_2}{(\varphi_2, \varepsilon_1)}k_{21} & 1 & ... & 0 & 0 \\ ... & ... & ... & ... & ... \\ \dfrac{\varphi_n}{(\varphi_n, \varepsilon_1)}k_{n1} & \dfrac{\varphi_n}{(\varphi_n, \varepsilon_2)}k_{n2} & ... & \dfrac{\varphi_n}{(\varphi_n, \varepsilon_{n-1})}k_{n,n-1} & 1 \end{Vmatrix}, \quad (5.15)$$

where

$$k_{ij} = \begin{cases} 0, & (\varphi_i, \varepsilon_j) = \varphi_j, \\ k_{ij0} + k_{ij1}x + ... + k_{ijh_{ij}}x^{h_{ij}}, & (\varphi_i, \varepsilon_j) \neq \varphi_j, \end{cases}$$

$$h_{ij} = \deg \frac{(\varphi_i, \varepsilon_j)}{\varphi_j} - 1, \ i = 2, 3, ..., n, \ j = 1, ..., n-1, \ i > j,$$

where $k_{ijl}$ are parameters, $i = 2, 3, ..., n$, $j = 1, ..., n-1$.

Denote by $\dfrac{f}{g}(\alpha)$ the value of the fraction of polynomials $f(x)$, $g(x)$ on the element $\alpha \in F$.

**Property 5.6.** *If*

$$\frac{\varphi_i}{\varphi_j}(\alpha) = 0, \quad \frac{\varepsilon_i}{\varepsilon_j}(\alpha) \neq 0, \ i > j,$$

*then*

$$\frac{\varphi_i}{(\varphi_i, \varepsilon_j)}(\alpha)k_{ij}(\alpha) \neq 0.$$

**181**

**Proof.** Suppose that $k_{ij}(x) = 0$. Therefore, $(\varphi_i, \varepsilon_j) = \varphi_j$. The equality

$$\frac{\varepsilon_i}{\varepsilon_j} = \frac{\varepsilon_i \varphi_i (\varphi_i, \varepsilon_j)}{\varepsilon_j \varphi_i \varphi_j} = \frac{\varphi_i}{\varphi_j} \frac{(\varepsilon_i \varphi_i, \varepsilon_i \varepsilon_j)}{\varphi_i \varepsilon_j}$$

is fulfilled. It follows that $\dfrac{\varepsilon_i}{\varepsilon_j}(\alpha) = 0$, which contradicts the condition of assertion.

Since $\dfrac{\varphi_i}{(\varphi_i, \varepsilon_j)}\left|\dfrac{\varepsilon_i}{\varepsilon_j}\right.$ and $\dfrac{\varepsilon_i}{\varepsilon_j}(\alpha) \neq 0$, we get $\dfrac{\varphi_i}{(\varphi_i, \varepsilon_j)}(\alpha) \neq 0$. Consequently,

$$\frac{\varphi_i}{(\varphi_i, \varepsilon_j)}(\alpha) k_{ij}(\alpha) \neq 0. \qquad \square$$

**Property 5.7.** *If $k_{ij}(x) = 0, i > j$, then*

$$k_{ij}(x) = k_{i-1.j}(x) = ... = k_{j+1.j}(x) = 0. \qquad (5.16)$$

**Proof.** According to the assumption $(\varphi_i, \varepsilon_j) = \varphi_j$. Since $\varphi_t | \varphi_i$, $t = i - 1, i - 2, ..., j + 1$, then $(\varphi_t, \varepsilon_j) | (\varphi_i, \varepsilon_j)$. On the other hand $\varphi_j | \varepsilon_j$ and $\varphi_i | \varphi_t$. That is $\varphi_j | (\varphi_t, \varepsilon_q)$. Hence, $(\varphi_i, \varepsilon_j) | (\varphi_t, \varepsilon_q)$. Therefore, $(\varphi_t, \varepsilon_q) = \varphi_q$, $t = i - 1$, $i - 2, ..., j + 1$. And this is equivalent to the fulfillment of equality (5.16). $\square$

Consider the product of matrices

$$V(\mathrm{E}, \Phi) T(x) = U(x),$$

where $T(x) \in \mathbf{G}_{\mathrm{E}}$. Denote by $U_i(x)$ the submatrix of the $U(x)$ obtained by removing of the first $i$th rows and the first $i$th columns of the matrix $U(x)$, $i = 1, ..., n - 1$. According to the Binet-Cauchy formula

$$\det U_i(x) = \sum_j |V_{ij}(\mathrm{E}, \Phi)||T_{ji}(x)| + \det T_i(x),$$

where $\sum_j |V_{ij}(\mathrm{E}, \Phi)||T_{ji}(x)|$ are the sum of the products of all possible of the $(n - i)$th order minors of the matrix $V_{ij}(\mathrm{E}, \Phi)$ built on the last rows of the matrix $V(\mathrm{E}, \Phi)$, except minor of a unitriangular matrix for the corresponding minors of the same order of the matrix $T(x)$. $T_i(x)$ is obtained by the crossing of the first $i$th rows and the first $i$th columns of the matrix $T(x)$, $i = 1, ..., n-1$. Denote by $V_{ij}(\mathrm{E}, \Phi, \alpha)$ the matrix obtained from $V_{ij}(\mathrm{E}, \Phi)$, by carrying out the change $x$ for $\alpha \in F$.

**Lemma 5.8.** *In order that*

$$\det U_i(\alpha) = \sum_j |V_{ij}(\mathrm{E}, \Phi, \alpha)||T_{ji}(\alpha)| + \det T_i(\alpha) = 0, \qquad (5.17)$$

*it is necessary and sufficient that all additions of this sum be equal to zero.*

**Proof. Sufficiency** is obvious.

**Necessity.** To prove that it is sufficient to note that minor $|V_{ij}(\mathrm{E}, \Phi, \alpha)|$ are the sum of the product of elements in $F$ and parameters $k_{ijl}, \ i = 2, 3, ..., n,$ $j = 1, ..., n - 1$. Herewith, the set of parameters in each such minor is not repeated in any other. $\quad\square$

**Theorem 5.15.** *The equality*

$$\left( \frac{\varphi_{i+1}}{\varphi_i}(x), \det U_i(x) \right) = 1, \ i = 1, ..., n - 1$$

*is fulfilled.*

**Proof.** Suppose for some $i = r$ we have

$$\left( \frac{\varphi_{r+1}}{\varphi_r}(x), \det U_r(x) \right) = \delta(x).$$

Let $\alpha$ be the root of the polynomial $\delta(x)$. Denote by

$$E_r(x) = \left\| \begin{matrix} \frac{\varepsilon_{r+1}}{\varepsilon_1}(x) & ... & \frac{\varepsilon_{r+1}}{\varepsilon_r}(x) \\ \frac{\varepsilon_{r+2}}{\varepsilon_1}(x) & ... & \frac{\varepsilon_{r+2}}{\varepsilon_r}(x) \\ ... & ... & ... \\ \frac{\varepsilon_n}{\varepsilon_1}(x) & ... & \frac{\varepsilon_n}{\varepsilon_r}(x) \end{matrix} \right\|, \quad F_r(x) = \left\| \begin{matrix} \frac{\varphi_{r+1}}{\varphi_1}(x) & ... & \frac{\varphi_{r+1}}{\varphi_r}(x) \\ \frac{\varphi_{r+2}}{\varphi_1}(x) & ... & \frac{\varphi_{r+2}}{\varepsilon_r}(x) \\ ... & ... & ... \\ \frac{\varphi_n}{\varphi_1}(x) & ... & \frac{\varphi_n}{\varphi_r}(x) \end{matrix} \right\|.$$

Suppose that $\frac{\varepsilon_{r+1}}{\varepsilon_r}(\alpha) = 0$. By Property 2.4, $E_r(\alpha) = \mathbf{0}$. That is, the matrix $T(\alpha)$ has the form

$$T(\alpha) = \left\| \begin{matrix} * & * \\ \mathbf{0} & T_r(\alpha) \end{matrix} \right\|.$$

By Lemma 5.8, the equality (5.17) is fulfilled if and only if $|T_r(\alpha)| = 0$. Hence, $|T(\alpha)| = 0$, which contradicts the invertibility of the matrix $T(x)$.

Suppose that the matrix $E_r(\alpha)$ contains no zeros. According to the assumption $\frac{\varphi_{r+1}}{\varphi_r}(\alpha) = 0$. By virtue of Property 2.4, we get $F_r(\alpha) = \mathbf{0}$. In view of Property 5.6, the matrix

$$\left\| \begin{matrix} \mu_{r+1.1}(\alpha) & ... & \mu_{r+1.r}(\alpha) \\ \mu_{r+2.1}(\alpha) & ... & \mu_{r+2.r}(\alpha) \\ ... & ... & ... \\ \mu_{n1}(\alpha) & ... & \mu_{nr}(\alpha) \end{matrix} \right\|, \quad \mu_{ij} = \frac{\varphi_i}{(\varphi_i, \varepsilon_j)}(\alpha) k_{ij}(\alpha),$$

contains no zero elements. By Lemma 5.8, the equality (5.17) is fulfilled if and only if all minors of $T_{jr}(\alpha)$ and $T_r(\alpha)$ are equal to zeros. It follows that $|T(\alpha)| = 0$ that is a contradiction.

Suppose that the matrix $E_r(\alpha)$ contains zeros. Let

$$\frac{\varepsilon_p}{\varepsilon_q}(\alpha) = 0, \quad \frac{\varepsilon_p}{\varepsilon_{q+1}}(\alpha) \neq 0, \quad \frac{\varepsilon_{p-1}}{\varepsilon_q}(\alpha) \neq 0.$$

By Corollary 2.4, the fraction $\dfrac{\varepsilon_p}{\varepsilon_q}$ lies on the first subdiagonal of the matrix $T(x)$. The matrix $E_r(\alpha)$ contains only one such element, namely $\dfrac{\varepsilon_{r+1}}{\varepsilon_r}(\alpha)$. That is $(p,q) = (r+1, r)$. Hence, $\dfrac{\varepsilon_{r+1}}{\varepsilon_r}(\alpha) = 0$, but this case we have considered above. It means that the zeros in the matrix $E_r(\alpha)$ can be accommodated in this way

$$a)\ E_r(\alpha) = \left\| \mathbf{0} \quad \divideontimes \right\|, \ \ b)\ E_r(\alpha) = \left\| \begin{matrix} \divideontimes \\ \mathbf{0} \end{matrix} \right\|, \ \ c)\ E_r(\alpha) = \left\| \begin{matrix} \mathbf{0} & \divideontimes \\ \mathbf{0} & \mathbf{0} \end{matrix} \right\|.$$

Consider the case $a)$. Without loss of generality of the proof, we implement one on the seventh-order matrices.

Let $n = 7$, $r = 3$ and the matrix $E_3(\alpha)$ has the form

$$E_3(\alpha) = \left\| \begin{matrix} 0 & * & * \\ 0 & * & * \\ 0 & * & * \\ 0 & * & * \end{matrix} \right\|,$$

where all $*$ are not zeros. Based on Property 2.5, the matrix $T(\alpha)$ has the form

$$T(\alpha) = \left\| \begin{matrix} \cdot & \cdot & \cdot & \cdot & \cdot & \cdot & \cdot \\ \hline 0 & \cdot & \cdot & & & & \\ 0 & \cdot & \cdot & & & & \\ 0 & t_{42} & t_{43} & & W & & \\ 0 & t_{52} & t_{53} & & & & \\ 0 & t_{62} & t_{63} & & & & \\ 0 & t_{72} & t_{73} & & & & \end{matrix} \right\|, \ \ t_{ij} \neq 0.$$

Since $F_3(\alpha) = \mathbf{0}$, by Property 5.6, the matrix $V(\mathrm{E}, \Phi, \alpha)$ has the form

$$V(\mathrm{E}, \Phi, \alpha) = \left\| \begin{matrix} 1 & & & & & & & \mathbf{0} \\ \cdot & 1 & & & & & & \\ \cdot & \cdot & 1 & & & & & \\ \hline \cdot & \mu_{42} & \mu_{43} & 1 & & & & \\ \cdot & \mu_{52} & \mu_{53} & \cdot & 1 & & & \\ \cdot & \mu_{62} & \mu_{63} & \cdot & \cdot & 1 & & \\ \cdot & \mu_{72} & \mu_{73} & \cdot & \cdot & \cdot & 1 & \end{matrix} \right\|, \ \ \mu_{\mathrm{ij}} \neq 0.$$

All the maximal order minors that are built on the last four rows and the last six columns of the matrix $V(\mathrm{E}, \Phi, \alpha)$, are not zeros and contain variables that are missing from all others. Therefore, equality (5.17) is fulfilled when all maximum order minors of the matrix $W$ is equal to zero. Then, by virtue of Corollary 3.2, $|T(\alpha)| = 0$ is a contradiction.

Case $b$). Let

$$E_3(\alpha) = \begin{Vmatrix} * & * & * \\ * & * & * \\ 0 & 0 & 0 \\ 0 & 0 & 0 \end{Vmatrix}.$$

Based on Property 2.5, the matrix $T(\alpha)$ has the form

$$T(\alpha) = \begin{Vmatrix} \cdot & \cdot & \cdot & \cdot & \cdot & \vline & \cdot & \cdot \\ \cdot & \cdot & \cdot & \cdot & \cdot & \vline & \cdot & \cdot \\ \cdot & \cdot & \cdot & \cdot & \cdot & \vline & \cdot & \cdot \\ \hline \cdot & \cdot & \cdot & \cdot & \cdot & \vline & \cdot \\ \cdot & \cdot & \cdot & \cdot & \cdot & \vline & \cdot & \cdot \\ 0 & 0 & 0 & 0 & 0 & \vline & \cdot & \cdot \\ 0 & 0 & 0 & 0 & 0 & \vline & \cdot & \cdot \end{Vmatrix} = \begin{Vmatrix} t_{ij} \end{Vmatrix}_1^7.$$

The matrix $V(\mathrm{E}, \Phi, \alpha)$ has the form

$$V(\mathrm{E}, \Phi, \alpha) = \begin{Vmatrix} 1 & & & & & & \mathbf{0} & \\ \cdot & 1 & & & & & & \\ \cdot & \cdot & 1 & & & & & \\ \hline \mu_{41} & \mu_{42} & \mu_{43} & 1 & & & & \\ \mu_{51} & \mu_{52} & \mu_{53} & \cdot & 1 & & & \\ \cdot & \cdot & \cdot & \cdot & \cdot & 1 & & \\ \cdot & \cdot & \cdot & \cdot & \cdot & \cdot & 1 \end{Vmatrix}, \quad \mu_{\mathrm{ij}} \neq 0.$$

All maximal order minors of the matrix

$$\begin{Vmatrix} \mu_{41} & \mu_{42} & \mu_{43} & 1 & 0 & \vline & 0 & 0 \\ \mu_{51} & \mu_{52} & \mu_{53} & \cdot & 1 & \vline & 0 & 0 \\ \cdot & \cdot & \cdot & \cdot & \cdot & \vline & 1 & 0 \\ \cdot & \cdot & \cdot & \cdot & \cdot & \vline & \cdot & 1 \end{Vmatrix},$$

containing its last two columns are nonzeros. Then equality (5.17) is satisfied when all corresponding maximal order minors of the matrix $T(\alpha)$ are zeros. All such minors have the form

$$\det \begin{Vmatrix} t_{i4} & t_{i5} & \vline & t_{i6} & t_{i7} \\ t_{j4} & t_{j5} & \vline & t_{j6} & t_{j7} \\ \hline 0 & 0 & \vline & t_{66} & t_{67} \\ 0 & 0 & \vline & t_{76} & t_{77} \end{Vmatrix}, \quad 1 \leq i < j \leq 5.$$

**185**

If

$$\det \begin{Vmatrix} t_{66} & t_{67} \\ t_{76} & t_{77} \end{Vmatrix} = 0,$$

then $|T(\alpha)| = 0$ is a contradiction. Let

$$\det \begin{Vmatrix} t_{66} & t_{67} \\ t_{76} & t_{77} \end{Vmatrix} \neq 0.$$

Then all the second-order minors of the $T(\alpha)$ constructed on its 4 and 5 columns, are zeros. That is, in this case $|T(\alpha)| = 0$.

Case $c$). Let

$$E_3(\alpha) = \begin{Vmatrix} 0 & 0 & * \\ 0 & 0 & * \\ 0 & 0 & 0 \\ 0 & 0 & 0 \end{Vmatrix}.$$

According to Property 2.5, the matrix $T(\alpha)$ has the form

$$T(\alpha) = \begin{Vmatrix} \cdot & \cdot & \cdot & \cdot & \cdot & \cdot & \cdot \\ 0 & \cdot & \cdot & \cdot & \cdot & \cdot & \cdot \\ 0 & 0 & \cdot & \cdot & \cdot & \cdot & \cdot \\ 0 & 0 & * & \cdot & \cdot & \cdot & \cdot \\ 0 & 0 & * & \cdot & \cdot & \cdot & \cdot \\ 0 & 0 & 0 & 0 & 0 & \cdot & \cdot \\ 0 & 0 & 0 & 0 & 0 & \cdot & \cdot \end{Vmatrix} = \begin{Vmatrix} t_{ij} \end{Vmatrix}_1^7.$$

The matrix $V(\mathrm{E}, \Phi, \alpha)$ has the form

$$V(\mathrm{E}, \Phi, \alpha) = \begin{Vmatrix} 1 & & & & & & 0 \\ \cdot & 1 & & & & & \\ \cdot & \cdot & 1 & & & & \\ \cdot & \cdot & \mu_{43} & 1 & & & \\ \cdot & \cdot & \mu_{53} & \cdot & 1 & & \\ \cdot & \cdot & \cdot & \cdot & \cdot & 1 & \\ \cdot & \cdot & \cdot & \cdot & \cdot & \cdot & 1 \end{Vmatrix}, \quad \mu_{\mathrm{ij}} \neq 0.$$

Consider the submatrix of the matrix $V(\mathrm{E}, \Phi, \alpha)$:

$$\begin{Vmatrix} \mu_{43} & 1 & 0 & 0 & 0 \\ \mu_{53} & \cdot & 1 & 0 & 0 \\ \cdot & \cdot & \cdot & 1 & 0 \\ \cdot & \cdot & \cdot & \cdot & 1 \end{Vmatrix}.$$

All maximal order minors containing its last two columns are nonzeros. Therefore, the corresponding minors of the matrix $T(\alpha)$ are zeros. They have the form

$$\det \left\| \begin{array}{cc|cc} t_{i4} & t_{i5} & t_{i6} & t_{i7} \\ t_{j4} & t_{j5} & t_{j6} & t_{j7} \\ \hline 0 & 0 & t_{66} & t_{67} \\ 0 & 0 & t_{76} & t_{77} \end{array} \right\|, \quad 3 \le i < j \le 5.$$

Further considerations are the same as in the case $b$). The theorem is proved. $\square$

## 5.7. Monic divisors of polynomial matrices

In 1978 it was announced [55], and in 1980 it was published the article of P. Kazimirsky [56] giving an exhaustive answer to the questions, which for many years interested algebraists, namely, it indicated the necessary and sufficient conditions under which a matrix polynomial over the field of complex numbers has a monic divisor. A year later, his monograph [57] was published. In this article the scientist outlined the basic principles of the theory of the decomposition of matrix polynomials into multipliers. It is safe to say that these were the heyday of the theory of factorization of polynomial matrices. In a rather short period of time dozens of articles were published in which, under one or another limitation, was solved the mentioned problem. This subsection proposed a modern approach to its solution.

Let $F$ be an algebraically closed field of characteristic zero, and $A(x)$ be a nonsingular $n \times n$ matrix over $F[x]$. There are invertible matrices $P(x)$ and $Q(x)$ such that

$$P(x)A(x)Q(x) = \mathrm{diag}(\varepsilon_1(x), ..., \varepsilon_n(x)) = \mathrm{E}(x)$$

is the Smith form of the matrix $A(x)$. Let $\Phi(x) = \mathrm{diag}(\varphi_1(x), ..., \varphi_n(x))$ be a $d$-matrix, moreover $\Phi(x)|\mathrm{E}(x)$, $\deg \det \Phi(x) = nr$. Consider the determining matrix $V(\mathrm{E}, \Phi)$. Denote by $F(k)$ the transcendental extension of the field $F$ by attachment of all parameters $k_{ijl}$ appearing in the matrix $V(\mathrm{E}, \Phi)$.

**Theorem 5.16.** *The matrix $A(x)$ has a left monic divisor with the Smith form $\Phi(x)$ if and only if the matrix*

$$\left(V(\mathrm{E}, \Phi)P(x)\right)^{-1}\Phi(x)$$

*is right regularized over $F(k)[x]$.*

**Proof. Necessity.** Let $A(x) = B(x)C(x)$, where $B(x)$ is a monic matrix polynomial with the Smith form $\Phi(x)$. According to Corollary 4.4, all divisors of the matrix $A(x)$ form a set

$$\left(\mathbf{L}(\mathrm{E}, \Phi)P(x)\right)^{-1}\Phi(x)\mathrm{GL}_n(F[x]).$$

Therefore, the matrix $B(x)$ can be written as

$$B(x) = (L(x)P(x))^{-1} \Phi(x)K(x),$$

where $L(x) \in \mathbf{L}(\mathrm{E}, \Phi)$, $K(x) \in \mathrm{GL}_n(F[x])$. By virtue of Theorem 5.5,

$$L(x) = H(x)V_0(x)S(x),$$

where $H(x) \in \mathbf{G}_\Phi$, $V_0(x) \in \mathbf{V}(\mathrm{E}, \Phi)$, $S(x) \in U_n^{up}(F[x]) \subset \mathbf{G}_\mathrm{E}$. Then

$$B(x) = (L(x)P(x))^{-1} \Phi(x)K(x) = (H(x)V_0(x)S(x)P(x))^{-1} \Phi(x)K(x).$$

By Property 2.2, the set of all left transforming matrices of the matrix $A(x)$ has the form $\mathbf{P}_{A(x)} = \mathbf{G}_\mathrm{E}P$. Thus, $S(x)P(x) = P_0(x)$ is a left transforming matrix of the matrix $A(x)$. Then

$$B(x) = (H(x)V_0(x)P_0(x))^{-1} \Phi(x)K(x) = P_0^{-1}(x)V_0^{-1}(x)H^{-1}(x)\Phi(x)K(x).$$

Since

$$H^{-1}(x)\Phi(x) = \Phi(x)H_1(x),$$

where $H_1(x) \in \mathrm{GL}_n(F[x])$, then

$$B(x) = P_0^{-1}(x)V_0^{-1}(x)\Phi(x)H_1(x)K(x) = (V_0(x)P_0(x))^{-1} \Phi(x)K_1(x).$$

The matrix $B(x)$ is monic. Consequently, the matrix $(V_0(x)P_0(x))^{-1} \Phi(x)$ is right regularized. To carry out the change of polynomial coefficients of the matrix $V_0(x)$ for corresponding parameters $k_{ijl}$, we get $(V(\mathrm{E}, \Phi)P_0(x))^{-1} \Phi(x)$ is right regularized over $F(k)[x]$.

Let $P_1(x) \in \mathbf{P}_{A(x)}$. We show that the matrix

$$D(x) = (V(\mathrm{E}, \Phi)P_1(x))^{-1}\Phi(x)$$

is also right regularized. Since $P_1(x) = N(x)P_0(x)$, where $N(x) \in \mathbf{G}_\mathrm{E}$, then

$$D(x) = (V(\mathrm{E}, \Phi)P_1(x))^{-1} \Phi(x) = (V(\mathrm{E}, \Phi)N(x)P_0(x))^{-1} \Phi(x) =$$
$$= ((V(\mathrm{E}, \Phi)N(x)) P_0(x))^{-1} \Phi(x).$$

According to Theorems 5.15 and 2.8, there is a matrix $T(x) \in \mathbf{G}_\Phi$ such that $T(x)V(\mathrm{E}, \Phi)N(x)$ is the lower unitriangular matrix over $F(k)[x]$. By Lemma 5.3, there is a matrix $T_1(x) \in \mathbf{G}_\Phi$ such that

$$T_1(x)T(x)V(\mathrm{E}, \Phi)N(x) = V_1(\mathrm{E}, \Phi) \in \mathbf{V}(\mathrm{E}, \Phi).$$

Then

$$D(x) = (T_1(x)T(x)\left(V(\mathrm{E}, \Phi)N(x)\right)P_0(x))^{-1} T_1(x)T(x)\Phi(x) =$$
$$= ((T_1(x)T(x)V(\mathrm{E}, \Phi)N(x)) P_0(x))^{-1} \Phi(x)(\tilde{T}(x)\tilde{T}_1(x)) =$$
$$= (V_1(\mathrm{E}, \Phi)P_0(x))^{-1} \Phi(x)T_2(x).$$

Redenoted parameters $k_{ijl}$ in the matrix $V_1(\mathrm{E},\Phi)$ for $k'_{ijl}$, we get that the matrix $(V(\mathrm{E},\Phi)P_1(x))^{-1}\Phi(x)T_2(x)$, and therefore $(V(\mathrm{E},\Phi)P_1(x))^{-1}\Phi(x)$ is right regularized over $F(k)[x]$.

**Sufficiency.** By assumption, the matrix

$$(V(\mathrm{E},\Phi)P(x))^{-1}\Phi(x) = B(x)$$

is right regularized over $F(k)[x]$. That is, there is $U(x) \in \mathrm{GL}_n(F(k)[x])$ such that

$$B(x)U(x) = Ix^r + B_{r-1}x^{r-1} + ... + B_1x + B_0 = D(x).$$

It is obvious that $A(x) \in M_n(F(k)[x])$. According to Corollary 4.4, all left divisors of the matrix $A(x)$ with the Smith form $\Phi(x)$ form a set

$$(\mathbf{L}(\mathrm{E},\Phi)P(x))^{-1}\Phi(x)\mathrm{GL}_n(F(k)[x]).$$

Since $V(\mathrm{E},\Phi) \in \mathbf{L}(\mathrm{E},\Phi)$, then $D(x)$ is a left divisor of the matrix $A(x)$ :

$$A(x) = D(x)C(x).$$

By Theorem 5.12, the matrix $B(x)$ is right regularized if and only if

$$\det M_{P(x)\|I\ \ Ix\ \ ...\ \ Ix^{r-1}\|}(\Phi) = f(k_{210}, ..., k_{n,n-1,h_{n,n-1}}, x) \neq 0.$$

Since $F$ is an algebraically closed field so $F$ contains an infinite number of elements. Therefore, there are elements $p_{210}, ..., p_{n,n-1,h_{n,n-1}}, p_{nn}$, such that

$$f(p_{210}, ..., p_{n,n-1,h_{n,n-1}}, p_{nn}) \neq 0.$$

Then the matrix $\overline{D}(x)$, which is obtained from the matrix $D(x)$ by changing of variables $k_{210}, ..., k_{n,n-1,h_{n,n-1}}, x$ for corresponding elements $p_{210}, ..., p_{n,n-1,h_{n,n-1}}, p_{nn}$ from the field $F$, and will be the desired monic divisor of the matrix $A(x)$. $\qquad\square$

We apply the obtained results to the search for solutions of one-sided matrix equations.

**Example 5.6.** Consider the equation

$$X^2 = \begin{Vmatrix} 0 & 1 & 0 \\ 0 & 0 & 0 \\ 0 & 0 & 0 \end{Vmatrix}. \tag{5.18}$$

The polynomial matrix

$$A(x) = Ix^2 - \begin{Vmatrix} 0 & 1 & 0 \\ 0 & 0 & 0 \\ 0 & 0 & 0 \end{Vmatrix} = \begin{Vmatrix} x^2 & -1 & 0 \\ 0 & x^2 & 0 \\ 0 & 0 & x^2 \end{Vmatrix}$$

corresponds to this equation. By the generalized Bezout theorem, the matrix $B$ is the root of the equation (5.18) if and only if

$$A(x) = (Ix - B)C(x).$$

That is, the problem of solving (5.18) is equivalent to finding the left monic divisors of the matrix $A(x)$. We reduce the matrix $A(x)$ to its Smith form:

$$\underbrace{\begin{Vmatrix} 1 & 0 & 0 \\ 0 & 0 & 1 \\ x^2 & 1 & 0 \end{Vmatrix}}_{P(x)} A(x) \underbrace{\begin{Vmatrix} 0 & 0 & 1 \\ -1 & 0 & x^2 \\ 0 & 1 & 0 \end{Vmatrix}}_{Q(x)} = \begin{Vmatrix} 1 & 0 & 0 \\ 0 & x^2 & 0 \\ 0 & 0 & x^4 \end{Vmatrix} = \mathrm{E}(x).$$

Potential divisors have Smith forms:

$$\Phi(x) = \mathrm{diag}\,(1, 1, x^3), \quad \Phi_1(x) = \mathrm{diag}\,(1, x, x^2).$$

First we treat divisors with the Smith form $\Phi(x)$. In this case

$$V(\mathrm{E}, \Phi) = \begin{Vmatrix} 1 & 0 & 0 \\ 0 & 1 & 0 \\ 0 & ax^2 + bx & 1 \end{Vmatrix},$$

where $a, b$ are parameters. Then the matrix $A(x)$ has the desired divisor if and only if the matrix $(V(\mathrm{E}, \Phi)P(x))^{-1}\Phi(x)$ is right regularized. Since

$$\det M_{V(\mathrm{E}, \Phi)P(x)}(\Phi) = \det M_{\| x^2 \ 1 \ ax^2 + bx \|}(x^3) = \det \underbrace{\begin{Vmatrix} 0 & 1 & 0 \\ 0 & 0 & b \\ 2 & 0 & 2a \end{Vmatrix}}_{T} = 2b \neq 0,$$

according to Theorems 5.16, 5.12, 5.13, the matrix $A(x)$ has a left monic divisor of the form

$$B = T^{-1} \begin{Vmatrix} 0 & 0 & 0 \\ 1 & 0 & 0 \\ 0 & 2 & 0 \end{Vmatrix} T = \begin{Vmatrix} 0 & ab^{-1} & b \\ 1 & 0 & 0 \\ 0 & b^{-1} & 0 \end{Vmatrix},$$

where $b \neq 0$. Substituting $ab^{-1}$ by parameter $c$, we get that all matrices of the set

$$\left\{ \begin{Vmatrix} 0 & c & b \\ 0 & 0 & 0 \\ 0 & b^{-1} & 0 \end{Vmatrix} \mid c \in F, \ b \in F \backslash \{0\} \right\} \tag{5.19}$$

are solutions of equation (5.18). Moreover, according to Theorem 5.9, the set (5.19) consists of all roots of the equation (5.18) with the Smith form $\Phi(x)$.

The matrix $(V(\mathrm{E}, \Phi_1)P(x))^{-1}\Phi_1(x)$ is not right regulated. So the set (5.19) consists of all roots of equation (5.18). $\diamondsuit$

## 5.8. Structure of g.c.d. of matrices

In this subsection, we continue the study of the greatest common matrix divisor, started in subsection 1.5. The main focus will be on the study of its structure, that is, its Smith form and transforming matrices. We begin this analysis with the second order matrices over Bezout rings of stable range 1.5.

Let $A$ be a $2 \times 2$ matrix. There are such invertible matrices $P_A, Q_A$, that $P_A A Q_A = \mathrm{E}$, where $\mathrm{E} = \mathrm{diag}(\varepsilon_1, \varepsilon_2)$, $\varepsilon_1 | \varepsilon_2$. Consider the group $\mathbf{G}_{\mathrm{E}}$ of all invertible matrices of the form

$$\left\|\begin{array}{cc} h_{11} & h_{12} \\ \dfrac{\varepsilon_2}{\varepsilon_1} h_{21} & h_{22} \end{array}\right\|.$$

For convenience, we will also sometimes denote this group as $\mathbf{G}_{\frac{\varepsilon_2}{\varepsilon_1}}$.

Consider the matrix $B = P_B^{-1} \Delta Q_B^{-1}$, where $\Delta = \mathrm{diag}(\delta_1, \delta_2)$, $\delta_1 | \delta_2$. By Property 2.2 the sets of left transforming matrices $\mathbf{P}_B, \mathbf{P}_A$ form adjacent classes $\mathbf{P}_B = \mathbf{G}_\Delta P_B, \mathbf{P}_A = \mathbf{G}_{\mathrm{E}} P_A$. Let's investigate the set $\mathbf{P}_B \mathbf{P}_A^{-1}$.

**Lemma 5.9.** *Let $P_B P_A^{-1} = \|s_{ij}\|_1^2 = S$. The element*

$$((\varepsilon_2, \delta_2), s_{21}[\varepsilon_1, \delta_1])$$

*is invariant with respect to the choice of transforming matrices $P_B$ and $P_A$.*

**Proof.** Let at least one of the matrices $A$, $B$ be nonsingular, and let $F_A$ and $F_B$ be their left transforming matrices. There are such $H_A \in \mathbf{G}_{\frac{\varepsilon_2}{\varepsilon_1}}$ and $H_B \in \mathbf{G}_{\frac{\delta_2}{\delta_1}}$, that

$$F_A = H_A P_A, \ F_B = H_B P_B.$$

Denote $F_B F_A^{-1} = \|s'_{ij}\|_1^2$. Consider the product of matrices

$$F_B F_A^{-1} = H_B P_B (H_A P_A)^{-1} = H_B P_B P_A^{-1} H_A^{-1} = H_B S H_A^{-1},$$

where $S = P_B P_A^{-1}$. Set $H_B S = \|k_{ij}\|_1^2$. Since

$$H_B = \left\|\begin{array}{cc} h_{11} & h_{12} \\ \dfrac{\delta_2}{\delta_1} h_{21} & h_{22} \end{array}\right\|,$$

then

$$k_{21} = \left\|\begin{array}{cc} \dfrac{\delta_2}{\delta_1} h_{21} & h_{22} \end{array}\right\| \left\|\begin{array}{c} s_{11} \\ s_{21} \end{array}\right\| = \dfrac{\delta_2}{\delta_1} h_{21} s_{11} + h_{22} s_{21}.$$

Consider

$$\mu = ((\varepsilon_2, \delta_2), k_{21}[\varepsilon_1, \delta_1]) = \left((\varepsilon_2, \delta_2), \left(\dfrac{\delta_2}{\delta_1} h_{21} s_{21} + h_{22} s_{21}\right)[\varepsilon_1, \delta_1]\right) =$$

$$= \left((\varepsilon_2, \delta_2), \dfrac{\delta_2}{\delta_1}[\varepsilon_1, \delta_1] h_{21} s_{21} + h_{22} s_{21}[\varepsilon_1, \delta_1]\right).$$

From the inclusion
$$\frac{\delta_2[\varepsilon_1,\delta_1]}{\delta_1(\varepsilon_2,\delta_2)} \in R,$$
implies that
$$(\varepsilon_2,\delta_2)\Big|\frac{\delta_2}{\delta_1}[\varepsilon_1,\delta_1].$$
Hence,
$$\mu = ((\varepsilon_2,\delta_2), h_{22}s_{21}[\varepsilon_1,\delta_1]) = \tau\left(\frac{(\varepsilon_2,\delta_2)}{\tau}, \frac{[\varepsilon_1,\delta_1]}{\tau}h_{22}s_{21}\right) =$$
$$= \tau\left(\frac{(\varepsilon_2,\delta_2)}{\tau}, h_{22}s_{21}\right),$$
where $\tau = ((\varepsilon_2,\delta_2),\ [\varepsilon_1,\delta_1])$. Since
$$\frac{\delta_2\left((\varepsilon_2,\delta_2),[\varepsilon_1,\delta_1]\right)}{\delta_1(\varepsilon_2,\delta_2)} = \frac{(\delta_2(\varepsilon_2,\delta_2),\delta_2[\varepsilon_1,\delta_1])}{\delta_1(\varepsilon_2,\delta_2)} \in R,$$
then
$$\frac{(\varepsilon_2,\delta_2)}{\tau}\Big|\frac{\delta_2}{\delta_1}.$$
From the invertibility of the matrix $H_B$ it follows that
$$\left(h_{22},\frac{\delta_2}{\delta_1}\right) = 1.$$
Therefore,
$$\left(h_{22},\frac{(\varepsilon_2,\delta_2)}{\tau}\right) = 1.$$
Consequently,
$$\mu = \tau\left(\frac{(\varepsilon_2,\delta_2)}{\tau}, s_{21}\right) = ((\varepsilon_2,\delta_2), s_{21}(\varepsilon_2,\delta_2), s_{21}[\varepsilon_1,\delta_1]) = ((\varepsilon_2,\delta_2), s_{21}[\varepsilon_1,\delta_1]).$$
So
$$\left((\varepsilon_2,\delta_2), k_{21}[\varepsilon_1,\delta_1]\right) = \left((\varepsilon_2,\delta_2), s_{21}[\varepsilon_1,\delta_1]\right).$$

Consider $SH_A^{-1} = \|t_{ij}\|_1^2$. Since
$$H_A^{-1} = \left\|\begin{matrix} v_{11} & v_{12} \\ \dfrac{\varepsilon_2}{\varepsilon_1}v_{21} & v_{22} \end{matrix}\right\|,$$
then
$$t_{21} = \left\|s_{21} \quad s_{22}\right\|\left\|\begin{matrix} \nu_{11} \\ \dfrac{\varepsilon_2}{\varepsilon_1}\nu_{21} \end{matrix}\right\| = s_{21}\nu_{11} + \frac{\varepsilon_2}{\varepsilon_1}\nu_{21}s_{22}.$$
Noting that
$$\frac{(\varepsilon_2,\delta_2)}{((\varepsilon_2,\delta_2),[\varepsilon_1,\delta_1])}\Big|\frac{\varepsilon_2}{\varepsilon_1}$$
and thinking as above, we get
$$((\varepsilon_2,\delta_2), t_{21}[\varepsilon_1,\delta_1]) = ((\varepsilon_2,\delta_2), s_{21}[\varepsilon_1,\delta_1]).$$

By virtue of the associativity of the ring $M_2(R)$, we conclude our consideration of this case.

Assume that

$$A \sim \left\| \begin{matrix} \varepsilon_1 & 0 \\ 0 & 0 \end{matrix} \right\|, \quad B \sim \left\| \begin{matrix} \delta_1 & 0 \\ 0 & 0 \end{matrix} \right\|$$

and $F_A \in \mathbf{P}_A$, $F_B \in \mathbf{P}_B$. There are such matrices

$$H_A = \left\| \begin{matrix} e_1' & v_{12} \\ 0 & e_2' \end{matrix} \right\|, \quad H_B = \left\| \begin{matrix} e_1 & h_{12} \\ 0 & e_2 \end{matrix} \right\|,$$

that $F_A = H_A P_A$, $F_B = H_B P_B$. Set $F_B F_A^{-1} = \|s_{ij}'\|_1^2 = S'$. To prove Lemma it is necessary to show that $s_{21}$ and $s_{21}'$ are different by the invertible element of the ring $R$. Consider the product of matrices

$$S' = F_B F_A^{-1} = H_B P_B (H_A P_A)^{-1} = H_B P_B P_A^{-1} H_A^{-1} = H_B S H_A^{-1},$$

where $S = P_B P_A^{-1}$. Write these matrices explicitly, that is

$$\left\| \begin{matrix} s_{11}' & s_{12}' \\ s_{21}' & s_{22}' \end{matrix} \right\| = \left\| \begin{matrix} e_1 & h_{12} \\ 0 & e_2 \end{matrix} \right\| \left\| \begin{matrix} s_{11} & s_{12} \\ s_{21} & s_{22} \end{matrix} \right\| \left\| \begin{matrix} e_1'^{-1} & * \\ 0 & e_2'^{-1} \end{matrix} \right\| =$$
$$= \left\| \begin{matrix} s_{11}'' & s_{12}'' \\ s_{21}u & s_{22}'' \end{matrix} \right\|,$$

where $u = e_2 e_2'^{-1}$ is an invertible element of the ring $R$. Therefore, $s_{21}' = s_{21}u$, which is to be proved. □

**Lemma 5.10.** *Let $a, b \in R$. In order an invertible matrix $S = \|s_{ij}\|_1^2$ could be written as $S = L_a L_b$, where $L_a \in \mathbf{G}_a$, $L_b \in \mathbf{G}_b$, it is necessary and sufficient $s_{21} = (a, b)t$.*

**Proof. Necessity.** Let $S = L_a L_b$, where $L_a \in \mathbf{G}_a$, $L_b \in \mathbf{G}_b$. Let us write down these matrices in an expanded form

$$S = \left\| \begin{matrix} s_{11} & s_{12} \\ s_{21} & s_{22} \end{matrix} \right\| = \left\| \begin{matrix} l_{11} & l_{12} \\ al_{21} & l_{22} \end{matrix} \right\| \left\| \begin{matrix} h_{11} & h_{12} \\ bh_{21} & h_{22} \end{matrix} \right\| = L_a L_b.$$

It follows that

$$\left\| \begin{matrix} s_{11} & s_{12} \\ s_{21} & s_{22} \end{matrix} \right\| = \left\| \begin{matrix} l_{11}' & l_{12}' \\ al_{21}h_{11} + bl_{22}h_{22} & l_{22}' \end{matrix} \right\|.$$

This yields

$$s_{21} = al_{21}h_{11} + bl_{22}h_{22}.$$

It means that $s_{21} = (a, b)t$.

**Sufficiency.** If $a = b = 0$, then $s_{21} = 0$. Consequently, $S \in \mathbf{G}_a = \mathbf{G}_b$ and the proof is obvious.

If $a = 0$ and $b \neq 0$, then $s_{21} = bt$. Therefore $S \in \mathbf{G}_b$.

Consider the case $a, b \neq 0$. Invertibility of the matrix $S$ implies that

$$(s_{21}, s_{22}) = 1.$$

By Theorem 1.9, there are $k_{12}$, $k_{22}$ such that

$$s_{21}k_{12} + s_{22}k_{22} = 1,$$

where $(k_{22}, b) = 1$. Hence, $(k_{12}b, k_{22}) = 1$. It follows that there are such $k_{11}$, $k_{21}$, that

$$k_{11}k_{22} - bk_{21}k_{12} = 1.$$

This means that the matrix $K_1 = \left\| \begin{matrix} k_{11} & k_{12} \\ bk_{21} & k_{22} \end{matrix} \right\|$ belongs to the group $\mathbf{G}_b$. Then

$$SK_1 = \left\| \begin{matrix} s_{11} & s_{12} \\ (a,b)t & s_{22} \end{matrix} \right\| \left\| \begin{matrix} k_{11} & k_{12} \\ bk_{21} & k_{22} \end{matrix} \right\| = \left\| \begin{matrix} g_{11} & g_{12} \\ (a,b)g_{21} & 1 \end{matrix} \right\| = S_1.$$

Hence

$$\left\| \begin{matrix} 1 & -g_{12} \\ 0 & 1 \end{matrix} \right\| S_1 = H_1 S_1 = \left\| \begin{matrix} q_{11} & 0 \\ (a,b)g_{21} & 1 \end{matrix} \right\| = S_2,$$

where $H_1 \in \mathbf{G}_a$. The matrix $S_2$ is invertible, thus $q_{11}$ is invertible element of the ring $R$. Then

$$\left\| \begin{matrix} q_{11}^{-1} & 0 \\ 0 & 1 \end{matrix} \right\| S_2 = H_2 S_2 = \left\| \begin{matrix} 1 & 0 \\ (a,b)g_{21} & 1 \end{matrix} \right\|,$$

where $H_2 \in \mathbf{G}_a$. Since $H_1, H_2 \in \mathbf{G}_a$ then $H_3 = H_2 H_1 \in \mathbf{G}_a$ i

$$\left\| \begin{matrix} 1 & 0 \\ (a,b)g_{21} & 1 \end{matrix} \right\| = H_3 S K_1.$$

There are $u$ and $v$ in the ring $R$ such that

$$(a,b)g_{21} = (au + bv)g_{21} = aug_{21} + bvg_{21}.$$

Consider the matrices

$$H_4 = \left\| \begin{matrix} 1 & 0 \\ aug_{21} & 1 \end{matrix} \right\| \in \mathbf{G}_a, \quad K_2 = \left\| \begin{matrix} 1 & 0 \\ bvg_{21} & 1 \end{matrix} \right\| \in \mathbf{G}_b.$$

Then $H_3 S K_1 = H_4 K_2$, that is

$$S = (H_3^{-1} H_4)(K_2 K_1^{-1}).$$

Noting that $H_3^{-1} H_4 \in \mathbf{G}_a$, $K_2 K_1^{-1} \in \mathbf{G}_b$ we are finishing the proof. $\square$

**Lemma 5.11.** *Let*

$$A = P_A^{-1} \left\| \begin{matrix} \varepsilon_1 & 0 \\ 0 & \varepsilon_2 \end{matrix} \right\| Q_A^{-1}, \quad B = P_B^{-1} \left\| \begin{matrix} \delta_1 & 0 \\ 0 & \delta_2 \end{matrix} \right\| Q_B^{-1}$$

*are matrices over* $R$ *and*

$$D = P_D^{-1} \left\| \begin{matrix} \varphi_1 & 0 \\ 0 & \varphi_2 \end{matrix} \right\| Q_D^{-1}, \ \ T = P_T^{-1} \left\| \begin{matrix} \gamma_1 & 0 \\ 0 & \gamma_2 \end{matrix} \right\| Q_T^{-1},$$

*moreover*

$$A = DA_1, \ \ B = DB_1, \ \ A = TA_2, \ \ B = TB_2.$$

*If* $\gamma_1 | \varphi_1 = (\varepsilon_1, \delta_1)$ *and* $\gamma_2 | \varphi_2$, *then* $D = TN$.

**Proof.** According to Theorem 4.2 we get

$$A = DA_1 \Rightarrow P_D = L_A P_A, \text{ where } L_A \in \mathbf{G}_{\frac{\varphi_2}{(\varphi_2, \varepsilon_1)}},$$
$$B = DB_1 \Rightarrow P_D = L_B P_B, \text{ where } L_B \in \mathbf{G}_{\frac{\varphi_2}{(\varphi_2, \delta_1)}}.$$

For similar reasons we have

$$A = TA_2 \Rightarrow P_T = K_A P_A, \text{ where } K_A \in \mathbf{G}_{\frac{\gamma_2}{(\gamma_2, \varepsilon_1)}},$$
$$B = TB_2 \Rightarrow P_T = K_B P_B, \text{ where } K_B \in \mathbf{G}_{\frac{\gamma_2}{(\gamma_2, \delta_1)}}.$$

Consider the product

$$P_T P_D^{-1} = K_A P_A P_A^{-1} L_A^{-1} = K_A L_A^{-1} =$$
$$= \underbrace{\left\| \begin{matrix} k_{11} & k_{12} \\ \dfrac{\gamma_2}{(\gamma_2, \varepsilon_1)} k_{21} & k_{22} \end{matrix} \right\|}_{K_A} \underbrace{\left\| \begin{matrix} h_{11} & h_{12} \\ \dfrac{\varphi_2}{(\varphi_2, \varepsilon_1)} h_{21} & h_{22} \end{matrix} \right\|}_{L_A^{-1}}.$$

That is

$$P_T P_D^{-1} = \left\| \begin{matrix} l_{11} & l_{12} \\ \left( \dfrac{\gamma_2}{(\gamma_2, \varepsilon_1)}, \dfrac{\varphi_2}{(\varphi_2, \varepsilon_1)} \right) l_{21} & l_{22} \end{matrix} \right\|.$$

Taking into account Property 1.6, we receive

$$\left( \frac{\gamma_2}{(\gamma_2, \varepsilon_1)}, \frac{\varphi_2}{(\varphi_2, \varepsilon_1)} \right) = \frac{(\gamma_2, \varphi_2)}{(\gamma_2, \varphi_2, \varepsilon_1)}.$$

Since $\gamma_2 | \varphi_2$ then

$$\frac{(\gamma_2, \varphi_2)}{(\gamma_2, \varphi_2, \varepsilon_1)} = \frac{\gamma_2}{(\gamma_2, \varepsilon_1)}.$$

Consequently,

$$P_T P_D^{-1} = \left\| \begin{matrix} l_{11} & l_{12} \\ \dfrac{\gamma_2}{(\gamma_2, \varepsilon_1)} l_{21} & l_{22} \end{matrix} \right\| = \left\| \begin{matrix} m_{11} & m_{12} \\ m_{21} & m_{22} \end{matrix} \right\| = M. \tag{5.20}$$

On the other hand

$$P_T P_D^{-1} = K_B P_B P_B^{-1} L_B^{-1} = K_B L_B^{-1}.$$

Similarly, we show that

$$P_T P_D^{-1} = \left\| \begin{matrix} v_{11} & \nu_{12} \\ \dfrac{\gamma_2}{(\gamma_2, \delta_1)} \nu_{21} & \nu_{22} \end{matrix} \right\| = \left\| \begin{matrix} m_{11} & m_{12} \\ m_{21} & m_{22} \end{matrix} \right\| = M. \qquad (5.21)$$

From equalities (5.20) and (5.21) it follows that $\dfrac{\gamma_2}{(\gamma_2, \varepsilon_1)} \,|\, m_{21}$ and $\dfrac{\gamma_2}{(\gamma_2, \delta_1)} \,|\, m_{21}$. Therefore

$$\left[ \frac{\gamma_2}{(\gamma_2, \varepsilon_1)}, \frac{\gamma_2}{(\gamma_2, \delta_1)} \right] \,|\, m_{21}.$$

By Property 1.8

$$\left[ \frac{\gamma_2}{(\gamma_2, \varepsilon_1)}, \frac{\gamma_2}{(\gamma_2, \delta_1)} \right] = \frac{\gamma_2}{(\gamma_2, (\varepsilon_1, \delta_1))}.$$

Since $(\varepsilon_1, \delta_1) = \varphi_1$ then

$$\frac{\gamma_2}{(\gamma_2, (\varepsilon_1, \delta_1))} = \frac{\gamma_2}{(\gamma_2, \varphi_1)}.$$

Hence, $M \in \mathbf{G}_{\frac{\gamma_2}{(\gamma_2, \varphi_1)}}$. According to the condition of our statement the matrix $\mathrm{diag}(\gamma_1, \gamma_2)$ is a divisor of the matrix $\mathrm{diag}(\varphi_1, \varphi_2)$. Based on Theorem 4.2, we get $D = TN$. Lemma is proved. $\qquad \square$

**Lemma 5.12.** *Suppose that $A$ and $B$ are singular matrices and $P_B P_A^{-1} =$*
$= \left\| \begin{matrix} e_1 & s \\ 0 & e_2 \end{matrix} \right\|$. *Then*

*1) for any other matrices $P_A' \in \mathbf{P}_A$ and $P_B' \in \mathbf{P}_B$ the equality*

$$P_A' P_B'^{-1} = \left\| \begin{matrix} s_{11}' & s_{12}' \\ 0 & s_{22}' \end{matrix} \right\|$$

*is fulfilled;*

*2) $\mathbf{P}_A \cap \mathbf{P}_B \neq \varnothing$.*

**Proof.** 1). Since $A$ and $B$ are singular matrices, consequently $\mathbf{G}_\mathrm{E}, \mathbf{G}_\Delta$ coincide with the group of invertible upper triangular matrices.

Let $P_A'$ and $P_B'$ are other left transforming matrices of the matrices $A$ and $B$. There are such $H_A \in \mathbf{G}_\mathrm{E}$ and $H_B \in \mathbf{G}_\Delta$, that

$$P_A' = H_A P_A, \quad P_B' = H_B P_B.$$

Consider the product of matrices

$$P_A' P_B'^{-1} = H_A P_A (H_B P_B)^{-1} = H_A P_A P_B^{-1} H_B^{-1} = H_A S H_B^{-1}.$$

Since $H_A$ and $H_B^{-1}$ are upper triangular matrices, then $H_A S H_B^{-1}$ will also be an upper triangular matrix.

2). Noting that

$$P_B P_A^{-1} = \left\| \begin{matrix} e_1 & 1 \\ 0 & e_2 \end{matrix} \right\| = H \in \mathbf{G}_\mathrm{E} = \mathbf{G}_\Delta$$

we get $P_B = H P_A$, where $H \in \mathbf{G}_\mathrm{E}$. Since $\mathbf{P}_A = \mathbf{G}_\mathrm{E} P_A$, we have $P_B \in \mathbf{P}_A$. That is $P_B \in \mathbf{P}_A \cap \mathbf{P}_B$. $\qquad \square$

**Theorem 5.17.** *Let $R$ be a Bezout ring of stable range of* 1.5. *And let*

$$A \sim \mathrm{E} = \mathrm{diag}(\varepsilon_1,\, \varepsilon_2), \ \ B \sim \Delta = \mathrm{diag}(\delta_1,\, \delta_2),$$

$P_B P_A^{-1} = \|s_{ij}\|_1^2,\ P_B \in \mathbf{P}_B,\ P_A \in \mathbf{P}_A.$ *Then*

$$(A, B)_l = (L_A P_A)^{-1}\Phi = (L_B P_B)^{-1}\Phi,$$
*where*
$$\Phi = \mathrm{diag}(\varphi_1,\, \varphi_2) = \mathrm{diag}((\varepsilon_1, \delta_1),\ ((\varepsilon_2, \delta_2), s_{21}[\varepsilon_1, \delta_1])),$$

*and the matrices $L_A$ and $L_B$ satisfy the equality $L_B^{-1} L_A = P_B P_A^{-1}$ and belong to groups $\mathbf{G}_{\frac{\varphi_2}{(\varphi_2,\varepsilon_1)}}$ and $\mathbf{G}_{\frac{\varphi_2}{(\varphi_2,\delta_1)}}$, respectively.*

**Proof.** Note that, by Lemma 5.9, the element $((\varepsilon_2, \delta_2), s_{21}[\varepsilon_1, \delta_1])$, and, hence, the matrix $\Phi$, do not depend on the choice of the transforming matrices $P_A$ and $P_B$. In this case, according to the Corollary 2.15 the matrix $(A, B)_l$ has the Smith form $\Phi$.

At first, consider the case when at least one of the matrices $A$, $B$ is non-singular, or $A$, and $B$ are singular, moreover $s_{21} \neq 0$.

We show that the matrix $P_B P_A^{-1}$ can be written as

$$P_B P_A^{-1} = MN, \tag{5.22}$$
*where*
$$M \in \mathbf{G}_{\frac{\varphi_2}{(\varphi_2,\delta_1)}},\ \ N \in \mathbf{G}_{\frac{\varphi_2}{(\varphi_2,\varepsilon_1)}},\ \ \varphi_2 = ((\varepsilon_2, \delta_2), s_{21}[\varepsilon_1, \delta_1]).$$

Using Property 1.7 we have

$$\left(\frac{\varphi_2}{(\varphi_2, \delta_1)}, \frac{\varphi_2}{(\varphi_2, \varepsilon_1)}\right) = \frac{\varphi_2}{(\varphi_2, [\varepsilon_1, \delta_1])} = \frac{((\varepsilon_2, \delta_2), s_{21}[\varepsilon_1, \delta_1])}{((\varepsilon_2, \delta_2), [\varepsilon_1, \delta_1])} =$$
$$= \left(\frac{(\varepsilon_2, \delta_2)}{((\varepsilon_2, \delta_2), [\varepsilon_1, \delta_1])}, \frac{[\varepsilon_1, \delta_1]}{((\varepsilon_2, \delta_2), [\varepsilon_1, \delta_1])} s_{21}\right) = \left(\frac{(\varepsilon_2, \delta_2)}{((\varepsilon_2, \delta_2), [\varepsilon_1, \delta_1])}, s_{21}\right) = \mu.$$

Since $\mu | s_{21}$, based on Lemma 5.10 the matrix $P_B P_A^{-1}$ can be represented as (5.22). From this equality it follows that

$$M^{-1} P_B = N P_A.$$

Denoting $M^{-1} = L_B$, $N = L_A$, we obtain

$$(L_A P_A)^{-1}\Phi = (L_B P_B)^{-1}\Phi = D.$$

Since $\Phi | \mathrm{E}$, $\Phi | \Delta$, and taking into account that

$$L_A \in \mathbf{G}_{\frac{\varphi_2}{(\varphi_2,\varepsilon_1)}} = \mathbf{L}(\mathrm{E}, \Phi) \neq \varnothing,\ \ L_B \in \mathbf{G}_{\frac{\varphi_2}{(\varphi_2,\delta_1)}} = \mathbf{L}(\Delta, \Phi) \neq \varnothing,$$

by Theorem 4.2, we receive that the matrix $D$ is a left common divisor of the matrices $A$ and $B$.

Suppose that $T = P_T^{-1} \Gamma Q_T^{-1}$, where $\Gamma = \mathrm{diag}(\gamma_1, \gamma_2)$ is the another left common divisor of the matrices $A$ and $B$. The Smith form of the matrix $T$ divides the Smith form of the matrix $(A, B)_l$, which according to Corollary 2.15, is the matrix $\Phi$. Hence, $\Gamma | \Phi$. Consequently, $\gamma_1 | \varphi_1 = (\varepsilon_1, \delta_1)$ and $\gamma_2 | \varphi_2$. By Lemma 5.11, the matrix $T$ is a left divisor of the matrix $D$. Therefore, $D$ is the left g.c.d. of matrices $A, B$.

Finally, consider the case $s_{21} = 0$ and

$$A \sim \mathrm{E} = \mathrm{diag}(\varepsilon_1, 0), \ \ B \sim \Delta = \mathrm{diag}(\delta_1, 0).$$

Taking into account Lemma 5.12, we receive $\mathbf{P}_A \cap \mathbf{P}_B \neq \varnothing$. Let $U \in \mathbf{P}_A \cap \mathbf{P}_B$. It means that the matrices $A$ and $B$ can be written as

$$A = U^{-1} \mathrm{E} Q_A^{-1}, \ \ B = U^{-1} \Delta Q_B^{-1}.$$

Consider the matrix

$$D = U^{-1} \mathrm{diag}((\varepsilon_1, \delta_1), 0) = U^{-1} \Phi.$$

Since

$$A = \left( U^{-1} \left\| \begin{matrix} (\varepsilon_1, \delta_1) & 0 \\ 0 & 0 \end{matrix} \right\| \right) \left( \left\| \begin{matrix} \dfrac{\varepsilon_1}{(\varepsilon_1, \delta_1)} & 0 \\ 0 & 0 \end{matrix} \right\| Q_A^{-1} \right) = DA_1,$$

$$B = \left( U^{-1} \left\| \begin{matrix} (\varepsilon_1, \delta_1) & 0 \\ 0 & 0 \end{matrix} \right\| \right) \left( \left\| \begin{matrix} \dfrac{\delta_1}{(\varepsilon_1, \delta_1)} & 0 \\ 0 & 0 \end{matrix} \right\| Q_B^{-1} \right) = DB_1,$$

then $D$ is a left common divisor of the matrices $A$ and $B$.

Let $T = P_T^{-1} \Gamma Q_T^{-1}$, where $\Gamma = \mathrm{diag}(\gamma_1, \gamma_2)$, $\gamma_1 | \gamma_2$ is another common left divisor of the matrices $A$ and $B$. Thus, $\Gamma | \mathrm{E}$ and $\Gamma | \Delta$. It follows that $\gamma_1 | \varepsilon_1$ and $\gamma_1 | \delta_1$. That is $\gamma_1 | (\varepsilon_1, \delta_1)$. Consequently, $\Gamma | \Phi$. By Lemma 5.11 the matrix $T$ is a left divisor of the matrix $D$. So $D$ is the left g.c.d. of $A$ and $B$. The theorem is proved. □

The following result concerns left g.c.d. of higher order matrices with some restrictions on their Smith forms.

Let $R$ be an elementary divisor ring. And let

$$A \sim \mathrm{E} = \mathrm{diag}(\varepsilon_1, \varepsilon_2, ..., \varepsilon_n), \ \ B \sim \Delta = \mathrm{diag}(1, ..., 1, \delta)$$

are nonsingular matrices over $R$. Set $P_B P_A^{-1} = S = \|s_{ij}\|_1^n$, where $P_B \in \mathbf{P}_B$, $P_A \in \mathbf{P}_A$.

**Lemma 5.13.** *The element*

$$((\varepsilon_n, \delta), s_{n1} \varepsilon_1, ..., s_{n.n-1} \varepsilon_{n-1})$$

*is invariant under the choice of the transforming matrices $P_B$ and $P_A$.*

**Proof.** Let $F_A$ and $F_B$ are other left transforming matrices of the matrices $A$ and $B$. Then there exist $H_A \in \mathbf{G}_E$ and $H_B \in \mathbf{G}_\Delta$ such that

$$F_A = H_A P_A, \quad F_B = H_B P_B.$$

Consider the product of matrices

$$F_B F_A^{-1} = H_B P_B (H_A P_A)^{-1} = H_B P_B P_A^{-1} H_A^{-1} = H_B S H_A^{-1},$$

where $S = P_B P_A^{-1}$. Denote $H_B S = \|k_{ij}\|_1^n$. By virtue of Theorem 2.7, the matrix $H_B$ has the form

$$H_B = \left\| \begin{array}{cccc} h_{11} & ... & h_{1.n-1} & h_{1n} \\ ... & ... & ... & ... \\ h_{n-1.1} & ... & h_{n-1.n-1} & h_{n-1.n} \\ \delta h_{n1} & ... & \delta h_{n.n-1} & h_{nn} \end{array} \right\|.$$

Thus,

$$k_{ni} = \left\| \begin{array}{cccc} \delta h_{n1} & ... & \delta h_{n.n-1} & h_{nn} \end{array} \right\| \left\| \begin{array}{c} s_{1i} \\ \vdots \\ s_{n-1.i} \\ s_{ni} \end{array} \right\| =$$

$$= \delta(h_{n1}s_{1i} + ... + h_{n.n-1}s_{n-1.i}) + h_{nn}s_{ni} = \delta l_i + h_{nn}s_{ni}, \quad i = 1, ..., n-1.$$

Consider

$$((\varepsilon_n, \delta), \varepsilon_1 k_{n1}, ..., \varepsilon_{n-1}k_{n.n-1}) =$$

$$= ((\varepsilon_n, \delta), \delta \varepsilon_1 l_1 + \varepsilon_1 h_{nn}s_{n1}, ..., \delta \varepsilon_{n-1}l_{n-1} + \varepsilon_{n-1}h_{nn}s_{n.n-1}) = d.$$

Since $(\varepsilon_n, \delta) | \delta$, we have

$$d = ((\varepsilon_n, \delta), \varepsilon_1 h_{nn}s_{n1}, ..., \varepsilon_{n-1}h_{nn}s_{n.n-1}) =$$

$$= ((\varepsilon_n, \delta), h_{nn}(\varepsilon_1 s_{n1}, ..., \varepsilon_{n-1}s_{n.n-1})).$$

It implies out of the invertibility of the matrix $H_B$ that $(h_{nn}, \delta) = 1$ and, therefore,

$$(h_{nn}, (\varepsilon_n, \delta)) = 1.$$

It means that

$$((\varepsilon_n, \delta), k_{n1}\varepsilon_1, ..., k_{n.n-1}\varepsilon_{n-1}) = ((\varepsilon_n, \delta), s_{n1}\varepsilon_1, ..., s_{n.n-1}\varepsilon_{n-1}).$$

Set $SH_A^{-1} = \|t_{ij}\|_1^n$. Since $H_A^{-1} \in \mathbf{G}_E$ taking into account Theorem 2.7, the matrix $H_A^{-1}$ has the form

$$H_A^{-1} = \left\| \begin{array}{ccccc} v_{11} & v_{12} & ... & v_{1.n-1} & v_{1n} \\ \dfrac{\varepsilon_2}{\varepsilon_1}v_{21} & v_{22} & ... & v_{2.n-1} & v_{2n} \\ ... & ... & ... & ... & ... \\ \dfrac{\varepsilon_n}{\varepsilon_1}v_{n1} & \dfrac{\varepsilon_n}{\varepsilon_2}v_{n2} & ... & \dfrac{\varepsilon_n}{\varepsilon_{n-1}}v_{n.n-1} & v_{nn} \end{array} \right\|.$$

Therefore,

$$t_{ni} = \left\| s_{n1} \quad s_{n2} \quad ... \quad s_{nn} \right\| \left\| \begin{matrix} v_{1i} \\ \vdots \\ v_{ii} \\ \dfrac{\varepsilon_{i+1}}{\varepsilon_i} v_{i+1.i} \\ \vdots \\ \dfrac{\varepsilon_n}{\varepsilon_i} v_{ni} \end{matrix} \right\| = s_{n1} v_{1i} + ... +$$

$$+ s_{ni} v_{ii} + s_{n.i+1} \frac{\varepsilon_{i+1}}{\varepsilon_i} v_{i+1.i} + ... + s_{nn} \frac{\varepsilon_n}{\varepsilon_i} v_{ni}, \quad i = 1, ..., n-1.$$

Consider

$$((\varepsilon_n, \delta), t_{n1}\varepsilon_1, ..., t_{n.n-1}\varepsilon_{n-1}) =$$

$$= ((\varepsilon_n, \delta), s_{n1}\varepsilon_1 v_{11} + s_{n2}\varepsilon_2 v_{21} + ... + s_{n.n-1}\varepsilon_{n-1}v_{n-1.1} + s_{nn}\varepsilon_n v_{n1}, ...$$

$$..., s_{n1}\varepsilon_{n-1}v_{1.n-1} + s_{n2}\varepsilon_{n-1}v_{2.n-1} + ... + s_{n.n-1}\varepsilon_{n-1}v_{n-1.n-1} + s_{nn}\varepsilon_n v_{n.n-1}).$$

Since $((\varepsilon_n, \delta), s_{n1}\varepsilon_1, ..., s_{n.n-1}\varepsilon_{n-1})$ is the divisor of all terms, we conclude that

$$((\varepsilon_n, \delta), s_{n1}\varepsilon_1, ..., s_{n.n-1}\varepsilon_{n-1}) \,|\, ((\varepsilon_n, \delta), t_{n1}\varepsilon_1, ..., t_{n.n-1}\varepsilon_{n-1}).$$

On the other hand, $S = \|t_{ij}\|_1^n H_A$. As a result of similar reasoning, we get

$$((\varepsilon_n, \delta), t_{n1}\varepsilon_1, ..., t_{n.n-1}\varepsilon_{n-1}) \,|\, ((\varepsilon_n, \delta), s_{n1}\varepsilon_1, ..., s_{n.n-1}\varepsilon_{n-1}).$$

This yields

$$((\varepsilon_n, \delta), t_{n1}\varepsilon_1, ..., t_{n.n-1}\varepsilon_{n-1}) = ((\varepsilon_n, \delta), s_{n1}\varepsilon_1, ..., s_{n.n-1}\varepsilon_{n-1}).$$

By virtue of the associativity of the ring $M_n(R)$, the proof is complete. □

**Lemma 5.14.** *Let*

$$\frac{\varphi}{(\varphi, \varepsilon_i)} = \mu_i,$$

*where* $\varphi = ((\varepsilon_n, \delta), s_{n1}\varepsilon_1, ..., s_{n.n-1}\varepsilon_{n-1})$. *Then* $\mu_i | s_{ni}$, $i = 1, ..., n-1$.

**Proof.** Since

$$\mu_i = \frac{((\varepsilon_n, \delta), s_{n1}\varepsilon_1, ..., s_{n.n-1}\varepsilon_{n-1})}{((\varepsilon_n, \delta), s_{n1}\varepsilon_1, ..., s_{n.i-1}\varepsilon_{i-1}, s_{ni}\varepsilon_{ni}, ..., s_{n.n-1}\varepsilon_{n-1}, \varepsilon_i)} =$$

$$= \frac{\left(((\varepsilon_n, \delta), s_{n1}\varepsilon_1, ..., s_{n.i-1}\varepsilon_{i-1}), \varepsilon_i \left(s_{ni}, \frac{\varepsilon_{i+1}}{\varepsilon_i}s_{n.i+1}, ..., s_{n.n-1}\frac{\varepsilon_{n-1}}{\varepsilon_i}\right)\right)}{((\varepsilon_n, \delta), s_{n1}\varepsilon_1, ..., s_{n.i-1}\varepsilon_{i-1}, \varepsilon_i)} =$$

$$= \left(\frac{((\varepsilon_n, \delta), s_{n1}\varepsilon_1, ..., s_{n.i-1}\varepsilon_{i-1})}{((\varepsilon_n, \delta), s_{n1}\varepsilon_1, ..., s_{n.i-1}\varepsilon_{i-1}, \varepsilon_i)}, \frac{\varepsilon_i}{((\varepsilon_n, \delta), s_{n1}\varepsilon_1, ..., s_{n.i-1}\varepsilon_{i-1}, \varepsilon_i)} \times \right.$$

$$\left. \times \left(s_{ni}, \frac{\varepsilon_{i+1}}{\varepsilon_i}s_{n.i+1}, ..., s_{n.n-1}\frac{\varepsilon_{n-1}}{\varepsilon_i}\right)\right) =$$

$$= \left(\frac{((\varepsilon_n, \delta), s_{n1}\varepsilon_1, ..., s_{n.i-1}\varepsilon_{i-1})}{((\varepsilon_n, \delta), s_{n1}\varepsilon_1, ..., s_{n.i-1}\varepsilon_{i-1}, \varepsilon_i)}, \left(s_{ni}, \frac{\varepsilon_{i+1}}{\varepsilon_i}s_{n.i+1}, ..., s_{n.n-1}\frac{\varepsilon_{n-1}}{\varepsilon_i}\right)\right),$$

we get $\mu_i | s_{ni}$, $i = 1, ..., n-1$. The lemma is proved. □

**Theorem 5.18.** *Let R be an elementary divisor ring, and let*

$$A \sim E = \mathrm{diag}(\varepsilon_1, \varepsilon_2, ..., \varepsilon_n), \quad B \sim \Delta = \mathrm{diag}(1, ..., 1, \delta)$$

*are invertible matrices. Assume that $P_B P_A^{-1} = \|s_{ij}\|_1^n$, $P_B \in \mathbf{P}_B$, $P_A \in \mathbf{P}_A$.
Then*

$$(A, B)_l = P_B^{-1} \Phi,$$

*where* $\Phi = \mathrm{diag}\,(1, ..., 1, \varphi)$, $\varphi = ((\varepsilon_n, \delta), \varepsilon_1 s_{n1}, ..., \varepsilon_{n-1} s_{n.n-1})$.

**Proof.** Note that, by Lemma 5.13, the element

$$((\varepsilon_n, \delta), s_{n1}\varepsilon_1, ..., s_{n.n-1}\varepsilon_{n-1}),$$

and, hence, the matrix $\Phi$, do not depend on the choice of the transforming
matrices $P_A$ and $P_B$.

Since $\varphi \,|\, \delta$ we get

$$B = (P_B^{-1}\Phi) \left( \mathrm{diag}\left(1, ..., 1, \frac{\delta}{\varphi}\right) Q_B^{-1} \right),$$

i.e., the matrix $D = P_B^{-1}\Phi$ is a left divisor of the matrix $B$.

On the other hand, $\varphi \,|\, \varepsilon_n$. This yields $\Phi|$E. By Lemma 5.14, the element
$\dfrac{\varphi}{(\varphi, \varepsilon_i)}$ is a divisor of $s_{ni}$, $i = 1, ..., n-1$. Consequently, $P_B P_A^{-1} = S \in \mathbf{L}(\mathrm{E}, \Phi)$.
Taking into account Theorem 4.2 we receive $D = P_B^{-1}\Phi$ is a left divisor of
the matrix $A$. Thus, the matrix $D = P_B^{-1}\Phi$ is a left common divisor of the
matrices $A, B$.

Let $T = P_T^{-1}\Gamma Q_T^{-1}$ be another left common divisor of the matrices $A, B$.
Since the Smith form of the matrix $T$ is a divisor of the Smith forms of the
matrices $A, B$, then the matrix $\Gamma$ has the form

$$\Gamma = \mathrm{diag}\,(1, ..., 1, \gamma),$$

where $\gamma \,|\, \varepsilon_n$ та $\gamma \,|\, \delta$. It means that

$$\gamma \,|\, (\varepsilon_n, \delta). \tag{5.23}$$

Based on Theorem 4.2, we have

$$P_T = K_A P_A,$$

where $K_A \in \mathbf{L}(\mathrm{E}, \Gamma)$,

$$P_T = K_B P_B,$$

where $K_B \in \mathbf{L}(\Delta, \Gamma)$. Thus,

$$K_A P_A = K_B P_B \Rightarrow P_B P_A^{-1} = K_B^{-1} K_A.$$

**201**

The matrix $K_B^{-1} K_A$ has the form

$$K_B^{-1} K_A = \left\| \begin{matrix} v_{11} & ... & v_{1.n-1} & v_{1n} \\ ... & ... & ... & ... \\ v_{n-1.1} & ... & v_{n-1.n-1} & v_{n-1.n} \\ \gamma v_{n1} & ... & \gamma v_{n.n-1} & v_{nn} \end{matrix} \right\| \times$$

$$\times \left\| \begin{matrix} u_{11} & ... & u_{1.n-1} & u_{1n} \\ ... & ... & ... & ... \\ u_{n-1.1} & ... & u_{n-1.n-1} & u_{n-1.n} \\ \dfrac{\gamma}{(\gamma, \varepsilon_1)} u_{n1} & ... & \dfrac{\gamma}{(\gamma, \varepsilon_{n-1})} u_{n.n-1} & u_{nn} \end{matrix} \right\| =$$

$$= \left\| \begin{matrix} l_{11} & ... & l_{1.n-1} & l_{1n} \\ ... & ... & ... & ... \\ l_{n-1.1} & ... & l_{n-1.n-1} & l_{n-1.n} \\ \dfrac{\gamma}{(\gamma, \varepsilon_1)} l_{n1} & ... & \dfrac{\gamma}{(\gamma, \varepsilon_{n-1})} l_{n.n-1} & l_{nn} \end{matrix} \right\| = P_B P_A^{-1} = \| s_{ij} \|_1^n .$$

Thus,

$$\frac{\gamma}{(\gamma, \varepsilon_i)} \,|\, s_{ni}, \;\; i = 1, ..., n-1.$$

This yields

$$\gamma \,|\, s_{ni}(\gamma, \varepsilon_i) = (\gamma s_{ni}, \varepsilon_i s_{ni}), \;\; i = 1, ..., n-1.$$

It follows that $\gamma \,|\, \varepsilon_i s_{ni}, i = 1, ..., n-1$. In view of divisibility (5.23), we obtain

$$\gamma \,|\, ((\varepsilon_n, \delta), \varepsilon_1 s_{n1}, ..., \varepsilon_{n-1} s_{n.n-1}) = \varphi.$$

Therefore, $\Gamma \,|\, \Phi$.

Since $P_T = K_B P_B$ and $P_D = P_B$ we get

$$P_T P_D^{-1} = K_B P_B P_B^{-1} = K_B \in \mathbf{L}(\Delta, \Gamma).$$

In view of the fact that $\mathbf{L}(\Delta, \Gamma) = \mathbf{L}(\Phi, \Gamma)$, and Theorem 4.2, we conclude that the matrix $T$ is a left divisor of the matrix $D = P_B^{-1} \Phi$. Thus, the matrix $P_B^{-1} \Phi$ is the left greatest common divisor of the matrices $A$ and $B$. The theorem is proved. $\square$

**Corollary 5.5.** *The matrices $A$ and $B$ are left relatively prime if and only if*

$$((\varepsilon_n, \delta), s_{n1} \varepsilon_1, ..., s_{n.n-1} \varepsilon_{n-1}) = 1. \qquad \square$$

**Corollary 5.6.** *The sets of transforming matrices $(A, B)_l$ and $B$ are connected by the relations:*

*1)* $\mathbf{P}_{(A,B)_l} = \mathbf{G}_\Phi P_B.$

*2)* $\mathbf{P}_B \subseteq \mathbf{P}_{(A,B)_l}.$

**Proof.** Equality 1) follows from the Property 2.2.

Since

$$\mathbf{P}_B = \mathbf{G}_\Delta P_B, \;\; \mathbf{P}_{(A,B)_l} = \mathbf{G}_\Phi P_B,$$

and $\mathbf{G}_\Delta \subseteq \mathbf{G}_\Phi$, we find $\mathbf{P}_B \subseteq \mathbf{P}_{(A,B)_l}$. $\square$

## 5.9. Structure of l.c.m. of matrices

In this subsection we establish the relations between Smith forms of the matrices $A$ and $B$ and Smith forms of their right least common multiple and the corresponding relationships between the transforming matrices of these matrices.

Let $R$ be a Bezout ring of stable range 1.5.

**Lemma 5.15.** *Let*

$$M = P_M^{-1} \operatorname{diag}(\omega_1, \omega_2) Q_M^{-1}, \quad F = P_F^{-1} \operatorname{diag}(\tau_1, \tau_2) Q_F^{-1}$$

*be right common multiples of matrices*

$$A = P_A^{-1} \operatorname{diag}(\varepsilon_1, \varepsilon_2) Q_A^{-1}, \quad B = P_B^{-1} \operatorname{diag}(\delta_1, \delta_2) Q_B^{-1}.$$

*If $\omega_1 | \tau_1$ and $[\varepsilon_2, \delta_2] = \omega_2 | \tau_2$, then $F = MN$.*

**Proof.** According to Theorem 4.2,

$$P_A = L_M P_M,$$

where $L_M \in \mathbf{G}_{\frac{\varepsilon_2}{(\varepsilon_2, \omega_1)}}$,

$$P_B = L_{M_1} P_M,$$

where $L_{M_1} \in \mathbf{G}_{\frac{\delta_2}{(\delta_2, \omega_1)}}$. It follows that $P_M = L_M^{-1} P_A$ and $P_M = L_{M_1}^{-1} P_B$.

Similarly we have

$$P_A = L_F P_F,$$

where $L_F \in \mathbf{G}_{\frac{\varepsilon_2}{(\varepsilon_2, \tau_1)}}$,

$$P_B = L_{F_1} P_F,$$

where $L_{F_1} \in \mathbf{G}_{\frac{\delta_2}{(\delta_2, \tau_1)}}$. Hence, $P_F = L_F^{-1} P_A$ and $P_F = L_{F_1}^{-1} P_B$. Then

$$P_F P_M^{-1} = L_F^{-1} P_A P_A^{-1} L_M = L_F^{-1} L_M =$$

$$= \underbrace{\left\| \begin{matrix} k_{11} & k_{12} \\ \dfrac{\varepsilon_2}{(\varepsilon_2, \tau_1)} k_{21} & k_{22} \end{matrix} \right\|}_{L_F^{-1}} \underbrace{\left\| \begin{matrix} h_{11} & h_{12} \\ \dfrac{\varepsilon_2}{(\varepsilon_2, \omega_1)} h_{21} & h_{22} \end{matrix} \right\|}_{L_M},$$

i.e.,

$$P_F P_M^{-1} = \left\| \begin{matrix} l_{11} & l_{12} \\ \left( \dfrac{\varepsilon_2}{(\varepsilon_2, \tau_1)}, \dfrac{\varepsilon_2}{(\varepsilon_2, \omega_1)} \right) l_{21} & l_{22} \end{matrix} \right\|.$$

By Property 1.7,

$$\left( \dfrac{\varepsilon_2}{(\varepsilon_2, \tau_1)}, \dfrac{\varepsilon_2}{(\varepsilon_2, \omega_1)} \right) = \dfrac{\varepsilon_2}{(\varepsilon_2, [\tau_1, \omega_1])}.$$

Since $\omega_1 | \tau_1$, then

$$\dfrac{\varepsilon_2}{(\varepsilon_2, [\tau_1, \omega_1])} = \dfrac{\varepsilon_2}{(\varepsilon_2, \tau_1)}.$$

It means that

$$P_F P_M^{-1} = \left\| \begin{matrix} l_{11} & l_{12} \\ \dfrac{\varepsilon_2}{(\varepsilon_2, \tau_1)} l_{21} & l_{22} \end{matrix} \right\| = \left\| \begin{matrix} g_{11} & g_{12} \\ g_{21} & g_{22} \end{matrix} \right\| = G. \qquad (5.24)$$

On the other hand,

$$P_F P_M^{-1} = L_{F_1}^{-1} P_B P_B^{-1} L_{M_1}^{-1} = L_{F_1}^{-1} L_{M_1}^{-1}.$$

Similarly, we show that

$$P_F P_M^{-1} = \left\| \begin{matrix} v_{11} & v_{12} \\ \dfrac{\delta_2}{(\delta_2, \tau_1)} v_{21} & v_{22} \end{matrix} \right\| = \left\| \begin{matrix} g_{11} & g_{12} \\ g_{21} & g_{22} \end{matrix} \right\| = G. \qquad (5.25)$$

Taking into account (5.24) and (5.25), we receive $\dfrac{\varepsilon_2}{(\varepsilon_2, \tau_1)} \mid g_{21}$ and $\dfrac{\delta_2}{(\delta_2, \tau_1)} \Big| g_{21}$. Consequently,

$$\left[ \dfrac{\varepsilon_2}{(\varepsilon_2, \tau_1)}, \dfrac{\delta_2}{(\delta_2, \tau_1)} \right] \Big| g_{21}.$$

According to Property 1.9

$$\left[ \dfrac{\varepsilon_2}{(\varepsilon_2, \tau_1)}, \dfrac{\delta_2}{(\delta_2, \tau_1)} \right] = \dfrac{[\varepsilon_2, \delta_2]}{([\varepsilon_2, \delta_2], \tau_1)}.$$

Since $[\varepsilon_2, \delta_2] = \omega_2$, then

$$\dfrac{[\varepsilon_2, \delta_2]}{([\varepsilon_2, \delta_2], \tau_1)} = \dfrac{\omega_2}{(\omega_2, \tau_1)}.$$

Thus, $G \in \mathbf{G}_{\frac{\omega_2}{(\omega_2, \tau_1)}}$. Since $\mathrm{diag}(\omega_1, \omega_2)$ is divisor of matrix $\mathrm{diag}(\tau_1, \tau_2)$, and also taking into account the Theorem 4.2, we obtain that $F = MN$. $\qquad \square$

**Theorem 5.19.** *Let*

$$A \sim \mathrm{E} = \mathrm{diag}(\varepsilon_1, \varepsilon_2), \quad B \sim \Delta = \mathrm{diag}(\delta_1, \delta_2),$$

$P_B P_A^{-1} = \|s_{ij}\|_1^2$, $P_B \in \mathbf{P}_B$, $P_A \in \mathbf{P}_A$. *Then*
*1) if at least one of the matrices $A$, $B$ is not nonsingular or $A$ and $B$ are singular, and moreover $s_{21} \neq 0$, then*

$$[A, B]_r = (L_M P_A)^{-1} \Omega = (L_{M_1} P_B)^{-1} \Omega,$$

*where*

$$\Omega = \left\| \begin{matrix} \omega_1 & 0 \\ 0 & \omega_2 \end{matrix} \right\| = \left\| \begin{matrix} \dfrac{[\varepsilon_1, \delta_1](\varepsilon_2, \delta_2)}{((\varepsilon_2, \delta_2), s_{21}[\varepsilon_1, \delta_1])} & 0 \\ 0 & [\varepsilon_2, \delta_2] \end{matrix} \right\|,$$

*and matrices $L_M, L_{M_1}$ satisfy the equality $L_{M_1}^{-1} L_M = P_B P_A^{-1}$ and belong to groups $\mathbf{G}_{\frac{\varepsilon_2}{(\varepsilon_2, \omega_1)}}$ and $\mathbf{G}_{\frac{\delta_2}{(\delta_2, \omega_1)}}$, respectively.*

2) *if matrices* $A, B$ *are singular, and moreover* $s_{21} = 0$, *then* $\mathbf{P}_A \cap \mathbf{P}_B \neq$ $\neq \varnothing$ *and*
$$[A, B]_r = P^{-1}\Omega,$$
*where*
$$\Omega = \left\| \begin{matrix} [\varepsilon_1, \delta_1] & 0 \\ 0 & 0 \end{matrix} \right\|, \;\; P \in \mathbf{P}_A \cap \mathbf{P}_B.$$

**Proof.** 1). By Lemma 5.9, the element $s_{21}[\varepsilon_1, \delta_1]$), and, hence, the matrix $\Omega$, do not depend on the choice of the transforming matrices $P_A$ and $P_B$.

Show that the matrix $P_B P_A^{-1}$ can be written as

$$P_B P_A^{-1} = KT, \tag{5.26}$$

where
$$K \in \mathbf{G}_{\frac{\delta_2}{(\delta_2, \omega_1)}}, \;\; T \in \mathbf{G}_{\frac{\varepsilon_2}{(\varepsilon_2, \omega_1)}}, \;\; \omega_1 = \frac{[\varepsilon_1, \delta_1] \cdot (\varepsilon_2, \delta_2)}{((\varepsilon_2, \delta_2), s_{21}[\varepsilon_1, \delta_1])}.$$

Using Property 1.6, we get

$$\left( \frac{\varepsilon_2}{(\varepsilon_2, \omega_1)}, \frac{\delta_2}{(\delta_2, \omega_1)} \right) = \frac{(\varepsilon_2, \delta_2)}{(\varepsilon_2, \delta_2, \omega_1)} = \frac{((\varepsilon_2, \delta_2), s_{21}[\varepsilon_1, \delta_1])}{((\varepsilon_2, \delta_2), [\varepsilon_1, \delta_1])} =$$
$$= \left( \frac{(\varepsilon_2, \delta_2)}{((\varepsilon_2, \delta_2), [\varepsilon_1, \delta_1])}, \frac{[\varepsilon_1, \delta_1]}{((\varepsilon_2, \delta_2), [\varepsilon_1, \delta_1])} s_{21} \right) = \left( \frac{(\varepsilon_2, \delta_2)}{((\varepsilon_2, \delta_2), [\varepsilon_1, \delta_1])}, s_{21} \right) = \mu.$$

Hence, $\mu | s_{21}$. Based on Lemma 5.10, the matrix $P_B P_A^{-1}$ can be represented as (5.26). Therefore,
$$K^{-1} P_B = T P_A.$$

Set $K^{-1} = L_{M_1}$, $T = L_M$, we get

$$(L_M P_A)^{-1} \Omega = (L_{M_1} P_B)^{-1} \Omega = M.$$

Since $\mathrm{E}|\Omega$, $\Delta|\Omega$, and $L_M^{-1} \in \mathbf{G}_{\frac{\varepsilon_2}{(\varepsilon_2, \omega_1)}}$, $L_{M_1}^{-1} \in \mathbf{G}_{\frac{\delta_2}{(\delta_2, \omega_1)}}$ on the basis of Theorem 4.2, we conclude that the matrix $M$ is the right common multiple of the matrices $A$ and $B$.

Let $F = P_F^{-1} \Upsilon Q_F^{-1}$, where $\Upsilon = \mathrm{diag}(\tau_1, \tau_2)$, $\tau_1 | \tau_2$ — the other right common multiple of the matrices $A$ and $B$. That is $F = A A_2$, $F = B B_2$. Thus, $\mathrm{E}|\Upsilon$ and $\Delta|\Upsilon$. Since $\varepsilon_2 | \tau_2$ and $\delta_2 | \tau_2$, then $[\varepsilon_2, \delta_2] = \omega_2 | \tau_2$. Since $\varepsilon_1 | \tau_1$ and $\delta_1 | \tau_1$, then $[\varepsilon_1, \delta_1] | \tau_1$, i.e., $\tau_1 = [\varepsilon_1, \delta_1] x$. Moreover, $P_A = K_A P_F$, where $K_A \in \mathbf{G}_{\frac{\varepsilon_2}{(\varepsilon_2, \tau_1)}}$ and $P_B = K_B P_F$, where $K_B \in \mathbf{G}_{\frac{\delta_2}{(\delta_2, \tau_1)}}$. It means that $P_F = K_A^{-1} P_A$ and $P_F = K_B^{-1} P_B$. Consequently,
$$K_B K_A^{-1} = P_B P_A^{-1}.$$

The matrix $K_B K_A^{-1}$ has the form

$$K_B K_A^{-1} = \left\| \begin{matrix} u_{11} & u_{12} \\ \dfrac{\delta_2}{(\delta_2, \tau_1)} u_{21} & u_{22} \end{matrix} \right\| \left\| \begin{matrix} v_{11} & v_{12} \\ \dfrac{\varepsilon_2}{(\varepsilon_2, \tau_1)} v_{21} & v_{22} \end{matrix} \right\| =$$

$$= \left\| \begin{matrix} l_{11} & l_{12} \\ \dfrac{(\varepsilon_2, \delta_2)}{((\varepsilon_2, \delta_2), \tau_1)} l_{21} & l_{22} \end{matrix} \right\| = \left\| \begin{matrix} l_{11} & l_{12} \\ \dfrac{z}{(z, tx)} l_{21} & l_{22} \end{matrix} \right\| = P_B P_A^{-1} = \|s_{ij}\|_1^2,$$

where $z = (\varepsilon_2, \delta_2)$, $t = [\varepsilon_1, \ \delta_1]$. That is $s_{21} = \dfrac{z}{(z, tx)} l_{21}$. Thinking similarly to the proof of Theorem 5.17, we show that $\omega_1 | \tau_1$. Therefore, $\Omega | \Upsilon$. By Lemma 5.15, the matrix $M$ is a left divisor of the matrix $F$. Thus, $M$ is the right least common multiple of the matrices $A$ and $B$.

2). By Lemma 5.12, $s_{21} = 0$ and $\mathbf{P}_A \cap \mathbf{P}_B \neq \varnothing$. Let $U \in \mathbf{P}_A \cap \mathbf{P}_B$. This means that the matrices $A$ and $B$ can be written as

$$A = U^{-1} \mathrm{E} Q_A^{-1}, \ \ B = U^{-1} \Delta Q_B^{-1}.$$

Consider the matrix

$$M = U^{-1} \left\| \begin{matrix} [\varepsilon_1, \delta_1] & 0 \\ 0 & 0 \end{matrix} \right\| = U^{-1} \Omega.$$

Since

$$M = \left( U^{-1} \left\| \begin{matrix} \varepsilon_1 & 0 \\ 0 & 0 \end{matrix} \right\| Q_A^{-1} \right) \cdot \left( Q_A \left\| \begin{matrix} \dfrac{[\varepsilon_1, \delta_1]}{\varepsilon 1} & 0 \\ 0 & 0 \end{matrix} \right\| \right) =$$

$$= \left( U^{-1} \left\| \begin{matrix} \delta_1 & 0 \\ 0 & 0 \end{matrix} \right\| Q_B^{-1} \right) \cdot \left( Q_B \left\| \begin{matrix} \dfrac{[\varepsilon_1, \delta_1]}{\delta_1} & 0 \\ 0 & 0 \end{matrix} \right\| \right),$$

then $M$ is a right common multiple of the matrices $A$ and $B$.

Suppose that $F = P_F^{-1} \Gamma Q_F^{-1}$ is another right common multiple of $A$ and $B$. Therefore, invariant factors of the matrix $F$ are multiples of corresponding invariant factors of the matrices $A$ and $B$. Second invariant factors of these matrices are zeros. Hence, the matrix $\Gamma$ has the form $\Gamma = \mathrm{diag}(\gamma_1, \ 0)$. In addition, $\varepsilon_1 | \gamma_1$ and $\delta_1 | \gamma_1$. Thus $[\varepsilon_1, \delta_1] | \gamma_1$, i.e., $\gamma_1 = [\varepsilon_1, \delta_1] \alpha$. Hence, $\Omega | \Gamma$. Since $F = AF_1$, then $P_A = U = LP_F$, where $L \in \mathbf{L}(\Gamma, \mathrm{E})$ is a group of invertible upper triangular matrices. Noting that $\mathbf{L}(\Gamma, \mathrm{E}) = \mathbf{L}(\Gamma, \Omega)$, we get $U = LP_F$, where $L \in \mathbf{L}(\Gamma, \Omega)$. Taking into account Lemma 5.9, we receive that the matrix $M$ is a left divisor of the matrix $F$. Therefore, $M$ is the right least common multiple of the matrices $A$ and $B$. The theorem is proved. $\square$

**Corollary 5.7.** *If the matrices $A$ and $B$ are singular and $s_{21} \neq 0$, then*

$$M = [A, B]_r = \mathbf{0}. \hspace{3cm} \square$$

Further studies are directed to the study of the right least common multiple of nonsingular matrices with Smith forms

$$\mathrm{E} = \mathrm{diag}(1, ..., 1, \varepsilon), \Delta = \mathrm{diag}(1, ..., 1, \delta)$$

over Bezout rings of stable range 1.5.

**Lemma 5.16.** *Let $A \sim E, B \sim \Delta$ and $P_B \in \mathbf{P}_B$, $P_A \in \mathbf{P}_A$, $P_B P_A^{-1} = = \|s_{ij}\|_1^n$. The element*

$$((\varepsilon, \delta), s_{n1}, s_{n2}, ..., s_{n.n-1})$$

*is invariant with respect to the choice of transforming matrices $P_A$ and $P_B$.*

**Proof.** Let $F_A$ and $F_B$ be other left transforming matrices of the matrices $A$ and $B$. There are $H_A \in \mathbf{G}_E$ and $H_B \in \mathbf{G}_\Delta$ such that $F_A = H_A P_A$, $F_B = = H_B P_B$. Consider the product of matrices

$$F_B F_A^{-1} = H_B P_B (H_A P_A)^{-1} = H_B P_B P_A^{-1} H_A^{-1} = H_B S H_A^{-1},$$

where $S = P_B P_A^{-1}$. In view of the structure of the matrices $H_B$, $H_A^{-1}$ (see Theorem 2.7) and the associativity of the ring $M_n(R)$, the proof of this assertion is evident. $\qquad\square$

**Lemma 5.17.** *Let*

$$\mu = \frac{(\varepsilon, \delta)}{((\varepsilon, \delta), \omega_{n-1})}, \;\; \omega_{n-1} = \frac{(\varepsilon, \delta)}{((\varepsilon, \delta), s_{n1}, s_{n2}, ..., s_{n.n-1})}.$$

Then $\mu \,|\, (s_{n1}, s_{n2}, ..., s_{n.n-1})$.

**Proof.** Since

$$\mu = \frac{(\varepsilon, \delta)}{\left((\varepsilon, \delta), \frac{(\varepsilon, \delta)}{((\varepsilon,\delta), s_{n1}, s_{n2}, ..., s_{n.n-1})}\right)} = \frac{(\varepsilon, \delta) \left((\varepsilon, \delta), s_{n1}, s_{n2}, ..., s_{n.n-1}\right)}{\left((\varepsilon, \delta) \left((\varepsilon, \delta), s_{n1}, s_{n2}, ..., s_{n.n-1}\right), (\varepsilon, \delta)\right)} =$$

$$= ((\varepsilon, \delta), (s_{n1}, s_{n2}, ..., s_{n.n-1})),$$

we find $\mu \,|(s_{n1}, s_{n2}, ..., s_{n.n-1})$. $\qquad\square$

**Lemma 5.18.** *Let $S = \|s_{ij}\|_1^n \in \mathrm{GL}_n(R)$ and*

$$\Omega = \mathrm{diag}(1, ..., 1, \omega_{n-1}, \omega_n),$$

*where $\omega_{n-1} | \omega_n$. In order that there exist invertible matrices*

$$L_A = \left\| \begin{matrix} & & * & & \\ \varepsilon l_{n1} & ... & \varepsilon l_{n.n-2} & \dfrac{\varepsilon}{(\varepsilon, \omega_{n-1})} l_{n.n-1} & l_{nn} \end{matrix} \right\|, \qquad (5.27)$$

$$L_B = \left\| \begin{matrix} & & * & & \\ \delta p_{n1} & ... & \delta p_{n.n-2} & \dfrac{\delta}{(\delta, \omega_{n-1})} p_{n.n-1} & p_{nn} \end{matrix} \right\| \qquad (5.28)$$

*such that $S L_A = L_B$, it is necessary and sufficient that*

$$\frac{(\varepsilon, \delta)}{((\varepsilon, \delta), \omega_{n-1})} \bigg| (s_{n1}, ..., s_{n.n-1}).$$

**Proof. Necessity.** Set
$$\frac{\varepsilon}{(\varepsilon, \omega_{n-1})} = a, \quad \frac{\delta}{(\delta, \omega_{n-1})} = b.$$
It is clear that
$$L_A \in \mathbf{G}_{\mathrm{diag}(1, \ldots, 1, a)} = \mathbf{G}_a \quad \text{and} \quad L_B \in \mathbf{G}_{\mathrm{diag}(1, \ldots, 1, b)} = \mathbf{G}_b.$$
Thus,
$$S = L_B L_A^{-1} \in \mathbf{G}_b \mathbf{G}_a.$$
Hence, we conclude that $(a, b) | s_{ni}, i = 1, \ldots, n-1$. In view of Property 1.6, we write
$$(a, b) = \left( \frac{\varepsilon}{(\varepsilon, \omega_{n-1})}, \frac{\delta}{(\delta, \omega_{n-1})} \right) = \frac{(\varepsilon, \delta)}{((\varepsilon, \delta), \omega_{n-1})}.$$
Thus,
$$\frac{(\varepsilon, \delta)}{((\varepsilon, \delta), \omega_{n-1})} \mid (s_{n1}, \ldots, s_{n.n-1}).$$

**Sufficiency.** First, we consider the case where the matrix $S$ has the form
$$S = \left\| \begin{array}{ccccc} 1 & 0 & \ldots & 0 & 0 \\ 0 & 1 & \ldots & 0 & 0 \\ \ldots & \ldots & \ldots & \ldots & \ldots \\ 0 & 0 & \ldots & 1 & 0 \\ s_{n1} & s_{n2} & \ldots & s_{n.n-1} & 1 \end{array} \right\|.$$
Let $(a, b) = \alpha$. Then there exist elements $t_1, t_2$ such that
$$a t_1 + b t_2 = \alpha.$$
According to the assumption of Lemma, $\alpha \mid (s_{n1}, \ldots, s_{n.n-1})$. So
$$(s_{n1}, \ldots, s_{n.n-1}) = \alpha \beta.$$
there exists an invertible matrix $U$ such that
$$\left\| s_{n1} \quad \ldots \quad s_{n.n-1} \right\| U = \left\| 0 \quad \ldots \quad 0 \quad \alpha \beta \right\|.$$
Then
$$S \left\| \begin{array}{c|c} U & \mathbf{0} \\ \hline 0 \ \ldots \ 0 \ \ -a\beta t_1 & 1 \end{array} \right\| = \left\| \begin{array}{c|c} U & \mathbf{0} \\ \hline 0 \ \ldots \ 0 \ \ b\beta t_2 & 1 \end{array} \right\|.$$
Which was to prove.

Now let $S = \|s_{ij}\|_1^n$ be an arbitrary matrix from $\mathrm{GL}_n(R)$. By Corollary 2.6, there exist $H_1 \in \mathbf{G}_\Delta$ and $T \in \mathbf{G}_\mathrm{E}$ such that
$$H_1 S T = \left\| \begin{array}{cccc|c} 1 & 0 & & 0 & 0 \\ c_{21} & 1 & & 0 & 0 \\ \ldots & \ldots & & & \ldots \\ c_{n-1.1} & \ldots & c_{n-1.n-2} & 1 & 0 \\ \hline c_{n1} & c_{n1} & & c_{n.n-1} & 1 \end{array} \right\| = \left\| \begin{array}{c|c} C_{11} & \mathbf{0} \\ \hline C_{21} & 1 \end{array} \right\|.$$

Consider the matrix
$$H_2 = \left\|\begin{array}{c|c} C_{11}^{-1} & \mathbf{0} \\ \hline \mathbf{0} & 1 \end{array}\right\|.$$

Then
$$H_2 H_1 ST = \left\|\begin{array}{cccc} 1 & 0 & ... & 0 \\ ... & ... & ... & ... \\ 0 & & 1 & 0 \\ c_{n1} & ... & c_{n.n-1} & 1 \end{array}\right\| = C.$$

Note that $H_2 \in \mathbf{G}_\Delta$. Hence, $H_2 H_1 \in \mathbf{G}_\Delta$. It is obvious that $\alpha \,|\, (c_{n1}, ..., c_{n.n-1})$. According to the assertion proved above, there exists the matrix $L'_A$ of the form (5.27) and a matrix $L'_B$ of the form (5.28) such that $CL'_A = L'_B$, i.e.,

$$HSTL'_A = L'_B.$$

Thus,
$$S(TL'_A) = H^{-1}L'_B.$$

In view of Property 4.1, the matrices $TL'_A = L_A$, $H^{-1}L'_B = L_B$ also have the form (5.27) and (5.28), respectively. The proof is completed. □

**Lemma 5.19.** *Let*

$$\mathrm{E} = \mathrm{diag}(1, ..., 1, \varepsilon), \quad \Delta = \mathrm{diag}(1, ..., 1, \delta), \quad U = \|u_{ij}\|_1^n \in \mathrm{GL}_n(R).$$

*Then*
$$\mathrm{E}U\Delta \sim \mathrm{diag}(1, ..., 1, \gamma_{n-1}, \gamma_n). \tag{5.29}$$

**Proof.** Consider the product

$$\mathrm{E}U\Delta = \left\|\begin{array}{ccc|c} u_{11} & ... & u_{1.n-1} & \delta u_{1n} \\ ... & ... & ... & ... \\ u_{n-1.1} & ... & u_{n-1.n-1} & \delta u_{n-1.n} \\ \hline \varepsilon u_{n1} & ... & \varepsilon u_{n.n-1} & \varepsilon\delta u_{nn} \end{array}\right\| = \left\|\begin{array}{cc} U_{11} & U_{12} \\ U_{21} & U_{22} \end{array}\right\|.$$

Since the matrix $U_{11}$ has the order $n-1$ and is a submatrix of an invertible matrix of the order $n$ by Proposition 3.7, we conclude that

$$U_{11} \sim \mathrm{diag}(\underbrace{1, ..., 1}_{n-2}, u).$$

This implies equivalence (5.29). □

**Theorem 5.20.** *Let $A, B$ be nonsingular $n \times n$ matrices over a Bezout ring of stable range 1.5, moreover*

$$A \sim \mathrm{diag}(1, ..., 1, \varepsilon) = \mathrm{E}, \quad B \sim \mathrm{diag}(1, ..., 1, \delta) = \Delta,$$

$P_B \in \mathbf{P}_B$, $P_A \in \mathbf{P}_A$, $P_B P_A^{-1} = \|s_{ij}\|_1^n$. *Then*

$$[A, B]_r = (L_A P_A)^{-1}\Omega = (L_B P_B)^{-1}\Omega,$$

*where*

$$\Omega = \mathrm{diag}\left(1, ...1, \frac{(\varepsilon, \delta)}{((\varepsilon, \delta), s_{n1}, s_{n2}, ..., s_{n.n-1})}, [\varepsilon, \delta]\right),$$

*and the matrices $L_A$ and $L_B$ belong to the sets $\mathbf{L}(\Omega, \mathrm{E})$ and $\mathbf{L}(\Omega, \Delta)$, respectively, and satisfy the equality $(P_B P_A^{-1})L_A = L_B$.*

**Proof.** According to Lemma 5.16, the element

$$((\varepsilon, \delta), s_{n1}, s_{n2}, ..., s_{n.n-1}),$$

and, hence, the matrix $\Omega$ are independent under the choice of the transforming matrices $P_A$ and $P_B$. Consider a matrix

$$\left\| A \quad B \right\| = \left\| P_A^{-1} \mathrm{E} Q_A^{-1} \quad P_B^{-1} \Delta Q_B^{-1} \right\|.$$

Then

$$P_B \left\| A \quad B \right\| \left\| \begin{matrix} Q_A & 0 \\ 0 & Q_B \end{matrix} \right\| = \left\| (P_B P_A^{-1})\mathrm{E} \quad \Delta \right\|.$$

By Corollary 2.6, we choose the matrix $P_B P_A^{-1}$ in a lower unitriangular form. Then

$$\left\| (P_B P_A^{-1})\mathrm{E} \quad \Delta \right\| \sim \left\| \begin{matrix} 1 & & 0 & 0 \\ & \ddots & & \\ * & & 1 & 0 \\ s_{n1} & ... & s_{n.n-1} & \varepsilon \end{matrix} \right\| \begin{matrix} 1 & & 0 & 0 \\ & \ddots & & \\ 0 & & 1 & 0 \\ 0 & ... & 0 & \delta \end{matrix} \right\| \sim$$

$$\sim \mathrm{diag}(1, ..., 1, ((\varepsilon, \delta), s_{n1}, s_{n2}, ..., s_{n.n-1})).$$

Taking into account Theorem 1.10, we receive

$$(A, B)_l \sim \mathrm{diag}(1, ..., 1, ((\varepsilon, \delta), s_{n1}, s_{n2}, ..., s_{n.n-1})).$$

Based on Theorem 1.20, write the matrix $[A, B]_r$ in the form

$$[A, B]_r = BUA_1,$$

where

$$A = (A, B)_l A_1, \quad U \in \mathrm{GL}_n(R).$$

The matrix $A$ has one invariant factor that is not equal to 1. It follows that, its right divisor $A_1$ also has the same number of invariant factors that different from identity. According to Lemma 5.19, $[A, B]_r$ has no more than two invariant factors different from identity:

$$[A, B]_r \sim \mathrm{diag}(1, ..., 1, \omega_{n-1}, \omega_n).$$

It follows from Theorem 1.18,

$$\det A \det B = \det(A, B)_l \det[A, B]_r,$$

i.e.,

$$\det[A, B]_r = \frac{\det A \det B}{\det(A, B)_l} = \frac{\varepsilon\delta}{((\varepsilon, \delta), s_{n1}, s_{n2}, ..., s_{n.n-1})} = \omega_{n-1}\omega_n.$$

The matrices $A, B$ are the left divisors of the matrix $[A, B]_r$. Therefore, their Smith forms divide the matrix $\mathrm{diag}(1, ..., 1, \omega_{n-1}, \omega_n)$. So $[\varepsilon, \delta] \,|\, \omega_n$. On the other hand, it follows from the structure of the matrices of the sets $\mathbf{L}(\Omega, \mathrm{E})$, $\mathbf{L}(\Omega, \Delta)$ that there are no restrictions on the invariant factor $\omega_n$. Hence,

$$\omega_{n-1} = \frac{\varepsilon\delta(\varepsilon, \delta)}{\varepsilon\delta\left((\varepsilon, \delta), s_{n1}, s_{n2}, ..., s_{n.n-1}\right)} = \frac{(\varepsilon, \delta)}{((\varepsilon, \delta), s_{n1}, s_{n2}, ..., s_{n.n-1})}.$$

This yields

$$[A, B]_r \sim \mathrm{diag}\left(1, ...1, \frac{(\varepsilon, \delta)}{((\varepsilon, \delta), s_{n1}, s_{n2}, ..., s_{n.n-1})}, [\varepsilon, \delta]\right) = \Omega.$$

By Property 1.6

$$\left(\frac{\varepsilon}{(\varepsilon, \omega_{n-1})}, \frac{\delta}{(\delta, \omega_{n-1})}\right) = \left(\frac{(\varepsilon, \delta)}{((\varepsilon, \delta), \omega_{n-1})}\right) = \mu.$$

In view of Lemma 5.17, we conclude that

$$\mu \,|\, (s_{n1}, s_{n2}, ..., s_{n.n-1}).$$

According to Lemma 5.18, there exist matrices $L_A \in L(\Omega, \mathrm{E})$, $L_B \in L(\Omega, \Delta)$ such that

$$P_B P_A^{-1} L_A = L_B.$$

It follows that

$$L_A^{-1} P_A \Omega = L_B^{-1} P_B \Omega = M.$$

Since $\mathrm{E}|\Omega$ and $\Delta|\Omega$, in view of Theorem 4.2, the matrix $M$ is a right common multiple of the matrices $A$ and $B$ .

Let $N$ be the right least common multiple of the matrices $A$ and $B$. According to just proved above it follows that $N \sim \Omega$. So, $N = P_N^{-1} \Omega Q_N^{-1}$. Then the matrix

$$M = L_A^{-1} P_A \Omega = P_M^{-1} \Omega$$

is a right multiple of the matrix $N$: $M = N N_1$. According to Theorem 4.2, this is equivalent to $P_N = L P_M$, where $L \in \mathbf{L}(\Omega, \Omega)$. Based on Property 4.6, we get

$$\mathbf{L}(\Omega, \Omega) = \mathbf{G}_\Omega.$$

By virtue of Theorem 4.5, the matrices $M$ and $N$ are right associated. Thus, the matrix $M$ is the right least common multiple of the matrices $A$ and $B$. The theorem is proved. $\qquad\square$

# INVARIANTS OF PRIMITIVE MATRICES WITH RESPECT TO ZELISKO GROUP ACTION

Solving some factorization problems, it becomes necessary to describe all unassociated matrices with the fixed Smith form. An example of this is the description of all left by right unassociated divisors of the matrix

$$A(x) = \left\| \begin{matrix} x^3 & 0 \\ 0 & x^3 \end{matrix} \right\|$$

with the Smith form

$$\Phi(x) = \left\| \begin{matrix} x & 0 \\ 0 & x^2 \end{matrix} \right\|$$

over polynomial ring $F[x]$, where $F$ is field. In this case we have to search all matrices with the Smith form $\Phi(x)$. These will be all matrices of the form $P^{-1}(x)\Phi(x)Q^{-1}(x)$, where $P(x), Q(x)$ are invertible matrices. Indeed

$$A(x) = \big(P^{-1}(x)\text{diag}(x, x^2)Q^{-1}(x)\big)\big(Q(x)\text{diag}(x^2, x)P(x)\big).$$

To establish the fact that matrices are right associates, the Hermite form is usually used. However, if we use this form as a generator of unassociated matrices, then only non-associative matrices with a given determinant can be described. So the set

$$\left\{ \left\| \begin{matrix} 1 & 0 \\ ax^2 + bx + c & x^3 \end{matrix} \right\|, \left\| \begin{matrix} x & 0 \\ ax + b & x^2 \end{matrix} \right\|, \left\| \begin{matrix} x^2 & 0 \\ a & x \end{matrix} \right\|, \left\| \begin{matrix} x^3 & 0 \\ 0 & 1 \end{matrix} \right\| \right\}, \ a, b, c \in F,$$

is the set of all right unassociated of the second order matrices with $\det \Phi(x) = x^3$. Herewith, matrices with the Smith form $\Phi(x)$ will belong to the set

$$\left\{ \left\| \begin{matrix} x & 0 \\ ax + b & x^2 \end{matrix} \right\|, \left\| \begin{matrix} x^2 & 0 \\ a & x \end{matrix} \right\| \right\},$$

where $a, b, c \in F$. However, this set also contains matrices which have no Smith form $\Phi(x)$. In particular, such matrices are

$$\left\|\begin{matrix} x & 0 \\ 1 & x^2 \end{matrix}\right\|, \left\|\begin{matrix} x^2 & 0 \\ 1 & x \end{matrix}\right\|.$$

This simple example has already shown that the Hermite form is a "rough tool" for solving this "delicate" problem. Therefore, there is a need to construct such a canonical form with respect to one-sided transformations, the very appearance of which would speak of the Smith form of the matrix.

Over an elementary divisor ring $R$ any matrix with the Smith form $\Phi$ can be represented as $P^{-1}\Phi Q^{-1}$, where $P, Q \in \mathrm{GL}_n(R)$. Thus, by the right transformations from $\mathrm{GL}_n(R)$, this matrix is reduced to the form $P^{-1}\Phi$. Since the matrix $\Phi$ is invariant by such transformations, it is natural to construct desired form as $P^{-1}\Phi$. Based on Property 2.2, each matrix in the adjacent class $\mathbf{G}_\Phi P$ can be selected as a matrix $P$. That is, the form $P^{-1}\Phi$ is not canonical. Therefore, the task is equivalent to finding "canonical" matrices in the class $\mathbf{G}_\Phi P$. The research of the last two sections of the monograph is devoted to this task.

## 6.1. Φ-rod of column and its properties

Let $P = \|p_{ij}\|_1^n$ be an invertible matrix and $\Phi = \mathrm{diag}\,(\varphi_1, ..., \varphi_n)$ is a nonsingular $d$-matrix. Let us introduce the following notation: $P^m$ is a matrix consisting of the last $m$ rows of the matrix $P$, $1 \le m < n$, $P^m_{i_1, ..., i_k}$ is a matrix consisting of $i_1, ..., i_k$ columns of the matrix $P^m$, $1 \le k \le n$, $1 \le i_1 < ... < i_k \le n$,

$$\delta^m_{i_1, ..., i_k} = \langle P^m_{i_1, ..., i_k} \rangle, \quad \Delta^m_{i_1, ..., i_k} = \langle V^m_{i_1, ..., i_k} \rangle,$$

where $V = HP$.

**Theorem 6.1.** *If $H \in \mathbf{G}_\Phi$, then*

$$\left(\delta^m_{i_1, ..., i_k}, \frac{\varphi_{n-m+1}}{\varphi_{n-m}}\right) = \left(\Delta^m_{i_1, ..., i_k}, \frac{\varphi_{n-m+1}}{\varphi_{n-m}}\right),$$

$m = 1, ..., n-1$.

**Proof.** Let $m < k \le n$ and $\mu$ be an arbitrary minor of the order $m$ of the matrix $P^m_{i_1, ..., i_k}$. The minor of the matrix $V^m_{i_1, ..., i_k}$ which similarly constructed denote by $\nu$. If $k = m$, then the correctness of our assertion is proved in Lemma 2.1. Hence,

$$\left(\mu, \frac{\varphi_s}{\varphi_{s-1}}\right) = \left(\nu, \frac{\varphi_s}{\varphi_{s-1}}\right).$$

It follows that

$$\left(\delta^m_{i_1, ..., i_k}, \frac{\varphi_s}{\varphi_{s-1}}\right) = \left(\Delta^m_{i_1, ..., i_k}, \frac{\varphi_s}{\varphi_{s-1}}\right).$$

Finally, let's consider the case $1 \leq k < m$. Let $\mu_1, ..., \mu_t$ are all minors of the $m$-th order of the matrix $P^m$, containing the submatrix $P^m_{i_1, ..., i_k}$. Denote by $\nu_1, ..., \nu_t$ the corresponding minors of the matrix $V^m$. According to the above

$$\left( \mu_i, \frac{\varphi_s}{\varphi_{s-1}} \right) = \left( \nu_i, \frac{\varphi_s}{\varphi_{s-1}} \right).$$

$i = 1, ..., t$. By Proposition 3.4

$$(\mu_1, ..., \mu_t) = \delta^m_{i_1, ..., i_k}, \ (\nu_1, ..., \nu_t) = \Delta^m_{i_1, ..., i_k}.$$

Thereby

$$\left( \delta^m_{i_1, ..., i_k}, \frac{\varphi_s}{\varphi_{s-1}} \right) = \left( \mu_1, ..., \mu_t, \frac{\varphi_s}{\varphi_{s-1}} \right) = \left( \left( \mu_1, \frac{\varphi_s}{\varphi_{s-1}} \right), ..., \left( \mu_t, \frac{\varphi_s}{\varphi_{s-1}} \right) \right) =$$
$$= \left( \left( \nu_1, \frac{\varphi_s}{\varphi_{s-1}} \right), ..., \left( \nu_t, \frac{\varphi_s}{\varphi_{s-1}} \right) \right) = \left( \nu_1, ..., \nu_t, \frac{\varphi_s}{\varphi_{s-1}} \right) = \left( \Delta^m_{i_1, ..., i_k}, \frac{\varphi_s}{\varphi_{s-1}} \right).$$

The theorem is proved. $\qquad \square$

Let $R$ be an elementary divisor ring and $\Phi = \mathrm{diag}(\varphi_1, ..., \varphi_n)$ is a nonsingular $d$-matrix over $R$. Each invertible matrix consists of primitive columns. Therefore, studying the action of $\mathbf{G}_\Phi$ on invertible matrices it is natural to begin by studying of its actions on primitive columns.

Denote by $\Phi_1 = I_n$,

$$\Phi_i = \mathrm{diag} \left( \frac{\varphi_i}{\varphi_1}, \frac{\varphi_i}{\varphi_2}, ..., \frac{\varphi_i}{\varphi_{i-1}}, \underbrace{1, ..., 1}_{n-i+1} \right), \ \ i = 2, ..., n.$$

**Definition 6.1.** *Let* $\mathbf{a} = \| a_1 \ ... \ a_n \|^T$ *be a primitive column and*

$$\Phi_i \mathbf{a} = \left\| \frac{\varphi_i}{\varphi_1} a_1 \ ... \ \frac{\varphi_i}{\varphi_{i-1}} a_{i-1} \ \ a_i \ ... \ \ a_n \right\|^T \sim \| \delta_i \ 0 \ ... \ 0 \|^T,$$

$i = 1, ..., n$. *The column* $\| \delta_1 \ ... \ \delta_n \|^T$ *is called* $\Phi$-**rod** *of the column* $\mathbf{a}$ *(in notations* $R_\Phi(\mathbf{a})$*).*

**Theorem 6.2.** *If* $H \in \mathbf{G}_\Phi$, *then* $R_\Phi(\mathbf{a}) = R_\Phi(H\mathbf{a})$.

**Proof.** According to Lemma 3.5

$$\Phi_i H \mathbf{a} \overset{l}{\sim} \Phi_i \mathbf{a} \overset{l}{\sim} \| \delta_i \ 0 \ ... \ 0 \|^T,$$

$i = 2, ..., n$. It is also obvious that

$$I_n H \mathbf{a} \overset{l}{\sim} I_n \mathbf{a} \overset{l}{\sim} \| 1 \ 0 \ ... \ 0 \|^T.$$

Consequently, $R_\Phi(\mathbf{a}) = R_\Phi(H\mathbf{a})$. $\qquad \square$

It follows from this Theorem that the Φ-rod is invariant under the actions of $\mathbf{G}_\Phi$ at column $\mathbf{a}$.

Note that $\delta_i | a_i$, $i = 1, ..., n$, and using Theorem 6.2, we get.

**Corollary 6.1.** *If $R_\Phi(\mathbf{a}) = \|\delta_1 \ ... \ \delta_n\|^T$, then the column $\mathbf{a}$ can be written as follows:* $\mathbf{a} = \|\delta_1 b_1 \ \delta_2 b_2 \ ... \ \delta_n b_n\|^T$. $\qquad\qquad\square$

Since $I_n\mathbf{a} \sim \|1 ... 0\|^T$, then $\delta_1 = 1$. We show the relationship between elements of the Φ-rod of the column $\mathbf{a}$.

**Property 6.1.** *Equalities are fulfilled:*

$$\delta_i = \left(\frac{\varphi_i}{\varphi_{i-1}}, \frac{a_i}{\delta_{i-1}}, \frac{a_{i+1}}{\delta_{i-1}}, ..., \frac{a_n}{\delta_{i-1}}\right)\delta_{i-1}, \ \ i = 2, ..., n.$$

**Proof.** It is easy to see that

$$\delta_i = \left(\frac{\varphi_i}{\varphi_1}a_1, ..., \frac{\varphi_i}{\varphi_{i-1}}a_{i-1}, a_i, ..., a_n\right) =$$

$$= \left(\frac{\varphi_i}{\varphi_{i-1}}\left(\frac{\varphi_{i-1}}{\varphi_1}a_1, ..., \frac{\varphi_{i-1}}{\varphi_{i-2}}a_{i-2}\right), a_{i-1}, a_i, ..., a_n\right) =$$

$$= \delta_{i-1}\left(\frac{\varphi_i}{\varphi_{i-1}}\frac{\left(\frac{\varphi_{i-1}}{\varphi_1}a_1, ..., \frac{\varphi_{i-1}}{\varphi_{i-2}}a_{i-2}\right)}{\delta_{i-1}}, \frac{a_{i-1}}{\delta_{i-1}}, ..., \frac{a_n}{\delta_{i-1}}\right) =$$

$$= \delta_{i-1}\left(\frac{\varphi_i}{\varphi_{i-1}}, \frac{a_i}{\delta_{i-1}}, ..., \frac{a_n}{\delta_{i-1}}\right),$$

$i = 2, ..., n.$ $\qquad\qquad\square$

**Corollary 6.2.** *The elements $\delta_i$, $i = 2, ..., n$, of the Φ-rod of the column $\mathbf{a}$ satisfy the conditions*

*1)* $\delta_1 | \delta_2 | ... | \delta_n$,

*2)* $\dfrac{\delta_i}{\delta_{i-1}} \Big| \dfrac{\varphi_i}{\varphi_{i-1}}$, $i = 2, ..., n.$ $\qquad\qquad\square$

Since $\delta_1 = 1$, every Φ-rod is a primitive column. However, not every primitive column will be a Φ-rod of some primitive column. This will only be done under the conditions stated in the following statement.

**Property 6.2.** *Column $\|1 \ \tau_2 \ ... \ \tau_n\|^T$ is the Φ-rod of some primitive column $\|a_1 \ ... \ a_n\|^T$ if and only if*

*1)* $\tau_i \Big| \dfrac{\varphi_i}{\varphi_1}$, $i = 2, ..., n;$

*2)* $\dfrac{\tau_i}{\tau_{i-1}} \Big| \dfrac{\varphi_i}{\varphi_{i-1}}$, $i = 2, ..., n.$

**Proof. Necessity.** Let

$$R_\Phi\left(\|a_1 ... a_n\|^T\right) = \|1\tau_2 ... \tau_n\|^T.$$

Then

$$\tau_i = \left(\frac{\varphi_i}{\varphi_1}a_1, ..., \frac{\varphi_i}{\varphi_{i-1}}a_{i-1}, a_i, ..., a_n\right),$$

$i = 2, ..., n$. According to Property 1.10

$$\tau_i = \left(\frac{\varphi_i}{\varphi_1}, \frac{\varphi_i}{\varphi_1}a_2, ..., \frac{\varphi_i}{\varphi_{i-1}}a_{i-1}, a_i, ..., a_n\right).$$

Therefore $\tau_i \Big| \dfrac{\varphi_i}{\varphi_1}$, $i = 2, ..., n$.

Condition 2) received in Corollary 6.2.

**Sufficiency.** Consider a primitive column

$$\tau = \|1 \ \tau_2 \ ... \ \tau_n\|^T,$$

elements of which satisfy equalities 1) and 2). Let us show that the $\Phi$-rod of this column is the same as the column itself. Let

$$R_\Phi(\tau) = \|\delta_1 \ ... \ \delta_n\|^T.$$

Then

$$\delta_i = \left(\frac{\varphi_i}{\varphi_1}\tau_1, \frac{\varphi_i}{\varphi_2}\tau_2, ..., \frac{\varphi_i}{\varphi_{i-1}}\tau_{i-1}, \tau_i, ..., \tau_n\right) =$$

$$= \left(\frac{\varphi_i}{\varphi_1}, \frac{\varphi_i}{\varphi_2}\tau_2, ..., \frac{\varphi_i}{\varphi_{i-1}}\tau_{i-1}, \tau_i, ..., \tau_n\right).$$

From condition 2) it follows that $\tau_i \,|\, \tau_{i+1} \,|\, ... \,|\, \tau_n$. According to condition 1) $\tau_i \Big| \dfrac{\varphi_i}{\varphi_1}$. Hence,

$$\left(\frac{\varphi_i}{\varphi_1}, \tau_i, \tau_{i+1}, ..., \tau_n\right) = \tau_i.$$

That is

$$\delta_i = \left(\frac{\varphi_i}{\varphi_2}\tau_2, ..., \frac{\varphi_i}{\varphi_{i-1}}\tau_{i-1}, \tau_i\right).$$

Since $\tau_i \Big| \dfrac{\varphi_i}{\varphi_{i-1}}\tau_{i-1}$, and

$$\frac{\varphi_i}{\varphi_{j-1}}\tau_{j-1} = \left(\frac{\varphi_i}{\varphi_j}\tau_j\right)\left(\frac{\varphi_j\tau_{j-1}}{\varphi_{j-1}\tau_j}\right) = \left(\frac{\varphi_i}{\varphi_j}\tau_j\right)\left(\frac{\varphi_j}{\varphi_{j-1}}\Big/\frac{\tau_j}{\tau_{j-1}}\right),$$

$2 \leqslant j - 1 \leqslant i - 1$, then in the sequence

$$\frac{\varphi_i}{\varphi_2}\tau_2, \ ..., \ \frac{\varphi_i}{\varphi_{i-1}}\tau_{i-1}, \ \tau_i$$

each previous element is divided into next. So, $\delta_i = \tau_i$, $i = 1, ..., n$. It means that $R_\Phi(\tau) = \tau$. $\qquad\square$

By Corollary 6.1 each primitive column **a** can be written as

$$\mathbf{a} = \|\delta_1 b_1 \ \ \delta_2 b_2 \ \ ... \ \ \delta_n b_n\|^T,$$

where

$$\|\delta_1 \ \ \delta_2 \ \ ... \ \ \delta_n\|^T = R_\Phi(\mathbf{a}).$$

Since $\delta_i, \ i = 1, ..., n$ are invariants regarding the action of the group $\mathbf{G}_\Phi$, the transformations from this group elements of the column **a** can be replaced at most by $\delta_i$ or zero. Let us transform a primitive column **a** by actions from $\mathbf{G}_\Phi$ to a simplest form.

**Theorem 6.3.** *Let $R$ be a Bezout ring of stable range* 1.5. *If*

$$R_\Phi(\mathbf{a}) = \|\delta_1 \ \ \delta_2 \ \ ... \ \ \delta_n\|^T$$

*is $\Phi$-rod of a primitive column* $\mathbf{a} = \|a_1 \ \ a_2 \ \ ... \ \ a_n\|^T$, *then in the group* $\mathbf{G}_\Phi$ *there exists a matrix $H$, such that*

$$H\mathbf{a} = \|b \ \ \delta_2 \ \ ... \ \ \delta_n\|^T, \ \ (b, \delta_n) = 1.$$

**Proof.** According to Theorem 1.9, there are $u_1, ..., u_n$, such that

$$\frac{\varphi_n}{\varphi_1} a_1 u_1 + ... + \frac{\varphi_n}{\varphi_{n-1}} a_{n-1} u_{n-1} + a_n u_n = \delta_n,$$

where

$$\left( u_n, \frac{\varphi_n}{\varphi_1} \right) = 1.$$

Since

$$\frac{\varphi_n}{\varphi_{n-1}} \ \bigg| \ \frac{\varphi_n}{\varphi_{n-2}} \ \bigg| ... \bigg| \ \frac{\varphi_n}{\varphi_1},$$

based on Property 1.10, we get

$$\left( \frac{\varphi_n}{\varphi_1} u_1, ..., \frac{\varphi_n}{\varphi_{n-1}} u_{n-1}, u_n \right) = \left( \frac{\varphi_n}{\varphi_1}, \frac{\varphi_n}{\varphi_2} u_2, ..., \frac{\varphi_n}{\varphi_{n-1}} u_{n-1}, u_n \right) =$$

$$= \left( \frac{\varphi_n}{\varphi_2} u_2, ..., \frac{\varphi_n}{\varphi_{n-1}} u_{n-1}, \left( u_n, \frac{\varphi_n}{\varphi_1} \right) \right) = 1.$$

By Theorem 1.1, the row

$$\left\| \frac{\varphi_n}{\varphi_1} u_1 \ \ ... \ \ \frac{\varphi_n}{\varphi_{n-1}} u_{n-1} \ \ u_n \right\|$$

can be complemented to an invertible matrix $H_n$ of the form

$$H_n = \left\| \begin{array}{ccccc} u_{11} & u_{12} & ... & u_{1.n-1} & u_{1n} \\ 0 & u_{22} & ... & u_{2.n-1} & u_{2n} \\ ... & ... & ... & ... & ... \\ 0 & ... & 0 & u_{n-1.n-1} & u_{n-1.n} \\ \frac{\varphi_n}{\varphi_1} u_1 & \frac{\varphi_n}{\varphi_2} u_2 & ... & \frac{\varphi_n}{\varphi_{n-1}} u_{n-1} & u_n \end{array} \right\|,$$

which belongs to the group $\mathbf{G}_\Phi$. Then

$$H_n \|a_1 \;\; a_2 \;\; ... \;\; a_n\|^T = \|b_1 \;\; ... \;\; b_{n-1} \;\; \delta_n\|^T.$$

Based on Theorem 6.2, this column has $\Phi$-rod $\|\delta_1 \;\; \delta_2 \;\; ... \;\; \delta_n\|^T$. It means that

$$\Phi_{n-1} H_n \|a_1 \;\; a_2 \;\; ... \;\; a_n\|^T \sim \|\delta_{n-1} \;\; 0 \;\; ... \;\; 0\|^T$$

So there are such elements $v_1, ..., v_n$, such that

$$\frac{\varphi_{n-1}}{\varphi_1} b_1 v_1 + ... + \frac{\varphi_{n-1}}{\varphi_{n-2}} b_{n-2} v_{n-2} + b_{n-1} v_{n-1} + \delta_n v_n = \delta_{n-1}.$$

According to Theorem 1.9, we choose these elements so that

$$(v_1, ..., v_{n-1}) = 1,$$

where

$$\left( v_{n-1}, \frac{\varphi_{n-1}}{\varphi_1} \right) = 1.$$

As in the previous case,

$$\left( \frac{\varphi_{n-1}}{\varphi_1} v_1, ..., \frac{\varphi_{n-1}}{\varphi_{n-2}} v_{n-2}, v_{n-1} \right) = 1.$$

It follows from Theorem 1.1 that in the group $\mathbf{G}_\Phi$ there is a matrix $H_{n-1}$ with the last two rows

$$\left\| \begin{matrix} \dfrac{\varphi_{n-1}}{\varphi_1} v_1 & ... & \dfrac{\varphi_{n-1}}{\varphi_{n-2}} v_{n-2} & v_{n-1} & v_n \\ 0 & ... & 0 & 0 & 1 \end{matrix} \right\|.$$

Then

$$H_{n-1} H_n \|a_1 \;\; a_2 \;\; ... \;\; a_n\|^T = \|c_1 \;\; ... \;\; c_{n-2} \;\; \delta_{n-1} \;\; \delta_n\|^T.$$

Continuing this process, we get such a matrix $H_2 ... H_n \in \mathbf{G}_\Phi$, that

$$H_2 ... H_n a = \|d \;\; \delta_2 \;\; ... \;\; \delta_n\|^T.$$

If $(d, \delta_n) = 1$, the theorem is proved. If this condition is not fulfilled, then with the primitiveness of the column $\mathbf{a}$, and the fact that $\delta_2 \,|\, \delta_3 \,|\, ... \,|\, \delta_n$, it follows that $(d, \delta_2) = 1$. Consequently, $(d, \delta_2, \delta_n) = 1$. Since $\delta_n \neq 0$ and the ring $R$ has stable range of 1.5, there is such $r$, that $(d + r\delta_2, \delta_n) = 1$. Thus in the group $\mathbf{G}_\Phi$ exists the matrix

$$H_1 = \left\| \begin{matrix} 1 & r \\ 0 & 1 \end{matrix} \right\| \oplus I_{n-2},$$

such that

$$H_1 ... H_n \mathbf{a} = H\mathbf{a} = \|b \;\; \delta_2 \;\; ... \;\; \delta_n\|^T,$$

where $b = d + r\delta_2$. The theorem is proved. $\qquad\square$

If a Bezout ring does not have stable range 1.5 Theorem 6.3, in general, does not hold.

**Example 6.1.** Let

$$R = \left\{ a + b_1 x + b_2 x^2 + ... | a \in \mathbb{Z}, \ b_i \in \mathbb{Q}, i \in \mathbb{N} \right\}$$

and

$$\Phi = \left\| \begin{matrix} 1 & 0 \\ 0 & 7x \end{matrix} \right\|, \quad \mathbf{a} = \left\| \begin{matrix} 1 \\ 5x \end{matrix} \right\|.$$

Then

$$R_\Phi(\mathbf{a}) = \left\| \begin{matrix} 1 \\ x \end{matrix} \right\|.$$

However, the group $\mathbf{G}_\Phi$ does not contain such matrix $H$, that

$$H\mathbf{a} = \left\| \begin{matrix} * \\ x \end{matrix} \right\|.$$

Suppose this is not the case. Then the group $\mathbf{G}_\Phi$ contains such a matrix

$$H = \left\| \begin{matrix} h_{11} & h_{12} \\ 7xh_{21} & h_{22} \end{matrix} \right\|,$$

that

$$H\mathbf{a} = \left\| \begin{matrix} * \\ x \end{matrix} \right\|.$$

It follows that

$$7xh_{21} + 5xh_{22} = x.$$

That is

$$7h_{21} + 5h_{22} = 1.$$

By Lemma 1.8

$$\left\| \begin{matrix} h_{21} & h_{22} \end{matrix} \right\| = \left\| \begin{matrix} -2 - 5r & 3 + 7r \end{matrix} \right\|,$$

where $r \in R$ and satisfies the condition

$$\left\| \begin{matrix} 7x(-2 - 5r) & 3 + 7r \end{matrix} \right\| \sim \left\| \begin{matrix} 1 & 0 \end{matrix} \right\|.$$

The element $x$ is relatively prime only with units of the ring $R$, which have the form

$$\pm 1 + b_1 x + b_2 x^2 + ... \, .$$

Consequently,

$$3 + 7r = \pm 1 + b_1 x + b_2 x^2 + ... \, .$$

That is

$$r_1 = -\frac{4}{7} + c_1 x + c_2 x^2 + ...$$

or

$$r_2 = -\frac{2}{7} + d_1 x + d_2 x^2 + ...,$$

which are not elements of the ring $R$. $\diamondsuit$

The elements $\delta_i$ of $\Phi$-rod of the column $\mathbf{a}$ are always not zeros. It means the row
$$H \, \|a_1 \, ... \, a_n\|^T = \|b \, \delta_2 \, ... \, \delta_n\|^T.$$
contains no zeros. However, if, for example,
$$\Phi = \mathrm{diag}(\varphi_1, \varphi_2, \varphi_3), \ \mathbf{a} = \left\| 1 \quad \frac{\varphi_2}{\varphi_1} \quad \frac{\varphi_3}{\varphi_1} \right\|^T,$$
then $R_\Phi(\mathbf{a}) = \mathbf{a}$ and in the group $\mathbf{G}_\Phi$ there exists such matrix
$$H = \left\| \begin{matrix} 1 & 0 & 0 \\ -\dfrac{\varphi_2}{\varphi_1} & 1 & 0 \\ 0 & -\dfrac{\varphi_3}{\varphi_2} & 1 \end{matrix} \right\|,$$
that
$$H\mathbf{a} = \|1 \ \ 0 \ \ 0\|^T.$$
So the question is to find the conditions under which the elements of the column $\mathbf{a}$ by transforming from the group $\mathbf{G}_\Phi$ can be replaced to zero.

**Theorem 6.4.** *Let $R$ be a Bezout ring of stable range* 1.5 *and*
$$R_\Phi(\mathbf{a}) = \|\delta_1 \ \ \delta_2 \ \ ... \ \ \delta_n\|^T.$$
*The group $\mathbf{G}_\Phi$ contains the matrix $K$ such that*
$$K\mathbf{a} = \|b \ \ \delta_2 \ \ ... \ \ \delta_k \ \ 0 \ \ ... \ \ 0\|^T, \ \ 1 \le k < n,$$
*(if $k = 1$, then $b = 1$), if and only if*
$$R_\Phi(\mathbf{a}) = \left\| \delta_1 \ \ ... \ \ \delta_{k-1} \ \ \delta_k \ \ \frac{\varphi_{k+1}}{\varphi_k}\delta_k \ \ \frac{\varphi_{k+2}}{\varphi_k}\delta_k \ \ ... \ \ \frac{\varphi_n}{\varphi_k}\delta_k \right\|^T.$$

**Proof. Necessity.** Acts $\mathbf{G}_\Phi$ do not change the $\Phi$-rod of a column. Therefore
$$R_\Phi\left( \|b \ \ \delta_2 \ \ ... \ \ \delta_k \ \ 0 \ \ ... \ \ 0\|^T \right) = \|\delta_1 \ \ \delta_2 \ \ ... \ \ \delta_n\|^T.$$
Property 6.1 implies that
$$\frac{\delta_i}{\delta_{i-1}} = \left( \frac{\varphi_i}{\varphi_{i-1}}, 0, ..., 0 \right) = \frac{\varphi_i}{\varphi_{i-1}},$$
$i = k+1, k+2, ..., n$. Hence,
$$\delta_{k+1} = \frac{\varphi_{k+1}}{\varphi_k}\delta_k,$$
$$\delta_{k+2} = \frac{\varphi_{k+2}}{\varphi_{k+1}}\frac{\varphi_{k+1}}{\varphi_k}\delta_k = \frac{\varphi_{k+2}}{\varphi_k}\delta_k,$$
$$.....................................................$$
$$\delta_n = \frac{\varphi_n}{\varphi_{n-1}}\frac{\varphi_{n-1}}{\varphi_{n-2}}...\frac{\varphi_{k+1}}{\varphi_k}\delta_k = \frac{\varphi_n}{\varphi_k}\delta_k.$$

**Sufficiency.** Let $k > 1$. By Theorem 6.3, without loss of generality, we can assume that

$$\mathbf{a} = \left\| b \ \ \delta_2 \ \ ... \ \ \delta_k \ \ \frac{\varphi_{k+1}}{\varphi_k}\delta_k \ \ \frac{\varphi_{k+2}}{\varphi_k}\delta_k \ \ ... \ \ \frac{\varphi_n}{\varphi_k}\delta_k \right\|^T,$$

where

$$\left( b, \frac{\varphi_n}{\varphi_k}\delta_k \right) = 1.$$

Then

$$\underbrace{\left( I_{k-1} \oplus \left\| \begin{matrix} 1 & & & & 0 \\ -\dfrac{\varphi_{k+1}}{\varphi_k} & 1 & & & \\ & & \ddots & \ddots & \\ 0 & & & -\dfrac{\varphi_n}{\varphi_{n-1}} & 1 \end{matrix} \right\| \right)}_{K} \mathbf{a} =$$

$$= \left\| b \ \ \delta_2 \ \ ... \ \ \delta_k \ \ 0 \ \ ... \ \ 0 \right\|^T,$$

where $K \in \mathbf{G}_\Phi$.

Let $k = 1$. Then

$$R_\Phi(\mathbf{a}) = \left\| 1 \ \ \frac{\varphi_2}{\varphi_1} \ \ \frac{\varphi_3}{\varphi_1} \ \ ... \ \ \frac{\varphi_n}{\varphi_1} \right\|^T.$$

According to Theorem 6.3, there exists the matrix $K_1 \in \mathbf{G}_\Phi$ such that

$$K_1\mathbf{a} = \left\| b \ \ \frac{\varphi_2}{\varphi_1} \ \ \frac{\varphi_3}{\varphi_1} \ \ ... \ \ \frac{\varphi_n}{\varphi_1} \right\|^T.$$

where

$$\left( b, \frac{\varphi_n}{\varphi_1} \right) = 1.$$

It follows that

$$\left( b, \frac{\varphi_2}{\varphi_1} \right) = 1.$$

There are $u, v \in R$ such that

$$bu + \frac{\varphi_2}{\varphi_1}v = 1.$$

Then

$$\left\| \begin{matrix} u & v & 0 & ... & 0 \\ -\dfrac{\varphi_2}{\varphi_1} & b & 0 & ... & 0 \\ 0 & -\dfrac{\varphi_3}{\varphi_1} & 1 & & 0 \\ & & \ddots & \ddots & \\ 0 & 0 & & -\dfrac{\varphi_n}{\varphi_{n-1}} & 1 \end{matrix} \right\| K_1\mathbf{a} = \left\| 1 \ \ 0 \ \ ... \ \ 0 \right\|^T.$$

The theorem is proved. $\qquad\qquad\square$

Now let $\delta_k$ be the first lower element of the matrix $K\mathbf{a}$ from Theorem 6.4, for which the condition
$$\frac{\delta_k}{\delta_{k-1}} = \frac{\varphi_k}{\varphi_{k-1}}$$
does not hold. This means that by transforming from the group $\mathbf{G}_\Phi$ the element $\delta_k$ cannot be replaced by zero. However, for other elements of this column it is possible. Let $\delta_t$, $2 \leqslant t < k$, is the first lower element of the matrix $K\mathbf{a}$, which satisfies the condition
$$\left(\frac{\varphi_t}{\varphi_{t-1}} \middle/ \frac{\delta_t}{\delta_{t-1}}, \frac{\delta_{t+1}}{\delta_t}\right) = 1.$$

**Theorem 6.5.** *Let $R$ be a Bezout ring of stable range 1.5 and*
$$K\mathbf{a} = \|b \ \ \delta_2 \ \ ... \ \ \delta_k \ \ 0 \ \ ...0\|^T$$
*is a primitive column, moreover*
$$R_\Phi(K\mathbf{a}) = \|\delta_1 \ \ \delta_2 \ \ ... \ \ \delta_n\|^T.$$

*There is a matrix $L$ in $\mathbf{G}_\Phi$ $L$, such that*
$$LK\mathbf{a} = \|b \ \ \delta_2 \ \ ... \ \ \delta_{s-1} \ \ 0 \ \ ... \ \ 0 \ \ \delta_{t+1} \ \ \delta_{t+2} \ \ ... \ \ \delta_k \ \ 0 \ \ ... \ \ 0\|^T, \qquad (6.1)$$
$2 \leqslant t < k$, *if and only if conditions*
$$\left(\frac{\varphi_i}{\varphi_{i-1}} \middle/ \frac{\delta_i}{\delta_{i-1}}, \frac{\delta_{i+1}}{\delta_i}\right) = 1, \qquad (6.2)$$
$i = t, t-1, ..., s$, *satisfied.*

**Proof. Necessity.** Since
$$R_\Phi(K\mathbf{a}) = R_\Phi(LK\mathbf{a}),$$
then
$$\frac{\delta_i}{\delta_{i-1}} = \left(\frac{\varphi_i}{\varphi_i}, \frac{\delta_{t+1}}{\delta_{i-1}}, \frac{\delta_{t+2}}{\delta_{i-1}}, ..., \frac{\delta_k}{\delta_{i-1}}\right) = \left(\frac{\varphi_i}{\varphi_{i-1}}, \frac{\delta_{t+1}}{\delta_{i-1}}\right),$$
$i = t, t-1, ..., s$. It means that
$$\left(\frac{\varphi_i}{\varphi_{i-1}} \middle/ \frac{\delta_i}{\delta_{i-1}}, \frac{\delta_{i+1}}{\delta_i}\right) = 1, \ \ i = t, t-1, ..., s.$$

**Sufficiency.** There are $u_t, v_t$ in $R$ such that
$$u_t \frac{\varphi_t}{\varphi_{t-1}} \middle/ \frac{\delta_t}{\delta_{t-1}} + v_t \frac{\delta_{t+1}}{\delta_t} = 1.$$

Consider the matrix
$$L_t = I_{t-2} \oplus \left\| \begin{matrix} 1 & 0 & 0 \\ \dfrac{\varphi_t}{\varphi_{t-1}} u_t & -1 & v_t \\ 0 & 0 & 1 \end{matrix} \right\| \oplus I_{n-t-1} \in \mathbf{G}_\Phi.$$

At the position $t$ in the column $L_t K \mathbf{a}$ there is an element

$$u_t \frac{\varphi_t}{\varphi_{t-1}} \delta_{t-1} - \delta_t + v_t \delta_{t+1} = \delta_t \left( u_t \frac{\varphi_t}{\varphi_{t-1}} \bigg/ \frac{\delta_t}{\delta_{t-1}} + v_t \frac{\delta_{t+1}}{\delta_t} - 1 \right) = 0.$$

Consequently,

$$L_t K a = \| b \ \ \delta_2 \ \ ... \ \ \delta_{t-1} \ \ 0 \ \ \delta_{t+1} \ \ \delta_{t+2} \ \ ... \ \ \delta_k \ \ ... \ \ 0 \|^T.$$

Again, there are such elements $u_{t-1}, v_{t-1}$ in $R$ such that

$$u_{t-1} \frac{\varphi_{t-1}}{\varphi_{t-2}} \bigg/ \frac{\delta_{t-1}}{\delta_{t-2}} + v_{t-1} \frac{\delta_{t+1}}{\delta_{t-1}} = 1.$$

Consider the matrix

$$L_{t-1} = I_{t-3} \oplus \left\| \begin{matrix} 1 & 0 & 0 & 0 \\ \dfrac{\varphi_{t-1}}{\varphi_{t-2}} u_{t-1} & -1 & 0 & v_{t-1} \\ 0 & 0 & 1 & 0 \\ 0 & 0 & 0 & 1 \end{matrix} \right\| \oplus I_{n-t-1} \in \mathbf{G}_\Phi.$$

Then

$$L_{t-1} L_t K \mathbf{a} = \| b \ \ \delta_2 \ \ ... \delta_{t-2} \ \ 0 \ \ 0 \ \ \delta_{t+1} \ \ \delta_{t+2} \ \ ... \ \ \delta_k \ \ 0 \ \ ... \ \ 0 \|^T.$$

Continuing this process, we reduce the column $\mathbf{a}$ to the form (6.1).

Separately, due to its specificity, consider the case when $s = 2$. Based on previous considerations, there exists $\overline{L} \in \mathbf{G}_\Phi$ such that

$$\overline{L} K \mathbf{a} = \| b \ \ \delta_2 \ \ 0 \ \ 0 \ \ \delta_{t+1} \ \ \delta_{t+2} \ \ ... \ \ \delta_k \ \ 0 \ \ ... \ \ 0 \|^T,$$

where

$$\left( \frac{\varphi_2}{\varphi_1} \bigg/ \frac{\delta_2}{\delta_1}, \frac{\delta_{t+1}}{\delta_2} \right) = 1.$$

Since $(b, \delta_n) = 1$ then

$$\left( b, \frac{\delta_{t+1}}{\delta_2} \right) = 1.$$

Thus

$$\left( b \frac{\varphi_2}{\varphi_1} \bigg/ \frac{\delta_2}{\delta_1}, \frac{\delta_{t+1}}{\delta_2} \right) = 1.$$

This means that there are such elements $u_2, v_2 \in R$ that

$$u_2 b \frac{\varphi_2}{\varphi_1} \bigg/ \frac{\delta_2}{\delta_1} + v_2 \frac{\delta_{t+1}}{\delta_2} = 1.$$

Then the matrix

$$L_2 = \left\| \begin{matrix} 1 & 0 & 0 & ... & 0 & 0 \\ \dfrac{\varphi_2}{\varphi_1} u_2 b & -1 & 0 & ... & 0 & v_2 \\ 0 & 0 & 1 & & & 0 \\ ... & ... & ... & ... & ... & ... \\ 0 & 0 & 0 & & 1 & 0 \\ 0 & 0 & 0 & ... & 0 & 1 \end{matrix} \right\| \oplus I_{n-t-1} \in \mathbf{G}_\Phi$$

will satisfy the equality

$$L_2 \overline{L} K \mathbf{a} = L K \mathbf{a} = \| b \ \ 0 \ \ 0 \ \ \delta_{t+1} \ \ \delta_{t+2} \ \ ... \ \ \delta_k \ \ 0 \ \ ... \ \ 0 \|^T.$$

The theorem is proved. □

If the element $\delta_{s-1}$ of the matrix $LK\mathbf{a}$ from Theorem 6.5 does not satisfy the condition (6.2) then considering the elements that stand above $\delta_{s-1}$ again we are looking for one element $\delta_l$ for which

$$\left( \frac{\varphi_l}{\varphi_{l-1}} \middle/ \frac{\delta_l}{\delta_{l-1}}, \frac{\delta_{l+1}}{\delta_l} \right) = 1,$$

and repeat the reasoning of Theorem 6.5. Having looked at all elements of the column $\mathbf{a}$, eventually they will be left unchanged, or replace them with zero.

We summarize the results obtained by pre-entering a number of notations.

Denote by $(k, ..., l)$ descending sequence of natural numbers $k, k-1, ..., l$. Herewith $(k, ..., k)$ is a one-element, and $(, ..., )$ is an empty sequence.

Let's compare all $\Phi$-rod $\| \delta_1 \delta_2 ... \delta_n \|^T$ a series of sequences

$$(i_1, ..., j_1), (i_2, ..., j_2), ..., (i_p, ..., j_p),$$
$$n \geqslant i_1 \geqslant j_1 > i_2 \geqslant j_2 > ... > i_p \geqslant j_p \geqslant 2,$$

by such a rule:

1a) If

$$\frac{\delta_s}{\delta_{s-1}} = \frac{\varphi_s}{\varphi_{s-1}}, \tag{6.3}$$

$s = n, n-1, ..., q$, herewith $s = q-1$ the condition (6.3) is not satisfied, then $(i_1, ..., j_1) = (n, ..., q)$.

1b) If condition (6.3) is not satisfied for $s = n$ then $(i_1, ..., j_1) = (, ..., )$.

2) In the descending sequence $j_1 - 2, j_1 - 3, ..., 2$ (if $(i_1, ..., j_1) = (, ..., )$ we put $j_1 = n + 1$) we find the first element $t$ for which

$$\left( \frac{\varphi_t}{\varphi_{t-1}} \middle/ \frac{\delta_t}{\delta_{t-1}}, \frac{\delta_{t+1}}{\delta_t} \right) = 1,$$

and we put $i_2 = t$. Let

$$\left( \frac{\varphi_r}{\varphi_{r-1}} \middle/ \frac{\delta_r}{\delta_{r-1}}, \frac{\delta_{t+1}}{\delta_r} \right) = 1, \tag{6.4}$$

$r = t, t-1, ..., l$, wherein for $r = l-1$ the condition (6.4) is not satisfied. Then $j_2 = l$.

3) All other pairs $(i_\mu, ..., j_\mu)$ are given by analogy with the pair $(i_2, ..., j_2)$, and only then consider numeric sequence $j_{\mu-1} - 2, j_{\mu-1} - 3, ..., 2$.

**Theorem 6.6.** *Let $R$ be a Bezout ring of stable range* 1.5 *and* $\mathbf{a}$ *is a primitive column with $\Phi$-rod $\| \delta_1 \ \ \delta_2 \ \ ... \ \ \delta_n \|^T$. Let*

$$(i_1, ..., j_1), (i_2, ..., j_2), ..., (i_p, ..., j_p)$$

be a set of sequences numbers corresponding to the Φ-rod $\|\delta_1 \ \ \delta_2 \ \ ... \ \ \delta_n\|^T$. There is such a matrix $H \in \mathbf{G}_\Phi$ that

$$H \|a_1 \ ...a_n\|^T = \|b \ \ \lambda_2 \ \ ... \ \ \lambda_n\|^T,$$

where $(b, \delta_n) = 1$, and $\lambda_i = 0$ if $i$ belongs to any of the numeric sequences, $\lambda_i = \delta_i$ in all other cases. □

## 6.2. Φ-skeleton of matrices and its properties

Let $R$ be an elementary divisor ring and $\Phi = \mathrm{diag}(\varphi_1, ..., \varphi_n)$ be a nonsingular $d$-matrix over $R$. Recall that

$$\Phi_i = \mathrm{diag}\left(\frac{\varphi_i}{\varphi_1}, \frac{\varphi_i}{\varphi_2}, ..., \frac{\varphi_i}{\varphi_{i-1}}, \underbrace{1, ..., 1}_{n-i+1}\right), \ \ i = 2, ..., n,$$

where $\Phi_1 = I_n$. Let $P \in \mathrm{GL}_n(R)$ and

$$\Phi_i P \overset{l}{\sim} \begin{Vmatrix} \sigma_{i1} & 0 & ... & 0 \\ d_{21}^i & \sigma_{i2} & ... & 0 \\ ... & ... & ... & ... \\ d_{n1}^i & ... & d_{n,n-1}^i & \sigma_{in} \end{Vmatrix} = \mathrm{triang}\,(\sigma_{i1}, \ ..., \ \sigma_{in})$$

are left Hermite forms of matrices $\Phi_i P$, $i = 1, ..., n$.

**Definition 6.2.** *The matrix* $S_\Phi(P) = \|\sigma_{ij}\|_1^n$ *is called* Φ-**skeleton** *of the matrix* $P$.

**Theorem 6.7.** *If* $H \in \mathbf{G}_\Phi$, *then* $S_\Phi(P) = S_\Phi(HP)$.

**Proof.** By Lemma 3.5, $\Phi_i P \overset{l}{\sim} \Phi_i HP$, $i = 2, ..., n$. That is, left Hermite forms of the matrices $\Phi_i P$ and $\Phi_i HP$ coincide, $i = 2, ..., n$. In addition, left Hermite forms of the matrices $\Phi_1 P = I_n P$ and $\Phi_1 HP = I_n HP$ is the identity matrix. □

**Corollary 6.3.** *If matrices* $A = P_A^{-1}\mathrm{E}Q_A^{-1}$, $B = P_B^{-1}\Phi Q_B^{-1}$ *are right associates, then* $S_\Phi(P_A) = S_\Phi(P_B)$.

**Proof.** According to Theorem 4.5, the matrices $A$ and $B$ are right associates if and only if $P_B = HP_A$, where $H \in \mathbf{G}_\Phi$. Then by virtue of Theorem 6.7, $S_\Phi(P_A) = S_\Phi(P_B)$. □

Let us investigate the properties of the Φ-skeleton of a matrix. Denote by

$$\Phi^k = \mathrm{diag}\left(\underbrace{\varphi, ..., \varphi}_{k}, \underbrace{1, ..., 1}_{n-k}\right), \ \ \varphi \neq 0, \ \ k = 1, ..., n-1.$$

**Lemma 6.1.** *If* $P \in \mathrm{GL}_n(R)$ *and*

$$\Phi^k P \overset{l}{\sim} \mathrm{triang}\,(\alpha_{k1}, \ ..., \ \alpha_{kn}), \tag{6.5}$$

*then the following conditions are hold:*

*1)* $\varphi^{k-1} | \alpha_{k\,i_1} \ldots \alpha_{k\,i_{n-1}}$ *for all* $1 \leqslant i_1 < \ldots < i_{n-1} \leqslant n$, $k = 1, \ldots, n-1$;

*2)* $\alpha_{k1} \ldots \alpha_{kn} = \varphi^k e_k$, *where* $e_k \in U(R)$, $k = 1, \ldots, n-1$;

*3)* $\alpha_{ki} | \varphi$, $k = 1, \ldots, n-1$, $i = 1, \ldots, n$.

**Proof.** The matrix

$$\mathrm{diag}\Big(1, \ldots, 1, \underbrace{\varphi, \ldots, \varphi}_{k}\Big).$$

is the Smith form of the matrix $\Phi^k P$. Therefore, $\big\langle \Phi^k P \big\rangle_{n-1} = \varphi^{k-1}$. The Smith forms of the matrices $\Phi^k P$, $\mathrm{triang}\,(\alpha_{k1}, \ldots, \alpha_{kn})$ are coincide. It follows that each $(n-1)$ order minor of the matrix $\mathrm{triang}\,(\alpha_{k1}, \ldots, \alpha_{kn})$ is divided into $\varphi^{k-1}$. In particular,

$$\varphi^{k-1} \left\| \begin{matrix} \alpha_{k\,i_1} & 0 & 0 \\ * & \ddots & 0 \\ * & * & \alpha_{k\,i_{n-1}} \end{matrix} \right\| = \alpha_{k\,i_1} \ldots \alpha_{k\,i_{n-1}},$$

$1 \leqslant i_1 < \ldots < i_{n-1} \leqslant n$.

The equivalence (6.5) implies that

$$\det \Phi^k P = \alpha_{k1} \ldots \alpha_{kn} \varepsilon_k,$$

where $\varepsilon_k \in U(R)$. On the other hand,

$$\det \Phi^k P = \varphi^k \det P = \varphi^k e,$$

where $e = \det V \in U(R)$. Hence,

$$\alpha_{k1} \ldots \alpha_{kn} = \varphi^k e \varepsilon_k^{-1} = \varphi^k e_k,$$

where $e_k = e \varepsilon_k^{-1} \in U(R)$. Consequently,

$$\alpha_{ki} = \frac{\varphi^k}{\alpha_{k1} \ldots \alpha_{k,\,i-1} \alpha_{k,\,i+1} \ldots \alpha_{kn} e_k^{-1}}.$$

Since

$$\varphi^{k-1} \mid \alpha_{k1} \ldots \alpha_{k,\,i-1} \alpha_{k,\,i+1} \ldots \alpha_{kn},$$

then

$$\alpha_{k1} \ldots \alpha_{k,\,i-1} \alpha_{k,\,i+1} \ldots \alpha_{kn} = \varphi^{k-1} s_i.$$

So

$$\alpha_{ki} = \frac{\varphi^k}{\varphi^{k-1} s_i e_k^{-1}} = \frac{\varphi}{s_i e_k^{-1}}.$$

That is $\alpha_{ki} | \varphi$, $k = 1, \ldots, n-1$, $i = 1, \ldots, n$. $\qquad \square$

**Theorem 6.8.** *If* $\| \sigma_{ij} \|_1^n$ *is* $\Phi$-*skeleton of the matrix* $P$ *and*

$$\Delta_{ij} = \frac{\sigma_{ij}}{\sigma_{i-1,j}},$$

*then the following conditions are hold:*

*1)* $\left(\dfrac{\varphi_i}{\varphi_{i-1}}\right)^{i-2} \Big| \Delta_{ij_1} \dots \Delta_{ij_{n-1}},$ *for all* $1 \leqslant j_1 < \dots < j_{n-1} \leqslant n$, $i = 2, \dots, n$,

*2)* $\Delta_{i1} \dots \Delta_{in} = \left(\dfrac{\varphi_i}{\varphi_{i-1}}\right)^{i-1} e_i,\ e_i \in U(R),\ i = 2, \dots, n$,

*3)* $\Delta_{ij} \Big| \dfrac{\varphi_i}{\varphi_{i-1}},\ i = 2, \dots, n,\ j = 1, \dots, n$,

*4)* $\left(\dfrac{\varphi_{i-1}}{\varphi_{i-2}}, \dfrac{\sigma_{i.i-1} \dots \sigma_{in}}{\left(\sigma_{i.i-1} \dots \sigma_{in}, \frac{\varphi_i}{\varphi_{i-1}}\right)}\right) \Big| \sigma_{i-1.i-1} \dots \sigma_{i-1.n},\ i = 3, \dots, n$,

*5)* $\prod\limits_{j=1}^{n} \sigma_{ij} = \dfrac{\varphi_i}{\varphi_1} \dots \dfrac{\varphi_i}{\varphi_{i-1}} e_i,\ e_i \in U(R),\ i = 2, \dots, n$,

*6)* $\sigma_{11} = \dots = \sigma_{1n} = 1$.

**Proof.** Write the matrix $\Phi_i$, $1 < i \leqslant n$ in the form

$$\Phi_i = \operatorname{diag}\left(\dfrac{\varphi_i}{\varphi_{i-1}}, \dots, \dfrac{\varphi_i}{\varphi_{i-1}}, \underbrace{1, \dots, 1}_{n-i+1}\right) \times$$

$$\times \operatorname{diag}\left(\dfrac{\varphi_{i-1}}{\varphi_1}, \dots, \dfrac{\varphi_{i-1}}{\varphi_{i-2}}, 1, \dots, 1\right) = \Phi^i \Phi_{i-1}.$$

According to the definition of the $\Phi$-skeleton, there exists an invertible matrix $U_{i-1}$ such that

$$U_{i-1} \Phi_{i-1} P = \operatorname{triang}(\sigma_{i-1,1}, \dots, \sigma_{i-1,n}) = D_{i-1}.$$

Hence,

$$\Phi_i P = \Phi^i \Phi_{i-1} P = (\Phi^i U_{i-1}^{-1})(U_{i-1} \Phi_{i-1} P) = \Phi^i U_{i-1}^{-1} D_{i-1}.$$

Let $\operatorname{triang}(\alpha_{i1}, \dots, \alpha_{in})$ be left Hermite form of the matrix $\Phi^i U_{i-1}^{-1}$. Then there is $S \in \operatorname{GL}_n(R)$ such that

$$S \Phi^i U_{i-1}^{-1} = \operatorname{triang}(\alpha_{i1}, \dots, \alpha_{in}).$$

So

$$S \Phi_i P = \operatorname{triang}(\alpha_{i1}, \dots, \alpha_{in}) D_{i-1} = \operatorname{triang}(\alpha_{i1}\sigma_{i-1.1}, \dots, \alpha_{in}\sigma_{i-1.n}).$$

On the other hand, $\Phi_i P \overset{l}{\sim} \operatorname{triang}(\sigma_{i1}, \dots, \sigma_{in})$. Therefore,

$$\operatorname{triang}(\sigma_{i1}, \dots, \sigma_{in}) \overset{l}{\sim} \operatorname{triang}(\alpha_{i1}\sigma_{i-1.1}, \dots, \alpha_{in}\sigma_{i-1.n}).$$

It follows that $\sigma_{ij}$ and $\alpha_{ij}\sigma_{i-1,j}$ are associated elements of the ring $R$. Then

$$\alpha_{ij} = \dfrac{\sigma_{ij}}{\sigma_{i-1,j}} e_{ij} = \Delta_{ij} e_{ij},\ e_{ij} \in U(R),$$

and the first three conditions of Theorems follow from Lemma 6.1.

Denote by $U_{i-1}$ a matrix which consists of the last $n - i + 2$ columns of the matrix $P$, and by $P_{i-1}$ a submatrix of an order $n - i + 2$ which is contained in the lower right corner of the matrix $P$, $3 \leqslant i \leqslant n$. From the definition of the $\Phi$-skeleton of the matrix $P$ it follows that $\sigma_{i.i-1} \ldots \sigma_{in} = \langle \Phi_i U_{i-1} \rangle_{n-1}$. Since $\dfrac{\varphi_i}{\varphi_{i-1}} \det V_{i-1}$ is one of such minor, then

$$\sigma_{i.i-1} \ldots \sigma_{in} \ \left| \ \frac{\varphi_i}{\varphi_{i-1}} \det P_{i-1}. \right.$$

Hence,

$$\frac{\sigma_{i.i-1} \ldots \sigma_{in}}{\left( \sigma_{i.i-1} \ldots \sigma_{in}, \frac{\varphi_i}{\varphi_{i-1}} \right)} \ \left| \ \det P_{i-1}. \right. \tag{6.6}$$

Since

$$\Phi_{i-1} = \mathrm{diag}\left( \frac{\varphi_{i-1}}{\varphi_{i-2}}, \ldots, \frac{\varphi_{i-1}}{\varphi_{i-2}}, \underbrace{1, \ldots, 1}_{n-i+2} \right) \Phi_{i-2},$$

then all maximal order minors of the matrix $\Phi_{i-1} U_{i-1}$, except $\det P_{i-1}$, are multiples to $\dfrac{\varphi_{i-1}}{\varphi_{i-2}}$. In view of (6.6), we have

$$\delta_{i-1} = \left( \frac{\varphi_{i-1}}{\varphi_{i-2}}, \frac{\sigma_{i,i-1} \ldots \sigma_{in}}{\left( \sigma_{i,i-1} \ldots \sigma_{in}, \frac{\varphi_i}{\varphi_{i-1}} \right)} \right)$$

is a divisor of $\det P_{i-1}$. Thus, $\delta_{i-1}$ is a divisor of $\sigma_{i-1,i-1} \ldots \sigma_{i-1,n}$.

The equality 5) follows from the definition of the $\Phi$-skeleton of the matrix $A$.

The matrix $\Phi_1 V = I V = V$ is invertible. Therefore, it has the Hermite form the identity matrix $I$. Thus, $\sigma_{11} = \ldots = \sigma_{1n} = 1$. The theorem is proved. $\square$

# Chapter 7

# ONE-SIDED EQUIVALENCE OF MATRICES

*By one-sided transformations of matrices we constructed a canonical form of the form $P^{-1}\Phi$.*

## 7.1. Nonassociative matrices with standard $\Phi$-skeletons

Let $R$ be an elementary divisor ring and $P = \|p_{ij}\|_1^n$ is a matrix over $R$. Denote by

$$P_j = \begin{Vmatrix} p_{1j} & p_{1.j+1} & ... & p_{1n} \\ ... & ... & ... & ... \\ p_{nj} & p_{n.j+1} & ... & p_{nn} \end{Vmatrix},$$

$j = 2, ..., n$. Recall that

$$\Phi_i = \text{diag}\left(\frac{\varphi_i}{\varphi_1}, \frac{\varphi_i}{\varphi_2}, ..., \frac{\varphi_i}{\varphi_{i-1}}, \underbrace{1, ..., 1}_{n-i+1}\right), \quad i = 2, ..., n.$$

**Theorem 7.1.** *Let $P$ be an invertible matrix and $\| \sigma_{ij} \|_1^n$ be its $\Phi$-skeleton. The equation*

$$x\Phi_i P_j = \|a_{ij} \quad p_{i.j+1} \quad ... \quad p_{in}\| \tag{7.1}$$

*has a solution if and only if*

$$p_{ij} \equiv a_{ij} (\text{mod}\,\sigma_{ij}), \quad i = 2, ..., n, j = 1, ..., n.$$

**Proof. Necessity.** Let $\Gamma_j^i = \text{diag}(\gamma_1^i, \gamma_2^i, ..., \gamma_{n-j+1}^i)$ be the Smith form of the matrix $\Phi_i P_j$. The definition of $\Phi$-skeleton implies that

$$\Phi_i P_j \overset{l}{\sim} \text{triang}(\sigma_{ij}, \sigma_{i.j+1}, ..., \sigma_{in}).$$

That is $\sigma_{ij}\sigma_{i.j+1}...\sigma_{in} = \langle\Phi_i P_j\rangle$. On the other hand, $\gamma_1^i\gamma_2^i...\gamma_{n-j+1}^i = \langle\Phi_i P_j\rangle$. Hence,

$$\gamma_1^i\gamma_2^i...\gamma_{n-j+1}^i = \sigma_{ij}\sigma_{i.j+1}...\sigma_{in}. \tag{7.2}$$

**229**

The matrix

$$
\left\|
\begin{array}{cccc}
\dfrac{\varphi_i}{\varphi_1}p_{1j} & \dfrac{\varphi_i}{\varphi_1}p_{1.j+1} & \ldots & \dfrac{\varphi_i}{\varphi_1}p_{1n} \\
\ldots & \ldots & \ldots & \ldots \\
\dfrac{\varphi_i}{\varphi_{i-1}}p_{i-1.j} & \dfrac{\varphi_i}{\varphi_{i-1}}p_{i-1.j+1} & \ldots & \dfrac{\varphi_i}{\varphi_{i-1}}p_{i-1.n} \\
p_{ij} & p_{i.j+1} & \ldots & p_{in} \\
\ldots & \ldots & \ldots & \ldots \\
p_{nj} & p_{n.j+1} & \ldots & p_{nn} \\
a_{ij} & p_{i.j+1} & \ldots & p_{in}
\end{array}
\right\|
$$

be an augmented matrix of equation (7.1). By subtracting from the last row of this matrix its $i$th row we get

$$
\left\|
\begin{array}{c|ccc}
\dfrac{\varphi_i}{\varphi_1}p_{1j} & \dfrac{\varphi_i}{\varphi_1}p_{1.j+1} & \ldots & \dfrac{\varphi_i}{\varphi_1}p_{1n} \\
\ldots & \ldots & \ldots & \ldots \\
\dfrac{\varphi_i}{\varphi_{i-1}}p_{i-1.j} & \dfrac{\varphi_i}{\varphi_{i-1}}p_{i-1.j+1} & \ldots & \dfrac{\varphi_i}{\varphi_{i-1}}p_{i-1.n} \\
p_{ij} & p_{i.j+1} & \ldots & p_{in} \\
\ldots & \ldots & \ldots & \ldots \\
p_{nj} & p_{n.j+1} & \ldots & p_{nn} \\
\hline
a_{ij}-p_{ij} & 0 & \ldots & 0
\end{array}
\right\|
=
\left\|
\begin{array}{cc}
* & \Phi_i P_{j+1} \\
a_{ij}-p_{ij} & \mathbf{0}
\end{array}
\right\|.
$$

From the definition of the $\Phi$-skeleton, it follows that there is an invertible matrix $V_{i.j+1}$ such that

$$
V_{i.j+1}\Phi_i P_{j+1} =
\left\|
\begin{array}{cccc}
\mathbf{0} & \mathbf{0} & & \mathbf{0} \\
\sigma_{i.j+1} & 0 & & 0 \\
* & \sigma_{i.j+1} & & 0 \\
& & \ddots & \\
* & * & & \sigma_{in}
\end{array}
\right\|.
$$

Therefore,

$$
(V_{i.j+1}\oplus 1)
\left\|
\begin{array}{cc}
* & \Phi_i P_{j+1} \\
a_{ij}-p_{ij} & \mathbf{0}
\end{array}
\right\|
=
$$

$$
=
\left\|
\begin{array}{c|cccc}
 & \mathbf{0} & \mathbf{0} & & \mathbf{0} \\
 & \sigma_{i.j+1} & 0 & & 0 \\
* & * & \sigma_{i.j+1} & & 0 \\
 & & & \ddots & \\
 & * & * & & \sigma_{in} \\
\hline
a_{ij}-p_{ij} & & & \mathbf{0} &
\end{array}
\right\|
= L_{ij}.
$$

Obviously, the augmented matrix and the matrix $L_{ij}$ are left associates. According to Theorem 2.6, equation (7.1) has a solution if and only if the invariant factors of the augmented matrix are associates with the corresponding invariant factors of the matrix $\Gamma_j^i$. It follows that

$$\langle L_{ij} \rangle = \gamma_1^i \gamma_2^i ... \gamma_{n-j+1}^i.$$

In view of (7.2), we have

$$\langle L_{ij} \rangle = \sigma_{ij} \sigma_{i,j+1} ... \sigma_{in}.$$

Consequently, $\sigma_{ij}\sigma_{i.j+1} ... \sigma_{in}$ is a divisor of all maximal order minors of the matrix $L_{ij}$. In particular,

$$\sigma_{ij}\sigma_{i.j+1} ... \sigma_{in} | (a_{ij} - p_{ij})\sigma_{i.j+1} ... \sigma_{in}.$$

That is $\sigma_{ij} | (a_{ij} - p_{ij})$.

**Sufficiency.** Let $a_{ij} = p_{ij} + r\sigma_{ij}$. According to the definition of Φ-skeleton, we get

$$V_i \operatorname{diag}\left(\frac{\varphi_i}{\varphi_1}, ..., \frac{\varphi_i}{\varphi_{i-1}}, 1, ..., 1\right) P = \operatorname{triang}(\sigma_{i1}, \sigma_{i2}, ..., \sigma_{in}),$$

where $V_i \in \mathrm{GL}_n(R)$. Let $\|v_{j1}^i \quad v_{j2}^i \quad ... \quad v_{jn}^i\|$ be the $j$th row of the matrix $V_i$. Then the equalities

$$\|r \quad 1\| \left\|\begin{matrix} v_{j1}^i & ... & v_{j.i-1}^i & v_{ji}^i & v_{j.i+1}^i & ... & v_{jn}^i \\ 0 & ... & 0 & 1 & 0 & ... & 0 \end{matrix}\right\| \Phi_i P_j =$$

$$= \|r \quad 1\| \left\|\begin{matrix} \sigma_{ij} & 0 & ... & 0 \\ p_{ij} & p_{i.j+1} & ... & p_{in} \end{matrix}\right\| = \|p_{ij} + r\sigma_{ij} \quad p_{i.j+1} \quad ... \quad p_{in}\|.$$

are fulfilled. Consequently,

$$\|r \quad 1 \quad \| \left\|\begin{matrix} v_{j1}^i & ... & v_{j.i-1}^i & v_{ji}^i & v_{j.i+1}^i & ... & v_{jn}^i \\ 0 & ... & 0 & 1 & 0 & ... & 0 \end{matrix}\right\| =$$

$$= \|rv_{j1}^i \quad ... \quad rv_{j.i-1}^i \quad rv_{ji}^i + 1 \quad rv_{j.i+1}^i \quad ... \quad rv_{jn}^i\|$$

is the desired solution of equation (7.1). The theorem is proved. $\square$

Let $B$ be a nonsingular matrix with the Smith form $\Phi = \operatorname{diag}(\varphi_1, ..., \varphi_n)$, and $B = P^{-1}\Phi Q^{-1}$. Suppose that the set of its left transforming matrices $\mathbf{P}_B$ contains a lower unitriangular matrix $P_0$. It follows that

$$S_\Phi(P_0) = \left\|\begin{matrix} 1 & 1 & ... & 1 & 1 \\ \dfrac{\varphi_2}{\varphi_1} & 1 & ... & 1 & 1 \\ ... & ... & ... & ... & ... \\ \dfrac{\varphi_n}{\varphi_1} & \dfrac{\varphi_n}{\varphi_2} & ... & \dfrac{\varphi_n}{\varphi_{n-1}} & 1 \end{matrix}\right\| = F(\Phi).$$

By Theorem 6.7, all left transforming matrices of $B$ has the Φ-skeleton $F(\Phi)$. Therefore, there will be no confusion if we identify the Φ-skeleton of the matrix $P$ with the $\Phi$-skeleton of the matrix $B$. So we will assume $S_\Phi(B) = S_\Phi(P)$, where $P \in \mathbf{P}_B$.

Denote by $\mathbf{T}(\Phi)$ the set of lower unitriangular matrices of the form

$$\left\| \begin{array}{ccccc} 1 & 0 & ... & 0 & 0 \\ t_{21} & 1 & ... & 0 & 0 \\ ... & ... & ... & ... & ... \\ t_{n1} & t_{n2} & ... & t_{n,n-1} & 1 \end{array} \right\|,$$

where $t_{ij} \in K\left(\dfrac{\varphi_i}{\varphi_j}\right)$, $i > j$. Consider $S_n$ be the symmetric group of degree $n$ and

$$\tau = \begin{pmatrix} 1 & ... & n \\ i_1 & ... & i_n \end{pmatrix} \in S_n.$$

Let $E(\tau) = \|\delta_{kj}\|_1^n$, where

$$\delta_{kj} = \begin{cases} 1, & j = i_k, \\ 0, & j \neq i_k, \end{cases}$$

$k = 1, ..., n$, is matrix presentation of $\tau$. It is known that

$$E(\sigma)E(\tau) = E(\tau\sigma).$$

We say that the matrix $B$ has a **standard** $\Phi$ skeleton if

$$S_\Phi(B) = F(\Phi)E(\tau), \ \ \tau \in S_n.$$

**Theorem 7.2.** *The set $\mathbf{T}(\Phi)$ consists of representatives of the various left cosets of $\mathrm{GL}_n(R)$ at $\mathbf{G}_\Phi$.*

**Proof.** Consider the matrix $\Psi = \varphi_n I_n$. Then $\mathbf{V}(\Psi, \Phi) = \mathbf{T}(\Phi)$. By Theorem 5.6, this set consists of representatives of the various left cosets of $\mathrm{GL}_n(R)$ at $\mathbf{G}_\Phi$. $\qquad \square$

Theorems 7.2 and 5.6, imply that $\mathbf{T}^{-1}(\Phi)\Phi$ consists of right unassociated matrices with the $\Phi$-skeleton $F(\Phi)$. To show that the inverse statement is correct, we establish a number of auxiliaries facts.

**Lemma 7.1.** *Let $S_\Phi(P) = \|\sigma_{ij}\|_1^n$, moreover $\sigma_{nk} = 1$, $1 \leqslant k \leqslant n$. Then there is $H \in \mathbf{G}_\Phi$ such that $\|0 \ ... \ 0 \ 1\|^T$ is the $k$th column of $HP$.*

**Proof.** Let $P = \|p_{ij}\|_1^n$ and

$$\Phi_n \|p_{1k} \ ... \ p_{nk}\|^T =$$

$$= \mathrm{diag}\left(\frac{\varphi_n}{\varphi_1}, ..., \frac{\varphi_n}{\varphi_{n-1}}, 1\right)\| p_{1k} \ ... \ p_{nk}\|^T \overset{l}{\sim} \|\alpha \ 0 \ ...0\|^T.$$

From the definition of the $\Phi$-skeleton it follows that

$$\Phi_n \left\| \begin{array}{ccc} p_{1.k+1} & ... & p_{1n} \\ ... & ... & ... \\ p_{n.k+1} & ... & p_{nn} \end{array} \right\| = \Phi_n P_{k+1} \overset{l}{\sim} \mathrm{triang}\,(\sigma_{n.k+1}, ..., \sigma_{nn}).$$

It means that $\sigma_{n.k+1} \ldots \sigma_{nn}| \langle \Phi_n P_{k+1} \rangle$. Since $\alpha| \left\langle \Phi_n \, \| p_{1k} \; \ldots \; p_{nk} \|^T \right\rangle$, then

$$\alpha \sigma_{n.k+1} \ldots \sigma_{nn} \left\langle \Phi_n \begin{Vmatrix} p_{1k} & p_{1,k+1} & \ldots & p_{1n} \\ \ldots & \ldots & \ldots & \ldots \\ p_{nk} & p_{n,k+1} & \ldots & p_{nn} \end{Vmatrix} \right\rangle = \Phi_n P_k \quad .$$

Noting that

$$\langle \Phi_n P_k \rangle = \sigma_{nk} \sigma_{n.k+1} \ldots \sigma_{nn} = \sigma_{n.k+1} \ldots \sigma_{nn},$$

we get

$$\alpha \sigma_{n.k+1} \ldots \sigma_{nn} | \sigma_{n.k+1} \ldots \sigma_{nn}.$$

Therefore $\alpha \in U(R)$. Since $\alpha$ is g.c.d. of the matrix $\Phi_n \, \| p_{1k} \; \ldots \; p_{nk} \|^T$ elements and is determined up to associativity, then $\alpha = 1$. By Lemma 5.4, there is $H \in \mathbf{G}_\Phi$ such that

$$H \, \| p_{1k} \; \ldots \; p_{nk} \|^T = \| 0 \; \ldots \; 0 \; 1 \|^T. \qquad \square$$

**Lemma 7.2.** *Let $C$ be $((n-1) \times k)$ matrix $(n-1) \geqslant k$, moreover*

$$C \overset{l}{\sim} \operatorname{triang}(\gamma_1, \gamma_2, \ldots, \gamma_k),$$

*where* $\gamma_1 \gamma_2, \ldots \gamma_k \neq 0$. *If*

$$\left\langle \begin{matrix} & C & \\ a_1 & a_2 & \ldots & a_k \end{matrix} \right\rangle = \gamma_1 \gamma_2 \ldots \gamma_k,$$

$a_i \in R$, $i = 1, \ldots, k$, *then*

$$\begin{Vmatrix} & C & \\ a_1 & a_2 & \ldots & a_k \end{Vmatrix} \overset{l}{\sim} \operatorname{triang}(\gamma_1, \gamma_2, \ldots, \gamma_k).$$

**Proof.** There is $V \in \operatorname{GL}_{n-1}(R)$ such that

$$VC = \operatorname{triang}(\gamma_1, \gamma_2, \ldots, \gamma_k) = \begin{Vmatrix} \mathbf{0} & \mathbf{0} & \ldots & \mathbf{0} \\ \gamma_1 & 0 & \ldots & 0 \\ * & \gamma_2 & & 0 \\ & & \ddots & \\ * & * & & \gamma_k \end{Vmatrix}.$$

Then

$$(V \oplus I_1) \begin{Vmatrix} & C & \\ a_1 & a_2 & \ldots & a_k \end{Vmatrix} = \begin{Vmatrix} \mathbf{0} & \mathbf{0} & \ldots & \mathbf{0} \\ \gamma_1 & 0 & \ldots & 0 \\ * & \gamma_2 & & 0 \\ & & \ddots & \\ * & * & & \gamma_k \\ a_1 & a_2 & \ldots & a_k \end{Vmatrix}.$$

This matrix contains a maximal order minor $\gamma_1 \gamma_2 \ldots \gamma_{k-1} a_k$. By assumption, $\gamma_1 \gamma_2 \ldots \gamma_k | \gamma_1 \gamma_2 \ldots \gamma_{k-1} a_k$. It follows that $a_k = \gamma_k a_k'$. So the relation

$$
\left\|
\begin{array}{cccc}
& I_{n-1} & & \mathbf{0} \\
0 & \ldots & 0 & -a_k' & 1
\end{array}
\right\|
\left\|
\begin{array}{cccc}
\mathbf{0} & \mathbf{0} & \ldots & \mathbf{0} \\
\gamma_1 & 0 & \ldots & 0 \\
* & \gamma_2 & & 0 \\
& & \ddots & \\
* & * & & \gamma_k \\
a_1 & a_2 & \ldots & a_k
\end{array}
\right\| =
$$

$$
= \left\|
\begin{array}{cccc}
\mathbf{0} & \mathbf{0} & \ldots & \mathbf{0} \\
\gamma_1 & 0 & \ldots & 0 \\
& \ddots & & \\
* & & \gamma_{k-1} & 0 \\
* & & * & \gamma_k \\
b_1 & \ldots & b_{k-1} & 0
\end{array}
\right\|
\overset{l}{\sim} \operatorname{triang}(\gamma_1, \gamma_2, \ldots, \gamma_k).
$$

is fulfilled. This matrix contains a minor $\gamma_1 \ldots \gamma_{k-2} b_{k-1} \gamma_k$, which is a multiple of $\gamma_1 \ldots \gamma_k$. Therefore,

$$
b_{k-1} = \gamma_{k-1} b_{k-1}'
$$

and

$$
\left\|
\begin{array}{ccccc}
& I_{n-1} & & & \mathbf{0} \\
0 & \ldots & 0 & -b_{k-1}' & 0 & 1
\end{array}
\right\|
\left\|
\begin{array}{cccc}
\mathbf{0} & \mathbf{0} & \ldots & \mathbf{0} \\
\gamma_1 & 0 & \ldots & 0 \\
& \ddots & & \\
* & & \gamma_{k-1} & 0 \\
* & & * & \gamma_k \\
b_1 & \ldots & b_{k-1} & 0
\end{array}
\right\| =
$$

$$
= \left\|
\begin{array}{cccc}
\mathbf{0} & \mathbf{0} & \mathbf{0} & \mathbf{0} \\
\gamma_1 & 0 & 0 & 0 \\
& \ddots & & \\
* & \gamma_{k-2} & 0 & 0 \\
* & * & \gamma_{k-1} & 0 \\
* & * & * & \gamma_k \\
c_1 & \ldots & c_{k-2} & 0 & 0
\end{array}
\right\|.
$$

Continuing this process, we find an invertible matrix $L$ such that

$$
L \left\|
\begin{array}{cccc}
\mathbf{0} & \mathbf{0} & \ldots & \mathbf{0} \\
\gamma_1 & 0 & \ldots & 0 \\
* & \gamma_2 & & 0 \\
& & \ddots & \\
* & * & & \gamma_k \\
a_1 & a_2 & \ldots & a_k
\end{array}
\right\| =
\left\|
\begin{array}{cccc}
\mathbf{0} & \mathbf{0} & \ldots & \mathbf{0} \\
\gamma_1 & 0 & \ldots & 0 \\
* & \gamma_2 & & 0 \\
& & \ddots & \\
* & * & & \gamma_k \\
0 & 0 & \ldots & 0
\end{array}
\right\|.
$$

Noting that the resulting matrix is left associated to the matrix

$$\left\| \begin{matrix} & C & & \\ a_1 & a_2 & ... & a_k \end{matrix} \right\|.$$

we complete the proof. □

Denote by $\bar{b}_s$ a column of height $n-1$ of the form

$$\bar{b}_s = \left\| \underbrace{0 \ ... \ 0 \ 1}_{s} \ * \ ... \ * \right\|^T,$$

$1 \leqslant s \leqslant n-1$. Consider the matrix

$$B = \left\| \bar{b}_{j_1} \ ... \ \bar{b}_{j_k} \right\|,$$

$1 \leqslant j_1, ..., j_k \leqslant n-1$, where $j_l \neq j_t$ herewith $l \neq t$.

**Lemma 7.3.** *If*

$$\operatorname{diag}\left( \frac{\varphi_n}{\varphi_1}, ..., \frac{\varphi_n}{\varphi_{n-1}} \right) B \overset{l}{\sim} \operatorname{triang}\left( \frac{\varphi_n}{\varphi_{j_1}}, ..., \frac{\varphi_n}{\varphi_{j_k}} \right),$$

*moreover*

$$\Phi_n \left\| \begin{matrix} & B & \\ a_1 & ... & a_k \end{matrix} \right\| \overset{l}{\sim} \operatorname{triang}\left( \frac{\varphi_n}{\varphi_{j_1}}, ..., \frac{\varphi_n}{\varphi_{j_k}} \right),$$

$a_i \in R, \ i = 1, ..., k,$ *then there is* $H \in \mathbf{G}_\Phi$ *such that*

$$H \left\| \begin{matrix} & B & \\ a_1 & ... & a_k \end{matrix} \right\| = \left\| \begin{matrix} B \\ \mathbf{0} \end{matrix} \right\|. \tag{7.3}$$

**Proof.** Let $j_s$ be the least indices in the set $j_1, ..., j_k$. Reasoning as in the proof of Lemma 7.2, we get

$$a_s = \frac{\varphi_n}{\varphi_{j_s}} a'_s.$$

Then the matrix

$$H_s = \left\| \begin{matrix} & & & I_{n-1} & & & & \mathbf{0} \\ 0 & ... & 0 & -\dfrac{\varphi_n}{\varphi_{j_s}} a'_s & 0 & ... & 0 & 1 \end{matrix} \right\|,$$

containing the element $-\dfrac{\varphi_n}{\varphi_{j_s}} a'_s$ at the position $(n, j_s)$, will belong to the group $\mathbf{G}_\Phi$, and fulfills the condition

$$H_s \left\| \begin{matrix} \bar{b}_{j_s} \\ a_s \end{matrix} \right\| = \left\| \begin{matrix} \bar{b}_{j_s} \\ 0 \end{matrix} \right\|.$$

From remaining indexes we choose a smallest one, namely $j_t$. If

$$H_s \left\| \begin{matrix} \bar{b}_{j_t} \\ a_t \end{matrix} \right\| = \left\| \begin{matrix} \bar{b}_{j_t} \\ b_t \end{matrix} \right\|,$$

we show that
$$b_t = \frac{\varphi_n}{\varphi_{j_t}} b'_t$$

and construct the matrix
$$H_t = \left\| \begin{matrix} & I_{n-1} & & \mathbf{0} \\ 0 & \ldots & 0 & -\dfrac{\varphi_n}{\varphi_{j_t}} b'_t & 0 & \ldots & 1 \end{matrix} \right\|,$$

in which the element $-\dfrac{\varphi_n}{\varphi_{j_t}} b'_t$ is at the position $(n, j_t)$. Then
$$H_t \left\| \begin{matrix} \overline{b}_{j_t} \\ b_t \end{matrix} \right\| = \left\| \begin{matrix} \overline{b}_{j_t} \\ 0 \end{matrix} \right\|.$$

Since $j_t < j_s$, then the $j_t$th row of the matrix
$$\left\| \begin{matrix} \overline{b}_{j_t} & \overline{b}_{j_s} \\ b_t & 0 \end{matrix} \right\|$$

has the form $\| 1 \ 0 \|$. Thus,
$$H_t H_s \left\| \begin{matrix} \overline{b}_{j_t} & \overline{b}_{j_s} \\ a_t & a_s \end{matrix} \right\| = H_t \left\| \begin{matrix} \overline{b}_{j_t} & \overline{b}_{j_s} \\ b_t & 0 \end{matrix} \right\| = \left\| \begin{matrix} \overline{b}_{j_t} & \overline{b}_{j_s} \\ 0 & 0 \end{matrix} \right\|,$$

moreover $H_t H_s \in \mathbf{G}_\Phi$. Continuing such transformations of the matrix
$$\left\| \begin{matrix} & B & \\ a_1 & \ldots & a_k \end{matrix} \right\|,$$

we get the matrix $H \in \mathbf{G}_\Phi$ which will satisfy equality (7.3). Lemma is proved. $\square$

Let $P = \|p_{ij}\|_1^n$ be an invertible matrix and $S_\Phi(P) = \|\sigma_{ij}\|_1^n$ be its $\Phi$-skeleton. By removing in this matrices $n$th row and $r$th column, we get the matrices $P_{nr}$, and $S_\Phi(P)_{nr}$, respectively.

**Lemma 7.4.** *If*
$$\|p_{1r} \ldots p_{nr}\|^T = \|0 \ldots 0 \ 1\|^T,$$

$1 \leqslant r \leqslant n$, i $\sigma_{1r} = \ldots = \sigma_{nr} = 1$, then
$$S_{\mathrm{diag}\,(\varphi_1, \ldots, \varphi_{n-1})}(P_{nr}) = S_\Phi(P)_{nr}. \tag{7.4}$$

**Proof.** Suppose that $r = n$. By a definition of the $\Phi$-skeleton,
$$\Phi_i P = \Phi_i \left\| \begin{matrix} P_{nn} & \mathbf{0} \\ * & 1 \end{matrix} \right\| \overset{l}{\sim} \mathrm{triang}\,(\sigma_{i1}, \ldots, \sigma_{i.n-1}, 1),$$

$i = 1, \ldots, n$. Note that for any $n \times k$ matrix of the form
$$\left\| \begin{matrix} U & \mathbf{0} \\ * & 1 \end{matrix} \right\|,$$

$n \geqslant k$, the equality

$$\left\langle \begin{matrix} U & \mathbf{0} \\ * & 1 \end{matrix} \right\rangle = \langle U \rangle$$

holds. It follows that the statement is correct for $r = n$.

Consider the case $1 \leqslant r < n$. Since $\sigma_{ir} = 1$, then g.c.d. of maximal order minors of matrices

$$\Phi_i \left\| \begin{matrix} 0 & p_{1.r+1} & ... & p_{1n} \\ ... & ... & ... & ... \\ 0 & p_{n-1.r+1} & ... & p_{n-1.n} \\ 1 & p_{n.r+1} & ... & p_{nn} \end{matrix} \right\|, \quad \Phi_i \left\| \begin{matrix} p_{1.r+1} & ... & p_{1n} \\ ... & ... & ... \\ p_{n-1.r+1} & ... & p_{n-1,n} \\ p_{n.r+1} & ... & p_{nn} \end{matrix} \right\|,$$

$i = 1, ..., n-1$ coincide and equal to $\sigma_{i.r+1} ... \sigma_{in}$. According to the remark made above, g.c.d. of maximal order minors of the first matrix is equal to

$$\left\langle \bar{\Phi}_i \left\| \begin{matrix} p_{1.r+1} & ... & p_{1n} \\ ... & ... & ... \\ p_{n-1.r+1} & ... & p_{n-1.n} \end{matrix} \right\| \right\rangle,$$

where

$$\bar{\Phi}_i = \text{diag}\left( \frac{\varphi_i}{\varphi_1}, ..., \frac{\varphi_i}{\varphi_{i-1}}, \underbrace{1, ..., 1}_{n-i} \right).$$

That is

$$\left\langle \bar{\Phi}_i \left\| \begin{matrix} p_{1.r+1} & ... & p_{1n} \\ ... & ... & ... \\ p_{n-1.r+1} & ... & p_{n-1.n} \end{matrix} \right\| \right\rangle = \sigma_{i.r+1} ... \sigma_{in}.$$

By Lemma 7.2,

$$\bar{\Phi}_i \left\| \begin{matrix} p_{1.r+1} & ... & p_{1n} \\ ... & ... & ... \\ p_{n-1.r+1} & ... & p_{n-1.n} \end{matrix} \right\| \overset{l}{\sim} \text{triang}\,(\sigma_{i.r+1}, ..., \sigma_{in}). \tag{7.5}$$

It is obvious that g.c.d. of maximal order minors of the matrices

$$\Phi_i \left\| \begin{matrix} p_{1.r-1} & 0 & p_{1,r+1} & ... & p_{1n} \\ ... & ... & ... & ... & ... \\ p_{n-1.r-1} & 0 & p_{n-1.r+1} & ... & p_{n-1.n} \\ p_{n.r-1} & 1 & p_{n.r+1} & ... & p_{nn} \end{matrix} \right\|,$$

$$\bar{\Phi}_i \left\| \begin{matrix} p_{1.r-1} & p_{1.r+1} & ... & p_{1n} \\ ... & ... & ... & ... \\ p_{n-2.r-1} & p_{n-2.r+1} & ... & p_{n-2,n} \\ p_{n-1.r-1} & p_{n-1.r+1} & ... & p_{n-1.n} \end{matrix} \right\|$$

**237**

coincide and equal to $\sigma_{i.r-1}\sigma_{i.r+1}...\sigma_{in}$. In view of (7.5), we have

$$\bar{\Phi}_i \left\| \begin{matrix} p_{1.r-1} & p_{1.r+1} & ... & p_{1n} \\ ... & ... & ... & ... \\ p_{n-2.r-1} & p_{n-2.r+1} & ... & p_{n-2,n} \\ p_{n-1.r-1} & p_{n-1.r+1} & ... & p_{n-1.n} \end{matrix} \right\| \overset{l}{\sim} \text{triang}\,(\sigma_{i.r-1}, \sigma_{i.r+1}, ..., \sigma_{in}).$$

Thinking further, we have

$$\bar{\Phi}_i \left\| \begin{matrix} p_{11} & ... & p_{1.r-1} & p_{1.r+1} & ... & p_{1n} \\ ... & ... & ... & ... & ... & ... \\ p_{n-1.1} & ... & p_{n-1.r-1} & p_{n-1.r+1} & ... & p_{n-1.n} \end{matrix} \right\| \overset{l}{\sim}$$

$$\overset{l}{\sim} \text{triang}\,(\sigma_{i1}, ..., \sigma_{i.r-1}, \sigma_{i.r+1}, ..., \sigma_{in}),$$

$i = 1, ..., n - 1$. It means that equality (7.4) is fulfilled. $\qquad \square$

The following Theorem summarizes the obtained results.

**Theorem 7.3.** *The set* $\mathbf{T}^{-1}(\Phi)\Phi$ *consists of all right nonassociated matrices with the* $\Phi$-*skeleton* $F(\Phi)$.

**Proof.** By virtue of Theorem 7.2, the set $\mathbf{T}(\Phi)$ consists of representatives of the various left cosets of $\mathrm{GL}_n(R)$ at $\mathbf{G}_\Phi$. Then according to Theorem 4.5, matrices from $\mathbf{T}^{-1}(\Phi)\Phi$ are right unassociated, and have the $\Phi$-skeleton $F(\Phi)$.

Now let the matrix $B = P^{-1}\Phi Q^{-1}$ has the $\Phi$-skeleton $F(\Phi) = \|\sigma_{ij}\|_1^n$. Let us show that $B \overset{r}{\sim} T^{-1}\Phi$, where $T \in \mathbf{T}(\Phi)$.

Suppose that $n = 2$. Since $\mathbf{P}_B = \mathbf{G}(\Phi)P$ and $\sigma_{22} = 1$, by Lemma 7.1, the set $\mathbf{P}_B$ contains a matrix of the form

$$P' = \left\| \begin{matrix} e & 0 \\ a & 1 \end{matrix} \right\|,$$

where $e \in U(R)$. Write the element $a$ in the form

$$a = \alpha + \frac{\varphi_2}{\varphi_1}t,$$

where $\alpha \in K\left(\dfrac{\varphi_2}{\varphi_1}\right)$. Then

$$\underbrace{\left\| \begin{matrix} e^{-1} & 0 \\ -\dfrac{\varphi_2}{\varphi_1}e^{-1}t & 1 \end{matrix} \right\|}_{H} \left\| \begin{matrix} e & 0 \\ a & 1 \end{matrix} \right\| = \left\| \begin{matrix} 1 & 0 \\ \alpha & 1 \end{matrix} \right\| = T \in \mathbf{T}(\Phi).$$

This yields $HP' = T$, where $H \in \mathbf{G}_\Phi$. By Theorem 4.5, it means that $B \overset{r}{\sim} T^{-1}\Phi$.

Suppose that the Theorem holds for all $k < n$. As above, we choose in the set $\mathbf{P}_B$ a transforming matrix of the form

$$P = \left\| \begin{matrix} U_{n-1} & \mathbf{0} \\ * & 1 \end{matrix} \right\|.$$

By virtue of Lemma 7.4,

$$S_{\operatorname{diag}(\varphi_1, \, ..., \, \varphi_{n-1})}(U_{n-1}) = F(\operatorname{diag}(\varphi_1, \, ..., \, \varphi_{n-1})).$$

By induction, there is $G_{n-1} \in \mathbf{G}_{\operatorname{diag}(\varphi_1, \, ..., \, \varphi_{n-1})}$ such that

$$G_{n-1}U_{n-1} = T_{n-1} \in \mathbf{T}(\operatorname{diag}(\varphi_1, \, ..., \, \varphi_{n-1})).$$

Therefore,

$$(G_{n-1} \oplus I_1)P = \left\| \begin{matrix} & T_{n-1} & & \mathbf{0} \\ a_1 & ... & a_{n-1} & 1 \end{matrix} \right\| = P_n,$$

where

$$a_{n-1} = \alpha_{n.n-1} + \frac{\varphi_n}{\varphi_{n-1}}t_{n.n-1}, \;\; \alpha_{n.n-1} \in K\left(\frac{\varphi_n}{\varphi_{n-1}}\right).$$

It is obvious that $G_{n-1} \oplus I_1 \in \mathbf{G}_\Phi$. Then

$$\left\| \begin{matrix} & I_{n-1} & & \mathbf{0} \\ 0 & ... & 0 & -\frac{\varphi_n}{\varphi_{n-1}}t_{n.n-1} & 1 \end{matrix} \right\| P_n =$$

$$= H_{n-1}P_n = \left\| \begin{matrix} & T_{n-1} & & \mathbf{0} \\ b_1 & ... & b_{n-2} & \alpha_{n.n-1} & 1 \end{matrix} \right\| = P_{n-1},$$

where

$$b_{n-2} = \alpha_{n.n-2} + \frac{\varphi_n}{\varphi_{n-2}}t_{n.n-2}, \;\; \alpha_{n.n-2} \in K\left(\frac{\varphi_n}{\varphi_{n-2}}\right).$$

Again, we get that

$$\left\| \begin{matrix} & I_{n-1} & & & \mathbf{0} \\ 0 & ... & 0 & -\frac{\varphi_n}{\varphi_{n-2}}t_{n.n-2} & 0 & 1 \end{matrix} \right\| P_{n-1} =$$

$$= H_{n-2}P_{n-1} = \left\| \begin{matrix} & T_{n-1} & & \mathbf{0} \\ c_1 & ... & \alpha_{n.n-2} & \alpha_{n.n-1} \end{matrix} \right\|.$$

Continuing the described process, we obtain a matrix

$$H = H_1 ... H_{n-1}(G_{n-1} \oplus I_1) \in \mathbf{G}_\Phi$$

such that

$$HP = \left\| \begin{matrix} & T_{n-1} & & \mathbf{0} \\ \alpha_{n1} & ... & \alpha_{n.n-2} & \alpha_{n.n-1} \end{matrix} \right\| \in \mathbf{T}(\Phi).$$

The Theorem is proved. $\qquad\square$

**Corollary 7.1.** *The set* $\mathbf{P}_B$ *contains a lower unitriangular matrix if and only if* $S_\Phi(B) = F(\Phi)$. $\qquad\square$

Now we turn to the description of unassociated matrices with standard $\Phi$-skeletons.

**Lemma 7.5.** *Let $P \in \mathrm{GL}_n(R)$ and $\tau \in S_n$. If $S_\Phi(P) = F(\Phi)E(\tau)$, then*

$$S_\Phi(PE^{-1}(\tau)) = F(\Phi).$$

**Proof.** Let

$$\tau = \begin{pmatrix} 1 & ... & n \\ i_1 & ... & i_n \end{pmatrix}.$$

Then all elements of the $i_n$th column of the matrix $F(\Phi)E(\tau)$ is equal to 1. By removing the $i_n$th column and the last row of the matrix $P$, we get the matrix $\bar{p}_{i_n}$. By Lemma 7.1, we can assume that

$$\bar{p}_{i_n} = \|0 \ ... \ 0 \ 1\|^T.$$

Denote by $P_{ni_n}$ the matrix which cut down from $P$ by removing of its last row and $i_n$th column. According to Lemma 7.4,

$$S_{\mathrm{diag}\,(\varphi_1, ..., \varphi_{n-1})}(P_{ni_n}) = S_\Phi(P)_{ni_n}.$$

Let us leave unchanged the numbering of the columns of the $P_{ni_n}$ and $S_\Phi(P)_{ni_n}$ as it was in the matrices $P$ and $S_\Phi(B)$, respectively. Denote by $\bar{u}_{i_{n-1}}$ the $i_{n-1}$th column of $P_{ni_n}$. All elements of the $i_{n-1}$th column of $S_\Phi(P)_{ni_n}$ are 1. It follows that there is $H_{n-1} \in \mathbf{G}_{\mathrm{diag}\,(\varphi_1, ..., \varphi_{n-1})}$ such that

$$H_{n-1}\bar{u}_{i_{n-1}} = \|0 \ ... \ 0 \ 1\|^T.$$

Therefore, the matrix consisting of $i_{n-1}$th and $i_n$th columns of the matrix $(H_{n-1} \oplus I_1)P$ has the form

$$\left\| \begin{matrix} \mathbf{0} & \\ 1 & 0 \\ * & 1 \end{matrix} \right\|.$$

Since $H_{n-1} \oplus I_1 \in \mathbf{G}_\Phi$, then

$$S_\Phi((H_{n-1} \oplus I_1)P) = S_\Phi(P) = F(\Phi)E(\tau).$$

Similarly there is $H_{n-2} \in \mathbf{G}_{\mathrm{diag}\,(\varphi_1, ..., \varphi_{n-2})}$ such that

$$(H_{n-2} \oplus I_2)\bar{w}_{i_{n-2}} = \|0 \ ... \ 0 \ 1 \ * \ *\|^T,$$

where $\bar{w}_{i_{n-2}}$ is the $i_{n-2}$th column of $(H_{n-1} \oplus I_1)P$. Continuing this process, we find in the group $\mathbf{G}_\Phi$ matrices of the form $H_{n-i} \oplus I_i$, where $H_{n-i} \in \mathbf{G}_{\mathrm{diag}\,(\varphi_1, ..., \varphi_{n-i})}$, $i = 1, ..., n-2$, for which $i_s$th columns of the matrices

$$(H_2 \oplus I_{n-2}) ... (H_{n-1} \oplus I_1)P$$

have the form

$$\left\| \underbrace{0 \ ... \ 0 \ 1}_{i_s} \ * \ * \right\|^T = \bar{t}_{i_s},$$

$s = 1, ..., n$. Consequently, the matrix $T = \left\| \bar{t}_{i_1} \ ... \ \bar{t}_{i_n} \right\|$ is the lower unitriangular matrix. Therefore, $S_\Phi(T) = F(\Phi)$. Noting that

$$\tau^{-1} = \begin{pmatrix} i_1 & ... & i_n \\ 1 & ... & n \end{pmatrix},$$

and $E(\tau^{-1}) = E^{-1}(\tau)$, we get

$$T = \underbrace{(H_2 \oplus I_{n-2}) ... (H_{n-1} \oplus I_1)}_{H} PE^{-1}(\tau),$$

where $H \in \mathbf{G}_\Phi$. Hence,

$$F(\Phi) = S_\Phi(T) = S_\Phi(HPE^{-1}(\tau)) = S_\Phi(PE^{-1}(\tau)).$$

Lemma is proved. □

**Remark**. The permutation of columns of an invertible matrix does not follow to a similar permutation of the $\Phi$-skeleton of an obtained matrix.

**Example 7.1.** Let $R = \mathbb{Z}$,

$$\Phi = \left\| \begin{matrix} 1 & 0 \\ 0 & \varphi \end{matrix} \right\|, \quad P = \left\| \begin{matrix} 1 & 0 \\ 1 & 1 \end{matrix} \right\|^{-1},$$

where $\varphi \notin \{0, \ \pm 1\}$, and

$$\tau = \begin{pmatrix} 1 & 2 \\ 2 & 1 \end{pmatrix}, \quad E(\tau) = \left\| \begin{matrix} 0 & 1 \\ 1 & 0 \end{matrix} \right\|.$$

However, the $\Phi$-skeletons of the matrices $P$ and $PE(\tau)$ coincide:

$$S_\Phi(P) = S_\Phi(PE(\tau)) = \left\| \begin{matrix} 1 & 1 \\ \varphi & 1 \end{matrix} \right\|. \qquad \diamond$$

The units located at main diagonals of matrices from $\mathbf{T}(\Phi)$, are called **diagonal units**. In every matrix of this set the elements at the position to the left of diagonal units we replace for zeros. The resulting set denote by $\mathbf{T}_\tau(\Phi)$.

**Theorem 7.4.** *The set $\mathbf{T}_\tau^{-1}(\Phi)\Phi$ consists of all right unassociated matrices with the $\Phi$-skeleton $F(\Phi)E(\tau)$.*

**Proof.** Let $T \in \mathbf{T}_\tau(\Phi)$, where

$$\tau = \begin{pmatrix} 1 & ... & n \\ i_1 & ... & i_n \end{pmatrix}, \quad \tau^{-1} = \begin{pmatrix} 1 & ... & n \\ j_1 & ... & j_n \end{pmatrix}.$$

Then

$$\left\| 1 \ ... \ n \right\| E(\tau) = \left\| j_1 \ ... \ j_n \right\|.$$

Therefore, the $s$th column of the matrix $\mathrm{triang}\,(1, ..., 1)E(\tau)$ is the $j_s$th column of the matrix $\mathrm{triang}\,(1, ..., 1)$, $s = 1, ..., n$. Consequently, $\bar{t}_s$ which is the $s$th column of the matrix $T$ has the form

$$\bar{t}_s = \left\| \underbrace{0 \ ... \ 0 \ 1}_{j_s} *...* \right\|^T.$$

From the method of construction of $\mathbf{T}_\tau(\Phi)$ it follows that

$$\bar{t}_n = \left\| \underbrace{0 \; ... \; 0 \; 1}_{j_n}, 0...0 \right\|^T.$$

Therefore,

$$\Phi_n \bar{t}_n \overset{l}{\sim} \left\| 0 \; ... \; 0 \; \frac{\varphi_n}{\varphi_{j_n}} \right\|^T,$$

where

$$\Phi_n = \operatorname{diag}\left(\frac{\varphi_n}{\varphi_1}, ..., \frac{\varphi_n}{\varphi_{n-1}}, 1\right) = \operatorname{diag}\left(\frac{\varphi_n}{\varphi_1}, ..., \frac{\varphi_n}{\varphi_{n-1}}, \frac{\varphi_n}{\varphi_n}\right).$$

The matrix $\Phi_n \left\| \bar{t}_{n-1} \; \bar{t}_n \right\|$ contains only one nonzero minor: $\pm \dfrac{\varphi_n}{\varphi_{j_{n-1}}} \dfrac{\varphi_n}{\varphi_{j_n}}$. Consequently,

$$\Phi_n \left\| \bar{t}_{n-1} \; \bar{t}_n \right\| \overset{l}{\sim} \operatorname{triang}\left(\frac{\varphi_n}{\varphi_{j_{n-1}}}, \frac{\varphi_n}{\varphi_{j_n}}\right).$$

Continuing similar considerations, we obtain that

$$\Phi_n T \overset{l}{\sim} \operatorname{triang}\left(\frac{\varphi_n}{\varphi_{j_1}}, ..., \frac{\varphi_n}{\varphi_{j_n}}\right).$$

That is, the last rows of the matrices $S_\Phi(T)$ and $F(\Phi)E(\tau)$ are coincide and equal to

$$\left\| \frac{\varphi_n}{\varphi_{j_1}} \; ... \; \frac{\varphi_n}{\varphi_{j_n}} \right\|.$$

Obviously, the rest of the rows of the matrices $F(\Phi)E(\tau)$ and $S_\Phi(T)$ are coincide. So $S_\Phi(T) = F(\Phi)E(\tau)$.

On the contrary, let the matrix $B = P^{-1}\Phi Q^{-1}$ has the $\Phi$-skeleton $F(\Phi)E(\tau)$. That is

$$S_\Phi(B) = S_\Phi(P) = F(\Phi)E(\tau).$$

According to Lemma 7.5,

$$S_\Phi(PE^{-1}(\tau)) = F(\Phi).$$

From Theorem 7.3 it follows that there is $H_n \in \mathbf{G}_\Phi$ such that

$$H_n PE^{-1}(\tau) = \operatorname{triang}(1, ..., 1).$$

Hence,

$$H_n P = \operatorname{triang}(1, ..., 1)E(\tau).$$

By Theorem 6.7,

$$F(\Phi)E(\tau) = S_\Phi(P) = S_\Phi(H_n P) = S_\Phi(\operatorname{triang}(1, ..., 1)E(\tau)). \tag{7.6}$$

The columns of the matrix $H_n P$ we number as following:

$$\bar{a}_s = \left\| \underbrace{0\ ...\ 0\ 1}_{s} * ... * a_s \right\|^T$$

is the $s$th column, $s = 1, ..., n$. Then

$$H_n P = \left\| \bar{a}_{j_1} \quad ... \quad \bar{a}_{j_n} \right\|.$$

From equality (7.6) it follows that the last row of the matrix $S_\Phi(H_n P)$ is

$$\left\| \frac{\varphi_n}{\varphi_{j_1}} \quad ... \quad \frac{\varphi_n}{\varphi_{j_n}} \right\|.$$

Let $j_l = n$. Then the matrix $H_n P$ has the form

$$H_n P = \left\| \begin{matrix} C & & \mathbf{0} & & B & \\ a_1 & ... & a_{l-1} & 1 & a_{l+1} & ... & a_n \end{matrix} \right\|.$$

According to the definition of $\Phi$-skeleton we get

$$\Phi_n \left\| \bar{a}_{j_l} \quad \bar{a}_{j_{l+1}} \quad ... \quad \bar{a}_{j_n} \right\| \overset{l}{\sim} \text{triang}\left(1, \frac{\varphi_n}{\varphi_{j_{l+1}}}, ..., \frac{\varphi_n}{\varphi_{j_n}}\right).$$

Taking into account the structure of the matrix $H_n P$, we have

$$\text{diag}\left(\frac{\varphi_n}{\varphi_1}, ..., \frac{\varphi_n}{\varphi_{n-1}}\right) B \overset{l}{\sim} \text{triang}\left(\frac{\varphi_n}{\varphi_{j_{l+1}}}, ..., \frac{\varphi_n}{\varphi_{j_n}}\right).$$

Let $j_k = n-1$. If $a_k$ is on the right of a diagonal unit, then by Lemma 7.3, there is $H_{n-1} \in \mathbf{G}_\Phi$ such that in the matrix $H_{n-1} H_n P$ the element $a_k$ is replaced by zero. If element $a_k$ is on the left of a diagonal unit, then using the methods that where used to prove Theorem 7.3, we construct the following matrix $H_{n-1}$ such that in the matrix $H_{n-1} H_n P$ the element $a_k$ is replaced by $\alpha \in K\left(\dfrac{\varphi_n}{\varphi_{n-1}}\right)$. Reasoning similarly, we find $L_n \in \mathbf{G}_\Phi$ such that in the matrix $L_n P$ on the right of a diagonal unit there will stand zeros and on the left are elements from $K\left(\dfrac{\varphi_n}{\varphi_i}\right)$, $1 \leqslant i \leqslant n-1$. We remove the last row and the first column in the matrix $L_n P$. Denote the resulting matrix by $(L_n P)_{nl}$. Using Lemma 7.4 and thinking as above, we will find $L_{n-1} \in \mathbf{G}_{diag\,(\varphi_1, ..., \varphi_{n-1})}$ such that in the last row of the matrix $L_{n-1}(L_n P)_{nl}$ on the right of a diagonal unit is zero, and on the left are elements from $K\left(\dfrac{\varphi_{n-1}}{\varphi_i}\right)$, $1 \leqslant i \leqslant n-2$. So the elements of the last two rows of the matrix $(L_{n-1} \oplus I_1) L_n P$, where $(L_{n-1} \oplus I_1) L_n \in \mathbf{G}_\Phi$ satisfy the conditions of the Theorem. Continuing the described process, we construct the matrix $H \in \mathbf{G}_\Phi$, such that $HP = T \in \mathbf{T}_\tau(\Phi)$. That is $B \overset{r}{\sim} T^{-1}\Phi$. The theorem is proved. $\qquad\square$

Denote by $S_n^\Phi$ permutations of order $n$ by which we get different matrices of the form $F(\Phi)E(\tau)$, and denote

$$\text{Stand}\,(\Phi) = \underset{\tau \in S_n^\Phi}{\cup}\, T_\tau(\Phi).$$

**Theorem 7.5.** *The set* $\text{Stand}\,(\Phi)^{-1}\Phi$ *consists of all right unassociated matrices with standard $\Phi$-skeletons.* $\qquad\square$

**Example 7.2.** Let $R = \mathbb{Z}$, $\Phi = \text{diag}\,(1,3,6)$. Since $S_n^\Phi = S_n$, then

$$\text{Stand}(\Phi) = \left\{ \begin{Vmatrix} 1 & 0 & 0 \\ a & 1 & 0 \\ b & c & 1 \end{Vmatrix}, \begin{Vmatrix} 1 & 0 & 0 \\ a & 0 & 1 \\ b & 1 & 0 \end{Vmatrix}, \begin{Vmatrix} 0 & 1 & 0 \\ 1 & 0 & 0 \\ c & b & 1 \end{Vmatrix} \right\} \bigcup$$

$$\bigcup \left\{ \begin{Vmatrix} 0 & 1 & 0 \\ 0 & a & 1 \\ 1 & 0 & 0 \end{Vmatrix}, \begin{Vmatrix} 0 & 0 & 1 \\ 1 & 0 & 0 \\ c & 1 & 0 \end{Vmatrix}, \begin{Vmatrix} 0 & 0 & 1 \\ 0 & 1 & 0 \\ 1 & 0 & 0 \end{Vmatrix} \right\},$$

where $a \in \{0,1,2\,\}$, $b \in \{0,1,...,5\}$, $c \in \{0,1\}$, and $\text{Stand}\,(\Phi)^{-1}\Phi$ is the set of all right unassociated matrices with standard $\Phi$-skeletons.

Note that the set $\text{Stand}\,(\Phi)$ does not exhaust all representatives of the left cosets of the group $\text{GL}_n(R)$ at $\mathbf{G}_\Phi$. Indeed, this set does not contain representatives of cosets with $\Phi$-skeletons:

$$\begin{Vmatrix} 1 & 1 & 1 \\ 1 & 1 & 3 \\ 2 & 2 & 3 \end{Vmatrix}, \begin{Vmatrix} 1 & 1 & 1 \\ 1 & 3 & 1 \\ 2 & 3 & 2 \end{Vmatrix}, \begin{Vmatrix} 1 & 1 & 1 \\ 3 & 1 & 1 \\ 3 & 2 & 2 \end{Vmatrix}.$$

Examples of matrices with such $\Phi$-skeletons are, in particular,

$$\left\{ \begin{Vmatrix} 0 & 0 & 1 \\ 0 & 1 & 0 \\ 1 & 0 & 3 \end{Vmatrix}, \begin{Vmatrix} 0 & 1 & 0 \\ 0 & 0 & 1 \\ 1 & 3 & 0 \end{Vmatrix}, \begin{Vmatrix} 0 & 1 & 0 \\ 0 & 0 & 1 \\ 1 & 2 & 0 \end{Vmatrix} \right\}^{-1} \Phi. \qquad \Diamond$$

Denote by $W(\Phi)$ all representatives of the left cosets of the group $\text{GL}_n(R)$ at $\mathbf{G}_\Phi$.

**Example 7.3.** Let $R = \mathbb{Q}[x]$, $\Phi = \text{diag}\,(1, 1, x^2 + 1)$. Then

$$F(\Phi) = \begin{Vmatrix} 1 & 1 & 1 \\ 1 & 1 & 1 \\ x^2 + 1 & x^2 + 1 & 1 \end{Vmatrix}.$$

Hence,

$$S_n^\Phi = \left\{ \begin{pmatrix} 1 & 2 & 3 \\ 1 & 2 & 3 \end{pmatrix}, \begin{pmatrix} 1 & 2 & 3 \\ 1 & 3 & 2 \end{pmatrix}, \begin{pmatrix} 1 & 2 & 3 \\ 3 & 2 & 1 \end{pmatrix} \right\}.$$

That is, in addition to $F(\Phi)$, there are two standard $\Phi$-skeletons:

$$\left\| \begin{array}{cccc} 1 & 1 & 1 & \\ 1 & 1 & 1 & v \\ x^2+1 & 1 & x^2+1 & \end{array} \right\|, \quad \left\| \begin{array}{ccc} 1 & 1 & 1 \\ 1 & 1 & 1 \\ 1 & x^2+1 & x^2+1 \end{array} \right\|.$$

The polynomial $x^2 + 1$ is irreducible in the ring $\mathbb{Q}[x]$. Consequently, there are no other $\Phi$-skeletons of a matrix with the Smith form $\Phi$. Therefore,

$$W(\Phi) = \text{Stand}\,(\Phi) =$$
$$= \left\{ \left\| \begin{array}{ccc} 1 & 0 & 0 \\ 0 & 1 & 0 \\ b(x) & c(x) & 1 \end{array} \right\|, \left\| \begin{array}{ccc} 1 & 0 & 0 \\ 0 & 0 & 1 \\ b(x) & 1 & 0 \end{array} \right\|, \left\| \begin{array}{ccc} 0 & 0 & 1 \\ 0 & 1 & 0 \\ 1 & 0 & 0 \end{array} \right\| \right\},$$

where $b(x), c(x) \in \{mx + n | m, n \in \mathbb{Q}\}$. Hence, $\text{Stand}\,(\Phi)^{-1}\Phi$ is the set of all right unassociated matrices with the Smith form $\Phi$. $\diamond$

The regularity established in this example can be transferred to a wider class of matrices.

**Corollary 7.2.** *The set* $\text{Stand}\,(\Phi)^{-1}\Phi$, *where*

$$\Phi = \text{diag}\,(1, ..., 1, \underbrace{\varphi, ..., \varphi}_{i}),$$

*and $\varphi$ is irreducible elements of $R$, $1 \leqslant i < n$, is the set of all right unassociated matrices with the Smith form $\Phi$.* $\square$

## 7.2. Normal form of disappear matrices with respect to one-sided transformations

Let $R$ be a Bezout ring of stable rank 1.5 and $P = \|p_{ij}\|_1^n$ be an invertible matrix over $R$. A matrix is called **disappear** if it has a Smith form:

$$\Phi = \text{diag}\,(1, ..., 1, \varphi), \quad \varphi \neq 0, \quad \varphi \notin U(R).$$

In this case the group $\mathbf{G}_\Phi$ consists of all invertible matrices of the form

$$\left\| \begin{array}{cccc} h_{11} & ... & h_{1.n-1} & h_{1n} \\ ... & ... & ... & ... \\ h_{n-1.1} & ... & h_{n-1.n-1} & h_{n-1.n} \\ \varphi h_{n1} & ... & \varphi h_{n.n-1} & h_{nn} \end{array} \right\|.$$

Consequently,

$$\Phi_1 = ... = \Phi_{n-1} = I, \quad \Phi_n = \text{diag}(\varphi, ..., \varphi, 1).$$

There is an invertible matrix $V$ such that

$$V\Phi_n P = \left\| \begin{array}{ccccc} \gamma_1 & 0 & & 0 & 0 \\ c_{21} & \gamma_2 & & 0 & 0 \\ \vdots & & \ddots & & \\ c_{n-1.1} & c_{n-1.2} & & \gamma_{n-1} & 0 \\ c_{n1} & c_{n2} & ... & c_{n.n-1} & \gamma_n \end{array} \right\| \tag{7.7}$$

is the Hermite form of the matrix $\Phi_n P$. So the $\Phi$-skeleton of $P$ is

$$V\Phi_n P = \begin{Vmatrix} 1 & 1 & ... & 1 \\ ... & ... & ... & ... \\ 1 & 1 & ... & 1 \\ \gamma_1 & \gamma_2 & ... & \gamma_n \end{Vmatrix} = S_\Phi(P).$$

**Lemma 7.6.** *The equalities*
$$\gamma_i = \frac{\varphi \nu_i}{\nu_{i+1}},$$
*where*
$$\nu_i = (p_{ni}, p_{n.i+1}, ..., p_{n.n-1}, p_{nn}, \varphi), \quad i = 1, ..., n, \quad \nu_{n+1} = \varphi,$$
*hold.*

**Proof.** Consider the matrix

$$P_i = \begin{Vmatrix} p_{1i} & p_{1.i+1} & ... & p_{1.n-1} & p_{1n} \\ ... & ... & ... & ... & ... \\ p_{n-1.i} & p_{n-1.i+1} & ... & p_{n-1.n-1} & p_{n-1.n} \\ p_{ni} & p_{n.i+1} & ... & p_{n.n-1} & \gamma_n \end{Vmatrix}.$$

By Proposition 3.8, the matrix $\Phi_n P_i$, $i = 1, ..., n-1$, has the Smith form

$$\Phi_n P_i \sim \begin{Vmatrix} \nu_i & 0 & & 0 \\ 0 & \varphi & & 0 \\ & & \ddots & \\ 0 & 0 & & \varphi \\ \mathbf{0} & \mathbf{0} & \mathbf{0} & \mathbf{0} \end{Vmatrix} = \begin{Vmatrix} \nu_i \oplus \varphi I_{n-i} \\ \mathbf{0} \end{Vmatrix},$$

where $\nu_i = (p_{ni}, p_{n.i+1}, ..., p_{nn}, \varphi)$. Hence,

$$\langle \Phi_n P_i \rangle = \det(\nu_i \oplus \varphi I_{n-i}) = \nu_i \varphi^{n-i}.$$

On the other hand, from equality (7.7) it follows that

$$\Phi_n P_i \overset{l}{\sim} \begin{Vmatrix} \mathbf{0} & \mathbf{0} & ... & \mathbf{0} & \mathbf{0} \\ \gamma_i & 0 & ... & 0 & 0 \\ c_{i+1.i-1} & \gamma_{i+1} & & 0 & 0 \\ & & \ddots & & \\ c_{n-1.i} & c_{n-1.i+1} & & \gamma_{n-1} & 0 \\ c_{ni} & c_{n.i+1} & ... & c_{n.n-1} & \gamma_n \end{Vmatrix}.$$

That is

$$\langle \Phi_n P_i \rangle = \det \begin{Vmatrix} \gamma_i & 0 & ... & 0 & 0 \\ c_{i+1.i} & \gamma_{i+1} & & 0 & 0 \\ & & \ddots & & \\ c_{n-1.i} & c_{n-1.i+1} & & \gamma_{n-1} & 0 \\ c_{ni} & c_{n.i+1} & ... & c_{n.n-1} & \gamma_n \end{Vmatrix} = \gamma_i \gamma_{i+1} ... \gamma_n.$$

So $\gamma_i\gamma_{i+1}\dots\gamma_n = \nu_i\varphi^{n-i}$. Similarly, we show that

$$\gamma_{i+1}\gamma_{i+2}\dots\gamma_n = \nu_{i+1}\varphi^{n-i-1}.$$

Then

$$\frac{\gamma_i\gamma_{i+1}\dots\gamma_n}{\gamma_{i+1}\gamma_{i+2}\dots\gamma_n} = \gamma_i = \frac{\nu_i\varphi^{n-i}}{\nu_{i+1}\varphi^{n-i-1}} = \frac{\varphi\nu_i}{\nu_{i+1}}.$$

Lemma is proved. □

In view of (7.7), we have

$$\gamma_n = (\varphi p_{1n},\dots,\varphi p_{n-1.n},p_{nn}) = (\varphi(p_{1n},\dots,p_{n-1.n}),p_{nn}).$$

Since $(p_{1n},\dots,p_{n-1.n},p_{nn}) = 1$, then $\gamma_n = (\varphi,p_{nn})$. It means that $\|1 \ \dots \ 1 \ \gamma_n\|^T$ is $\Phi$-rod of the column $\|p_{1n} \ p_{2n} \ \dots \ p_{nn}\|^T$. According to Theorem 6.3, there is $H \in \mathbf{G}_\Phi$ such that

$$H\|p_{1n} \ p_{2n} \ \dots \ p_{nn}\|^T) = \|p'_{1n} \ \dots \ p'_{n-1.n} \ \gamma_n\|^T.$$

Consider the matrix $HP$. By Theorem 6.7, $S_\Phi(P) = S_\Phi(HP)$. Therefore, it can be assumed that the element $p_{nn}$ of the matrix $P$ is equal to $\gamma_n$.

The following statement shows how the elements of the matrix $P$ are changed by action from $\mathbf{G}_\Phi$.

**Theorem 7.6.** *Let* $P = \|p_{ij}\|_1^n$ *be an invertible matrix over $R$ moreover* $p_{nn} = \gamma_n|\varphi$, *and*

$$S_\Phi(P) = \begin{Vmatrix} 1 & 1 & \dots & 1 \\ \dots & \dots & \dots & \dots \\ 1 & 1 & \dots & 1 \\ \gamma_1 & \gamma_2 & \dots & \gamma_n \end{Vmatrix}.$$

*If $p_{ni} \neq 0$, then the matrix equation*

$$x\Phi_n P_i = \|a_i \ \ p_{n.i+1} \ \ \dots \ \ p_{n.n-1} \ \ \gamma_n\|, \tag{7.8}$$

*where*

$$P_i = \begin{Vmatrix} p_{1i} & p_{1.i+1} & \dots & p_{1.n-1} & p_{1n} \\ \dots & \dots & \dots & \dots & \dots \\ p_{n-1.i} & p_{n-1.i+1} & \dots & p_{n-1.n-1} & p_{n-1.n} \\ p_{ni} & p_{n.i+1} & \dots & p_{n.n-1} & \gamma_n \end{Vmatrix},$$

$1 \leqslant i \leqslant n - 1$, *is solvable and the group $\mathbf{G}_\Phi$ contains a matrix of the form*

$$\begin{Vmatrix} & & * & \\ \varphi\overline{x}_1 & \dots & \varphi\overline{x}_{n-1} & \overline{x}_n \end{Vmatrix}, \tag{7.9}$$

*where $\|\overline{x}_1 \ \dots \ \overline{x}_{n-1} \ \overline{x}_n\|$ is the solution of the equation (7.8), if and only if conditions hold*

*1) $a_i \equiv p_{ni}(\mathrm{mod}\gamma_i)$,*

*2) $(a_i,p_{n.i+1},\dots,p_{n.n-1},\gamma_n) = (p_{ni},p_{n.i+1},\dots,p_{n.n-1},\gamma_n)$.*

*Or conditions equivalent to them:*

*3)* $(a_i, \varphi) = (p_{ni}, \varphi)$,

*4)* $a'_i \equiv p'_{ni} \left( \mod \dfrac{\varphi}{[(p_{ni}, \varphi), (p_{n.i+1}, ..., p_{n.n-1}, \gamma_{nn})]} \right)$, *where $a'_i$ and $p'_{ni}$ is a quotient of division $a_i, p_{ni}$ by $(p_{ni}, \varphi)$,*

$$[(p_{ni}, \varphi), (p_{n.i+1}, ..., p_{n.n-1}, \gamma_{nn})] = \frac{(p_{ni}, \varphi)(p_{n.i+1}, ..., p_{n.n-1}, \gamma_{nn})}{((p_{ni}, \varphi), (p_{n.i+1}, ..., p_{n.n-1}, \gamma_{nn}))}.$$

**Proof. Necessity**. According to Theorem 7.1, condition 1) is necessary and sufficient to solvability of equation (7.8).

Suppose that the group $\mathbf{G}_\Phi$ contains a matrix $H_i$ of the form (7.9). Then

$$H_i P = \left\| \begin{matrix} & & & & * & & & \\ b_{n1} & ... & b_{n.i-1} & a_{ni} & p_{n.i+1} & ... & p_{n.n-1} & \gamma_n \end{matrix} \right\| 1.$$

By Theorem 6.1, $(a_i, \varphi) = (p_{ni}, \varphi)$. Since $\gamma_n | \varphi$, then

$$(a_i, p_{n.i+1}, ..., p_{n.n-1}, \gamma_n) = (a_i, p_{n.i+1}, ..., p_{n.n-1}, (\gamma_n, \varphi)) =$$

$$= ((a_i, \varphi), p_{n.i+1}, ..., p_{n.n-1}, \gamma_n) = ((p_{ni}, \varphi), p_{n.i+1}, ..., p_{n.n-1}, \gamma_n) =$$

$$= (p_{ni}, p_{n.i+1}, ..., p_{n.n-1}, \gamma_n).$$

**Sufficiency**. Consider the last row of the matrix $P$:

$$\left\| \begin{matrix} p_{n1} & p_{n2} & ... & p_{n.n-1} & \gamma_n \end{matrix} \right\|.$$

Denote by

$$(p_{nj}, p_{n.j+1}, ..., p_{n.n-1}, \gamma_n) = \nu_j, \ \ j = 1, ..., n-1.$$

Let $p_{n.k_1}, p_{n.k_2}, ..., p_{n.k_p}$ be nonzero elements of the last row of the matrix $P$, $1 \le k_1 < k_2 < ... < k_p \le n-1$. Let $i = k_s$, $k_1 \le k_s \le k_p$, $i \ne 1$. From conditions 1) and 2) it follows that $a_i = p_{ni} + \gamma_i r$ и $(a_i, \nu_{i+1}) = \nu_i$. Hence, $(p_{ni} + \gamma_i r, \nu_{i+1}) = \nu_i$. That is

$$\left( \frac{p_{ni}}{\nu_i} + \frac{\gamma_i}{\nu_i} r, \frac{\nu_{i+1}}{\nu_i} \right) = 1.$$

By Lemma 7.6, $\gamma_i = \dfrac{\varphi \nu_i}{\nu_{i+1}}$. Therefore,

$$\left( \frac{p_{ni}}{\nu_i} + \frac{\varphi}{\nu_{i+1}} r, \frac{\nu_{i+1}}{\nu_i} \right) = 1. \tag{7.10}$$

It follows from the definition of the $\Phi$-skeleton of the matrix $P$ that there exists an invertible matrix $V$ such that

$$V\Phi_n P = \begin{Vmatrix} \gamma_1 & 0 & & 0 & 0 & 0 & & 0 & 0 \\ c_{21} & \gamma_2 & & 0 & 0 & 0 & & 0 & 0 \\ \cdots & & \ddots & & & & & & \\ c_{i-1.1} & c_{i-1.2} & & \gamma_{i-1} & 0 & 0 & & 0 & 0 \\ c_{i1} & c_{i2} & \cdots & c_{i.i-1} & \gamma_i & 0 & & 0 & 0 \\ c_{i+1.1} & c_{i+1.2} & \cdots & c_{i+1.i-1} & c_{i+1.i} & \gamma_{i+1} & & 0 & 0 \\ \cdots & \cdots & \cdots & \cdots & \cdots & & \ddots & & \\ c_{n-1.1} & c_{n-1.2} & \cdots & c_{n-1.i-1} & c_{n-1.i} & c_{n-1.i+1} & & \gamma_{n-1} & 0 \\ p_{n1} & p_{n2} & \cdots & p_{n.i-1} & p_{ni} & p_{n.i+1} & \cdots & p_{n.n-1} & \gamma_n \end{Vmatrix}.$$

The matrix $\Phi_n P$ has the Smith form $\operatorname{diag}(1, \varphi, ..., \varphi)$. The product of the first two invariant factors of the matrix $\Phi_n P$ is equal to g.c.d. of its second-order minors. So $\varphi$ is the divisor of all second-order minors of the matrix $\Phi_n P$. The matrix $V\Phi_n P$ is equivalent to the matrix $\Phi_n P$. Thus, all second-order minors of the matrix $V\Phi_n P$ are divisible by $\varphi$. In particular,

$$\varphi \ \Big| \ \begin{vmatrix} c_{ij} & \gamma_i \\ p_{nj} & p_{ni} \end{vmatrix}, \ \ j = k_1, ..., k_{s-1},$$

which is equivalent to

$$\frac{\varphi}{\nu_i} \ \Big| \ \begin{vmatrix} c_{ij} & \dfrac{\varphi}{\nu_{i+1}} \\ p_{nj} & \dfrac{p_{ni}}{\nu_i} \end{vmatrix}.$$

Consequently,

$$\frac{\nu_{i+1}}{\nu_i} \ \Big| \ \begin{vmatrix} c_{ij} & \dfrac{\varphi}{\nu_{i+1}} \\ p_{nj} & \dfrac{p_{ni}}{\nu_i} \end{vmatrix}.$$

Taking into account (7.10), and Proposition 1.9, we get

$$\left(p_{ni} + rc_{ij}, \frac{\nu_{i+1}}{\nu_i}\right) = \left(p_{ni}, c_{ij}, \frac{\nu_{i+1}}{\nu_i}\right), \ \ j = k_1, ..., k_{s-1}. \tag{7.11}$$

Consider the matrix

$$\left(I_{i-1} \oplus \begin{Vmatrix} 1 & 0 & & 0 & 0 \\ 0 & 1 & & 0 & 0 \\ & & \ddots & & \\ 0 & 0 & & 1 & 0 \\ r & 0 & & 0 & 1 \end{Vmatrix}\right) V\Phi_n P =$$

$$= \begin{Vmatrix} \gamma_1 & & 0 & 0 & 0 & & 0 & 0 \\ & \ddots & & & & & & \\ c_{i-1.1} & & \gamma_{i-1} & 0 & 0 & & 0 & 0 \\ c_{i1} & \dots & c_{i.i-1} & \gamma_i & 0 & & 0 & 0 \\ c_{i+1.1} & \dots & c_{i+1.i-1} & c_{i+1.i} & \gamma_{i+1} & & 0 & 0 \\ & & & & & \ddots & & \\ c_{n-1.1} & \dots & c_{n-1.i-1} & c_{n-1.i} & c_{n-1.i+1} & & \gamma_{n-1} & 0 \\ p_{n1}+rc_{i1} & \dots & p_{n.i-1}+rc_{i.i-1} & p_{ni}+r\gamma_i & p_{n.i+1} & \dots & p_{n.n-1} & \gamma_n \end{Vmatrix} =$$

$$= V_{i-1}\Phi_n P.$$

Set

$$(p_{n.k_1}+rc_{i.k_1}, ..., p_{n.k_{s-1}}+rc_{i.k_{s-1}}, p_{ni}+r\gamma_i, p_{n.i+1}, ..., p_{n.n-1}, \gamma_n) = \delta_i.$$

Under condition 2)

$$(p_{ni}+r\gamma_i, p_{n.i+1}, ..., p_{n.n-1}, \gamma_n) = (p_{ni}, p_{n.i+1}, ..., p_{n.n-1}, \gamma_n).$$

In addition, using equality (7.11), we have

$$\left(\delta_i, \frac{\nu_{i+1}}{\nu_i}\right) = \left(\left(p_{n.k_1}+rc_{i.k_1}, \frac{\nu_{i+1}}{\nu_i}\right), ..., \left(p_{n.k_{s-1}}+rc_{i.k_{s-1}}, \frac{\nu_{i+1}}{\nu_i}\right),\right.$$

$$\left.(p_{ni}+r\gamma_i, p_{n.i+1}, ..., p_{n.n-1}, \gamma_n)\right) =$$

$$= \left(p_{n.k_1}, c_{i.k_1}, ..., p_{n.k_{s-1}}, c_{i.k_{s-1}}, p_{ni}, p_{n.i+1}, ..., p_{n.n-1}, \gamma_n, \frac{\nu_{i+1}}{\nu_i}\right) =$$

$$= \left((p_{n.k_1}, ..., p_{n.k_{s-1}}, p_{ni}, p_{n.i+1}, ..., p_{n.n-1}, \gamma_n), (c_{i.k_1}, ..., c_{i.k_{s-1}}), \frac{\nu_{i+1}}{\nu_i}\right).$$

Since $\begin{Vmatrix} p_{n1} & p_{n2} & \dots & p_{n.n-1} & \gamma_n \end{Vmatrix}$ is the last row of the invertible matrix $P$, then

$$1 = (p_{n1}, ..., p_{n.i-1}, p_{ni}, p_{n.i+1}, ..., p_{n.n-1}, \gamma_n) =$$

$$= (p_{n.k_1}, ..., p_{n.k_{s-1}}, p_{ni}, p_{n.i+1}, ..., p_{n.n-1}, \gamma_n).$$

Hence,

$$\left(\delta_i, \frac{\nu_{i+1}}{\nu_i}\right) = 1.$$

By definition

$$(p_{n.i+1}, p_{n.i+2}, ..., p_{n.n-1}, \gamma_n) = \nu_{i+1}.$$

There is an invertible matrix $Q$ such that

$$\begin{Vmatrix} p_{n.i+1} & p_{n.i+2} & \dots & p_{n.n-1} & \gamma_n \end{Vmatrix} Q = \begin{Vmatrix} \nu_{i+1} & 0 & \dots & 0 \end{Vmatrix}.$$

Then

$$V\Phi_n P \left\| \begin{matrix} I_i & \mathbf{0} \\ \mathbf{0} & Q \end{matrix} \right\| = \left\| \begin{matrix} \gamma_1 & 0 & & 0 & 0 & 0 & & 0 & 0 \\ c_{21} & \gamma_2 & & 0 & 0 & 0 & & 0 & 0 \\ \dots & & \ddots & & & & & & \\ c_{i-1.1} & c_{i-1.2} & & \gamma_{i-1} & 0 & 0 & & 0 & 0 \\ c_{i1} & c_{i2} & \dots & c_{i.i-1} & \gamma_i & 0 & & 0 & 0 \\ c_{i+1.1} & c_{i+1.2} & \dots & c_{i+1.i-1} & c_{i+1.i} & * & * & * & * \\ \dots & \dots & \dots & \dots & \dots & \dots & \dots & \dots & \dots \\ c_{n-1.1} & c_{n-1.2} & \dots & c_{n-1.i-1} & c_{n-1.i} & * & * & * & * \\ p_{n1} & p_{n2} & \dots & p_{n.i-1} & p_{ni} & \nu_{i+1} & 0 & \dots & 0 \end{matrix} \right\|.$$

This matrix has the Smith form $\operatorname{diag}(1, \varphi, ..., \varphi)$. Since

$$\left\| \begin{matrix} c_{ij} & 0 \\ u_{nj} & \nu_{i+1} \end{matrix} \right\|$$

are submatrices of a resulting matrix, $j = 1, ..., i-1$, then

$$\varphi \left| \begin{matrix} c_{ij} & 0 \\ u_{nj} & \nu_{i+1} \end{matrix} \right| = c_{ij}\nu_{i+1}.$$

Therefore, $\dfrac{\varphi}{\nu_{i+1}} \Big| c_{ij}$. That is $c_{ij} = \dfrac{\varphi}{\nu_{i+1}} c'_{ij}$, $j = 1, ..., i-1$. Considering that $\gamma_i = \dfrac{\varphi}{\nu_{i+1}}\nu_i$, we get

$$\left(\delta_i, \frac{\varphi}{\nu_{i+1}}\right) = \left(p_{n1} + r\frac{\varphi}{\nu_{i+1}}c'_{i1}, ..., p_{n.i-1} + r\frac{\varphi}{\nu_{i+1}}c'_{i.i-1},\right.$$

$$\left. p_{ni} + r\frac{\varphi}{\nu_{i+1}}\nu_i, p_{n.i+1}, ..., p_{n.n-1}, \gamma_n, \frac{\varphi}{\nu_{i+1}}\right) =$$

$$= \left((p_{n1}, ..., p_{n.i-1}, p_{ni}, p_{n.i+1}, ..., p_{n.n-1}, \gamma_n), \frac{\varphi}{\nu_{i+1}}\right) = 1.$$

Since

$$\left(\delta_i, \frac{\nu_{i+1}}{\nu_i}\right) = 1,$$

then

$$\left(\delta_i, \frac{\varphi}{\nu_{i+1}}\frac{\nu_{i+1}}{\nu_i}\right) = \left(\delta_i, \frac{\varphi}{\nu_i}\right) = 1.$$

Write the matrix $V\Phi_n P$ in a block form:

$$V\Phi_n P = \left\| \begin{matrix} \gamma_1 & & 0 & 0 & & 0 \\ & \ddots & & 0 & & \\ c_{i-1.1} & \dots & \gamma_{i-1} & 0 & & 0 \\ c_{i1} & \dots & c_{i.i-1} & \gamma_i & & 0 \\ \dots & \dots & \dots & & \ddots & \\ p_{n1} & \dots & p_{n.i-1} & p_{ni} & \dots & \gamma_n \end{matrix} \right\| = \left\| \begin{matrix} B_{11} & \mathbf{0} \\ B_{21} & B_{22} \end{matrix} \right\|.$$

Since

$$\left\| \begin{matrix} \mathbf{0} \\ B_{22} \end{matrix} \right\| \overset{l}{\sim} \Phi_n P_i,$$

then the Smith forms of matrices

$$\left\| \begin{matrix} \mathbf{0} \\ B_{22} \end{matrix} \right\|, \quad \Phi_n P_i$$

are coincide. According to Proposition 3.8,

$$\Phi_n P_i \sim \left\| \begin{matrix} \nu_i \oplus \varphi I_{n-i} \\ \mathbf{0} \end{matrix} \right\|.$$

Hence,

$$B_{22} \sim \mathrm{diag}(\nu_i, \underbrace{\varphi, ..., \varphi}_{n-i}).$$

By virtue of Theorem 3.2,

$$B_{11} \sim \mathrm{diag}\left( \frac{\varphi}{\nu_i}, \varphi, ..., \varphi \right).$$

It means that

$$\langle B_{11} \rangle_1 = \frac{\varphi}{\nu_i}.$$

Since

$$\left( \delta_i, \frac{\varphi}{\nu_i} \right) = \left( (p_{n1} + rc_{i1}, ..., p_{ni} + r\gamma_i, p_{n.i+1}, ..., p_{n.n-1}, \gamma_n), \frac{\varphi}{\nu_i} \right) = 1,$$

we conclude that g.c.d. of elements of the matrix

$$\left\| \begin{matrix} \gamma_1 & & 0 & 0 & 0 & ... & 0 & 0 \\ & \ddots & & & & & & \\ c_{i-1.1} & & \gamma_{i-1} & 0 & 0 & ... & 0 & 0 \\ p_{n1} + rc_{i1} & ... & p_{n.i-1} + rc_{i.i-1} & p_{ni} + r\gamma_i & p_{n.i+1} & ... & p_{n.n-1} & \gamma_n \end{matrix} \right\| = T_i$$

is equal to 1. This matrix has maximal rank. So by Theorem 2.21, there are $s_1, s_2, ..., s_{i-1}$ such that

$$\left\| s_1 \quad s_2 \quad ... \quad s_{i-1} \quad 1 \right\| T_i =$$
$$= \left\| l_1 \quad ... \quad l_{i-1} \quad p_{ni} + r\gamma_i \quad p_{n.i+1} \quad ... \quad p_{n.n-1} \quad \gamma_n \right\| \sim$$
$$\sim \left\| 1 \quad 0 \quad ... \quad 0 \right\|.$$

So the row

$$\left\| l_1 \quad ... \quad l_{i-1} \quad p_{ni} + r\gamma_i \quad p_{n.n+1} \quad ... \quad p_{n.n-1} \quad \gamma_n \right\|$$

is primitive. It follows that the last row of the matrix

$$
\begin{Vmatrix}
1 & & 0 & 0 & ... & 0 & 0 \\
 & \ddots & & & & & \\
0 & & 1 & 0 & ... & 0 & 0 \\
0 & ... & 0 & 1 & & 0 & 0 \\
 & & & & \ddots & & \\
0 & ... & 0 & 0 & & 1 & 0 \\
s_1 & ... & s_{i-1} & 0 & ... & 0 & 1
\end{Vmatrix} V_{i-1} \Phi_n P =
$$

$$
=
\begin{Vmatrix}
\gamma_1 & & 0 & 0 & 0 & & 0 & 0 \\
 & \ddots & & & & & & \\
c_{i-1.1} & & \gamma_{i-1} & 0 & 0 & & 0 & 0 \\
c_{i1} & ... & c_{i.i-1} & \gamma_i & 0 & & 0 & 0 \\
c_{i+1.1} & ... & c_{i+1.i-1} & c_{i+1.i} & \gamma_{i+1} & & 0 & 0 \\
 & & & & & \ddots & & \\
c_{n-1.1} & ... & c_{n-1.i-1} & c_{n-1.i} & c_{n-1.i+1} & & \gamma_{n-1} & 0 \\
l_1 & ... & l_{i-1} & p_{ni} + r\gamma_i & p_{n.i+1} & ... & p_{n.n-1} & \gamma_n
\end{Vmatrix}
\qquad (7.12)
$$

is primitive. Denote by $\left\| t_1 \quad t_2 \quad ... \quad t_n \right\|$ the last row of the matrix

$$
\begin{Vmatrix}
1 & & 0 & 0 & ... & 0 & 0 \\
 & \ddots & & & & & \\
0 & & 1 & 0 & ... & 0 & 0 \\
0 & ... & 0 & 1 & & 0 & 0 \\
 & & & & \ddots & & \\
0 & ... & 0 & 0 & & 1 & 0 \\
s_1 & ... & s_{i-1} & 0 & ... & 0 & 1
\end{Vmatrix} V_{i-1}.
$$

From equality (7.12) it follows that

$$
\left\| t_1 \quad t_2 \quad ... \quad t_n \right\| \Phi_n P_i = \left\| p_{ni} + r\gamma_i \quad p_{n.i+1} \quad ... \quad p_{n.n-1} \quad \gamma_n \right\|.
$$

Therefore, $\left\| t_1 \quad t_2 \quad ... \quad t_n \right\|$ is the solution of equation (7.8). In addition, since

$$
\left\| t_1 \quad t_2 \quad ... \quad t_n \right\| \Phi_n P =
$$
$$
= \left\| l_1 \quad ... \quad l_{i-1} \quad p_{ni} + r\gamma_i \quad p_{n.i+1} \quad ... \quad p_{n.n-1} \quad \gamma_n \right\|,
$$

then

$$
\left\| t_1 \quad t_2 \quad ... \quad t_n \right\| \Phi_n =
$$
$$
= \left\| l_1 \quad ... \quad l_{i-1} \quad p_{ni} + r\gamma_i \quad p_{n.i+1} \quad ... \quad p_{n.n-1} \quad \gamma_n \right\| P^{-1}.
$$

The row

$$
\left\| l_1 \quad ... \quad l_{i-1} \quad p_{ni} + r\gamma_i \quad p_{n.i+1} \quad ... \quad p_{n.n-1} \quad \gamma_n \right\| P^{-1}
$$

**253**

is primitive. Consequently, $\left\| t_1 \quad t_2 \quad ... \quad t_n \right\| \Phi_n$ be a primitive row. From Theorem 1.1, (see p. 24) it follows that the group $\mathbf{G}_\Phi$ contains a matrix of the form

$$\left\| \begin{matrix} & & * & \\ \varphi t_1 & ... & \varphi t_{n-1} & t_n \end{matrix} \right\|.$$

That is, if $i > 1$ the theorem is proved.

Let $i = 1$. We show that the equation

$$x\Phi_n P = \left\| p_{n1} + r\gamma_1 \quad p_{n2} \quad ... \quad p_{n.n-1} \quad \gamma_n \right\|, \tag{7.13}$$

where

$$(p_{n1} + r\gamma_1, p_{n2}, ..., p_{n.n-1}, \gamma_n) = (p_{n1}, p_{n2}, ..., p_{n.n-1}, \gamma_n) = 1 \tag{7.14}$$

has the desired solution.

In view of (7.14), by Theorem 7.1, equation (7.13) is solvable. So there is a row $\left\| \overline{x}_1 \quad ... \quad \overline{x}_n \right\|$ such that

$$\left\| \overline{x}_1 \quad ... \quad \overline{x}_n \right\| \Phi_n P = \left\| p_{n1} + r\gamma_1 \quad p_{n2} \quad ... \quad p_{n.n-1} \quad \gamma_n \right\|.$$

It follows that

$$\left\| \overline{x}_1 \quad ... \quad \overline{x}_n \right\| \Phi_n = \left\| p_{n1} + r\gamma_1 \quad p_{n2} \quad ... \quad p_{n.n-1} \quad \gamma_n \right\| P^{-1}.$$

The right side of this equality is a primitive row. Consequently, $\left\| \overline{x}_1 \quad ... \quad \overline{x}_n \right\| \Phi_n$ is a primitive row. Therefore, as above, the group $\mathbf{G}_\Phi$ contains a desired matrix.

Now we show that conditions 1), 2) are equivalent to conditions 3), 4). Suppose that the conditions 1) and 2) are fulfilled.

Reasoning similarly as above, we construct the matrix $H \in \mathbf{G}_\Phi$ such that

$$HP = \left\| \begin{matrix} & & & * & & & \\ * & ... & * & a_i & p_{n.i+1} & ... & p_{n.n-1} & \gamma_n \end{matrix} \right\|.$$

By virtue of Theorem 6.1,

$$(a_i, \varphi) = (p_{ni}, \varphi), \quad i = 1, ..., n-1.$$

Set
$$a_i = (p_{ni}, \varphi)a_i', \quad p_{ni} = (p_{ni}, \varphi)p_{ni}'.$$

Write condition 1) in the form

$$(p_{ni}, \varphi)a' \equiv (p_{ni}, \varphi)p_{ni}' \left( \mathrm{mod}\, \frac{\varphi(p_{ni}, ..., p_{n.n-1}, \gamma_n)}{(p_{n.i+1}, ..., p_{n.n-1}, \gamma_n)} \right).$$

Hence,

$$a' \equiv p_{ni}' \left( \mathrm{mod}\, \frac{\varphi(p_{ni}, ..., p_{n.n-1}, \gamma_n)}{(p_{ni}, \varphi)(p_{n.i+1}, ..., p_{n.n-1}, \gamma_n)} \right).$$

Since $\gamma_n | \varphi$, then

$$((p_{ni}, \varphi), (p_{n.i+1}, ..., p_{n.n-1}, \gamma_n)) =$$
$$= (p_{ni}, p_{n.i+1}, ..., p_{n.n-1}, (\gamma_n, \varphi)) =$$
$$= (p_{ni}, p_{n.i+1}, ..., p_{n.n-1}, \gamma_n).$$

Note that

$$\frac{(p_{ni}, \varphi)(p_{n.i+1}, ..., p_{n.n-1}, \gamma_n)}{((p_{ni}, \varphi), (p_{n.i+1}..., p_{n.n-1}, \gamma_n))} = [(p_{ni}, \varphi), (p_{n.i+1}..., p_{n.n-1}, \gamma_n)],$$

we get

$$a_i' \equiv p_{ni}' \left( \mathrm{mod} \frac{\varphi}{[(p_{ni}, \varphi), (p_{n.i+1}, ..., p_{n.n-1}, \gamma_{nn})]} \right).$$

It is easy to check that the reversed reasoning will be correct. The theorem is proved. $\qquad\square$

**Corollary 7.3.** *Let* $P = \|p_{ij}\|_1^n$ *be an invertible matrix. The group* $\mathbf{G}_\Phi$ *contains a matrix* $H$ *such that*

$$HP = \left\| \begin{matrix} & * & \\ s_{n1} & ... & s_{n.n-1} & \gamma_n \end{matrix} \right\|,$$

*where* $\gamma_n | \varphi$ *and* $s_{ni} \in K(\varphi)$, $i = 1, ..., n-1$.

**Proof**. Since $\Phi = \mathrm{diag}(1, ..., 1, \varphi)$, then

$$S_\Phi(P) = \left\| \begin{matrix} 1 & 1 & ... & 1 \\ ... & ... & ... & ... \\ 1 & 1 & ... & 1 \\ \gamma_1 & \gamma_2 & ... & \gamma_n \end{matrix} \right\|,$$

moreover

$$R_\Phi(P) = (\|p_{1n}...p_{nn}\|^T) = \|1...1\gamma_n\|^T.$$

By Theorem 6.3, there is $H_n \in \mathbf{G}_\Phi$ such that

$$H_n P = \left\| \begin{matrix} & * & \\ v_{n1} & ... & v_{n.n-1} & \gamma_n \end{matrix} \right\|.$$

In addition, according to Theorem 6.1, , $(v_{ni}, \varphi) = (p_{ni}, \varphi)$, $i = 1, ..., n-1$. If $v_{n.n-1} = 0$, then consider the element $v_{n.n-2}$. Let $v_{n.n-1} \neq 0$ and

$$s_{n.n-1} \equiv v_{n.n-1} \, (\mathrm{mod} \, \varphi),$$

where $s_{n.n-1} \in K(\varphi)$. Consequently,

$$s_{n.n-1} = v_{n.n-1} + r\varphi = v_{n.n-1} + \left( r \frac{\varphi}{\gamma_{n-1}} \right) \gamma_{n-1}.$$

That is

$$s_{n.n-1} \equiv v_{n.n-1} \, (\mathrm{mod} \, \gamma_{n-1}).$$

Noting that $\gamma_n | \varphi$ we have

$$(s_{n.n-1}, \gamma_n) = (v_{n.n-1} + r\varphi, \gamma_n) = (v_{n.n-1}, \gamma_n).$$

By Theorem 7.6, there is $H_{n-1} \in \mathbf{G}_\Phi$ such that

$$H_n U = \left\| \begin{matrix} & & & * \\ v'_{n1} & ... & v'_{n.n-1} & s_{n.n-1} & \gamma_n \end{matrix} \right\|.$$

By a similar scheme we consistently replace the elements at positions $(n, n-2), (n, n-3), ..., (n1)$, for their representatives from $K(\varphi)$. The proof is complete. $\square$

**Lemma 7.7.** *Suppose that invertible $n \times n$ matrices $U$ and $V$ have the same last row. Then there is $H \in \mathbf{G}_\Phi$ such that $HV = U$.*

**Proof.** Since

$$U = \left\| \begin{matrix} u_{11} & ... & u_{1n} \\ ... & ... & ... \\ u_{n-1.1} & ... & u_{n-1.n} \\ u_{n1} & ... & u_{nn} \end{matrix} \right\|, \quad V = \left\| \begin{matrix} v_{11} & ... & v_{1n} \\ ... & ... & ... \\ v_{n-1.1} & ... & v_{n-1.n} \\ u_{n1} & ... & u_{nn} \end{matrix} \right\|,$$

then

$$UV^{-1} = \left\| \begin{matrix} * & ... & * & * \\ ... & ... & ... & ... \\ * & ... & * & * \\ 0 & ... & 0 & 1 \end{matrix} \right\| = H.$$

It is obvious that $H \in \mathbf{G}_\Phi$. $\square$

We introduce the following notations:

$D(\varphi)$ is the set of all unassociated divisors of $\varphi$ except $\varphi$ itself;

$M_\varphi(\mu) = \{a \in K(\varphi) | (a, \varphi) = \mu\}$;

$M'_\varphi(\mu)$ is a quotient of $M_\varphi(\mu)$ by $\mu$;

$M'_\varphi(\mu, \gamma)$ is the set of coset representatives of $M'_\varphi(\mu)$ modulo $\dfrac{\varphi}{[\mu, \gamma]}$, $\gamma \in$

$\in D(\varphi)$. Note that if $a \in M'_\varphi(\mu)$, then a representative of coset which contain $a$ is

$$a_1 = a + r\frac{\varphi}{[\mu, \gamma]} \in K(\varphi)$$

where

$$\left( a_1, \frac{\varphi}{\mu} \right) = 1.$$

Let $P = \|p_{ij}\|_1^n$ be an invertible matrix and $p_{nk}$ is the on the right first of its last row item that is not divisible by $\varphi$, $2 \leqslant k \leqslant n$. Denote by

$$\mu_i = (p_{ni}, \varphi), \quad i = 1, ..., n.$$

That is, the last row of the matrix $P$ has the form

$$\|\mu_1 v_1 \ \dots \ \mu_k v_k \ \varphi v_{k+1} \ \dots \ \varphi v_n\|.$$

**Theorem 7.7.** *In the group* $\mathbf{G}_\Phi$ *there exists a matrix $H$ such that*

$$HP = \left\|\begin{matrix} \mathbf{0} & \bar{I}_{n-k} \\ L_k & \mathbf{0} \end{matrix}\right\|,$$

*where*

$$L_k = \left\|\begin{matrix} f_{k-1} & 0 & 0 & \dots & 0 & f_k \\ f_{k-2} & 0 & 0 & \dots & 1 & 0 \\ \dots & \dots & \dots & \dots & \dots & \dots \\ f_2 & 0 & 1 & \dots & 0 & 0 \\ f_1 & 1 & 0 & \dots & 0 & 0 \\ \mu_1 s_1 & \mu_2 s_2 & \mu_3 s_3 & \dots & \mu_{k-1} s_{k-1} & \mu_k \end{matrix}\right\|, \qquad (7.15)$$

*moreover if $\mu_i = \varphi$, then $s_i = 0$, if $\mu_i \in D(\varphi)$, then $s_i \in M'_\varphi(\mu_i, (\mu_{i+1}, ..., \mu_k))$, $i = 1, ..., k-1$, $f_j \in R$ , $j = 1, ..., k$, where $\bar{I}_{n-k}$ is a matrix of order an $n-k$ of the form*

$$\left\|\begin{matrix} 0 & & 1 \\ & \ddots & \\ 1 & & 0 \end{matrix}\right\|.$$

*In addition, the matrix $HP$ is a single matrix in the class $\mathbf{G}_\Phi P$ up to complement of its last row to an invertible matrix.*

**Proof.** If $k < n$, then $\mu_n = \varphi$. By Theorem 6.3, there is $H_n \in \mathbf{G}_\Phi$ such that

$$H_n P = \left\|\begin{matrix} p'_{11} & \dots & p'_{1.n-1} & p'_{1n} \\ \dots & \dots & \dots & \dots \\ p'_{n-1.1} & \dots & p'_{n-1.n-1} & p'_{n-1.n} \\ p'_{n1} & \dots & p'_{n.n-1} & \varphi_n \end{matrix}\right\|.$$

Let

$$(p'_{1n}, ..., p'_{n-1.n}) = \delta_n.$$

There is an invertible matrix $S_n$ such that

$$S_n \left\|\begin{matrix} p'_{1n} \\ p'_{2n} \\ \dots \\ p'_{n-1.n} \end{matrix}\right\| = \left\|\begin{matrix} \delta_n \\ 0 \\ \dots \\ 0 \end{matrix}\right\|.$$

Since

$$1 = (p'_{1n}, p'_{2n}, ..., p'_{n-1.n}, \varphi) = (\delta_n, \varphi),$$

there are $u_n, v_n$ such that

$$\delta_n u_n + \varphi v_n = 1.$$

Then

$$
\begin{Vmatrix}
u_n & 0 & ... & 0 & v_n \\
0 & 1 & ... & 0 & 0 \\
... & ... & \ddots & ... & ... \\
0 & 0 & ... & 1 & 0 \\
-\varphi & 0 & ... & 0 & \delta_n
\end{Vmatrix} (S_n \oplus 1) H_n P =
$$

$$
=
\begin{Vmatrix}
u_n & 0 & ... & 0 & v_n \\
0 & 1 & ... & 0 & 0 \\
... & ... & \ddots & ... & ... \\
0 & 0 & ... & 1 & 0 \\
-\varphi & 0 & ... & 0 & \delta_n
\end{Vmatrix}
\begin{Vmatrix}
* & ... & * & \delta_n \\
* & ... & * & 0 \\
... & ... & ... & ... \\
* & ... & * & 0 \\
* & ... & * & \varphi
\end{Vmatrix} =
$$

$$
=
\begin{Vmatrix}
* & ... & * & v_{1.n-1} & 1 \\
* & ... & * & v_{2.n-1} & 0 \\
* & ... & * & v_{n-1.n-1} & 0 \\
* & ... & * & v_{n.n-1} & 0
\end{Vmatrix} = P_n.
$$

Let $n - 1 > k$. The matrices

$$
\begin{Vmatrix}
u_n & 0 & ... & 0 & v_n \\
0 & 1 & ... & 0 & 0 \\
... & ... & \ddots & ... & ... \\
0 & 0 & ... & 1 & 0 \\
-\varphi & 0 & ... & 0 & \delta_n
\end{Vmatrix}, \ \ (S_n \oplus 1)
$$

belong to the group $\mathbf{G}_\Phi$. By virtue of Theorem 6.1, we receive $(v_{n.n-1}, \varphi) = \varphi$. That is $v_{n.n-1} = \varphi t_{n-1}$. Let $(v_{2.n-1}, ..., v_{n-1.n-1}) = \delta_{n-1}$. There is an invertible matrix $S_{n-1}$ such that

$$
S_{n-1}
\begin{Vmatrix}
v_{2.n-1} \\
v_{3.n-1} \\
... \\
v_{n-1.n-1}
\end{Vmatrix}
=
\begin{Vmatrix}
\delta_{n-1} \\
0 \\
... \\
0
\end{Vmatrix}.
$$

Then

$$
(1 \oplus S_{n-1} \oplus 1) P_n =
\begin{Vmatrix}
* & ... & * & v_{1.n-1} & 1 \\
* & ... & * & \delta_{n-1} & 0 \\
* & ... & * & 0 & 0 \\
... & ... & ... & ... & ... \\
* & ... & * & 0 & 0 \\
* & ... & * & \varphi t_{n-1} & 0
\end{Vmatrix}.
$$

The last two columns of this matrix is a primitive matrix. So $(\delta_{n-1}, \varphi t_{n-1}) = 1$. Thus, there are $u_{n-1}, v_{n-1}$ such that

$$
\delta_{n-1} u_{n-1} + \varphi t_{n-1} v_{n-1} = 1.
$$

Then

$$\begin{Vmatrix} 1 & 0 & 0 & ... & 0 & 0 \\ 0 & u_{n-1} & 0 & ... & 0 & v_{n-1} \\ 0 & 0 & 1 & ... & 0 & 0 \\ ... & ... & & \ddots & & ... \\ 0 & 0 & 0 & ... & 1 & 0 \\ 0 & -\varphi t_{n-1} & 0 & ... & 0 & \delta_{n-1} \end{Vmatrix} (1 \oplus S_{n-1} \oplus 1)P_n =$$

$$= \begin{Vmatrix} * & ... & * & v_{1.n-1} & 1 \\ * & ... & * & 1 & 0 \\ * & ... & * & 0 & 0 \\ * & ... & * & 0 & 0 \\ ... & ... & ... & ... & ... \\ * & ... & * & 0 & 0 \end{Vmatrix}.$$

Left-multiply this matrix by

$$\begin{Vmatrix} 1 & -v_{1.n-1} & 0 & ... & 0 \\ 0 & 1 & 0 & ... & 0 \\ 0 & 0 & 1 & & 0 \\ ... & ... & & \ddots & \\ 0 & 0 & 0 & & 1 \end{Vmatrix},$$

we get a matrix with two last columns:

$$\begin{Vmatrix} \bar{I}_2 \\ \mathbf{0} \end{Vmatrix}.$$

Continuing this process, we find a matrix $H'_k \in \mathbf{G}_\Phi$ such that

$$H'_k P = \begin{Vmatrix} & \begin{matrix} * \\ \vdots \\ * \\ d \\ 0 \\ \vdots \\ 0 \\ \mu_k v'_k \end{matrix} & \begin{matrix} \bar{I}_{n-k} \\ \hline \\ \mathbf{0} \end{matrix} \end{Vmatrix}.$$

This matrix is invertible. It follows that $(\mu_k v'_k, d) = 1$. Since $(\mu_k v'_k, \varphi) = \mu_k$, then

$$(\mu_k v'_k, \varphi d) = \mu_k.$$

So there are $a, b$ such that $\mu_k v'_k a + \varphi db = \mu_k$, moreover, according to Theorem 1.9, $(a, \varphi) = 1$. Since $(a, b) = 1$, then $(a, \varphi b) = 1$. Therefore, $af + \varphi bg = 1$ for some $f, g \in R$. Then

$$\left( I_{n-k} \oplus \left\| \begin{array}{ccccc} f & 0 & ... & 0 & -g \\ 0 & 1 & ... & 0 & 0 \\ ... & ... & \ddots & ... & ... \\ 0 & 0 & ... & 1 & 0 \\ \varphi b & 0 & ... & 0 & a \end{array} \right\| \right) H'_k P =$$

$$= \left\| \begin{array}{cccc|c} * & ... & * & * & \bar{I}_{n-k} \\ \hline * & ... & * & * & \\ c_{k1} & ... & c_{k.k-1} & \mu_k & \mathbf{0} \end{array} \right\| = \left\| \begin{array}{cc} B_k & \bar{I}_{n-k} \\ C_k & \mathbf{0} \end{array} \right\| = P_k.$$

The matrix $P_k$ is invertible. It follows that $C_k$ is invertible matrix. Then

$$\underbrace{\left\| \begin{array}{cc} I_{n-k} & -B_k C_k^{-1} \\ \mathbf{0} & I_k \end{array} \right\|}_{U} P_k = \left\| \begin{array}{cc} \mathbf{0} & \bar{I}_{n-k} \\ C_k & \mathbf{0} \end{array} \right\|,$$

where $U \in \mathbf{G}_\Phi$.

Let $F = \operatorname{diag}(1, ..., 1, \varphi)$ be a matrix of the order $k$. Let us show that by transformation from the group $\mathbf{G}_F$ the matrix $C_k$ can be reduced to the form (7.15).

According to Corollary 7.3, we can assume that $c_{ki} \in K(\varphi)$, $i = 1, ..., k-1$. Hence, if $\varphi | c_{k.k-1}$, then $c_{k.k-1} = 0$. Let $\varphi \nmid c_{k.k-1}$ and

$$c_{k.k-1} = \mu_{k-1} c'_{k-1} \in K_\varphi,$$

where $\mu_{k-1} \in D(\varphi)$. It follows that $c'_{k-1} \in M'_\varphi(\mu_{k-1})$. Let $s_{k-1} \in M'_\varphi(\mu_{k-1}, \mu_k)$, moreover

$$c'_{k-1} \equiv s_{k-1} \left( \operatorname{mod} \frac{\varphi}{[\mu_{k-1}, \mu_k]} \right).$$

By virtue of Theorem 7.6, and Corollary 7.3, there is $H_{k-1} \in \mathbf{G}_F$ such that

$$H_{k-1} C_k = \left\| \begin{array}{ccccc} & & & * & \\ * & ... & * & d_{k-2} & \mu_{k-1} s_{k-1} & \mu_k \end{array} \right\|,$$

where $d_{k-2} \in K(\varphi)$.

According to the described scheme, by transformation from $\mathbf{G}_F$, $C_k$ will be reduced to the form

$$C = \left\| \begin{array}{cccc} & & * & \\ \mu_1 s_1 & ... & \mu_{k-1} s_{k-1} & \mu_k \end{array} \right\|,$$

where

$$s_i \in M'_\varphi(\mu_i, (\mu_{i+1}, ..., \mu_k)), i = 1, ..., k-1.$$

By Theorem 1.2, the primitive last row of this matrix can be complemented to an invertible matrix of the form

$$L_k = \begin{Vmatrix} f_{k-1} & 0 & 0 & ... & 0 & f_k \\ f_{k-2} & 0 & 0 & ... & 1 & 0 \\ ... & ... & ... & ... & ... & ... \\ f_2 & 0 & 1 & ... & 0 & 0 \\ f_1 & 1 & 0 & ... & 0 & 0 \\ \mu_1 s_1 & \mu_2 s_2 & \mu_3 s_3 & ... & \mu_{k-1} s_{k-1} & \mu_k \end{Vmatrix}.$$

Based on Lemma 7.7, there is $D \in \mathbf{G}_F$ such that $DC = L_k$. Therefore, the group $\mathbf{G}_F$ contains a matrix $K$ such that $KC_k = L_k$. Then

$$\underbrace{\begin{Vmatrix} I_{n-k} & \mathbf{0} \\ \mathbf{0} & K \end{Vmatrix}}_{V} \begin{Vmatrix} \mathbf{0} & \bar{I}_{n-k} \\ C_k & \mathbf{0} \end{Vmatrix} = \begin{Vmatrix} \mathbf{0} & \bar{I}_{n-k} \\ L_k & \mathbf{0} \end{Vmatrix}.$$

where $V \in \mathbf{G}_\Phi$.

Suppose that $U \in \mathbf{G}_\Phi P$ and has the form

$$U = \begin{Vmatrix} \mathbf{0} & \bar{I}_{n-l} \\ N_l & \mathbf{0} \end{Vmatrix},$$

where

$$N_l = \begin{Vmatrix} f'_{l-1} & 0 & 0 & ... & 0 & f'_l \\ f'_{l-2} & 0 & 0 & ... & 1 & 0 \\ ... & ... & ... & ... & ... & ... \\ f'_2 & 0 & 1 & ... & 0 & 0 \\ f'_1 & 1 & 0 & ... & 0 & 0 \\ \mu'_1 s'_1 & \mu'_2 s'_2 & \mu'_3 s'_3 & ... & \mu'_{l-1} s'_{l-1} & \mu'_l \end{Vmatrix},$$

where $\mu'_i$ and $s'_i$ defined by the same rules as $\mu_i$ and $s_i$. The matrices $HP$ and $U$ are representatives of the coset $\mathbf{G}_\Phi P$. It follows that there is $T \in \mathbf{G}_\Phi$ such that $TU = HP$. Therefore, the $\Phi$-rods of corresponding columns of $U$ and $HP$ are coincide. Consequently,

$$\mu_i = \mu'_i, \quad i = 1, ..., k, \quad \mu_{k+1} = ... = \mu_n = \varphi.$$

So $l = k$. Since $TU = HP$ we have

$$T \begin{Vmatrix} & & * & & \\ \mu_{k-1} s'_{k-1} & \mu_k & 0 & ... & 0 \end{Vmatrix} = \begin{Vmatrix} & & * & & \\ \mu_{k-1} s_{k-1} & \mu_k & 0 & ... & 0 \end{Vmatrix}.$$

From the definition method of elements $s'_{k-1}, s_{k-1}$ and by Theorem 7.6, we get that $s'_{k-1} = s_{k-1}$. Again, considering the equality

$$T \begin{Vmatrix} & & & * & & \\ \mu_{k-2} s'_{k-2} & \mu_{k-1} s_{k-1} & \mu_k & 0 & ... & 0 \end{Vmatrix} =$$
$$= \begin{Vmatrix} & & & * & & \\ \mu_{k-2} s_{k-2} & \mu_{k-1} s_{k-1} & \mu_k & 0 & ... & 0 \end{Vmatrix},$$

we have $s'_{k-2} = s_{k-2}$. Reasoning similarly, we receive $s'_i = s_i$, $i = 1, ..., k$. Using Lemma 7.7, we complete the proof. $\square$

We fix the element $\gamma \in D(\varphi)$. Introduce the following notations:

$K_\varphi(\gamma) = \bigcup\limits_{\mu \in D(\varphi)} \mu M'_\varphi(\mu, \gamma) \bigcup \{0\}$,

$\bar{K}_\varphi(\gamma) = \{a \in K_\varphi(\gamma) | (a, \gamma) = 1\}$,

$N_k(\gamma)$ be the set of all rows of the form

$$\| a_1\ a_2\ ...\ a_{n-k-1}\ \gamma\ \underbrace{0\ ...\ 0}_{k} \|,$$

where

$$a_1 \in \bar{K}_\varphi(a_2, ..., a_{n-k-1}, \gamma),$$

$$a_i \in K_\varphi(a_{i+1}, ..., a_{n-k-1}, \gamma),\ \ i = 2, ..., n-k-1,\ \ 0 \leqslant k \leqslant n-2,$$

$$N(\varphi) = \bigcup\limits_{\gamma \in D(\varphi)} \bigcup\limits_{k=0}^{n-2} N_k(\gamma) \bigcup \{\| 1\ \ 0\ \ ...\ \ 0 \|\}.$$

Complement each row of $N(\varphi)$ to invertible matrices. The obtained set denote by $N_n(\varphi)$.

**Theorem 7.8.** *The set $N_n(\varphi)$ consists of representatives of all left cosets of a factor group* $\mathrm{GL}_n(R)$ *at* $\mathbf{G}_\Phi$.

**Proof** follows from Theorems 7.6 and 7.7. The more stringent restrictions, imposed on the choice of $a_1$, is conditioned by a primitivity of

$$\begin{Vmatrix} a_1 & a_2 & ... & a_{n-k-1} & \gamma & 0 & ... & 0 \end{Vmatrix}$$

$\square$

**Theorem 7.9.** *The set $N_n^{-1}(\varphi)\Phi$ consists of all right unassociated matrices with the Smith form $\Phi$.*

**Proof.** According to Theorem 4.5, the matrices $A$ and $B$ with the Smith form $\Phi = \mathrm{diag}\,(1, ..., 1, \varphi)$, are right associates if and only if $\mathbf{P}_A = \mathbf{P}_B$. The set $\mathbf{P}_A$ is a left coset of $\mathrm{GL}_n(R)$ at $G_\Phi$. So by Theorem 7.8, the assertion is true. $\square$

**Example 7.4.** Let $R = \mathbb{Z}$ and $\varphi = 30$. Then

$$D(\varphi) = \{1, 2, 3, 5, 6, 10, 15\},\ K(30) = \{0, 1, ..., 29\}.$$

We construct the set $K_{30}(3)$. First let us describe the sets $M_{30}(\mu)$, $M'_{30}(\mu)$, where $\mu \in D(\varphi)$:

$M_{30}(1) = \{1, 7, 11, 13, 17, 19, 23, 29\}$,
$M'_{30}(1) = M_{30}(1)$.

$M_{30}(2) = \{2, 4, 8, 14, 16, 22, 26, 28\}$,
$M'_{30}(2) = \{1, 2, 4, 7, 8, 11, 13, 14\}$.

$M_{30}(5) = \{5, 25\}$,
$M'_{30}(5) = \{1, 5\}$.

$M_{30}(6) = \{6, 12, 18, 24\}$,
$M'_{30}(6) = \{1, 2, 3, 4\}$.

$M_{30}(10) = \{10, 20\}$,
$M_{30}(10) = \{1, 2\}$.

$M_{30}(15) = \{15\}$,
$M_{30}(15) = \{1\}$.

Now we find the sets $M'_{30}(\mu, 3)$, where $\mu \in D(\varphi)$. To construct the set $M'_{30}(1, 3)$, it is necessary to brake down the set $M'_{30}(1)$ into the congruence class of a modulo

$$\frac{30}{[1, 3]} = 10.$$

Representatives are those elements that are relatively prime with

$$\frac{\varphi}{\mu} = \frac{30}{1} = 30.$$

Since

$$M_{30}(1) = \{\{1, 11\}, \{7, 17\}, \{13, 23\}, \{19, 29\}\},$$

then

$$M'_{30}(1, 3) = \{1, 7, 13, 19\}.$$

Despite the fact note

$$13 \equiv 3 (\mathrm{mod} 10),$$

the element 3 cannot be a representative of the class $\{13, 23\}$ because

$$(3, 30) \neq 1.$$

The simple calculations show that

$$M'_{30}(2, 3) = \{1, 2, 4, 8\}, \quad M'_{30}(3, 3) = \{1, 3, 7, 9\},$$

$$M'_{30}(5, 3) = \{1\}, \quad M'_{30}(6, 3) = \{1, 2, 3, 4\},$$

$$M'_{30}(10, 3) = \{1\}, \quad M'_{30}(15, 3) = \{1\}.$$

Hence,

$$K_{30}(3) = \{0, \ \{1, 7, 13, 19\}, \ 2\{1, 2, 4, 8\},$$

$$3\{1, 3, 7, 9\}, 5, \ 6\{1, 2, 3, 4\}, 10, 15\}. \qquad \diamond$$

We apply the obtained results to describe the left divisors of matrices. Let $A$ be an arbitrary $n \times n$ matrix with the Smith form

$$\mathrm{E} = \mathrm{diag}\,(\varepsilon_1, \dots, \varepsilon_t, 0, \dots, 0), \ \varepsilon_t \neq 0, \ t \leq n.$$

And let $\Phi = \mathrm{diag}\,(1, ..., 1, \varphi)|E$. Denote by $W(\mathrm{E}, \Phi)$ a subset of $N_{n \times n}(\varphi)$, which consists of matrices, the last row of which has the form

$$\left\| \frac{\varphi}{(\varphi, \varepsilon_1)} l_1 \quad ... \quad \frac{\varphi}{(\varphi, \varepsilon_t)} l_t \quad l_{t+1} \quad ... \quad l_n \right\|.$$

If the matrix

$$S = \left\| \frac{\varphi}{(\varphi, \varepsilon_1)} l_1 \quad ... \quad \overset{*}{\frac{\varphi}{(\varphi, \varepsilon_t)}} l_t \quad l_{t+1} \quad ... \quad l_n \right\|$$

is a representative of coset $\mathbf{G}_\Phi S$, then by Theorem 6.1, all matrices of this class has the form

$$\left\| \frac{\varphi}{(\varphi, \varepsilon_1)} v_1 \quad ... \quad \overset{*}{\frac{\varphi}{(\varphi, \varepsilon_t)}} v_t \quad v_{t+1} \quad ... \quad v_n \right\|.$$

**Theorem 7.10.** *If $\varphi$ be an irreducible element of R, then the set $\boldsymbol{W}(\mathrm{E}, \Phi)$ consists of matrices*

$$\left\| \begin{array}{ccccc|c} & & I_{n-1} & & & 0 \\ \hline 0 & ... & 0 \; a_s & ... & a_{n-1} & I_1 \end{array} \right\|,$$

$$\left\| \begin{array}{ccccc|c} & & I_{n-2} & & & 0 \\ \hline & & 0 & & & \\ 0 & ... & 0 \; a_s & ... & a_{n-2} & I_2 \end{array} \right\|, ...,$$

$$\left\| \begin{array}{cccc|c} & I_s & & & 0 \\ \hline & 0 & & & \\ 0 & ... & 0 & a_s & I_{n-s} \end{array} \right\|, \left\| \begin{array}{c|c} I_{s-1} & 0 \\ \hline 0 & I_{n-s+1} \end{array} \right\|,$$

*where $s$ is the number of the first invariant factor of the matrix $A$, which is divided into $\varphi$, $a_i \in K(\varphi)$, $i = s, s+1, ..., n-1$. Moreover, $(\boldsymbol{W}(\mathrm{E}, \Phi)P)^{-1}\Phi$ is the set of all irreducible divisors of the matrix $A = P^{-1}\mathrm{E}Q^{-1}$ with the Smith form $\Phi$.*

**Proof.** Since

$$\frac{\varphi}{(\varphi, \; \varepsilon_i)} = \varphi, \;\; i = 1, ..., s-1,$$

then the set $W(\mathrm{E}, \Phi)$ consists of matrices from $N_n(\varphi)$, which have the last row of the form

$$\left\| 0 \quad ... \quad 0 \quad l_s \quad l_{s+1} \quad ... \quad l_n \right\|.$$

Since $D(\varphi) = \{1\}$ it follows that

$$N(\varphi) = \bigcup_{k=0}^{n-1} N_k(1),$$

where $N_k(1)$ consists of all rows of the form

$$\|a_1\ a_2\ ...\ a_{n-k-1}\ 1\ \underbrace{0\ ...\ 0}_{k}\|,$$

where $a_i \in K(\varphi)$, $i = 1, ...n - k - 1$.

Selecting from the set $N(\varphi)$ rows with the desired structure and complementing them to the invertible matrices, we obtain the set $W(\mathrm{E}, \Phi)$.

To complete, it is enough to note that all irreducible matrices have the Smith form of the form diag $(1,...,1,\varphi)$, where $\varphi$ is an irreducible element in $R$ and using Theorem 4.6. □

## 7.3. Application of the obtained results to the solution of matrix unilateral equations

Example 5.6 (see p. 189) demonstrates that, in some cases, to find all solutions of the matrix equation it is sufficient to use the determining matrix. This will be if potential Smith forms of monic divisors either satisfy the requirements of Theorem 5.9, or there are no divisors with such forms. Otherwise the determining matrix generates only a part of the desired solutions. In this case, we use the notion of a standard $\Phi$-skeleton to extend the set of solutions. Let us illustrate this with specific examples.

Let $F$ be an algebraically closed field of characteristic zero. Consider the equations: $X^2 = \mathbf{0}$ and $X^2 = X$, where $X$ is $3 \times 3$ matrix. Their solutions are all third order nilpotent and idempotent matrices. Skeptics will say that solutions of this equations have been known long ago. Of course, they will be right. Indeed, nilpotent third order matrices are

$$\mathbf{0},\ \ T^{-1} \begin{Vmatrix} 0 & 0 & 0 \\ 0 & 0 & 0 \\ 0 & 1 & 0 \end{Vmatrix} T,$$

and idempotent are

$$\mathbf{0},\ \ I,\ \ T^{-1} \begin{Vmatrix} 0 & 0 & 0 \\ 0 & 0 & 0 \\ 0 & 0 & 1 \end{Vmatrix} T, T^{-1} \begin{Vmatrix} 1 & 0 & 0 \\ 0 & 1 & 0 \\ 0 & 0 & 0 \end{Vmatrix} T,$$

where $T$ is an arbitrary nonsingular matrix. However, such form of solutions has serious drawbacks. First, there are repeated ones. Second, much more seriously, the description of the roots in this form is reduced to the description of all third order invertible matrices. But this is another unsolved problem.

**Example 7.5.** Let us describe all solutions of the equation

$$X^2 = \mathbf{0}, \tag{7.16}$$

where $X$ be the third order matrix.

The matrix equation (7.16) corresponds to a polynomial matrix $A(x) =$ $= Ix^2$, which has the Smith form $E(x) = Ix^2$. It is obvious that $P_A = I$. Our task is to find all left unassociated divisors of the matrix polynomial $Ix^2$, having determinant degree 3 and separate among them monic. The potential Smith forms of these divisors are matrices

$$\Phi_1 = Ix, \quad \Phi_2 = \operatorname{diag}(1, x, x^2).$$

By Theorem 5.1, the matrix $\Phi_1$ corresponds to a single divisor $Ix = Ix - \mathbf{0}$. That is, the solution of the equation (7.16) is the matrix $X_1 = \mathbf{0}$.

According to Theorem 5.3, all matrices with the Smith form $\Phi_2$ are divisors of the matrix polynomial $A(x)$. Therefore, it is necessary to describe all right unassociated matrices with the Smith form $\Phi_2$. We classify them by $\Phi$-skeletons. The standard $\Phi_2$-skeletons are:

$$\begin{Vmatrix} 1 & 1 & 1 \\ x & 1 & 1 \\ x^2 & x & 1 \end{Vmatrix}, \quad \begin{Vmatrix} 1 & 1 & 1 \\ x & 1 & 1 \\ x^2 & 1 & x \end{Vmatrix}, \quad \begin{Vmatrix} 1 & 1 & 1 \\ 1 & x & 1 \\ x & x^2 & 1 \end{Vmatrix},$$

$$\begin{Vmatrix} 1 & 1 & 1 \\ 1 & 1 & x \\ x & 1 & x^2 \end{Vmatrix}, \quad \begin{Vmatrix} 1 & 1 & 1 \\ 1 & x & 1 \\ 1 & x^2 & x \end{Vmatrix}, \quad \begin{Vmatrix} 1 & 1 & 1 \\ 1 & 1 & x \\ 1 & x & x^2 \end{Vmatrix}.$$

The set of invertible matrices with standard $\Phi_2$-skeletons are:

$$\operatorname{Stand}(\Phi_2) =$$

$$= \left\{ \begin{Vmatrix} 1 & 0 & 0 \\ a & 1 & 0 \\ bx + c & d & 1 \end{Vmatrix}, \begin{Vmatrix} 1 & 0 & 0 \\ a & 0 & 1 \\ bx + c & 1 & 0 \end{Vmatrix}, \begin{Vmatrix} 0 & 1 & 0 \\ 1 & 0 & 0 \\ d & bx + c & 1 \end{Vmatrix} \right\} \bigcup$$

$$\bigcup \left\{ \begin{Vmatrix} 0 & 0 & 1 \\ 1 & 0 & 0 \\ d & 1 & 0 \end{Vmatrix}, \begin{Vmatrix} 0 & 1 & 0 \\ 0 & a & 1 \\ 1 & 0 & 0 \end{Vmatrix}, \begin{Vmatrix} 0 & 0 & 1 \\ 0 & 1 & 0 \\ 1 & 0 & 0 \end{Vmatrix} \right\},$$

where $a, b, c, d \in F$. Consequently, the set of all matrices with standard $\Phi_2$-skeletons are $\operatorname{Stand}^{-1}(\Phi_2)\Phi_2$. Let us choose from this set all right regularized matrices. The matrix $\Phi_2$ corresponds to

$$J(\Phi_2) = \begin{Vmatrix} 0 & 0 & 0 \\ 0 & 0 & 0 \\ 0 & 1 & 0 \end{Vmatrix}.$$

Then

$$M_{\left\|\begin{smallmatrix} 1 & 0 & 0 \\ a & 1 & 0 \\ bx+c & d & 1 \end{smallmatrix}\right\|}(\Phi_2) = \left\|\begin{matrix} a & 1 & 0 \\ c & d & 1 \\ b & 0 & 0 \end{matrix}\right\| = T_2.$$

Hence,

$$X = T_2^{-1} J(\Phi_2) T_2 = b^{-1} \left\|\begin{matrix} c & d & 1 \\ -ac & -ad & -a \\ (da-c)c & (da-c)d & da-c \end{matrix}\right\|.$$

So the set

$$\mathbf{X}_1 = \left\{ b^{-1} \left\|\begin{matrix} c & d & 1 \\ -ac & -ad & -a \\ (da-c)c & (da-c)d & da-c \end{matrix}\right\| \right\},$$

where $a, c, d \in F$, $b \in F \backslash \{0\}$, is the set of nilpotent matrices generated by the determining matrix $V(\mathrm{E}, \Phi_2)$.

Similarly, we obtain the following sets of solutions:

$$\mathbf{X}_2 = \left\{ b^{-1} \left\|\begin{matrix} -c & -1 & 0 \\ c^2 & c & 0 \\ ac & a & 0 \end{matrix}\right\| \right\}, \quad \mathbf{X}_3 = \left\{ b^{-1} \left\|\begin{matrix} 0 & 0 & 0 \\ -d & -c & -1 \\ cd & c^2 & c \end{matrix}\right\| \right\}.$$

Since

$$\det M_{\left\|\begin{smallmatrix} 0 & 0 & 1 \\ 1 & 0 & 0 \\ d & 1 & 0 \end{smallmatrix}\right\|}(\Phi_2) = \det M_{\left\|\begin{smallmatrix} 0 & 1 & 0 \\ 0 & a & 1 \\ 1 & 0 & 0 \end{smallmatrix}\right\|}(\Phi_2) =$$

$$= \det M_{\left\|\begin{smallmatrix} 0 & 0 & 1 \\ 0 & 1 & 0 \\ 1 & 0 & 0 \end{smallmatrix}\right\|}(\Phi_2) = 0,$$

then the matrices from the set

$$\left\{ \left\|\begin{matrix} 0 & 0 & 1 \\ 1 & 0 & 0 \\ d & 1 & 0 \end{matrix}\right\|, \left\|\begin{matrix} 0 & 1 & 0 \\ 0 & a & 1 \\ 1 & 0 & 0 \end{matrix}\right\|, \left\|\begin{matrix} 0 & 0 & 1 \\ 0 & 1 & 0 \\ 1 & 0 & 0 \end{matrix}\right\| \right\}^{-1} \Phi_2$$

are not right regularized and not give solutions.

The matrices with $\Phi$-skeletons

$$\left\|\begin{matrix} 1 & 1 & 1 \\ 1 & x & 1 \\ x & x & x \end{matrix}\right\|, \left\|\begin{matrix} 1 & 1 & 1 \\ 1 & 1 & x \\ x & x & x \end{matrix}\right\|$$

correspond to the next sets of solutions, respectively:

$$\mathbf{X}_4 = \left\{ \left\|\begin{matrix} 0 & 0 & 0 \\ m & 0 & 0 \\ n & 0 & 0 \end{matrix}\right\| \right\},$$

$m \in F\backslash\{0\}$, $n \in F$,

$$\mathbf{X}_5 = \left\{ \begin{Vmatrix} 0 & 0 & 0 \\ 0 & 0 & 0 \\ m & n & 0 \end{Vmatrix} \right\},$$

$m, n \in F$, where $m, n$ simultaneously do not equal to zero. $\diamondsuit$

**Example 7.6.** Let us describe all solutions of the equation

$$X^2 = X, \tag{7.17}$$

where $X$ be third order matrix.

The matrix equation (7.17) corresponds to a polynomial matrix $A(x) = Ix^2 - Ix$ with Smith form $\mathrm{E}(x) = I(x^2 - x)$. Hence, $P_A = I$.

Our task is to find all left unassociated divisors of the matrix polynomial $I(x^2 - x)$, having determinant degree 3 and separate among them monic. The potential Smith forms of these divisors are matrices

$$\Phi_1 = Ix, \quad \Phi_2 = I(x-1),$$

$$\Phi_3 = \mathrm{diag}\,(1, x, x(x-1)), \quad \Phi_4 = \mathrm{diag}(1, x-1, x(x-1)).$$

The matrices $\Phi_1$ and $\Phi_2$ correspond to single solutions

$$\mathbf{X}_1 = \mathbf{0}, \quad \mathbf{X}_2 = I.$$

The solutions with Smith form $\Phi_3$ of its characteristic matrices we classify by $\Phi_3$-skeletons:

$$\begin{Vmatrix} 1 & 1 & 1 \\ x & 1 & 1 \\ x(x-1) & x-1 & 1 \end{Vmatrix} \Rightarrow \mathbf{X}_3 = \left\{ b^{-1} \begin{Vmatrix} b+c & d & 1 \\ -(b+c)a & -da & -a \\ (b+c)(da-c) & d(da-c) & da-c \end{Vmatrix} \right\},$$

$$\begin{Vmatrix} 1 & 1 & 1 \\ x & 1 & 1 \\ x(x-1) & 1 & x-1 \end{Vmatrix} \Rightarrow \mathbf{X}_4 = \left\{ -b^{-1} \begin{Vmatrix} -(b+c) & -1 & 0 \\ (b+c)c & c & 0 \\ (b+c)a & a & 0 \end{Vmatrix} \right\},$$

$$\begin{Vmatrix} 1 & 1 & 1 \\ 1 & x & 1 \\ x-1 & x(x-1) & 1 \end{Vmatrix} \Rightarrow \mathbf{X}_5 = \left\{ -b^{-1} \begin{Vmatrix} 0 & 0 & 0 \\ -d & -b-c & -1 \\ dc & (b+c)c & c \end{Vmatrix} \right\},$$

where $a, c, d \in F$, $b \in F\backslash\{0\}$,

$$\begin{Vmatrix} 1 & 1 & 1 \\ 1 & 1 & x \\ x-1 & x-1 & x \end{Vmatrix} \Rightarrow \mathbf{X}_6 = \left\{ \begin{Vmatrix} 0 & 0 & 0 \\ 0 & 0 & 0 \\ m & n & 1 \end{Vmatrix} \right\},$$

$$\begin{Vmatrix} 1 & 1 & 1 \\ 1 & x & 1 \\ x-1 & x & x-1 \end{Vmatrix} \Rightarrow \mathbf{X}_7 = \left\{ \begin{Vmatrix} 0 & 0 & 0 \\ m & 1 & 0 \\ nm & n & 0 \end{Vmatrix} \right\},$$

$$\left\| \begin{matrix} 1 & 1 & 1 \\ x & 1 & 1 \\ x & x-1 & x-1 \end{matrix} \right\| \Rightarrow \mathbf{X}_8 = \left\{ \left\| \begin{matrix} 1 & 0 & 0 \\ m & 0 & 0 \\ n & 0 & 0 \end{matrix} \right\| \right\},$$

where $m, n \in F$.

Finally, we schematically write the solutions of equation (7.17), with the Smith form $\Phi_4$ of their characteristic matrices:

$$\left\| \begin{matrix} 1 & 1 & 1 \\ x-1 & 1 & 1 \\ x(x-1) & x & 1 \end{matrix} \right\| \Rightarrow \mathbf{X}_9 = \left\{ -b^{-1} \left\| \begin{matrix} -b-c & -d & -1 \\ ac & ad-b & a \\ (c+b-da)c & (c+b-da)d & c-da \end{matrix} \right\| \right\},$$

$$\left\| \begin{matrix} 1 & 1 & 1 \\ x-1 & 1 & 1 \\ x(x-1) & 1 & x \end{matrix} \right\| \Rightarrow \mathbf{X}_{10} = \left\{ b^{-1} \left\| \begin{matrix} b+c & 1 & 0 \\ -c(b+c) & -c & 0 \\ -ac & -a & b \end{matrix} \right\| \right\},$$

$$\left\| \begin{matrix} 1 & 1 & 1 \\ 1 & x-1 & 1 \\ x & x(x-1) & 1 \end{matrix} \right\| \Rightarrow \mathbf{X}_{11} = \left\{ b^{-1} \left\| \begin{matrix} b & 0 & 0 \\ d & b+c & 1 \\ -(b+c)d & -(b+c)c & -c \end{matrix} \right\| \right\},$$

where $a, c, d \in F$, $b \in F \backslash \{0\}$,

$$\left\| \begin{matrix} 1 & 1 & 1 \\ x-1 & 1 & 1 \\ x-1 & x & x \end{matrix} \right\| \Rightarrow \mathbf{X}_{12} = \left\{ \left\| \begin{matrix} 0 & 0 & 0 \\ m & 1 & 0 \\ n & 0 & 1 \end{matrix} \right\| \right\},$$

$$\left\| \begin{matrix} 1 & 1 & 1 \\ 1 & x-1 & 1 \\ x & x-1 & x \end{matrix} \right\| \Rightarrow \mathbf{X}_{13} = \left\{ \left\| \begin{matrix} 1 & 0 & 0 \\ m & 0 & 0 \\ -nm & n & 1 \end{matrix} \right\| \right\},$$

$$\left\| \begin{matrix} 1 & 1 & 1 \\ 1 & 1 & x-1 \\ x & x & x-1 \end{matrix} \right\| \Rightarrow \mathbf{X}_{14} = \left\{ \left\| \begin{matrix} 1 & 0 & 0 \\ 0 & 1 & 0 \\ m & n & 0 \end{matrix} \right\| \right\},$$

where $m, n \in F$. $\diamond$

# NOTATIONAL CONVENTIONS

| | |
|---|---|
| $R$ | — commutative finitely generated principal ideal ring without zero divisors (commutative Bezout domain) |
| $U(R)$ | — group of all units in the ring $R$ |
| $F$ | — field |
| $F[x]$ | — ring of polynomials with coefficients in $F$ |
| $\mathbb{C}$ | — field of complex numbers |
| $\mathbb{C}^n$ | — set of all $n \times 1$ matrices over $\mathbb{C}$ |
| $M_n(R)$ | — ring of all $n \times n$ matrices with entries from the ring $R$ |
| $M_{m \times n}(R)$ | — set of all $m \times n$ matrices with entries from the ring $R$ |
| $\mathrm{GL}_n(R)$ | — group of all invertible $n \times n$ matrices over the ring $R$ |
| $\mathbf{G}_\Phi$ | — sugroup of group $\mathrm{GL}_n(R)$ consists of all matrices $H$ such that $H\Phi = \Phi F$, where $F \in \mathrm{GL}_n(R)$ |
| $\mathbf{L}(\mathrm{E}, \Phi)$ | — set of all invertible matrix $H$ such that $H\mathrm{E} = \Phi S$, where $S \in M_n(R)$ |
| $\mathbf{W}(\mathrm{E}, \Phi)$ | — set of representatives of the left cosets of $\mathbf{G}_\Phi$ in $\mathbf{L}(\mathrm{E}, \Phi)$ |
| $a_i \mid a_j$ | — $a_i$ divides $a_j$ |
| $(a_1, a_2, ..., a_n)$ | — greatest common divisor of elements $a_1, a_2, ..., a_n$ |
| $[a_1, a_2, ..., a_n]$ | — least common multiple of elements $a_1, a_2, ..., a_n$ |
| $A^T$ | — transpose of the matrix $A$ |
| $\det A, |A|$ | — determinant of the matrix $A$ |
| $\langle A \rangle$ | — greatest common divisor of maximal order minors of matrix $A$ |
| $\langle A \rangle_i$ | — greatest common divisor of $i$th order minors of matrix $A$ |
| $I$ | — identity matrix |
| $I_n$ | — $n \times n$ identity matrix |

| | | |
|---|---|---|
| $\mathbf{0}$ | — | zero matrix |
| $\mathbf{0}_{m \times n}$ | — | $m \times n$ zero matrix |
| $\mathrm{diag}(a_1, a_2, ..., a_n)$ | — | diagonal matrix with elements $a_1, a_2, ..., a_n$ at main diagonal (may be rectangular) |
| $d$-matrix | — | $\mathrm{diag}(a_1, a_2, ..., a_n)$, where $a_i | a_{i+1}$, $i = 1, ..., n-1$ |
| $\mathrm{triang}(a_1, a_2, ..., a_n)$ | — | lower triangular matrix with elements $a_1, a_2, ..., a_n$ at main diagonal (square) |
| $A \oplus B$ | — | direct sum of matrices $A$ and $B$, i.e., $A \oplus B = \begin{Vmatrix} A & \mathbf{0} \\ \mathbf{0} & B \end{Vmatrix}$ |
| $\mathbf{P}_A$ | — | set of left transforming matrices of matrix $A$ to its Smith form |
| $K(f)$ | — | set of cosets representatives of factor ring $R/Rf$ |

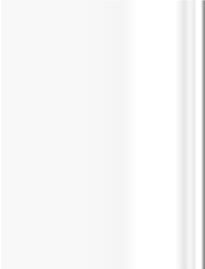

**274**

Монографія присвячена дослідженню арифметики кілець матриць над певними класами комутативних областей скінченно породжених головних ідеалів. Основну увагу зосереджено на побудові теорії розкладності матриць на множники. Висвітлюється тісний зв'язок факторизовності матриць із певними властивостями підгруп повної лінійної групи та спеціальною нормальною формою матриць стосовно односторонньої еквівалентності. Ґрунтовно вивчаються властивості матриць над кільцями стабільного рангу 1,5.

Для спеціалістів з теорії кілець, лінійної алгебри та студентів і аспірантів.

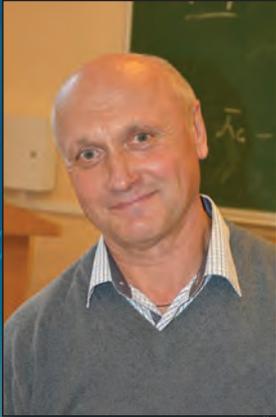

**Volodymyr P. Shchedryk**
Doctor of Physical and Mathematical Sciences,
Leading Researcher at the Pidstryhach Institute
of Applied Problems of Mechanics and Mathematics
of the National Academy of Sciences of Ukraine.

Volodymyr P. Shchedryk was born in Lypovets
of Vinnytsia region (Ukraine), in 1980 graduated
from the Vinnytsia State Pedagogical University,
and in 1988 defended his Ph.D. thesis at the Leningrad
State University (Russia).

In 2012 V. Shchedryk received a Doctor of Science
degree from the Institute of Mathematics
of the National Academy of Sciences of Ukraine.